\documentclass{amsart}
\usepackage{fullpage} 

\let\oldtocsection=\tocsection
\let\oldtocsubsection=\tocsubsection
\let\oldtocsubsubsection=\tocsubsubsection
 
\renewcommand{\tocsection}[3]{\hspace{0em}\oldtocsection{#1}{#2}{\bf#3}}
\renewcommand{\tocsubsection}[2]{\hspace{2.5em}\oldtocsubsection{#1}{#2}}
\renewcommand{\tocsubsubsection}[2]{\hspace{4em}\oldtocsubsubsection{#1}{#2}}

\usepackage{float}
\usepackage{marginnote}
\usepackage{amsmath, amscd, amssymb,amsfonts,amsthm,amscd,graphicx,psfrag,epsfig,mathrsfs}
\usepackage{youngtab}
\usepackage{graphicx}
\usepackage[all]{xy}
\usepackage{array}
\usepackage{amsthm}
\usepackage{stmaryrd}
\usepackage[titletoc,toc]{appendix}
\usepackage{stmaryrd}
\usepackage{mathtools}

\makeatletter
\def\subsection{\@startsection{subsection}{2}%
  \z@{.5\linespacing\@plus.7\linespacing}{.1\linespacing}%
  {\normalfont\bf}}
\makeatother

\usepackage[
colorlinks=true,%
linkcolor=blue,%
citecolor=blue,%
filecolor=blue,%
menucolor=blue,%
urlcolor=blue]{hyperref}

\usepackage{tikz}
\usepackage{tikz-cd}

\usetikzlibrary{arrows,chains,shapes,matrix,positioning,scopes} 
\usetikzlibrary{decorations.markings}

\tikzstyle arrowstyle=[scale=1.1]
\tikzstyle directed=[postaction={decorate,decoration={markings,
    mark=at position 1 with {\arrow[arrowstyle]{latex}}}}]
\tikzstyle reverse directed=[postaction={decorate,decoration={markings,
    mark=at position .45 with {\arrowreversed[arrowstyle]{latex};}}}]    
\tikzstyle arrowstyle=[scale=1.1]
\tikzstyle directed=[postaction={decorate,decoration={markings,
    mark=at position 1 with {\arrow[arrowstyle]{latex}}}}]
\tikzstyle reverse directed=[postaction={decorate,decoration={markings,
    mark=at position .45 with {\arrowreversed[arrowstyle]{latex};}}}]
    
\usetikzlibrary{decorations.pathreplacing,decorations.markings,decorations.pathreplacing,arrows,shapes}   
\tikzset{
  -z>/.style={
    decoration={
      show path construction,
      lineto code={
        \path (\tikzinputsegmentfirst) -- (\tikzinputsegmentlast) coordinate[pos=#1] (mid);
        \draw (\tikzinputsegmentfirst) -- (mid);
        \draw[double distance=1.5pt, arrows = {- Latex[length=0pt 2 0]}] (mid) -- (\tikzinputsegmentlast);      }
    },decorate
  }, -z>/.default=.5,
  z->/.style={
    decoration={
      show path construction,
      lineto code={
          \path (\tikzinputsegmentfirst) -- (\tikzinputsegmentlast) coordinate[pos=#1] (mid);
                \draw[double distance=1.5pt] (\tikzinputsegmentfirst) -- (mid);
                \draw[decoration={markings, mark=at position 1 with {\arrow[scale=1.2]{latex}}},postaction={decorate}] (mid) -- (\tikzinputsegmentlast);
      }
    },decorate
  }, z->/.default=.5,
  --z>/.style={
    decoration={
      show path construction,
      lineto code={
        \path (\tikzinputsegmentfirst) -- (\tikzinputsegmentlast) coordinate[pos=#1] (mid);
        \draw [dashed](\tikzinputsegmentfirst) -- (mid);
        \draw [double distance=1.5pt, dashed, arrows = {- Latex[length=0pt 2 0]}] (mid) -- (\tikzinputsegmentlast);      }
    },decorate
  }, --z>/.default=.5,
    z-->/.style={
    decoration={
      show path construction,
      lineto code={
          \path (\tikzinputsegmentfirst) -- (\tikzinputsegmentlast) coordinate[pos=#1] (mid);
                \draw[double distance=1.5pt, dashed] (\tikzinputsegmentfirst) -- (mid);
                \draw[dashed, decoration={markings, mark=at position 1 with {\arrow[scale=1.2]{latex}}},postaction={decorate}] (mid) -- (\tikzinputsegmentlast);
      }
    },decorate
  }, z-->/.default=.5,
   x->/.style={
    decoration={
      show path construction,
      lineto code={
          \path (\tikzinputsegmentfirst) -- (\tikzinputsegmentlast) coordinate[pos=#1] (mid);
                \draw[double distance=2pt] (\tikzinputsegmentfirst) -- (mid);
                \draw[decoration={markings, mark=at position 1 with {\arrow[scale=1.2]{latex}}},postaction={decorate}] (\tikzinputsegmentfirst) -- (\tikzinputsegmentlast);
      }
    },decorate
  }, x->/.default=.5,
  -x>/.style={
    decoration={
      show path construction,
      lineto code={
        \path (\tikzinputsegmentfirst) -- (\tikzinputsegmentlast) coordinate[pos=#1] (mid);
        \draw[double distance=2pt, shorten >=5pt] (mid) -- (\tikzinputsegmentlast); 
        \draw[decoration={markings, mark=at position 1 with {\arrow[scale=1.5]{latex}}},postaction={decorate}] (\tikzinputsegmentfirst) -- (\tikzinputsegmentlast);}
    },decorate
  }, -x>/.default=.5,
   x-->/.style={
    decoration={
      show path construction,
      lineto code={
          \path (\tikzinputsegmentfirst) -- (\tikzinputsegmentlast) coordinate[pos=#1] (mid);
                \draw[double distance=2pt, dashed] (\tikzinputsegmentfirst) -- (mid);
                \draw[dashed, decoration={markings, mark=at position 1 with {\arrow[scale=1.2]{latex}}},postaction={decorate}] (\tikzinputsegmentfirst) -- (\tikzinputsegmentlast);
      }
    },decorate
  }, x-->/.default=.5,
  --x>/.style={
    decoration={
      show path construction,
      lineto code={
        \path (\tikzinputsegmentfirst) -- (\tikzinputsegmentlast) coordinate[pos=#1] (mid);
        \draw[double distance=2pt, dashed, shorten >=5pt] (mid) -- (\tikzinputsegmentlast); 
        \draw[dashed, decoration={markings, mark=at position 1 with {\arrow[scale=1.4]{latex}}},postaction={decorate}] (\tikzinputsegmentfirst) -- (\tikzinputsegmentlast);}
    },decorate
  }, --x>/.default=.5,
}

\tikzset{
    partial ellipse/.style args={#1:#2:#3}{
        insert path={+ (#1:#3) arc (#1:#2:#3)}
    }
}

\newtheorem{theorem}{Theorem}[section]

\newtheorem{lemma}[theorem]{Lemma}
\newtheorem{proposition}[theorem]{Proposition}
\newtheorem{corollary}[theorem]{Corollary}

\newtheorem{conjecture}[theorem]{Conjecture}
\newtheorem{theorem-definition}[theorem]{Theorem-Definition}
\newtheorem{theorem-construction}[theorem]{Theorem-Construction}
\newtheorem{lemma-definition}[theorem]{Lemma--Definition}
\newtheorem{lemma-construction}[theorem]{Lemma--construction}
\newtheorem{definition}[theorem]{Definition}

\theoremstyle{definition}

\newtheorem{remark}[theorem]{Remark}
\newtheorem{example}[theorem]{Example}

\newcommand{\cal}{\mathcal}
\newcommand{\bg}{\begin{equation}\begin{gathered}}
\newcommand{\eg}{\end{gathered}\end{equation}}

\newcommand{\F}{\mathbb{F}}
\newcommand{\rI}{{\rm I}}
\newcommand{\rF}{{\rm F}}
\newcommand{\rE}{{\rm E}}
\newcommand{\g}{\mathfrak{g}}

\newcommand{\n}{\mathfrak{n}}

\newcommand{\p}{\mathfrak{p}}

\newcommand{\R}{\mathbb{R}}

\newcommand{\C}{\mathbb{C}}

\newcommand{\Z}{\mathbb{Z}}

\newcommand{\hlra}{\lhook\joinrel\longrightarrow}
\newcommand{\bF}{{\bf F}}
\newcommand{\bE}{{\bf E}}
\newcommand{\bK}{{\bf K}}
\newcommand{\X}{{\rm X}}
\newcommand{\old}[1]{}
\newcommand{\lr}{{\rm L}}

\newcommand{\D}{{\mathbb D}}
\newcommand{\Q}{{\mathbb Q}}
\newcommand{\T}{{\mathbb T}}

\newcommand{\A}{{\rm A}}
\newcommand{\B}{{\rm B}}
\newcommand{\G}{{\rm G}}
\newcommand{\U}{{\rm U}}
\renewcommand{\H}{{\rm H}}
\renewcommand{\P}{{\mathbb P}}

\newcommand{\bS}{{\mathbb{S}}}
\newcommand{\lms}{\longmapsto}
\newcommand{\lra}{\longrightarrow}
\newcommand{\hra}{\hookrightarrow}

\newcommand{\be}{\begin{equation}}
\newcommand{\ee}{\end{equation}}
\newcommand{\bt}{\begin{theorem}}
\newcommand{\et}{\end{theorem}}
\newcommand{\bd}{\begin{definition}}
\newcommand{\ed}{\end{definition}}
\newcommand{\bp}{\begin{proposition}}
\newcommand{\ep}{\end{proposition}}

\newcommand{\blc}{\begin{lemma-definition}}
\newcommand{\elc}{\end{lemma-definition}}
\newcommand{\btc}{\begin{theorem-construction}}
\newcommand{\etc}{\end{theorem-construction}}
\newcommand{\bl}{\begin{lemma}}
\newcommand{\el}{\end{lemma}}
\newcommand{\bc}{\begin{corollary}}
\newcommand{\ec}{\end{corollary}}
\newcommand{\bcon}{\begin{conjecture}}
\newcommand{\econ}{\end{conjecture}}
\newcommand{\la}{\label}

\input xy
\xyoption{all}

\begin{document}

\title{ Quantum geometry of moduli spaces of local systems  and representation theory}
\author{Alexander B. Goncharov, Linhui Shen}
\date{}
\maketitle


\begin{abstract} Let $\G$ be a split semi-simple adjoint group over $\Q$, and $\bS$ a {\it colored decorated surface}, given by an oriented  surface with punctures,  special boundary points, and a specified collection of boundary intervals. 

We introduce a  moduli space ${\mathscr P}_{\G, \bS}$ parametrizing $\G$-local systems on  $\bS$   with certain  boundary data, and prove that it carries a  cluster Poisson structure, equivariant under the action of  a discrete cluster modular  group  ${\Gamma}_{\G, \bS}$. 
The group ${\Gamma}_{\G, \bS}$ contains  
 the mapping class group  of $\bS$,   the product of  Weyl groups   over  punctures,   the product of generalized   braid groups   
over  boundary components, and the group of outer automorphisms of $\G$.  
 We prove  that  the dual   moduli   space ${\mathscr A}_{\G', \bS}$   
  carries a ${\Gamma}_{\G, \bS}$-equivariant cluster $K_2$-structure,  so that  
  the pair $({\mathscr A}_{\G', \bS}, {\mathscr P}_{\G, \bS})$ is a cluster ensemble. 
  This  generalizes the works  of V. Fock $\&$ the first author, and of I. Le. 
  \vskip 0.5mm
  
 The main theorem has many applications, including the following three, established in the paper. 

 \begin{itemize}
  \item Combining   with the  quantization of cluster Poisson varieties,  
we get a ${\Gamma}_{\G, \bS}$-equivariant {\it cluster quantization} of  ${\mathscr P}_{\G, \bS}$. It is given by a $\ast$-algebra $\mathcal{A}_\hbar({\mathscr P}_{\G, \bS})$    
 and  its  ${\Gamma}_{\G, \bS}$-equivariant principal series  $\ast$-representation in a Hilbert space.  
The Langlands modular duality is  built in the construction.

  \item  When $\bS$ is   punctured disc with two special   points,  we get  a geometric construction of the principal series $\ast$--representations of  the  
Langlands modular double of the quantum group ${\cal U}_q(\frak{g})$.

\item   Using the cluster Poisson structure of    ${\mathscr P}_{\G, \bS}$, we proved that its Donaldson-Thomas    transformation is  
a cluster transformation, and compute it explicitly. Combing this  with the  Gross-Hacking-Keel-Kontsevich construction, we  finally prove the
Duality Conjecture for ${\mathscr P}_{\G, \bS}$, including   a canonical linear basis in ${\cal O}({\mathscr P}_{\G, \bS})$. 
  \end{itemize}

 For  a genus $g$ surface $\bS$  with  $n>0$ punctures, we conjecture a canonical nondegenerate continuous pairing between  the  two infinite dimensional  vector bundles on the   space  $\widetilde {\cal M}_{g,n}$:
\vskip .5mm

\begin{enumerate}  \item  The bundle of  conformal blocks in the Liouville/Toda  Conformal Field Theory. It is a   bundle of {\it discrete} vector spaces,  given by  the coinvariants of   oscillatory   representations of 
the $W$-algebra of ${\mathfrak{g}^\vee}$.
 \vskip .5mm
 
\item   The bundle of {\it topological} vector spaces,    provided by the  cluster quantization of   ${\mathscr P}_{\G, \bS}$.  
\vskip .5mm
 The second bundle in the genus zero case is  identified by Modular Functor Conjecture with:
\vskip .5mm
\item  Tensor product invariants of  the principal series representations of the  quantum group ${\cal U}_q(\frak{g})$.
\end{enumerate}
\vskip 0.5mm

 The Modular Functor Conjecture is claimed for ${\rm PGL_2}$ by Teschner and for ${\rm PGL_m}$ by Schrader-Shapiro.  

\vskip .5mm

The conjectural  relation 1) $\leftrightarrow$ 3) between  representations of $W-$algebras and  principal series representations of quantum groups is a continuous analog 
of  the Kazhdan-Lusztig Theory.  \end{abstract}

\tableofcontents

\section{Summary}

\subsubsection{Main results.} Let $\G$ be a split semi-simple adjoint algebraic group over $\Q$. Let $\bS$ be a {\it colored decorated surface}, defined as an oriented   surface with  punctures and corners on the boundary, equipped  with a specified collection of boundary intervals,  called {\it colored boundary intervals}, see  Section \ref{2.1.1}.  

We introduce and quantize the  moduli space ${\mathscr P}_{\G, \bS}$ parametrizing $\G$-local system on  
  $\bS$ with some boundary data, see Definition \ref{DEFP2}. We achieve this  by proving our first main result,    Theorem \ref{MTH}: 
  \be \la{ONE}
\mbox{\it The moduli space    ${\mathscr P}_{\G, \bS}$ carries a canonical $\Gamma_{\G, \bS}-$equivariant cluster Poisson structure}. 
\ee
The  discrete symmetry group $\Gamma_{\G, \bS}$ is   
 defined in  Section \ref{2.1.6}. It contains the mapping class group of $\bS$, the product of  Weyl groups   over  punctures of $\bS$,  
 the product of generalized   braid groups   
over  boundary components of $\bS$, and the group   of outer automorphisms of $\G$. \vskip 1mm

When  $\bS$ varies,   moduli spaces  ${\mathscr P}_{\G, \bS}$ behave nicely under cutting and gluing 
of $\bS$, see \ref{1.2.3}-\ref{1.2.5}. Quantized moduli spaces provide an {\it algebraic-geometric avatar} of an extended TQFT\footnote{TQFT stands for topological quantum field theory},
 depending on a   number $\hbar\in \C$. \vskip 1mm
 
Applying to this non-linear and  non-commutative  object linearization procedures of various flavors, traditionally called  {\it quantizations},  and perceived as linear representations 
of the deformed  algebras of functions,  we get  more familiar  linear objects. They are  either traditional TQFT's or  their   continuous variants. 
In the continuous case,   simple objects are parametrized by spaces rather than   discrete sets.   \vskip 1mm

We stress  that the quantized moduli spaces  ${\mathscr P}_{\G, \bS}$ themselves, rather than their linear realizations, are the primary objects. 
We support this idea by the following three examples. \vskip 2mm

First, the  moduli spaces ${\mathscr P}_{\G, \bS}$ recover traditional quantum objects, such as quantum groups, in a simple and geometrically transparent way. 
The key features of   quantum groups, such as the Hopf algebra structure,   the   Chevalley generators,  the Cartan group, the R-matrix and the braided monoidal  structure on the category of representations, etc., are  built in the geometry of ${\mathscr P}_{\G, \bS}$. Thanks to  (\ref{ONE}),  the very  existence of these features no longer looks like a surprise and does not require   calculations. \vskip 1mm

  We stress that our approach leads to  genuine quantum objects, depending on the Planck constant   $\hbar \in \C$,  rather than   deformation quantization objects,  depending on the  formal parameter $\hbar$. \vskip 1mm
  
  Second,  
     quantum group  representations  of various flavors,  finite/infinite  dimensional, are naturally built into a much richer structure, given by   TQFT's of   algebraic  or analytic flavors. 
   Both  representations   
   and   morphisms between them   reveal themself as   objects of the same nature:  
   vector spaces    with  extra data, provided 
by   quantizations of   the    spaces   ${\mathscr P}_{\G, \bS}$. In particular, they are   modules over   non-commutative $\ast-$algebras ${\cal O}_q({\mathscr P}_{\G, \bS})$    
 quantizing   algebras   of  regular   functions on   ${\mathscr P}_{\G, \bS}$,   or  their Langlands modular doubles ${\cal A}_\hbar({\mathscr P}_{\G, \bS})$.   The operations of 
  tensor product  and its decomposition   into irreducibles 
    are incorporated  into the properties  of representations of these $\ast-$algebras    under the gluing and cutting of surfaces $\bS$.   We conclude that: \vskip 1mm
    
     {\it The  representation theory of quantum groups becomes a part of the representation theory of these algebras}. 
     
     \vskip 1mm
    This is a new phenomenon  in the representation theory, even for finite dimensional representations.   \vskip 1mm
    
    The resulting TQFT's are deeply related to Mathematical Physics. 
    In Section \ref{TQFT}   we discuss two examples, discrete and continuous, related to the Wess-Zumino-Witten and Liouville-Toda theories respectively.     
    \vskip 2mm

Third, let us turn on a complex structure,  starting from a punctured Riemann surface $\Sigma$. Consider 
the stack ${\rm Loc}_{\rm DR}(\G, \Sigma)$ of meromorphic $\G-$bundles with connections on $\Sigma$, with any singularities  
  at  the punctures. The Riemann-Hilbert correspondence identifies it with the Betti stack 
${\rm Loc}_{\rm B}(\G, S)$, which depends only on the  underlying topological surface   $S$, and  parametrises    
 monodromy of the connection   and  the Stokes data at  the punctures. 
Stacks ${\mathscr P}_{\G, \bS}$   are  generic  strata of ${\rm Loc}_{\rm B}(\G, S)$, and their slight generalization   describes all strata of the stack ${\rm Loc}_{\rm B}(\G, S)$. Stacks ${\rm Loc}_{\rm DR}(\G, \Sigma)$   carry 
a Poisson analog of  hyperk\"ahler structure. \\

Next, let  $\G'$ be the universal cover   of $\G$. Recall the moduli spaces ${\mathscr A}_{\G', \bS}$ introduced in \cite{FG03a},  see  Definition \ref{DEFP3}.  We prove in Theorem \ref{MTHa} our second main result: 
 \be \la{TWO}
 \begin{split}
 & \mbox{\it The moduli space  ${\mathscr A}_{\G', \bS}$ carries a canonical   $\Gamma_{\G, \bS}-$equivariant  cluster $K_2-$structure. }\\
 &\mbox {\it The pair  $({\mathscr A}_{\G', \bS}, {\mathscr P}_{\G, \bS})$ forms a $\Gamma_{\G, \bS}-$equivariant  cluster ensemble, in the sense of \cite{FG03b}.}\\
 \end{split}
 \ee  
 
 These results has numerous applications, including the following ones, which we briefly discuss below. 

  \subsubsection{Application I: Donaldson-Thomas transformations for the spaces ${\mathscr P}_{\G, \bS}$.} 
 The Donaldson-Thomas (DT) transformation is a remarkable transformation of a cluster Poisson variety ${\mathscr X}$, defined by Kontsevich and Soibelman \cite{KS}.  
 It encapsulates all information about   motivic Donaldson-Thomas invariants of the related Calabi-Yau categories, which  generalize the 
  classical numerical Donaldson-Thomas invariants. 
  However in general the DT-transformation is only a formal map, given by a collection of formal power series, and thus does not produce an actual transformation of  
  the space ${\mathscr X}$. \vskip 1mm
 
 We prove in Theorem \ref{CDT} that, barring a few of special cases,  
 \be \la{THREEDT}
 \begin{split}
 & \mbox{\it The DT-transformation of  the space ${\mathscr P}_{\G, \bS}$ is a cluster Poisson transformation.     }\\
 &\mbox{\it  We calculate it  explicitly as the geometric transformation ${\mathcal D}{\mathcal T}_{\G, \bS}$ from Definition \ref{KSDT}.}\\ 
 \end{split}
 \ee 
 This is a technically challenging result. For    $\G={\rm PGL}_m$ it was proved in \cite{GS16}, and required 40 pages to accomplish.  
 The proof for general $\G$ in Section \ref{DT=cluster}   is much more transparent, and takes just 2 pages plus  a few pictures. 
The dramatic simplification is achieved by using  our uniform treatment of the cluster Poisson 
 structure on the space ${\mathscr P}_{\G, \bS}$,  
 and  demonstrates clearly its greater flexibility and advantages even   for ${\rm PGL}_m$. 
 
\medskip

 \subsubsection{Intermezzo: Tropical points  in Higher Teichmuller theory.}  Before we proceed to the next major applications, let us recall the following. 
 We say that $\mathscr{X}$ is a positive space if it carries a collection of rational coordinate systems 
 such that the transition functions between any two of them are subtraction free rational functions. In this case the set of points 
 $\mathscr{X}({\cal F})$ with values in any semifield ${\cal F}$, e.g. $\R_{>0}$ and the semifields  $\Z^t, \Q^t, \R^t$ of integral, rational or real tropical numbers are well defined, 
 see \cite[Section 1]{FG03b}. Any cluster variety is a positive space. 
 So,  thanks to  (\ref{TWO}), we can define the Higher Teichmuller space ${\mathscr P}_{\G, \bS}(\R_{>0})$ and  the set of unbounded measured $\G-$laminations 
 ${\mathscr P}_{\G, \bS}(\R^t)$. Similarly, we can define the sets ${\mathscr A}_{\G', \bS}(\Z^t)$  and ${\mathscr A}_{\G', \bS}(\R^t)$  of integral and real tropical points of the space ${\mathscr A}_{\G', \bS}$. When $\G'={\rm SL}_2$, we recover the  integral and measured bounded laminations in the classical Teichmuller theory,  
 see \cite[Section 12]{FG03a}, \cite{FG05b}. The    group $\Gamma_{\G, \bS}$ acts on all these spaces.

   \medskip
   
 \subsubsection{Application II: Duality Conjectures and canonical linear basis in ${\cal O}({\mathscr P}_{\G, \bS})$.} \la{1.0.3}

 The result \eqref{THREEDT} is important not  only because it allows  to calculate motivic DT-invariants. 
 It also completes the proof of    Duality Conjectures \cite{FG03a} for the moduli space ${\mathscr P}_{\G, \bS}$. Namely,   denote by  ${\cal O}({\mathscr P}_{\G, \bS})$   the algebra of
  regular functions on the   moduli space ${\mathscr P}_{\G, \bS}$.  Let  $\G^\vee$ be the Langlands dual group. Theorem \ref{CBA} asserts that
 \be \la{FOURCB}
 \begin{split}
 & \mbox{\it There is a canonical linear basis in the vector space ${\cal O}({\mathscr P}_{\G, \bS})$, parametrised by the set 
${\mathscr A}_{\G^\vee, \bS}(\Z^t)$ }\\
 &\mbox{\it of the integral tropical points 
of the space ${\mathscr A}_{\G^\vee, \bS}$. This parametrization is $\Gamma_{\G, \bS}-$equivariant.  }\\ 
 \end{split}
 \ee  
We emphasize that   the canonical basis  in (\ref{FOURCB}) does not refer to cluster structures.  Indeed, the algebra ${\cal O}({\mathscr P}_{\G, \bS})$ is determined just by the stack ${\mathscr P}_{\G, \bS}$. And the positive structure on   ${\mathscr A}_{\G^\vee, \bS}$,   defined in \cite{FG03a},  is sufficient 
  to define its tropical points. However its proof   uses crucially the main results (\ref{ONE}) and (\ref{TWO}).   \vskip 1mm
  
  Here are the key steps of the proof. 

\begin{itemize}

\item   By  (\ref{ONE}) $\&$ (\ref{TWO}),   spaces ${\mathscr P}_{\G, \bS}$ and ${\mathscr A}_{\G^\vee, \bS}$ carry $\Gamma_{\G, \bS}-$equivariant Langlands dual cluster structures. 
\vskip 1mm

\item  Thanks to this,  we can apply the Gross-Hacking-Keel-Kontsevich construction  \cite{GHKK}. It provides for each cluster Poisson coordinate system  a 
  collection of   a priori  infinite formal Laurent series, parametrised by the 
set  ${\mathscr A}_{\G^\vee, \bS}(\Z^t)$ and related by cluster Poisson transformations. 
\vskip 1mm\item  The crucial application of  (\ref{THREEDT}) is that these Laurent series 
 are  in fact  Laurent polynomials in any cluster Poisson coordinate system on   ${\mathscr P}_{\G, \bS}$. The latter form an algebra, denoted by 
  ${\cal O}^{\rm cl}({\mathscr P}_{\G, \bS})$. Yet a priori these functions may not be regular on   the space ${\mathscr P}_{\G, \bS}$, i.e. may have singularities. 
\vskip 1mm
\item Finally, ${\cal O}({\mathscr P}_{\G, \bS}) = {\cal O}^{\rm cl}({\mathscr P}_{\G, \bS})$ by    \cite[Theorem 1.1]{S20}.

\end{itemize}

Linear bases in the space ${\cal O}({\mathscr P}_{\G, \bS})$, parametrised by the set 
${\mathscr A}_{\G^\vee, \bS}(\Z^t)$, but not $\Gamma_{\G, \bS}-$equivariant, were constructed in \cite{GS13}.  
\vskip 1mm

Our next goal is to take the full advantage of the fact that  elements of the canonical basis are actual functions, rather than just formal power series. In particular, 
we can evaluate them on the positive locus. This, plus the strong convexity properties of the canonical basis,  allow us to introduce several  kinds of limits, conjectured in \cite[Section 4]{FG03b} and established in Section \ref{SECT2.4a}, which all go under the roof of tropical limits.  They  have important applications in the Higher Teichmuller theory.  Let us discuss this. 

 \medskip

\subsubsection{Application III: Cluster convexity, length functions, and  Earthquake flow in Higher Teichmuller theory}   
  
 \bd \la{DEFCONV} Given a positive space $\mathscr{X}$, 
  a function $f$ on $\mathscr{X}(\R_{>0})$ or   $\mathscr{X}(\R^t)$ is  {\em convex}, if for any positive coordinate system $\{X_i\}$, the function $f(x)$  is a convex function in the logarithmic coordinates $\{x_i:=\log X_i\}$. 
 A function  on the set  $\mathscr{X}(\Z^t)$ or $\mathscr{X}(\R^t)$ 
  is convex if it is convex in any tropicalised positive coordinate system.\ed

We  define for any cluster variety the  {\it cluster convexity} as the convexity   for the cluster atlas. 
\vskip 2mm

   
   Theorem \ref{CBA}, summarised in (\ref{FOURCB}), 
can be restated as the existence a canonical 
$\Gamma_{\G, \bS}-$equivariant pairing
\be \la{CPa}
\mathbb{I}:  {\mathscr A}_{\G^\vee, \bS}(\Z^t)\times {\mathscr P}_{\G, \bS}(\C)\lra \C.
 \ee
Indeed, given a point $l\in  {\mathscr A}_{\G^\vee, \bS}(\Z^t)$, the function $\mathbb{I}(l, \cdot)$ on  $ {\mathscr P}_{\G, \bS}(\C)$ is   the 
regular function representing the canonical basis element parametrised by $l$. 
 Key properties of the canonical basis provide convexity properties, which imply the existence of   tropical limits of two flavors,  
 leading to the following  canonical    $\Gamma_{\G, \bS}-$equivariant maps, see Section \ref{SECT2.4a}:
  
\vskip 1mm
 (i) The  \underline{cluster length function}: 
  \be \la{SLim}
\begin{split}
&{\rm I}(l; e^x):  {\mathscr A}_{\G^\vee, \bS}(\R^t)\times {\mathscr P}_{\G, \bS}(\R_{>0})\lra \R. \\
\end{split}
\ee

\vskip 1mm
(ii) The \underline{cluster intersection pairing}:
 \be \la{SLime}
\begin{split}
&{\cal I}(l, m): {\mathscr A}_{\G^\vee, \bS}(\R^t)\times {\mathscr P}_{\G, \bS}(\R^t)\lra \R.\\
\end{split}
  \ee
We prove that both the cluster length function and the cluster intersection pairing  are 
   cluster convex in each of the two variables. 
  The cluster intersection pairing  is homogeneous of the degree $1$ in the first variable. 
  
  \vskip 1mm

Recall that the space ${\mathscr P}_{\G, \bS}$ is Poisson. So, given a tropical point $l \in  {\mathscr A}_{\G^\vee, \bS}(\R^t)$, we can view the cluster length function   ${\rm I}(l, e^x)$  as a Hamiltomian 
on the higher Teichm\"uller space ${\mathscr P}_{\G, \bS}(\R_{>0})$. 
Its   Hamiltonian flow  for the time $1$, depending on $l$, provides  the 
 analog of the Earthquake map:
 \be \la{EM}
{\cal E}_{\G, \bS}:  {\mathscr A}_{\G^\vee, \bS}(\R^t)\times  {\mathscr P}_{\G, \bS}(\R_{>0}) \lra  {\mathscr P}_{\G, \bS}(\R_{>0}). 
\ee
 
We conjecture that for any point $x \in  {\mathscr P}_{\G, \bS}(\R_{>0}) $, the Earthquake map ${\cal E}_{\G, \bS}$ provides an isomorphism
\be \la{EMI}
 {\cal E}^x_{\G, \bS}:  {\mathscr A}_{\G^\vee, \bS}(\R^t)\times  \{x\} \stackrel{\sim}{\lra}  {\mathscr P}_{\G, \bS}(\R_{>0}).  
 \ee

 For the group $\G={\rm PGL}_2$ most of these results  
are well known  in the classical Teichmuller theory. Namely, the cluster length function is the length of the geodesic representative of a measured lamination 
 in a given hyperbolic metric. The cluster intersection pairing is the intersection pairing between the bounded and unbounded laminations.  The 
 map (\ref{EM}) is Thurston's Earthquake map, and the map (\ref{EMI}) is an isomorphism by Thurston's Earthquake theorem. The cluster 
 convexity of the length function ${\rm I}(l; e^x)$ is an alternative to  Kerckhoff's theorem \cite{Ke}, see Section \ref{2.4.2a}. A conjectural analog of Kerckhoff's theorem is  Conjecture \ref{KCON}. 
 \medskip

\subsubsection{The main tool: cluster quantization.} \la{1.1.2} The key  application  of  cluster Poisson varieties  is the  {\it cluster Poisson quantization}.  
   Let us briefly recall 
its main features. \vskip 1mm

Given a cluster Poisson variety  ${\mathscr X}$ and a   number $q \in \C^*$,   there is the non-commutative  $q-$deformation  ${\cal O}_q({\mathscr X})$  of the   algebra of   functions on ${\mathscr X}$ \cite{FG03b}[Section 3]. It provides a functor from cluster Poisson varieties:
\be \la{qdeff}
{\mathscr X} \lra {\cal O}_q(\mathscr X).
\ee 
In particular, the   cluster modular group $\Gamma_{\mathscr X}$   acts  by automorphisms of ${\cal O}_q({\mathscr X})$. \vskip 1mm

Recall   the Langlands dual cluster Poisson variety  ${\mathscr X}^\vee$ \cite[Section 1.2.10]{FG03b}, see also Section \ref{2.4.6}. 

 Then, given  a    Planck constant $\hbar \in \C$, there is the {\it cluster Langlands  modular double}:
\be \la{LMD}
\mathcal{A}_\hbar({\mathscr X}):= {\cal O}_q({\mathscr X})\otimes_\C {\cal O}_{q^\vee}({\mathscr X}^\vee), \qquad  q = e^{i\pi \hbar}, \qquad q^\vee = e^{i\pi/ \hbar}.
\ee  
If $\hbar >0$ or $|\hbar |=1$, the algebra $\mathcal{A}_\hbar({\mathscr X})$ carries   a $\ast$-algebra structure, such that 
the   anti-involution   $\ast$ preserves the factors if $\hbar >0$, and swipes them if $|\hbar|=1$, see Sections \ref{1.1} and \ref{Sec8a}. \vskip 1mm

 The   cluster quantization  of  ${\mathscr X}$
provides  a $\ast-$representation of the $\ast-$algebra $\mathcal{A}_\hbar({\mathscr X})$    by unbounded operators in a Hilbert space 
${\cal H}_{\mathscr X}$, and a triple of spaces defined in \cite[Section 5]{FG07}, see also Section \ref{S17.2.10}:
 \be \la{SHS1}
{\cal S}_{\mathscr X} \subset {\cal H}_{\mathscr X} \subset {\cal S}^*_{\mathscr X}.
\ee
\vskip 2mm

Referring to (\ref{SHS1}) as  just vector spaces is a gross oversimplifying. In reality we define a collection, usually infinite, of vector spaces,  depending on a data of combinatorial 
nature, called a {\it quiver} ${\bf q}$, see Definition \ref{DQUIV}:
 \be \la{SHS1*}
{\cal S}_{{\mathscr X}, {\bf q}} \subset {\cal H}_{{\mathscr X}, {\bf q}}  \subset {\cal S}^*_{{\mathscr X}, {\bf q}} .
\ee
The algebra  $\mathcal{A}_\hbar({\mathscr X})$ acts by unbounded operators in each of the Hilbert spaces ${\cal H}_{{\mathscr X}, {\bf q}}$.
The subspace ${\cal S}_{{\mathscr X}, {\bf q}}$ is   the  common domain of the definition of operators from $\mathcal{A}_\hbar({\mathscr X})$, i.e.  the largest  subspace of  ${\cal H}_{{\mathscr X}, {\bf q}}$  
where the algebra $\mathcal{A}_\hbar({\mathscr X})$ acts -  the analog of the classical Schwarz space. We call it the {\it cluster Schwartz space}. 
 The action of  $\mathcal{A}_\hbar({\mathscr X})$ plus the Hilbert   metric  provide  each of the spaces  ${\cal S}_{{\mathscr X}, {\bf q}}$ with  a structure of a topological Fr\'echet space. 
 The space  ${\cal S}^*_{{\mathscr X}, {\bf q}}$ is the topological dual of ${\cal S}_{{\mathscr X}, {\bf q}}$. We call it  the space of {\it cluster distributions}.  \vskip 1mm
 
 The cluster spaces (\ref{SHS1*}) and the actions of the algebra $\mathcal{A}_\hbar({\mathscr X})$ are related as follows. For each pair of quivers related by so-called {\it cluster transformation} 
 ${\bf q} \to {\bf q}'$ we defined a unitary isomorphism 
 \be \la{INTERT}
 {\cal I}_{{\bf q} \to {\bf q}'}: {\cal H}_{{\mathscr X}, {\bf q}}  \stackrel{\sim}{\lra} {\cal H}_{{\mathscr X}, {\bf q'}} .
\ee
  The intertwiners provide isomorphisms between each of the three vector spaces in the triple (\ref{SHS1*}). So in addition to (\ref{INTERT}) we have isomorphisms 
  of topological vector space
 \be \la{INTERT1}
   {\cal I}_{{\bf q} \to {\bf q}'}: {\cal S}_{{\mathscr X}, {\bf q}} \stackrel{\sim}{\lra} {\cal S}_{{\mathscr X}, {\bf q'}},  \quad  {\cal I}_{{\bf q} \to {\bf q}'}: {\cal S}^*_{{\mathscr X}, {\bf q'}}  \stackrel{\sim}{\lra} {\cal S}^*_{{\mathscr X}, {\bf q}}.
\ee
  Furthermore, for each cluster transformation we have a canonical isomorphism  of $\ast-$algebras 
    \be \la{INTERT2}
 {\bf i}_{{\bf q} \to {\bf q}'}: {\cal A}_{{\mathscr X}}  \stackrel{\sim}{\lra} {\cal A}_{{\mathscr X}} .
\ee   
The intertwiners (\ref{INTERT1}) intertwine the algebra isomorphisms  ${\bf i}_{{\bf q} \to {\bf q}'}$ and the actions $\rho_{\bf q}$ and $\rho_{\bf q'}$ of the algebra ${\cal A}_{{\mathscr X}}$ on the 
topological vector spaces ${\cal S}_{{\mathscr X}, {\bf q}}$ and ${\cal S}_{{\mathscr X}, {\bf q'}}$: 
$$
\rho_{\bf q'}({\bf i}_{{\bf q} \to {\bf q}'}(a)) \cdot \ {\cal I}_{{\bf q} \to {\bf q}'}   (a)v=  {\cal I}_{{\bf q} \to {\bf q}'}  \ \rho_{\bf q}(a) \cdot v\ \ \ \ \forall a \in {\cal A}_{{\mathscr X}} \ \ \ \forall v \in {\cal S}_{{\mathscr X}, {\bf q}} .
$$

   The intertwiner  assigned to the composition ${\bf q} \to {\bf q}''$ of cluster transformation  ${\bf q} \to {\bf q}'$ and ${\bf q}' \to {\bf q}''$ is the composition of the corresponding intertwiners:
 $$
  {\cal I}_{{\bf q} \to {\bf q}'}: =   {\cal I}_{{\bf q} \to {\bf q}'}\circ {\cal I}_{{\bf q'} \to {\bf q}''}.
  $$

  There is a subclass of cluster transformations ${\bf q} \to {\bf q}$, called {\it trivial cluster transformations} - the ones which act trivially on the algebra $\mathcal{A}_\hbar({\mathscr X})$.  
  The corresponding intertwiners acts by multiplication by a unitary scalar. The cluster transformations considered modulo the trivial ones form the {\it cluster modular groupoid}.  
 Its objects are clusters, 
  and morphisms are cluster transformations modulo the trivial ones. The fundamental group of the cluster modular gropoid is called the {\it cluster modular group} $\Gamma_{\mathscr X}$.
  \vskip 2mm
  
To avoid abuse of terminology, we refer to the triple (\ref{SHS1})   as a triple of {cluster vector spaces}. 
The very notion of a {\it cluster vector space} assumes that we have a collection of vector spaces related by the intertwiners.    
 \vskip 2mm

 We conclude that  the   cluster modular group $\Gamma_{\mathscr X}$   acts by automorphisms   of the  $\ast-$algebra $\mathcal{A}_\hbar({\mathscr X})$. 
  The intertwiners provide  a unitary projective representation of the   group $\Gamma_{\mathscr X}$ on the cluster Hilbert space ${\cal H}_{\mathscr X}$, preserving the cluster   triple (\ref{SHS1}). 
So we get a 
 $\Gamma_{\mathscr X}-$equivariant $\ast-$representation of $\mathcal{A}_\hbar({\mathscr X})$,  called 
the   {\it principal series representation}, see Theorem \ref{MTHCQ}. It has been defined for $\hbar >0$ in  \cite{FG07}, and extended to the case $|\hbar|=1$ in  Section \ref{Sec8a}.


\medskip 

\subsubsection{Application III: quantization of  moduli spaces ${\mathscr P}_{\G, \bS}$.} \la{1.1.2}

The major application of  the result (\ref{ONE}) is the cluster Poisson quantization of  the  moduli space ${\mathscr P}_{\G, \bS}$.   
The algebra ${\cal O}({\mathscr P}_{\G, \bS})$ of regular functions on ${\mathscr P}_{\G, \bS}$  is isomorphic to the cluster Poisson algebra  ${\cal O}^{\rm cl}({\mathscr P}_{\G, \bS})$   \cite{S20}. 
So the algebra ${\cal O}_q({\mathscr P}_{\G, \bS})$, defined via (\ref{qdeff}) as the $q-$deformation of  the algebra ${\cal O}^{\rm cl}({\mathscr P}_{\G, \bS})$, is in fact a 
 $q-$deformation of the algebra   ${\cal O}({\mathscr P}_{\G, \bS})$. 
Let   $ \G_{\rm ad}^\vee$  be the adjoint group of $\G^\vee$. Following (\ref{LMD}), we introduce 
the {\it Langlands modular double} of the algebra  ${\cal O}_q({\mathscr P}_{\G, \bS})$   by setting 
 \be \la{AH}
 \mathcal{A}_\hbar({\mathscr P}_{\G, \bS}):= {\cal O}_q({\mathscr P}_{\G, \bS})\otimes_\C{\cal O}_{q^\vee}({\mathscr P}_{\G_{\rm ad}^\vee,  \bS} ), \qquad q= e^{i\pi \hbar}, \qquad q^\vee:= e^{i \pi / \hbar}.
\ee
We stress that the  Langlands modular duality  is built into   the construction:
$$
  \mathcal{A}_\hbar({\mathscr P}_{\G, \bS}) = {\mathcal A}_{\hbar^\vee}({\mathscr P}_{\G_{\rm ad}^\vee,  \bS} ), \qquad  \hbar^\vee:= 1/\hbar.
$$
 The   cluster Poisson quantization provides a $\Gamma_{\G, \bS}-$equivariant  $\ast-$representation of the $\ast-$algebra  $\mathcal{A}_\hbar({\mathscr P}_{\G, \bS})$ in the Hilbert space $\mathcal{H}_{\G, \bS}$. 
 
\medskip 

\subsubsection{Application IV: quantization of   generalized character varieties.}  In the case when the decorated surface $\bS$ is a surface $S$ with $n$ punctures but without boundary, 
the moduli space ${\mathscr P}_{\G, S}$ is closely related to the traditional  Betti moduli space ${\rm Loc}_{\G, S}$ of $\G-$local systems on $S$, also known as  the character variety. 
However unless $n=0$, these two moduli spaces are different. Their algebras of regular functions are related as follows: 
\be \la{ONEr}
\mathcal{O}({\rm Loc}_{\G, S}) = {\cal O}({\mathscr P}_{\G, S})^{W^n}. 
\ee 

The space  ${\rm Loc}_{\G, S}$ does not carry any cluster structure. However using the cluster nature of the action of the group $W^n$ on 
the moduli space ${\mathscr P}_{\G, S}$  (\ref{ONE}) and   relation (\ref{ONEr}), we can apply the cluster quantization machine. \vskip 1mm

Generalising this,  we introduce a new moduli space ${\rm Loc}_{\G, \bS}$ for any decorated surface $\bS$. 
One has  
$$
\mathcal{O}({\rm Loc}_{\G, \bS}) = {\cal O}({\mathscr P}_{\G, \bS})^{W^n}. 
$$ 
 Using crucially   (\ref{ONE}) and the cluster nature of the $W^n-$action   on ${\mathscr P}_{\G, \bS}$, 
we  quantize  ${\rm Loc}_{\G, \bS}$ and introduce its  Langlands modular double: 
 \be \la{AH+}
 \begin{split}
&\mathcal{O}_q({\rm Loc}_{\G, \bS}) := {\cal O}_q({\mathscr P}_{\G, \bS})^{W^n}; \\
& \mathcal{A}_\hbar({\rm Loc}_{\G, \bS}):= {\cal O}_q({\rm Loc}_{\G, \bS})\otimes_\C{\cal O}_{q^\vee}({\rm Loc}_{\G_{\rm ad}^\vee, \bS}).\\
\end{split}
 \ee

The cluster quantization of   the space  ${\mathscr P}_{\G, \bS}$ provides  a triple of  spaces, see (\ref{SHS1}), which are denoted  by   
 \be \la{SHS1a}
{\cal S}_{\G, \bS} \subset {\cal H}_{\G, \bS} \subset {\cal S}^*_{\G, \bS}.
\ee
The $\ast-$algebra $\mathcal{A}_\hbar({\rm Loc}_{\G, \bS})$ acts on the  cluster Schwarz   subspace ${\cal S}_{\G, \bS}$.  
The Weyl group $W^n$ acts   by unitary intertwiners on the Hilbert space ${\cal H}_{\G, \bS}$, preserving the cluster Schwarz  subspace ${\cal S}_{\G, \bS}$. 

\medskip
 
  \subsubsection{A hypothetical link to the quantum geometric Langlands.} \la{1.0.7}The cluster quantisation of the    space  ${\mathscr P}_{\G, \bS}$   can be viewed as a analytic incarnation  of the quantum geometric Langlands correspondence. 

Indeed, by Lemma \ref{6.9.03.11},  for any cluster Poisson variety ${\mathscr  X}$ 
 there is a canonical isomorphism of algebras
 $${\cal O}_q({\mathscr  X})^{\rm op} = {\cal O}_{q^{-1}}({\mathscr  X}).
 $$
For the moduli spaces related to surfaces this gives  
 $$
 {\cal O}_{q^\vee}({\mathscr P}_{\G_{\rm ad}^\vee, \bS})^{\rm op} =   {\cal O}_{(q^\vee)^{-1}}({\mathscr P}_{\G_{\rm ad}^\vee, \bS}).
  $$ 
 Therefore the  cluster Schwarz  space ${\cal S}_{\G, \bS}$  is a 
  $\Gamma_{\G, \bS}-$equivariant  bimodule
\be \la{bimod}
 {\cal S}_{\G, \bS}\ \ \in \ \  {\cal O}_{q}({\mathscr P}_{\G,  \bS})- \mbox{\rm mod} -   {\cal O}_{(q^\vee)^{-1}}({\mathscr P}_{\G_{\rm ad}^\vee, \bS}).
  \ee
 Passing to the $W^n-$coinvariants,   and using definition (\ref{AH+}), we arrive at the  $\Gamma_{\G, \bS}/W^n-$equivariant  bimodule  
 \be \la{bimodS}
 ({\cal S}_{\G, \bS})_{W^n}\ \ \in \ \   {\cal O}_{q}({\rm Loc}_{\G,  \bS}) - \mbox{\rm mod} -   {\cal O}_{(q^\vee)^{-1}}({\rm Loc}_{\G_{\rm ad}^\vee, \bS}).
  \ee

 It is tempting to think about it as of the kernel for the quantum analytic Langlands correspondence. \vskip 1mm
 
 It would be very interesting to relate this     to Gaitsgory's    quantum geometric Langlands  conjecture \cite{G}. It says, in the unramified set-up, 
 for a complex projective curve $\Sigma$, that 
   the category of ${\cal D}_\hbar({\rm Bun}_\G)-$modules is  equivalent to the category of ${\cal D}_{-\frac{1}{\hbar}}({\rm Bun}_{\G^\vee})-$modules. 
  Here ${\cal D}_\hbar({\rm Bun}_\G)$ refers to the differential operators on the $\hbar$'s power of a natural line bundle on ${\rm Bun}_\G$. 
  We stress that in both cases the equivalence interchanges 
  $$
  \hbar \leftrightarrow -\frac{1}{\hbar} \ \ \ \,\mbox{and} \ \ \ \G \leftrightarrow \G^\vee.
  $$
  To compare with  bimodule (\ref{bimodS}), let $\Sigma$ be now   a genus $g$ complex curve  with $n>0$ punctures.  
 Then the quantum geometric Langlands studies  a family    of  categories over the moduli space ${\cal M}_{g,n}$, related by the quantum geometric Langlands duality equivelences.

 Let $S$ be a topological surface of the same topology as $\Sigma$. We suggest that the $\Gamma_{S}-$equivariant\footnote{When $\bS=S$ is a surface with punctures, the group $\Gamma_{\G, \bS}$ does not have the braid group factors, and   $\Gamma_{\G, \bS} / W^n =  \Gamma_S \times {\rm Out}(\G)$, where $\Gamma_S$  is the mapping class group  of $S$.}  
  bimodule   (\ref{bimodS})   describes the analytic realization of the family of quantum Langlands duality equivalences   over  ${\cal M}_{g,n}$. 
  \vskip 2mm
  
  
   However there seems to be an apparent issue. The quantum  geometric Langlands equivalence is uniquely determined by the complex curve $\Sigma$:  
  moving the $\Sigma$ around a loop $\gamma$ in the moduli space ${\cal M}_{g,n}$ we come back to the same equivalence. On the other hand, 
  the action of the mapping class group $\Gamma_S$ on  bimodule (\ref{bimodS}) is the key part of cluster quantization: 
  the homotopy class of the loop  $\gamma$ acts on  bimodule (\ref{bimodS}). 
  
  This issue is resolved as follows, see also Section \ref{6.0.8} for a similar 
  discussion. Let ${\cal K}_{\G, \Sigma}$ stands for the kernel providing the quantum geometric Langlands correspondence for the  curve $\Sigma$. 
  We imagine   an equivalence,  depending on a combinatorial data on $S$ -   the quiver ${\bf q}$  which we use to construct the bimodule  
  ${\cal S}_{{\G, S}_{\bf q}}$, e.g. a triangulation of $S$ plus some data:
    $$
  \Phi_{\Sigma, {\bf q}}: \mbox{Analytic realization of} \ {\cal K}_{\G, \Sigma}  \lra  \mbox{bimodule} \ {{\cal S}_{\G, S}}_{\bf q}.
   $$
  Moving around the loop $\gamma$ we move the combinatorial data as well, and therefore get naturally an equivalence    $$
 \Phi_{\Sigma, \gamma({\bf q})}: \mbox{Analytic realization of} \ {\cal K}_{\G, \Sigma}  \lra  \mbox{bimodule} \ {{\cal S}_{\G, S}}_{\gamma({\bf q})}.
   $$
 The bimodules     ${{\cal S}_{\G, S}}_{\bf q}$ and    
  ${{\cal S}_{\G, S}}_{\gamma({\bf q})}$   are not identical, but rather  intertwined  by the   intertwiner ${\cal I}_\gamma$ acting on  the underlying vector space and by the 
  action of the element $\gamma$ on the algebra $ \mathcal{A}_\hbar({\rm Loc}_{\G, \bS})$.  
    
\subsubsection{Application V:  Poisson Lie groups via moduli spaces.}  \la{sect1.5}
Consider    colored decorated surfaces   on Fig. \ref{R7}.

  \begin{figure}[ht]
\begin{tikzpicture}[scale=1]
\draw[dashed, red, thick] (0,0.5) circle (.5cm);
\node[blue] at (0,0.5) {$\circ$};
\node at (0,1) {$\bullet$};
\node at (0,0) {$\bullet$};

\begin{scope}[xshift=-2cm]
\node at (4,0) {$\bullet$};
\node at (4,1) {$\bullet$};
\node at (5,0) {$\bullet$};
\node at (5,1) {$\bullet$};
\draw[thick] (4,0) -- (5,0);
\draw[thick] (4,1)--(5,1);
\draw[dashed, red, thick]  (4,1) -- (4,0);
\draw[dashed, red, thick]  (5,1) --(5,0);
\end{scope}

\begin{scope}[xshift=1cm]
\node at (4,0) {$\bullet$};
\node at (4,1) {$\bullet$};
\node at (5,0) {$\bullet$};
\node at (5,1) {$\bullet$};
\draw[thick] (4,0) -- (5,0);
\draw[thick] (4,1)--(5,1);
\draw[dashed, red, thick]  (4,1) -- (4,0);
\draw[dashed, red, thick]  (5,1) --(5,0);
\node[blue] at (4.5, 0.5) {$\circ$};
\end{scope}

\end{tikzpicture}
\caption{Three colored decorated surfaces, related to the Poisson Lie groups $\G^*$, $\G$ and $\D(\G)$. }
\label{R7}
\end{figure}

The related moduli spaces provide the three basic Poisson Lie groups related to the Lie algebra $\mathfrak{g}$. \vskip 1mm
\vskip 1mm

1. Let $\bS=\odot$ be a decorated surface given by a punctured disc with two special boundary points. 
The  Poisson moduli space ${\rm Loc}_{\G, \odot}$ has a natural map $$
\mu_{\rm out}: {\rm Loc}_{\G, \odot} \lra \H
$$ to the Cartan group $\H$ of $\G$, called the outer monodromy map. Denote by ${\mathscr L}_{\G, \odot}$ 
its fiber over the unit $e \in \H$. We prove that the moduli space ${\mathscr L}_{\G, \odot}$ has a natural Poisson Lie group structure, and identify it with the dual Poisson Lie group $\G^*$, see Section \ref{SEC1.7}:
\be \la{3ms1}
{\mathscr L}_{\G, \odot} = \G^*.
\ee 

\vskip 1mm

2. The Poisson Lie group $\G$ is canonically identified with the moduli space ${\rm Loc}_{\G, \square}$ assigned to a rectangle $\square$ with two opposite colored sides,   see Section \ref{SECT4.2}:
\be \la{3ms2}
{\rm Loc}_{\G, \square} = \G.
\ee 

\vskip 1mm

3. The Drinfeld double $\D(\G)$ of the Poisson Lie group  $\G$ is identified with  the moduli space ${\rm Loc}_{\G, \boxed{\cdot}}$ assigned to a punctured 
rectangle $\boxed{\cdot}$ with two opposite colored sides, see Section \ref{SECT4.2}:
\be \la{3ms3}
{\rm Loc}_{\G, \boxed{\cdot}} = \D(\G).
\ee 
 
More examples, including 
the cluster Poisson structure of  Grothendick's resolution $\widehat \G \lra \G$,  the dual Poisson Lie group $\D(\G)^*$,  and  the Heisenberg double $\D_{\rm H}(\G)$ of the Poisson Lie group $\G$,  see  in Section \ref{PLG}. 

\medskip

\subsubsection{Application VI: Quantum groups and their principal series of representations.}

Recall the $q-$deformation functor  $ 
{\mathscr X} \lra {\cal O}_q(\mathscr X) 
$  from (\ref{qdeff}). 
Applied to  the  cluster Poisson varieties  associated   by   theorem  (\ref{ONE}) with  the   moduli spaces  discussed Section \ref{PLG},  we arrive at   the corresponding 
quantum groups, understood as Hopf algebras. We stress that the Hopf algebra structure on the algebras ${\cal O}_q(\mathscr X)$ is immediate due to the following:

 \vskip 1mm {\it The composition law for the corresponding moduli spaces is given by  the gluing maps, sometimes followed by the cutting maps. These   are cluster Poisson maps, and thus 
 lead to maps of quantum algebras by (\ref{qdeff})}.  
 \vskip 1mm
 
The   cluster quantization given by the Application IV  delivers immediately the principal series $\ast-$representations of the Langlands modular doubles of the corresponding quantum groups, see Sections \ref{SECT3.4}, \ref{SEC1.8}  and   \ref{SSEECC11.3}. In particular, we get the principal series  $\ast-$representations for     the following quantum groups, where ${\rm D}(\mathfrak{g})$ is the Drinfeld double of the Lie bialgebra $\mathfrak{g}$: 
\be
\U_q(\mathfrak{g}), \quad {\cal O}_q(\G), \quad  {\U}_q({\rm D}(\mathfrak{g})), \quad {\cal O}_q({\rm D}(\mathfrak{g})).
\ee
Here $\U_q(\mathfrak{g})$ is the Hopf algebra quantizing the universal enveloping algebra of $\mathfrak{g}$, ${\cal O}_q(\G)$ is the dual Hopf algebra, quantizing regular functions on $\G$, and so on. 

 \medskip

 \subsubsection{Applications VII: quantization of ${\rm Loc}_{\G, \bS}$ and representation theory $\&$ Mathematical Physics.} \la{1.0.8}

Let $S$ be an oriented surface   with $n$ punctures. 
 The cluster quantization   of the space ${\rm Loc}_{\G, S}$ 
 provides   projective representations of the mapping class group $\Gamma_S$ of $S$  in the triple of spaces similar to  (\ref{SHS1}):
 \be  \la{Gtriple}
{\cal S}({\rm Loc}_{\G, S}) \subset {\cal H}({\rm Loc}_{\G, S}) \subset {\cal S}^*({\rm Loc}_{\G, S}).
\ee 
Conjecture \ref{GENDKL}   relates 
  quantized moduli spaces ${\rm Loc}_{\G, S}$  to the representation theory 
  of $W$-algebras and  
      Toda Conformal Field Theory for the group  $\G$. The Toda theory for $sl_2$  is the Liouville theory. \vskip 2mm

Precisely, denote by ${\rm CB}_{\G, \Sigma}$   the space of {\it conformal blocks} in the Toda    theory, defined as    the coinvariants of the  tensor product of  oscillatory   representations of 
the $W$-algebra of ${\mathfrak{g}^\vee}$,  assigned to the punctures of a variable Riemann surface $\Sigma$. It is an     infinite dimensional {\it discrete} vector space. 
Let  $S$  be the topological surface underlying $\Sigma$.  Conjecture \ref{GENDKL}  predicts a canonical nondegenerate continuous pairing 
\be \la{CBPar}
{\cal C}: {\rm CB}_{\G, \Sigma} \otimes {\cal S}({\rm Loc}_{\G, S})  \lra \C.
\ee
Since the pairing is non-degenerate, it tells   the size of the space of conformal blocks.\footnote{One can not have an isomorphism 
${\rm CB}_{\G, \Sigma} = {\cal H}({\rm Loc}_{\G, S})$ since these vector spaces are of    different nature: one is discrete, another is topological.} The continuous   pairing (\ref{CBPar}) is     the same thing as the canonical map
\be \la{CBPar2}
{\rm C}: {\rm CB}_{\G, \Sigma} \lra {\cal S}^*({\rm Loc}_{\G, S}).
\ee
The space ${\rm CB}_{\G, \Sigma}$ has a distinguished line ${\cal L}_{\G, S}$, generated by  the highest weight vector ${l}_{\G, \Sigma}$. Therefore we arrive at a family of one dimensional subspaces, depending on 
a Riemann surface $\Sigma$:  
\be \la{CBPar3}
{\rm C} ({\cal L}_{\G, \Sigma}) \subset  {\cal S}^*({\rm Loc}_{\G, S}).
\ee 
  The vector ${\rm C} ({l}_{\G, \Sigma})$ is the conformal block. Its regularised Hilbert norm 
  is conjectured to be the {\it correlation function}  for the Toda theory.   \vskip 2mm

 The space ${\cal H}({\rm Loc}_{\G, S})$  is  identified by the Modular Functor Conjecture with the  tensor product invariants of  the principal series representations of the  
 modular double of the quantum group ${\cal U}_q(\frak{g})$. \vskip 1mm

 The conjectural  relation between  representations of $W-$algebras and  principal series representations of quantum groups is a continuous analog 
of  the Kazhdan-Lusztig Theorem \cite{KL}.  Let us explain this. 

\medskip

 \subsubsection{Application VIII:   Finite dimensional  ${\mathscr U}_q(\mathfrak{g})-$modules  and WZW theory from ${\cal O}_q({\cal P}_{\G, \bS})-$modules.} 
 The category of  finite dimensional ${\mathscr U}_q(\mathfrak{g})-$modules    contains the subcategory of  class one modules. The latter is a  deformation of the category of finite dimensional  $\mathfrak{g}-$modules. 
 Any irreducible   ${\mathscr U}_q(\mathfrak{g})-$module  is a tensor product of a class one module and one of the $2^r$ one dimensional    
 ${\mathscr U}_q(\mathfrak{g})-$modules. 
To detect the subcategory of class one ${\mathscr U}_q(\mathfrak{g})-$modules, we consider  a Harish-Chandra-like pair $({\mathscr U}_q(\mathfrak{g}), {\rm H}, \tau)$, where $\rm H$ is the Cartan group of $\G$, and $\tau$ is its natural action by automorphisms  ${\mathscr U}_q(\mathfrak{g})$. Compatible representations of such a pair, see   Definition \ref{D5.10}, 
are excatly the  class one  ${\mathscr U}_q(\mathfrak{g})-$modules, see Section \ref{5.3.3}. \vskip 1mm

This construction has a   natural  generalizing. Namely, we   consider the triples 
$$
( {\cal O}_q({\mathscr P}_{\G, \bS}), {\rm H}^{c_\bS}, \tau_\bS)
$$
where $c_\bS$ is the number of colored boundary intervals on $\bS$, and    $\tau_\bS$ is an action of the group ${\rm H}^{c_\bS}$ by automorphisms of the algebra 
${\cal O}_q({\mathscr P}_{\G, \bS})$.  The ${\cal O}_q({\mathscr P}_{\G, \bS})-$modules  with a compatible action of the   group ${\rm H}^{c_\bS}$, see Definition \ref{def5.11},  are called   {\it enhanced  ${\cal O}_q({\mathscr P}_{\G, \bS})-$modules}. 

The action of the group ${\rm H}^{c_\bS}$ commutes with the $W^n-$action,  thus preserving the subalgebra ${\cal O}_q({\rm Loc}_{\G, \odot})$  in (\ref{AH}). 
This allows us to define the category of enhanced   ${\cal O}_q({\rm Loc}_{\G, \bS})-$modules. 
When $\bS=\odot$ is a punctured disc with two special points on Figure \ref{R7},  it contains the subcategory of enhanced  ${\cal O}_q({\mathscr L}_{\G, \odot})-$modules.

\bt   
For  any semi-simple Lie algebra $\mathfrak{g}$, Theorem \ref{UEAB} provides a monoidal functor 
\be\la{mcl}
\{\mbox{enhanced   ${\cal O}_q({\mathscr L}_{\G, \odot})-$modules}\} \lra \{\mbox{class one ${\mathscr U}_q(\mathfrak{g})-$modules}\}.
\ee
In the simply-laced case, it is  an equivalence thanks to \cite{S22}. Conjecturally, this is always the case. 
\et

The monoidal category  on the left in (\ref{mcl})  is naturally built into  the genus zero sector of an  
algebraic-geometric TQFT, given by enhanced ${\cal O}_q({\rm Loc}_{\G, \bS})-$modules $V_{\G, \bS; \lambda}$ of certain type. Let us elaborate on this. \\

Let $V_\lambda$ be the finite dimensional   ${\mathscr U}_q(\mathfrak{g})-$module with the highest weight $\lambda$. 
Recall  the multitplicity space   
\be \la{mspi}
\Bigl(V_{\lambda_1}\otimes \ldots \otimes V_{\lambda_n}\Bigr)^{{\mathscr U}_q(\mathfrak{g})}. 
\ee
 It has a natural action of the  braid group $\Gamma_{S^2_{(n)}}$  
  for a sphere with $n$ punctures $S^2_{(n)}$.  \vskip 2mm
 
  Recall that given a split torus $\T$, we denote by $X_*(\T)$ the group of cocharacters ${\rm Hom}({\Bbb G}_m, \T)$. 
      
       \bcon   \la{MFK2i}  $\bS $ be a sphere with $n$ punctures, and $q\in \C^\times$ is not a root of unity. Then 
       
   1.     The  
   category of enhanced finite dimensional    ${\cal O}_q({{\cal P}_{\G, \bS}})-$modules  is semi-simple. 

The 
irreducible modules  ${\cal V}_{\G, \bS, \alpha}$ 
are parametrised by collections of  integral coweights  at the punctures:
\be \la{insi}
\alpha   = (\alpha _1, \ldots , \alpha_n)  \in X_*({\rm H}^n).
\ee
The center  ${\cal O}({\rm H}^n)$   acts 
on  ${\cal V}_{\G, \bS, \lambda}$ by the character     
$\alpha(q) \in {\rm H}^n(\C)$.  
The   $W^n-$action on   ${\cal O}_q({{\cal P}_{\G, \bS}})$ provides  intertwiners 
\be
{\cal I}_w: {\cal V}_{\G, \bS, \alpha} \lra {\cal V}_{\G, \bS, w(\alpha)}, \qquad w \in W^n.
\ee 

2. Let $V_{\G, \bS, \alpha}$ be the ${\cal O}_q({{\rm Loc}_{\G, \bS}})-$module given by the restriction of   ${\cal V}_{\G, \bS, \alpha}$ 
to  the subalgebra ${\cal O}_q({{\rm Loc}_{\G, \bS}})$.  
 The simple  enhanced ${\cal O}_q({{\rm Loc}_{\G, \bS}})-$modules  
are precisely the ones $V_{\G, \bS, \lambda}$,    parametrised by     dominant coweights:
\be \la{ins1}
\lambda   = (\lambda _1, \ldots , \lambda_n)  \in X^+_*({\rm H}^n).
\ee

 3.  There is an isomorphism of $\Gamma_{\bS}-$modules 
\be \la{mspia}
V_{\G, \bS,   \lambda_1, ..., \lambda_n} = \Bigl(V_{\lambda_1}\otimes \ldots \otimes V_{\lambda_n}\Bigr)^{{\U}_q(\mathfrak{g})}. 
\ee 
 \econ

The Kazhdan-Lusztig theory \cite{KL} identifies  the  $\Gamma_{S^2_{(n)}}-$modules    (\ref{mspia})
 with the spaces of conformal blocks for the WZW theory, indexed by  (\ref{ins1}). So the following three finite dimensional $\Gamma_{S^2_{n}}-$modules   
 are isomorphic:
  \be \la{mspii}
(\mbox{WZW conformal blocks})_\lambda   =  \mbox{multiplicity spaces} ~ 
\Bigl(V_{\lambda_1}\otimes \ldots \otimes V_{\lambda_n}\Bigr)^{{\U}_q(\mathfrak{g})}=\mbox{${\cal O}_q({{\cal P}_{\G, \bS}})-$modules} ~V_{\G, \bS,   \lambda}. 
\ee 

  The new   feature here is that the vector spaces  
$V_{\G, \bS; \lambda}$ and hence the multiplicity spaces and  WZW conformal blocks   should carry an  action of the   algebra  ${\cal O}_q({\mathscr P}_{\G, \bS})$, equivariant under  the action of the  group $\Gamma_\bS$. \vskip 1mm

Let $S$ be a punctured surface   of  positive genus. Then there are no analogs of 
the multiplicity spaces in (\ref{mspii}).  The space of conformal blocks is   defined, carries an action of the mapping class group $\Gamma_S$, but, unless $\lambda =0$,   is   infinite dimensional. We conjecture that there are irreducible ${\cal O}_q({\mathscr P}_{\G, S})-$modules $V_{\G, S, \lambda}$, parametrised by    (\ref{ins1}), and isomorphic to the   WZW conformal blocks as $\Gamma_S-$modules, see Conjecture \ref{5.13}.  

\medskip

\subsubsection{The irregular case. } \la{5.3.5} 
 We expect the same is true in the most general set-up. Given a punctured   topological surface $S$, there is the de Rham stack ${\cal M}_{\rm DR}(S; \beta)$ 
 of framed   meromorphic  $\G-$bundles with connections on Riemann surfaces  of topological type $S$, with 
  irregular singularities at the punctures $p$ of  given  type $\beta = \{\beta_p\}$. Here $\beta_p$ is  an element of the cyclic envelope of the positive braid semigroup  of $\mathfrak{g}$. 
 Denote by ${\cal P}_{S, {\beta}}$ the moduli space  parametrising such Riemann surfaces   with  extra data describing irregular  singular points   \cite[Section 11.1]{GK}.   
 So there is a   map ${\cal P}_{S, {\beta}}\lra {\cal M}_{g,n}$ forgetting the extra data. 
  The de Rham stack ${\cal M}_{\rm DR}(S; \beta)$
  is fibered over   ${\cal P}_{S, {\beta}}$. Next, there is the constant Betti stack ${\cal M}_{\rm B}(S;\beta) \lra{\cal P}_{S, {\beta}}$ of framed Stokes data of   type $\beta$. The Riemann-Hilbert correspondence gives a complex analytic equivalence of  the two complex stacks over ${\cal P}_{S, {\beta}}$:
 \be \la{rs22xxys} 
 \begin{gathered}
     \xymatrix{
   {\cal M}_{\rm DR}(S, {\beta}) \ar[dr]   \ar[rr]^{ {\cal R}{\cal H} } & & \ar[dl]{\cal M}_{\rm B}(S, {\beta})  \\
         &    {\cal P}{\cal M}_{S, {\beta}}     &   }
 \end{gathered}
  \ee  
As explained in \cite[Section 11.4]{GK},   the main result (\ref{ONE})  implies that, barring the case  when $S$ is a disc and $\beta$ is rather degenerate, the Betti stack ${\cal M}_{\rm B}(S;\beta)$ has a cluster Poisson structure  equivariant under the action of the 
 generalized mapping class group 
 $$
 \Gamma_{S, \beta} := \pi_1({\cal P}{\cal M}_{S, {\beta}}  ).
 $$
  Therefore we get   $\Gamma_{S, \beta}-$equivariant quantized algebra ${\cal O}_q({\cal M}_{\rm B}(S; \beta))$. \vskip 2mm
  
  \bcon There are    spaces of $\beta-$irregular WZW conformal blocks. They are    $\Gamma_{S, \beta} -$modules,  isomorphic as $\Gamma_{S, \beta} -$modules to certain irreducible ${\cal O}_q({\cal M}_{\rm B}(C; \beta))-$modules 
 $V_{\G, S; \beta;  \lambda}$.
 \econ
 
When $\bS$ is a  colored decorated surface, each boundary component $\pi$ of $\bS$ describes an irregularity type $\beta_\pi$. The stack ${\cal P}_{\G, \bS}$ is a   stratum in the Betti stack ${\cal M}_{\rm B}(S; \beta)$, where $\beta = \{\beta_\pi\}$, and $S$ is a surface  obtained by shrinking   boundary components  of $\bS$ into   punctures  \cite[Section 11.5]{GK}. Denote by ${\rm H}^\circ$ the  subset of regular  elements in ${\rm H}$. The Weyl group $W$ acts freely on ${\rm H}^\circ$, and $\pi_1( {\rm H}^\circ(\C)/W)$ is the Artin 
 braid group ${\Bbb B}_{\mathfrak{g}}$ of $\mathfrak{g}$:  
 $$
\pi_1( {\rm H}^\circ(\C)/W) = {\Bbb B}_{\mathfrak{g}}. 
$$ 
Assume for simplicity  that each boundary component of $\bS$ carries an even number of special points. Denote by $n$ the number of boundary components on $\bS$. 
Then the moduli space ${\cal P}{\cal M}_{S, {\beta}} $ 
 is   the fibration 
 \be \la{fide}
({\rm H}^\circ(\C)/W)^n \lra {\cal P}{\cal M}_{S, {\beta}}  \lra {\cal M}_{g,n}.
   \ee 
 Therefore the  group ${\Bbb B}_{\mathfrak{g}}^n$  acts by automorphisms of the Betti stack ${\cal M}_{\rm B}(S, {\beta})$. This is  the braid group action we see in Section  \ref{2.1.6} and 
 Section   \ref{SECT9.4}. \vskip 1mm
 
 In particular, the stack ${\cal P}_{\G, \odot}$ is   a  stratum in the stack of Stokes data for the meromorphic $\G-$bundles with connections on ${\Bbb C}{\rm P}^1$ with a regular singularity at $0$ and an order two  irregular singularity at $\infty$, described by the  cyclic braid element $w_0^2$. The braid group ${\Bbb B}_{\mathfrak{g}}$ acts by its automorphisms. The   WZW conformal blocks   for the reduced stack ${\cal L}_{\G, \odot}$  are identified, as  ${\Bbb B}_{\mathfrak{g}}-$modules, with finite dimensional class one representations $V_\lambda$ of ${\mathscr U}_q(\mathfrak{g})$. The braid group action coincides with the one defined by Soibelman \cite{So1} under the name quantum Weyl group, and by Lusztig \cite{L},  see Section \ref{SSECC8}.\vskip 2mm
 
 For the similar reduced stack ${\cal L}_{\G, \odot_n}$, related to an $n-$punctured disc $\odot_n$ with two special boundary points,  see Section \ref{TQFT1},   the  
 de Rham stack parametrises $\G-$bundles with meromorphic connections on ${\Bbb C}\P^1$ with regular singularities at the punctures $z_1, ..., z_n$, and an order two pole at $\infty$ of the irregularity   type $w_0^2$. The   WZW conformal blocks   are identified with 
  tensor products $V_{\lambda_1} \otimes \ldots \otimes V_{\lambda_n}$ of  irreducible 
   finite dimensional class one  ${\mathscr U}_q(\mathfrak{g})-$modules. The algebra 
  ${\mathscr U}_q(\mathfrak{g})^{\otimes n}$ acts on them. 
  So in this case 
  \be \la{brN1}
  {\cal O}_q\Bigl({\cal M}_{\rm B}(S^2-\{z_1, ..., z_n, \infty\}; w_0^2)\Bigr) = {\U}_q(\mathfrak{g})^{\otimes n}.
    \ee
   In agreement with (\ref{fide}), the wild mapping class group  $\Gamma_{\C-\{z_1, ..., z_n\}; w_0^2} $  is the product of the  mapping class group $\Gamma_{\C-\{z_1, ..., z_n\}}$, which can interchange the punctures $z_1, ..., z_n$,   and  the  braid groups ${\mathbb B}_{\mathfrak{g}}$:
   \be
  \Gamma_{\C-\{z_1, ..., z_n\}; w_0^2} =   \Gamma_{\C-\{z_1, ..., z_n\}} \times {\mathbb B}_{\mathfrak{g}} .  
   \ee
The   group $\Gamma_{\C-\{z_1, ..., z_n\}; w_0^2}$ acts by automorphisms of each of the sides in (\ref{brN1}), and the group ${\mathbb B}_{\mathfrak{g}}$ acts diagonally. The isomorphism (\ref{brN1}) intertwines these two actions.

\medskip

 Section  \ref{1.0.8}   discusses the continuous analog of this story. It is elaborated in Sections \ref{TQFT1} - \ref{SEC1.8}. The variety of connections with the  representation theory and Mathematical Physics are discussed in Section \ref{sec1.9}. We expect all of them to hold in the most general irregular set-up, while for simplicity we mostly focus 
 on the case of regular singularities.

\subsubsection{\it Acknowledgments.} This paper   continues   series of works of the first author with Vladimir Fock, and contains some key ideas 
conceived more then a decade ago,   but not documented at that time. In particular, the idea that the Poisson Lie group $\G^*$ can be 
interpreted as an appropriately defined moduli space of $\G$-local systems on a punctured disc with two 
special points on the boundary is due to V.  Fock and the first author. The first author  thanks Vladimir Fock for the joy of collaboration. 
\vskip 2mm

The authors are very grateful to the Referee, who read the paper very carefully, and made many comments and remarks, incorporated in to the text, which greatly improved the exposition. 

\vskip 2mm
The work of the first author was
supported by the  NSF grants  DMS-1564385,  DMS-1900743, and  DMS-2153059. The work of the second author was supported by the NSF grant DMS-2200738. 
The final version of the paper was prepared when the  authors visited 
the MPI (Bonn) and IHES    (Bures sur Yvette) during the Summer of 2018, and finished when the first author stayed at  IHES in February of 2019, and 
at the Hebrew University of Jerusalem in March of 2019.  
The hospitality and support of MPI, IHES, and HUJI under the European
Research Council (ERC) grant N 669655,  are gratefully acknowledged. The first author was Simons Fellow  in 2019 and 2023. 
The generous support of the Simons Foundation is gratefully acknowledged. We are grateful to David Kazhdan for many fruitful conversations, and to Joerg Teschner for the comments 
on the first draft of the paper. 

\medskip

 \section{Introduction}  
 
\medskip
  \subsection{Cluster Poisson moduli spaces assigned to colored decorated surfaces}  \la{Sc1.2}

\subsubsection{Colored decorated surfaces.} \la{2.1.1}
A  {\it decorated surface} $\bS$ is an oriented topological surface  with   {\it punctures} inside, and 
a finite number of {\it special points} on the boundary, considered modulo isotopy, see Figure \ref{amalgam11}. 
We assume that each boundary component has at least one special point,  and 
    the total number   of punctures and special points   is positive. 
Consider   the punctured boundary
\be \la{PBDSS}
\widehat  \partial \bS:= \partial \bS - \{\mbox{\rm special points}\}.
\ee
 {\it Boundary intervals}  are  closures of  connected components of the punctured boundary. They are orientated by   the orientation of $\partial \bS$. 
For each special point $s$, there is a unique boundary interval ${\rm I}_s$ starting at $s$.  \vskip 1mm

A  {\it colored decorated surface} $\bS$ is a decorated surface with an additional data: a subset (possibly empty) 
$$
 \{\mbox{\it colored boundary intervals}\} \subset \{\mbox{\it boundary intervals}\}.
 $$
 \vskip 1mm

An {\it ideal triangulation} of $\bS$ is a triangulation with vertices at the special points and punctures. 

  \begin{figure}[ht]
\begin{center}
\begin{tikzpicture}[scale=.5]
 \draw(-1,-.5) arc (3.3:70:3);
  \draw (1,-.5) arc (176.7:110:3);
   \draw (1.5,2.5) arc (-20:-160:1.6 and .6);
   \draw (0.75,1) arc (-60:-120:1.5);
   \draw (0.6,.92) arc (60:120:1.2);
   \draw[dashed] (1, -.5) arc (0:180:1 and .4);
      \draw (1, -.5) arc (0:-180:1 and .4);
       \draw[red, thick, dashed] (-1.5, 2.5) arc (0:360:1 and .4);
       \draw[red, thick, dashed] (1.5, 2.5) arc (180:540:1 and .4);
       \node  at (0, -1.3) {\footnotesize $\B$};
       \node[red] at (0, -0.9) {\footnotesize $\bullet$};
       \node[blue] at (0, 0.3) {\footnotesize $\bullet$};
       \node  at (-0.3, .3) {\footnotesize $\B$};
        \node[blue] at (0.9, 1.6) {\footnotesize $\bullet$};
         \node  at (0.6, 1.6) {\footnotesize $\B$};
        \node[red] at (-3.5, 2.5) {\footnotesize $\bullet$};
         \node  at (-3.8, 2.5) {\footnotesize $\B$};
        \node[red] at (3.5, 2.5) {\footnotesize $\bullet$};
         \node  at (3.8, 2.5) {\footnotesize $\B$};
         \node[red] at (-1.9, 2.8) {\footnotesize $\bullet$};
          \node[red] at (1.9, 2.8) {\footnotesize $\bullet$};
          \node  at (-1.8, 3.1) {\footnotesize $\B$};
          \node  at (1.8, 3.1) {\footnotesize $\B$};
\end{tikzpicture}
\end{center}
\caption{A colored decorated surface $\bS$ with 2 punctures (blue),  5 special points (red), 4 colored boundary intervals (red dashed), and 1 non-colored boundary interval (black). 
The moduli space ${\mathscr P}_{\G, \bS}$  parametrizes $\G$-local systems  ${\cal L}$ on $\bS$ 
 with   invariant flags $\B$ near all punctures  and   special points, and pinnings over colored boundary segments.}
\label{amalgam11}
\end{figure}

Our construction of the cluster structures   of  various moduli spaces related to the moduli spaces of $\G-$local systems on surfaces 
 uses systematically the {\it amalgamation} of cluster  varieties \cite{FG05}, recalled in Section \ref{sec2}, 
which allows us to assemble more complicated 
cluster varieties from the simpler ones. We shall pick an ideal triangulation  of $\bS$,  and then assemble 
a  moduli  space assigned to $\bS$  from  the elementary ones assigned to the triangles.  Note that the  spaces ${\mathscr X}_{\G, \bS}$,  defined in \cite{FG03a} for    surfaces with corners,   are  not suitable for this purpose. For example, the 
 space   ${\mathscr X}_{ {\rm PGL_2}, t}$ for a triangle $t$    is just a point. But we can not build a space from finitely many points. 
This raises a key question: 
\be \la{CORRPS}
\mbox{\it What is the ``correct" Poisson moduli space assigned to a triangle, or a surface with corners?}
\ee
The answer is given by    moduli spaces  ${\mathscr P}_{\G, \bS}$ assigned to  {colored decorated surfaces   $\bS$}, which we define   now. 
  
\subsubsection{The moduli space ${\mathscr P}_{\G, \bS}$.} \la{1.1.3}   Recall the flag variety ${\cal B}$ and the decorated flag variety ${\mathcal A}:=\G/\U$, 
where $\U$ is a maximal unipotent subgroup in $\G$, see Section \ref{2.1.2}.   There is a canonical projection $\pi:{\mathcal A} \to{\cal B}$. The flag $\pi(\A)$ is called the {\it underlying flag} for a decorated flag $\A$. 
The $\G$-orbits $(\A_1, \A_2)$ on the space of pairs of decorated flags are described by the following two invariants, see Lemma \ref{LD2.2}:
\be \la{INVWH}
w(\A_1, \A_2) \in W, \ \ \ \ h(\A_1, \A_2)\in \H. 
\ee
We say that a pair $(\A_1, \A_2)$ is generic if $w(\A_1, \A_2)$ is the longest element $w_0\in W$. \vskip 1mm

A reduction of a $\G$-local system   ${\cal L}$ to a Borel subgroup near a   point $p$ is the same thing as 
 a  flat section   of  the local system ${\cal L}_{\cal B}:= {\cal L}\times_{\G} {\cal B}$  of flag varieties   near $p$.  Abusing terminology, we call it 
  a {\it flag at   $p$}, and denote  by $\B_p$. The point $p$  can be a puncture. Then we say that $\B_p$ is {\it an invariant flag near $p$}. A decorated flag  $\A_p$ at   $p$ is flat sections of the local system 
${\cal L}_{\cal A}:= {\cal L}\times_{\G} {\cal A}$ of decorated flags near $p$. \vskip 1mm

Let us now introduce pinnings, which play the crucial role in the paper. 
\vskip 1mm

  We start with a $\G-$local system ${\cal L}$ on an oriented interval  ${\rm I}$,    with a pair of flags $(\B_s, \B_t)$ 
      at the fibers of ${\cal L}$ at the endpoints $s,t$ of ${\rm I}$. We assume that moving the flags to the same fiber  we get a generic pair of flags,  
      and that the order $(s,t)$ provides the orientation of ${\rm I}$.       
        Note that ${\cal L}$  is trivial, but not  trivialised.  
        
        Assuming the group $\G$ is adjoint, a  {\it pinning} $p_{\rm I}$   of the data $({\cal L}, (\B_s, \B_t))$ on an oriented interval ${\rm I}$        is a pair of decorated flags  $(\A_s, \A_t)$ over the flags $(\B_s, \B_t)$ such that         \be
       h (\A_s, \A_t)=1.
        \ee 
      See Definition \ref{DEFPIN} and Lemma \ref{PIN1} for more details about the pinnings.   \vskip 1mm

   The key fact about pinnings is that  pinning $p_{\rm I}$  provides a trivialization 
   of   ${\cal L}$,   
   which identifies the flags $(\B_s, \B_t)$ with the  standard flags $(\B_-, \B_+)$ in $\G$. \vskip 1mm
   
   Changing  the orientation of ${\rm I}$,  the pinning $p=  (\A_s, \A_t)$ gives rise to   the {\it opposite pinning}  $p^*:= (\A_t, \A_s)$, see Definition \ref{opposite.pinning.2020}. 
   This operation is evidently reflexive: $p^{**}=p$. \vskip 1mm
   
    Given a colored decorated surface $\bS$ we  assume that the boundary is orientated by the orientation of $\bS$. The pinnings are assigned 
  to the canonically oriented boundary intervals, unless   stated otherwise. \vskip 1mm
  
      \bd \la{DEFP2} Let $\G$ be a split semi-simple adjoint algebraic group over $\Q$, and $\bS$ a colored decorated surface. 
 The moduli space   ${\mathscr P}_{\G, \bS}$  parametrizes data $({\cal L}, \beta, p)$, where 
 
\begin{itemize}

\item  ${\cal L}$ is $\G$-local system    on $\bS$. 
 
\vskip 1mm

 \item  $\beta$ is a framing on  ${\cal L}$, given by   invariant flags     $\B_p$ 
  near   punctures $p$, and   invariant flags     $\B_s$ 
  near   special points $s$, such that 
   \be \la{KAs}
 \mbox{for each colored boundary interval ${\rm I}$, the pair of flags  at its ends   is generic}.\footnote{ Note that a decorated flag $\A_s$  at a point $s$ determines a flag  $\pi(\A_s)$  at   $s$.} 
 \ee
\vskip 1mm  
\item  $p$ is a collection of pinnings $\{p_{\rm I}\}$  on   colored decorated boundary intervals ${\rm I}$ of $\bS$. 
\end{itemize}
   \ed


   
%
      \vskip 2mm

 When $\bS$ has no colored    intervals,   the  space ${\mathscr P}_{\G, \bS}$ coincides with the one  ${\mathscr X}_{\G, \bS}$,  
introduced in \cite{FG03a}. \vskip 2mm

For each non-colored boundary interval ${\rm J}$, the pair of flags $(\B_s, \B_t)$ at its ends   gives rise to an invariant 
$$
w_{\rm J}:= \omega(\B_s, \B_t) \in W,
$$
called  the 
$w-$distance   between $\B_s$ and $\B_t$, see (\ref{INVWH}). The moduli space ${\mathscr P}_{\G, \bS}$ is stratified by the data $\{w_{\rm J}\}$. \vskip 1mm

In particular, given a positive simple root $\alpha$, there are two moduli spaces ${\mathscr P}_{\G, \alpha}$ and ${\mathscr P}_{\G, \overline \alpha}$, 
 related to an oriented triangle with   a single non-colored "short" side ${\rm J}$, where  
$w_{\rm J}=s_\alpha$ is the simple reflection at $\alpha$. These moduli spaces, and their variants with  "short pinnings" at the side ${\rm J}$,  are our main building blocks, see   Section \ref{pcs.sec}. 
To explain how we build the moduli spaces from them we need the crucial operation of {\it gluing}.     \vskip 2mm 


\subsubsection{Gluing maps.}  \la{1.2.3}

Each colored boundary interval ${\rm I}$ of $\bS$ carries a pinning $p_{\rm I}$. The pinnings allow  to define the 
   {\it gluing map}. Its very existence  is 
the key advantage of the space ${\mathscr P}_{\G, \bS}$.  Let us elaborate on this.   \vskip 2mm

Given a decorated
surface $\bS$, possibly disconnected, we choose two   colored boundary intervals ${\rm I}_1$ and ${\rm I}_2$. 
Gluing  
the intervals ${\rm I}_1$ and ${\rm I}_2$ so that their orientations are opposite to each other, 
we get an oriented  
surface
 $$
\bS':= \bS/({\rm I}_1 \sim {\rm I}_2), 
$$  
see Figure \ref{fga5}.  The images of   punctures,  special  points, and colored boundary intervals  on $\bS$ define   a  colored decoration of $\bS'$. 
\begin{figure}[ht]
\epsfxsize 200pt
\center{
\begin{tikzpicture}[scale=0.8]
\draw (1,-1) -- (0,0);
\draw (0,0) -- (1,1);
\draw[dashed, thick, red] (1,1) -- (1,-1);    
\draw (1.3,-1) -- (2.3,0) -- (1.3, 1);
\draw[red, thick, dashed] (1.3,1) -- (1.3, -1); 
\node[black] at (0,0) {{\tiny $\bullet$}};
\node[black] at (1,-1) {\tiny $\bullet$};
\node[black] at (1,1) {\tiny $\bullet$};
\node[black] at (1.3,-1) {\tiny $\bullet$};
\node[black] at (2.3,0) {\tiny $\bullet$};
\node[black] at (1.3,1) {\tiny $\bullet$};
\draw[ultra thick, -latex] (3,0) -- (4,0);
\begin{scope}[shift={(1,0)}];
\draw (5,-1) -- (4,0);
\draw (4,0) -- (5,1);
\draw (5,1) -- (6,0) -- (5,-1);    
\draw[dashed] (5,-1) -- (5,1); 
\node  at (4,0) {\tiny $\bullet$};
\node  at (5,1) {\tiny $\bullet$};
\node  at (6,0) {\tiny $\bullet$};
\node  at (5,-1) {\tiny $\bullet$};
\end{scope}
\end{tikzpicture}
 }
\caption{Gluing a quadrangle out of two triangles with two vertical colored sides.}
\label{fga5}
\end{figure}

A point of the  moduli space ${\mathscr P}_{\G, \bS }$ provides   a $\G$-local system ${\mathscr L}$ with  pinnings  on 
  the boundary intervals ${\rm I}_1, {\rm I}_2$.   
As   explained above, the pinnings  trivialize the restrictions of the local system ${\mathscr L}$ to  the boundary intervals ${\rm I}_1, {\rm I}_2$.  
Therefore there is a unique isomorphism of   the restrictions of the local system ${\mathscr L}$ 
to the  intervals ${\rm I}_1$ and ${\rm I}_2$ which identify    the   pinning on ${\rm I}_1$ with the \underline{opposite pinning} on ${\rm I}_2$. This way we get a new $\G-$local system 
${\cal L}'$ on the decorated surface $\bS'$, see Figure \ref{fga11}.     It inherits    from ${\cal L}$ the framing and the pinnings on the colored intervals of $\bS'$.    So we get a gluing map, see  Lemma \ref{GLUINGM}: 
 \be \la{AMA}
 \begin{split}
 \gamma_{{\rm I}_1, {\rm I}_2}:  {\mathscr P}_{\G, \bS }   &\lra {\mathscr P}_{\G, \bS'}.\\
 ({\cal L}, \mbox{\rm boundary data}) &\lms ({\cal L}', \mbox{\rm induced boundary data}). 
\end{split}
\ee

  The gluing maps allow   to assemble the  space ${\mathscr P}_{\G, \bS}$  
  from the   ones assigned to the   most basic decorated surface: the  triangle, that is   a disc with three special   points given by the vertices.  
We start with an { ideal triangulation} ${\cal T}$ of   $\bS$. 
 Any decorated surface $\bS$, except  a bigon and once-punctured disc with a single special point, admits an ideal triangulation ${\cal T}$  without self-folded triangles. Any two 
   such triangulations are related by a sequence of flips.   
A triangle $t$  
 gives rise to a moduli space  ${\mathscr P}_{\G, t}$.
    Gluing the surface $\bS$ from the triangles $t$ of   ${\cal T}$, and using  maps (\ref{AMA}) for each pair of the glued sides, we get the gluing map:
    \be  \la{GMAT}
\gamma_{\rm T}: 
\prod_{t \in {\cal T}} {\mathscr P}_{\G, t} \lra {\mathscr P}_{\G, \bS}.
\ee
 
\begin{figure}[ht]
\epsfxsize 200pt
\center{
\begin{tikzpicture}[scale=1]
\draw (1,-1) -- (0,0) -- (1,1);    
\draw (1.7,-1) -- (2.7,0) -- (1.7,1); 
\draw[dashed, thick, red, arrows={Latex[length=2mm 3 0]-}] (1,1) -- (1,-1);
\draw[dashed, thick, red, arrows={Latex[length=2mm 3 0]-}] (1.7,-1) -- (1.7,1);
\node[blue] at (0.8,0.2) {\small $p$};
\node[blue] at (1.9,-0.2) {\small $q$};
\begin{scope}[shift={(3.5,0)}]
\draw (1.3,-1) -- (0.3,0) -- (1.3,1);    
\draw (2,-1) -- (3,0) -- (2,1); 
\draw[dashed, thick, red, arrows={Latex[length=2mm 3 0]-}] (1.3,1) -- (1.3,-1);
\draw[dashed, thick, red, arrows={-Latex[length=2mm 3 0]}] (2,-1) -- (2,1);
\node[blue] at (1.1,0.2) {\small $p$};
\node[blue] at (2.2,-0.1) {\small $q^\ast$};
\end{scope}
\begin{scope}[shift={(4,0)}]
\draw (5,-1) -- (4,0) -- (5,1) -- (6,0) -- (5,-1);    
\draw[thick] (5,-1) -- (5,1); 
\node[blue] at (7,0) {$p=q^\ast$};
\end{scope}
\end{tikzpicture}
 }
\caption{Gluing framed local systems with pinnings on two triangles.}
\label{fga11}
\end{figure}

\subsubsection{Cutting maps.} \la{cuttingmap} Given a   non-boundary segment $\sigma$ connecting two marked points on $\bS$, we can cut 
the surface $\bS$ along  $\sigma$, getting new colored decorated surface  $\bS-\sigma$. The two new boundary segments on  $\bS-\sigma$ are non-colored.  
Then there is the corresponding cutting map: 
\be
 {\rm C}_\sigma: {\mathscr P}_{\G, \bS} \lra {\mathscr P}_{\G, \bS-\sigma}.
\ee    
It is  obtained by   restricting  
framed $\G-$local system ${\cal L}$ with pinnings from $\bS$ to   the  surface $\bS-\sigma$.  The invariant flags at the special points/punctures and the pinnings are inherited from ${\cal L}$.

  \subsubsection{A better variant, and the cutting maps.} \la{1.2.5} Given a simple loop $\alpha$ on $\bS$ we can  cut 
the surface $\bS$ along $\alpha$, getting new colored decorated surfaces $\bS-\alpha$.  However there is an issue with the   cutting maps for the moduli spaces, 
since there are no   
 invariant flags nearby the two new boundary loops on $\bS-\alpha$. To handle  this, we introduce a more flexible variant of colored decorated surfaces and related moduli spaces ${\mathscr P}_{\G, \bS}$.  \vskip 1mm
 
 We alter the  definition of a colored decorated surface, distinguishing now the {\it punctures} from the {\it holes}.
 
 \bd A  {\it fully  decorated surface} $\bf S$ is an oriented  surface  with   {\em punctures} inside, and 
a finite number of {\it special points} on the  boundary, modulo isotopy. 
Boundary components without  special points are called {\em holes}. 
 \ed

Defining the   moduli space ${\mathscr P}_{\G, \bf S}$, we  require   invariant flags near   holes, but   not   near
 the punctures.       We do not specify any invariant flag near the punctures. 
 
 \bd \la{DEFP2*} Let  $\bf S$ a fully   decorated surface. 
 The stack   ${\mathscr P}_{\G, \bf S}$  parametrizes data $({\cal L}, \beta, p)$, where 
 
\begin{itemize}

\item  ${\cal L}$ is $\G$-local system    on $\bf S$. 
 
\vskip 1mm

 \item  $\beta$ = $\{$invariant flags     $\B_h$ 
on holes, and  flags     $\B_s$ 
  near   special points $s$, satisfying condition (\ref{KAs})$\}$.
  
\vskip 1mm

\item  $p = \{p_{\rm I}\}$  are pinnings on   colored decorated boundary intervals ${\rm I}$ of $\bf S$. 
\end{itemize}
   \ed

So for a surface $\bf S$ is a surface $S$ which has    punctures,   but no holes and boundary components, 
then   ${\mathscr P}_{\G, \bf S}$ is   the  stack ${\rm Loc}_{\G, S}$ of $\G-$local systems on $S$.  
 The stack ${\mathscr P}_{\G, \bf S}$ from Definition \ref{DEFP2*} is the right generalization of  stacks of $\G-$local systems to decorated surfaces. 
 
 In the presence of  punctures,   the  moduli spaces   from Definition \ref{DEFP2*}  may not even be unirational. Therefore, although  they are Poisson, they 
   no longer carry a cluster Poisson structure.    However, thanks to the cluster nature of the Weyl group actions at the punctures  on the  moduli spaces ${\mathscr P}_{\G, \bS}$ from Definition \ref{DEFP2}, see Theorems \ref{MTH} and \ref{Th8.5A}, the moduli spaces ${\mathscr P}_{\G, \bf S}$  can be quantized, see Section \ref{sqmsls}.  \\

Given a simple loop $\alpha$ on $\bf S$, we define the fully colored surface ${\bf S}-\alpha$ by requiring    the two new boundary components to be punctures. This  allows to 
define the cutting maps
\be
{\rm C}_\alpha: {\mathscr P}_{\G, \bf S} \lra {\mathscr P}_{\G, \bf S-\alpha}.
\ee

\vskip 1mm

  Just to avoid  overloading  the statements of   theorems, 
   we stick to   
  moduli spaces   ${\mathscr P}_{\G, \bS}$ from Definition \ref{DEFP2}.  
  
  Yet the variant  distinguishing  punctures and  holes is certainly the   preferable one. 
  
\subsubsection{Discrete symmetries of the moduli space ${\mathscr P}_{\G, \bS}$.}  \la{2.1.6}  Recall the braid group ${\mathbb B}_{\mathfrak g}$ assigned to the Lie algebra ${\mathfrak g}$ of the group $\G$, see Section \ref{SEC12.1b}. The Weyl group $W$ of ${\mathfrak g}$ has a canonical set theoretic  embedding   $\mu$  into the braid group  
${\mathbb B}_{\mathfrak g}$. 
 Let $w\to w^*$  be the braid group automorphism  
 given by the conjugation by $\mu(w_0)$. Since $\mu(w_0)^2$ lies in the center of ${\mathbb B}_{\mathfrak g}$,  this is an involution. 
  
  Consider the subgroup of stable points of the involution given by the conjugation by  $\mu(w_0)$:
\be \la{SBG}
{\mathbb B}^\ast_{\mathfrak g} = \{ w \in {\mathbb B}_{\mathfrak g} ~|~ ww_0 = w_0w\} \subset {\mathbb B}_{\mathfrak g}.
\ee
\be
 \mbox{\it Let $d_\pi$ be the number of special points on a boundary component $\pi$ of a decorated surface $\bS$.}
 \ee
  We associate with  a boundary component  $\pi$ of $\bS$ a  group  
 \be
\la{BR*}
\begin{split}
&{\mathbb B}^{(\pi)}_{\mathfrak g} =   \left\{  \begin{array}{ll}    {\mathbb B}_{\mathfrak g}     & \mbox{if  $d_\pi$   is even,}  \\
      {\mathbb B}^\ast_{\mathfrak g} & \mbox{if   $d_\pi$   is odd.} \\
   \end{array}\right.\\
   \end{split}
\ee
  Equivalently,   ${\mathbb B}^{(\pi)}_{\mathfrak g}$ is the subgroup of ${\mathbb B}_{\mathfrak g}$  stable under the action of the involution 
  $\mu(w_0)^{d_\pi}$. Since $\mu(w_0)^{2}$  is in the center, for even $d_\pi$ we get the whole group ${\mathbb B}_{\mathfrak g}$.
  
     \vskip 2mm

The following discrete  groups act (birationally in the second case) on  the  moduli space 
 ${\mathscr P}_{\G, \bS}$: 
 
 \begin{itemize}
 
   \item  The mapping class group $\Gamma_\bS$  of $\bS$ acts by automorphisms of the moduli space ${\mathscr P}_{\G, \bS}$. It essentially goes back to \cite{FG03a}.

\vskip 1mm

  \item  For each   puncture on $\bS$, there is a birational action  of the  Weyl group $W$ of $\G$ on ${\mathscr P}_{\G, \bS}$. The definition of the Weyl group action goes back to \cite[Section 1.2]{FG03a}. We recall its definition and prove its cluster nature in Section \ref{Sec5.2}.

\vskip 1mm

  \item  For each boundary component $\pi$ of $\bS$, we define in Section \ref{SEC5.1} an action of the  braid group ${\mathbb B}^{(\pi)}_{\mathfrak g}$ by automorphisms of    the moduli space 
  ${\mathscr P}_{\G, \bS}$.  
  
\vskip 1mm

  \item The group ${\rm Out}(\G):= {\rm Aut}(\G) / {\rm Inn}(\G)$  of the outer automorphisms of the group $\G$ acts by automorphisms of the moduli space ${\mathscr P}_{\G, \bS}$.  It essentially goes back to \cite{GS16}.  \end{itemize}

 The mapping class group $\Gamma_\bS$   permutes   punctures and       boundary components  with the same number of special points. 
So there is a canonical map
\be \la{CANP}
\Gamma_\bS \lra {\rm Perm} \Bigl(\{ \mbox{punctures of $\bS$}\} \times \{\mbox{connected components of the boundary of $\bS$}\}  \Bigr).
\ee
Consider the following semi-direct product, where    
   $\Gamma_\bS$ acts   via map (\ref{CANP}), and   ${\rm Out (\G})$ acts naturally:   
\be \la{GGBG}
\Gamma_{\G, \bS}:= (\Gamma_\bS \times {\rm Out} (\G)) \ltimes \Bigl (\prod_{\{\mbox{components $\pi$ of $\partial\bS$}\}}  
{\mathbb B}^{(\pi)}_{\mathfrak g} \times \prod_{\{\mbox{punctures of $\bS$}\}} W\Bigr ).
\ee
The   group $\Gamma_{\G, \bS}$   acts on the space ${\mathscr P}_{\G, \bS}$.

\vskip 2mm \paragraph{\it The subgroup $\Gamma'_{\G, \bS}$.} 
For technical reasons, we consider a subgroup 
$\widetilde {\mathbb B}^{\ast}_{\mathfrak g} \subset {\mathbb B}^\ast_{\mathfrak g}$   generated by the elements $\mu(w)$ where $w\in W$ and $w^* = w$.  
It is  very likely that $\widetilde {\mathbb B}^{\ast}_{\mathfrak g} = {\mathbb B}^{\ast}_{\mathfrak g}$. 
  We set  
  \be
\la{BR*}
\begin{split}
&\widetilde {\mathbb B}^{(\pi)}_{\mathfrak g} =   \left\{  \begin{array}{ll}    \widetilde {\mathbb B}^{\ast}_{\mathfrak g}    & \mbox{if $d_\pi$     is odd,}  \\
      {\mathbb B}_{\mathfrak g} & \mbox{if   $d_\pi$   is even.} \\
   \end{array}\right.\\
   \end{split}
\ee 
Let  $ \Gamma'_{\G, \bS}$ be the subgroup of
 $\Gamma_{\G, \bS}$ given by  semi-direct product (\ref{GGBG}) where   ${\mathbb B}^{(\pi)}_{\mathfrak g}$ is replaced by  $\widetilde {\mathbb B}^{(\pi)}_{\mathfrak g}$. \vskip 2mm
 
 To state the main result, we have to recall first what is a cluster Poisson variety and its $q-$deformation. We start with the most basic example. 
\subsubsection{Quantum tori} \la{2.1.5}
Let $\Lambda$ be a lattice  with a skew symmetric form  
$\langle \ast, \ast\rangle: \Lambda \times \Lambda \to \Z$. It 
gives rise to a split  algebraic torus ${\rm T}_{\Lambda}:=  {\rm Hom}(\Lambda, {\mathbb G}_m)$ equipped with a Poisson bracket
\be \la{LOGCa}
\{X_\lambda, X_\mu\}:= 2\langle \lambda, \mu\rangle X_\lambda X_\mu, \qquad  \forall \lambda, \mu \in \Lambda.
\ee
Here $X_\lambda$ is the character of  ${\rm T}_{\Lambda}$ assigned to   $\lambda$. 
One quantizes ${\rm T}_{\Lambda}$, getting  a {\it quantum torus algebra} ${\cal O}_{q}({\rm T}_\Lambda)$, given by 
  a free  $\Z[q, q^{-1}]$-module    with a basis  $X_{\lambda}$,  $\lambda \in \Lambda$, and  the 
product
\be \la{QTA}
X_{\lambda}X_{\mu} = 
q^{\langle \lambda, \mu\rangle}X_{\lambda + \mu}. 
\ee
 The algebra ${\cal O}_{q}({\rm T}_\Lambda)$ is the group algebra of the discrete Heisenberg group 
assigned to $\Lambda$ and the form $\langle\ast, \ast \rangle$, where $q$ is a generator of the center. Its quasi-classical limit $q\to 1$ recovers   Poisson algebra \eqref{LOGCa}. 

 \subsubsection{Cluster Poisson varieties.} A space ${\mathscr X}$   with a cluster Poisson structure   carries an atlas of  log-canonical rational coordinate systems. Any two of them
 are related by specific birational transformations, called 
{\it cluster Poisson transformations}. The latter  are   compositions of monomial transformations and 
   {\it cluster Poisson mutations}. 
    Here {\bf q}re detailed account. We start from the notion of a quiver. 
    
\begin{definition} \la{DQUIV}
{\it A {\it quiver} 
 $ {\bf q}$ is a datum 
$
\Bigl(\Lambda, (\ast, \ast), \{e_i\}, \{g_j\}, \{d_i\}\Bigr), 
$ where: 
\vskip 1mm
i)  $\Lambda$ is a lattice,
 $\{e_i\}$   a basis of $\Lambda$, and 
$\{g_j\}\subset \{e_i\}$  a subset of {\it frozen} basis vectors, 
\vskip 1mm

ii)  $(\ast, \ast)$ a $\frac{1}{2}\Z$-valued 
bilinear form on $\Lambda$, with 
$(e_i, e_j) \in \Z$ 
unless both  
$e_i$, $e_j$ are frozen,  
\vskip 1mm

iii)  $d_i\in \Z_{>0}$ are multipliers such that 
$\langle e_i, e_j\rangle := (e_i, e_j)d^{-1}_j$ is skew-symmetric.} 

\end{definition}

 A quiver ${\bf q}$ is the same thing as  
  a {\it marked quiver}.   
  Indeed, the  vertices of the quiver are  given by the basis vectors $e_i$;  they   
  are marked by the integers $d_j$.  There are $(e_i, e_j)$ arrows  
   $i\to j$ if $( e_i, e_j)>0$, and $(e_j, e_i)$ arrows   $j\to i$ otherwise. Note that we allow {\it half-arrows} between the frozen vertices.

The   lattice $\Lambda_{\bf q}$ of a quiver ${\bf q}$ gives rise to a   torus ${\rm T}_{\bf q}:= {\rm Hom}(\Lambda_{\bf q}, {\Bbb G}_m)$. 
The       skew-symmetric $\frac{1}{2}\Z$-valued  bilinear form $\langle \ast, \ast\rangle_{\bf q}$ on the lattice $\Lambda_{\bf q}$  
provides   log-canonical Poisson structure (\ref{LOGCa}) on the torus.

 A cluster Poisson structure on a space ${\mathscr X}$ is given by   a  collection   of  Poisson 
tori ${\rm T}_{\bf q}$ assigned to certain quivers ${\bf q}$, and  birational isomorphisms
\be  \la{ic1}
i_{\bf q}: {\rm T}_{\bf q} \lra {\mathscr X}.
\ee 
 
 The    basis $(e_1, ..., e_n)$ of  $\Lambda_{\bf q}$  provides, via the   map $i_{\bf q}$,   a ra{\bf q}tional cluster Poisson coordinate system 
$(X_{e_1}, \ldots, X_{e_n})$ on  the space  ${\mathscr X}$. 
 For each basis vector $e_k$ there is  an elementary transformation of quivers $\mu_{e_k}: {\bf q} \to {\bf q'}$, called {\it quiver mutation}, see (\ref{12.12.04.2a}), which give rise to birational 
Poisson transformations 
\be \la{MUTa}
{{\mu}_{e_k}}: {\rm T}_{\bf q} \lra {\rm T}_{{\bf q}'},
\ee
intertwining the maps $i_{\bf q}$ in (\ref{ic1}), that is making the following diagram commute:
 \begin{displaymath}
    \xymatrix{
        {\rm T}_{\bf q} \ar[r]^{i_{{\bf q}}} \ar[d]_{{\mu}_{e_k}} &{\mathscr X}         \ar[d]^{{=}} \\
         {\rm T}_{{\bf q}'} \ar[r]^{i_{{\bf q}'}} &  {\mathscr X}}
         \end{displaymath}

Any two cluster coordinate systems ${\bf q}_1$ and ${\bf q}_2$ of the  cluster Poisson atlas are related by a   sequence of quiver mutations and quiver isomorphisms. 
Their composition is called a {\it cluster Poisson transformation}, and denoted   by $\mu_{{\bf q}_1 \to {\bf q}_2}$.

\subsubsection{$q$-deformations of  cluster Poisson varieties.} Denote by ${\cal F}_q({\rm T}_{\bf q})$ the fraction field of the  quantum torus  algebra ${\cal O}_q({\rm T}_{\bf q})$. 
Mutation maps (\ref{MUTa}) can be quantized, providing    isomorphisms 
\be \nonumber
{{\mu}^*_{e_k}}: {\cal F}_q({\rm T}_{\bf q'}) \lra {\cal F}_q({\rm T}_{\bf q}).
\ee
So any cluster Poisson transformation $\mu_{{\bf q}_1 \to {\bf q}_2}$ gives rise to  a quantum cluster transformation 
\be \la{MUTT1a}
  \mu^*_{{\bf q}_1 \to {\bf q}_2}: {\cal F}_q({\rm T}_{\bf c_2})  \stackrel{\sim}{\lra} {\cal F}_q({\rm T}_{\bf c_1}).  
\ee 
  
  \bd The  algebra ${\cal O}_q({\mathscr X})$ consists  of all  elements $f\in {\cal O}_q({\rm T}_{\bf q})$  which remain 
 Laurent polynomials  over $\Z[q, q^{-1}]$ under any quantum cluster transformation (\ref{MUTT1a}). 
 
 For $q=1$ we get the algebra of universally Laurent polynomials on ${\mathscr X}$, 
 denoted by    ${\cal O}^{\rm cl}({\mathscr X})$.
 \ed
So   for each quiver ${\bf q}$ of the cluster atlas there is a canonical  embedding  of  algebras 
\be \la{astA1}
\begin{split}
&i^*_{\bf q}: {\cal O}_q({\mathscr X})  \subset {\cal O}_q({\rm T}_{\bf q}).\\
\end{split}
\ee
Cluster transformations identify   subalgebras (\ref{astA1}). So there are commutative diagrams
 \begin{displaymath}
    \xymatrix{
        {\cal O}_q({\mathscr X})  \ar[r]^{\mu^*_{{\bf q}_1 \to {\bf q}_2}} \ar[d]_{{i}^*_{\bf c_2}} &{\cal O}_q({\mathscr X})           \ar[d]^{{i}^*_{\bf c_1}} \\
          {\cal F}_q({\rm T}_{\bf c_2}) \ar[r]^{\mu^*_{{\bf q}_1 \to {\bf q}_2}} &    {\cal F}_q({\rm T}_{\bf c_1})}
         \end{displaymath}

This structure is best described by a {\it cluster modular groupoid} ${\rm M}_{\mathscr X}$. 
   It is a category, whose objects  are cluster Poisson coordinate systems ${\bf q}$ on ${\mathscr X}$, and morphisms are  cluster Poisson transformations. The  
    fundamental group $ \pi_1({\rm M}_{\mathscr X}, {\bf q})$ is  called the {\it cluster modular group ${\Gamma}_{\mathscr X}$}. 
    So a quantum deformation of   ${\mathscr X}$ is   a functor from {the cluster modular groupoid} ${\rm M}_{\mathscr X}$ 
    to the category of quantum torus algebras and birational isomorphisms between them. \vskip 2mm

  A priori, two different algebras  are assigned to an algebraic variety/stack    ${\mathscr X}$, equipped with a cluster Poisson structure.  One is the 
algebra ${\cal O}({\mathscr X})$ of regular functions  on the variety/stack  ${\mathscr X}$. The other is the algebra 
${\cal O}^{\rm cl}({\mathscr X})$  of { universally Laurent polynomials}. The $q$-deformation of ${\cal O}^{\rm cl}({\mathscr X})$ is the algebra 
${\cal O}_q({\mathscr X})$. This notation is unambiguous: there is no alternative definition of ${\cal O}_q({\mathscr X})$.

\subsubsection{The main result.} 
  To state Theorem \ref{MTH}, 
  we   exclude  a few pairs $(\G, \bS)$:  
  \vskip 2mm
 \noindent
1. We always assume that the total number  $\mu$ of marked points on $\bS$ is positive,  
and    \be \la{COND1} 
 \mbox{\it Exclude discs with $\leq 2$ special points $\&$ punctured discs with $1$ special point if $\G$ of type ${\rm A_1, A_2}$.} 
 \ee
 2. Talking about   cluster nature of the Weyl group action:    
 \be \la{COND2} 
 \mbox{\it If $\G$ is of type $\A_1$, we exclude surfaces $\bS$ with a single puncture and no boundary components.}
 \ee 
3. Talking about   cluster nature of   braid group actions we  sometimes assume for simplicity  that:
 \be \la{COND3} 
 \mbox{\it The total number of punctures and connected components of $\bS$ is at least two}.
 \ee
 
\bt \la{MTH}
Let $\G$  be any   split semi-simple adjoint   algebraic group over $\Q$. Let $\bS$ be a colored 
decorated surface with $n$ punctures and $\mu>0$, excluding the cases (\ref{COND1}). Then

\begin{enumerate}

\vskip 1mm
\item The   space ${\mathscr P}_{\G, \bS}$ has a $\Gamma_\bS \times {\rm Out}(\G)$-equivariant  cluster Poisson  structure. 

\vskip 1mm

\item Excluding the cases (\ref{COND2}),  the   group $W^n$ acts by cluster Poisson transformations. 

\vskip 1mm

\item Assuming  conditions (\ref{COND3}), for each boundary component $\pi$ of $\bS$, the   braid group $\widetilde {\mathbb B}^{(\pi)}_{\mathfrak g}$ acts   on the space  ${\mathscr P}_{\G, \bS}$ by   quasi-cluster transformations,  
and   on  ${\mathscr X}_{\G, \bS}$  by  cluster transformations.\footnote{We   have a proof without assuming condition (\ref{COND3}). 
But it is not used  in this paper, and thus will be presented later on.}  
\end{enumerate}

So, assuming   (\ref{COND2}) and   (\ref{COND3}), the group $\Gamma'_{\G, \bS}$ maps to  the cluster modular groups of ${\mathscr X}_{\G, \bS}$ and ${\mathscr P}_{\G, \bS}$. 
\et

  To prove   Theorem \ref{MTH} we  construct  a  $\Gamma_\bS\times {\rm Out}(\G)$-equivariant cluster Poisson  atlas   on  the   space 
  ${\mathscr P}_{\G, \bS}$, with   a finite number of coordinate systems modulo the action of $\Gamma_\bS$. 
 The construction  works for    all $\G$ at once. 
It is new  even for ${\rm PGL}_m$,  where it describes explicitly a larger and much more flexible class of cluster Poisson coordinate systems than in \cite{FG03a}.  
For example, the 
cluster nature of  the group ${\rm Out}(\G)$  action   was established for ${\rm PGL}_m$ in \cite[Section 9]{GS16}  
by rather elaborate arguments, and takes 10 pages. The present proof, for any group $\G$,  
 is obvious from  the  definition of the cluster atlas on the space ${\mathscr P}_{\G, \bS}$.  This   shows 
 that even for ${\rm PGL}_m$,   the current approach is much more advanced  than   the one   in \cite{FG03a}. 

   \vskip 2mm
  
Cluster Poisson coordinates   are  rational functions on the moduli space ${\mathscr P}_{\G, \bS}$. Note that a cluster structure on a variety is completely determined by exhibiting a single 
cluster, that is a  quiver together with a   coordinate system consisting of functions, labeled by the   vertices of the quiver. 

 We give a very simple and transparent definition of   Poisson clusters assigned to ideal triangulations of $\bS$. 
It  takes just a couple of pages:  the cluster coordinates are defined  in Section \ref{CSEC4.3}, and the quivers  
in  Section \ref{sec.3.1}. 
However the proof of the key fact that  these clusters belong to the same cluster atlas is rather non-trivial.  
It takes Sections \ref{pcs.sec} - \ref{CSEC7} to accomplish.

 \vskip 2mm
The proof of   the cluster nature of the braid group action is  transparent for the space ${\mathscr X}_{\G, \bS}$, see Section \ref{SEC5.1}, but becomes quite demanding and involved 
 for the space   ${\mathscr P}_{\G, \bS}$ due to the presence of the frozen coordinates. We present the proof  in Sections \ref{SECT9.4} - \ref{sect10.4}.

 \vskip 2mm 
  Theorem \ref{MTH}    
  gives interesting examples of   cluster modular groups enlarging   mapping class groups.  Cluster nature of the group $ \Gamma'_{\G, \bS}$ action is   crucial for   applications.  
  Here are the first one.

  \bd The algebra ${\cal O}_q({\mathscr P}_{\G, \bS})$ consists  of all  elements $f\in {\cal O}_q({\rm T}_{\bf q})$   remaining  
 Laurent polynomials  for any quantum cluster Poisson transformation   for the atlas fromTheorem \ref{MTH}. 
\ed

\bt  \la{Th8.5A}  The    group $\Gamma'_{\G, \bS}$   acts by automorphisms of  the algebra  ${\cal O}_q({\mathscr P}_{\G, \bS})$  given by quantum quasi-cluster transformations. 
They are quantum cluster transformations on the  subgroup $\Gamma_\bS \times W^n$. 
\et

 Theorem \ref{Th8.5A} is proved in Section \ref{Sec8.4}.

\subsubsection{Construction of the cluster Poisson structure on the   space ${\mathscr P}_{\G, \bS}$.} \la{1.2.9} It  goes in three steps.

  \begin{enumerate}
  
\item   For a  colored triangle $t$,  we introduce   cluster Poisson coordinates on  ${\mathscr P}_{\G, t}$ by a simple explicit construction in Section \ref{CSEC4.3}, 
 amalgamating   elementary cluster Poisson varieties ${\mathscr P}_{\G, \alpha}$ assigned to simple positive 
roots. The construction depends on a choice of a reduced 
  decomposition of the  element $w_0$   into a product of   generators, and   a choice of an angle of the triangle $t$. We prove:\vskip 1mm
  
  {\it Cyclic invariance.} The resulting cluster structure does  not depend neither on the choice of the reduced decomposition of $w_0$, nor on the angle of the triangle   $t$. \vskip 1mm

   \item  For each side ${\rm I}$ of the triangle $t$ there is a  canonical   action   of the Cartan group $\H$ on  ${\mathscr P}_{\G, t}$,  see (\ref{ALa}).  
    It provides      
  a    free action of   $\H^3$ on ${\mathscr P}_{\G, t}$.   From the cluster Poisson perspective, the $\H^3-$orbits are  the ``frozen fibers", and   the quotient ${\mathscr P}_{\G, t}/\H^3$ is  the ``non-frozen base"  for  ${\mathscr P}_{\G, t}$,  see (\ref{FRX}). The base is    the   space  ${\rm Conf}^\times_3({\cal B})$ of generic  
  $\G-$orbits on the space of triples of flags: 

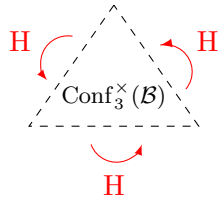
\[
\begin{tikzpicture}[scale=1.2]
\draw[red, thick, dashed] (90:1) -- (200:1) -- (340:1) -- (90:1);    
\node at (0,0) {{\footnotesize ${\rm Conf}_3^\times(\mathcal{B})$}};
\node[red] at (-90:1) { ${\rm H}$};
\node[red] at (30:1.2) { ${\rm H}$};
\node[red] at (150:1.2) {${\rm H}$};
 \draw[red, -latex] (245:.6) arc (200:340:.3);
  \draw[red, -latex] (125:.8) arc (80:220:.3);
  \draw[red, directed] (5:.8) arc (-40:100:.3);
\end{tikzpicture}
\]

If the triangle $t$ has non-colored sides, we just drop the pinnings at those sides. 

We assemble  a cluster Poisson structure on  
  ${\mathscr P}_{\G, \bS}$  
  by amalgamating   the ones on    ${\mathscr P}_{\G, t}$ over    
  triangles  $\{t\}$ of an ideal triangulation ${\cal  T}$ of $\bS$.   For   triangles $t_1, t_2$  sharing a colored edge,  
  the amalgamation ${\mathscr P}_{\G, t_1} \ast {\mathscr P}_{\G, t_2}$ amounts to  the quotient of ${\mathscr P}_{\G, t_1} \times{\mathscr P}_{\G, t_2}$ by the antidiagonal action of $\H$:      
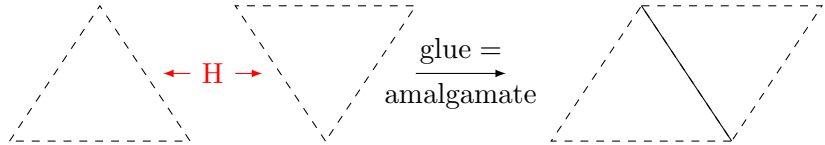
\[
\begin{tikzpicture}[scale=1.2]
\draw (0,0) -- (1,1.5);
\draw[red, thick, dashed] (1,1.5)  -- (2,0) ;
\draw  (2,0)-- (0,0);    
\draw[red, thick, dashed] (3.5,0) -- (2.5,1.5);
\draw(2.5,1.5) -- (4.5,1.5) -- (3.5,0);  
\draw[red, latex-] (1.7, 0.75) -- (2, 0.75); 
\draw[red, -latex] (2.5, 0.75) -- (2.8, 0.75); 
\node[red] at (2.25, 0.75) { ${\rm H}$};
\draw  (7,1.5) -- (9,1.5) -- (8,0) -- (7,1.5)--(6,0)--(8,0); 
\draw  (8,0) -- (7,1.5);       
 \draw[-latex] (4.5, 0.75) -- (5.5, 0.75);
\node at (5,1) {$\footnotesize{\rm glue =}$};
\node at (5,0.5) {$\footnotesize{\rm amalgamate}$};
\end{tikzpicture}
\]
We prove that the resulting cluster structure does  not depend on a choice of  ${\cal T}$ thanks to:
  
  {\it Flip invariance.} For a quadrangle $q$, the  cluster   structure  
  of the  spaces  ${\mathscr P}_{\G, q}$       
 does not depend on the choice of a diagonal cutting  $q$ into   triangles. \vskip 1mm
 
  \item

   We show that the   gluing  maps  (\ref{GMAT}) are compatible with the cluster amalgamation maps, thus establishing   
   a cluster Poisson structure on the   space ${\mathscr P}_{\G, \bS}$. 
          
  \end{enumerate}

 This implies that the     mapping class group $\Gamma_\bS$ acts on   ${\mathscr P}_{\G, \bS}$   
   by cluster Poisson  transformations. 
   
   Indeed, take an element $\gamma \in {\Gamma}_\bS$, and a generic point $[{\mathscr L}] \in {\mathscr P}_{\G, \bS}$.
      Since our cluster atlas is $\Gamma_\bS$-equivariant,  the  coordinates of  $[{\mathscr L}]$  
   in  a cluster coordinate system ${\bf q}$ assigned to an ideal triangulation ${\rm T}$  are the same as the  coordinates of  $\gamma[{\mathscr L}]$ for the cluster coordinate system 
   $\gamma({\bf q})$ assigned to $\gamma({\rm T})$.        Let us  
   relate  triangulations $\gamma({\rm T})$ and ${\rm T}$ be a sequence of flips.    Then, thanks to the cyclic and flip invariance, 
   there is a cluster transformation   expressing the  coordinates of   $\gamma[{\mathscr L}]$  in the cluster coordinate system  $\gamma({\bf q})$ 
    via the ones for $ {\bf q}$.     This just means that  $\gamma$ acts by  a cluster Poisson transformation.

 

  \medskip
     
  \subsection{Cluster structure of the dual moduli space ${\mathscr A}_{\G, \bS}$} 
   \medskip
   
   Cluster varieties  come in pairs: ${\mathscr A}$ $\&$ ${\mathscr X}$. 
Here  ${\mathscr X}$   
is a {cluster Poisson variety}. The 
 cluster variety ${\mathscr A}$ is its dual. It   is a geometric incarnation of a cluster algebra \cite{FZI}: the algebra of universally Laurent polynomials  ${\cal O}^{\rm cl}({\mathscr A})$ is the 
 upper cluster algebra \cite{BFZ}. 
 The cluster variety 
${\mathscr A}$ carries a dual structure, which includes a class in $K_2({\mathscr A})$. To stress the difference between the two, we  refer to it as 
a {\it cluster $K_2$-variety}. Cluster varieties ${\mathscr A}$ $\&$ ${\mathscr X}$ are related in many   ways, and form a {\it cluster ensemble} \cite{FG03b}.

\subsubsection{The main result.} The moduli space ${\mathscr A}_{\G, \bS}$ was introduced in \cite{FG03a}. It parametrizes twisted $\G$-local systems on $\bS$   equipped with some extra data 
  near the  punctures and special  points, see Definition \ref{DEFP3}. 
This data forces the  monodromy of the local systems around the punctures to be   unipotent.  

\vskip 2mm
 {\bf Convention.} {\it We  abuse notation by assuming that the group $\G$ in  ${\mathscr P}_{\G, \bS}$ is always adjoint, while the group $\G$ in  ${\mathscr A}_{\G, \bS}$ is always simply-connected, unless stated otherwise. }
 \vskip 2mm

The spaces ${\mathscr A}_{\G, \bS}$ and ${\mathscr P}_{\G, \bS}$ are dual to each other 
  in several ways, and related by a   map
\be \la{MPP}
p: {\mathscr A}_{\G, \bS} \lra {\mathscr P}_{\G, \bS}.
\ee

 Given a puncture   on $\bS$, there is a geometric action of  the   Weyl group   on  the space    ${\mathscr A}_{{\G}, \bS}$  
  introduced in \cite{FG03a} for ${\rm SL_2}$, and in  \cite[Section 6]{GS16} in general.  
The geometric action of the braid group $\B^{(\pi)}_\G$ on the space ${\mathscr A}_{{\G}, \bS}$ is defined essentially the same way  as for the space ${\mathscr P}_{{\G}, \bS}$. 

Any   
 cluster transformation acts simultaneously on the cluster Poisson and cluster $K_2-$spaces which form a cluster ensemble. 
 Since we already defined actions of the Weyl and braid groups on the cluster Poisson space 
 ${\mathscr P}_{{\G}, \bS}$, the $K_2-$flavor of the same cluster transformation provides us  cluster actions of the Weyl and braid groups by cluster $K_2-$transformations 
 of the space ${\mathscr A}_{{\G}, \bS}$.
 
 So, as soon as we   prove  that the spaces ${\mathscr A}_{{\G}, \bS}$ and  ${\mathscr P}_{{\G}, \bS}$ form a cluster ensemble, we get two a priori different 
  actions of the Weyl and braid group on the space ${\mathscr A}_{{\G}, \bS}$:  the geometric and the cluster. 
 
 \bt \la{MTHa}
 Under the same assumptions as in Theorem \ref{MTH}, 

\begin{enumerate}

\item  The  space ${\mathscr A}_{\G, \bS}$ has a $({\rm Out}(\G)\times \Gamma_\bS)\ltimes W^n$-equivariant structure of a cluster $K_2$-variety. 
\vskip 1mm

\item Cluster   structures on the spaces ${\mathscr A}_{\G, \bS}$ and ${\mathscr P}_{\G, \bS}$ are dual to each other:  
the pair   
$({\mathscr A}_{\G, \bS}, {\mathscr P}_{\G, \bS})$, related by the map   (\ref{MPP}), has   a $\Gamma'_{\G, \bS}$-equivariant cluster ensemble structure. 
 
\vskip 1mm

\item The cluster variety structures on the  spaces ${\mathscr A}_{\G^\vee, \bS}$ and ${\mathscr P}_{\G, \bS}$ are Langlands dual. 
\vskip 1mm

 \item For any  boundary component $\pi$, the  braid group $\widetilde \B^{(\pi)}_\G$ acts       by    quasi-cluster transformations of ${\mathscr A}_{\G, \bS}$.

\end{enumerate}
\et

 The dual analog of the gluing map (\ref{GMAT}) for the  space ${\mathscr A}_{\G, \bS}$ is the {\it restriction map}:
$$
{\rm Res}_{\rm T}: {\mathscr A}_{\G, \bS} \lra \prod_{t\in {\rm T}} {\mathscr A}_{\G, t}.
$$
 It is injective, and  allows to decompose  the  space ${\mathscr A}_{\G, \bS}$  
into the   ones assigned to the triangles.

\vskip 2mm

We give a very simple  definition of  cluster $K_2-$coordinates  on the moduli space ${\mathscr A}_{\G, \bS}$  assigned to ideal triangulations of $\bS$ 
  in Section \ref{SeCC2.4A}. They   are  regular functions, labeled by the vertices of quivers from 
  Section \ref{sec.3.1}. An alternative way to define these functions is given in Section \ref{CSEC12.3}. 
Just like in the Poisson case, proving that these cluster coordinates are related by  cluster $K_2$-transformations  is not easy, 
and takes Sections \ref{pcs.sec} and \ref{sec.conf3} to accomplish. This   
 implies $\Gamma_\bS$-equivariance of the   cluster structure.

It remains to show that the cluster and geometric actions of the 
Weyl and braid group coincide. For the Weyl group it follows immediately from \cite{GS16}. For the braid group it is   similar to the Poisson case. 
\vskip 2mm

The very existence of cluster coordinates on the moduli spaces ${\mathscr P}_{{\rm G}, { \bS}}$  and ${\mathscr A}_{{\rm G}, { \bS}}$ is a mystery: 
we do not understand yet why such reach and rigid structures on these moduli spaces should exist. 

\vskip 2mm
 The key feature of the pair of spaces   $({\mathscr A}_{\G, \bS}, {\mathscr P}_{\G, \bS})$, which enters to all   dualities,  is    this: 
\be \la{DIMEQ}
  {\rm dim}~{\mathscr A}_{\G, \bS} = {\rm dim}~{\mathscr P}_{\G, \bS}.
\ee
If $\bS$ has special points,  property (\ref{DIMEQ}) fails for the pair   $({\mathscr A}_{\G, \bS},  {\mathscr X}_{\G, \bS})$, see   (\ref{PGSX}).    
 This is why, see   (\ref{CORRPS}), 
\be \la{slogan}
 \mbox{\it  The right  Poisson moduli space  for the pair $(\G, \bS)$ is the  space ${\mathscr P}_{\G, \bS}$, rather than   ${\mathscr X}_{\G, \bS}$.}
\ee

\medskip
  
\subsection{Applications to DT-transformations and  Duality Conjectures}

\medskip 

In \cite[Definition 1.5]{GS16} we introduced a  transformation 
${\mathcal D}{\mathcal T}_{\G, \bS}$ of  ${\mathscr X}_{\G, \bS}$. Let us extend it to   ${\mathscr P}_{\G, \bS}$. 

\bd \la{KSDT} The transformation ${\mathcal D}{\mathcal T}_{\G, \bS}$ of the moduli space ${\mathscr P}_{\G, \bS}$ is defined as follows: 
 \be \la{DTC}
{\mathcal D}{\mathcal T}_{\G, \bS}:= \ast \circ \Bigl(\prod_{\mbox{ \rm components of $\partial \bS$}}\mu(w_0)\Bigr) \circ \Bigl(\prod_{\mbox{\rm punctures of $\bS$}}  {w}_0\Bigr).
\ee
\ed
It is interesting to note that it is  the product of three   different   incarnations of the element $w_0$: 

\begin{itemize}

\item The   element $\ast \in {\rm Out}(\G)$ arises from  the involution $\ast$ of the Dynkin diagram acting on the simple positive roots by $\alpha_{i^*}:= -w_0(\alpha_i)$. 

\vskip 1mm\item The braid group elements $\mu(w_0)\in {\Bbb B}_{\mathfrak g}^{(\pi)}$  for each of the boundary components $\pi$ of $\bS$. 

\vskip 1mm\item The Weyl group elements $w_0$   for each of the punctures of $\bS$. 

\end{itemize}

For the space ${\mathscr X}_{\G, \bS}$ the element $\mu(w_0)$ reduces to  the ``shift  by one" of    special points on each boundary component. 
This is the way it was defined in \cite{GS16}. Definition \ref{KSDT} looks   more natural.

\bt \la{CDT} Let $\bS$ be an admissible colored decorated surface, that is we exclude surfaces (\ref{COND1})-(\ref{COND3}).  Then  

\vskip 1mm

a) For any   adjoint group $\G$, the transformation ${\mathcal D}{\mathcal T}_{\G, \bS}$ is a quasi-cluster Poisson transformation. 

\vskip 1mm

b) The map ${\mathcal D}{\mathcal T}_{\G, \bS}$  is  the Donaldson-Thomas transformation for the  cluster    variety  
${\mathscr P}_{\G, \bS}$.\et

Theorem \ref{CDT} is proved in Section \ref{DT=cluster}. The part  a) follows immediately from 
Theorem \ref{MTH}.  \vskip 1mm

For $\G = {\rm PGL}_m$, Theorem \ref{CDT}  can be deduced  from the similar result in \cite{GS16} for the  space 
${\mathscr X}_{\G, \bS}$.  \vskip 1mm

Since  by the part a) of Theorem \ref{CDT} the transformation ${\mathcal D}{\mathcal T}_{\G, \bS}$ is cluster Poisson,  it remains to verify the    criteria \cite[Definition 1.15]{GS16} to prove
Theorem \ref{CDT}. This is done in Section \ref{DT=cluster}. \vskip 2mm

As    explained in Section \ref{1.0.3}, Theorem \ref{CDT}, combined with the work \cite{GHKK} and using the work \cite{S20}, prove  Duality Conjectures for the   space ${\mathscr P}_{\G, \bS}$ from \cite[Section 12.4]{FG03a}, stated in the cluster set-up in \cite[Section 4]{FG03b}.  Precisely, we get the following.

\bt \la{CBA} Under the same assumptions as in Theorem \ref{CDT}:
\vskip 1mm
\begin{enumerate}

\item  There is a canonical linear basis ${\Bbb I}(l) = {\Bbb I}_{\cal A}(l)$ on the vector space ${\cal O}({\mathscr P}_{\G, \bS})$, parametrised by the elements $l$ of the set 
${\mathscr A}_{\G^\vee, \bS}(\Z^t)$ of integral tropical points $l$ 
of the space ${\mathscr A}_{\G^\vee, \bS}$. This parametrization is $\Gamma_{\G, \bS}-$equivariant.  
\vskip 1mm

\item In any cluster Poisson coordinate system on ${\mathscr P}_{\G, \bS}$ the  elements ${\Bbb I}(l)$ of the canonical basis are Laurent polynomials with positive integral coefficients.

\vskip 1mm

\item One has the product formula
\be \la{PRFO}
{\Bbb I}(l_1) \cdot {\Bbb I}(l_2)  = \sum_{l_3} c(l_1, l_2, l_3){\Bbb I}(l_3).
 \ee
Here the coefficients $c(l_1, l_2, l_3)$ are positive integers, and the sum over $l_3 \in {\mathscr A}_{\G^\vee, \bS}(\Z^t)$ is finite. 

\vskip 1mm

\item   Let $(a_1, ..., a_n)$ be the tropical 
coordinates of a point $l\in {\mathscr A}_{\G^\vee, \bS}(\Z^t)$ in a cluster ${\bf q}$. Then the highest term\footnote{for the partial order on monomials $X^a = \prod_{i\in I}X_i^{a_i}$ such that $X^a\geq X^b$ if $a_i \geq b_i$ for all $i\in I$.} 
of ${\Bbb I}(l)$ in the cluster Poisson coordinates $X_i$ for the cluster ${\bf q}$ is $\prod_iX_i^{a_i}$.

 \end{enumerate}

\et

Property (4) tells how  elements of the canonical basis determine   integral tropical points of the dual space. \\

Here is an important  special case of Theorem \ref{CBA}, which we  prove directly in Section \ref{SECT13.4}.
 
 \bt \la{4.20.24.1} 
 For any   $\G-$local system ${\cal L}$ on $\bS$, any finite dimensional representation $V$ of $\G$, 
 and any loop $\alpha$ on $\bS$, the function ${\rm Tr}{\rm Mon}_\alpha({\cal L}_V)$ given by the trace of the monodromy 
 along $\alpha$ in the local system ${\cal L}_V= {\cal L}\times_\G V$ is a 
 Laurent polynomial with positive integral coefficients in any cluster Poisson coordinates on ${\mathscr P}_{\G, \bS}$.
 \et

 The functions ${\rm Tr}{\rm Mon}_\alpha({\cal L}_V)$ for the  fundamental representations $V$  and simple loops $\alpha$ 
should belong to the canonical basis. 
Moreover, for any multicurve $\sum n_j\alpha_j$,   the product 
$
\prod_j{\rm Tr}{\rm Mon}_{\alpha_j}({\cal L}_{V_j})^{n_j}
$  where  $V_j$ are    fundamentals   should have the same property.  But they give only a  part of the canonical basis if ${\rm rk}(\G)>1$. 
 \medskip

   \subsection{Applications to higher Teichm\"uller theory}\la{SECT2.4a}
\medskip

In Section \ref{SECT2.4a} we show that the canonical basis together with its  positivity properties from Theorem \ref{CBA} is the backbone of the various notions and  results in higher Teichmuller theory, including the $\G-$analogs of the length function and its cluster convexity, the intersection pairing. We use the  latter  to define the $\G-$analogs of Thurston's 
  earthquake and horocycles flows.

  \bt  \la{AHTT} Under the same assumptions as in Theorem \ref{CDT}:
  
  \begin{enumerate}

\item There is    a $\Gamma_{\G, \bS}-$equivariant \underline{cluster length function}
 \be \la{SLim}
\begin{split}
&{\rm I}(l; e^x):  {\mathscr A}_{\G^\vee, \bS}(\R^t)\times {\mathscr P}_{\G, \bS}(\R_{>0})\lra \R.\\ 
\end{split}
\ee
It is   cluster
   convex, in the sense of Definition \ref{DEFCONV}, in  each of the  variables:
 \be \la{SLime24a}
\begin{split}
& {\rm I}(l_1; e^x) +  {\rm I}(l_2; e^x) \geq   2\cdot {\rm I}((l_1+l_2)/2; e^x).\\
&  {\rm I}(l; e^{x_1}) +  {\rm I}(l; x_2) \geq   2\cdot {\rm I}(l; e^{(x_1+x_2)/2}).\\
\end{split}
  \ee
 
  \item There is  a canonical    $\Gamma_{\G, \bS}-$equivariant \underline{cluster intersection pairing}
 \be \la{SLime}
\begin{split}
&{\cal I}: {\mathscr A}_{\G^\vee, \bS}(\R^t)\times {\mathscr P}_{\G, \bS}(\R^t)\lra \R.\\
\end{split}
  \ee
It is homogeneous of the degree $1$ in the first variable, and 
  cluster convex in  both variables:
 \be \la{SLime24}
\begin{split}
&{\cal I}(al; x) = a {\cal I}(l; x)\ \ \ \ \ \ \forall a\in \R_{>0}.\\
& {\cal I}(l_1; x) +  {\cal I}(l_2; x) \geq   2\cdot {\cal I}((l_1+l_2)/2; x).\\
&  {\cal I}(l; x_1) +  {\cal I}(l;  x_2) \geq   2\cdot {\cal I}(l; (x_1+x_2)/2).\\
\end{split}
  \ee
\vskip 1mm

  \end{enumerate}
  
    \et
  
 \begin{proof}  (1) Let us restrict the elements ${\Bbb I}_l$ of  canonical basis to the positive locus $\{X_i = e^{x_i}>0\}$, and consider the following limit where $l \in  {\mathscr A}_{\G^\vee, \bS}(\Q^t)$ 
 and $Cl\in {\mathscr A}_{\G^\vee, \bS}(\Z^t)$: 
 \be \la{SLimQ}
\begin{split}
& {\rm I}(l;  e^x):={\lim}_{C \to \infty}\frac{\log {\Bbb I}(C\cdot l)( e^x)}{C}.\\ 
&{\rm I}(l; e^x):  {\mathscr A}_{\G^\vee, \bS}(\Q^t)\times {\mathscr P}_{\G, \bS}(\R_{>0})\lra \R.\\
 \end{split}
\ee
Here the limit is over  integers $C= \prod^n_{i=1} p_i^{a_i}$, where $p_i$ are the primes ordered by their value, and $a_i, n \to \infty$.\vskip 2mm

Let us show that this limit exists.

 Thanks to  property (4) in Theorem \ref{CBA},  the leading coefficient    in the product formula (\ref{PRFO}) is 
\be \la{LICO}
 c(l_1, l_2, l_1+l_2)=1.
\ee

 Therefore  we have 
  \be \la{B}
     {\Bbb  I}(l_1; e^x)  {\Bbb  I}(l_2; e^x)  = {\Bbb I}(l_1+l_2; e^x) + \mbox{\rm other terms with positive coefficients}.
\ee  
So  
\be \la{4.23.24}
{\Bbb  I}(l_1; e^x)  {\Bbb  I}(l_2; e^x)  \geq  {\Bbb I}(l_1+l_2; e^x).
\ee
 In particular,  for any integer $m$  
 $$
 {\Bbb I}(ml; e^x)^{1/m} \leq  {\Bbb  I}(l; e^x).
 $$
So the  functions in the first line of (\ref{SLimQ}) is bounded from above, and passing from $C$ to $CC'$ we decrease the function standing under the limit sign. 
On the other hand, all these functions are bounded from above by the leading monomial of ${\Bbb  I}(l; e^x)$, provided by the property (4) of Theorem \ref{CBA} - see also (\ref{A}) below. 
  Therefore  the limit  exists.  Let us show that it is 
a convex function in $x$. 

Let $P(X)$ be a Laurent polynomial in $X = (X_1, ..., X_n)$ with positive  coefficients. Let $X_i = e^{x_i}$ and let $e^{x}:= (e^{x_1}, ..., e^{x_n})$. Then, since $e^x$ is a convex function, 
  $$
  P(e^x) \cdot P(e^y) \geq P(e^{{(x+y)}/2})^2.
  $$
So the function $\log P(e^x)$ is convex in any 
   logarithmic cluster Poisson coordinates $(x_1, ..., x_n)$: 
 $$
 \log P(e^x) + \log P(e^y) \geq 2 \log P(e^{{(x+y)}/2}).
 $$
 It remains to note that a limit of convex functions, which as we have seen exists,  is a convex function.

Since ${\rm I}(l; x)$ is evidently a 
 degree $1$ homogeneous function in $l$, inequality (\ref{4.23.24}) implies the convexity in $l$
 $$
 {\rm I}(l_1; x)  + {\rm I}(l_2; x)  \geq  2 \cdot {\rm I}((l_1+l_2)/2; x).
 $$   
 Therefore the function (\ref{SLimQ}) where $l$ is from the rational tropical locus can be extended by continuity to a function (\ref{SLim}) where $l$ is on the real tropical locus. 
 The obtained function is convex in both variables. 
  \vskip 2mm

 (2)  Property (2) of Theorem \ref{CBA} allows to define the tropicalization of the canonical basis:
  \be \la{SLimeaa}
\begin{split}
&{\cal I}: {\mathscr A}_{\G^\vee, \bS}(\Z^t)\times {\mathscr P}_{\G, \bS}(\R^t)\lra \R.\\
& {\cal I}(l,\bullet):= {\Bbb I}^t(l, \bullet), \ \ \ \   {\cal I}(l,x)= \lim_{C \to \infty}\frac{{\Bbb I}(l)( e^{Cx})}{C}.\\
\end{split}
  \ee 
     The obtained function is evidently convex in the $x-$variable, do to the property (2)  in Theorem \ref{CBA}. 
     
\bp \la{AHTT1}

The    function  ${\cal I}(l,u)$ is degree $1$ homogeneous in $l$. 
   \ep
   
     \begin{proof}   Let us fix a seed, which gives rise to a local chart $(x_1, …, x_n)$  and a scattering diagram. 

\bl
We have for any positive integer $m$
\be \la{A}
 {\Bbb  I}(ml; x_1, x_2,…, x_n) = {\Bbb I}(l; x_1^m, …, x_n^m) + \mbox{\rm other terms with positive coefficients}.
\ee
\el

\begin{proof}  Following Definition 3.1 and 3.3 of\cite{GHKK}   let $(\gamma(t), m_0)$ be a broken line. Then $(\gamma(2t), 2m_0) $
 is also a broken line. Therefore, if ${\rm Mono}(\gamma(t))$ is a term of ${\Bbb I}(l),$ then ${\rm Mono}(\gamma(2t))$ is a term of ${\Bbb  I}(2l)$.

The topicalization of ${\Bbb I}(l; x_1, ..., x_n)^m$ coincides with the topicalization of ${\Bbb I}(l; x_1^2, …, x_n^2)$. 
So  we get 
$$
 m {\Bbb I}^t (l; x_1, ..., x_n)\stackrel{(\ref{B})}{\geq} {\Bbb I}^t (ml; x_1, ..., x_n) \stackrel{(\ref{A})}{\geq}   {\Bbb I}^t(l; x_1^m, …, x_n^m) = m {\rm I}^t(l; x_1, ..., x_n).
 $$
This implies the desired equality:
$
m {\Bbb I}^t (l; x)= {\Bbb I}^t (ml; x). 
$
\end{proof}

       Therefore it can be extended by $l-$homogenuity to a function
   \be \la{SLimb}
 {\cal I}(l;x): {\mathscr A}_{\G^\vee, \bS}(\Q^t) \times 
  {\mathscr P}_{\G, \bS}(\R_{>0}) \lra \R. 
  \ee    

  The function  ${\rm I}(l;x)$ is convex in the $l$-variable.   
  Therefore it can be extended by continuity to a function 
    \be \la{SLimc}
 {\cal I}(l;x): {\mathscr A}_{\G^\vee, \bS}(\R^t) \times 
  {\mathscr P}_{\G, \bS}(\R_{>0}) \lra \R 
  \ee   
  
Proposition \ref{AHTT1} is proved. \end{proof}
Theorem \ref{AHTT} is proved.   \end{proof}

  \subsubsection{Cluster convexity in the part (1) of Theorem \ref{AHTT} and Kirckhoff's theorem.} \la{2.4.2a} When $\G = {\rm PGL}_2$, the function  ${\rm I}(l,u)$ is    the length of an integral geodesic lamination. 
 Indeed,  in this case  the canonical basis ${\Bbb I}(l)$ is parametrised by the  multicurves $l=\sum_j n_j\alpha_j$ on $\bS$, where $\{\alpha_j\}$ is a collection of non-intersecting simple loops 
 on $\bS$,  $n_j \in \Z$, and  $n_j >0$ if $\alpha_j$ is not a boundary loop  \cite[Section 12]{FG03a}. 
So ${\rm I}(l,u) $ is the length of the multicurve in the hyperbolic metric for the point $u$ of the Teichmuller space of $\bS$. 

The cluster convexity of the length function ${\rm I}(l,u) $ is a variant of Kerckhoff's theorem \cite{Ke}. 
Indeed, the cluster convexity  tells that the  function  
${\rm I}(l; u)$, as a function in $u$,  is convex in any logarithmic cluster coordinate system. Kerckhoff's theorem tells  that, 
for any Fenchel-Nielsen coordinate system, the 
 length function  ${\rm I}(l; u)$ is convex on the tajectories of the twist flows for this coordinate system.  Precisely, for a generic 
 multicurve $m = \sum_jm_j\beta_j$, the loops $\{\beta_j\}$  determine a pair of pants decomposition of the surface, and hence the corresponding Fenchel-Nielsen 
 coordinate system. The  trajectories of the  Hamiltonians  ${\rm I}(m,u) $ are the twist flows for 
 this Fenchel-Nielsen coordinate system.  \vskip 2mm

 \bcon \la{KCON} For any pair $l, m \in {\mathscr A}_{\G^\vee, \bS}(\R^t)$, the cluster length function  ${\rm I}(l; x)$ is convex on the  trajectories of the Hamiltonian ${\rm I}(m; x)$. 
\econ 

Conjecture \ref{KCON} for ${\rm PGL}_2$  is a  reformulation of Kerckhoff's theorem  via the canonical bases perspective. 
\vskip 2mm

Theorem \ref{4.20.24.1} implies immediately the following. 
  \bt In the assumptions of Theorem \ref{4.20.24.1}:
  
 i) The function $\log {\rm Tr}{\rm Mon}_\alpha({\cal L}_V)$ is  cluster convex  function 
 on ${\mathscr P}_{\G, \bS}(\R_{>0})$.
 
 ii) The following limit exists, providing   a cluster convex function on  ${\mathscr P}_{\G, \bS}(\R_{>0})$:
 $$
 {\rm I}(\alpha, V;  e^x):={\lim}_{n \to \infty}\frac{\log {\rm Tr}{\rm Mon}^n_\alpha({\cal L}_V)( e^x)}{n}.
 $$  
It is called the cluster length function 
 of the pair $(\alpha, V)$. 
 \et

 \subsubsection{The earthquake map.} Since the cluster length function ${\rm I}(l,x)$  is convex in the $x-$variable, it is continues, and  differentiable everywhere except  on a countable subset of hypersurfaces. 
 We can define the Hamiltonian flow for such a function. 
 The time 1 Hamiltonian flow  for the Hamiltonian    ${\rm I}(l,x)$ on the higher Teichm\"uller space ${\mathscr P}_{\G, \bS}(\R_{>0})$ is the 
 analog of the earthquake map
 $$
{\cal E}_{\G, \bS}:  {\mathscr A}_{\G^\vee, \bS}(\R^t)\times  {\mathscr P}_{\G, \bS}(\R_{>0}) \lra  {\mathscr P}_{\G, \bS}(\R_{>0}). $$

 \bcon \la{Earthq} For any point $x \in  {\mathscr P}_{\G, \bS}(\R_{>0}) $, the earthquake map ${\cal E}_{\G, \bS}$ provides an isomorphism
 $$
 {\cal E}^x_{\G, \bS}:  {\mathscr A}_{\G^\vee, \bS}(\R^t)\times  \{x\} \stackrel{\sim}{\lra}  {\mathscr P}_{\G, \bS}(\R_{>0}).  $$
 \econ 
 
 For the group $\G={\rm PGL}_2$ this is  Thurston's theorem    in the Teichmuller theory.

For any multicurve $\sum n_i \alpha_i$, with $n_i>0$, colored by representations $V_i$ of $\G$, there is the cluster  length function 
$
\sum_i  n_i {\rm I}(\alpha_i, V_i;  e^x). 
$ 
It is therefore a cluster convex function on ${\mathscr P}_{\G, \bS}(\R_{>0})$.  However if ${\rm rk} \G>1$, 
the supply of  functions obtained by their limits    is  much smaller than  needed  for 
Earthquake Conjecture \ref{Earthq}. 

  \subsubsection{The horocycle flow} Given a point $l \in  {\mathscr A}_{\G^\vee, \bS}(\R^t)$, the cluster intersection pairing provides  the function ${\cal I}(l, \cdot)$ on  the real Poisson tropical space 
${\mathscr P}_{\G, \bS}(\R^t)$.  The cluster intersection pairing ${\cal I}(l,m)$  is convex in the $x-$variable, and therefore it is continues, and  differentiable everywhere except  on a countable subset of hypersurfaces. 
 We can view this function as  a Hamiltonian. Its Hamiltonian flow is the $\G-$analog of the horocycle flow. Indeed, when $\bS$ is a closed surface and $\G = {\rm PGL}_2$, 
 this is a theorem of Papadopoulos \cite{Pa}. \vskip 2mm

\subsection{Quantization of moduli spaces of $\G$-local systems on decorated surfaces} \la{1.1}


 \medskip

\subsubsection{The result.}  \la{5.1.1}  
Given a     topological  surface $S$ and a group $\G$, the  {\it character variety}   ${\rm Loc}_{\G, S}$ 
parametrizes homomorphisms    $\pi_1(S) \to \G$,  considered modulo  $\G$-conjugation. 
 It is    the Betti version of the moduli space of $\G$-local systems on $S$.  
The mapping class group $\Gamma_S= {\rm Diff}(S) / {\rm Diff}_0(S)$ of $S$ acts  on ${\rm Loc}_{\G, S}$.  \vskip 1mm
 
 Now let $\G$ be a split reductive algebraic group over $\Q$ with connected center. 
 We define a generalization ${\rm Loc}_{\G, \bS}$ of the character variety for any colored decorated surface $\bS$ as the moduli space parametrising the data 
 obtained by forgetting the flags at the punctures (but not the decorated flags at the special points) in   Definition \ref{DEFP2}. Equivalently, it is just the moduli space ${\mathscr P}_{\G, {\bf S}}$ 
 from Definition \ref{DEFP2*} obtained by viewing $\bS$ as a fully colored decorated surface $\bf S$. Here are two examples. \vskip 1mm
 
\begin{enumerate}

\item  If $\bS = {\Bbb D}_m$ is a disc with $m$  boundary points and all boundary intervals colored,    we show in Lemma \ref{EQQ} that  we get   
the configuration space of $m$  decorated flags: 
 \be
 {\rm Loc}_{\G, {\Bbb D}_m} = {\rm Conf}_m({\cal A}_\G):= \G \backslash (\G/\U)^m. 
 \ee
 
\vskip 1mm

\item  If $\bS=\odot$ is the punctured disc with two special boundary points, 
and two colored sides, the moduli space  ${\rm Loc}_{\G, \odot}$ parametrises $\G-$local systems on the punctured disc 
 with reduction to $\B$ at the two boundary points. \vskip 2mm
 \end{enumerate}
 
 Let $n$ be the number of punctures on $\bS$.  We will quantize the  space ${\rm Loc}_{\G, \bS}$ in the following sense.   
\begin{itemize} 
\item Using Theorem \ref{MTH}, and in particular the cluster Poisson nature of the action of the  group $W^n$,  
we  define a non-commutative algebra  ${\cal O}_q({\rm Loc}_{\G, \bS})$   
   deforming  the algebra   ${\cal O}({\rm Loc}_{\G, \bS})$ 
   by setting 
 \be \la{AHLQ}
 {\cal O}_q({\rm Loc}_{\G, \bS}):= {\cal O}_{q}({\mathcal P}_{\G, \bS})^{W^n}, \qquad q\in \C^*.
\ee
 
\vskip 1mm

\item We consider the {\it Langlands modular double} of the algebra  ${\cal O}_q({\rm Loc}_{\G, \bS})$, defined as the algebra
 \be \la{AH}
 \mathcal{A}_\hbar({\rm Loc}_{\G, \bS}):= {\cal O}_q({\rm Loc}_{\G, \bS})\otimes_\C{\cal O}_{q^\vee}({\rm Loc}_{\G_{\rm ad}^\vee, \bS}), \qquad q= e^{i\pi \hbar}, \qquad q^\vee:= e^{i \pi / \hbar}.
\ee
Here  $ \G_{\rm ad}^\vee$ is the minimal quotient for the Langlands dual group with connected center.    

Note that  the Langlands modular duality  enters to the construction by the very definition:
$$
  \mathcal{A}_\hbar({\rm Loc}_{\G, \bS}) = {\mathcal A}_{\hbar^\vee}({\rm Loc}_{\G_{\rm ad}^\vee, \bS}), \qquad  \hbar^\vee:= 1/\hbar.
$$

\vskip 1mm 

\item We define an involutive antiholomorphic antiautomorphism   $\ast$ of  algebra $ \mathcal{A}_\hbar({\rm Loc}_{\G, \bS})$,  assuming 
\be \la{Condh}
\hbar+\hbar^{-1} +2 \geq 0\qquad \Longleftrightarrow \qquad \hbar \in \R_{>0},  \qquad\mbox{or}   \qquad|\hbar|=1. 
\ee

\item  We construct  a $\Gamma_{\G, \bS}$-equivariant  representation of the $\ast$-algebra  $ \mathcal{A}_\hbar({\rm Loc}_{\G, \bS})$,\footnote{We assume that $\G$ is simply-laced if $|\hbar|=1$.} given by:

\begin{itemize}
\vskip 1mm 

\item[i)]  a unitary   representation $ \gamma \mapsto \rho_\gamma$ of a central extension $\widehat \Gamma_{\G, \bS}$ of  $\Gamma_{\G, \bS}$ 
 in a Hilbert space  ${\cal H}_{\G, S}$; 

\vskip 1mm

\item[ii)] a   representation $\circ$ of the $\ast$-algebra  $\mathcal{A}_\hbar({\rm Loc}_{\G, S})$ by unbounded operators in a dense subspace ${\cal S}_{\G, S}$ of the   Hilbert space ${\cal H}_{\G, S}$ with  an inner product  $\langle \cdot, \cdot \rangle$, such that: 
  \be  \la{3.4ax}
\begin{split}
&{\rho}_\gamma(A \circ s) = \gamma(A)\circ \rho_\gamma(s), \qquad   \forall \gamma \in \widehat \Gamma_{\G, \bS},\qquad \forall A \in \mathcal{A}_\hbar({\rm Loc}_{\G, S}),\qquad  \forall s \in {\cal S}.\\
&\langle A s_1, s_2\rangle = \langle  s_1, \ast As_2\rangle. \\
\end{split}
\ee 
\end{itemize}
\end{itemize}

 Let us explain how all this   looks in the simplest, and yet fundamental example, when  $\G = {\Bbb G}_m$.
  
\subsubsection{An example: $\G = {\Bbb G}_m$.}      The space  ${\rm Loc}_{{\mathbb G}_m, S} = {\rm Hom}(H_1(S, \Z), {\Bbb G}_m)$ is a split algebraic torus. The intersection form on the lattice $\Lambda:= H_1(S, \Z)$ gives rise to the quantum torus algebra ${\cal O}_q({\rm T}_\Lambda)$, see Section \ref{2.1.5}. 
 The   mapping class group $\Gamma_S$ acts    by  automorphisms  of 
the lattice $H_1(S, \Z)$ and preserves the  intersection form.  
Consider the {\it modular double} of the quantum torus algebra: 
$$
{\cal A}_{\hbar}({\rm T}_{\Lambda}):= 
{\cal O}_{q}({\rm T}_{\Lambda})\otimes_\C {\cal O}_{q^\vee}({\rm T}_{\Lambda}), \qquad  q=e^{i\pi \hbar}, \qquad q^\vee=e^{i\pi/ \hbar}, \qquad\hbar \in \C^*.
$$

Denote by $Y_\lambda$ the  generators of the algebra ${\cal O}_{q^\vee}({\rm T}_{\Lambda})$.

\bd \la{DEFINV}For any    Planck constant $\hbar$ such that $\hbar \in\R$ or $|\hbar|=1$, the algebra ${\mathcal A}_{\hbar}({\rm T}_{\Lambda})$ 
has an antiholomorphic    involutive antiautomorphism $\ast$, 
 defined as follows:
\be \la{INVRUa}
\begin{split}
\ast_{\rm R}, ~\hbar \in \R_{}:\qquad&\ast_{\rm R} (X_\lambda) = X_\lambda, \qquad\ast_{\rm R} (Y_\lambda) = Y_\lambda, \qquad\ast_{\rm R}(q) = q^{-1}, \qquad \ast_{\rm R}(q^\vee) = {q^\vee}^{-1},\\
~\ast_{\lr} , ~|\hbar|=1: \qquad&\ast_{\lr} (X_\lambda) = Y_\lambda, \qquad \ast_{\lr} (Y_\lambda) = X_\lambda, \qquad \ast_{\lr}(q) = {q^\vee}^{-1}, \qquad \ast_{\lr}(q^\vee) = q^{-1}.\\
\end{split}
\ee 
\ed

The involution  $\ast_{\rm R}$   defines a real form  of the quantum torus algebra ${\cal O}_{q}({\rm T}_{\Lambda})$. 

The involution $\ast_{\lr}$ flips the two factors of ${\cal A}_{\hbar}({\rm T}_{\Lambda})$. It does not act on ${\cal O}_{q}({\rm T}_{\Lambda})$.

 \vskip 2mm
 
  The group ${\rm Sp}(\Lambda)$ of automorphisms of the lattice with the form 
  acts by  automorphisms of    algebras ${\cal O}_q({\rm T}_{\Lambda})$ and ${\cal A}_{\hbar}({\rm T}_{\Lambda})$.
Let $\widetilde {\rm Sp}(\Lambda)$ be its metaplectic cover. 
  Using the Weil representation of the 
metaplectic  group,  we get a $\widetilde {\rm Sp}(\Lambda)$-equivariant 
quantization of the torus  ${\rm T}_{\Lambda}$. 

\bt \la{Weil}  
 For any $\hbar$ such that $\hbar \in\R$ or $|\hbar|=1$, there is a  representation $ \circ  $  of 
the $\ast$-algebra ${\cal A}_{\hbar}({\rm T}_{\Lambda})$ 
in a dense   subspace ${\cal S}$ of a Hilbert space ${\cal H}$. It is equivariant under 
 a unitary representation $ \gamma \mapsto \rho_\gamma$ of the group the   group $\widetilde {\rm Sp}(\Lambda)$  in the Hilbert space  space ${\cal H}$:   
 \be  \la{Weil1}
\rho_\gamma(a\circ s) = 
\gamma(a)\circ \rho_\gamma(s)\qquad  
 \forall \gamma \in \widetilde {\rm Sp}(\Lambda), ~ \forall a \in {\cal A}_{\hbar}({\rm T}_{\Lambda}), ~ \forall s\in {\cal S}. 
\ee
\et
The space ${\cal S} $ is the maximal subspace in ${\cal H} $ where the algebra ${\cal A}_{\hbar}({\rm T}_{\Lambda})$ acts.
There is a   Freschet space structure 
on ${\cal S} $ for which the algebra ${\cal A}_{\hbar}({\rm T}_{\Lambda})$ acts by   continuous operators. Theorem \ref{Weil} is proved in Section \ref{SEC8.4}. 
In the case when $\hbar \in\R$ it was proved in \cite[Theorem 7.7]{FG07}.

It is  instructive to elaborate Theorem \ref{Weil} in the simplest   and at the same time most basic example.

\subsubsection{Basic example.} 
Let $\Lambda = \Z^2, \langle e_1, e_2\rangle =1$. Set 
$\hbar = \beta^2$. 
Let ${\cal H}:= L^2(\R)$ be the   Hilbert 
space  of functions $f(t)$. Consider the shift  operator $T_z(f)(t) := f(t+z)$. 
There are the following 
  operators acting on functions $f(t)$: 
\be \la{1.18.09.2}
\begin{split}
&\qquad X_1 := T_{2\pi i \beta},   \qquad  
 Y_1 := T_{2\pi i /\beta},\\
& \qquad  X_2 :=e^{\beta t}, \qquad  
 Y_2 :=e^{t/\beta}.\\
\end{split}
\ee
Let ${\cal S}\subset L^2(\R)$  be    the maximal subspace on which any polynomial of   operators  (\ref{1.18.09.2}) is defined, and its image lies in $L_2(\R)$. 
It consists of all functions $f(t)$ 
which for any real $c$ decay faster then $e^{ct}$, 
and whose Fourier transform has the same property.  The space ${\cal S}$  contains   functions  
$$
e^{-t^2+at+b}P(t), \qquad  \mbox{where $P(t)$ is an a polynomial and $a, b\in \C$}. 
$$
The operators  (\ref{1.18.09.2}) acting in ${\cal S}$ 
satisfy the following commutation relations 
$$
q^{-1} X_1 X_2 = q X_2 X_1, \qquad  
{q^\vee}^{-1} Y_1 Y_2 = 
q^\vee Y_2 Y_1, \qquad  
[ X_i,  Y_j]=0.
$$

This means that we get a representation of the  algebra ${\cal A}_{\hbar}({\rm T}_{\Lambda})$.   
One  checks that for any $\hbar$ such that $\hbar >0$  or $|\hbar|=1$ it is a  representation of the $\ast$-algebra ${\cal A}_{\hbar}({\rm T}_{\Lambda})$  for   the involution  (\ref{INVRUa}).

The unitary projective representation  
of ${\rm SL_2}(\Z)$ is obtained by the restriction of the Weil representation for the metaplectic cover of ${\rm SL}_2(\R)$, see \cite{LV}. It acts by operators
$\rho_\gamma$ given by  kernels $K_\gamma$:
$$
f(x) \lms \int_{\R}K_\gamma(x,y)f(x)dx, \qquad  \gamma \in {\rm SL_2}(\Z).
$$
For the  generators $T= \begin{pmatrix}  
      1 & 1 \\
      0  & 1 \\
 \end{pmatrix}$ and 
$S=  \begin{pmatrix}  
      0 & 1 \\
      -1 & 0 \\
 \end{pmatrix}$ of ${\rm SL_2}(\Z)$, we define the kernels by setting 
$$
K_T(x,y):=  e^{-x^2/4\pi i  }\delta(x-y), \qquad  
K_S(x,y):= \frac{1}{2\pi i}e^{xy/2\pi i}. 
$$  

The kernel 
$K_\gamma(x,y)$ for $\beta \not \in \Q$  is determined uniquely up to a constant 
by the system of difference equations arising from the equivariance condition (\ref{Weil1}). 
For example, for the kernel $K_S$
\be \nonumber
\begin{split} 
&K_S(x+2\pi i \beta,y) = e^{\beta y}K_S(x,y), \qquad  
K_S(x+2\pi i/\beta,y) = e^{y/ \beta}K_S(x,y)\\
&K_S(x,y+2\pi i \beta) = e^{\beta x}K_S(x,y), \qquad  
K_S(x,y+2\pi i/\beta) = e^{x/ \beta}K_S(x,y).\\
\end{split}
\ee
 
Although one of the factors of  the algebra $\mathcal{A}_\hbar({\rm T}_{\Lambda})$ does not have a limit when $\hbar \to 0$ or $\hbar \to \infty$, the kernels $K_\gamma(x,y)$ do not depend on $\hbar$, 
and define the metaplectic representation. 

\vskip 2mm
Recall the canonical $\ast$-representation of the Heisenberg algebra with the Planck constant $1$:
\be \la{}
\begin{split}
& p_1 := 2\pi i  \cdot \frac{\partial}{\partial t},  \qquad     p_2: =   t, 
\qquad  [ p_1,  p_2] = 2\pi i, \qquad  \ast p_1 = p_1 \qquad  \ast p_2 = p_2.
\end{split}
\ee
Rescaling it by $\beta$ and   $\beta^{-1}$, we get 
 two Heisenberg algebras   given 
by   logarithms of operators (\ref{1.18.09.2}):
\be \nonumber
\begin{split}
& x_1 := \beta p_1 = 2\pi i \beta\cdot \frac{\partial}{\partial t},  \qquad  x_2: = \beta p_2 = \beta t, 
\qquad  [ x_1,  x_2] = 2\pi i \beta^2, \\
& y_1 := \beta^{-1}p_1 = 2\pi i \beta^{-1} \cdot \frac{\partial}{\partial t},   \qquad  y_2: =   \beta^{-1} p_2 = \beta^{-1}t, 
\qquad  [ y_1,  y_2] = 2\pi i \beta^{-2}.\\
\end{split}
\ee

\subsubsection{The general quantum torus case.} Recall the Heisenberg Lie algebra ${\it Heis}_\Lambda$ assigned to a lattice $\Lambda$ with a form $\langle \ast, \ast\rangle$. 
It is a central extension of the abelian Lie algebra $\Lambda$ by $\Z$:
$$
0 \lra  \Z \lra {\it Heis}_\Lambda \lra \Lambda \lra 0.
$$
The commutator of  any   elements $\widetilde p_1, \widetilde p_2$ projecting to   $p_1, p_2 \in \Lambda$ is given by 
$[\widetilde p_1, \widetilde p_2]:=   \langle p_1, p_2\rangle$. \vskip 1mm
 
Suppose first that the form $\langle \ast, \ast\rangle$ is non-degenerate. Then there is the Heisenberg  representation of the 
Lie algebra ${\it Heis}_\Lambda$ with ${\rm 1} \in \Z \lms 2\pi i ~{\rm Id}$, and elements $p \in \Lambda$ acting by    self-adjoint operators $\widehat p$.    The exponents of the rescaled by $\beta$ and $\beta^{-1}$    operators $\widehat p$   generate   two commuting quantum torus algebras, 
providing a representation of the  algebra $\mathcal{A}_\hbar({\rm T}_{\Lambda})$: 
$$
X_p := e^{\beta \widehat p}, \qquad  Y_q := e^{\beta^{-1} \widehat q}, \qquad [X_p, Y_q]=0.
$$
If $\hbar$ satisfies condition (\ref{Condh}), we get a representation of the $\ast$-algebra $\mathcal{A}_\hbar({\rm T}_{\Lambda})$ from Definition \ref{DEFINV}. 

Choosing a decomposition $\Lambda = L \oplus L'$ into a sum of two Lagrangian subspaces we get a  realization of the representation of the 
Lie algebra ${\it Heis}_\Lambda$ in the Hilbert space   of semiforms on the vector space $L \otimes \R$. 
Namely, an element $q \in L$ acts as  a vector field $ 2\pi i \partial_q$, and an element $p\in L'$ acts as the operator of multiplication by the linear function 
$\langle p, \ast\rangle$.

 Hilbert spaces assigned to different  decompositions of $\Lambda$  are related by unitary intertwiners, 
providing a unitary projective representation of the metaplectic group $\widetilde {\rm Sp}(\Lambda)$, known as the Weil representation. \vskip 1mm

Now let $\Lambda$ be  lattice with an arbitrary form $\langle\ast, \ast\rangle$. 
Let $\Lambda_0\subset \Lambda$ be the kernel of the form.   Then there is a canonical map, whose fibers are   
the symplectic leaves of the Poisson structure on ${\rm T}_{\Lambda}$:
 $$
\mu: {\rm T}_{\Lambda} \to {\rm T}_{\Lambda_0}.
$$ 
The Hilbert space 
${\cal H}_\Lambda$ assigned to $\Lambda$ is the space of semiforms on   ${\rm T}_{\Lambda_0}(\R_{>0})$ 
  with values in    ${\cal H}_{\Lambda/\Lambda_0}$.
  
The subalgebra $\mu^*{\cal O}({\rm T}_{\Lambda_0})$   coincides with  the center of the algebra ${\cal O}_q({\rm T}_{\Lambda})$ unless  $\hbar  \in \Q$. 
If $\hbar  \in \Q$, the center is much bigger. So one can   decompose  ${\cal H}_\Lambda$   according to the characters  of the center. The characters $\lambda$ we get this way  are  paramatrized by the points 
   of the group $ {\rm T}_{\Lambda_0}(\R_{>0})$:
$$
{\cal H}_\Lambda= 
\int {\cal H}_{\lambda}d\lambda, \qquad \lambda \in {\rm T}_{\Lambda_0}(\R_{>0}).
$$

\subsubsection{Our strategy for an arbitrary $\G$.} 
One could hope to  find  a birational Poisson isomorphism between   ${\rm Loc}_{\G, \bS}$ and a Poisson torus 
${\rm T}_\Lambda$, and  then quantize  ${\rm T}_\Lambda$.  
However this idea faces two  serious issues:\vskip 1mm

1. In general  such a birational isomorphism    does not exist: 
the space   ${\rm Loc}_{\G, \bS}$    is not rational. 

\vskip 1mm 2. Although the mapping class group  $\Gamma_\bS$ acts on ${\rm Loc}_{{\Bbb G}_m, \bS}$ by monomial transformations, its   action  
  on  ${\rm Loc}_{\G, \bS}$ is  non-linear  if $\G \not = {\Bbb G}_m$. 
So given one  coordinate system, the action of   $\Gamma_\bS$  generates infinitely  many others, and we must 
relate their quantizations  to produce a $\Gamma_\bS$-equivariant quantization. This   requires a non-linear analog of the Weil representation. 

\vskip 2mm
We handle these   issues  using the following strategy. \vskip 1mm

1. We consider the moduli space ${\mathscr P}_{\G, \bS}$, which is birationally isomorphic to a torus.   
Forgetting  flags at  the punctures, we get  a   surjective map, which is a Galois $W^n-$cover over the generic point:
 \be \la{PRW}
\pi: {\mathscr P}_{\G, \bS} \lra {\rm Loc}_{\G, \bS}.
\ee

2.  By Theorem \ref{MTH}, the  space ${\mathscr P}_{\G, S}$ carries a $\Gamma_{\G, \bS}$-equivariant cluster Poisson structure.  
This allows to apply cluster Poisson quantization to quantize  the space ${\mathscr P}_{\G, S}$. Since $W^n \subset \Gamma_{\G, \bS}$, the projection   (\ref{PRW}) has cluster Poisson nature. 
This allows to  quantize  
  the  space ${\rm Loc}_{\G, S}$.

 \subsubsection{Cluster Poisson quantization.}  \la{2.4.6}  Given a quiver, its {\it Langlands dual quiver} is defined by changing the skew-symmetrizable form $(*,*)$ by setting 
 $(x, y):= -(y, x)$, and keeping the rest of the quiver data intact. Given a cluster Poisson space ${\mathscr X}$,  we denote by ${\mathscr X}^\vee$ the   
{\it  Langlands dual cluster Poisson variety}. 
Recall  {\it quantum modular double algebra } of ${\mathscr X}$: 
\be \la{CLAlg}
\mathcal{A}_\hbar({\mathscr X}):= {\cal O}_q({\mathscr X}) \otimes_\Z {\cal O}_{q^\vee}({\mathscr X}^\vee), \qquad  q=e^{i\pi \hbar}, \qquad q^\vee=e^{i\pi/ \hbar}, \qquad\hbar \in \C^*.
\ee
Using the map (\ref{astA1}) we get for each quiver ${\bf q}$ of the cluster atlas an embedding  of  algebras 
\be \la{astA1q}
\begin{split}
&\alpha^*_{\bf q}: \mathcal{A}_\hbar({\mathscr X})  \subset \mathcal{A}_\hbar({\rm T}_{\bf q}).\\
\end{split}
\ee

 The next two theorems describe   the cluster Poisson quantization. 
 \vskip 2mm
       
 According to Definition \ref{DEFINV}, any Planck constant $\hbar$  satisfying condition (\ref{Condh}) gives rise to an   involution $\ast$ of the algebra 
 $\mathcal{A}_\hbar({\rm T}_{\bf q})$. For each quiver ${\bf q}$ from the cluster atlas the embedding (\ref{astA1q}) induces a $\ast$-algebra structure 
 on     $\mathcal{A}_\hbar({\mathscr X})$. Theorem \ref{INVRU} tells that it does not depend on   ${\bf q}$ if $\hbar$ satisfies condition (\ref{Condh}).
  
   \bt \la{INVRU} For  any   $\hbar$   satisfying condition (\ref{Condh}) there is a unique involution $\ast$ of the algebra $\mathcal{A}_\hbar({\mathscr X})$ 
   such that for any quiver ${\bf q}$ from the cluster atlas the map (\ref{astA1q}) is a map of $\ast$-algebras.   
   \et

If  $\hbar \in \R$, the involution $\ast$ acts on each  factor in (\ref{CLAlg}).

In the   $|\hbar |=1$ case, the involution $\ast$ interchanges the   factors of the algebra $\mathcal{A}_\hbar({\mathscr X})$ in (\ref{CLAlg}). \vskip 1mm

For any cluster Poisson variety ${\mathscr X}$, there is a split algebraic torus ${\rm C}_{\mathscr X}$, 
called the {\it cluster Casimir torus},  defined in  \cite[Section 2.2]{FG03b},\footnote{We usually omit the adjective ``cluster".}  and a canonical surjective map 
\be \la{CLTX}
\mu:{\mathscr X} \lra {\rm C}_{\mathscr X}.
\ee
The subalgebra $\mu^*{\cal O}({\rm C}_{\mathscr X})$ is the center of the Poisson algebra ${\cal O}({\mathscr X})$   \cite[Section 3]{FG03b}.\vskip 2mm

For any cluster Poisson variety ${\mathscr X}$, the center of the algebra ${\cal O}_q({\mathscr X})$  contains the algebra of 
regular functions on the Cluster casimir torus  ${\rm C}_{{\mathscr X}}$. 
We prove in Section \ref{SECT20.4} that if $q$ is not a root of unity, then, under some mild condition on ${\mathscr X}$, which is satisfied if the Duality Conjectures are valid for ${\mathscr X}$, 
they coincide: 
\be \la{CCasT2}
{\rm Center}\ {\cal O}_q({\mathscr X}) = {\cal O}({\rm C}_{{\mathscr X}}) \ \ \ \ \ \ \mbox{if $q$ is not a root of unity, and ${\cal O}_q({\mathscr X})$ is ``large"}.
\ee


 Theorem \ref{MTHCQ} describes  the cluster Poisson quantization. 
      
   \bt \la{MTHCQ} Let ${\mathscr X}$ be a  cluster Poisson variety, and $\hbar$   a  Planck constant     satisfying condition (\ref{Condh}).  
    Then the    $\ast$-algebra $\mathcal{A}_\hbar({\mathscr X})$    from Theorem \ref{INVRU} has    a     series  of $\Gamma_{\mathscr X}$-equivariant representations,      
parametrized by   
  $\lambda \in {\rm C}_{\mathscr X}(\R_{>0})$,   
 and involving the following data.\footnote{We  assume that if $|\hbar|=1$ then ${\mathscr X}$ is simply-laced, i.e. ${\mathscr X} = {\mathscr X}^\vee$.} 
  \vskip 2mm
  
  i) For   each quiver ${\bf q}$: 
  
   \begin{itemize}
  
 \vskip 1mm  \item  A Hilbert space ${\cal H}_{{\bf q}, \lambda}$, assigned to   the   torus   $ {\rm T}_{\bf q}$ by Theorem \ref{Weil}.  
  
 \vskip 1mm \item A dense  subspace ${\cal S}_{{\bf q}, \lambda}\subset {\cal H}_{{\bf q}, \lambda}$, equipped with a  structure of a Frechet topological vector space.
  
 \vskip 1mm \item A  continuous representation $\rho_{{\bf q} }$ of the   $\ast$-algebra   $\mathcal{A}_\hbar({\mathscr X})$ in   ${\cal S}_{{\bf q}, \lambda}$.

 \end{itemize}
 \vskip 1mm
 
  ii)  For each cluster mutation ${\bf q}_1 \to {\bf q}_2$:  
 
  \begin{itemize}
  
 \vskip 1mm \item
  
  a   unitary operator ${\rm I}_{{\bf q}_1, {\bf q}_2}$ 
 intertwining\footnote{Note that  the second line in (\ref{3.4a})  makes  sense     only for $A \in \mathcal{A}_\hbar({\mathscr X}) $.}   representations 
 $\rho_{{\bf q}_1 }$ and $\rho_{{\bf q}_2 }$, and identifying   ${\cal S}_{{\bf q}, \lambda}$:    
  \be \la{3.4a}
\begin{split}
&{\cal I}_{{\bf q}_1\to  {\bf q}_2}: {\cal H}_{{\bf q}_1, \lambda} \stackrel{\sim}{\lra}  {\cal H}_{{\bf q}_2, \lambda},\qquad 
{\cal I}_{{\bf q}_1\to {\bf q}_2}: {\cal S}_{{\bf q}_1, \lambda} \stackrel{\sim}{\lra}  {\cal S}_{{\bf q}_2, \lambda},\\
&{\cal I}_{{\bf q}_1\to {\bf q}_2}(A) \circ  \rho_{\bf c_1}(s) = \rho_{\bf c_2} \circ \mu_{{\bf q}_1 \to {\bf q}_2}(A)(s), \qquad \forall A \in \mathcal{A}_\hbar({\mathscr X}),
\qquad \forall s \in {\cal S}_{{\bf q}, \lambda}.\\
\end{split}
\ee 
  \end{itemize}
   \vskip 1mm
   
   iii) Denote by ${\rm H}{\rm ilb}$ the category of Hilbert spaces with  morphisms given by unitary  operators, defined  up to a unitary scalar. 
Intertwiners (\ref{3.4a}) define  a functor from the cluster modular groupoid ${\rm M}_{\mathscr X}$:
   \be \la{3.4aac}
\begin{split}
 & {\cal Q}:  {\rm M}_{\mathscr X} \lra {\rm H}{\rm ilb}, \\
 &{\cal Q}({\bf q}) := {\cal H}_{{\bf q}, \lambda}, \qquad {\cal Q}({{\bf q}_1 \to {\bf q}_2}) := {\cal I}_{{\bf q}_1\to {\bf q}_2}\\
 \end{split}
\ee 
 Since $\pi_1({\rm M}_{\mathscr X})=\Gamma_{\cal X}$,  we get a unitary projective representation of   
the   group $ \Gamma_{\mathscr X}$ in the  space  ${\cal H}_{{\bf q}, \lambda}$. 
  \et

Theorem \ref{MTHCQ}  for $\hbar >0$ was  proved in \cite{FG07}.  
In   Section \ref{SEC8.4} we do the      
 $|\hbar|=1$ case.

\subsubsection{Quantization of the   space $ {\mathscr P}_{\G, \bS}$.}  \la{sqmsls} 

Since the space ${\mathscr P}_{\G, \bS}$ has a cluster Poisson structure,   
we have its cluster Casimir torus, denoted by ${\rm C}_{\mathscr{P}_{\G, \bS}}$. Our next goal is to explain how it looks like. 

\vskip 2mm  
 
      Let $\H_*$ be the coinvariants of the involution $\ast: h \lms w_0(h^{-1})$ of the Cartan group $\H$ of $\G$.  
  Set  
  \be
\la{BR**}
\begin{split}
&\H_{(\pi)} =   \left\{  \begin{array}{ll}    \H_*   & \mbox{if   $d_\pi$   is odd,}  \\
      \H & \mbox{if   $d_\pi$   is even.} \\
   \end{array}\right.\\
   \end{split}
\ee 
  
 \bd \la{TORID}  
 Let $\bS$ be a colored decorated surface without non-colored boundary intervals. 
 
 The  Casimir torus  ${\rm C}_{\G, \bS}$ is a split torus,  given by
  the product of the Cartan groups $\H$ over the punctures and    components $\pi$   with  even $d_\pi$,  and the groups $\H_\ast$ over  
 components $\pi$   with odd $d_\pi$:
    \be \la{DEF1.15a}
  {\rm C}_{\G, \bS}:=    \H^{ \{\mbox{\rm punctures}\}} \times \H^{\{\mbox{\rm $\pi~|~d_\pi$ is even}\}}
 \times  \H_\ast^{\{\mbox{\rm     $\pi~|~d_\pi$ is odd}\}}.   
 \ee
  \ed 
  
As explained  in Theorem \ref{TH1.16a} below, for the cluster Poisson variety describing  the  
cluster Poisson structure on the space ${\mathscr P}_{\G, \bS}$,  the corresponding cluster Casimir 
torus ${\rm C}_{\mathscr {P}_{\G, \bS}}$  is  related to the Casimir torus $ {\rm C}_{\G, \bS}$ from (\ref{DEF1.15a}) via an isogeny 
$i_{\rm C}: {\rm C}_{\mathscr X}\lra  {\rm C}_{\G, \bS}$, which induces an isomorphism on the real poisitive points. 
Furthermore, the projection $\mu$ from (\ref{CLTX})
 agrees with the projection  
   \be \la{CPmu}
\mu_\bS: {\mathscr P}_{\G, \bS} \lra {\rm C}_{\G, \bS}, 
\ee
defined  in  Theorem \ref{Th1.14} below, see (\ref{mubs-tcas}), by $\mu_\bS=i_{\rm C} \circ \mu$.  
  We  use crucially  the description of 
 the center of the algebra ${\cal O}_q({\mathscr P}_{\G, \bS})$ for generic values of $q$ from Theorem \ref{TH1.16a} in 
 terms of the Casimir torus $ {\rm C}_{\G, \bS}$. \vskip 2mm

Theorems \ref{MTHCQ},  \ref{MTH}  and \ref{TH1.16a}  immediately imply  our next main result.
 
 \bt \la{Th1.12} Let $\G$ be a semi-simple split adjoint group, and  $\bS$   a colored decorated surface as in Theorem \ref{MTH}. Then for any  $\hbar$ such that $\hbar \in \R_{>0}$,  or $\G$ is simply-laced and $|\hbar|=1$,  there is    a family of $\Gamma_{\G, \bS}$-equivariant $\ast-$representations of the $\ast$-algebra $\mathcal{A}_\hbar({\mathscr P}_{\G, \bS})$ by unbounded operators in  Hilbert spaces ${\cal H}({\mathscr P}_{\G, \bS})_\lambda$, parametrized 
 by  the positive points of the Casimir torus  $\lambda\in {\rm C}_{\G, \bS}(\R_{>0})$. We   refer to it  as the {\it principal series of representations}.
It includes   unitary projective representations of   the group $ \Gamma_{\G, \bS}$   in    ${\cal H}({\mathscr P}_{\G, \bS})_\lambda$.   \et

\subsubsection{Quantization of the moduli space $ {\rm Loc}_{\G, \bS}$.} Projection (\ref{PRW}) gives rise to an embedding 
$$
\pi^*: {\cal O}({\rm Loc}_{\G, \bS}) \hlra {\cal O}({\mathscr P}_{\G, \bS})^{W^n}.
$$

Recall    the  two a priori different algebras  assigned to the moduli space ${\mathscr P}_{\G, \bS}$: the  
algebra of regular functions ${\cal O}({\mathscr P}_{\G, \bS})$, and the algebra 
${\cal O}^{\rm cl}({\mathscr P}_{\G, \bS})$  of { universally Laurent polynomials}.

\bt \la{MONal} \cite{S20} Unless  $\bS$ has  no special points and $1$ puncture,  the map $\pi^*$ induces an isomorphism
\be \la{MONbl}
 {\cal O}({\rm Loc}_{\G, \bS}) ={\cal O}^{\rm cl}({\mathscr P}_{\G, \bS})^{W^n}.
\ee
\et

By Theorem \ref{MTH},  the    
 group $W^n$ acts on the space ${\mathscr P}_{\G, S}$ by  cluster Poisson transformations.     This allows   to define  the algebra ${\cal O}_q({\rm Loc}_{\G, \bS})$ as the $W^n-$invariants 
 of the algebra ${\cal O}_q({\mathscr P}_{\G, \bS})$, 
 and  introduce its Langlands modular double  $\mathcal{A}_\hbar({\rm Loc}_{\G, \bS})$ as in (\ref{AH}). Theorem \ref{MONal} justifies the notation.  
 \be
 \begin{split}
 &{\cal O}_q({\rm Loc}_{\G, \bS}):= {\cal O}_q({\mathscr P}_{\G, \bS})^{W^n}.\\
 &\mathcal{A}_\hbar({\rm Loc}_{\G, \bS}):= {\cal O}_q({\rm Loc}_{\G, \bS})\otimes_\C{\cal O}_{q^\vee}({\rm Loc}_{\G^\vee, \bS}).\\
 \end{split}
 \ee
 
Theorem \ref{MTHCQ}  provides the  principal series 
  representations $\rho_\lambda$,  where  $\lambda \in \H(\R_{>0})^n$,  of the $\ast$-algebra $\mathcal{A}_\hbar({\mathscr P}_{\G, \bS})$. 
 Denote by ${\rm Res}_{{\rm Loc}}\rho_\lambda$ it restriction to  the   $\ast$-algebra $\mathcal{A}_\hbar({\rm Loc}_{\G, \bS})$.  \vskip 1mm

 Projection (\ref{CPmu}) intertwines the actions of the group $W^n$ on ${\mathscr P}_{\G, \bS}$ and on 
  ${\rm C}_{\G, \bS}(\R_{>0})$.   Since the action of   $W^n$ on   ${\mathscr P}_{\G, \bS}$ is cluster,  by Theorem  \ref{MTHCQ}
   it gives rise to unitary  intertwiner operators ${\cal I}_w$,   $w \in W^n$. They establish  unitary equivalence of 
   representations  of the  $\ast$-algebra $\mathcal{A}_\hbar({\rm Loc}_{\G, \bS})$ on the same $W^n$-orbit: 
\be \nonumber
 {\cal I}_w:  {\rm Res}_{{\rm Loc} }\rho_{\lambda} \stackrel{\sim}{\lra} {\rm Res}_{{\rm Loc} }\rho_{w(\lambda)}.
\ee
 Indeed, by the very definition, $\mathcal{A}_\hbar({\rm Loc}_{\G, \bS})$ lies in the algebra of $W^n$-invariants of $\mathcal{A}_\hbar({\mathscr P}_{\G, S})$. 

\vskip 1mm
{\it This way we get the quantized moduli space ${\rm Loc}_{\G, \bS}$ with all the properties listed in Section \ref{5.1.1} }.

 \medskip 
 
  \subsection{The modular functor conjecture.}  \medskip
  
  In this Section we elaborate on the Conjecture stated in \cite[Conjecture 6.2]{FG07}.

  Below $\bS$ is an arbitrary colored decorated surface.

\subsubsection{The algebra side of the Modular Functor conjecture} \la{sec5.3.3} 
Given an oriented loop $\alpha$ on $\bS$, the semisimple part of the monodromy around $\alpha$ defines a regular map 
\be
{\rm m}_\alpha: {\mathscr P}_{\G, \bS} \lra {\rm H}/W.
\ee

\bcon \la{QMONC}  There is a canonical q-deformation of the monodromy map ${\rm m}_\alpha$, providing a map of algebras, called the quantum monodromy map:
\be \la{QMas}
{\rm M}^*_\alpha: {\cal O}({\rm H}/W) \lra {\cal O}_q({\mathscr P}_{\G, \bS}).
\ee

\econ

We have an explicit construction  when $\G$ is of type $\A_n$, which  we will discuss on another occasion. \\

  Let $\alpha$ be a simple non-trivial  loop on a colored decorated surface $\bS$. Let  
$$
\bS_\alpha:= \bS - \alpha
$$
 be the surface obtained by cutting $\bS$ along the loop $\alpha$. 
Denote by $\alpha_{\pm}$ the new boundary loops of  $\bS_\alpha$, oriented by the orientation of $\bS_\alpha$. 
There is a canonical map given by the monodromies around the oriented loops $\alpha_{\pm}$:
\be \la{qMmP}
{\rm m}_{\alpha_{\pm}}: {\mathscr P}_{\G, \bS_\alpha} \lra {\rm H}  \times {\rm H}. 
\ee


\bl
The map ${\rm m}_{\alpha_{\pm}}$ in (\ref{qMmP})   gives rise to a commutative subalgebra  
\be \la{comsa}
{\rm M}_{\alpha_{\pm}}^*{\cal O}({\rm H} \times {\rm H})  \subset {\cal O}_q({\mathscr P}_{\G, \bS_\alpha}).
 \ee
 \el
 
 \begin{proof} The moduli space ${\mathscr P}_{\G, \bS_\alpha}$ has a cluster Poisson structure. In any cluster Poisson coordinate system, the map $m_{\alpha_{\pm}}$ is given by a monomial map.  Its monomial components lie in   the Poisson center, and described by the general cluster Poisson construction \cite{FG03b}. In particular, they are immediately quantized, 
 and the result does not depend on the cluster coordinate system. 
 \end{proof}
 
Consider  the pull back of the ideal ${\cal I}_\alpha$ defining the antidiagonal $(h, h^{-1})$ in $({\rm H}  \times {\rm H})$: 
 $$
 {\cal I}_\alpha \subset {\rm M}^*_{\alpha_\pm}{\cal O}({\rm H}  \times {\rm H}).
 $$ 
There is an  algebra quantizing the subspace of the moduli space ${\mathscr P}_{\G, \bS_\alpha}$ defined by the condition that 
 the monodromies around the   oriented loops $\alpha_+$ and $\alpha_-$ are inverse to each other:
  \be \la{132}
  {\cal O}_q({\mathscr P}_{\G, \bS_\alpha})/{\cal I}_\alpha. 
  \ee 
  The group $W \times W$ acts on ${\rm H}  \times {\rm H}$, and it   acts   by cluster transformations on ${\cal O}_q({\mathscr P}_{\G, \bS_\alpha})$. 
 Embedding (\ref{comsa}) is $W \times W-$equivariant. 
  The  diagonal subgroup $W_\Delta \subset W \times W$ preserves  the ideal  ${\cal I}_\alpha$, and hence acts on the quotient (\ref{132}). 
 Consider the subalgebra of $W_\Delta-$invariants:
  \be
   \Bigl( {\cal O}_q({\mathscr P}_{\G, \bS_\alpha})/{\cal I}_\alpha\Bigr)^{W_\Delta}.  
   \ee
 Modular Functor Conjecture describes how   algebras $ {\cal O}_q({\mathscr P}_{\G, \bS})$ behave  under the cutting of $\bS$ by a loop $\alpha$. 
  
    \bcon \la{koskin} There is a canonical isomorphism of algebras
  \be \la{MFKiso}
  \Bigl( {\cal O}_q({\mathscr P}_{\G, \bS_\alpha})/{\cal I}_\alpha\Bigr)^{W_\Delta}  =   {\cal O}_q({\mathscr P}_{\G, \bS})^{{\rm M}^*_\alpha {\cal O}({\rm H}/W)}.
      \ee
  
  \econ
  
 \subsubsection{The analytic side of the Modular Functor conjecture.} 

The principal series representations 
assigned to $\bS$ are parametrized by  $\lambda \in {\rm C}_{\G, \bS}(\R_{>0})$. The ones 
assigned to $\bS_\alpha$ are parametrized by the triples  $$
(\lambda, \chi_-, \chi_+) \in {\rm C}_{\G, \bS}(\R_{>0}) \times \H(\R_{>0})\times \H(\R_{>0}).
$$ 
Here $\chi_\pm$   are the parameters at the   punctures $p_{\pm}$.  By Theorem \ref{Th1.12}, representations attached to  elements $(\lambda, w_1(\chi_-), w_2(\chi_+))$ were $w_1, w_2 \in W$ 
 are canonically equivalent. The equivalences are given by the unitary intertwiners corresponding to the elements   $(w_1, w_2)$.
 
The restriction to $ \bS_\alpha$ provides a  map of the moduli spaces and a dual map of the $\ast$-algebras:
\be \la{MK2}
\begin{split}
&{\rm Res}_\alpha: {\mathscr P}_{\G,   \bS} \lra {\mathscr P}_{\G,   \bS_\alpha}.\\
&{\rm Res}^*_\alpha: \mathcal{A}_\hbar({\mathscr P}_{\G,   \bS_\alpha}) \lra \mathcal{A}_\hbar({\mathscr P}_{\G,   \bS}).\\
\end{split}
\ee


\begin{conjecture} \label{MK} Given a decorated surface $\bS$, the principal series 
representation ${\cal H}({\mathscr P}_{\G,   \bS})_\lambda$ 
has a  natural $\Gamma_{\G,   \bS_\alpha}$-equivariant 
decomposition into an integral of Hilbert spaces
\begin{equation} \label{MK1}
{\cal H}({\mathscr P}_{\G,   \bS})_\lambda \stackrel{\sim}{=} 
\int_{\chi}{\cal H}({\mathscr P}_{\G,   \bS_\alpha})_{\lambda, \chi, \chi^{-1}} d\mu_\chi, \qquad \qquad \chi\in \H/W(\R_{>0}).
\end{equation} 
Here ${\cal H}({\mathscr P}_{\G,   \bS_\alpha})_{\lambda, \chi, \chi^{-1}} $ is the principal series representation 
assigned to   $\bS_\alpha$.

The isomorphism (\ref{MK1}) is an isomorphism of $\mathcal{A}_\hbar({\mathscr P}_{\G,   \bS_\alpha})$-modules in the following sense:\footnote{The algebra $\mathcal{A}_\hbar({\mathscr P}_{\G,   \bS'})$ acts by unbounded operators,  so   map (\ref{MK1})  is not   literally a map of $\mathcal{A}_\hbar({\mathscr P}_{\G,   \bS_\alpha})$-modules.} 

\noindent
It intertwines the  induced via  map (\ref{MK2})  action of the $\ast$-algebra $\mathcal{A}_\hbar({\mathscr P}_{\G,   \bS_\alpha})$ on the Schwartz space on the left  with its 
 natural action   on the Schwartz spaces on the right.
\end{conjecture}

 The   modular functor conjecture  for ${\rm PGL_2}$ was claimed   by Teschner \cite{T}. A proof  for ${\rm PGL}_m$ was   carried out by Schrader  $\&$ Shapiro in  \cite{SS2}-\cite{SS4}.  
Although we consider   spaces ${\mathscr P}_{\G, \bS}$ rather then ${\mathscr X}_{\G, \bS}$, the arguments \cite{SS2}-\cite{SS4}  should not require essential changes.

 \medskip
\subsection{Historical comments and the architecture of the paper} 
\medskip

 \subsubsection{Historical comments.} Some parts of Theorems \ref{MTH} $\&$ \ref{MTHa} were known. Here is a brief survey.
 
 \vskip 1mm
 
 A $\Gamma_\bS$-equivariant positive atlas  on   the moduli spaces  ${\mathscr A}_{{\G}, \bS}$  and  ${\mathscr X}_{{\G}, \bS}$ was   
was defined  in \cite{FG03a} using total positivity \cite{L1}. This is sufficient to develop many aspects of the Higher Teichmuller theory, e.g. the Higher Teichmuller and lamination spaces,  but not to quantize it. 
 
A $\Gamma_\bS$-equivariant cluster atlases  on   the moduli spaces  ${\mathscr A}_{{\rm SL_m}, \bS}$  and  ${\mathscr P}_{{\rm PGL_m}, \bS}$ were    
defined  in \cite{FG03a}. The construction relies   on the existence of the {\it special cluster coordinate system} on the triangle in both cases, which 
 is a unique feature of the ${\rm A}_m$-case. The cyclic invariance holds for it on the nose.

The cluster nature  of the geometric Weyl group action for the type ${\rm A}$ was established \cite[\S 7-8]{GS16}. Using crucially
   the technique  developed there,  proving the cluster nature 
  of the Weyl group for general $\G$ become a   straightforward exercise, executed in Section \ref{Sec5.2}. 
  
\vskip 1mm
Breaking the symmetry,   pick  an angle of   $t$  and identify   the  space  ${\mathscr P}_{\G, t}$ with an open part of $\B \times \H$, see Lemma \ref{1.2.09.1}. 
Then the action of the group $\H^3$ on ${\mathscr P}_{\G, t}$ from Section \ref{1.2.9}  is given by the left and right actions of $\H$ on $\B$, and the action of $\H$ on itself. The cluster Poisson structure on the Borel group $\B$ was defined in \cite{FG05}. Its quiver $ {\bf q}_\B$ is   a large subquiver
 of the coveted quiver $ {\bf q}_{\G, t}$ for  ${\mathscr P}_{\G, t}$. 
 
It remains    to add to the quiver  $ {\bf q}_\B$ the  extra $r$ frozen vertices,    half-integral arrows between them, and the arrows connecting them with 
  $ {\bf q}_\B$. This    is done  in Section \ref{pcs.sec}. In contrast with the type ${\rm A}$, it is very unclear why the mutation class of the quiver ${\bf q}_{\G, t}$ is  preserved by  rotations of the triangle $t$. We prove this in Section \ref{sec.conf3}.  The cyclic invariance of the cluster Poisson structure of  the non-frozen base 
  ${\rm Conf}^\times_3({\cal B})$ is relatively easy. The difficult part is to handle the  frozen fibers. 
 The frozen variables are essential:  we use them to amalgamate   spaces ${\mathscr P}_{\G, t}$. \vskip 1mm
  
 Generalising \cite{FG03a} and using  the amalgamation technique   \cite{FG05},  Ian Le in tour de force works \cite{IL1}, \cite{IL2}  
 established cluster structures on the   spaces ${\mathscr A}_{\G, \bS}$  and 
 ${\mathscr X}_{\G, \bS}$  for the type $D_m$, and then by using folding, for all classical groups and $\G_2$.    However  these cases were treated by case by case study. 
 I. Le elaborated   very interesting specific cluster coordinate systems  for each case.
 
 Our   approach is systematic,   works for all   groups at once,   
 and  much   more flexible even for the type ${\rm A}$. 
 
 \vskip 1mm
 Cluster nature of the  Weyl group action on  spaces ${\mathscr A}_{\G, \bS}$ for  types ${\rm B,C,D}$  was  shown in \cite{IIO}.
  
  \vskip 1mm
  Another  quiver   related to the space ${\rm Conf}_3({\cal A})$ was studied by Jiarui Fei  \cite{F} from a very different  perspective.  See \cite[Appendix A]{F} for a list of such quivers.
  The arrows containing  unfrozen vertices in Fei's quivers coincide with ours. However  the arrows between the frozen vertices are different, e.g.
   there are no half-integral arrows in Fei's  quivers.   
Therefore one can not use them to define quivers for  surfaces other than a triangle.  
The construction of Fei's quivers were from the \emph{Auslander-Reiten} theory of the category of projective presentations of Dynkin quivers. In this context,  quivers with half integral arrows are not defined. 
 It would be very interesting to understand the relation between the two stories.

\subsubsection{The architecture of the paper.}  
The paper can be subdivided into four parts.

\begin{enumerate}  

\item {\bf  Quantum Geometry of Moduli Spaces.} It consists of Sections \ref{SEC1} - \ref{sec1.9}.   \vskip 1mm

   In Sections \ref{SEC1} - \ref{KFMS} we present   main definitions    and some of the main constructions,  and state the results. 
  We provide  proofs  when they are simple. 
  The   proofs requiring elaborate arguments appear later in the text.  We do not duplicate    later on any of the discussions in Sections \ref{SEC1} -\ref{KFMS}. 
  
  \vskip 1mm
In Section \ref{PLG} we realize the Poisson Lie group $\G$, the dual Poisson Lie group $\G^*$, the Grothendieck resolution $\widehat G \lra \G$, 
  the Drinfeld double $\D(\G)$ and its   dual $\D(\G)^*$, the Heisenberg double $\D_{\rm H}(\G)$,  and the  dressing transformations 
via    moduli spaces 
of $\G-$local systems on decorated surfaces. 

\vskip 1mm
In Section \ref{TQFT1} we explain that the quantized  moduli spaces of $\G$-local systems on decorated surfaces  provide  a nonlinear geometric avatar of  TQFT. 

One of its linearizations relates to infinite-dimensional $\ast-$representations of quantum groups. 

In Sections \ref{SEC1.8} and  \ref{SC2.7} we present applications to the representation theory of quantum groups. 

In Section  \ref{SECT2.8a} we define   integrable systems related to the moduli  spaces of $\G-$local systems. 

Section \ref{SEC1.8} is elaborated   further in Section \ref{SSEECC11.3}. 
 Sections \ref{SC2.7} - \ref{SECT2.8a} are elaborated in Section \ref{Sec8}. 
 
 \vskip 1mm

 Section \ref{sec1.9}  presents our Main Conjecture, which relates the cluster quantization of   moduli spaces ${\rm Loc}_{\G, S}$ with   
 each of the following three objects: 
 
 (i) coinvariants of oscillatory representations of $W-$algebras,  
 
 (ii) principal series of $\ast-$representations of quantum groups,

 (iii) conformal blocks for Toda theories. \vskip 2mm

\item {\bf The Cluster Structure.}  It  consists of Sections \ref{SSECC3} - \ref{SEC5}.  \vskip 1mm

In Section \ref{SSECC3} we give an elementary description of  the cluster Poisson structure for the moduli space ${\mathscr P}_{{\rm PGL_m, \bS}}$ by amalgamation of the 
moduli spaces assigned to the  triangles of a triangulation of $\bS$. It serves an introduction to   Sections \ref{Clcoord} - \ref{CSEC7}. The general case is much more complicated: it takes few  pages to do the $\A_m$ case, and almost 50 to do the general $\G$ case.   \vskip 1mm

 In Sections \ref{Clcoord}-\ref{SEC5} we establish the cluster Poisson/$K_2-$structure of the moduli space ${\mathscr P}_{{\G, \bS}}$/${\mathscr A}_{{\G, \bS}}$, and prove their equivariance under the action of the mapping class group of $\bS$.   \vskip 1mm 
 

In Section \ref{Clcoord} we introduce a collection of  coordinate systems on each of the   
spaces ${\mathscr P}_{\G,  \bS}$ and ${\mathscr A}_{\G,  \bS}$. We start from  a triangle $t$.  
Given an angle $a_t$ of the triangle, and a reduced decomposition  of the 
longest element $w_0$ of the Weyl group, we introduce   coordinates on  the  spaces ${\mathscr P}_{\G,  t}$ and ${\mathscr A}_{\G,  t}$.   
Given a triangulation ${\cal T}$ of $\bS$, 
we assemble    spaces ${\mathscr P}_{\G,  \bS}$ and ${\mathscr A}_{\G,  \bS}$ 
by  gluing   the spaces ${\mathscr P}_{\G,  t}$ and ${\mathscr A}_{\G, t}$ assigned to   triangles $t$ of ${\cal T}$. 
So  given  a combinatorial datum:
 \be \la{data}
\begin{split}
&\mbox{(A triangulation ${\cal T}$ of $\bS$, an angle $a_t$ for each triangle $t$ of the triangulation, }  \\
&\mbox{and,  for each of the   angles, a reduced decomposition    of $w_0$).}\\
\end{split}
\ee
we  get   coordinate systems  on the  spaces ${\mathscr P}_{\G,  \bS}$ and ${\mathscr A}_{\G,  \bS}$ 
by the amalgamation  of the ones on the triangles. 
We review the amalgamation construction   in Section \ref{sec2}.   \vskip 1mm

In Section \ref{SECT4} we give a representation theoretic description of   cluster coordinates on the 
space ${\mathscr A}_{\G,  t}$. This leads to a coordinate free definition of   cluster coordinates on  ${\mathscr A}_{\G,  \bS}$.
Yet the explicit definition of Section \ref{Clcoord} is indispensable in the proofs.   \vskip 1mm

After Section  \ref{Clcoord} and \ref{SECT4} two crucial things remains to be done. 

First,  in Section \ref{pcs.sec} we introduce the quiver  assigned to the combinatorial datum (\ref{data}). We show  that its  vertices match the coordinates on each of the spaces 
${\mathscr P}_{\G,  \bS}$ and ${\mathscr A}_{\G,  \bS}$ introduced in   Section \ref{Clcoord}. 
This way we get a cluster coordinate system assigned to the datum (\ref{data}).  \vskip 1mm

Second,  in Section \ref{sec.conf3} - \ref{CSEC7} we prove that the cluster coordinate systems assigned to different data (\ref{data}) are related by cluster transformations. 
This implies that the obtained cluster structures are equivariant under the action of the mapping class group of $\bS$.  \vskip 1mm

In Section \ref{SEC5} we define the braid group actions on the  spaces ${\mathscr P}_{\G,  \bS}$ and ${\mathscr A}_{\G,  \bS}$. We show that the Weyl, the braid and 
the outer automorphism   groups act by cluster transformations of   ${\mathscr P}_{\G,  \bS}$ and ${\mathscr A}_{\G,  \bS}$.   \vskip 2 mm

\item {\bf  Applications}. It consists of Sections \ref{DT=cluster} - \ref{Sec8}.  \vskip 1mm

In Section \ref{DT=cluster} we prove Theorem \ref{CDT} 
on the structure of the Donaldson-Thomas transformations for the moduli spaces ${\mathscr P}_{\G, \bS}$.  \vskip 1mm

 In Section \ref{SSEECC11.3} we present a geometric cluster realization of  quantum groups. We prove that the braid group action coincides with Lusztig's braid group action. 
Combined with the results 
 \cite{FG07} and Sections \ref{Sec8a} - \ref{Sec.quant.s}, it delivers a construction of the principal series of $\ast-$representations of the modular double of the quantum group, together 
with the unitary intertwiners corrrsponding to the Weyl group action, which are   quantum analogs of the Gelfand-Graev intertwiners \cite{GG73}.  \vskip 1mm

 Section \ref{Sec8}   develops realizations of  the  principal series   $\ast-$representations of   the spaces ${\mathscr P}_{\G, \bS}$, defined by presenting the decorated surface $\bS$ as a topological double of another decorated surface. We apply this to construct realizations of the principal series representations of the quantum groups, their 
 multiplicity spaces, and the analog of the regular regular representation for the quantum group. \vskip 2mm

\item {\bf   Tools and Supplements}.  It consists of Sections \ref{Sec8a} - \ref{Appen12}. \vskip 1mm
In Section \ref{Sec8.1} we recall the quantum double construction, which was used in Section \ref{Sec8}.

In Sections \ref{Sec8.2a}  - \ref{Sec8.2} we discuss the principal series $\ast-$representations of the cluster quantum double, which we   used in Section \ref{Sec8}.

In Section \ref{SEC8.4} we discuss the principal series $\ast-$representations of   cluster Poisson varieties. 

The crucial new feature of Sections  \ref{Sec8.2a}  - \ref{SEC8.4}  is the $|\hbar|=1$ case.  \vskip 1mm

In Section \ref{Sec.quant.s} we prove the Quantum Lift Theorem. It is the main tool   to prove that   certain rational function lift canonically to   the algebra of quantum cluster regular functions.  
We use it crucially in Section \ref{SSEECC11.3}.  It is independent of the rest of the paper. \vskip 1mm

In Section \ref{Appen12} we recall our main working horse - the  amalgamation of cluster varieties, and discuss   quasicluster transformations.

 \end{enumerate}


\medskip

 \section{Basic constructions} \la{SEC1} 
 
\medskip
  
In Section \ref{SECC1.2} we recall basics  on split reductive groups and  discuss decorated flags and pinnings. 
  
In Section \ref{Sec2.2}    we  revisit the  moduli spaces ${\mathscr P}_{\G, \bS}$  and  the gluing maps.

In Section \ref{POTEN} we recall from \cite{GS13} the construction of the potential functions. 
 
In Section \ref{EXAMPLES} we discuss the principal affine space and the potentials functions for the groups of type $\A_n$.

 

   \subsection{Basics of split semi-simple algebraic groups} \la{SECC1.2}

\medskip
\subsubsection{The Cartan group of $\G$.} Take a Borel subgroup 
$\B$ of $\G$. 
The Cartan group of $\B$ is the quotient $\B/\U$, where $\U$ is the maximal unipotent subgroup inside $\B$. 
The Cartan group of another Borel subgroup $\B'$ is canonically isomorphic to the Cartan group of $\B$. Indeed, the subgroups $\B$ and $\B'$ are conjugated, and $\B$ acts trivially on $\B/\U$. Thus the group $\G$ determines {\it the Cartan group} $\H$. It is equipped with the root data. 
A maximal split torus  in $\G$ is called a {\it Cartan subgroup}. It is isomorphic to the Cartan group without  canonical root data.

\subsubsection{Flags and decorated flags.} \la{2.1.2} The {\it flag variety} ${\cal B}$ parametrizes   Borel subgroups   in $\G$. The group $\G$ acts transitively on  ${\cal B}$. A choice of a Borel subgroup $\B$ defines an isomorphism ${\cal B} = \G/\B$. A pair of flags $(\B_1, \B_2)$  is  {\it generic} if it belongs to the unique open $\G$-orbit in ${\cal B}\times {\cal B}$. Equivalently, $\B_1\cap \B_2$ is a Cartan subgroup of $\G$. 
Assigning to a Borel subgroup $\B$ its unipotent radical $\U$, we get a bijection between   Borel subgroups and maximal unipotent subgroups in $\G$. The Borel subgroup $\B$ assigned to $\U$ is the normaliser of $\U$. 
Let  ${\mathcal A}:=\G/\U$ be the {\it decorated flag variety}, also known as the {\it principal affine space}. 
Let $\A \in {\cal A}$ denote a decorated flag.
There is a natural projection 
\be\pi: {\cal A}\lra {\cal B}, \hskip 7mm \A \lms \B.
\ee
The Cartan group ${\rm H}$ acts on ${\cal A}$ on the right. 
The  fibers of the map $\pi$ are    ${\rm H}$-torsors.
A pair $(\A_1, \A_2)$ of decorated flags is said to be {\it generic} if its underlying pair of flags $(\B_1, \B_2)$ is generic. 

\subsubsection{The root data.} Let  ${\rm I}=\{1, \ldots, r\}$ be the set of vertices of the Dynkin diagram, and   $C_{ij}$  the Cartan matrix.   
Let $\alpha_1, \ldots, \alpha_r$ be the simple positive roots, $\alpha_1^\vee, \ldots, \alpha^\vee_r$ the  simple positive coroots, and  $\Lambda_1, \ldots, \Lambda_r$  the fundamental weights, understood as     characters / cocharacters of the Cartan group ${\rm H}$:
\[
\alpha_i:\H \lra {\Bbb G}_m, \hskip 5mm \alpha_i^\vee:{\Bbb G}_m \lra \H, \hskip 5mm \Lambda_j: \H \lra {\Bbb G}_m; \hskip 10mm \Lambda_j\circ \alpha_i^\vee =  {\rm Id}^{\delta_{ij}}_{{\Bbb G}_m}, 
 \hskip 5mm  \alpha_i \circ \alpha_j^\vee = {\rm Id}^{C_{ij}}_{{\Bbb G}_m}.
\]

\subsubsection{Pinnings in $\G$.} \la{3.1.4a}
Let $(\B, \B^-)$ be a generic pair of flags. 
Let ${\rm U}:=[ \B, \B ]$ and ${\rm U}^-:=[\B^-, \B^-]$ be maximal unipotent subgroups. 
Let $x_i: {\mathbb A}^1 \rightarrow {\rm U}$  be a unipotent subgroup associated to the simple root $\alpha_i$. 
Equivalently, the choices of $x_1, \ldots, x_r$ determine an additive isomorphism
\be 
\la{additive.isomorphism.5.21}
\begin{split}
&(\chi_1, \chi_2,\ldots, \chi_r):~\U/[\U, \U] \stackrel{\sim}{\lra} {\mathbb A}^r,\\
&\chi_i(x_j(a)) =  \left\{  \begin{array}{ll}  
      a & \mbox{if } i=j \\
      0 & \mbox{if } i\neq j. \\
   \end{array}\right.\\
   \end{split}
\ee
Let $y_i: {\mathbb A}^1 \rightarrow {\rm U}^- $ be a unipotent subgroup associated to the simple root $-\alpha_i$. 

\bd \la{DEFPIN} The datum  $p = (\B, \B^-, x_i,  y_i; i\in {\rm I})$ is called a {\it pinning} over $(\B, \B^-)$ if it gives rise to a homomorphism $\gamma_i : {\rm SL}_2 \to \G$ for each $i\in {\rm I}$ such that 
\be
\la{5.21.2016.ss}
\gamma_i   \begin{pmatrix}  
      1 & a \\
      0 & 1 \\
   \end{pmatrix}= x_i(a),  
 \hskip 7mm \gamma_i \begin{pmatrix}  
      1 & 0 \\
      a & 1 \\
   \end{pmatrix}= y_i(a), \hskip 7mm 
   \gamma_i \begin{pmatrix}  
      a & 0 \\
      0 & a^{-1} \\
   \end{pmatrix}=\alpha_i^\vee(a).
\ee
\ed
We sometimes fix a pinning $p$ and refer to it as {\it a standard pinning}. This choice will be useful for concrete calculation. 

\subsubsection{The Weyl group of $\G$.} \la{3.1.5} A pinning provides two ways to lift elements of the Weyl group $ W$ to $\G$.  
Namely, let $s_i \,(i\in {\rm I})$ be the simple reflections generating  $W$. Set 
\be
\overline{s}_i:=y_i(1)x_i(-1)y_i(1), \hskip 7mm \overline{\overline{s}}_i:=x_i(1)y_i(-1)x_i(1).
\ee
The elements $\{ \overline{s}_i \}$ and $\{\overline{\overline{s}}_i\}$ satisfy the braid relations.
Thus, we can associate with each $w\in W$ two representatives $\overline{w}, \overline{\overline{w}} \in \G$ such that for any reduced decomposition $w=s_{i_1}\ldots s_{i_k}$ one has $\overline{w}=\overline{s}_{i_1}\ldots \overline{s}_{i_k}$ and $ \overline{\overline{w}}=\overline{\overline{s}}_{i_1}\ldots\overline{ \overline{s}}_{i_k}.$ 

Let $w_0$ be the longest element of the Weyl group. The canonical central element $s_\G$ is the image of $   \begin{pmatrix}  
      -1 & 0 \\
      0  & -1 \\
 \end{pmatrix}$ under a principal embedding ${\rm SL}_2 \hookrightarrow \G$.
We have $s_\G= \overline{w}_0^2 = \overline{\overline{w}}_0^2$. 

Recall the involution $i \mapsto i^*$ of ${\rm I}$ such that
$
\alpha_{i^*}^\vee = -w_0(\alpha_i^\vee).
$ For $w\in W$, we set $w^\ast= w_0 w w_0^{-1}$. A reduced decomposition $w=s_{i_1}\ldots s_{i_k}$ induces a reduced decomposition $w^\ast=s_{i_1^\ast}\ldots s_{i_k^\ast}$.

\subsubsection{Decorated flags} \la{SECC1.3}

Fix a standard pinning. Let $[{\rm U}]\in \mathcal{A}$ be the left coset of ${\rm U}$ in ${\rm G}$ containing $1$ and let $[{\rm U}^{-}]= \overline{w}_0\cdot [{\rm U}]$.

 \blc \la{LD2.2} Consider the configuration space
$$
{\rm Conf}_2({\cal A}):= \G \backslash {\cal A}^2. 
$$ 
For each pair $\A_1, \A_2 \in {\cal A}$, there  exist a unique $h \in {\rm H}$ and a unique  $w\in W$ such that
\[
(\A_1, \A_2) = (h\cdot [\U],  \overline{w} \cdot [\U]) \hskip 2mm  \in {\rm Conf}_2({\cal A}).
\]
Define the {\it basic invariants} 
\be \la{99}
 h(\A_1, \A_2):= h, \hskip 7mm w(\A_1, \A_2):= w. 
\ee
A pair of decorated flags  $(\A_1$, $\A_2)$ is generic if and only if $w(\A_1, \A_2)=w_0$. 
\elc
\begin{proof}
Recall the Bruhat decomposition
$
\G= \bigsqcup_{w \in W} \U {\rm H} \overline{w} \U.
$ 
Therefore
\[
{\rm Conf}_2({\cal A})=\G \backslash (\G/ {\rm U})^2\stackrel{\sim}{=}  \U \backslash \G / \U  = \bigsqcup_{w \in W} {\rm H} \overline{w}
\]
consists of $|W|$-many copies of ${\rm H}$.
The rest follows easily.
\end{proof}

Similarly, we define $h_-(\A_1, \A_2):= h_-$ and $w_-(\A_1, \A_2):=w_-$ such that
\[
(\A_1, \A_2)=(h_- \cdot [{\rm U}^-], ~\overline{w}_-\cdot [{\rm U}^-]).
\]
It is easy to show that
\be
h(\A_1, \A_2) = w_0(h_-(\A_1, \A_2)), \hskip 7mm w(\A_1, \A_2) = w_-(\A_1, \A_2)^*.
\ee

The invariant  $w(\A_1, \A_2)$ is called the $w$-distance and depends only on the underlying flags.   We denote  by $w(\B_1, \B_2)$ the $w$-distance of  a  pair $(\B_1, \B_2)$, and by  $w_-(\B_1, \B_2)$ the $w_-$-distance.

\bl 
\la{lem.part.flag.decom}
Suppose $u, v\in W$ are such that the length 
$l(uv) = l(u) + l(v).$ Then
\begin{itemize}
\vskip 1mm \item[1.] If $w(\B_1, \B_2)= uv$, then there is a unique flag $\B'$ such that  
\[
w(\B_1, \B')= u, \hskip 7mm w(\B', \B_2) =v.
\]
\vskip 1mm \item[2.] Conversely, if $w(\B_1, \B')= u$ and $w(\B', \B_2) = v$, then $w(\B_1, \B_2)=uv$.
\end{itemize}
\el 

 Lemma \ref{lem.part.flag.decom} is well-known.  The following Corollary is a direct consequence of Lemma \ref{lem.part.flag.decom}.
\bc
\la{USFX}
Let  $(\B_l, \B_r)$ be a pair of   flags with
$w(\B_l, \B_r)=w.$   Every reduced word ${\bf i}=(i_1, \ldots, i_p)$ of $w$ gives rise to a unique ${\bf i}$-chain of  flags  
$
\{\B_l=\B_0, \, \B_1,\, \B_2, \ldots,\, \B_p=\B_r\} 
$ such that  
$w(\B_{k-1}, \B_{k})= s_{i_k}$ for  $k\in \{1,\ldots, p\}$, as   illustrated below:
\begin{center}
\begin{tikzpicture}[scale=0.8]
\node [circle,draw=red,fill=red,minimum size=3pt,inner sep=0pt,label=below:{\small $\B_l={\B}_0\qquad ~$}] (a) at (0,0) {};
 \node [circle,draw=red,fill=white,minimum size=3pt,inner sep=0pt,label=below:{\small $~{\B}_1$}]  (b) at (1,0){};
 \node [circle,draw=red,fill=white,minimum size=3pt,inner sep=0pt,label=below:{\small $~{\B}_2$}] (c) at (2,0) {};
 \node [circle,draw=red,fill=white,minimum size=3pt,inner sep=0pt,label=below:{\small $~{\B}_3$}] (d) at (3,0) {};
  \node [circle,draw=red,fill=white,minimum size=3pt,inner sep=0pt,label=below:{\small $~{\B}_{4}$}] (e) at (4,0) {};
  \node [circle,draw=red,fill=red,minimum size=3pt,inner sep=0pt,label=below:{\small $\qquad \qquad ~{\B}_{p}=\B_r$}] (f) at (5.5,0) {};
  \node[blue] at (0.5, 0.2) {\small $s_{i_1}$};
   \node[blue] at (1.5, 0.2) {\small $s_{i_2}$};
    \node[blue] at (2.5, 0.2) {\small $s_{i_3}$};
     \node[blue] at (3.5, 0.2) {\small $s_{i_4}$};
 \foreach \from/\to in {a/b, b/c, c/d, e/f, d/e}
                   \draw[thick] (\from) -- (\to);
   \node at (4.725, -.4) {$\cdots$};         
    \node at (4.725, .2) {$\cdots$};       
 \end{tikzpicture}
\end{center}
\ec

\bp \la{P2.5}
\la{idenity.basic.invariants.tt}
1. If $w(\A_1, \A_2)=u$, then 
\bg  \la{lem.16.13.05.hhh.3.26}
h(\A_1\cdot h_1, \A_2\cdot h_2) =  h_1 u(h_2^{-1}) \, h(\A_1, \A_2), \\
 h(\A_1, \A_2)\, \overline{\overline{u}}= {\overline{u}}\, h(\A_2, \A_1)^{-1}.
 \eg

2.
Assume $w(\A_1, \A_2)=u$ and  $w(\A_2, \A_3)=v$. If $l(uv)=l(u)+l(v)$, then 
\be  
\la{chain.formula.for.cartan.5.8}
h(\A_1, \A_3) =h(\A_1, \A_2) \cdot u \left(h(\A_2, \A_3)\right).
\ee

3. We have
\[
h(x_i(b)\cdot [\U^-], [\U^-]) =\alpha_{i^*}^\vee(b).
\]
\ep
\begin{proof} Part 1 is by definition. Part 2 follows from Lemma \ref{lem.part.flag.decom} and Part 1. Note that
\[
   \begin{pmatrix} 
      1 & b \\
      0 & 1 \\
   \end{pmatrix}  =    \begin{pmatrix} 
      1 & 0 \\
      b^{-1} & 1 \\
   \end{pmatrix}      \begin{pmatrix} 
      0 & 1 \\
      -1 & 0 \\
   \end{pmatrix}    \begin{pmatrix} 
      b^{-1} & 0 \\
      0 & b \\
   \end{pmatrix}
    \begin{pmatrix} 
      1 & 0 \\
      b^{-1} & 1 \\
   \end{pmatrix}.
\]
Therefore 
\be
\la{basic.exchange.formula.x--y}
x_{i}(b) = y_{i}(b^{-1}) {\overline{\overline{s}}_{ i}} \alpha_{i}^\vee (b^{-1}) y_{i}(b^{-1}).
\ee
Note that $y_i(b^{-1})$ can be absorbed in to $[\U^-]$. So 
$
\big(x_i(b)\cdot [\U^-], [\U^-]\big) =\big(\alpha_{i}^\vee(b^{-1})\cdot [\U^-], \overline{s}_{ i} \cdot [\U^-] \big). 
$
Therefore $h_-= \alpha_{i}^\vee(b^{-1})$ and $h= \alpha_{i^*}^\vee(b)$.
\end{proof}

\subsubsection{Decorated flags for the adjoint group.}

 An additive character $\psi: {\rm U} \rightarrow \mathbb{G}_a$ is   {\it non-degenerate}  
if the stabilizer of the Cartan group $\H$, acting on $\psi$ via $h\circ \psi(u):= \psi(h^{-1}uh)$, is the center of $\G$. 

Assume   that the center of $\G$ is trivial. Then the   space $\mathcal{A}_\G = \G/\U$ can be defined as a moduli space.

\bd \la{12.22.08.1} 
A {\em decorated flag} for an adjoint group  $\G$ is a pair 
$(\U, \psi)$ where $\U$ is  a maximal   
unipotent subgroup and 
 $\psi: \U \rightarrow \mathbb{G}_a$ 
  a non-degenerate additive character. 
The  space ${\cal A}_{\G}$  
parametrises decorated flags. 
Forgetting $\psi$ and taking the unique Borel subgroup containing $\U$, we get the map $\pi: {\cal A}_\G\to{\cal B}_\G$.  
\ed

\vskip 2mm \paragraph{\bf Remark.} 
The group $\G$ acts by conjugation on the pairs $(\U, \psi)$.
 This action is transitive. Since the center of $\G$ is trivial,  
 the stabilizer of  such a pair $(\U, \psi)$ is $\U$. 
 Hence we obtain an isomorphism 
$
i_{\chi}: {\cal A}_{\G}\stackrel{\sim}{\lra} {\rm G/U}. 
$ 
This isomorphism is not canonical because specifying a point in $ {\cal A}_{\G}$ requires choosing a non-degenerate character $\psi$, or better,   a {pinning}.  
Therefore, when we write ${\cal A}_{\rm G} = {\rm G}/{\rm U}$,  we abuse notation, 
keeping in mind a choice of a pinning. The coset  $[\U]$ in ${\rm G}/{\rm U}$ refers  to the standard maximal unipotent subgroup for a given pinning. 

\vskip 2mm

\subsubsection{Pinnings revisited.} In Section \ref{3.1.4a} we recalled the notion of {\it pinnings in a group $\G$}.  Let us  introduce now the crucial for this paper notion of 
a {\it pinning} for an adjoint group $\G$, which we use throughout the paper, and which is equivalent to the notion of  a {pinning in a group $\G$}.\footnote{For  any  $\G$, not necessarily adjoint, one can still use Definition \ref{PIN}.   
The  defined this way pinnings  form a 
$\G-$torsor. }

 \bd \la{PIN} Let $\G$ be any split semi-simple  adjoint group. 
 A pinning   is an ordered generic pair of decorated flags 
 $(\A_1, \A_2)$    such that 
 $$
 h(\A_1, \A_2)=1.
 $$
 \ed
 
 Denote by ${\mathcal P}$ the moduli space of pinnings. 
Then   ${\mathcal P}$ 
 is a left principal homogeneous $\G$-space: the group $\G$  acts simply transitively from the left  on the  pairs $(\A_1, \A_2)$ such that  $h(\A_1, \A_2)=1$. 
  
 There is a  surjective projection of $\G$-spaces  onto the subvariety of generic pairs of flags \be \la{AABB+}
 {\mathcal P} \to {\cal B}\times {\cal B}, \ \ \ \ (\A_1, \A_2) \lms (\B_1, \B_2):= (\pi(\A_1), \pi(\A_2)).
 \ee
 
By (\ref{lem.16.13.05.hhh.3.26}), the Cartan group $\H$ acts simply transitively on each fiber  of this map  by 
\be \la{H-action}
h: (\A_1, \A_2 )\mapsto \left(\A_1\cdot h, \A_2\cdot  w_0(h)\right).
\ee

 We often refer to a pinning  $(\A_1, \A_2)$ as a {\it pinning over a generic pair of flags $(\B_1, \B_2)$}.

\bd
\label{opposite.pinning.2020}
Let $p=(\A_1, \A_2)$ be a pinning over a  pair $(\B_1, \B_2)$. Then the opposite pinning  $p^\ast$  is 
$$
p^*:= (\A_2, \A_1\cdot s_\G).
$$
\ed 
By Proposition \ref{P2.5}, (\ref{chain.formula.for.cartan.5.8}),  the map $p \mapsto p^*$ intertwines the action of an element $h\in \H$   on $p$ with the action of $w_0(h)$   on $p^\ast$.
So  reversing the order of the decorated flags we indeed get  a new { pinning}.  
\be
p^*:=(\A_2, \A_1\cdot s_{\G}), \quad h(\A_2, \A_1\cdot s_\G)=1.
\ee

 
 \bl \la{PIN1} 
 A pinning  over a generic pair $(\B_1, \B_2)$ is the same  as 
 a pair $(\A_1, \B_2)$, where 
   $\A_1$ is a decorated flag over $\B_1$.  
 \el

\begin{proof}    Every $(\A_1, \A_2)$ provides a pair  $(\A_1, \pi(\A_2))$. Conversely,  
given a generic pair  $(\A_1, \B_2)$,  there is a unique $\A_2$ over $\B_2$ such that 
$h(\A_1, \A_2)=1$. 
\end{proof}

\bl \la{twopinnings} If the group $\G$ is adjoint, the pinning in $\G$ from Section \ref{3.1.4a} is the same thing as a pinning from Definition \ref{PIN}.
\el

\begin{proof} The equivalence between  pinnings from Section \ref{3.1.4a} and   pairs $(\A_1, \B_2)$ from Lemma \ref{PIN1} is seen as follows. 
 Let ${\rm U}\subset \B_1$ be the maximal unipotent subgroup. A choice of  pinning data $x_i$ is equivalent to a choice of a non-degenerate character $\psi$ such that $\psi(x_i(a))=a$. 
The choice of $y_i$ is determined by  $x_i$. \end{proof} 
\vskip 2mm

\medskip

\subsection{Moduli spaces ${\mathscr P}_{\G,  \bS}$ of framed $\G$-local systems  with pinnings on $\bS$, and their variants} \la{Sec2.2}
  
\medskip  In Sections \ref{Sec2.2} - \ref{SECTT4.3} the group $\G$ has trivial center. We denote by $\G'$ its universal cover.

\subsubsection{Moduli spaces revisited.}  \la{SEC3.2.1} Given a boundary interval ${\rm I}$, the flags at its ends  
 are ordered by the     boundary orientation  of $\bS$.  
Transporting them   to a common point of ${\rm I}$ we get a pair of flags. Assume that it is generic. Then a 
  {\it pinning  $p_{\rm I}$ at a boundary interval ${\rm I}$} is a pinning  over  this ordered  pair of flags. \vskip 2mm

\begin{definition} \label{2A} 
Let $\bS$ be a \underline{colored} decorated surface. 
\vskip 1mm

i) The moduli space ${\mathscr P}_{\G,  \bS}$ parametrizes data $({\cal L}, \beta; p )$ where 

\begin{itemize} 

\item ${\cal L}$ is $\G$-local system on $\bS$,

\vskip 1mm \item 
$\beta$  is a framing on ${\cal L}$, given by invariant flags near all special points  and punctures.  We  assume  that 
 \be \la{KA}
 \mbox{for each colored boundary interval ${\rm I}$, the pair of flags  at its ends   is generic}.\footnote{ Note that a decorated flag $\A_s$  at a point $s$ determines a flag  $\pi(\A_s)$  at   $s$.} 
 \ee 
\vskip 1mm

 \item $p =\{p_{\rm I}\}$ is a collection of pinnings assigned to the colored boundary intervals ${\rm I}$, 
 which are oriented according to the boundary orientation induced by the orientation of $\bS$. 
\end{itemize} 
\vskip 1mm

ii) The   moduli space ${\rm Loc}_{\G,\bS}$ parametrise the similar data, where we drop   invariant flags at the punctures.
 \ed

The moduli space ${\mathscr P}_{\G,  \bS}$ and ${\rm Loc}_{\G,  \bS}$ assigned to a colored decorated surface $\bS$ has a natural stratification, obtained by assigning to each non-colored 
boundary interval ${\rm J}$ of an element $w_{\rm J}$, and postulating that the two flags at the end of the boundary interval are in the $w-$position $w_{\rm J}$:
\be
\begin{split}
&{\mathscr P}_{\G,  \bS}= \coprod_{w_{\rm J} \in W}{\mathscr P}_{\G,  \bS; w_{\rm J}}, \qquad \quad 
 {\rm Loc}_{\G,  \bS}= \coprod_{w_{\rm J} \in W}{\rm Loc}_{\G,  \bS; w_{\rm J}}.\\
\end{split}
\ee

 Denote by   $T\bS-\{0\}$ the tangent bundle  to $\bS$ with the zero section removed.  Recall 
the punctured boundary $\widehat \partial \bS$ of $\bS$, see (\ref{PBDSS}). Recall the central element $s_{\G'}$, see Section \ref{3.1.5}. One has $s_{\G'}^2=1$. 

  \bd  \la{DEFP3}    Let $\bS$ be a \underline{colored} decorated surface. The moduli space ${\mathscr A}_{  \G',  \bS}$ parametrizes   $ \G'$-local systems ${\cal L}$ on 
  $T\bS-\{0\}$  with the monodromy $s_{  \G'}$ 
around a generator of  $\pi_1(T_x\bS-\{0\})$,  and a flat section of the 
local system ${\cal L}_{\cal A}:= {\cal L}\times_{ \G'} {\cal A}$ near each special point or puncture. 

Moreover, for each \underline{non-colored} boundary interval ${\rm J}$ with the boundary points $(s,t)$ we have $h(\A_s, \A_t)=1$, where $\A_s$ is a decorated flag in the fiber  over the tangent vector at $s$, which follows the orientation of $\partial \bS$.
  \ed  
   
  \begin{figure}[ht]
\begin{center}
\begin{tikzpicture}[scale=.5]
 \draw (-1,-.5) arc (3.3:70:3);
  \draw (1,-.5) arc (176.7:110:3);
   \draw (1.5,2.5) arc (-20:-160:1.6 and .6);
   \draw (0.75,1) arc (-60:-120:1.5);
   \draw (0.6,.92) arc (60:120:1.2);
   \draw[dashed] (1, -.5) arc (0:180:1 and .4);
      \draw (1, -.5) arc (0:-180:1 and .4);
       \draw[dashed] (-1.5, 2.5) arc (0:360:1 and .4);
       \draw[dashed] (1.5, 2.5) arc (180:540:1 and .4);
       \node[red] at (0, -1.3) {\footnotesize $\A$};
       \node[red] at (0, -0.9) {\footnotesize $\bullet$};
       \node[blue] at (0, 0.3) {\footnotesize $\bullet$};
       \node[red] at (-0.3, .3) {\footnotesize $\A$};
        \node[blue] at (0.9, 1.6) {\footnotesize $\bullet$};
         \node[red] at (0.6, 1.6) {\footnotesize $\A$};
        \node[red] at (-3.5, 2.5) {\footnotesize $\bullet$};
         \node[red] at (-3.8, 2.5) {\footnotesize $\A$};
        \node[red] at (3.5, 2.5) {\footnotesize $\bullet$};
         \node[red] at (3.8, 2.5) {\footnotesize $\A$};
         \node[red] at (-1.9, 2.8) {\footnotesize $\bullet$};
          \node[red] at (1.9, 2.8) {\footnotesize $\bullet$};
          \node[red] at (-1.8, 3.1) {\footnotesize $\A$};
          \node[red] at (1.8, 3.1) {\footnotesize $\A$};
\end{tikzpicture}
\end{center}
\caption{ The moduli space ${\mathscr A}_{\G, \bS}$ for a colored decorated surface  $\bS$, without non-colored boundary intervals,  parametrizes $\G$-local systems  ${\cal L}$ on $\bS$ 
 with an invariant    flag $\A$ near each special point or puncture.}
\label{amalgam11}
\end{figure}

\bl  \la{EQQ}
Let $\bS$ be  a decorated surface   without  non-colored  boundary intervals.     
The moduli space   ${\mathscr P}_{\G, \bS}$  parametrizes $\G$-local systems  ${\cal L}$ on $\bS$  with an invariant flag     $\B_p$ 
  near each puncture $p$, and   a decorated flag $\A_s$  near each special point $s$,  
so that the pairs of decorated flags at the ends of   boundary intervals are generic.     \el

\begin{proof}  
The   space ${\mathscr P}_{\G,  \bS}$ from the Lemma parametrizes   $\G$-local systems ${\cal L}$ on $\bS$ where instead of the pinnings   there are  
  decorated flags $\A_s$  at    all special points  $s$.   In  Definition \ref{2A} we see a pair of decorated flags $(\A_s^-, \A_s^+)$ near each special point $s$ above the same flag $\B_s$, and the boundary is oriented by  $\A_s^-\to  \A_s^+$. There is a bijection between these  data,   illustrated on Figure \ref{amalgam3}, such that  
\be \la{AAaa}
 \A_s= \A_s^+.
\ee
 \begin{figure}[ht]
\epsfxsize80pt
\center{ \begin{tikzpicture}[scale=.6]
\draw[]  (0,0) -- (6,0);
\node[label=below:{$\A_s$}] at (1,0) {$\bullet$};
\node[label=below:{$\A_t$}] at (5,0) {$\bullet$};
\draw[]  (10,0) -- (16,0);
\node[label=above:{$\B_s$}] at (11,0) {$\bullet$};
\node[label=above:{$\B_t$}] at (15,0) {$\bullet$};
\node[label=below:{$\A_s^-$}] at (10.3,0) {};
\node[label=below:{$\A_s^+$}] at (11.7,-0.1) {};
\node[label=below:{$\A_t^-$}] at (14.3,0) {};
\node[label=below:{$\A_t^+$}] at (15.7,-0.1) {};
\draw[dashed, latex-latex] (7,0) -- (9,0);
\draw[latex-latex] (12.5,-.7) -- (13.5, -.7);
\node at (13, -.4) {1};
\end{tikzpicture}}
 \caption{A generic pair of decorated flags $(\A_s, \A_t)$ is the same thing as two pairs of decorated flags $(\A_s^-, \A_s^+), (\A_t^-, \A_t^+)$ such that $h(\A_s^+, \A_t^-)=1$, 
 and    $(\A_i^-, \A_i^+)$ decorate the same  flag $\B_i$. }
\label{amalgam3}
\end{figure} 
 Namely, for each special point $s$, we move 
 the  decorated flag $\A_s$  a bit  {\it to the right}, i.e.  
following the boundary orientation, getting a decorated flag  $\A_s^+$.   Let  $t$ be the endpoint of the     boundary interval  ${\rm I}_s$ starting   at $s$. 
Then there is unique decorated flag  $\A_t^-$, placed a bit to the left of  $t$, such  that 
 $h(\A_s^+, \A_t^-)=1 \in \H$, and the underlying flags  for   $\A^-_t$  and  $\A_t$ coincide. The  pair  
 $(\A_s^+, \A^-_t)$ is a {pinning} at the boundary interval ${\rm I}_s$. 
\end{proof}

A few remarks are in order  
  \begin{enumerate}

\vskip 1mm\item  If $\bS$ has neither punctures, nor non-colored boundary segments, the moduli spaces ${\mathscr A}_{\G', \bS}$ and ${\mathscr P}_{\G, \bS}$ are almost isomorphic. Indeed, via Lemma \ref{EQQ},   the projection $\G' \to \G$ provides a finite cover 
$$
\xi: {\mathscr A}_{\G', \bS} \to {\mathscr P}'_{\G, \bS},
$$
  where ${\mathscr P}'_{\G,  \bS}$ is the moduli space obtained by dropping  generality conditions on the pairs of flags at   boundary intervals.     
In this case the space ${\mathscr A}_{\G', \bS}$ carries, in addition to the $K_2-$structure,  a \underline{canonical} $\Gamma_\bS-$equivariant 
Poisson structure. This mechanism provides   very important cluster $K_2-$varieties with canonical Poisson structures.   An example,  is given by  the group $\G'$, see Section \ref{SEC4.4}. The  $2$-form provided by the $K_2-$structure and the Poisson structure are typically degenerate, and thus are not inverses of  each other.  \vskip 1mm

\vskip 1mm\item Given a colored boundary interval ${\rm I}$ on   $\bS$,   making the  interval ${\rm I}$ non-colored we get a new colored decorated 
surface $\bS'$. Then the dimensions of the related moduli  spaces drop by ${\rm rk}  \G$:
\be
\begin{split}
&{\rm dim}{\mathscr P}_{\G, \bS'} ={\rm dim}{\mathscr P}_{\G, \bS} - {\rm rk}  \G;\\
& {\rm dim}{\mathscr A}_{\G', \bS'} ={\rm dim}{\mathscr A}_{\G', \bS} - {\rm rk}  \G.\\
\end{split}
\ee
Indeed, in the first case we remove the pinning $p_{\rm I}$, and in the second we impose the condition $h(\A_s, \A_t)=1$, where $s,t$ are the left and right 
boundary vertices of the ${\rm I}$. 

For any colored decorated surface $\bS$ we have 
\be
\begin{split}
&{\rm dim}{\mathscr P}_{\G, \bS} = {\rm dim}{\mathscr A}_{\G', \bS}.\\
\end{split}
\ee

\vskip 1mm \item For any colored decorated surface $\bS$ the moduli spaces ${\mathscr A}_{\G', \bS}$ and ${\mathscr P}_{\G, \bS}$ are dual to each other. In particular Theorem \ref{MTHa} 
is valid for any such $\bS$. For example,   the pair  $({\mathscr A}_{\G', \bS}, {\mathscr P}_{\G, \bS})$ gives rise to a cluster ensemble. 
The proof for general  $\bS$ is reduced to Theorem \ref{MTHa}. Indeed, the boundary intervals carry frozen variables; 
 dropping them reduces the clusters, but does not affect the mutations we use. 

\vskip 1mm\item  The better way to picture a colored decorated surface $\bS$ is  the {\it open string format}. We  
expand  special points  to  
{\it black boundary  intervals},   call the colored 
intervals between them {\it red boundary intervals}, the non-colored boundary intervals just {\it boundary intervals}, and   expand punctures to  {\it black boundary circles}. Then a boundary component   is either 
 a polygon with    black/red/non-colored sides,  or  a black circle.  
The    boundary is   used as follows: 

\begin{itemize}

\item
We put  a  flat section of the local systems ${\cal L}_{\cal A}$ or ${\cal L}_{\cal B}$,  
over each black    interval/circle.

\vskip 1mm\item
The pinnings are at the red boundary intervals. 
We glue 
surfaces along the red intervals. 

\vskip 1mm\item For the ${\cal A}-$space, the pair of flags $(\A_s, \A_t)$ over each non-colored interval satisfies $h(\A_s, \A_t)=1$. 
\end{itemize}

Yet to simplify the pictures we  stick to presenting decorated surfaces using marked points. \vskip 1mm
 
\vskip 1mm\item    For each boundary interval ${\rm I}$, the    action  (\ref{H-action}) of the group $\H$  on   pinnings at ${\rm I}$ gives 
    a free  action 
\be \la{ALa}
    \tau_{\rm I}: {\H} \times{\mathscr P}_{\G, \bS} \lra {\mathscr P}_{\G, \bS}.
\ee
If the group $\G$ is adjoint, and the pinning $p_{\rm I}$ at ${\rm I}$ is given by a pair $(\A_1, \A_2)$, then, according to (\ref{H-action}),  the action of an element $h \in {\rm H}$ is given by altering this pinning 
$$
\tau_{\rm I}(h): (\A_1, \A_2) \lms (\A_1\cdot h, \A_2\cdot w_0(h)),
$$ and keeping intact the rest of the data. 
  \end{enumerate}

  \subsubsection{Example.} A triangle $t$ topologically is a disc with three special  points.    
The   moduli space ${\mathscr P}_{\G, t}$  parametrises 
triples of decorated flags $(\A_1, \A_2, \A_3)$ that are pairwisely generic, or what is the same by Lemma \ref{EQQ}, triples of flags with pinnings, considered modulo the diagonal $\G$-action: 
\be \la{12.31.08.4}
{\mathscr P}_{\G, t}:= \left\{(\A_1, \A_2, \A_3)\right\}\left/\G\right.  \stackrel{\ref{EQQ}}{=} 
\left\{(\B_1, \B_2, \B_3; p_{12}, p_{23}, p_{31})\right\}\left/\G\right..
\ee
Here $\B_1, \B_2, \B_3$ are pairwisely generic, and $p_{ij}$ is a pinning 
over $(\B_i, \B_j)$.   We picture 
an element of $\mathscr{P}_{\G, t}$  by oriented triangles, 
with the flags assigned to vertices, 
and arrows showing the pinnings, see Figure \ref{fga32}. 
Reversing an arrow means replacing the corresponding 
pinning by the opposite one in Definition \ref{opposite.pinning.2020}.

\begin{figure}[ht]
\epsfxsize 200pt
\center{
\begin{tikzpicture}[scale=0.5]
\tikzset{>=open triangle 45}
\foreach \count in {1,2,3}
{
\node at (-30+120*\count:2.5) {\small ${\A}_\count$};
\draw [directed]  (210-120*\count:2cm) -- (330-120*\count:2cm);
}
\node[blue] at (150:1.5) {\small };
\node[blue] at (270:1.5) {\small  };
\node[blue] at (30:1.5) {\small };
\end{tikzpicture} \qquad \qquad \qquad 
\begin{tikzpicture}[scale=0.5]
\tikzset{>=open triangle 45}
\foreach \count in {1,2,3}
{
\node at (-30+120*\count:2.5) {\small ${\B}_\count$};
\draw [directed]  (210-120*\count:2cm) -- (330-120*\count:2cm);
}
\node[blue] at (150:1.5) {\small $p_{12}$};
\node[blue] at (270:1.5) {\small $p_{23}$};
\node[blue] at (30:1.5) {\small $p_{31}$};
\end{tikzpicture}
 }

\caption{Two incarnations of the  moduli space ${\mathscr P}_{\G, t}$.}
\label{fga32}
\end{figure}

Let us pick a standard pinning $p=(\B, \B^-, x_i, y_i; i\in {\rm I})$ and set $\B_\ast = \B \cap \B^- w_0 \B^-$. 
\bl \la{1.2.09.1}
Given a triangle $t$, a choice of its vertex $v$   
provides an isomorphism 
$$
i_v: {\mathscr P}_{\G, t} \stackrel{\sim}{\lra} \B_\ast \times \H.
$$ 
\el

\begin{proof}  
Every $\G$-orbit in \eqref{12.31.08.4} contains a unique representative such that $p_{12}$ is the standard pinning $p$ chosen. The choice of a third flag $\B_3$ and a pinning $p_{31}$ over $(\B_3, \B_1)$ one-to-one corresponds to a choice of $b\in \B_\ast$ such that 
$p_{31}= b \cdot p_{12}^\ast$.
Here $b\in \B_\ast$ since the pair $(\B_2, \B_3)$ is required to be generic. The choice of $p_{23}$ corresponds to  $h\in {\rm H}$ such that its associated character $\psi(p_{23})= h \cdot \psi(p_{12}^\ast)$.
\end{proof}

\subsubsection{Gluing revisited.}  
Gluing a pair of colored boundary intervals ${\rm I}_1, {\rm I}_2$ on $\bS$, we get   a new colored   decorated
surface $\bS'$, see Section \ref{1.1.3}.

    \bl \la{GLUINGM} There exists  a canonical  regular dominant map, which we call the gluing map: 
    \be \la{AMAX}
\gamma_{{\rm I}_1, {\rm I}_2}: {\mathscr P}_{\G, \bS }   \lra {\mathscr P}_{\G, \bS'}.
\ee 
Its 
  fibers are principal $\H$-bundles. The   map  (\ref{AMAX}) admits  a factorization 
 $
{\mathscr P}_{\G, \bS }  \lra {\mathscr P}_{\G, \bS} /\H \hlra {\mathscr P}_{\G, \bS'}.   
 $
Here the first map is the quotient by the   action  $(\tau_{\rm I_1}, \tau_{\rm I_2^{op}})$ of  $\H$ in (\ref{ALa}).  
The second     is   injective.
\el 
 
 \begin{proof}  
 We have already explained in Section \ref{1.1.3} how to define the gluing map. A few more comments are necessary. 

  Since in the   adjoint group $\G$ we have $ \overline w_0^2=1$, 
  identifying the pinnings $p_{{\rm I}_1}$ and $p^*_{{\rm I}_2}$ amounts to  the same result as identifying the pinnings $p_{{\rm I}_1}^\ast $ and $p_{{\rm I}_2}$.
    
  Since  pinnings with the same underlying pair of flags form an $\H$-torsor, we get an  $\H$-torsor over the image. 
Any framed $\G$-local system on $\bS'$ whose restriction to the edge ${\rm I}$ is given by a generic pair of flags lies in the image of the gluing map. 
 \end{proof}

Given an {ideal triangulation} ${\cal T}$ of a decorated surface $\bS$, we can cut    $\bS$ into a collection 
 $\{t\}$ of   ideal triangles of   ${\cal T}$. Conversely,   a pair  $(\bS, {\cal T})$    is just  a surface glued from triangles. 
   The   {\it edges} of the triangulation ${\cal T}$ are the   internal edges, that is, the ones  we use to cut   $\bS$ into triangles. They are   
  identified with the   pairs of glued sides of   triangles. The unglued sides of   triangles are identified with   boundary intervals of $\bS$. 
  For example, an ideal triangulation of a rectangle has one edge and four boundary intervals. 
 
    The gluing maps allow us to 
 assemble the   space ${\mathscr P}_{\G, \bS}$ from the   spaces ${\mathscr P}_{\G, t}$ assigned to the triangles.  
This way  we arrive at the following key result.

\begin{theorem} \label{1.8.09.1}  Let $\bS$ be a decorated surface, and ${\cal T}$ an ideal triangulation  of $\bS$. 
Then 
the   space ${\mathscr P}_{\G, \bS}$ is birationally equivalent to the one obtained by gluing   the 
  spaces ${\mathscr P}_{\G, t}$ assigned to the triangles 
$t$ of  ${\cal T}$.  Moreover, there is a  canonical Zariski open embedding 
\be \la{GM}
\begin{split}
& \gamma_{\cal T}: \prod_{t}  {\mathscr P}_{\G, t}\left /{\rm H}^{\{\mbox{{\small \rm edges of} ${\cal T}$}\}}\right.  \stackrel{ }{\hlra}{\mathscr P}_{\G, \bS}.\\
\end{split}
\ee
 In particular, the moduli space ${\mathscr P}_{\G, \bS}$ is rational.   
\end{theorem}

\begin{proof} The first  claim follows from (\ref{GM}). 
The moduli space ${\mathscr P}_{\G, t}$ is rational by Lemma \ref{1.2.09.1}. 
The gluing preserves rationality. This implies the rationality of 
the space ${\mathscr P}_{\G, \bS}$. \end{proof}

We describe the cluster structure of ${\mathscr P}_{{\rm G,} \bS}$ via the  amalgamation  for   $\G={\rm PGL_2}$  in   Section \ref{SEC5.100}, and 
for $\G ={\rm PGL_m}$ in Section \ref{SSECC3.2}. 
Section \ref{SSECC3}  is independent of the rest of the paper.

 \medskip

\subsection{The potential functions $ {\cal W}_{\alpha_i} $} \la{POTEN} \medskip

Let us recall the definition of the potential function \cite[Section 1.1.1]{GS13}. 
Consider an open subvariety  
\be \la{ABA}
{\rm Conf}^*({\cal A}, {\cal B}, {\cal B}) \subset {\rm G} \left\backslash ({\cal A}\times {\cal B} \right.
\times  {\cal B})
\ee 
parametrizing $\G$-orbits   of triples $ \{{\A},  {\rm B}_{1}, {\rm B}_{2}\}$ such that  $\{{\A}, {\rm B}_1\}$ and $\{{\A}, {\rm B}_2\}$ are generic. 

Let us define for each simple positive root $\alpha_i$   a regular function, called {\it partial potential}
$$
{\cal W}_{\alpha_i}: {\rm Conf}^*({\cal A}, {\cal B},  {\cal B})   \lra {\mathbb A}^1. 
$$
Given such a triple $ \{{\A},  {\rm B}_{1}, {\rm B}_{2}\}$, 
there is a unique $u$ in the stabilizer  ${\rm U}_{{\B}}$ of $\A$ such that
$
\{ {\A},{\rm B}_{2}\} = u \cdot \{ {\A}, {\rm B}_{1}\}. 
$

Next,  the decorated flag $\A $ provides a   character  
$$
\psi_{\A}: {\rm U}_{\B} \lra {\mathbb A}^1.
$$
Indeed, for the  adjoint group    $\A = ({\rm U}_{\B}, \psi_{\A})$ by Definition \ref{12.22.08.1}. 
Otherwise we use the  projection ${\cal A}_{\widetilde{\G}} \to {\cal A}_{\G}$  to get the character. 
Denote by $\mathfrak{u}$ the Lie algebra of $\U$. The $\H$-module 
$\mathfrak{u}/[\mathfrak{u}, \mathfrak{u}]$ is a direct sum of  one dimensional $\H$-modules corresponding to 
simple positive roots. So any additive character   of $\U$ is decomposed canonically  into  a sum of characters   parametrised by   positive simple roots $\alpha_i$. In our case:
$$
\psi_{\A} =  \sum_{i \in {\rm I}} \psi_{\A, \alpha_i}.
$$

\bd The  function ${\cal W}_{\alpha_i}$ assigned to a simple positive root $\alpha_i$ is given by 
\be \la{POTF}
{\cal W}_{\alpha_i}({\A},  {\rm B}_{1}, {\rm B}_{2}):= \psi_{\A, \alpha_i}(u).  
\ee
\ed

The function ${\cal W}_{\alpha_i}$ is invariant under diagonal $\G$-action, and descends to the space  \eqref{ABA}. 
 \medskip 
 
 \subsection{The principal affine space and pinnings: examples} \la{EXAMPLES}
 
\medskip

\subsubsection{$\G = {\rm SL}_2$} Take a two dimensional vector space $V_2$ with a volume form $\omega$. Then ${\rm SL}_2 = {\rm Aut}(V_2, \omega)$, ${\cal A}_{{\rm SL}_2}= V_2-\{0\}$, and ${\cal B}_{{\rm SL}_2}={\Bbb P}(V_2)$. 
A configuration in ${\rm Conf}^\ast({\cal A}, {\cal B},{\cal B})_{{\rm SL}_2} $ consists of a nonzero vector $v$ and two lines ${\rm L}_1, {\rm L}_2$ not containing $v$.  The angle invariant $u$ of $(v, {\rm L}_1, {\rm L}_2)$ is an unipotent element of ${\rm SL}_2$ that fixes $v$ and maps ${\rm L}_1$ to ${\rm L}_2$. We denote it by
\be
\la{example.sl2.unipotent.potential}
u =   \begin{pmatrix} 
      1 & b \\
      0 & 1 \\
   \end{pmatrix}  \in {\rm U}_{{\rm SL}_2}.
\ee 
Let  $l_1\in {\rm L}_1$ and $l_2\in {\rm L}_2$ be vectors such that
$
\omega (v \wedge l_1) = \omega (v \wedge l_2) =1.
$
Then  
$ 
l_2-l_1 = b v$, where $b$ is in \eqref{example.sl2.unipotent.potential}. Note that
$b=\omega (l_2\wedge l_1)$.  
\begin{figure}[ht]
\epsfxsize80pt
\label{potential.sl2ax}
\center{ \begin{tikzpicture}[scale=0.6]
\draw[]  (0.25,-1) -- (-0.75,3);
\draw[] (-0.5,-1) -- (1.5,3);                   
\draw[directed,thick] (0,0) -- (2,0); 
\draw[directed,thick] (0,0) -- (1,2);
\draw[directed,thick] (0,0) -- (-0.5,2);
\draw[directed,  blue ] (-0.5,2) -- (1,2);
\node [label=below:{$v$}] at (1, 0.2) {};
\node [blue] at (0.25, 2.4) {$bv$};
\node [label=left: {$l_1$}] at (0,1) {};
\node [label=above: {$L_1$}] at (-1,2.1) {};
\node [label=left: {$l_2$}] at (1.75,1) {};
\node [label=above: {$L_2$}] at (1.85,2.1) {};
 \end{tikzpicture}}
\end{figure}
The {\it potential} on ${\rm Conf}^\ast ({\cal A}, {\cal B},{\cal B})_{{\rm SL}_{2}} $ is a regular function
\begin{align}
{\cal W}:~ &{\rm Conf}^\ast ({\cal A}, {\cal B},{\cal B})_{{\rm SL}_2}  \lra {\Bbb A}^1,\nonumber\\
\la{sl2.potential.6.3}
(v, {\rm L}_1, {\rm L}_2) &\lms b:=\frac{\omega(l_1\wedge l_2)}{\omega(v\wedge l_1) ~\omega(v\wedge l_2)}, \hskip 4mm   l_1\in {\rm L}_1,   l_2\in {\rm L}_2.
\end{align}
It does not depend on non-zero vectors 
 $l_1\in {\rm L}_1$ and $l_2\in {\rm L}_2$ chosen. 
 
 \bl 
\la{6.30.17.43.hh}
The potential function is additive: 
$$
  {\cal W}(v, {\rm L}_1, {\rm L}_2)+ {\cal W}(v, {\rm L}_2, {\rm L}_3)={\cal W}(v, {\rm L}_1, {\rm L}_3) .
$$
\el
\begin{proof} It follows directly from $l_2-l_1=bv$.
\end{proof} 
 
\subsubsection{$\G = {\rm PGL_2}$} \la{3.4.2} Let  $V_2$ be a 
two dimensional vector space.  The principal affine  space    ${\cal A}_{\rm GL_2}$
parametrizes   pairs 
$(v, \omega)$ where $v \in V_2-\{0\}$, and $\omega \in \wedge^2V_2^*-\{0\}$. The   group ${\mathbb G}_m$ acts on ${\cal A}_{\rm GL_2}$  by 
$$
(v, \omega) \lms (tv, t^{-2}\omega).
$$ 
The principal affine  space    ${\cal A}_{\rm PGL_2}$ parametrizes
 the orbits $[v, \omega]$ of ${\mathbb G}_m$   on ${\cal A}_{\rm GL_2}$.\vskip 1mm
 
Let $\U_v$ be the unipotent subgroup stabilising $v$. 
An element    $[v, \omega] \in {\cal A}_{\rm PGL_2}$ provides a map 
$$
\lambda_{v, \omega}: {\mathbb G}_a \to \U_v, \qquad  t\lms \lambda_{v, \omega}(t), \qquad  
\lambda_{v, \omega}(t): l \lms l+ t\omega(l, v)v.
$$

\vskip 2mm The   moduli space ${\rm Conf}_2({\cal A}_{\rm PGL_2}) $ parametrizes pairs $([v_1, \omega_1], [v_2, \omega_2])$ modulo the left diagonal action of 
${\rm PGL_2}$. There is a canonical map to the Cartan group $\H$ of ${\rm PGL_2}$:
\be \la{98}
h: {\rm Conf}_2({\cal A}_{\rm PGL_2}) \lra \H={\Bbb G}_m, \qquad h([v_1, \omega_1], [v_2, \omega_2]) := \omega_1(v_1, v_2) \omega_2(v_1, v_2). 
\ee
The Cartan group $\H$   acts on ${\cal A}_{\rm PGL_2}$  on the right by
$$
[v, \omega] \lra  [tv, t^{-1} \omega] = [v, t \omega].
$$
Composing with the $\H$-action, the   map (\ref{98}) is rescaled by $t$. \vskip 1mm

{\it The potential ${\cal W}_\alpha$ for ${\rm PGL}_2$.} It   
is a function on the triples $([v, \omega],  {\rm L}_1, {\rm L}_2)$, where ${\rm L}_1, {\rm L}_2$ are lines in $V_2$:
\be \la{98++}
\begin{split}
&{\cal W}_\alpha: {\rm PGL_2} \backslash \left( {\cal A}_{\rm PGL_2} \times {\cal B} \times   {\cal B}\right) \lra {\Bbb A}^1,\\
&{\cal W}_\alpha= \frac{\omega(l_1, l_2)}{\omega(v, l_1)\omega(v, l_2)},\hskip 4mm   l_1\in {\rm L}_1, \   l_2\in {\rm L}_2.\\
\end{split}
\ee
It is invariant under the action
 $(v, \omega) \lms (tv, t^{-2}\omega)$.\vskip 2mm
 
 Before we proceed further, let us recall the cross-ratio. Set ${\mathbb P}^1:= P(V_2)$. 
Denote by $x'$ a vector in $V_2$ projecting to a point $x$. 
Pick  a volume form $\omega_2$  in $V_2$.  Recall the cross ratio of four points $x_1, x_2, x_3, x_4 \in {\mathbb P}^1$:
\be \la{CRRA}
r(x_1, x_2, x_3, x_4):= 
\frac{\omega (x'_1,x'_2)\omega(x'_3,x'_4)}{\omega(x'_1,x'_4)\omega(x'_2,x'_3)}.
\ee
The cross ratio has the following properties:
\[
r(x_1, x_2, x_3, x_4) = r(x_2, x_3, x_4, x_1)^{-1} = r(x_4, x_3, x_2, x_1)^{-1}.
\]

\vskip 2mm
{\it Pinnings for ${\rm PGL}_2$.}  A pinning $\pi$ in the ${\rm PGL_2}$ can be defined in one of the following  equivalent ways:

\begin{enumerate}

\vskip 1mm\item A pair of decorated flags $(\A_1, \A_2)$ in generic position such that $h(\A_1, \A_2)=1$.

\vskip 1mm \item   A projective frame in $V_2$, that is a frame $(e_1, e_2)$ considered up to rescaling $(e_1, e_2) \lra (te_1, te_2)$.

\vskip 1mm\item An ordered triple of non-collinear vectors $(e_1, e_2, e_p)$ in $V_2$,  considered up to a rescaling  
$(e_1, e_2, e_p)\lra (te_1, te_2, te_p)$, such that 
$$
e_1+e_2+e_p=0.
$$

\item  An ordered triple of distinct points $\pi = (x_1, x_2; p)$ 
on ${\mathbb P}^1$.


\end{enumerate}

The equivalence of these descriptions is seen as follows. 


\vskip 1mm  $1) \rightarrow 2)$.  Given a pair $\A_1= [v_1, w_1]$ and $\A_2=[v_2, w_2]$ such that $h(\A_1, \A_2)=1$ for map (\ref{98}),  we define
$$
e_1:= w_1(v_1,v_2) v_1, \ \ e_2:= v_2.
$$

\vskip 1mm $2) \rightarrow 1)$. A pair of non-collinear vectors $(e_1, e_2)$ provides a  volume form $\omega$ such that $\omega(e_1, e_2) =1$. 
We set $\A_1 := [e_1, \omega]$, $\A_2 :=[e_2, \omega]$. Then  $h(\A_1, \A_2)=1$. Rescaling $(e_1, e_2)\to (te_1, te_2)$ leads to equivalent 
$\A_1$, $\A_2$.

\vskip 1mm $2) \leftrightarrow 3)$. 
Obvious.   

\vskip 1mm $3) \leftrightarrow 4)$.  Given a triple of distinct points  $(x_1, x_2; p)$, there is a unique up to rescaling triple of  vectors $(e_1, e_2, e_p)$ 
with $e_1+e_2+e_p=0$ so that the points  $(x_1, x_2; p)$ are  projectivisations of the vectors $(e_1, e_2, e_p)$.


\vskip 2mm

{\it Dual pinnings for ${\rm PGL}_2$}. Accordingly, there are several ways to express the definition of the dual pinning.

\begin{enumerate}

\vskip 1mm\item Given a pinning  $(\A_1, \A_2)$   such that $h(\A_1, \A_2)=1$, the dual pinning is $(\A_2, \A_1)$, by Definition \ref{opposite.pinning.2020}.

\vskip 1mm\item   For a pinning given by a  projective frame $(e_1, e_2)$ in $V_2$, the dual pinning is  $(e_2, -e_1)$.


\vskip 1mm\item  For a pinning $(x_1, x_2; p)$, the dual pinning is 
$ (x_2, x_1; p^*)$
 where $p^*\in {\mathbb P}^1$   the unique point  such that 
\be \la{12.21.08.1x++}
r(x_1,p, x_3,p^*)=1.
\ee

\end{enumerate}

So a pair of decorated flags $\A_1= [v_1, w_1]$ and $\A_2=[v_2, w_2]$, such that $h(\A_1, \A_2)=1$ for the map (\ref{98}), provides  a pinning $\pi = ([v_1], [v_2]; [p])$ 
 and the dual pinning $\pi^* = ([v_2], [v_1]; [p^*])$ where 
 $$
p= w_1(v_1,v_2) v_1 + v_2, \ \ \ p^*= w_2(v_2, v_1) v_2 + v_1.
$$
 
  \vskip 2mm
  
  \subsubsection{$\G = {\rm PGL_2}$ revisited.} Assuming we work over an algebraically closed field, here is a simpler variant. Take a two dimensional vector space $V_2$ with a non-zero area form $\omega$. 
  Then the principal affine space ${\cal A}_{{\rm PGL}_2}$ parametrises  vectors $v \in V_2 - \{0\}$ considered up to equivalence $v \sim -v$. We denote by $[v]$ the equivalence class. 
  The potential is given by the formula (\ref{98++}). The canonical map (\ref{98}) is 
  \be \la{98a}
h: {\rm Conf}_2({\cal A}_{\rm PGL_2}) \lra {\Bbb G}_m, \qquad h([v_1], [v_2]) := \omega(v_1, v_2)^2. 
\ee 
The equivalence between pinnings and triples of distinct points on ${\rm P}^1$ is seen as follows. Take a pinning understood as a pair $([v_1], [v_2])$ such that $\omega(v_1, v_2) = \pm 1$. 
Pick a representative $v_1$ of $[v_1]$. Then there is a unique $v_2$ such that $\omega(v_1, v_2) =1$. We set $p= -v_1-v_2$, getting the triple of points obtained by projectivisation of 
$(v_1, v_2, p)$. The dual pinning is given by $(v_2, -v_1)$. So $p^* = -v_2+v_1$. 

  \vskip 2mm
\subsubsection{$\G= {\rm PGL_m}$}  Let   $V_m$ be an $m$-dimensional   space.     The  space ${\cal A}_{\rm GL_m}$ parametrizes 
   pairs   
$(\rF_{\bullet}, \{v_i\})$,  where 
$$
 \rF_{\bullet} = \{0 \subset  \rF_0 \subset 
\rF_1 \subset \rF_2 \subset \ldots \subset \rF_{m} = V_m\}, \qquad
v_i \in \rF_i/\rF_{i-1} - \{0\}, \qquad i=1, ..., m.
$$
Decorated flags for ${\rm PGL_m}$ 
are  the orbits of the multiplicative group ${\mathbb G}_m$ acting by 
$v_i \lms tv_i$ on these pairs. For $m=2$ we recover the definition  above:
the form $\omega$ is the dual to $v_1\wedge v_2$. 

Let $\U$ be the subgroup stabilising a decorated flag. 
Then there are natural maps  
$$
x_i: {\mathbb G}_a \to \U/[\U,\U], \qquad  x_i(t): { \rF_{i+1}/\rF_{i-1}\lra \rF_{i+1}/\rF_{i-1}}, 
$$
$$
x_i(t): l\in { \rF_{i+1}/\rF_{i-1}} \lms 
\Bigl(l+ t\frac{l\wedge v_i\wedge \ldots  \wedge v_1}
{v_{i+1}\wedge v_{i}\wedge \ldots \wedge v_1}v_i\Bigr) ~{\rm mod} ~\rF_{i-1}.
$$

\vskip 2mm
{\it Pinnings for ${\rm PGL}_m$.}  A pinning $\pi$ in the ${\rm PGL_m}$ can be defined in one of the following  equivalent ways:

\begin{enumerate}

\vskip 1mm\item A pair of decorated flags $(\A_1, \A_2)$ in generic position such that $h(\A_1, \A_2)=1$.

\vskip 1mm\item   A projective frame in $V_m$, that is a frame $(e_1, \ldots , e_m)$  up to rescalings $(e_1, \ldots , e_m) \lra (te_1, \ldots , te_m)$.

\vskip 1mm\item A  set of vectors $(e_1, \ldots , e_m; e_p)$ in $V_m$,    up to rescalings  
$(e_1, \ldots , e_m; e_p)\lra (te_1, \ldots , te_m; te_p)$, where the first $m$  form a basis, and 
$$
e_1+\ldots + e_m+e_p=0.
$$

\vskip 1mm\item  A collection of  $(m+1)$ points $\pi = (x_1, \ldots , x_m; p)$ 
on ${\mathbb P}^{m-1}$, such that no $m$ of them lie in a hyperplane.


\end{enumerate}

Note that   given a generic collection of   points  $(x_1, \ldots , x_m; p)$, there is a unique up to rescaling collection of  vectors $(e_1, \ldots , e_m; e_p)$ lifting these points such that 
  $e_1+\ldots + e_m+e_p=0$.

\vskip 2mm

Given a pinning defined by  a  projective frame $(e_1, \ldots , e_m)$, the dual pinning is  given by
$$
((-1)^{m-1}e_m, (-1)^{m-2}e_{m-1}, \ldots ,  e_1).
$$ 
 
 \medskip
 
 \section{Key features and canonical functions on the moduli space ${\mathscr P}_{\G, \bS}$} \la{KFMS}

 \medskip

 In Section \ref{KFMS}   we introduce  a   collection of special regular functions on the moduli spaces ${\mathscr P}_{\G, \bS}$,   invariant under the mapping class group  action. 
These functions play the crucial role in our constructions. 
We stress that  most of them, most notably the outer monodromy and the   potential functions, 
 do not even exist for the spaces ${\mathscr X}_{\G, \bS}$. 
 
 All of these functions admit canonical quantum deformation, discussed in Section \ref{SECT3.4}. 
 This way they gives rise to the embeddings of the quantum groups to the algebras ${\cal O}_q({\mathscr P}_{\G, \bS})$. 
 
 Surprisingly, the potential functions play also a key role, as the Landau-Ginzburg potentials,  in the mirror symmetry conjectures \cite{GS13}, discussed briefly in Section \ref{SECTT4.3}. 
 
    \subsection{Key features of the moduli space ${\mathscr P}_{\G, \bS}$}   \la{SECCA1.5}

\medskip \bt \la{Th1.14} Let $\bS$ be a colored decorated surface. Then the moduli space ${\mathscr P}_{\G, \bS}$ has the following features.  
  \begin{enumerate} 
  
 \vskip 1mm
  \item For each \underline{puncture $p$}:
  \vskip 1mm
  
a)   A birational action of the  Weyl group $W$.
 
 b) A  projection onto the  Cartan group $\H$, provided by the monodromy around $p$:
\be \la{mup}
\mu_p:  {\mathscr P}_{\G, \bS}\lra \H.
\ee
It  intertwines the Weyl group action on the space ${\mathscr P}_{\G, \bS}$ from a) with its natural action  on $\H$.
  \vskip 2mm

 \item  For each \underline{special   point $s$, shared by two colored intervals}:
    \vskip 1mm
a)   
A collection of regular functions  on ${\mathscr P}_{\G, \bS}$,  
labeled by   simple positive roots  and called  {\em potentials}:
\be \la{FA}
  {\cal W}_{s,  i}\in {\cal O}({\mathscr P}_{\G, \bS}).
\ee 
 
b)  A  regular map onto the Cartan group $\H$:
\be \la{sspm}
\rho_s:  {\mathscr P}_{\G, \bS}\lra \H.
\ee
\vskip 2mm
\item For each \underline{colored boundary interval ${\rm I}$}, 
 an action  of the Cartan group $\H$ on the space ${\mathscr P}_{\G, \bS}$:
  \be \la{sspm+}
\tau_{\rm I}:  \H\times {\mathscr P}_{\G, \bS}\lra {\mathscr P}_{\G, \bS}.
\ee 

\vskip 2mm

 \item For each \underline{boundary component  $\pi$}, with a specified  reference special point $s$: 
   \vskip 1mm
   
    a)  An action of the  braid group $\B^{(\pi)}_{\mathfrak g}$, see (\ref{BR*}),   by regular transformations of the space ${\mathscr P}_{\G, \bS}$.  
 
 b) If all boundary intervals on $\pi$ are colored, a  projection onto the group $\H_{(\pi)}$ in (\ref{BR**}): 
\be  
\begin{split}
 \mu^\pi_{{\rm out}}: {\mathscr P}_{\G, \bS} \lra  \H_{(\pi)}.\\  
    \end{split}
\ee   
This projection is called the {\rm outer monodromy} around   $\pi$. 

 \vskip 2mm
\item  The actions of the Weyl,  braid,  mapping class        groups, and the  group ${\rm Out}(\G)$  provide an action of the group $\Gamma_{\G, \bS}$.       

\end{enumerate}
    \et
      
      We say that a colored decorated surface $\bS$ is {\it totally colored} if all boundary intervals of $\bS$ are colored.      
\bd 
Assuming that   $\bS$ is totally colored, 
 the monodromy and outer monodromy maps provide the projection $\mu_\bS$ of the space $ {\mathscr P}_{\G, \bS} $ onto the Casimir torus ${\rm C}_{\G, \bS}$ from  Definition \ref{TORID}:
\be  
\label{mubs-tcas}
\mu_\bS:= (\mu_p,  \mu^\pi_{{\rm out}}): {\mathscr P}_{\G, \bS} \lra {\rm C}_{\G, \bS} :=   \H^{\{\mbox{\rm punctures}\}}  \ \times \ \prod_{ \pi}
\H_{(\pi)}. \ee
Here the product is over all punctures $p$ and all boundary components $\pi$ of $\bS$. 
\ed

 \paragraph{\bf Comments.}  
1. Note  similarity between the statements 1a),  4a): they describe the action of the Weyl and the braid groups, assigned to the punctures $p$ 
and  boundary components $\pi$ with special points. 
Note also similarity   between the statements 1b),  4b): the projections in 1b) and 4b)  describe the components of the map $\mu_\bS$ for the punctures $p$ and 
totally colored boundary components $\pi$.

2. The map $\mu_\bS$ is the true analog of the monodromy map for the moduli space $\mathscr{P}_{\G, \bS}$ for any  totally colored decorated surface $\bS$. 
Indeed,  it describes a finite index subalgebra of the Poisson center for the moduli space  $\mathscr{P}_{\G, \bS}$, as well as of the 
center of the quantum algebrs ${\cal O}_q({\mathscr P}_{\G, \bS})$ for generic $q$, 
see Theorem \ref{TH1.16a}.

\begin{proof}  1a) The monodromy of a  generic  $\G$-local system around $p$  is a regular conjugacy class in $\G$.  
So  the  flags near $p$   invariant under the monodromy 
  form a principal homogeneous $W$-set. 
The Weyl group  acts by altering the invariant flag near $p$, keeping the rest of the data intact.   The flat section $\beta_p$ is one of them. So for each $w\in W$ we get a unique flat section $w\cdot \beta_p$.  

\begin{figure}[ht]
 \epsfxsize 200pt
\center{
\begin{tikzpicture}[scale=0.5]
\draw  (0,0) circle (20mm);
\node  at (15:3.2) {{\small $\B_s^r$}};
\node  at (200:3) {{\small $\A^-_s$}};
\node  at (300:3.2) {{\small $\B_s^l$}};
\node at (200:1.4) {{\small $s$}};
\node[red]  at (200:2) {{\small $\bullet$}};
\draw[  -latex  ] ([shift=(25:25mm)]0,0) arc (25:190:25mm);
\draw[  -latex] ([shift=(295:25mm)]0,0) arc (295:210:25mm);
    \end{tikzpicture}
 }
\caption{The triple $(\A^-_s, \B_s^r, \B_s^l)\in {\rm Conf}^*({\cal A}, {\cal B}, {\cal B})$.}
\label{pin12}
\end{figure}

1b) The semi-simple part of the monodromy of a $\G$-local system  lives in $\H/W$. The framing at $p$ allows to lift it to an element of $\H$, providing   the map 
(\ref{mup}). It intertwines the action of the Weyl group $W$ by its very definition. 
 \vskip 1mm
 
2a)  
Recall from Lemma \ref{EQQ} the pair of 
   decorated flags 
   \be \la{pdfl}
   (\A_s^-, \A_s^+)
   \ee near a special point $s$ shared by two colored intervals, sitting over the flag $\B_s$ at $s$. 
   
   Denote  by ${\rm I}_s^r$  the boundary interval containing $s$ and located to the right, i.e. following the boundary orientation,  of $s$. 
   This interval carries two flags at the ends. 
Transporting the rightmost   flag   
towards  the point $s$ along this   interval, we get a flag by $\B_s^r$ near $s$, see  Figure \ref{pin12}. Similarly, using the interval ${\rm I}_s^l$ sharing $s$ and located to the   left of  $s$, 
we get a flag $\B_s^l$ near $s$. So we get a triple 
\be \la{triple}
(\A^-_s, \B_s^r, \B_s^l) \in {\rm Conf}^*({\cal A}, {\cal B}, {\cal B}).
\ee
This construction works even if there is just one special point $s$ on a boundary component, and thus ${\rm I}_s^l= {\rm I}_s^r$.  
Namely, we transport the flag $\B_s$ against the boundary orientation around the boundary component, getting a flag $\B_s^r$ near  $s$, and   transport the same flag $\B_s$ along the boundary orientation, getting a flag $\B_s^l$ near $s$. 

Given a point of $ {\mathscr P}_{\G, \bS}$, we evaluate the   function  ${\cal W}_{\alpha_i}$ in (\ref{POTF})   
  on the   triple (\ref{triple}), getting our function:
 $$
 {\cal W}_{s, \alpha_i}:= {\cal W}_{\alpha_i}(\A^-_s, \B_s^r, \B_s^l).
  $$
  This function is regular. Indeed,  the pairs of flags $(\B_s, \B_s^r)$ and $(\B_s, \B_s^l)$ are generic since the corresponding boundary intervals were assumed to be colored. 
  
  Observe that we can similarly get a function  ${\cal W}_{s, \alpha_i}$ for any special point sharing just a single colored interval, but it is not going to be regular since the second boundary interval sharing $s$ is non-colored. \vskip 1mm

     2b)  
   We define the map $\rho_s$ in (\ref{sspm})   
    by using the pair of decorated flags (\ref{pdfl}), and setting
\be \la{rhoel}
\rho_s({\cal L};  \beta, p):= h(\A^+_s, \A^-_s) \in \H, \ \ \ \ \ \ ({\cal L};  \beta, p) \in {\mathscr P}_{\G, \bS}.
\ee
 \vskip 1mm
 
  3)  It is the action $\tau_{\rm I_s}$  from (\ref{ALa}),
    for the colored boundary interval ${\rm I}_s$ following $s$.  \vskip 1mm

4a) See Section \ref{SEC5}.  \vskip 1mm

4b) Let $s_1, ..., s_d$ be the special points on the  component $\pi$, ordered  following the orientation of the boundary, starting from the specified special point $s_1$. 
 Then, using the map $\rho_s$ in (\ref{sspm}),  we set
\be  \la{OUTM}
\begin{split}
\mu^\pi_{{\rm out}}:=    \left\{  \begin{array}{ll}    \rho_{s_1} \cdot \rho_{s_2}^* \cdot \ldots  \cdot \rho_{s_d}^*   & \mbox{if  $d$   is even,}  \\
      (\rho_{s_1} \cdot \rho_{s_2}^* \cdot \ldots  \cdot \rho_{s_d})_*  & \mbox{if   $d$  is odd.} \\
   \end{array}\right.\\
   \end{split}
\ee  
Here $h^*:= w_0(h^{-1})$, and $(-)_*$ denotes the projection onto the coinvariants $\H_*$, see (\ref{BR**}).    \vskip 1mm

 5) Clear by the very definition of the actions of the all four groups involved.   
\end{proof}
\vskip 1mm

\subsubsection{The Casimir torus and the center.} Our next  result  describes the cluster Casimir torus   ${\rm C}_{\mathscr{P}_{\G, \bS}}$  for  totally colored  $\bS$ 
via the Casimir torus  ${\rm C}_{\G, \bS}$ for  $\mathscr{P}_{\G, \bS}$ from Definition \ref{TORID}. 

  \vskip 2mm

\bt \la{TH1.16a} Let $\bS$ be a totally colored decorated surface. 


\vskip 1mm

1. There is a canonical isogeny from the   cluster Casimir torus  to the Casimir torus:
\be \la{ISOIC}
i_{\rm C}:   {\rm C}_{\mathscr{P}_{\G, \bS}}\lra   {\rm C}_{\G, \bS}.
\ee
It induces an isomorphism of the sets of the real positive points:
\be \la{ISORPP}
i_{\rm C}:   {\rm C}_{\mathscr{P}_{\G, \bS}}(\R_{>0}) \stackrel{\sim}{\lra}   {\rm C}_{\G, \bS}(\R_{>0}).
\ee

\vskip 1mm

2. 
 The center of the    Poisson algebra   ${\cal O}({\mathscr P}_{\G, \bS})$ is a finite dimensional module over the subalgebra $\mu_\bS^* {\cal O}({\rm C}_{\G, \bS})$. 
\vskip 1mm

3.  If $q$ is not a root of unity, then the center of the algebra ${\cal O}_q({\mathscr P}_{\G, \bS})$ is a finite dimensional module over its subalgebra  canonically isomorphic to the algebra $ {\cal O}({\rm C}_{\G, \bS})$:
\be \la{CCasT}
   {\rm Center}\ {\cal O}_q({\mathscr P}_{\G, \bS}) = \mbox{finite dimensional} \ {\cal O}({\rm C}_{\G, \bS})-\mbox{module},   \ \ \ \ \mbox{if $q$ is not a root of unity}.
\ee

  \et

\begin{proof} Here is the architecture of the proof.\vskip 2mm

(1) This is   proved in    Section \ref{SECT20.4}, see Section \ref{20.3.3p}. Note that the claim that the map (\ref{ISORPP}) is an isomorphism implies that the map (\ref{ISOIC}) is an isogeny. The key   step  is  Proposition \ref{14.3.1}, where we prove that   the subalgebra provided by the 
outer monodromy lies in both the Poisson and the quantum center.   
\vskip 2mm

(2),  (3). The Duality Conjectures 
are proved for all moduli spaces $\mathscr{P}_{\G, \bS}$ by Theorem \ref{CBA}. Therefore,  as we prove in Section \ref{SECT20.4}, we have  (\ref{CCasT2}). Now both (2) and (3) 
 are corollaries of (1). 
\end{proof}

Let us remark that the map (\ref{ISOIC}) could indeed have a non-trivial finite kernel. 
For example, it contains a cyclic order $m-1$ subgroup  when $\G = {\rm PGL}_m$ and $S$ is a surface with punctures, see Section \ref{20.3.4}. A complete description of the center see in 
Section \ref{20.3.4}.\vskip 2mm


 {\it The action/projection data.} In addition to the Casimir torus $ {\rm C}_{\G, \bS}$, which   describes the center of the Poisson algebra ${\cal O}({\mathscr P}_{\G, \bS})$, there are   two other important split tori related to the
pair $(\G, \bS)$: 
   \be \la{DEF1.15}
 \begin{split}
 &\widehat  {\rm C}_{\G, \bS}:=    \H^{ \{\mbox{\rm punctures}\}} \times \H^{\{\mbox{\rm special points}\}}.   \\
 &  {\rm F}_{\G, \bS}:=      \H^{ \{\mbox{\rm boundary intervals}\}}.\\
 \end{split}
 \ee 
 
 Using the open string format discussed in Section \ref{SEC3.2.1}, one can say that
    \be \la{DEF1.15}
 \begin{split}
 &\widehat  {\rm C}_{\G, \bS}:=    \H^{ \{\mbox{\rm \underline{black} boundary intervals/components}\}}, \\
 &{\rm F}_{\G, \bS}:=    \H^{ \{\mbox{\rm \underline{white} boundary intervals}\}}.\\
 \end{split}
 \ee

Forgetting   pinnings, we get  a canonical dominant projection 
\be \la{eta}
\eta: {\mathscr P}_{\G, \bS} \lra {\mathscr X}_{\G, \bS}.
\ee 
Indeed, pinnings exist for framed local systems from 
a Zariski open dense subset of ${\mathscr X}_{\G, \bS}$. 
 The torus  ${\rm F}_{\G, \bS}$ acts  simply transitively on the fibers  of $\eta$ via maps  (\ref{ALa}),  
providing  an open embedding
\be \la{PGSX}
\begin{split}
&{\mathscr P}_{\G, \bS}  \left/ {\rm F}_{\G, \bS} \right. \hlra {\mathscr X}_{\G, \bS}.\\
\end{split}
\ee

 The collection of maps $\{\rho_s\}$ in (\ref{sspm}), assigned to the special points of $\bS$,  provides   a  map 
 \be \la{rho*}
 \rho: {\mathscr P}_{\G, \bS} \lra \widehat {\rm C}_{\G, \bS}.
 \ee  
Together with the map $\eta$,   it provides a   double bundle    
 \begin{displaymath}
    \xymatrix{
       & {\mathscr P}_{\G, \bS}  \ar[dl]_{\rho }  \ar[dr]^{\eta} &       \\
    \widehat {\rm C}_{\G, \bS}    &   &  {\mathscr X}_{\G, \bS} }
         \end{displaymath}

 \bt \la{TH1.16}  

Algebras $ \rho^* {\cal O}(\widehat {\rm C}_{\G, \bS})$  and 
        $\eta^*{\cal O}({\mathscr X}_{\G, \bS})$ centralize each other  in the Poisson algebra ${\cal O}({\mathscr P}_{\G, \bS})$. 
     
Quantum  subalgebras $ \rho^* {\cal O}(\widehat {\rm C}_{\G, \bS})$  and 
        $\eta^*{\cal O}_q({\mathscr X}_{\G, \bS})$    centralize each other in ${\cal O}_q({\mathscr P}_{\G, \bS})$ for generic $q$. 
  \et

 To amplify   claim (\ref{slogan}), we  stress the duality between the action/projection data for the pair   $({\mathscr A}_{\G', \bS},  {\mathscr P}_{\G, \bS})$:

 \begin{enumerate}
 
 \item The torus $\widehat {\rm C}_{ \G', \bS}$ action on  ${\mathscr A}_{ \G', \bS}$ $~~ \longleftrightarrow ~~$ The  canonical projection $\rho: {\mathscr P}_{\G, \bS} \lra \widehat {\rm C}_{\G, \bS}$.  
   \vskip 1mm
   
   \item  The torus    ${\rm F}_{\G, \bS}$ action on  ${\mathscr P}_{\G, \bS}$ $~~ ~~\longleftrightarrow ~~$  The  canonical projection  $\tau: {\mathscr A}_{\G', \bS} \lra  {\rm F}_{ \G', \bS}$.  
  \end{enumerate} 
    
\medskip

\subsection{Canonical functions on moduli spaces  ${\mathscr P}_{\G, \bS}$, their quantization, and generalizations}  \la{SECT3.4} 
\medskip

The algebra   of   functions on the Cartan group $\H$ is generated 
by   simple   roots $\alpha_i$ and their inverses $\alpha^{-1}_i$.  

Thanks to  (\ref{FA}) and (\ref{sspm}), given a special point $s$ shared by two colored boundary intervals,   there are the following  regular functions     on   ${\mathscr P}_{\G, \bS}$:
\be \la{136}
{\cal W}_{s, \alpha_i}, \qquad {\cal K}^{\pm 1}_{s, \alpha_i}:=\rho_s^*\alpha^{\pm 1}_i.
\ee

 Recall the action of the Weyl group $W^n$ on $\mathscr{P}_{\G, \bS}$, where $n$ be the number of punctures on $\bS$. 
 
  \bt \la{QWE} Given a special point $s$ on a colored decorated surface $\bS$, shared by two colored boundary intervals, the functions ${\cal W}_{s,  i}$ and ${\cal K}^{\pm 1}_{s, i}$  give rise to   elements of the quantum algebra
\be \la{WH}
 {\bf W}_{s,  i}, ~{\bf K}^{\pm 1}_{s, i}~\in {\cal O}_q({\mathscr P}_{\G, \bS})^{W^n}.
\ee
\et
This means that for \underline{each} cluster Poisson coordinate system ${\bf c}$ on the space ${\mathscr P}_{\G, \bS}$ we have the following: 

\vskip 2mm

1.   There are elements  ${\rm W}^{{\bf c}}_{s,  i}, ~({{\rm K}^{\bf c}_{s, i}})^{\pm 1}$  
 of the   quantum torus algebra ${\cal O}_q({\rm T}_{{\bf c}})$ assigned to  ${\bf c}$, whose 
  $q \to 1$ limits   exist  and coincide with 
the   functions ${\cal W}_{s,  i}, ~{\cal K}_{s, i}^{\pm 1}$, expressed in the coordinates ${\bf c}$. 

\vskip 1mm 

2. Since both ${\bf K}_{s, i}$ and its inverse ${\bf K}_{s, i}^{-1}$ are  Laurent polynomials in every seed, they must be  Laurent {\it monomials} in every seed.

\vskip 1mm

3. For  {any} cluster Poisson transformation ${\bf c}_1 \to {\bf c}_2$, the corresponding   quantum cluster transformation 
  ${\rm Frac}~{\cal O}_q({\rm T}_{{\bf c}_2}) \lra {\rm Frac}~{\cal O}_q({\rm T}_{{\bf c}_1}) $ maps   elements 
${\rm W}^{{\bf c}_2}_{s,  i}, ~({\rm K}_{s, i}^{{\bf c}_2})^{\pm 1}$ to the   ones ${\rm W}^{{\bf c}_1}_{s,  i}, ~({\rm K}_{s, i}^{{\bf c}_1})^{\pm 1}$.

\vskip 2mm

 Theorem \ref{QWE}  is a central and  difficult result. Its proof uses the   Quantum Lift Theorem \ref{quantum.promotion.f}.
\vskip 1mm

Next, let us  glue colored boundary interval ${\rm I}_1$ and ${\rm I}_2$ on a colored decorated surface $\bS$, so that the special points $a \in {\rm I}_1$ and $b\in {\rm I}_2$  
are glued into a special point $s$, getting a new decorated surface $\bS'$. The gluing of  surfaces induces the map of the corresponding 
moduli spaces and their quantized algebras of functions: 
\be \la{GLM}
\begin{split}
&\gamma: \mathscr{P}_{\G, \bS} \lra \mathscr{P}_{\G, \bS'},\\
&\gamma^*: {\cal O}_q(\mathscr{P}_{\G, \bS'})^{W^n} \lra {\cal O}_q(\mathscr{P}_{\G, \bS})^{W^n}. \\
\end{split}
\ee 
So it is natural to ask the following two basic questions:
\be \la{QPSA}
\begin{split}
&\mbox{ \it What is the subalgebra of the algebra ${\cal O}_q({\mathscr P}_{\G, \bS})$  generated by   elements (\ref{WH})?}\\
&\mbox{ \it How the  elements (\ref{WH})  behave under the gluing of decorated surfaces?}\\
\end{split}
\ee
Theorem \ref{UEAB} below addresses   questions (\ref{QPSA}), revealing their connection to quantum groups.

\subsubsection{Relation to  quantum groups.} Recall  the        quantum universal enveloping algebra ${\U}_q(\mathfrak{g})$   of the Lie algebra $\mathfrak{g}$. 
It is defined as a Hopf algebra with the standard generators $\{E_i, F_i, K^\pm_i\}$, used e.g. in  the book \cite{L}. See Section \ref{SSECC11.1} for more details. 

We use  the generators    $\{\bE_i, \bF_i, \bK_i^{\pm}~|~ i=1,..., r\}$ of ${\U}_q({\mathfrak{g}})$,   
obtained by  rescaling the standard generators:
\be \la{RESCX}
\bE_i= q^{-\frac{1}{2}}(q-q^{-1}) E_i, \hskip 7mm \bF_i= q^{\frac{1}{2}}(q^{-1}-q) F_i, \hskip 7mm \bK_i =K_i. 
\ee
They satisfy  standard relations, see Section \ref{SSECC11.1}. 
There is an anti-involution  of  the algebra $\U_q(\mathfrak{g})$:
\be \la{ASST}
\ast: \U_q(\g) \lra \U_q(\g)^{\rm op}, \hskip 9mm \bE_i \lms \bE_i, \hskip 6mm \bF_i \lms \bF_i, \hskip 6mm \bK_i\lms \bK_i, \hskip 6mm q\lms q^{-1}.
\ee
 Rescaling (\ref{RESCX}) of the generators  of $\U_q(\mathfrak{g})$ is forced on us by the requirement to have an $\ast$-algebra structure 
on $\U_q(\mathfrak{g})$ with the involution $\ast$ such that   the generators are selfadjoint.  \\

A framing at a special point $s$ of   $\bS$ provides   a reduction of a $\G$-local system near $s$ 
to the Borel subgroup $\B$.  The Hopf subalgebra    ${\U}_q(\mathfrak{b})\subset {\U}_q(\mathfrak{g})$ assigned to  $\B$ is generated by   elements ${\bf E}_i, {\bf K}_i$. \vskip 1mm

Recall the action of the group $W^n$  on the algebra ${\cal O}_q({\mathscr P}_{\G, \bS})$ provided by Theorem \ref{Th8.5A}. \vskip 1mm

 Note that 
by Part 4) of Theorem \ref{TH1.16}, $\mu_{(\pi)}^*{\cal O}_{\H_{(\pi)}}$ lies in the center of the algebra ${\cal O}_q({\mathscr P}_{\G, \bS})$. 
So we can consider the quotient   ${\cal O}_q({\mathscr P}_{\G, \bS})_{\mu_{(\pi)}=1}$  of the algebra ${\cal O}_q({\mathscr P}_{\G, \bS})$ by the ideal given by the condition that the outer 
monodromy around the boundary component $\pi$ is equal to $1$. Note that this quotient does not depend on the choice of a reference special point on $\pi$.

\bt \la{UEAB} Let $\G$ be an adjoint split semi-simple algebraic group over $\Q$. Then for any colored decorated surface $\bS$:  

\begin{enumerate}   
 
\vskip 1mm \item For any special point $s$ on $\bS$, shared by    colored boundary intervals,   elements (\ref{WH}) provide   an injective map of algebras\footnote{In this paper we prove this assuming that the component $\pi$ containing $s$ has $>1$ special points.} 
\be \la{383aa}
\begin{split}
&\kappa_s: {\U}_q(\mathfrak{b}) \hlra {\cal O}_q({\mathscr P}_{\G, \bS})^{W^n}.\\
&{\bf E}_i \lms {\bf W}_{s,  i},\qquad~ \ \ \ {\bf K}_i \lms{\bf K}_{s, i}.\\
\end{split}
\ee

\vskip 1mm 
  
\item
 Then gluing map $\gamma^*$ in (\ref{GLM}) acts on the canonical elements (\ref{WH})  as follows:   
 \be \la{383}  \begin{split}
 &\gamma^{\ast} ({\bf K}_{s,i}) = {\bf K}_{a, i} \otimes  {\bf K}_{b, i},\\
  &\gamma^{\ast} ({\bf W}_{s, i}) = {\bf W}_{a, i} \otimes 1 + {\bf K}_{a, i} \otimes {\bf W}_{b, i}.\\
    \end{split}
  \ee 
So  maps $\kappa_*$ intertwine the coproduct $\Delta$ with the gluing map $\gamma^*$, providing a commutative diagram
  \begin{displaymath}
    \xymatrix{
       {\U}_q(\mathfrak{b}) \ar[r]^{\Delta \qquad} \ar[d]_{\kappa_s} &{\U}_q(\mathfrak{b}) \otimes {\U}_q(\mathfrak{b})    \ar[d]_{\kappa_a}^{\kappa_b} \\
         {\cal O}_q({\mathscr P}_{\G, \bS'})   \ar[r]^{\gamma^*}&  {\cal O}_q({\mathscr P}_{\G, \bS})  }
         \end{displaymath}
         
          \item For any boundary component $\pi$ of $\bS$ with exactly   two  special points  $e$ and $f$,  and assuming that both boundary intervals on $\pi$ are colored, 
           there is a canonical injective map   of $\ast$-algebras:  
\be \la{MKa}
\begin{split}
\kappa_{e,f}: ~&\U_q({\g}) \hra {\cal O}_q({\mathscr P}_{\G, \bS})_{\mu_{\rm out}=1}^{W^n};\\
&{\bf E_i} \lms {\bf W}_{e,  i}, \ \ \   {\bf F_{i}} \lms {\bf W}_{f,  i^*}, \ \ \ {\bf K_i} \lms  {\bf K}_{e, i}.\\
\end{split}
\ee

 \end{enumerate}
\et

 Theorems \ref{QWE} and  \ref{UEAB} are proved in Section \ref{SSEECC11.3}. 
 
\vskip 2mm

There are two ways to generalize Theorem \ref{UEAB}. 

First, the map $\kappa_s$ in (\ref{383aa}) can be defined starting from any root system $\Delta$, 
mapping the  quantum Borel algebra ${\U}_q(\mathfrak{b}_\Delta)$ for the root system $\Delta$ 
to the quantum algebra ${\cal O}_q(\mathscr{X}_{[b]})$ related to the cluster Poisson variety $\mathscr{X}_{[b]}$ assigned to any element $b$ 
of the positive braid semigroup ${\Bbb B}_\Delta^+$ of $\Delta$. 
We explain this  in Section \ref{SECT4.2.1}.
 
Second, assuming $\Delta$ is of finite type, we  describe the relations between the generators of the algebras 
$\kappa_s({\U}_q(\mathfrak{b}))$ for different special points $s$ on $\bS$. This  generalizes   part (3) of Theorem \ref{UEAB}, see  
 Section \ref{SECT4.2.2}. 
 \medskip
 
 \subsubsection{A generalization  to an arbitrary root system} \la{SECT4.2.1} Let $\Delta$ be the root system associated with any skewsymmetrizable 
 Cartan matrix ${\rm C}_{ij}$, where  $i,j \in {\rm I}$.   Then there is the associated braid group ${\Bbb B}_\Delta$, and its positive semigroup ${\Bbb B}_\Delta^+$. We  define,
  just as above, the quantum universal enveloping algebra ${\U}_q(\mathfrak{b}_\Delta)$ for the Borel subalgebra $\mathfrak{b}_\Delta$ of the Kac-Moody Lie algebra for  $\Delta$. For any   $b \in {\Bbb B}_\Delta^+$ there is  
 the  cluster Poisson variety $\mathscr{X}_{[b]}$, defined  in \cite{FG05} for the finite type $\Delta$, and in \cite{SW} in general, see Section \ref{SECT16.2}. 
 
 Next, let $b_1, b_2\in {\Bbb B}^+_\Delta$. Then we have the gluing map
 \be \la{GLMa}
\begin{split}
&\gamma_{b_1, b_2}: \mathscr{X}_{[b_1]} \times \mathscr{X}_{[b_2]}\lra \mathscr{X}_{[b_1b_2]},\\
&\gamma_{b_1, b_2}^*: {\cal O}_q(\mathscr{X}_{[b_1b_2]}) \lra {\cal O}_q(\mathscr{X}_{[b_1]}) \otimes  {\cal O}_q(\mathscr{X}_{[b_2]}).\\
\end{split}
\ee

 Generalizing Theorem \ref{QWE} and  parts (1)-(2) of Theorem \ref{UEAB}, we prove the following general result.
 
 \bt \la{QWEaa} Let $\Delta$ be  any   root system,  and  $b \in {\Bbb B}^+_\Delta$ any positive braid semigroup element. Then
 
 1) For each $i \in {\rm I}$, there are canonical elements 
 \be \la{WHa}
 {\bf W}_{b,  i}, ~{\bf K}^{\pm 1}_{b, i}~\in {\cal O}_q({\mathscr X}_{[b]}).
\ee 
 
 2)  There is   a canonical map of algebras  
\be
\begin{split}
&\kappa_b: {\U}_q(\mathfrak{b}) \hlra {\cal O}_q({\mathscr X}_{[b]}).\\
&{\bf E}_i \lms {\bf W}_{b,  i},\qquad~ \ \ \ {\bf K}_i \lms{\bf K}_{b, i}.\\
\end{split}
\ee

3)  For any $b_1, b_2\in {\Bbb B}_\Delta^+$, and any  $i \in {\rm I}$, the gluing map $\gamma_{b_1, b_2}^*$ in (\ref{GLMa}) acts  as follows:   
 \be \la{383}  \begin{split}
 &\gamma_{b_1, b_2}^{\ast} ({\bf K}_{b_1b_2,i}) = {\bf K}_{b_1, i} \otimes  {\bf K}_{b_2, i},\\
  &\gamma_{b_1, b_2}^{\ast} ({\bf W}_{b_1b_2, i}) = {\bf W}_{b_1, i} \otimes 1 + {\bf K}_{b_1, i} \otimes {\bf W}_{b_2, i}.\\
    \end{split}
  \ee 
So  maps $\kappa_*$ intertwine the coproduct $\Delta$ in ${\U}_q(\mathfrak{b})$ with the gluing map $\gamma_{b_1, b_2}^*$, providing a commutative diagram
  \begin{displaymath}
    \xymatrix{
       {\U}_q(\mathfrak{b}_\Delta) \ar[r]^{\Delta \qquad} \ar[d]_{\kappa_{b_1, b_2}} &{\U}_q(\mathfrak{b}_\Delta) \otimes {\U}_q(\mathfrak{b}_\Delta)    \ar[d]_{\kappa_{b_1}}^{\kappa_{b_2}} \\
         {\cal O}_q({\mathscr X}_{[b_1b_2]})   \ar[r]^{\gamma_{b_1, b_2}^*\ \ \ \ }&  {\cal O}_q({\mathscr X}_{[b_1]})  \otimes  {\cal O}_q({\mathscr X}_{[b_2]}) }
         \end{displaymath} \et
         
         Theorem \ref{QWEaa} is proved in Section \ref{SECT16.2}.
 
\subsubsection{Relations between the elements ${\bf W}_{s,i}, {\bf K}_{s,j}$ at different special points}\la{SECT4.2.2} Let $\Delta$ be a root system of finite type. Recall the involution $i \to i^*:=w_0(-\alpha_i)$ on the set ${\rm I}$ of the vertices of the Dynkin diagram of $\Delta$. 

\bt \la{5.25.24q} Let $\G$ be  as in Theorem \ref{UEAB}. Given any two special points $(s,t)$ on a decorated surface $
\bS$, connected by a single colored boundary interval, oriented  $s \to t$ by the boundary orientation, we have 
\be \la{5.25.24}
\begin{split}
&[{\bf W}_{s,i}, {\bf W}_{t, j^*}]=  
    \delta_{ij} (q^{-1}_i-q_i){\bf K}_{s,i}. \\
&{\bf K}_{s,i}{\bf W}_{t,j^*} = q^{-{\rm C}_{ij}}{\bf W}_{t,j^*}{\bf K}_{s,i}.\\
&{\bf K}_{t,i}{\bf W}_{s,j} = {\bf W}_{s,j}{\bf K}_{t,i}.\\
&{\bf K}_{s,i}{\bf K}_{t,j^*} = q^{-{\rm C}_{ij}}{\bf K}_{t,j^*}{\bf K}_{s, i}.\\
\end{split}
\ee
For  special points $s,t$ on different boundary components, the subalgebras $\kappa_s({\U}_q(\mathfrak{b}))$ and $\kappa_t({\U}_q(\mathfrak{b}))$ commute. 
\et

We prove Theorem \ref{5.25.24q}  in Section \ref{24.5.26.10.19}.
Note that in the part 3) of Theorem \ref{UEAB} there are just two special boundary points on the component  $\pi$, thus connected by two boundary intervals, 
while in    (\ref{5.25.24}) we have two special points connected by a single boundary interval. Considered together, they treat all possible cases for a pair of special 
boundary points. 

And they are  compatible. In the former case 
 the relation can be written as 
\be \la{5.25.24a}
\begin{split}
&[{\bf W}_{s,i}, {\bf W}_{t, j^*}]=  \delta_{ij}
    (q^{-1}_i-q_i)({\bf K}_{s,i}- {\bf K}_{t, j^*}).\\
\end{split}
\ee
It matches the quantum group relation
$[{\bf E}_i, {\bf F}_j]=\delta_{ij}(q_i^{-1}-q_i)({\bf K}_i- {\bf K}_i^{-1})$, see (\ref{Re1}),  since ${\bf K}_{t, j^*} = {\bf K}_{s, j}^{-1}$ 
after imposing the condition that  the outer monodromy is equal to $1$. 
\medskip

\subsubsection{Example.} Here is how all these relations look in the simplest case when $\bS$ is a triangle $t$, and $\G={\rm PGL}_2$. 
\be
{\bf W}_s = X_{e_3}, \ \ \ {\bf K}_s = X_{e_1+e_3}, \ \ \ {\bf W}_t = X_{e_1}, \ \ \ {\bf K}_t= X_{e_1+e_2}.
\ee
\begin{figure}
\begin{tikzpicture}[scale=0.8]
\draw (90:2)--(210:2)--(330:2)--(90:2);
\node[red] (a) at (30:1) {$\bullet$};
\node[red] (b) at (150:1) {$\bullet$};
\node[red] (c) at (270:1) {$\bullet$};
\node at (90:2.3) {$t$};
\node at (210:2.3) {$s$};
\node  at (30:1.3) {$e_2$};
\node  at (150:1.3) {$e_1$};
\node  at (270:1.3) {$e_3$};
\draw[-latex, blue, thick] (a)--(b);
\draw[-latex, blue, thick] (b)--(c);
\draw[-latex, blue, thick] (c)--(a);
\end{tikzpicture}
\caption{The quiver for $\mathscr{P}_{{\rm PGL}_2, t}$. Note that the boundary orientation  induced by the triangle orientation is shown as the clockwise orientation on the pictures. }
\label{24.5.26.9.31}
\end{figure}
As illustrated in Figure \ref{24.5.26.9.31}, the marked point $s$ precedes $t$ in the clockwise order. Then 
\[
[{\bf W}_s, {\bf W}_t] = X_{e_3}X_{e_1}- X_{e_1}X_{e_3} = (q^{-1}-q)X_{e_1+e_3} = (q^{-1}-q){\bf K}_s,
\]
\[
{\bf K}_s{\bf W}_t = X_{e_1+e_3}X_{e_1}=q^{-2}{\bf W}_t{\bf K}_s,
\]
\[
{\bf K}_t{\bf W}_s=X_{e_1+e_2}X_{e_3}= {\bf W}_s {\bf K}_t,
\]
\[
{\bf K}_s{\bf K}_t=X_{e_1+e_3}X_{e_1+e_2}=q^{-2}{\bf K}_t{\bf K}_s.
\]


\bcon \la{LS22} Let ${\rm D}_n^*$ be a punctured disc with $n$ special boundary points. Then for any $n>0$ the algebra 
\be \la{Dn}
{\cal O}_q(\mathscr{P}_{\G, {\rm D}_n^*})^W
\ee is isomorphic to the algebra 
generated by the elements ${\bf W}_{s, i}, {\bf K}_{s, i}$, assigned to the special boundary points $s$, and satisfying the relations (\ref{5.25.24}), as well as the Serre relations for each $s$. 
\econ

For $n=2$  Conjecture \ref{LS22} is proved for simply laced $\G$ in \cite{S22}. 

Consider the subalgebra of algebra (\ref{Dn}) given by  the condition that the outer monodromy $\mu_{\rm out}$ is equal to $1\in \H$:
\be \la{APN}
{\cal O}_q({\rm Loc}_{\G, {\rm D}_n^*}):= {\cal O}_q(\mathscr{P}_{\G, {\rm D}_n^*})_{\mu_{\rm out}=1}^W.
\ee
Then for each internal/external diagonal $\rm E$ cutting the convex $n-$gon ${\rm P}_n$ into  two polygons ${\rm P}_{a+1}$ and ${\rm P}_{b+1}$, where $a+b=n$, $a,b\geq 1$,
 we have the  map of algebras
\be
\Delta_{\rm E}: {\cal O}_q({\rm Loc}_{\G, {\rm D}_n^*}) \lra {\cal O}_q({\rm Loc}_{\G, {\rm D}_{a+1}^*})\otimes {\cal O}_q({\rm Loc}_{\G, {\rm D}_{b+1}^*}).
\ee
It is induced by the gluing  of two decorated surfaces ${\rm D}_{a+1}^*$ and ${\rm D}_{b+1}^*$ along a pair of sides into the decorated surface ${\rm D}_{a+b}^{**}$ with two punctures, 
and then encircling the punctures, and cutting out the obtained twice punctured disc. The part 3) of Theorem \ref{QWEaa} allows to calculate
 the effect of these maps on the generators, generalizing the formulas for the coproduct in the quantum group Hopf algebra. Given any pair of nonintersecting diagonals ${\rm E}, {\rm F}$ of 
 ${\rm P}_n$ we have $\Delta_{\rm F}\circ \Delta_{\rm E} = \Delta_{\rm F}\circ \Delta_{\rm E} $ - the analog of the coassociativity condition. 
 For example, specifying a side of the triangle, we provide   ${\cal O}_q({\rm Loc}_{\G, {\rm D}_3^*}) $ with a structure of the comodule over the Hopf algebra ${\cal O}_q({\rm Loc}_{\G, {\rm D}_2^*}) $. 
 
 It would be  interesting to study the representation theory of these algebras, e.g. finite dimensional ones. 
 For example,  how to build all representations using  representations of  algebras (\ref{APN})
 for $n=2$ and $n=3$ - the former are   ${\rm U}_q(\mathfrak{g})-$modules?

 \subsection{Landau-Ginzburg potentials  in mirror symmetry, and quantum group generators} \la{SECTT4.3}

\medskip For either the space ${\mathscr A}_{\G, \bS}$, or the space ${\mathscr P}_{\G, \bS}$ with all sides colored and due to Lemma \ref{EQQ}, or each special point $s$  carries a decorated   flag $\A_s$. 
We can  assign to $s$ a  regular function ${\cal W}_s$. It  is  decomposed as a sum: ${\cal W}_s = \sum_{i \in {\rm I}}{\cal W}_{s, i}$. \vskip 1mm
 
 Amazingly, the canonical functions ${\cal W}_{s, i}$ play a key role in two seemingly unrelated stories: 
 
\begin{enumerate}

\item  The quantization ${\rm W}_{s, i}$ of the functions ${\cal W}_{s, i}$  give rise to the quantum group  generators ${\bf E}_i, {\bf F}_i$. 
 
\vskip 1mm\item  Functions ${\cal W}_{s, i}$ are the Landau-Ginzburg potentials in the mirror symmetry conjectures \cite{GS13}. 

\end{enumerate}

 Let us elaborate on the mirror symmetry. 
 The sum of the functions ${\rm W}_{s,}$ over all special points $s$ on either moduli space ${\mathscr P}_{\G, \bS}$ or ${\mathscr A}_{\G, \bS}$  is denoted by ${\cal W}_{\rm S}$.  
For the space ${\mathscr A}_{\G, \bS}$,  each puncture $p$ also carries a decorated flag $\A_p$ and hence a potential ${\cal W}_p$. 
 The sum of the   potentials ${\cal W}_p$ over all punctures $p$ 
 is denoted by ${\cal W}_{\rm P}$.    \vskip 1mm
 
Let  $\ast \bS$ be the new decorated surface which coincides with $\bS$ as a surface, and  whose special points are the midpoints of  boundary segments of $\bS$. 
 Then, elaborating \cite[\S 10.2]{GS13}, we state, see   Figure \ref{amalgam5}:
 
 \bcon \la{MIRC} Let $\G$ be a split adjoint semi-simple   group over $\Q$, and $\bS$ a  decorated surface. Then

\begin{itemize} \vskip 1mm\item The moduli space $({\mathscr P}_{\G, \bS}, {\cal W}_{\rm S})$ is  mirror dual to the 
 moduli space  ${\mathscr A}_{\G^\vee, \ast \bS}$.

\vskip 1mm\item The moduli space $({\rm Loc}_{\G, \bS}, {\cal W}_{\rm S})$ is  mirror dual to the 
 Landau-Ginzburg model   $({\mathscr A}_{\G^\vee, \ast \bS},  {\cal W}_{\rm P})$.

\end{itemize} 
\econ

     \begin{figure}[ht]
     \begin{center}
\begin{tikzpicture}
 \draw[red, thick, dashed] (0,0) circle (1);
  \draw[blue] (0,0) circle (.1);
  \draw[red]  (4,0) circle (.1);
  \draw[thick]  (4,0) circle (1);
         \node[red] at (0, 1) {\footnotesize $\bullet$};
        \node[red] at (0, -1) {\footnotesize $\bullet$};
         \node[red] at (3, 0) {\footnotesize $\bullet$};
         \node[red] at (5, 0) {\footnotesize $\bullet$};
          \node[red] at (0, 1.4) {\footnotesize $\A_1$};
          \node[red] at (0, -1.4) {\footnotesize $\A_2$};
          \node[red] at (2.7, 0) {\footnotesize $\A_3$};
         \node[red] at (5.3, 0) {\footnotesize $\A_4$};
          \node[blue] at (0.4, 0) {\footnotesize $\B_p$};
           \node[red] at (4.4, 0) {\footnotesize $\A_p$};
\end{tikzpicture}
\end{center}
\caption{Mirror duality for the decorated surface $\odot$. On the left: the space ${\rm Loc}_{\G, \odot}$ with the potential 
${\cal W}_{\rm S}  = {\cal W}_{1} + {\cal W}_{2} $. Decorated
 flags $\A_1, \A_2$ are generic on each  arc. On the right: the space ${\mathscr A}_{\G'^\vee, \odot}$ with the potential 
${\cal W}_{p}$. No restrictions on  $\A_3, \A_4$. }
\label{amalgam5}
\end{figure}

  Conjecture \ref{MIRC} can be enhanced by incorporating the action/projection data from Section \ref{SECCA1.5}.3. Namely, when   one of the tori acts on any of the spaces, we consider the related  equivariant 
  ${\rm D}^b{\rm Coh}$/Fukaya category, while  for the projections $\mu_\bS$ and $\tau_\bS$ we consider the corresponding family of categories.

 
\medskip

\section{Poisson Lie groups as moduli spaces of $\G-$local systems on   surfaces} \la{PLG}

\medskip

\subsection{An overview} 

\medskip

In Section \ref{PLG} 
 we show that   classical Poisson Lie groups related to semi-simple Lie algebras $\mathfrak{g}$  can be   described as  moduli spaces related to certain 
colored decorated surfaces $\bS$ and the adjoint group $\G$. \vskip 1mm

Here are the universal features of this approach. \vskip 1mm

The group law is   given by the gluing map,   followed in some cases by the restriction map. 
Pinnings are crucial, since allow to glue the moduli spaces: the group law can not be defined 
without    pinnings. 

Since   the gluing and cutting maps,   as well as the action of the Weyl group $W$ at the punctures   are   Poisson maps,  the group law respects the  Poisson structure. 

The unit is  given by  the trivial $\G-$local system on $\bS$ with certain specific  pinnings at the colored   intervals. 

The inversion map is   induced by a   flip of the   surface, which alters the   orientation, and thus  the sign of the Poisson bivector. \vskip 1mm

The pinning provides a trivialization of the $\G-$local systems near the boundary interval where it sits. The identification of the moduli space with the related group is given by the action of the parallel transport  on pinnings. Here we see once again the crucial role pinnings play in the whole story. 
\vskip 1mm

The origin of the Poisson bracket is always the same - cluster. The cluster Poisson structure on the  moduli spaces allows to think about the Poisson structure on the  groups geometrically. 
It does not require  $r-$matrices. 
Different Poisson structures on the same group appear due to different geometric realizations of the   group. 
An example is given by the two  Poisson structures on the group $\G$. \vskip 1mm

Most importantly, the cluster quantization machine  allows us to quantize these Poisson Lie groups 
 in the strongest possible way, where $\hbar$ is  any complex number rather than a formal parameter, as was in the original Drinfeld-Jimbo approach. 
 The coproduct and the Hopf algebra structure results from the cluster Poisson nature of the gluing maps described by the cluster Poisson amalgamation,   followed by cluster projection. \vskip 1mm
 
 The universal quantum $R-$matrix $R$, which  lies in a completion  of the tensor
  square ${\cal H} \otimes {\cal H}$ of the Hopf algebra ${\cal H}$, provides an automorphism     of  the algebra ${\cal H} \otimes {\cal H}$:
\be \la{INDR}
 {\rm Ad}_R:  {\cal H} \otimes {\cal H} \lra {\cal H} \otimes {\cal H}, \qquad X \lms R X R^{-1}.
\ee
  The key point is that the  automorphism (\ref{INDR})   is induced by an element ${\cal R}$ of the cluster modular group of the underlying 
  cluster Poisson variety, which we call ${\cal R}$ the {\it geometric ${\cal R}-$matrix}. Moreover,  ${\cal R}$   is an automorphism of the underlying stack.  \vskip 1mm
  
  In Section \ref{TQFT} we show that  considerations of Section \ref{PLG} are specific examples  of the general TQFT-like structure on the moduli spaces.

\subsubsection{The general set-up.} 
Given a Poisson Lie group $\rm P$, there are three related Poisson Lie groups: the dual Poisson Lie group $\rm P^*$, the Drinfeld double $\D(\rm P)$ of   $\rm P$, 
 and  the Poisson Lie group   $\D(\rm P)^*$ dual to   $\D(\rm P)$. 

Let  $\mathfrak{p}:={\rm Lie} ~\rm P$, $\mathfrak{p}^*:={\rm Lie} ~{\rm P}^*$ and ${\rm D}_{\mathfrak{p}}:= {\rm Lie}~\D(\rm P)$. Then 
$({\rm D}_{\mathfrak{p}}; \mathfrak{p}, \mathfrak{p}^*)$ is a Manin triple. 

This means that there is a canonical  vector space isomorphism   ${\rm D}_{\mathfrak{p}} = \mathfrak{p}\oplus \mathfrak{p}^*$. It  provides the  non-degenerate symmetric bilinear form $\langle (e, f), (e',f')\rangle:= f(e')+f'(e)$ such that 
 $ \mathfrak{p}$ and $ \mathfrak{p}^*$ are isotropic and $\langle \cdot,  \cdot\rangle: \mathfrak{p}\otimes \mathfrak{p}^*\to \R$ is the natural pairing. There is a unique  Lie algebra structure 
 $ [\cdot, \cdot]_{{\rm D}_{\mathfrak{p}}}$ on ${\rm D}_{\mathfrak{p}}$ such that the   embeddings $\mathfrak{p} \hra  {\rm D}_{\mathfrak{p}}$ and $\mathfrak{p}^* \hra  {\rm D}_{\mathfrak{p}}$ 
 are Lie algebra maps, and the form $\langle \cdot,  \cdot\rangle$ is ${\rm D}_{\mathfrak{p}}-$invariant:
 \be
 [(e,f), (e', f')]_{{\rm D}_{\mathfrak{p}}}:= \Bigl([e,e'] _{\mathfrak{p}}+ {\rm ad}_{\mathfrak{p}^*}f(e') - {\rm ad}_{\mathfrak{p}^*}f'(e), ~~[f,f']_{\mathfrak{p}^*} + {\rm ad}_{\mathfrak{p}}e(f') - {\rm ad}_{\mathfrak{p}}e'(f)\Bigr).
 \ee

 \vskip 1mm

There are canonical maps of the Poisson Lie groups:
\be \la{PPDD}
\begin{split}
&\rm P \hra \D(\rm P), \quad \quad \rm P^* \hra \D(\rm P)\\
& \D(\rm P)^* \lra \rm P, \quad  \D(\rm P)^* \lra \rm P^*.\\
\end{split}
\ee
They  induce the following   maps, where the first map is   Poisson, and the second is a Lie group map: 
\be
\begin{split}
& \rm P \times \rm P^* \lra \D(\rm P), \\
&  \D(\rm P)^* \lra \rm P \times  \rm P^*.\\
\end{split}
\ee 
The  first map  is the composition $\rm P \times \rm P^* \lra \D(\rm P) \times \D(\rm P) \lra  \D(\rm P)$, where the first map is induced by the first line in (\ref{PPDD}), and the second   is the product. 
So it is a Poisson map,   a local isomorphism, but not a group map. 
The second map is a group map, but it is not   Poisson.\vskip 1mm

  There are   related quantum groups, that is Hopf algebras 
$$
\U_q(\mathfrak{p}),  \ \U_q(\mathfrak{p}^*), \ \U_q({\rm D}_{\mathfrak{p}}), \ \U_q({\rm D}^*_{\mathfrak{p}}).
$$
The Hopf algebra $\U_q({\rm D}_{\mathfrak{p}})$ is the quantum Drinfeld double of 
$\U_q(\mathfrak{p})$. \vskip 1mm

Below we consider {\it cluster Poisson Lie groups} $\rm P$ and $\rm P^*$, that is the Poisson Lie groups which  carry  cluster Poisson   structures, compatible with the group structure. Therefore they provide us algebras  ${\cal O}_q(\rm P^*)$ and ${\cal O}_q({\rm D}({\rm P})^*)$, which have the Hopf algebra structure with the coproduct induced by the group law. 

One should have isomorphisms of Hopf algebras, which needs to be established: 
\be \la{149}
\U_q(\mathfrak{p}) \stackrel{?}{=} {\cal O}_q(\rm P^*), \qquad {\U}_q({\rm D}_{\mathfrak{p}}) \stackrel{?}{=}  {\cal O}_q({\rm D}({\rm P})^*). 
\ee

Precisely, we consider the cluster Poisson Lie groups $\rm P= \B$ given by the Borel subgroup, and $\rm P =\G$, realizing them as moduli spaces of $\G-$local systems related to decorated surfaces. \vskip 1mm

  
 In Sections  \ref{SEC1.7} -  \ref{SECT4.2d}  we realize     via moduli spaces of $\G-$local systems 
  on  
 colored decorated surfaces $\bS$   various Poisson Lie groups and related spaces, including the following ones: 
 
 \begin{enumerate}

\item  The  Poisson Lie group $\G$.
 
\vskip 1mm \item The Poisson Lie Borel subgroup $\B \subset \G$. 
 
 \vskip 1mm\item  The  Poisson Lie group $\G^*$ dual to the one $\G$.

\vskip 1mm\item  The   Drinfeld double $\D(\G)$ of the  
 Poisson Lie group $\G$.
 
\vskip 1mm \item  The   dual Poisson Lie group $\D(\G)^*$ of the Drinfeld double $\D(\G)$.
 
\vskip 1mm\item The Grothendieck resolution  $\widehat \G\lra \G$ and its Poisson structure.\footnote{It provides the second Poisson structure on $\G$, non-compatible with the group law.}
 
\vskip 1mm \item  The  Heisenberg double ${\rm D}_{\rm H}(\G)$  of $\G$. 

\end{enumerate}
 
 All features of these Poisson Lie groups, including   canonical maps (\ref{PPDD})  and their quantum counterparts, are evident from the general features of the moduli spaces and their cluster Poisson structure. 
 However relating them to the classical definitions requires more work.

\medskip

\subsection{The moduli space ${\mathscr L}_{\G, \odot}$ $=$ the Poisson Lie group $\G^*$} \la{SEC1.7}

 \medskip
\subsubsection{The group $\G^*$.} \la{4.1}

Let us denote by  $\odot$ a colored decorated surface given by a  disc  with a puncture $p$,  two ordered special  points $s, t$, and two colored boundary intervals.  

The moduli space ${\mathscr P}_{\G, \odot}$ parametrizes 
the data illustrated on Figure \ref{gs1}. 
 \begin{figure}[ht]
 \epsfxsize 200pt
\center{
\begin{tikzpicture}[scale=0.5]
\begin{scope}[rotate=270]
\draw  (0,0) circle (20mm);
\node [blue, thick]  at (0,0) {\Large $\circ$};
\node [red]  at (2,0) {\small $\bullet$};
\node [red]  at (-2,0) {\small $\bullet$};
\node [label=right: {\small $\A_s^r$}, label=left:  {\small $\A_s^l$}, label=above:{\small $h_s$}] at (-2.9,0) {};
\node [label=right: {\small $\A_t^l$}, label=left:  {\small $\A_t^r$}, label=below:{ \small $h_t$}] at (2.7,0) {};
\draw[thick, latex- ] ([shift=(25:25mm)]0,0) arc (25:155:25mm);
\draw[thick,  -latex] ([shift=(-25:25mm)]0,0) arc (-25:-155:25mm);
\end{scope}
    \end{tikzpicture}
 }
\caption{Punctured disc with two special   points and pinnings $(\A_s^r, \A_t^l)$ and $(\A_t^r, \A_s^r)$.  
 }
\label{gs1}
\end{figure}

Recall the $\H-$invariants   at the special points $s, t$:  
\be \la{HSi}
\rho_{s}:= h(\A_s^r, \A_s^l) \in \H, \ \ \ \ \rho_{t}:= h(\A_t^r, \A_t^l) \in \H. 
\ee 
Recall the involution $h^*:= w_0(h^{-1})$ on the Cartan group,  
and  the {\it outer  monodromy map}:  
\begin{equation} \label{5}
\mu_{\rm out}: {\mathscr P}_{\G, \odot} \lra \H, \qquad \mu_{\rm out}:= \rho_{s}\rho_{t}^\ast.
\end{equation}

The  moduli space   ${\rm Loc}_{\G, \odot}$ 
is obtained from   ${\mathscr P}_{\G, \odot}$   by forgetting the   flag at the puncture.
Its $\H$-invariant is defined the same way:
\begin{equation} \label{5X}
\mu_{\rm out}: {\rm Loc}_{\G, \odot} \lra \H.
\end{equation}
\begin{definition} \label{8X} 
 The moduli space 
${\mathscr L}_{\G, \odot}$ is the fiber of the map (\ref{5X}) over the unit $e \in \H$. 
\ed

Unlike the map $\mu_{\rm out}$,  the space ${\mathscr L}_{\G, \odot}$ does not depend on the order of  special points $s_1, s_2$ on $\odot$.

Recall   the    Poisson-Lie group structure 
on $\G$. Its Drinfeld's dual is the Poisson-Lie group $\G^*$. 
As a group, $\G^*$  is a subgroup of   $\B_-\times \B_+$, given by the fiber  over $e$ of the canonical map $\B_-\times \B_+ \to \H$. 

Abusing notation, denote by ${\B}_\pm$  the standard Poisson Lie groups provided by the Borel subgroups $\B_\pm\subset \G$. 

 Recall the Drinfled double $\D(\B)$ of the Poisson Lie group ${\B}$. Denote by $\D(\B)^*$ its Poisson Lie dual group. 
 
\begin{theorem} \label{4}  
The moduli space 
 ${\rm Loc}_{\G, \odot}$  is a Poisson Lie group,    canonically isomorphic to   $\D(\B)^*$:
 \be \la{20aa}
 {\rm Loc}_{\G, \odot} =  \D(\B)^*.
 \ee
 The subspace ${\mathscr L}_{\G, \odot}\subset {\rm Loc}_{\G, \odot}$ is a Poisson Lie subgroup. 
 There is a  
 canonical Poisson  group isomorphism:
\be \la{20a}
\begin{split}
&\qquad {\mathscr L}_{\G, \odot} = \G^*.\\
\end{split}
\ee
 \end{theorem}

\begin{proof} 
  Let us show  that ${\rm Loc}_{\G, \odot}$ is a Poisson Lie group. The order of the special points determines the  
order of the two boundary intervals: we call the  one following $t$ clockwise, that is along the boundary orientation,  the left, and the second one the right. The 
identifications (\ref{20a})  depend on this choice.

\begin{figure}[ht]
\epsfxsize 200pt
\center{
\begin{tikzpicture}[scale=0.3]
\draw  (0,0) circle (20mm);
\draw [thick, blue]  (0,0) circle (4mm);
\node [red, thick] (a) at (0,2) {\small $\bullet$};
\node [red, thick] (b) at (0,-2) {\small $\bullet$};
\begin{scope}[shift={(5,0)}]
\draw  (0,0) circle (20mm);
\draw [thick, blue]  (0,0) circle (4mm);
\node [red] (a) at (0,2) {\small $\bullet$};
\node [red] (b) at (0,-2) {\small $\bullet$};
\end{scope}
\draw [thick, -latex] (8.5,0) -- (10.5,0);
\begin{scope}[shift={(14,0)}]
\draw [thick, blue]  (2,0) circle (4mm);
\draw [thick, blue]  (0,0) circle (4mm);
\draw ([shift=(90:20mm)]2,0) arc (90:-90:20mm);
\draw ([shift=(-90:20mm)]0,0) arc (-90:-270:20mm);
\draw[] (0,2) -- (2,2);
\draw[] (0,-2) -- (2,-2);
\draw[] (1,2) -- (1,-2);
\node [red] (a) at (1,2) {\small $\bullet$};
\node [red] (b) at (1,-2) {\small $\bullet$};
\end{scope}

\draw [thick, -latex] (19.5,0) -- (21.5,0);
\begin{scope}[shift={(25,0)}]
\draw [blue, thick]  (1,0) circle (7mm);
\draw ([shift=(90:20mm)]2,0) arc (90:-90:20mm);
\draw ([shift=(-90:20mm)]0,0) arc (-90:-270:20mm);
\draw (0,2) -- (2,2);
\draw (0,-2) -- (2,-2);
\node [red] (a) at (1,2) {\small $\bullet$};
\node [red] (b) at (1,-2) {\small $\bullet$};
\end{scope}
    \end{tikzpicture}
 }
 
 \caption{The gluing map, followed by cutting out the two  punctures.}
 \label{fga1}
\end{figure}

\subsubsection{The composition.} Given two  decorated surfaces $\odot$, 
let us glue the right boundary interval of the 
first with the left one of the second one, see Figure \ref{fga1}.  
We get a twice punctured disc with two special points $\odot \ast \odot$. 
Removing  a  disc 
containing the two punctures 
 we get an annuli  
 with two special points on the outer boundary. Shrinking the hole inside into a puncture we get   the decorated surface $\odot$.

 Using pinnings, we glue   framed $\G$-local systems given on each of the    two decorated  surfaces $\odot$, 
and then  restrict to the annulus,  see Figure \ref{fga2}. We get a composition map: 
\be \la{CMA}
{\rm Loc}_{\G, \odot} \times {\rm Loc}_{\G, \odot} \lra 
{\rm Loc}_{\G, \odot \ast \odot} \lra {\rm Loc}_{\G, \odot}.
\ee
It is evidently associative. It is a Poisson map since the gluing and   restriction maps in (\ref{CMA}) are Poisson. 

\begin{figure}[ht]
\epsfxsize 200pt
\center{
\begin{tikzpicture}[scale=0.35]
\draw   (0,0) circle (20mm);
\draw [blue, thick]  (0,0) circle (4mm);
\node [red, label=above:${\rm F}_+$] (a) at (0,2) {\small $\bullet$};
\node [red, label=below:${\rm F}_-$] (b) at (0,-2) {\small $\bullet$};
\node (g) at (2.6,0) {\small $p_2$};
\node (h) at (-2.6,0) {\small $p_1$};   
\begin{scope}[shift={(7,0)}]
\draw   (0,0) circle (20mm);
\draw [blue, thick]  (0,0) circle (4mm);
\node [red, label=above:${\rm F}_+$] (a) at (0,2) {\small $\bullet$};
\node [red, label=below:${\rm F}_-$] (b) at (0,-2) {\small $\bullet$};
\node  at (2.6,0) {\small $p_3$};
\node  at (-2.6,0) {\small $p_2$};  
\end{scope} 
\draw [thick, -latex] (12,0) -- (14,0);
\begin{scope}[shift={(18,0)}]
\draw   (0,0) circle (20mm);
\draw [blue, thick]  (0,0) circle (4mm);
\node [red, label=above:${\rm F}_+$] (a) at (0,2) {\small $\bullet$};
\node [red, label=below:${\rm F}_-$] (b) at (0,-2) {\small $\bullet$};
\node  at (2.6,0) {\small $p_3$};
\node  at (-2.6,0) {\small $p_1$};  
\end{scope}
    \end{tikzpicture}
 }
 
 \caption{The composition. Here and below we use the notation ${\rm F}$ for flags, rather than $\B$.}
 \label{fga2}
 \end{figure}

 \subsubsection{The unit.} It    is   
the  trivial $\G$-local system, equipped with pinnings $(p, p)$ where $p$ is a (standard) pinning, see Figure \ref{fga2id}.  
The outer monodromy is equal to $e$ since  $\rho_{s_\pm} = e$. Indeed, the pairs of decorated flags at each special point coincide by the construction. 
The composition with this object  does not change neither the $\G$-local system, nor the pinnings.

\begin{figure}[ht]
\epsfxsize 200pt
\center{
\begin{tikzpicture}[scale=0.35]
\draw   (0,0) circle (20mm);
\draw [blue, thick]  (0,0) circle (4mm);
\node [red, label=above:${\rm F}_+$] (a) at (0,2) {\small $\bullet$};
\node [red, label=below:${\rm F}_-$] (b) at (0,-2) {\small $\bullet$};
\draw[thick, latex- ] ([shift=(115:25mm)]0,0) arc (115:245:25mm);
\draw[thick,  latex-] ([shift=(65:25mm)]0,0) arc (65:-65:25mm);
\node (g) at (3,0) {\small $p$};
\node (h) at (-3,0) {\small $p$};   
    \end{tikzpicture}
 }
 
 \caption{The identity: on the picture both pinnings are given by a pinning $p: {\rm F}_- \lms {\rm F}_+$.}
 \label{fga2id}
 \end{figure}

\subsubsection{The inverse.} It is an involutive map of   moduli spaces, induced by   the disc rotation by $180^\circ$ 
around the diameter connecting the  two special points,   altering the disc orientation, see Figure \ref{fga2i}:
\be \la{inv}
{\rm Inv}: {\rm Loc}_{\G, \odot} \lra {\rm Loc}_{\G, \odot^\circ}. 
\ee


\begin{figure}[ht]
\epsfxsize 200pt
\center{
\begin{tikzpicture}[scale=0.35]
\draw   (0,0) circle (20mm);
\draw [blue, thick]  (0,0) circle (4mm);
\node [red, label=above:${\rm F}_+$] (a) at (0,2) {\small $\bullet$};
\node [red, label=below:${\rm F}_-$] (b) at (0,-2) {\small $\bullet$};
\node (g) at (2.6,0) {\small $q$};
\node (h) at (-2.6,0) {\small $p$};   
\draw [thick, -latex] (4,0) -- (6,0);
\begin{scope}[shift={(10,0)}]
\draw   (0,0) circle (20mm);
\draw [blue, thick]  (0,0) circle (4mm);
\node [red, label=above:${\rm F}_+$] (a) at (0,2) {\small $\bullet$};
\node [red, label=below:${\rm F}_-$] (b) at (0,-2) {\small $\bullet$};
\node  at (2.6,0) {\small $p$};
\node  at (-2.6,0) {\small $q$};  
\end{scope}
    \end{tikzpicture}
 }
  \caption{The inverse.}
 \label{fga2i}
 \end{figure}

\subsubsection{The  isomorphisms (\ref{20aa})}.  Denote the two flags 
at the  special points   
by ${\rm F}_+$ and ${\rm F}_-$, and  by $\A^l_\pm$ and $\A^r_\pm$ the two decorated flags near the flag ${\rm F}_\pm$, see Figure \ref{fga3}. 
 \begin{figure}[ht]
\epsfxsize 200pt
\center{
\begin{tikzpicture}[scale=0.35]
\draw   (0,0) circle (20mm);
\draw [blue, thick]  (0,0) circle (2.5mm);
\node [red, label=above:${\rm F}_+$] (a) at (0,2) {\small $\bullet$};
\node [red, label=below:${\rm F}_-$] (b) at (0,-2) {\small $\bullet$};
\node (i) at (1.7,2.3) {\small $ {\rm A}^r_+$};
\node (n) at (-1.6,2.28) {\small ${\rm A}^l_+$};  
\begin{scope}[shift={(9,0)}]
\draw  (0,0) circle (20mm);
\draw [blue, thick]  (0,0) circle (2.5mm);
\node [red, label=above:${{\rm F}}_+$] (a) at (0,2) {\small $\bullet$};
\node [red, label=below:${{\rm F}}_-$] (b) at (0,-2) {\small $\bullet$};  
\node  at (1.7,-2.2) {\small ${\A}^l_-$};
\node  at (-1.7,-2.2) {\small ${\A}^r_-$};  
\end{scope}
    \end{tikzpicture}
 }
 \caption{Identifying the  Lie group ${\rm Loc}_{{\rm G}, \odot}$ with the Lie group ${\rm B}_-\times {\rm B}_+$.}
 \label{fga3}
 \end{figure}  
 Moving 
 $\A_+^l$ and $\A_+^r$ to the flag $\F_-$ along the arcs connecting them with  ${\rm F}_-$ we get a triple 
 $\{\A_+^l,{\rm F}_-, \A_+^r\}$. Simiarly, we get a triple   $\{\A_-^l, {\rm F}_+, \A_-^r\}$.
   Since $\G$ acts freely and transitively on   pinnings,  there are unique $b_\pm \in \B_\pm$  such that 
\be \nonumber
\begin{split}
&\{{\rm F}_-, \A_+^r\} = b_- \cdot \{{\rm F}_-,  \A_+^l\},\\
&\{{\rm F}_+,  \A_-^l\} = b_+ \cdot \{{\rm F}_+,  \A_-^r\}.\\
\end{split}
\ee

Assigning to a point of ${\rm Loc}_{\G, \odot}$ the  pair 
$(b_-, b_+)$ we get a regular map 
\begin{equation} \label{31}
{f}: {\rm Loc}_{\G, \odot} \lra {\B}_- \times {\B}_+.
\end{equation}
The monodromy  of the local system   is given by the conjugacy class of $b_-^{-1}b_+$.    
 
The map (\ref{31}) is evidently a group homomorphism. 
To invert  the map (\ref{31}), pick a representative for every element in ${\rm Loc}_{\G, \odot}$ such that $p$ is a  standard  pinning. This data provides the restriction 
of the desired object of ${\rm Loc}_{\G, \odot}$ 
to the left  boundary interval. Using  elements $b_-, b_+$ we reconstruct the other pinning, getting a local system with the desired monodromy. 

Each of the two components of the map (\ref{31}) is a Poisson map by the definition of the Poisson structure on the space ${\mathscr P}_{\G, \odot}$, and hence ${\mathscr L}_{\G, \odot}$, 
 given in Section \ref{CSEC4.3} below.
 \be
{f}_-: {\rm Loc}_{\G, \odot} \lra  {\B}_-, \qquad 
  {f}_+: {\rm Loc}_{\G, \odot} \lra {\B}_+.
  \ee
Indeed, consider the map  
${\mathscr P}_{\G, \odot} \to {\B}_-$ is obtained by taking a  loop $\alpha_t$ around the puncture containing  $t$, and considering the  restriction map   to the triangle   given  by $\alpha_t$ and the two boundary intervals, see Figure \ref{amalgam99}. It is a Poisson map. 
\begin{figure}[ht]
\epsfxsize 200pt
\center{
\begin{tikzpicture}[scale=0.36]
\draw   (0,0) circle (20mm);
\draw [blue]  (0,0) circle (2mm);
\node [red] (a) at (0,2) {\small $\bullet$};
\node [red] (b) at (0,-2) {\small $\bullet$};
\draw[blue] (0,0.5) ellipse (6mm and 1.5cm);
\begin{scope}[shift={(10,0)}]
\draw  (0,0) circle (20mm);
\draw [blue]  (0,0) circle (2mm);
\node [red] (a) at (0,2) {\small $\bullet$};
\node [red] (b) at (0,-2) {\small $\bullet$};  
\draw[blue] (0,-0.5) ellipse (6mm and 1.5cm);
\end{scope}
    \end{tikzpicture}
 }
\caption{The maps $f_-$ and $f_+$.}
\label{amalgam99}
 \end{figure}  
It  evidently factorises as ${\mathscr P}_{\G, \odot} \to  {\rm Loc}_{\G, \odot}  \to {\B}_-$. Since the map ${\mathscr P}_{\G, \odot} \to  {\rm Loc}_{\G, \odot} $ is a Poisson finite cover at the generic point, 
the second map is Poisson. The map ${\mathscr P}_{\G, \odot} \to {\B}_+$ uses a similar loop $\alpha_s$ containing  $s$.

\subsubsection{The group ${\mathscr L}_{\G, \odot} $.}\la{thegr} 
The outer monodromy  (\ref{5})  
for the composition 
$ ({\cal L}, p )\circ ({\cal L}', p' ) $ 
is evidently the product of the outer monodromies 
of $({\cal L}, p )$ and $({\cal L}', p')$. 
In other words, there is a commutative diagram 
 \begin{displaymath}
    \xymatrix{
        {\rm Loc}_{\G, \odot} \times {\rm Loc}_{\G, \odot}  \ar[r]^{\qquad \circ} \ar[d]^{\mu_{\rm out}\qquad \qquad }_{\qquad \qquad \mu_{\rm out}} &{\rm Loc}_{\G, \odot}          \ar[d]^{\mu_{\rm out}} \\
          \H\times \H \ar[r]^{\qquad m} &    \H }
                \end{displaymath}
 where   the bottom  arrow is the multiplication. 
So the composition  induces 
a map 
\be \la{circm}
\circ: {\mathscr L}_{\G, \odot} \times {\mathscr L}_{\G, \odot} \lra 
{\mathscr L}_{\G, \odot}.  
\ee
To prove that 
${\mathscr L}_{\G, \odot}=\G^*$ it remains to check that 
under  identification (\ref{31})  the outer monodromy 
  is  $h_+h_-$, where 
$h_\pm$ is the Cartan part of $b_\pm$. 

\subsubsection{Poisson isomorphism property.} This can be proved by direct calculations with $r-$matrices, as was done for $\G$ in \cite{FG05}. 
On the other hand, we can just refer to \cite[Theorem 1.2]{S20}. Briefly, the argument there is this. Theorem \ref{UEAB} gives an algebra map ${\U}_q(\mathfrak{g}) \lra  {\cal O}_q({\mathscr L}_{\G, \odot})$.\footnote{The proof of Theorem \ref{UEAB} does not use Theorem \ref{4}.} 
The semi-classical limit $q\to1$ gives a Poisson algebra map
$$
\kappa: {\cal O}_q(\G^*) \lra  {\cal O}({\mathscr L}_{\G, \odot}).
$$
A priori, this map may not be surjective, see \cite[Remark 2.4]{S20}. 
In our case $\kappa$ is an isomorphism since there is  the inverse map  \cite[Formula (20)]{S20}. 
\end{proof}

\medskip
\subsection{The Poisson Lie group $\G$   $=$ the moduli space  of $\G-$local systems on a colored rectangle} \la{SECT4.2a}

 \medskip

\subsubsection{The Poisson Lie group $\G$.}  \la{PLGG}

Denote by  $\square$ the {\it colored rectangle}, defined as   
  a   disc  with four clockwise ordered special  points $s_1, \ldots, s_4$, and two colored boundary intervals $s_1s_2$ and $s_3s_4$.  
   ``Clockwise" means following the rectangle orientation. 

The moduli space ${\rm  Loc}_{\G, \square}$ parametrizes 
  $\G-$local systems on the rectangle  with two pinnings $p, q$ at the vertical sides, where $p$ is then pinning on the boundary interval $s_1s_2$. 
  Since any $\G-$local system on a rectangle is trivial, it is the same thing as isomorphism classes of pairs of pinnings $\{p, q\}$. 
 We stress that we  put no restriction to the pairs of  flags at the ends of the sides without pinnings. 
  
  Equivalently, it is the configuration space of quadruples of (decorated) flags $(\A_1, \B_2, \A_3, \B_4)$ 
at the vertices $s_1, \ldots,  s_4$, see Figure \ref{R5}. The pinnings are given by $p=(\A_1, \B_2)$ and $q = (\A_3, \B_4)$. Let $q^\ast$ be the opposite pinning of $q$.

  \begin{figure}[ht]
\begin{tikzpicture}[scale=1.5]
\begin{scope}[xshift=-.5cm]
\draw[thick] (0,0) -- (1,0);
\draw[thick] (0,1)--(1,1);
\draw[dashed,red, thick, -stealth] (0,0) -- (0,1);
\draw[dashed,red, thick, -stealth] (1,1) -- (1,0);
\node at (-.2, 0.5) {$p$};
\node at (1.2, 0.5) {$q$};
\end{scope}

\node at (1.3, 0.5) {$\sim$};
\begin{scope}[xshift=2cm]
\draw[thick] (0,0) -- (1,0);
\draw[thick] (0,1)--(1,1);
\draw[thick] (0,0) -- (0,1);
\draw[thick] (1,1) -- (1,0);
\node at (-.15, -.15) {${\rm A}_1$};
\node at (-.15, 1.15) {${\rm B}_2$};
\node at (1.15, -.15) {${\rm B}_4$};
\node at (1.15, 1.15) {${\rm A}_3$};
\end{scope}

\begin{scope}[xshift=6cm]
\draw[thick] (0,0) -- (1,0);
\draw[thick] (0,1)--(1,1);
\draw[dashed,red, thick, -stealth] (0,0) -- (0,1);
\draw[dashed,red, thick, -stealth] (1,1) -- (1,0);
\node at (-.2, 0.5) {$p$};
\node at (1.2, 0.5) {$q$};
\draw[black!40!green, -stealth, thick] (0.1,0.5) -- (0.9, 0.5);
\node at (0.5, 0.7) {$g$};
\end{scope}
\end{tikzpicture}
\caption{The moduli space ${\rm Loc}_{\G, \square}$, and the isomorphism ${\rm Loc}_{\G, \square}\stackrel{\sim}{=}\G$. }
\label{R5}
\end{figure}

 \bt  \la{Th4.3.3} The moduli space ${\rm Loc}_{\G, \square}$ is a Poisson Lie group, which is canonically isomorphic to 
 the  Poisson Lie group $\G$:
\be
{\rm Loc}_{\G, \square} \stackrel{\sim}{\lra}  \G.
\ee
\et

\begin{proof} The gluing map, followed by encircling the   middle special points, see Figure \ref{R6},  provides an  associative composition law
\be \la{COmm}
\circ:  {\rm Loc}_{\G, \square} \times {\rm Loc}_{\G, \square} \lra {\rm Loc}_{\G, \square}.
\ee
   \begin{figure}[ht]
\begin{tikzpicture}[scale=1.5]
\node at (4,0) {$\bullet$};
\node at (4,1) {$\bullet$};
\node at (5,0) {$\bullet$};
\node at (5,1) {$\bullet$};
\draw[thick] (4,0) -- (5,0);
\draw[thick] (4,1)--(5,1);
\draw[dashed, red, thick, stealth-]  (4,1) -- (4,0);
\draw[dashed, red, thick, -stealth]  (5,1) --(5,0);

\begin{scope}[xshift=2cm]
\node at (4,0) {$\bullet$};
\node at (4,1) {$\bullet$};
\node at (5,0) {$\bullet$};
\node at (5,1) {$\bullet$};
\draw[thick] (4,0) -- (5,0);
\draw[thick] (4,1)--(5,1);
\draw[dashed, red, thick, stealth-]  (4,1) -- (4,0);
\draw[dashed, red, thick, -stealth]  (5,1) --(5,0);
\end{scope}

\begin{scope}[xshift=5.2cm]
\draw[directed]  (2.7, 0.5) -- (3.5,0.5);
\node[red] at (3.1, 0.7) {glue};
\end{scope}

\begin{scope}[xshift=5.5cm]
\node at (3.7,0) {$\bullet$};
\node at (3.7,1) {$\bullet$};
\node at (5.3,0) {$\bullet$};
\node at (5.3,1) {$\bullet$};
\node at (4.5,0) {$\bullet$};
\node at (4.5,1) {$\bullet$};
\draw[thick] (3.7,0) -- (5.3,0);
\draw[thick] (3.7,1)--(5.3,1);
\draw[dashed, red, thick, stealth-]  (3.7,1) --(3.7,0);
\draw[dashed, red, thick, -stealth]  (5.3,1) --(5.3,0);
\draw[blue, dashed] (3.7,1) arc (255:285:3.2);
\draw[blue, dashed] (3.7,0) arc (105:75:3.2);
\end{scope}

\begin{scope}[xshift=8.4cm]
\draw[directed]  (2.7, 0.5) -- (3.5,0.5);
\node[blue] at (3.1, 0.7) {cut};
\end{scope}

\begin{scope}[xshift=8.5cm]
\node at (4,0) {$\bullet$};
\node at (4,1) {$\bullet$};
\node at (5,0) {$\bullet$};
\node at (5,1) {$\bullet$};
\draw[thick] (4,0) -- (5,0);
\draw[thick] (4,1)--(5,1);
\draw[dashed, red, thick, stealth-]  (4,1) -- (4,0);
\draw[dashed, red, thick, -stealth]  (5,1) --(5,0);
\end{scope}
\end{tikzpicture}
\caption{The composition law $\circ: {\rm Loc}_{\G, \square} \times {\rm Loc}_{\G, \square} \lra {\rm Loc}_{\G,\square}$ .}
\label{R6}
\end{figure}
 The unit is given by   $p=q^\ast$. 
 The inverse map is given by the reflection of the rectangle with respect to a  line parallel to the colored sides. So  ${\rm Loc}_{\G, \square}$ is a Lie group. 
There is an isomorphism of   Lie groups, see Figure \ref{R5}:
  \be \la{phi}
  \begin{split}
   \psi:  & {\rm Loc}_{\G, \square}\stackrel{\sim}{\lra} \G, \\
 & \{p, q\} \lms g, ~~~~~~~~\qquad q^\ast=gp, ~~g \in \G.\\  
 \end{split}
     \ee
  More precisely, we shall fix a representative in ${\rm Loc}_{\G, \square}$ such that $p$ is a standard pinning. Then there is a unique $g\in \G$ taking $p$ to $q^\ast$. Further details about the transporting $g$ can be found in \cite[Sec. 6.1]{GS13}.    
  
  The claim that $\psi$ is an isomorphism of Poisson Lie groups follows immediately from \cite{FG05}.
\end{proof} 

\medskip

\subsection{Simply-connected Poisson Lie group $\G'$ carries both cluster Poisson $\&$   $K_2-$structures} \la{SEC4.4}

\medskip

Recall  the canonical projection 
$$
\pi:\G'\to \G.
$$
Note that ${\mathscr A}_{\G', \square} = {\mathscr  P}_{\G', \square}$ by Lemma \ref{EQQ}. 
 
 \bt \la{BZFG} For the adjoint group  $\G$,  the moduli space ${\mathscr P}_{\G, \square}=\G$ carries a cluster Poisson structure. 

 For the simply connectedgroup $\G'$,  the moduli space ${\mathscr A}_{\G', \square} = {\mathscr P}_{\G', \square}=\G'$ carries a cluster $K_2-$structure. 
 
 The  two cluster structures are compatible via the map $\pi^*$ in the sense of Section \ref{sec18}. 
\et

\begin{proof} In each of the cases, the cluster structure is given by the cluster structure on the moduli spaces studied in Sections \ref{pcs.sec} - \ref{CSEC7}. 
For the ${\mathscr P}_{\G,  \square}$, it is the cluster structure studied in Section \ref{CSEC7}. 

For the space ${\mathscr A}_{\G', \square}$, it is the cluster structure on ${\rm Conf}^{w_0, w_0}({\cal A})$, see Definition \ref{def10.1} and Theorem \ref{main.result1.a.structure}. 
\vskip 1mm

The definition of compatibility is stated at the beginning of Section \ref{sec18.1}.

The map $\pi$ is finite  and surjective. Locally, when we fix a seed, we get a monomial map whose pull pack is
$\pi^* X_i = \prod A_j^{p_{ij}}$. 
The matrix $p_{ij}$ is invertible and satisfies Lemma \ref{Lem18.1} since $\pi$ is a finite  surjective map. 
Therefore, we get a compatible pair as in Section  \ref{sec18.1}.
\end{proof}

Therefore the simply-connected group $\G'$ has both the cluster $K_2-$structure and 
a Poisson structure $\pi_{\G'}:= \pi^*(\pi_\G)$,  induced   from the Poisson structure $\pi_\G$ on $\G$. 

This allows us to define quantum algebras: 
\be
\begin{split}
&{\cal O}_q(\G):= {\cal O}_q({\mathscr P}_{\G, \square}).\\
&{\cal O}^{\rm cl}_q(\G'):= {\cal O}^{\rm cl}_q({\mathscr A}_{\G', \square}).\\
\end{split}
\ee

Indeed,  the 
quantization of cluster Poisson varieties  in \cite{FG03b} provides us the algebra ${\cal O}_q(\G)$.

On the other hand,   the Poisson structure $\pi_{\G'}$ on $\G'$ is compatible with the cluster $K_2-$structure in the sense of \cite{BZ}, see also Section \ref{sec18}.  
Therefore it can be $q-$deformed in the sense of \cite{BZ}, providing a distinguished $q-$deformation 
$  {\cal O}_q^{\rm cl}({\mathscr A}_{\G', \square})$ of the upper cluster algebra $  {\cal O}^{\rm cl}({\mathscr A}_{\G', \square})$.  

Note that the  Poisson structure $\pi_{\G'}$   belongs to the family of Poisson structures compatible   with the cluster algebra structure on $\G'$.  
However it is  the only one relevant to the classical quantum group   deformation ${\cal O}_q(\G')$ of the algebra of regular functions on $\G'$.

In Section  \ref{sec18.3}  
we show that the  quantization in \cite{BZ} is compatible with the one from \cite{FG03b}. Therefore there is an injective map  of algebras
\be
{\cal O}_q({\mathscr P}_{\G, \square}) \hra  {\cal O}_q^{\rm cl}({\mathscr A}_{\G', \square}).
\ee

\bcon \la{GYC}
 The algebra ${\cal O}_q(\G')$ is isomorphic to the subalgebra  of the quantized upper cluster algebra ${\cal O}^{\rm cl}_q(\G')$, defined without inverting the frozen variables. 
 
 For the adjoint group $\G$, we shall have 
 \be
 {\cal O}_q(\G) = {\cal O}_q(\G') \cap \pi^* ({\cal O}_q(\G^{w_0, w_0})).
 \ee 
 \econ
 The work  of Goodearl-Yakimov \cite{GY} provides an  evidence for Conjecture \ref{GYC}. 

  For the simply-connected $\G'$, the frozen cluster ${\cal A}-$variables $\A_1, ..., \A_r$ and $\A'_1, ..., \A'_r$, assigned to the non-colored boundary intervals of $\square$, 
   naturally identify the $2r$ boundary divisors of $\G'$.

 \begin{figure}[ht]
\begin{tikzpicture}[scale=1.7]
\draw (0,0) -- (0,1) -- (1,1) -- (1,0) -- (0,0);
\draw (0,1) -- (1,0);
\node[red] at (-.1,-.1) {${0 \choose 1}$};
\node[red] at (-.1, 1.1) {${1 \choose 0}$};
\node[red] at (1.1, -.1) {${b \choose d}$};
\node[red] at (1.1, 1.1) {${a \choose c}$};
\node[blue] at (0.5, 0) {$\bullet$};
\node[blue] at (.5, .5) {$\bullet$};
\node[blue] at (.5, 1) {$\bullet$};
\draw[-stealth, thick, teal] (.5, .95) -- (.5, .55);
\draw[-stealth, thick, teal] (.5, .05) -- (.5, .45);
\node at (.6, 1.1) {$c$};
\node at (.6, .5) {$d$};
\node at (.6, -.1) {$b$};
\end{tikzpicture}
\caption{The cluster  coordinates on ${\rm SL}_2$.}
\label{R2000}
\end{figure}

{\bf Example}.  The cluster $K_2-$structure on ${\rm SL}_2$ is depicted on Figure \ref{R2000}. Here $   \begin{pmatrix}       a & b \\
      c & d \\
   \end{pmatrix} \in {\rm SL}_2$. It has two clusters, corresponding to the two triangulations of the square.
   The frozen variables $b, c$ are assigned to the horizontal = non-colored boundary intervals.  We have $ad=bc+1$, providing the exchange relation.\vskip 2mm

The projection $\pi: {\rm SL}_2 \rightarrow {\rm PGL}_2$ is described on Figure \ref{R21}. Here 
\[
\pi^*x=\frac{b}{d}, \qquad
\pi^*y=\frac{1}{bc}; \qquad
\pi^*z=\frac{c}{d}
\]
Note that the frozen variables $x,z$ are assigned to the vertical = colored boundary intervals.
\begin{figure}[hb]
\begin{tikzpicture}[scale=1.7]
\draw (0,0) -- (0,1) -- (1,1) -- (1,0) -- (0,0);
\draw (0,1) -- (1,0);
\node[blue] at (0.5, 0) {$\bullet$};
\node[blue] at (.5, .5) {$\bullet$};
\node[blue] at (.5, 1) {$\bullet$};
\draw[-stealth, thick, teal] (.5, .95) -- (.5, .55);
\draw[-stealth, thick, teal] (.5, .05) -- (.5, .45);
\node at (.6, 1.1) {$c$};
\node at (.6, .5) {$d$};
\node at (.6, -.1) {$b$};
\draw[-latex] (1.5,.5) -- (2.5,.5);
\begin{scope}[xshift=3cm]
\draw (0,0) -- (0,1) -- (1,1) -- (1,0) -- (0,0);
\draw (0,1) -- (1,0);
\node[blue] at (1, 0.5) {$\bullet$};
\node[blue] at (.5, .5) {$\bullet$};
\node[blue] at (0, .5) {$\bullet$};
\draw[-stealth, thick, teal] (.45, .5) -- (.05, .5);
\draw[-stealth, thick, teal] (.55, .5) -- (.95, .5);
\node at (-.08, .35) {$x$};
\node at (.5, .35) {$y$};
\node at (1.08, .35) {$z$};
\end{scope}
\end{tikzpicture}
\caption{Cluster  description of the projection $\pi: {\rm SL}_2 \lra {\rm PGL}_2$.}
\label{R21}
\end{figure}

\medskip
\subsection{The Drinfeld double $\D(\G)$ of   $\G$ via moduli spaces.} \la{SECT4.2} 

\medskip

Denote by  $\boxed{\cdot}$ a {\it punctured colored rectangle}, given by    
  a   punctured disc  with  four ordered special  points $s_1, \ldots, s_4$, and two  colored sides $s_1s_2$ and $s_3s_4$.  
The moduli space ${\mathscr P}_{\G, \boxed{\cdot}}$ parametrizes $\G-$local systems on the punctured disc   with an invariant flag   near  the puncture, 
and   pinnings $(p, q)$ at the colored sides. Forgetting the flag at the puncture, we get the moduli space ${\rm Loc}_{\G, \boxed{\cdot}}$, see 
  Figure \ref{R4}. 

  \begin{figure}[ht]
\begin{tikzpicture}[scale=1.5]
\draw[thick] (4,0) -- (5,0);
\draw[thick] (4,1)--(5,1);
\draw[dashed, red, thick, stealth-]  (4,1) arc (150:210:1);
\draw[dashed, red, thick, -stealth]  (5,1) arc (30:-30:1);
\node[blue] at (4.5, 0.5) {$\circ$};
\node at (3.7, 0.5) {$p$};
\node at (5.3, 0.5) {$q$};
\end{tikzpicture}
\caption{The moduli space  ${\rm Loc}_{\G, \boxed{\cdot}}$ for a punctured rectangle $\boxed{\cdot}$ with two colored  boundary intervals.}
\label{R4}
\end{figure}


Recall the Drinfled double $\D(\G)$ of the Poisson Lie group $\G$. There is a Lie group isomorphism:
\be
\delta: \D(\G) \stackrel{\sim}{\lra} \G \times \G.
\ee
However it is not a  Poisson isomorphism. 

\bt  The moduli space ${\rm Loc}_{\G, \boxed{\cdot}}$ is a   Poisson Lie group. 

It is    canonically isomorphic, as a Poisson Lie group,  to the Drinfeld double $\D(\G)$ of  $\G$:
\be
\varphi: {\rm Loc}_{\G, \boxed{\cdot}} \stackrel{\sim}{\lra} \D(\G).
\ee

The moduli spaces ${\mathscr L}_{\G, \odot}$ and ${\rm Loc}_{\G, \square}$ are Poisson Lie subgroups in the Poisson Lie group ${\rm Loc}_{\G, \boxed{\cdot}}$.
\et

\begin{proof}
 The gluing map, followed by encircling of the two punctures, as well as  the middle special points on the non-colored sides, see Figure \ref{R3}, provides an  associative group law:
\be \la{COmm1}
\circ: {\rm Loc}_{\G, \boxed{\cdot}} \times {\rm Loc}_{\G, \boxed{\cdot}} \lra {\rm Loc}_{\G, \boxed{\cdot}}.
\ee
   \begin{figure}[ht]
\begin{tikzpicture}[scale=1.5]
\node at (4,0) {$\bullet$};
\node at (4,1) {$\bullet$};
\node at (5,0) {$\bullet$};
\node at (5,1) {$\bullet$};
\draw[thick] (4,0) -- (5,0);
\draw[thick] (4,1)--(5,1);
\draw[dashed, red, thick, stealth-]  (4,1) arc (150:210:1);
\draw[dashed, red, thick, -stealth]  (5,1) arc (30:-30:1);
\node[blue] at (4.5, 0.5) {$\circ$};
\node at (3.7, 0.5) {$p_1$};
\node at (5.3, 0.5) {$q_1$};

\begin{scope}[xshift=2cm]
\node at (4,0) {$\bullet$};
\node at (4,1) {$\bullet$};
\node at (5,0) {$\bullet$};
\node at (5,1) {$\bullet$};
\draw[thick] (4,0) -- (5,0);
\draw[thick] (4,1)--(5,1);
\draw[dashed, red, thick, stealth-]  (4,1) arc (150:210:1);
\draw[dashed, red, thick, -stealth]  (5,1) arc (30:-30:1);
\node[blue] at (4.5, 0.5) {$\circ$};
\node at (3.7, 0.5) {$p_2$};
\node at (5.3, 0.5) {$q_2$};
\end{scope}

\begin{scope}[xshift=5.2cm]
\draw[directed]  (2.7, 0.5) -- (3.5,0.5);
\node[red] at (3.1, 0.7) {glue};
\end{scope}

\begin{scope}[xshift=5.5cm]
\node at (4,0) {$\bullet$};
\node at (4,1) {$\bullet$};
\node at (5,0) {$\bullet$};
\node at (5,1) {$\bullet$};
\node at (4.5,0) {$\bullet$};
\node at (4.5,1) {$\bullet$};
\draw[thick] (4,0) -- (5,0);
\draw[thick] (4,1)--(5,1);
\draw[dashed, red, thick, stealth-]  (4,1) arc (150:210:1);
\draw[dashed, red, thick, -stealth]  (5,1) arc (30:-30:1);
\node[blue] at (4.3, 0.5) {$\circ$};
\node[blue] at (4.7, 0.5) {$\circ$};
\node at (3.7, 0.5) {$p_1$};
\node at (5.3, 0.5) {$q_2$};
\draw[blue, dashed] (4,1) arc (245:295:1.25);
\draw[blue, dashed] (4,0) arc (115:65:1.25);
\draw[blue,dashed] (4.5, 0.5) ellipse (0.4cm and 0.2cm);
\end{scope}

\begin{scope}[xshift=8.4cm]
\draw[directed]  (2.7, 0.5) -- (3.5,0.5);
\end{scope}

\begin{scope}[xshift=8.5cm]
\node at (4,0) {$\bullet$};
\node at (4,1) {$\bullet$};
\node at (5,0) {$\bullet$};
\node at (5,1) {$\bullet$};
\draw[thick]  (4,0) -- (5,0);
\draw[thick]  (4,1)--(5,1);
\node[blue] at (3.0, 0.7) {cut};
\draw[dashed, red, thick, stealth-]  (4,1) arc (150:210:1);
\draw[dashed, red, thick, -stealth]  (5,1) arc (30:-30:1);
\node[blue] at (4.5, 0.5) {$\circ$};
\node at (3.7, 0.5) {$p_1$};
\node at (5.3, 0.5) {$q_2$};
\end{scope}
\end{tikzpicture}
\caption{The composition law $\circ:  {\rm Loc}_{\G, \boxed{\cdot}} \times {\rm Loc}_{\G, \boxed{\cdot}} \lra {\rm Loc}_{\G,\boxed{\cdot}}$ .}
\label{R3}
\end{figure}
 The unit is given by the trivial $\G-$local system on $\boxed{\cdot}$ with $p=q^\ast$. 
 The inverse map is given by the reflection of the rectangle with respect to the a   line parallel to  colored sides. So   ${\mathscr L}_{\G, \boxed{\cdot}}$ is a Lie group.

  \begin{figure}[ht]
\begin{tikzpicture}[scale=1.5]
\draw[thick] (4,0) -- (5,0);
\draw[thick] (4,1)--(5,1);
\draw[dashed, red, thick, stealth-]  (4,1) arc (150:210:1);
\draw[dashed, red, thick, -stealth]  (5,1) arc (30:-30:1);
\node[blue] at (4.5, 0.5) {$\circ$};
\node at (3.7, 0.5) {$p$};
\node at (5.3, 0.5) {$q$};
\draw[black!40!green, thick, -stealth]  (3.9,.6) arc (110:70:1.75);
\draw[black!40!green, thick, -stealth]  (3.9,0.4) arc (-110:-70:1.75);
\node at (4.5, 0.83) {$g_+$}; 
\node at (4.5, 0.17) {$g_-$}; 
\end{tikzpicture}
\caption{The isomorphism $\varphi:   {\rm Loc}_{\G, \boxed{\cdot}}\lra \G \times \G$. }
\label{R2}
\end{figure}
     
     There is a canonical isomorphism of Lie groups, illustrated on Figure \ref{R2}:  
  \be \la{psi1}
  \begin{split}
& \varphi:   {\rm Loc}_{\G, \boxed{\cdot}}\stackrel{\sim}{\lra} \G\times\G, \\
 & ({\mathscr L}, p, q) \lms (g_+, g_-)\\  
 \end{split}
     \ee
   Precisely, given a $\G-$local system ${\cal L}$ on the punctured rectangle together with   pinnings $\{p, q\}$, 
   the parallel transports of the pinning $p$  above/below the puncture towards the right side, see Figure \ref{R2}, provides new pinnings $p_{+}$ and $p_{-}$ on the right side. Then   
     $q^\ast=g_+p_+$ and $q^\ast=g_-q_-$ where $g_+, g_- \in \G$.    
     
     \medskip
     
        The composition map  (\ref{COmm1}) is Poisson. Indeed,  denote by $\boxdot_2$ the twice punctured rectangle with two opposite colored sides, pictured as the vertical sides. Then the gluing of colored decorated  surfaces 
 $$
 \boxed{\cdot}\times \boxed{\cdot} \lra \boxdot_2
 $$
 induces a map of the moduli spaces,  provided by gluing the  middle pinnings, and encircling the special points on the non-colored sides, illustrated by the first arrow on Figure \ref{R3}:

 \be \la{COmm}
{\mathscr P}_{\G, \boxed{\cdot}} \times {\mathscr P}_{\G, \boxed{\cdot}} \lra {\mathscr P}_{\G, \boxdot_2}.
\ee
\bl \la{4.5} The map (\ref{COmm}) is   a cluster Poisson projection map. 
\el

\begin{proof} We present the middle twice punctured colored square $\boxdot_2$ on 
Figure \ref{R2} by gluing   the internal rectangle $r_{\rm int}$ with the top and the bottom  triangles $t_{\rm top}$ and $t_{\rm bot}$:
\be
\boxdot_2 = r_{\rm int} \ast t_{\rm top} \ast t_{\rm bot}.
\ee
So the cluster Poisson structure on ${\mathscr P}_{\G, \boxdot_2}$ is the amalgamation of the ones on ${\mathscr P}_{\G, r_{\rm int}}$, ${\mathscr P}_{\G,  t_{\rm top}}$, and 
${\mathscr P}_{\G, t_{\rm bot}}$. Then   projection (\ref{COmm}) is described by forgetting the cluster coordinates related to   triangles $t_{\rm top}$ and $t_{\rm bot}$. 
 \end{proof}

Forgetting the flags at the punctures, we get a map 
\be \la{COmma}
{\rm Loc}_{\G, \boxed{\cdot}} \times {\rm Loc}_{\G, \boxed{\cdot}} \lra {\rm Loc}_{\G, \boxdot_2}.
\ee
It is a Poisson map at the generic point, since at the generic point it is obtained by taking the quotient under the action of the group $W \times W$, and this action   cluster Poisson. Therefore it is a Poisson map. 
Then encircling the two punctures  we get the composition map, see Figure \ref{R3}. 
This is a Poisson map since the restriction map is Poisson. 
   The inverse map is induced by the reflection   of the rectangle, which alters the orientation of the rectangle. 
   So it is an anti-Poisson map. So 
 ${\rm Loc}_{\G, \boxed{\cdot}}$ is a Poisson Lie group. 
    
       \medskip
        Let us recall the isomorphism of Poisson Lie groups defined in Section \ref{SEC1.7}:
     \be \la{isoG*}
\xi:  {\mathscr L}_{\G, \odot}\stackrel{\sim}{\lra} \G^*.     \ee
     
  \medskip
  
 The gluing maps on Figure \ref{R1} provide the $2^r:1$ cover  maps 
\be \la{psi}
  \begin{split}
  &  \gamma_\alpha:  {\rm Loc}_{\G, \square }\times  {\mathscr L}_{\G, {\odot}} \lra  {\rm Loc}_{\G, \boxed{\cdot}}, \\
   & \gamma_\beta:   {\mathscr L}_{\G, {\odot}}\times  {\rm Loc}_{\G, \square }\ \lra  {\rm Loc}_{\G, \boxed{\cdot}}.\\
       \end{split}
 \ee
  Same argument as above shows that they are Poisson maps. \vskip 1mm
 
Using isomorphisms (\ref{phi}), (\ref{psi1}), (\ref{isoG*}) we identify them with the Poisson $2^r:1$ cover maps 
 \be \la{psi2}
  \begin{split}
  &  \alpha:  \G\times  \G^*\lra  \G \times \G, \\
   & \beta:   \G^*\times  \G \lra  \G \times \G.\\
       \end{split}
 \ee  
They are not group maps. \vskip 1mm

Gluing with  trivial $\G-$local systems with pinnings $p=q$, representing the units  $e \in {\mathscr L}_{\G, {\odot}}  $ and 
 $e\in {\rm Loc}_{\G, \square } $, provide the canonical embeddings  (\ref{PPDD}) of Poisson Lie groups:
\be  
 \begin{split}
 &  \gamma_\alpha(\ast, e): \G=  {\rm Loc}_{\G, \square }  \hra  {\rm Loc}_{\G, \boxed{\cdot}} = \G \times \G, \\
 &\gamma_\beta(\ast, e): \G^*=   {\mathscr L}_{\G, {\odot}}   \hra  {\rm Loc}_{\G, \boxed{\cdot}}= \G \times \G. \\
       \end{split}
       \ee

 \begin{figure}[ht]
\begin{tikzpicture}[scale=1.5]
\begin{scope}[xshift=-.5cm]
\draw[thick] (0,0) -- (1,0);
\draw[thick] (0,1)--(1,1);
\draw[dashed,red, thick, -stealth] (0,0) -- (0,1);
\draw[dashed,red, thick, -stealth] (1,1) -- (1,0);
\node at (-.2, 0.5) {$p_1$};
\node at (1.2, 0.5) {$q_1$};
\end{scope}

\node at (1, 0.5) {$\times$};

\begin{scope}[xshift=-.2cm]
\draw[dashed, red, thick, stealth-]  (2,1) arc (120:240:0.6);
\draw[dashed, red, thick, -stealth]  (2,1) arc (60:-60:0.6);
\node[blue] at (2, 0.5) {$\circ$};
\node at (1.55, 0.5) {$p_2$};
\node at (2.5, 0.5) {$q_2$};
\end{scope}

\draw[thick] (4,0) -- (5,0);
\draw[thick] (4,1)--(5,1);
\draw[dashed, red, thick, stealth-]  (4,1) arc (150:210:1);
\draw[dashed, red, thick, -stealth]  (5,1) arc (30:-30:1);
\node[blue] at (4.5, 0.5) {$\circ$};

\begin{scope}[xshift=5.2cm]
\draw[dashed, red, thick, stealth-]  (2,1) arc (120:240:0.6);
\draw[dashed, red, thick, -stealth]  (2,1) arc (60:-60:0.6);
\node[blue] at (2, 0.5) {$\circ$};
\end{scope}

\node at (7.96, 0.5) {$\times$};

\begin{scope}[xshift=8.5cm]
\draw[thick] (0,0) -- (1,0);
\draw[thick] (0,1)--(1,1);
\draw[dashed,red, thick, -stealth] (0,0) -- (0,1);
\draw[dashed,red, thick, -stealth] (1,1) -- (1,0);
\end{scope}

\draw[directed] (2.7,0.5) -- (3.5, 0.5);
\node[red] at (3.1, 0.7) {glue};
\node at (3.1, 0.3) {$2^r:1$};

\begin{scope}[xshift=2.8cm]
\draw[directed]  (3.5, 0.5) -- (2.7,0.5);
\node[red] at (3.1, 0.7) {glue};
\node at (3.1, 0.3) {$2^r:1$};
\draw[directed] (3.5,-1.6) -- (2.7, -1.6);
\node at (3.1, -1.4) {$\beta$};
\end{scope}

\draw[directed] (1.4, -.3) -- (1.4, -1.3); 
\draw[directed] (4.5, -.3) -- (4.5, -1.3); 
\draw[directed] (7.6, -.3) -- (7.6, -1.3); 

\node at (1.4,-1.6) {$\G\times \G^\ast$};
\node at (4.5,-1.6) {$\G\times \G$};
\node at (7.6,-1.6) {$\G^\ast\times \G$};
\draw[directed] (2.7,-1.6) -- (3.5, -1.6);
\node at (3.1, -1.4) {$\alpha$};
\node at (1.6, -.8) {=};
\node at (7.8, -.8) {=};

\node at (3.72, 0.5) {$p_1$};
\node at (5.3, 0.5) {$q_2$};

\node at (6.75, 0.5) {$p_1$};
\node at (7.65, 0.5) {$q_1$};
\node at (8.3, 0.5) {$p_2$};
\node at (9.7, 0.5) {$q_2$};
\end{tikzpicture}
\caption{The geometric incarnation of the commutative diagram (\ref{Uhg}).}
\label{R1}
\end{figure}

We get a commutative diagram,  where the vertical maps are isomorphisms, and the horisontal ones are $2^r:1$ covers, also illustrated  it on  Figure \ref{R1}:   
  \be
\la{Uhg}
\begin{gathered}
 \xymatrix{     
 {\rm Loc}_{\G, \square }\times  {\mathscr L }_{\G, {\odot}}    \ar[r]^{~~~~~~~~ \gamma_{\alpha}}  \ar[d]^{\xi}_{\psi}   \ar[r]_{ } &   {\rm Loc}_{\G, \boxed{\cdot}} \ar[d]_{\sim}^{\varphi} &  \ar[l]^{ }_{\gamma_{\beta}~~~~}   \ar[d]_{\xi}^{\psi} {\mathscr L}_{\G, {\odot}}\times  {\rm Loc}_{\G, \square }\\
\G  \times  \G^* \ar[r]_{2^r:1}^{\alpha}&\G \times \G & \ar[l]^{2^r:1}_{\beta} \G^*\times  \G \\}
\end{gathered}
\ee  
 The maps $\psi$ and $\xi$ are isomorphisms of   Poisson Lie groups. The horizontal maps are Poisson since they are   gluing maps. 
The maps $\alpha$ and $\beta$ are   Poisson. Therefore the isomorphism $\varphi$ is a Poisson map. \end{proof}

\medskip
 \subsection{Cluster Poisson structure of   Grothendieck's resolution  $\widehat{\G} \lra \G$ via moduli spaces} \la{Grot}
\medskip

Consider the following well known  $2^r:1$  map
   \be \la{mpp1}
 \begin{split}
 & p': \G^*\lra  \G ,\\
 &(b_+, b_-) \in \G^* \lms b_+b_-^{-1} \in \G.\\
  \end{split}
   \ee
   
 Recall the  2-torsion subgroup $S_\G:= \H[2]$ of    Cartan group $\H$. The  map (\ref{mpp1})   can be described  as the composition 
 $$
 \G^* \lra \G^*/S_\G \hra \G
  $$
 of the quotient by the action of the   group $S_\G$, followed by an open embedding to $\G$. \vskip 1mm

 It is  known \cite{EL} that there is a unique a Poisson structure $\pi$ on   $\G$ such that the map $p'$ is   Poisson.   
 However $(\G, \pi)$ is not a Poisson Lie group: the  Poisson structure $\pi$  is not compatible   with the   product in $\G$. So we arrive at the maps of Poisson spaces
 \be
  \G^* \lra \G^*/S_\G \hra (\G, \pi).
   \ee

Our next goal is to give   an interpretation of these  results via the moduli spaces. \\
 
 Denote by ${\rm D}^*_1$ the punctured disc with two special points and one colored boundary segment, see Figure \ref{R19}. 
 It provides    moduli spaces    
 ${\mathscr P}_{\G, {\rm D}^*_1}$ and ${\rm Loc}_{\G, {\rm D}^*_1}$.  
      Forgetting the second pinning, we get a  Poisson map  
  \be \la{mpp}
 \begin{split}
 & p: {\mathscr L}_{\G, \odot} \lra  {\rm Loc}_{\G, {\rm D}^*_1}.\\
  \end{split}
     \ee

 \begin{figure}[ht]
\begin{tikzpicture}[scale=2]
\draw[dashed, red, thick, stealth-]  (1.97,1) arc (110:250:0.6);
\draw (2.03,1) arc (70:-70:0.6);
\node[blue] at (2, .48) {$\circ$};
\node at (2, -.12) {$\bullet$};
\node at (2, 1) {$\bullet$};
\node[red] at (1.4, .5) {$p$};
\draw[thick, black!40!green, -stealth] (1.6,.51) to[out=45, in=180] (1.9, 0.65) to[out=0, in=120] (2.12, 0.55) to[out=-60, in=60] (2.12, .41) to[out=-120, in=0] (1.9, 0.31) to[out=180, in=-45] (1.6, .45);
\node at (2.2, 0.5) {${\tiny g}$};
\end{tikzpicture}
\caption{The moduli space  ${\rm Loc}_{\G, {\rm D}_1^*}$ for a punctured disc with two special points and one  colored  boundary interval, and the isomorphism ${\rm Loc}_{\G, {\rm D}_1^*} = \G$.}
\label{R19}
\end{figure}

 Forgetting   the second pinning is equivalent to taking the quotient by the action of the group $S_\G$ on 
${\mathscr L}_{\G, \odot}$. Therefore  the map $p$  is the composition:  
\be \la{152}
 {\mathscr L}_{\G, \odot} \lra   {\mathscr L}_{\G, \odot}/S_\G \hra  {\rm Loc}_{\G, {\rm D}^*_1}.
 \ee
  Here the first map is  $2^r:1$, and the second  is a Zariski open embedding. \\

 Recall   the Grothendieck resolution $\widehat \G$ of the group $\G$.   

\bt i) There are canonical isomorphisms of moduli spaces:
 \be \la{cips}
 \begin{split}
 &{\rm Loc}_{\G, {\rm D}^*_1} = \G;\\ 
 & {\mathscr P}_{\G, {\rm D}^*_1} = \widehat \G.\\
 \end{split}
  \ee
  \vskip 1mm
   ii) 
    Combined with the isomorphism   ${\mathscr L}_{\G, \odot}  = \G^*$, we get  
   a commutative diagram of Poisson spaces:
   \begin{displaymath}
    \xymatrix{
         {\mathscr L}_{\G, \odot} \ar[r]^{} \ar[d]^{=}&{\mathscr L}_{\G, \odot}/S_\G\ar[r]^{\subset}\ar[d]^{=} &   {\rm Loc}_{\G, {\rm D}^*_1}\  \ar[d]^{=}\\
           \G^*  \ar[r]^{} &   \G^*/S_\G \ar[r]^{\subset}& (\G, \pi)  }
                \end{displaymath} 
In particular, the first isomorphism (\ref{cips})  provides    a Poisson structure   on  $\G$, which coincides with the one  $\pi$:
\be \la{323}
  {\rm Loc}_{\G, {\rm D}^*_1}  = (\G, \pi).  
\ee
\vskip 1mm

              iii)   The second isomorphism   provides  a   $W-$equivariant cluster Poisson variety structure on $\widehat \G$.
    
So  we get a commutative diagram of Poisson moduli spaces,    describing    Grothendieck's resolution $\widehat \G \to \G$: 
      \begin{displaymath}
    \xymatrix{
         {\mathscr P}_{\G, {\rm D}^*_1} \ar[r]^{} \ar[d]^{=}  &   {\rm Loc}_{\G, {\rm D}^*_1}  \ar[d]^{=}\\
           \widehat \G  \ar[r]^{} &   (\G, \pi) }
                \end{displaymath}

       \et
  
  \begin{proof} 
i) Look at  Figure \ref{R19}. Denote by $p'$ the pinning at the colored boundary interval obtained by the counterclockwise parallel transport of the pinning $p$ around the puncture. 
  Then $p'=gp$, where $g \in \G$. 
  This way we get the first isomorphism in (\ref{cips}). 
  
  To get the second isomorphism in (\ref{cips}), note that to get the moduli space $ {\mathscr P}_{\G, {\rm D}^*_1}$  
  we just add to the data parametrised by ${\rm Loc}_{\G, {\rm D}^*_1}$ a flag $\B$ near the puncture,  invariant under the monodromy, given by $g$. 
Now note that $\widehat \G$ parametrises pairs $(g, \B)$ where $g \in \G$ and $\B\subset \G$ is a Borel subgroup  such that $g\B g^{-1} = \B$.

 ii) This is evident: both projections are obtained by taking the quotient by the action of the group $S_\G$.    
 
 iii)  Follows from    the construction and Theorem \ref{MTH}.     \end{proof}

   \subsubsection{Dressing transformations.} 
Let us define 
 a Poisson Lie group action, which we call the {\it dressing action}: 
 \be \la{dres2}
      {\rm Loc}_{\G, \square} \times {\rm Loc}_{\G, {\rm D}_1^*}  \lra{\rm Loc}_{\G, {\rm D}_1^*}. 
\ee  
First,  glue the colored decorated surfaces $\square$ and ${\rm D}_1^*$ as shown on Figure \ref{R20a}. Then  cut the obtained colored decorate surface    by an arc   $\gamma$ connecting  the ends of the left colored interval, and containing the puncture on the left. Finally,   restrict to  
  the   punctured disc ${\rm D}_1^*$, obtianed by the cut. The   induced map of the moduli spaces  is the desired map (\ref{dres2}).

  \begin{figure}[ht]
\begin{tikzpicture}[scale=1.8]
\draw (1.7,1) --(2.3,1);
\draw (1.7,0) --(2.3,0);
\draw[dashed, red, thick, -stealth]  (1.7,0) --(1.7,.97);
\draw[dashed, red, thick, stealth-]  (2.3,0.03) --(2.3,1);
\node at (1.7, 1) {$\bullet$};
\node at (2.3,1) {$\bullet$};
\node at (1.7, 0) {$\bullet$};
\node at (2.3,0) {$\bullet$};

\begin{scope}[xshift=1.5cm, yshift=.05cm]
\draw[dashed, red, thick, stealth-]  (1.97,1) arc (110:250:0.6);
\draw (2.03,1) arc (70:-70:0.6);
\node[blue] at (2, .48) {$\circ$};
\node at (2, -.12) {$\bullet$};
\node at (2, 1) {$\bullet$};
\end{scope}

\draw[directed] (4.3, .5) -- (4.8, .5);
\node[red] at (4.55, .6) {glue};

\begin{scope}[xshift=3.5cm]
\draw (1.7,1) --(2.3,1);
\draw (1.7,0) --(2.3,0);
\draw (2.3,1) arc (90:-87:0.5);
\draw[dashed, red, thick, -stealth]  (1.7,0) --(1.7,.97);
\node[blue] at (2.0, 0.5) {$\circ$};
\node at (1.7, 1) {$\bullet$};
 \node at (2.3,1) {$\bullet$};
\node at (1.7, 0) {$\bullet$};
 \node at (2.3,0) {$\bullet$};
 \draw[blue] (1.7, 1) arc (90:-90:0.5);
\node[blue] at (2.3, 0.5) {${\tiny \gamma }$};
\end{scope}

\draw[directed] (6.8, .5) -- (7.3, .5);
\node[blue] at (7.05, .6) {cut $\gamma$};

\begin{scope}[xshift=6.1cm, yshift=.05cm]
\draw[dashed, red, thick, stealth-]  (1.97,1) arc (110:250:0.6);
\draw (2.03,1) arc (70:-70:0.6);
\node[blue] at (2, .48) {$\circ$};
\node at (2, -.12) {$\bullet$};
\node at (2, 1) {$\bullet$};
\end{scope}

\end{tikzpicture}
\caption{The   action  ${\rm Loc}_{\G, \square} \times {\rm Loc}_{\G, {\rm D}_1^*} \lra {\rm Loc}_{\G, {\rm D}_1^*}$.}
\label{R20a}
\end{figure}

                 \bl   The canonical  isomorphisms  of Poisson moduli spaces 
 $$
 {\mathscr L}_{\G, \square}  = \G, \qquad {\mathscr L}_{\G, {\rm D}^*_1}  = (\G, \pi) 
 $$  
provide a commutative   diagram of Poisson spaces, where on the bottom  the group $\G$ on itself by   conjugation: 
                   \begin{displaymath}
    \xymatrix{
         {\rm Loc}_{\G, \square} \times {\rm Loc}_{\G, {\rm D}_1^*}  \ar[r]^{\ \ \ \ ~~~~{\rm dressing}} \ar[d]^{=}&{\rm Loc}_{\G, {\rm D}_1^*} \ar[d]^{=}\\
          \G \times (\G, \pi)   \ar[r]^{\ \ {\rm conjugation}}  &  (\G, \pi)  }
                \end{displaymath}
 In particular the bottom action is  Poisson Lie group action of $\G$ on $(\G, \pi)$:
 \be \la{dressc}
 \G \times (\G, \pi) \lra (\G, \pi).
  \ee

  \el
 
 \begin{proof} 
  It is clear from Figure \ref{R20}     that the action (\ref{dres2})   amounts to the conjugation in the group $\G$.     
  The rest is also evident from Figure \ref{R20}. 
  \end{proof}    
    
                  \begin{figure}[ht]
\begin{tikzpicture}[scale=1.8]
\draw (1.7,1) --(2.3,1);
\draw (1.7,0) --(2.3,0);
\draw[dashed, red, thick, -stealth]  (1.7,0) --(1.7,.97);
\draw[dashed, red, thick, stealth-]  (2.3,0.03) --(2.3,1);
\node at (1.7, 1) {$\bullet$};
\node at (2.3,1) {$\bullet$};
\node at (1.7, 0) {$\bullet$};
\node at (2.3,0) {$\bullet$};
\draw[thick, black!40!green, -stealth] (1.7, 0.5) --(2.3,0.5);
\node[red] at (1.6, 0.55) {$p$};
\node[red] at (2.4, 0.55) {$q$};
\node at (2, 0.6) {$h$};

\begin{scope}[xshift=1.5cm, yshift=.05cm]
\draw[dashed, red, thick, stealth-]  (1.97,1) arc (110:250:0.6);
\draw (2.03,1) arc (70:-70:0.6);
\node[blue] at (2, .48) {$\circ$};
\node at (2, -.12) {$\bullet$};
\node at (2, 1) {$\bullet$};
\node[red] at (1.4, .5) {$r$};
\draw[thick, black!40!green, -stealth] (1.6,.51) to[out=45, in=180] (1.9, 0.65) to[out=0, in=120] (2.12, 0.55) to[out=-60, in=60] (2.12, .41) to[out=-120, in=0] (1.9, 0.31) to[out=180, in=-45] (1.6, .45);
\node at (2.2, 0.5) {${\tiny g}$};
\end{scope}

\draw[directed] (4.3, .5) -- (4.8, .5);
\node[red] at (4.55, .6) {glue};

\begin{scope}[xshift=3.5cm]
\draw (1.7,1) --(2.3,1);
\draw (1.7,0) --(2.3,0);
\draw (2.3,1) arc (90:-87:0.5);
\draw[dashed, red, thick, -stealth]  (1.7,0) --(1.7,.97);
\draw[dashed, red, thick, stealth-]  (2.3,0.03) --(2.3,1);
\node[blue] at (2.4, 0.5) {$\circ$};
\node at (1.7, 1) {$\bullet$};
\node at (2.3,1) {$\bullet$};
\node at (1.7, 0) {$\bullet$};
\node at (2.3,0) {$\bullet$};
\node[red] at (1.6, 0.5) {$p$};
\draw[black!40!green, -stealth, thick] (2.3, 0.7) arc (90:-90:0.2);
\draw[black!40!green, -stealth, thick] (1.7, 0.7) -- (2.29, 0.7);
\draw[black!40!green, stealth-, thick] (1.7, 0.3) -- (2.29, 0.3);
\node at (2, 0.8) {$h$};
\node at (2, 0.2) {$h^{-1}$};
\node at (2.6, 0.5) {$g$};
\end{scope}
\end{tikzpicture}
\caption{The   action  ${\rm Loc}_{\G, \square} \times {\rm Loc}_{\G, {\rm D}_1^*} \lra {\rm Loc}_{\G, {\rm D}_1^*}$  amounts to  conjugation $(h, g) \lms hgh^{-1}$.}
\label{R20}
\end{figure} 
     
  Similarly,  
  there is 
 a partially defined  Poisson Lie group   {\it dressing action} of the Poisson Lie group $\G^*$: 
 \be \la{dres22}
      {\rm Loc}_{\G, \square}/S_\G \times {\rm Loc}_{\G, \odot}  \lra{\rm Loc}_{\G, \square}/S_\G. 
\ee  
Namely, we glue the colored decorated surfaces $\square$ and $\odot$. Then  cut the obtained   surface    by the arc   $\gamma$ connecting  the ends of the left vertical boundary interval, and containing the puncture, and  restrict to  
  the   square. The   induced map of the moduli spaces  is the desired map (\ref{dres22}).

  \begin{figure}[ht]
\begin{tikzpicture}[scale=1.8]
\draw (1.7,1) --(2.3,1);
\draw (1.7,0) --(2.3,0);
\draw[dashed, red, thick, -stealth]   (1.7,0) --(1.7,.97);
\draw[dashed, red, thick, stealth-]  (2.3,0.03) --(2.3,1);
\node at (1.7, 1) {$\bullet$};
\node at (2.3,1) {$\bullet$};
\node at (1.7, 0) {$\bullet$};
\node at (2.3,0) {$\bullet$};

\begin{scope}[xshift=1.5cm, yshift=.05cm]
\draw[dashed, red, thick, stealth-]  (1.97,1) arc (110:250:0.6);
\draw[dashed, red, thick, -stealth]  (2.03,1) arc (70:-70:0.6);
\node[blue] at (2, .48) {$\circ$};
\node at (2, -.12) {$\bullet$};
\node at (2, 1) {$\bullet$};
\end{scope}

\draw[directed] (4.3, .5) -- (4.8, .5);
\node[red] at (4.55, .6) {glue};

\begin{scope}[xshift=3.5cm]
\draw (1.7,1) --(2.3,1);
\draw (1.7,0) --(2.3,0);
\draw[red, dashed, thick, -stealth] (2.3,1) arc (90:-87:0.5);
\draw[red, dashed, thick, -stealth] (1.7,0) --(1.7,.97);
\node[blue, dashed] at (2.0, 0.5) {$\circ$};
\node at (1.7, 1) {$\bullet$};
 \node at (2.3,1) {$\bullet$};
\node at (1.7, 0) {$\bullet$};
 \node at (2.3,0) {$\bullet$};
 \draw[blue] (1.7, 1) arc (90:-90:0.5);
\node at (2.3, 0.5) {${\tiny  \gamma }$};
\end{scope}

\draw[directed] (6.8, .5) -- (7.3, .5);
\node[blue] at (7.05, .6) {cut $\gamma$};

\begin{scope}[xshift=6.1cm, yshift=.05cm]
\draw (1.7,1) --(2.3,1);
\draw (1.7,0) --(2.3,0);
\draw[red, dashed, thick, -stealth] (1.7,0) --(1.7,.97);
 \draw[dashed, red, thick, stealth-]  (2.3,0.03) --(2.3,1);
\node at (1.7, 1) {$\bullet$};
 \node at (2.3,1) {$\bullet$};
\node at (1.7, 0) {$\bullet$};
 \node at (2.3,0) {$\bullet$};
\end{scope}

\end{tikzpicture}
\caption{The   action  ${\rm Loc}_{\G, \square}/S_\G \times {\rm Loc}_{\G, \odot} \lra {\rm Loc}_{\G, \square}/S_\G$.}
\label{R20b}
\end{figure}

{\bf Remark}. The classical dressing transformations \cite{STS} lead to partially defined and $1:2^r$ multi-valued ``actions" $a: \G \times \G^* \lra \G^*$ and $b: \G^* \times \G \lra \G$. 
 The issues with them are clear from Figures \ref{R20a} and \ref{R20b}. The  map  $a$ is not well defined since after cutting along the arc $\gamma$ we can not reconstruct the  pinning, and the pair of flags assigned to the arc $\gamma$ may not be generic. Same for the map $b$. There is a well defined analog of $a$  -  the dressing action (\ref{dressc}). 
 It would be interesting to find a compactification of the space ${\rm Loc}_{\G, \square}/S_\G$ where the dressing action really acts.

  Yet the partially defined multi-valued dressing transformations exist in a very general set-up. 

 \medskip
 
 \subsection{The Heisenberg double $\D_{\rm H}(\G)$ of   $\G$ via moduli spaces} \la{SECT4.2d} 
 
 \medskip
 
 Let us define a symplectic moduli space $\D_{\rm H}(\G)$, which we call the {\it Heisenberg double} of the Poisson Lie group $\G$. The Heisenberg double is not  a Poisson Lie group.
It is  closely related but   different than the symplectic space considered   in \cite{STS92}. 
The name   is justified by Conjecture \ref{HD} below. \vskip 1mm

Denote by ${\rm C}_{2,1}$ the colored decorated surface given by a 
cylinder with a pair of special points on each of the boundary circles,  two colored intervals on the left   circle, and one colored interval  on the right one. 
  Equivalently,  it is an annulus   
 with two pairs of special points,    two colored  boundary intervals  outside, and one    inside.  \vskip 1mm
 
  The  Poisson moduli space   ${\rm Loc}_{\G, {\rm C}_{2,1}}$  
  parametrizes $\G-$local systems on the cylinder with a pinning  at each of the three colored intervals. There are no conditions for the pair of flags over the non-colored  interval. 
 
  The  Poisson moduli space   ${\mathscr L}_{\G, {\rm C}_{2,1}}$  is the fiber of the map given by the outer monodromy around the boundary circle which carries two pinnings:
  \be
 \begin{split}
 &\mu_{\rm out}: {\rm Loc}_{\G, {\rm C}_{2,1}}\lra \rm H,\\
 & {\mathscr L}_{\G, {\rm C}_{2,1}}:= \mu_{\rm out}^{-1}(e).
  \end{split}
    \ee

\bt There are 
 actions of the Poisson Lie groups   ${\mathscr L}_{\G, \odot}$ and   ${\rm Loc}_{\G, \square}$ on   ${\mathscr L}_{\G, {\rm C}_{2,1}}$:
\be \la{abproj1}
 \begin{split}
  &\Phi: {\mathscr L}_{\G, \odot} \times  {\mathscr L}_{\G, {\rm C}_{2,1}} \lra    {\mathscr L}_{\G, {\rm C}_{2,1}},\\  
 &\Psi: {\mathscr L}_{\G, {\rm C}_{2,1}} \times {\rm Loc}_{\G, \square} \lra {\mathscr L}_{\G, {\rm C}_{2,1}}.  \\
 \end{split}
  \ee 
  \et

  \begin{proof} To define the action $\Phi$,     glue  the  surface $\odot$ 
   to the  annulus ${\rm C}_{2,1}$,  so that the only non-colored interval on the annulus is inside, and  cut     the puncture by an arc  ending at the    special points on the non-colored boundary interval, see Figure \ref{R17}.  
        \begin{figure}[ht]
\begin{tikzpicture}[scale=2]
\draw[dashed, red, thick, stealth-]  (1.97,1) arc (120:240:0.6);
\draw[dashed, red, thick, -stealth]  (2.03,1) arc (60:-60:0.6);
\node[blue] at (2, 0.5) {$\circ$};
\node at (2, -0.05) {$\bullet$};
\node at (2,1) {$\bullet$};

\begin{scope}[xshift=1cm]
\draw[dashed, red, thick, stealth-]  (1.97,1) arc (120:240:0.6);
\draw[dashed, red, thick, stealth-] (2.02,0.75) arc (60:-60:.28);
\draw[ thick] (1.98,0.75) arc (120:240:.28);
\draw[dashed, red, thick, -stealth]  (2.03,1) arc (60:-60:0.6);
\node at (2, -0.05) {$\bullet$};
\node at (2,1) {$\bullet$};
\node at (2, .75) {$\bullet$};
\node at (2, .25) {$\bullet$};
\end{scope}
\draw[directed] (3.6,0.5)--(4.2,0.5);
\node[red] at (3.9, 0.6) {glue};
\draw[directed] (5.5,0.5)--(6.3,0.5);
\node[blue] at (5.9,0.6) {encircle};
\node[blue] at (5.9,0.4) {puncture};
\draw[directed] (7.6,0.5)--(8.2,0.5);
\node at (7.9,0.6) {cut};

\begin{scope}[xshift=3cm]
\draw[dashed, red, thick, stealth-]  (1.97,1) arc (90:270:0.52);
\draw[dashed, red, thick, stealth-] (2.02,0.75) arc (60:-60:.28);
\draw[thick] (1.98,0.75) arc (120:240:.28);
\draw[dashed, red, thick, -stealth] (2.03,1) arc (60:-60:0.6);
\node at (2, -0.05) {$\bullet$};
\node at (2,1) {$\bullet$};
\node at (2, .75) {$\bullet$};
\node at (2, .25) {$\bullet$};
\node[blue] at (1.75, .5) {$\circ$};
\end{scope}

\begin{scope}[xshift=5cm]
\draw[dashed, red, thick, stealth-]  (1.97,1) arc (90:270:0.52);
\draw[dashed, red, thick, stealth-] (2.02,0.75) arc (60:-60:.28);
\draw[thick] (1.98,0.75) arc (120:240:.28);
\draw[red, dashed, thick, -stealth] (2.03,1) arc (60:-60:0.6);
\node at (2, -0.05) {$\bullet$};
\node at (2,1) {$\bullet$};
\node at (2, .75) {$\bullet$};
\node at (2, .25) {$\bullet$};
\node[blue] at (1.75, .5) {$\circ$};
\draw[blue, thick] (1.97,0.75) arc (80:280:.255);
\end{scope}

\begin{scope}[xshift=7cm]
\draw[dashed, red, thick, stealth-]  (1.97,1) arc (120:240:0.6);
\draw[dashed, red, thick, stealth-] (2.02,0.75) arc (60:-60:.28);
\draw[thick] (1.98,0.75) arc (120:240:.28);
\draw[red, dashed, thick, -stealth] (2.03,1) arc (60:-60:0.6);
\node at (2, -0.05) {$\bullet$};
\node at (2,1) {$\bullet$};
\node at (2, .75) {$\bullet$};
\node at (2, .25) {$\bullet$};
\end{scope}
\end{tikzpicture}
\caption{Gluing   from the left the decorated surface $\odot$ to the one ${\rm C}_{2,1}$, and   cutting out  the puncture.}
\label{R17}
\end{figure}     

    To define the 
    action $\Psi$,   glue  the    rectangle $\square$   to the   annulus  ${\rm C}_{2,1}$   
    along the only colored boundary interval there,  see Figure \ref{R18}, encircle the hole   
     by an arc ending at the external special points, and cut   the arc.
  
    \begin{figure}[ht]

\begin{tikzpicture}[scale=2]
\draw[thick]  (1.97,1) arc (120:240:0.6);
\draw[dashed, red, thick, -stealth] (2.02,0.75) arc (60:-60:.28);
\draw[dashed, red, thick, stealth-] (1.98,0.75) arc (120:240:.28);
\draw[dashed, red, thick, -stealth]  (2.03,1) arc (60:-60:0.6);
\node at (2, -0.05) {$\bullet$};
\node at (2,1) {$\bullet$};
\node at (2, .75) {$\bullet$};
\node at (2, .25) {$\bullet$};
\draw[directed] (3.6,0.5)--(4.2,0.5);
\node[red] at (3.9, 0.6) {glue};
\draw[directed] (5.5,0.5)--(6.3,0.5);
\node[blue] at (5.9,0.6) {encircle};
\draw[directed] (7.6,0.5)--(8.2,0.5);
\node at (7.9,0.6) {cut};

\begin{scope}[xshift=1cm]
\draw[thick] (1.7,1) --(2.3,1);
\draw[thick] (1.7,0) --(2.3,0);
\draw[dashed, red, thick, -stealth]  (1.7,0) --(1.7,.97);
\draw[dashed, red, thick, stealth-]  (2.3,0.03) --(2.3,1);
\node at (1.7, 1) {$\bullet$};
\node at (2.3,1) {$\bullet$};
\node at (1.7, 0) {$\bullet$};
\node at (2.3,0) {$\bullet$};
\end{scope}

\begin{scope}[xshift=3cm]
\draw[thick] (1.4, 1) -- (2,1);
\draw[thick] (1.4, -.05) -- (2,-.05);
\begin{scope}[xshift=-.35cm]
\draw[dashed, red, thick, -stealth] (2.02,0.75) arc (60:-60:.28);
\draw[dashed, red, thick, stealth-] (1.98,0.75) arc (120:240:.28);
\node at (2, .75) {$\bullet$};
\node at (2,0.25) {$\bullet$};
\end{scope}
\draw[dashed, red, thick, -stealth] (2.03,1) arc (60:-60:0.6);
 \node at (2, -0.05) {$\bullet$};
 \node at (2,1) {$\bullet$};
\node at (1.4, -0.05) {$\bullet$};
\node at (1.4,1) {$\bullet$};
\draw[thick] (1.4, 0) -- (1.4,0.97);
\end{scope}

\begin{scope}[xshift=5cm]
\draw[blue, thick]  (2, 0.0) arc (90:-90:-0.5);
\draw[thick] (1.4, 1) -- (2,1);
\draw[thick] (1.4, -.05) -- (2,-.05);
\draw[dashed, red, thick, -stealth] (2.02,0.75) arc (60:-60:.28);
\draw[dashed, red, thick, stealth-] (1.98,0.75) arc (120:240:.28);
\node at (2, .75) {$\bullet$};
\node at (2, .25) {$\bullet$};
\draw[dashed, red, thick, -stealth] (2.03,1) arc (60:-60:0.6);
\node at (2, -0.05) {$\bullet$};
\node at (2,1) {$\bullet$};
\node at (1.4, -0.05) {$\bullet$};
\node at (1.4,1) {$\bullet$};
\draw[thick] (1.4, 0) -- (1.4,0.97);
\end{scope}

\begin{scope}[xshift=6.8cm]
\draw[thick]  (1.97,1) arc (120:240:0.6);
\draw[dashed, red, thick, -stealth] (2.02,0.75) arc (60:-60:.28);
\draw[dashed, red, thick, stealth-] (1.98,0.75) arc (120:240:.28);
\draw[dashed, red, thick, -stealth]  (2.03,1) arc (60:-60:0.6);
\node at (2, -0.05) {$\bullet$};
\node at (2,1) {$\bullet$};
\node at (2, .75) {$\bullet$};
\node at (2, .25) {$\bullet$};
\end{scope}
\end{tikzpicture}
\caption{Gluing   from the right the  rectangle $\square$ to the surface ${\rm C}_{2,1}$, and   cutting.}
\label{R18}
\end{figure}

     
     Since the gluing and cutting maps are Poisson,   these are   Poisson group actions. 
    \end{proof} 
    
  
\bp 
There are  surjective Poisson projections: 
 \be  \la{pqprojb} 
 \begin{split}
  &{\rm P}:  {\mathscr L}_{\G, {\rm C}_{2,1}}\lra {\rm Loc}_{\G, \square}.\\
  &{\rm Q}:  {\mathscr L}_{\G, {\rm C}_{2,1}}\lra   {\mathscr L}_{\G, \odot}.\\
 \end{split}
 \ee
 \ep 
 
  \begin{proof} The   map  ${\rm P}$ is induced by  the restriction to a rectangle $R\subset{\rm C}_{2,1}$ with the vertices at the special points  and two colored sides,       shown on the left  of Figure \ref{R18a}. 

 \begin{figure}[ht]

\begin{tikzpicture}[scale=2]
\draw[draw=blue!.1, fill=blue, fill opacity=0.2] (2, 0.75) -- (2, 1) -- (2.03, 1) arc (60:-60:0.6) -- (2,0)--(2,.25)-- (2.02, 0.25)  arc (-60:60:.28) -- cycle;
\draw[dashed, red, thick, stealth-]  (1.97,1) arc (120:240:0.6);
\draw[dashed, red, thick, -stealth] (2.02,0.75) arc (60:-60:.28);
\draw[thick] (1.98,0.75) arc (120:240:.28);
\draw[dashed, red, thick, -stealth] (2.03,1) arc (60:-60:0.6);
\node at (2, -0.05) {$\bullet$};
\node at (2,1) {$\bullet$};
\node at (2, .75) {$\bullet$};
\node at (2, .25) {$\bullet$};
 \draw[directed] (2.5,0.5) -- (3, 0.5);
 \node at (2.75, 0.6) {cut};
 
 \begin{scope}[xshift=1.5cm]
 \draw[draw=blue!.1, fill=blue, fill opacity=0.2]  (1.7,1)--(2.3,1)--(2.3,0)--(1.7,0)--cycle;
\draw (1.7,1) --(2.3,1);
\draw (1.7,0) --(2.3,0);
\draw[dashed, red, thick, -stealth]  (1.7,0) --(1.7,.97);
\draw[dashed, red, thick, stealth-]  (2.3,0.03) --(2.3,1);
\node at (1.7, 1) {$\bullet$};
\node at (2.3,1) {$\bullet$};
\node at (1.7, 0) {$\bullet$};
\node at (2.3,0) {$\bullet$};
 \end{scope}
 
\begin{scope}[xshift=4cm]
 \draw[dashed, red, thick, stealth-]  (1.97,1) arc (120:240:0.6);
\draw[dashed, red, thick, -stealth] (2.02,0.75) arc (60:-60:.28);
\draw[thick] (1.98,0.75) arc (120:240:.28);
\draw[dashed, red, thick, -stealth] (2.03,1) arc (60:-60:0.6);
\draw[blue] (2, .5) ellipse (0.25cm and 0.35cm);
\node at (2, -0.05) {$\bullet$};
\node at (2,1) {$\bullet$};
\node at (2, .75) {$\bullet$};
\node at (2, .25) {$\bullet$};
 \draw[directed] (2.5,0.5) -- (3, 0.5);
 \node at (2.75, 0.6) {encircle};
\end{scope}

\begin{scope}[xshift=5.5cm]
 \draw[dashed, red, thick, stealth-]  (1.97,1) arc (120:240:0.6);
\draw[dashed, red, thick, -stealth] (2.03,1) arc (60:-60:0.6);
\node[blue] at (2, .5) {$\circ$};
\node at (2, -0.05) {$\bullet$};
\node at (2,1) {$\bullet$};
\end{scope}
\end{tikzpicture}

\caption{On the left: the rectangle $R$ on the  surface ${\rm C}_{2,1}$. On the right: the cutting loop.}
\label{R18a}
\end{figure}

 The   map  $\rm Q$ is induced by   cutting the  ${\rm C}_{2,1}$ by a simple loop, and restricting to the internal component,  shown on the right  of Figure \ref{R18a}.  
 Since  restriction maps  are  Poisson, these are Poisson maps.  
  \end{proof}  
 
  Composing  ${\rm P}$ with   the Dehn twist ${\cal D}$, rotating  one of the  boundary circle by $2\pi$, we get a map ${\rm P}':= {\cal D}\circ {\rm P}$. \vskip 1mm
   
The maps (\ref{pqprojb}) provide surjective projections of Poisson spaces
\be \la{abproj2}
 \begin{split}
&p: {\mathscr L}_{\G, \odot} \backslash {\mathscr L}_{\G, {\rm C}_{2,1}} \lra    {\rm Loc}_{\G, \square};   \\
&  q: {\mathscr L}_{\G, {\rm C}_{2,1}} / {\rm Loc}_{\G, \square} \lra 
{\mathscr L}_{\G, \odot}.  \\
  \end{split}
    \ee 


 \bd The Heisenberg double $\D_{\rm H}(\G)$ of the Poisson Lie group $\G$  is a symplectic moduli space 
 \be
 \begin{split}
 \D_{\rm H}(\G):=  {\mathscr L}_{\G, {\rm C}_{2,1}}.  
  \end{split}
   \ee
   \ed

 The Heisenberg double $\D_\H(\G)$  has the following 
 features.  
 
\begin{enumerate}

\vskip 1mm\item  The actions of the Poisson Lie groups $\G$ and $\G^*$, provided by the actions (\ref{abproj1}):
  \be \la{pqproj1}
 \begin{split}
  &\varphi: \G^* \times \D_{\rm H}(\G)\lra \D_{\rm H}(\G);\\
 & \psi: \D_{\rm H}(\G) \times \G\lra  \D_{\rm H}(\G).\\ 
 \end{split}
 \ee 
 
\vskip 1mm\item  The  surjective projections of  Poisson spaces, provided  by (\ref{abproj2}):
   \be \la{pqproj2}
 \begin{split}
  &\rm P:    \D_{\rm H}(\G)\lra \G;\\ 
   & \rm Q:   \D_{\rm H}(\G)   \lra  \G^*.\\
 \end{split}
 \ee 
  
 \vskip 1mm \item The  induced   surjective projections of  Poisson spaces 
   \be \la{pqproj2}
 \begin{split}
  &p:  \G^* \backslash \D_{\rm H}(\G)\lra \G;\\ 
   & q:   \D_{\rm H}(\G) / \G \lra  \G^*.\\
 \end{split}
 \ee 

\end{enumerate}

\subsubsection{Heisenberg double of the Hopf algebra ${\cal O}_q(\G)$ and ${\cal O}_q(\G^*)$.} Let  ${\cal H}$ and ${\cal H}^*$ be a pair of Hopf algebras with the coproducts
 $\Delta_{{\cal H}}$ and $\Delta_{{\cal H}^*}$. Let us assume that there is a   Hopf pairing
\be
 \begin{split}
 &\varphi (\cdot, \cdot ): {\cal H} \otimes_\C {\cal H}^* \lra \C, \\
 &\varphi ( xy, u) = \varphi ( x\otimes y, \Delta_{{\cal H}^*}(u)), \\
 & \varphi ( x, uv) = \varphi( \Delta_{{\cal H}}(x), u\otimes v). \\
 \end{split}
 \ee
 Then the vector   space ${\cal H}^* \otimes_\C {\cal H}$ carries a canonical structure of an associative algebra, such that ${\cal H}^* \otimes 1$ and $1 \otimes {\cal H}$ are subalgebras identified with ${\cal H}^*$ and ${\cal H}$. 
 Namely, given $x \in {\cal H}, u \in {\cal H}^*$, we set, following \cite{K97}:
 \be
 x \cdot u:= ({\rm Id} \otimes \varphi \otimes {\rm Id} )(\Delta_{{\cal H}^*}(u) \otimes \Delta_{{\cal H}}(x)). 
 \ee
 Here the map  $ {\rm Id} \otimes \varphi \otimes {\rm Id} $  is induced by applying the pairing to the two middle factors:
  \be
 {\rm Id} \otimes \varphi \otimes {\rm Id} : {\cal H}^* \otimes {\cal H}^* \otimes {\cal H} \otimes H \lra {\cal H}^* \otimes  {\cal H}  
   \ee

   The algebra ${\cal H}^* \otimes_\C {\cal H}$ is   called the Heisenberg double for the Hopf pairing  $\varphi$.  
 It was defined   in  \cite{RS90},  \cite{AF91}, \cite{STS92},  \cite{Nov92}    in the set up when ${\cal H}^*$ is the dual Hopf algebra to ${\cal H}$. \vskip 2mm

\bl The   projections $\rm P$ and $\rm Q$ in (\ref{pqproj2}) induce maps of  quantum algebras:
\be
 \begin{split}
 &{\rm P}^*: {\cal O}_q(\G) \lra {\cal O}_q(\D_{\rm H}(\G));\\
&{\rm Q}^*: {\U}_q(\mathfrak{g}) \lra {\cal O}_q(\D_{\rm H}(\G)).\\
 \end{split}
 \ee
 \el 
 
 \begin{proof} The map $\rm P$ is a cluster Poisson projection, and thus induces a map ${\rm P}^*$ of the q-deformed algebras. 
 
 The map ${\rm Q}^*$ is defined   on the generators ${\bf E}_i, {\bf F}_i, {\bf K}_i \in {\cal O}_q(\G^*) $  generating ${\U}_q(\mathfrak{g})$. 
 Equivalently, the map  ${\rm Q}^*$ coincides with the following canonical injective map, assigned by Theorem \ref{UEAB} to the two special points on the relevant boundary component of the decorated surface ${\rm C}_{2,1}$:
 $$
 \kappa:  {\U}_q(\mathfrak{g})  \hra {\cal O}_q({\cal L}_{\G, {\rm C}_{2,1}}).
  $$
   \end{proof}
  
 Therefore we get a map of vector spaces, induced by the product of the images of the maps ${\rm P}^*$ and ${\rm Q}^*$:
 \be \la{0987}
 {\rm P}^*\cdot {\rm Q}^* : {\cal O}_q(\G)\otimes_\C  {\U}_q(\mathfrak{g}) \lra {\cal O}_q(\D_{\rm H}(\G)).
 \ee

  \bcon \la{HD}
 i) There is  the canonical Hopf pairing between the Hopf algebras ${\cal O}_q(\G)$ and ${\U}_q(\mathfrak{g})$:
 \be \la{987654}
 \langle\cdot, \cdot \rangle: {\cal O}_q(\G)\otimes_\C  {\U}_q(\mathfrak{g})  \lra \C. 
 \ee
 
 ii) The corresponding map (\ref{0987}) is an isomorphism of vector spaces. 

  The algebra ${\cal O}_q(\D_{\rm H}(\G))$ for generic $q$ is   isomorphic to the Heisenberg double  for the Hopf pairing   (\ref{987654}).  
 \econ  
  
For the traditional quantum group ${\cal O}_q(\G)$ the pairing (\ref{987654}) is well known.  We use above the cluster Hopf algebra ${\cal O}_q(\G)$, which is conjecturally isomorphic to the traditional quantum group ${\cal O}_q(\G)$.   
 
\medskip

\subsection{The dual Poisson Lie  group $\D(\G)^*$ for the Drinfeld double $\D(\G)$.} \la{SECT4.2c}

\medskip

Denote by  $\odot_n$ the decorated   disc  with $n$  punctures and two  special   points, see Figure \ref{R11}. We assume  the punctures are numbered. It is convenient to think that the punctures are    on the axis connecting    special points. Then there are the   cluster Poisson moduli space  ${\cal P}_{\G, \odot_n}$, and the Poisson moduli space  
 ${\rm Loc}_{\G, \odot_n}$.

\begin{figure}[ht]
\begin{tikzpicture}[scale=2]
\draw[dashed, red, thick, stealth-]  (1.97,1) arc (120:240:0.6);
\draw[dashed, red, thick, -stealth]  (2.03,1) arc (60:-60:0.6);
\node[blue] at (2, 0.7) {$\circ$};
\node[blue] at (2, 0.3) {$\circ$};
\node[red] at (2, -0.05) {$\bullet$};
\node[red] at (2,1) {$\bullet$};
\end{tikzpicture}
\caption{The twice punctured disc $\odot_2$ with two special points, and two colored sides.}
\label{R11}
\end{figure}

\bp The moduli space ${\rm Loc}_{\G, \odot_n}$ is a Poisson Lie group.
 with the product map
\be \la{CoMPL}
{\rm Loc}_{\G, \odot_n} \times {\rm Loc}_{\G, \odot_n}\lra {\rm Loc}_{\G, \odot_n},
\ee
provided by   gluing the boundary intervals, and   encircling each pair of the punctures, see Figure \ref{R8a}. 
\ep 

\begin{proof}
\begin{figure}[ht]
\begin{tikzpicture}[scale=2]
\draw[dashed, red, thick, stealth-]  (1.97,1) arc (120:240:0.6);
\draw[dashed, red, thick, -stealth]  (2.03,1) arc (60:-60:0.6);
\node[blue] at (2, 0.7) {$\circ$};
\node[blue] at (2, 0.3) {$\circ$};
\node[red] at (2, -0.05) {$\bullet$};
\node[red] at (2,1) {$\bullet$};

\begin{scope}[xshift=1.5cm]
\draw[dashed, red, thick, stealth-]  (1.97,1) arc (120:240:0.6);
\draw[dashed, red, thick, -stealth]  (2.03,1) arc (60:-60:0.6);
\node[blue] at (2, 0.7) {$\circ$};
\node[blue] at (2, 0.3) {$\circ$};
\node[red] at (2, -0.05) {$\bullet$};
\node[red] at (2,1) {$\bullet$};
\end{scope}

\draw[directed,thick] (4.2,0.5)--(4.8,0.5);
\node[red] at (4.5, 0.7) {glue};
\draw[directed, thick] (6.2,0.5)--(6.8,0.5);
\node[blue] at (6.5, 0.7) {encircle};

\begin{scope}[xshift=3.5cm]
\draw[dashed, red, thick, stealth-]  (1.97,1) arc (120:240:0.6);
\draw[dashed, red, thick, -stealth]  (2.03,1) arc (60:-60:0.6);
\node[blue] at (1.9, 0.7) {$\circ$};
\node[blue] at (2.1, 0.7) {$\circ$};
\node[blue] at (1.9, 0.3) {$\circ$};
\node[blue] at (2.1, 0.3) {$\circ$};
\node[red] at (2, -0.05) {$\bullet$};
\draw[dashed, blue] (2, 0.7) ellipse (0.2cm and 0.1cm);
\draw[dashed, blue] (2, 0.3) ellipse (0.2cm and 0.1cm);
\node[red] at (2,1) {$\bullet$};
\end{scope}

\begin{scope}[xshift=5.5cm]
\draw[dashed, red, thick, stealth-]  (1.97,1) arc (120:240:0.6);
\draw[dashed, red, thick, -stealth]  (2.03,1) arc (60:-60:0.6);
\node[blue] at (2, 0.7) {$\circ$};
\node[blue] at (2, 0.3) {$\circ$};
\node[red] at (2, -0.05) {$\bullet$};
\node[red] at (2,1) {$\bullet$};
\end{scope}
\end{tikzpicture}
\caption{The composition law is induced by the map $\odot_2\times \odot_2 \lra \odot_2$.}
\label{R8a}
\end{figure}

The unit is   the trivial $\G-$local system with the   pinnings $(p,p)$. 

The inverse map is given by the rotation of the surface $\odot_n$ along the central vertical axis. \end{proof}

There is the outer monodromy map, defined just as for the decorated surface $\odot$: 
\be
\mu_{\rm out}: {\rm Loc}_{\G, \odot_n} \lra \H.
\ee
We denote by ${\mathscr L}_{\G, \odot_n}$ its fiber over the unit $e \in \H$: 
\be
{\mathscr L}_{\G, \odot_n} := \mu_{\rm out}^{-1}(e).
\ee

Just like in Section \ref{thegr}, the composition law (\ref{CoMPL}) induces a product map
\be \la{CoMPLa}
{\mathscr L}_{\G, \odot_n} \times {\mathscr L}_{\G, \odot_n}\lra {\mathscr L}_{\G, \odot_n}.
\ee

\bl There are group isomorphisms
\be \la{CoMPLb}
\begin{split}
&{\rm Loc}_{\G, \odot_n}  \stackrel{\sim}{\lra} \B \times \G^{n-1} \times\B.\\
&{\mathscr L}_{\G, \odot_n} \stackrel{\sim}{\lra}  \G^{n-1} \times \G^*.\\
\end{split}
\ee
\el 

\begin{proof} For $n=2$, it is obtained by  transporting   the left pinning to the right one along the three paths $\alpha_1, \alpha, \alpha_2$  on Figure \ref{R13}. 
It assigns  to a $\G-$local system with pinnings $(p,q)$ the element $(b_1, g, b_2)\in \B \times \G \times \B$ such that 
\be
q= b_1\alpha_1(p), \quad q=g\alpha(p), \quad q=b_2\alpha_2(p).
\ee
 The map (\ref{CoMPLb}) induces the second group isomorphism.

\begin{figure}[ht]
\begin{tikzpicture}[scale=2.5]
\draw[dashed, red, thick, stealth-]  (1.97,1) arc (120:240:0.6);
\draw[dashed, red, thick, -stealth]  (2.03,1) arc (60:-60:0.6);
\node[blue] at (2, 0.7) {$\circ$};
\node[blue] at (2, 0.3) {$\circ$};
\node[red] at (2, -0.05) {$\bullet$};
\node[red] at (2,1) {$\bullet$};
\draw[ thick, black!40!green, -stealth] (1.8, 0.8)--(2.2, 0.8);
\draw[thick, black!40!green, -stealth] (1.7, 0.5)--(2.3, 0.5);  
\draw[ thick, black!40!green, -stealth] (1.8, 0.2)--(2.2, 0.2); 
\node  at (2, 0.87) {$b_1$};
\node  at (2, 0.57) {$g$};
\node  at (2, 0.12) {$b_2$};
\end{tikzpicture}
\caption{The group isomorphism  ${\mathscr L}_{\G, \odot_2} \stackrel{\sim}{\lra}  \G \times \G^*$.}
\label{R13}
\end{figure}
Alternatively, there are two canonical projections of Poisson Lie groups:
\be \la{MPLGP}
\begin{split}
&{\mathscr L}_{\G, \odot_n} \lra {\rm Loc}_{\G, \square} = \G^{n-1}, \\
&{\mathscr L}_{\G, \odot_n} \lra {\mathscr L}_{\G, \odot} = \G^*.\\
\end{split}
\ee
The first is given by cutting the $n-$punctured disc along two arcs, each containing a special point and  encircling the nearby puncture, and restricting the $\G-$local systems, see Figure \ref{R10} on the left.

The second is given by encircling   punctures, and restricting   $\G-$local systems, see Figure \ref{R10} on the right.   
\begin{figure}[ht]
\begin{tikzpicture}[scale=2]
\begin{scope}[xshift=-.2cm]
\draw[blue] (2,.2) ellipse (0.1cm and 0.25cm);
\draw[blue] (2,.77) ellipse (0.1cm and 0.25cm);
\draw[dashed, red, thick, stealth-]  (1.97,1) arc (120:240:0.6);
\draw[dashed, red, thick, -stealth]  (2.03,1) arc (60:-60:0.6);
\node[blue] at (2, 0.7) {$\circ$};
\node[blue] at (2, 0.3) {$\circ$};
\node[red] at (2, -0.05) {$\bullet$};
\node[red] at (2,1) {$\bullet$};
\end{scope}

\draw[directed, thick] (1.8, -.3) -- (1.8, -.7);  
\draw[directed, thick] (5, -.3) -- (5, -.7);  
\node[blue] at (1.6, -.5) {cut};
\node[blue] at (4.8, -.5) {cut};
\node at (0.9, -1.25) {G~=};
\node at (5.5, -1.25) {$={\rm G}^\ast$};

\begin{scope}[xshift=-.2cm]
\draw[red,dashed,-stealth,thick] (1.5, -1.5) -- (1.5, -1);
\draw[red,dashed,stealth-,thick] (2.5, -1.5) -- (2.5, -1);
\draw[blue, dashed,thick] (1.5, -1.5) -- (2.5, -1.5);
\draw[blue, dashed, thick] (1.5, -1) -- (2.5, -1);
\end{scope}

\begin{scope}[xshift=3cm]
\draw[blue] (2,.5) ellipse (0.15cm and 0.3cm);
\draw[dashed, red, thick, stealth-]  (1.97,1) arc (120:240:0.6);
\draw[dashed, red, thick, -stealth]  (2.03,1) arc (60:-60:0.6);
\node[blue] at (2, 0.7) {$\circ$};
\node[blue] at (2, 0.3) {$\circ$};
\node[red] at (2, -0.05) {$\bullet$};
\node[red] at (2,1) {$\bullet$};
\end{scope}

\begin{scope}[xshift=4cm, scale=.5, yshift=-3cm]
\draw[dashed, red, thick, stealth-]  (1.97,1) arc (120:240:0.6);
\draw[dashed, red, thick, -stealth]  (2.03,1) arc (60:-60:0.6);
\node[blue] at (2, 0.5) {$\circ$};
\node[red] at (2, -0.05) {$\bullet$};
\node[red] at (2,1) {$\bullet$};
\end{scope}
\end{tikzpicture}
\caption{Cutting $\odot_2$ in two different ways provides maps to Poisson Lie groups $\G$ and $\G^*$.}
\label{R10}
\end{figure}

These maps are group maps by the construction. They are Poisson   since induced by the restriction of  $\G-$local systems with extra data.  
So projections (\ref{MPLGP}) are    Poisson Lie group maps.  \end{proof}

Although the map 
$$
{\mathscr L}_{\G, \odot_n} \lra  \G^{n-1}\times \G^*
$$
 is a group map,  it is not Poisson. \\

The Poisson Lie groups  ${\rm Loc}_{\G, \odot_n}$ and ${\mathscr L}_{\G, \odot_n}$  
  are the Fock-Rosly polyubles \cite{FR}. 
  
  Our approach  shows that 
  they give rise to  quantum Hopf algebras
\be
{\cal O}_q({\mathscr L}_{\G, \odot_n}), \qquad {\cal O}_q({\rm Loc}_{\G, \odot_n}).
\ee
\vskip 1mm

\bt \la{THD} The maps (\ref{MPLGP}) provide an isomorphism of Poisson Lie groups
\be
{\mathscr L}_{\G, \odot_2} \stackrel{\sim}{\lra} \D(\G)^*.
\ee
\et

The proof is left as an exercise. 

\subsubsection{The geometric ${\cal R}-$matrix.} Consider the decorated surface 
$$
\odot_{n,n}:= \odot_n\ast\odot_n
$$
 obtained by gluing two decorated 
surfaces $\odot_n$ along the   boundary intervals. 
It carries $n$ pairs of punctures. The briad group element $\beta_{n}$ acts by interchanging the punctures in each pair, see Figure \ref{R14}. 
\begin{figure}[ht]
\begin{tikzpicture}[scale=1.2]
\draw (0,0) ellipse (0.8cm and 1cm);
\node[red] at (0,1) {$\bullet$};
\node[red] at (0,-1) {$\bullet$};
\node[blue,thick] at (-.2, 0.4) {$\circ$};
\node[blue,thick] at (.2, .4) {$\circ$};
\node[blue,thick] at (-.2, -0.4) {$\circ$};
\node[blue,thick] at (.2, -.4) {$\circ$};
\draw[-stealth] (-.2,0.44) arc (140: 40:.28);
\draw[-stealth] (.2,0.36) arc (-40: -140:.28);
\draw[-stealth] (-.2,-.36) arc (140: 40:.28);
\draw[-stealth] (.2,-0.44) arc (-40: -140:.28);
\end{tikzpicture}
\caption{The geometric ${\cal R}-$matrix is induced by the braid group element $\beta_{2}$.}
\label{R14}
\end{figure}
It induces a cluster automorphism ${\cal R}_{n}$ of the moduli space ${\mathscr L}_{\G, \odot_{n,n}}$, which we call the {\it geometric ${\cal R}-$matrix}:
\be
{\cal R}_{n}: {\mathscr L}_{\G, \odot_{n,n}}\lra {\mathscr L}_{\G, \odot_{n,n}}.
\ee
It provides an automorphism of algebras
\be \la{AUTR22}
{\cal R}^*_{n}: {\cal O}_q({\mathscr L}_{\G, \odot_{n,n}}) \lra {\cal O}_q({\mathscr L}_{\G, \odot_{n,n}}).
\ee
To justify the name, recall that  there is the $R-$matrix   for the quantum double  ${\U}_q(\D_{\mathfrak{g}})$ of ${\U}_q(\mathfrak{g})$ \cite{D}: \be
{\rm R}_{2} \in { \U}_q(\D_{\mathfrak{g}}) \  \widehat \otimes \ \U_q(\D_{\mathfrak{g}}).
\ee

\subsubsection{The quantum double ${\U}_q(\D_{\mathfrak{g}})$ of ${\U}_q(\mathfrak{g})$.} According to the general principles, see (\ref{149}), the quantum double 
Hopf algebra ${\U}_q(\D_{\mathfrak{g}})$ should be realized as the quantum algebra of regular functions on the dual Poisson Lie group $ \D(\G)^*$ of the Poisson Lie group 
$ \D(\G)$. So Theorem \ref{THD} suggests the following Conjecture.

\bcon There is   an isomorphism of Hopf algebras 
\be \la{drdr}
 {\U}_q(\D_{\mathfrak{g}})  \stackrel{\sim}{\lra}  {\cal O}_q({\mathscr L}_{\G, \odot_2}).
   \ee
 The    conjugation by   ${\rm R}_{2}$ is intertwined  with 
automorphism   ${\cal  R}_{2}$ in (\ref{AUTR22}) by the  
  map  of Hopf algebras,  induced by (\ref{drdr}) and the $2^r:1$ gluing map:
  \be
\U_q(\D_{\mathfrak{g}})\otimes  \U_q(\D_{\mathfrak{g}})  \stackrel{\sim}{\lra}  
{\cal O}_q({\mathscr L}_{\G, \odot_2}) \otimes {\cal O}_q({\mathscr L}_{\G, \odot_2}) \hra  {\cal O}_q({\mathscr L}_{\G, \odot_2\ast \odot_2}).
\ee
  \econ
   
Here is what we can see already. 
By the part a) of Theorem \ref{UEAB},  there are canonical  injective maps of algebras, assigned to the  special points $s, t$ on the decorated surface $\odot_2$,  provided by the quantum lifts of the potentials and elements ${\bf K}_i$ at these points: 
\be \la{GRD0}
\kappa_s, \kappa_t:   {\U}_q(\mathfrak{b}) \hra {\cal O}_q({\mathscr L}_{\G, \odot_2}).
\ee
By the part a) of Theorem \ref{UEAB}, they are maps of Hopf algebras. 
Since the decorated surface $\odot_2$ has exactly two special points, by the part b) of Theorem \ref{UEAB} these maps generate  an injective  map of Hopf algebras
\be \la{GRD}
\kappa:   {\U}_q(\mathfrak{g}) \hra {\cal O}_q({\mathscr L}_{\G, \odot_2}).
\ee

\vskip 2mm
\paragraph{\bf An example.} 
\begin{figure}
\begin{tikzpicture}
\draw  (0,0) circle (2cm);
\node at (0,2) {$\bullet$};
\node at (0,-2) {$\bullet$};
\node at (0,1) {$\bullet$};
\node at (0,-1) {$\bullet$};
\draw plot [smooth] coordinates {(0,2) (-0.4, 1.3) (-0.4, 0.7) (0,0) (0.4, -0.7) (0.4, -1.3) (0, -2)};
\draw (-.06,2) --(-.06, 1);
\draw (.06,2) --(.06, 1);
\draw (-.06,-2) --(-.06,-1);
\draw (.06,-2) --(.06,-1);
\node at (0.06,1.2) {$\bowtie$};
\node at (0.06,-1.2) {$\bowtie$};
\node (A) at (7, 0) {1};
\node (B) at (8, 0) {2};
\node (C) at (9, 0) {3};
\node (E) at (7.5, -1) {5};
\node (F) at (7.5, -2) {7};
\node (G) at (8.5, 1) {4};
\node (H) at (8.5, 2) {6};
\foreach \from\to in {B/A, B/C, A/E, E/B, A/F, F/B, C/G, G/B, C/H, H/B}
\draw[-latex] (\from) -- (\to);
\end{tikzpicture}
\caption{The moduli space $\mathscr{P}_{{\rm PGL}_2, \odot_2}$ and its quiver.}
\label{quiver.o2.m2}
\end{figure}

Figure \ref{quiver.o2.m2} shows a tagged triangulation of the disk $\odot_2$ with two punctures and two marked boundary points,   and the corresponding quiver. We  impose the following condition on the outer monodromy
\[
\mu_{\rm out}=X_{v} =1, \qquad v = 2e_1+2e_3+e_2+e_4+e_5+e_6+e_7.
\]
The algebra $\mathcal{O}_q(\mathscr{L}_{{\rm PGL}_2, \odot_2})$ contains $\mathcal{O}_q({\rm PGL}_2)$ and $\mathcal{U}_q(\mathfrak{sl}_2)$ as Hopf subalgebras.

After the extension from ${\rm PGL}_2$ to ${\rm SL}_2$,  the Hopf subalgebra $\mathcal{O}_q({\rm SL}_2)$ is generated by
\[
{\bf A}= X_a, \qquad {\bf B}=X_b, \qquad {\bf C}=X_c, \qquad {\bf D}=X_d+ X_{-a}.
\]
Here 
\[
a=\frac{-e_1-e_2-e_3}{2}, \qquad b=\frac{-e_1-e_2+e_3}{2}, \qquad c=\frac{e_1-e_2-e_3}{2}, \qquad d=\frac{e_1-e_2+e_3}{2}.
\]
Classically, they correspond to the entries of the $2\times 2$ unimodular matrix   
$ \begin{bmatrix} 
      {\bf A} & {\bf B} \\
      {\bf C} & {\bf D} \\
  \end{bmatrix}.$

The Hopf subalgebra ${\U}_q(\mathfrak{sl}_2)$ is generated by
\[
{\bf E}= X_{v_1}+ X_{v_2}+ X_{v_3}+ X_{v_4}+ X_{v_5}, \qquad 
{\bf F}= X_{w_1}+ X_{w_2}+ X_{w_3}+X_{w_4}+ X_{w_5},
\]
\[
{\bf K}= X_{k}.
\]
Here
\[
v_1=e_1, \qquad v_2= e_1+e_2, \qquad v_3=e_1+e_2+e_4, \qquad v_4=e_1+e_2+e_6, \qquad v_5=e_1+e_2+e_4+e_6,
\]
\[
w_1= e_3, \qquad w_2= e_3+e_2, \qquad w_3= e_3+e_2+e_5, \qquad w_4=e_3+e_2+e_7, \qquad w_5=e_3+e_2+e_5+e_7,
\]
\[
k=e_1+e_2+e_4+e_6+e_3.
\]

A direct check shows that the generators satisfy the following relations:
\be \la{QDUB}
\left\{
\begin{array}{ll}
{\bf E}{\bf A} &= q^{-\frac{1}{2}}(q-q^{-1}) {\bf C} + q^{-1} {\bf A}{\bf E}, \\
{\bf E}{\bf B} &= q^{-\frac{1}{2}}(q-q^{-1}) {\bf D} + q^{-1} {\bf B}{\bf E}- q^{-\frac{1}{2}}(q-q^{-1}){\bf A}{\bf K} ,\\
{\bf E}{\bf C} &=q {\bf C}{\bf E}, \\
{\bf E}{\bf D} &= q{\bf DE}-q^{\frac{3}{2}}(q-q^{-1}) {\bf C}{\bf K};
\end{array}
\right.
\ee
\be
\left\{
\begin{array}{ll}
{\bf F}{\bf A} &= q^{-\frac{1}{2}}(q-q^{-1}) {\bf B} + q^{-1} {\bf A}{\bf F}, \\
{\bf F}{\bf B} &=q {\bf B}{\bf F}, \\
{\bf F}{\bf C} &= q^{-\frac{1}{2}}(q-q^{-1}) {\bf D} + q^{-1} {\bf C}{\bf F}- q^{-\frac{1}{2}}(q-q^{-1}){\bf A}{\bf K}^{-1} ,\\
{\bf F}{\bf D} &= q{\bf DF}-q^{\frac{3}{2}}(q-q^{-1}) {\bf B}{\bf K}^{-1};
\end{array}
\right.
\ee
\be
\left\{
\begin{array}{ll}
{\bf KA}&={\bf AK},\\
 {\bf K}{\bf B} &= q^{-2}{\bf B}{\bf K}, \\ 
 {\bf KC} & =q^2{\bf CK}, \\
  {\bf KD}&={\bf DK}.
  \end{array}
\right.
\ee

\subsubsection{More Poisson Lie semigroups.} 
Denote by ${\rm Cyl}_{1, 1; \circ}$ the decorated surface given by a punctured cylinder with a   special point on each of the two boundary components, see Figure \ref{R8}. Then there is the   cluster Poisson moduli space  ${\cal P}_{\G, {\rm Cyl}_{1,1; \circ}}$, and the Poisson moduli space  
 ${\rm Loc}_{\G, {\rm Cyl}_{1,1; \circ}}$. The Poisson moduli space  ${\rm Loc}_{\G, {\rm Cyl}_{1,1; \circ}}$ is a semigroup, see Figure \ref{R8a}. Note that the group identity does not exist. 
 
\begin{figure}[ht]
\begin{tikzpicture}
\draw[thick] (0,0) ellipse (0.2cm and 0.5cm);
\draw[thick, dashed] (2,.5) arc (90:270:0.2cm and 0.5cm);
\draw[thick] (2,.5) arc (90:-90:0.2cm and 0.5cm);
\draw[thick] (0,-.5) -- (2,-.5);
\draw[thick] (0,.5) -- (2,.5);
\node[red] at (.2,0) {$\bullet$};
\node[red] at (2.2,0) {$\bullet$};
\node[blue] at (1.1,0) {$\circ$};
\end{tikzpicture}
\caption{The punctured  cylinder with a special point  on each boundary component.}
\label{R8}
\end{figure}
\begin{figure}[ht]
\begin{tikzpicture}
\draw[thick] (0,0) ellipse (0.2cm and 0.5cm);
\draw[thick, dashed] (2,.5) arc (90:270:0.2cm and 0.5cm);
\draw[thick] (2,.5) arc (90:-90:0.2cm and 0.5cm);
\draw[thick] (0,-.5) -- (2,-.5);
\draw[thick] (0,.5) -- (2,.5);
\node[red] at (.2,0) {$\bullet$};
\node[red] at (2.2,0) {$\bullet$};
\node[blue] at (1.1,0) {$\circ$};

\begin{scope}[xshift=3cm]
\draw[thick] (0,0) ellipse (0.2cm and 0.5cm);
\draw[thick, dashed] (2,.5) arc (90:270:0.2cm and 0.5cm);
\draw[thick] (2,.5) arc (90:-90:0.2cm and 0.5cm);
\draw[thick] (0,-.5) -- (2,-.5);
\draw[thick] (0,.5) -- (2,.5);
\node[red] at (.2,0) {$\bullet$};
\node[red] at (2.2,0) {$\bullet$};
\node[blue] at (1.1,0) {$\circ$};
\end{scope}

\begin{scope}[xshift=3.2cm, yshift=-.5cm]
\draw[directed]  (2.4, 0.5) -- (3.8,0.5);
\node at (3.1, 0.7) {gluing};
\end{scope}

\begin{scope}[xshift=7.5cm]
\draw[thick] (0,0) ellipse (0.2cm and 0.5cm);
\draw[thick, dashed] (3,.5) arc (90:270:0.2cm and 0.5cm);
\draw[thick] (3,.5) arc (90:-90:0.2cm and 0.5cm);
\draw[thick] (0,-.5) -- (3,-.5);
\draw[thick] (0,.5) -- (3,.5);
\node[red] at (.2,0) {$\bullet$};
\node[red] at (3.2,0) {$\bullet$};
\node[blue] at (1,0) {$\circ$};
\node[blue] at (2.2,0) {$\circ$};
\draw[thick, dashed, green] (1.5,.5) arc (90:270:0.2cm and 0.5cm);
\draw[thick, green] (1.5,.5) arc (90:-90:0.2cm and 0.5cm);
\node[red] at (1.7,0) {$\bullet$};
\draw[blue, dashed] (1.6,0) ellipse (0.9cm and 0.2cm);
\end{scope}

\begin{scope}[xshift=9.1cm, yshift=-.5cm]
\draw[directed]  (2.1, 0.5) -- (4.1,0.5);
\node at (3.1, 0.7) {encircling};
\node at (3.1, 0.3) {punctures};
\end{scope}

\begin{scope}[xshift=14cm]
\draw[thick] (0,0) ellipse (0.2cm and 0.5cm);
\draw[thick, dashed] (2,.5) arc (90:270:0.2cm and 0.5cm);
\draw[thick] (2,.5) arc (90:-90:0.2cm and 0.5cm);
\draw[thick] (0,-.5) -- (2,-.5);
\draw[thick] (0,.5) -- (2,.5);
\node[red] at (.2,0) {$\bullet$};
\node[red] at (2.2,0) {$\bullet$};
\node[blue] at (1.1,0) {$\circ$};
\end{scope}
\end{tikzpicture}
\caption{The composition law via the gluing of the discs and encircling the pairs of the punctures.}
\label{R8a}
\end{figure}

\medskip
\section{Quantized  moduli spaces of $\G$-local systems   as a geometric avatar of  TQFT} \la{TQFT}

\medskip

In Section \ref{TQFT1} we explain how these constructions fit into the general framework of extended TQFT provided by the moduli spaces we consider. 
 
 Using this, in Section \ref{SEC1.8} we   
define the principal series of $\ast$-representations of the  quantum group modular double. The proofs and other details are deferred to Section \ref{SSEECC11.3}.

In Section \ref{SC2.7} we introduce canonical realizations of the principal series, the spaces of tensor product invariants, and  the regular representation.  The proofs and other details are postponed to Section \ref{Sec8}.

Combining Sections \ref{TQFT1}, \ref{SEC1.8},  \ref{SC2.7}, we see that the linearization of the geometric TQFT, given by the cluster quantization of the quantum 
 moduli spaces ${\mathscr P}_{\G, \bS}$, is the continuous analog of the braided monoidal category  provided by the principal series representations. 

\medskip

\subsection{Quantized  moduli spaces of $\G$-local systems on decorated surfaces as a TQFT} \la{TQFT1}

\medskip

  The moduli space  ${\mathscr P}_{\G, \bS}$ parametrizes   objects which   localize the \underline{topological} notion of a  
  $\G$-local system on a  
  surface $\bS$ using topologically trivial objects - triangles.  
The extra data, given by the   flags at the vertices and pinnings  at the sides, allow to build the moduli space ${\mathscr P}_{\G, \bS}$ 
from the elementary  ones ${\mathscr P}_{\G, t}$, and in particular to reconstruct a    $\G$-local system on $\bS$ from the data assigned to the triangles.  
 The localization on triangles persists on the quantum level. 
 
This contrasts with the standard technique of  pair of pants decompositions of surfaces, where the elementary objects - pair of pants - are topologically non-trivial. 
The Modular Functor Conjecture \ref{MK}  relates the two approaches.

The emerging structure  can be viewed as a continuous analog 
 of a TQFT. 
 Indeed, in the classical TQFT we assign to a circle a semisimple monoidal category ${\cal C}$ with simple objects $X_\alpha$, where 
 the set of indices $\{\alpha\}$ is discrete, and the tensor product decomposition is given by a finite sum, where $C_{\alpha, \beta}^\gamma$ are finite dimensional vector spaces: 
  $$
 X_\alpha \otimes_{\cal C} X_\beta = \sum_{\gamma}C_{\alpha, \beta}^\gamma \otimes X_\gamma. 
 $$
 Equivalently,  the spaces of invariants in the tensor product, called the {\it multiplicity spaces}, or {\it correlators}:
 $$
 \langle X_{\alpha_1} \otimes \ldots \otimes  X_{\alpha_n} \rangle_{\cal C} := {\rm Hom}_{\cal C}({\bf 1}_{\cal C}, X_{\alpha_1} \otimes \ldots \otimes  X_{\alpha_n})
 $$
  are finite dimensional for any $n \geq 2$,  and for any collection of simple objects $X_{\alpha_1},  \ldots , X_{\alpha_n}$. 
  \vskip 2mm
  
In our case the   simple objects are the principal series representations $V_\alpha$ of the modular double ${\cal A}_\hbar(\mathfrak{g})$ of 
the quantum group $\U_q(\mathfrak{g})$. They are 
parametrised by the points of an orbifold: 
\be \la{ALPHA}
\alpha \in \H/W(\R_{>0}).
\ee
As we will see momentarily, Modular Functor Conjecture \ref{MK} implies that  the multiplicity spaces are infinite dimensional and   given by 
\be \la{TQFTCL}
   \langle V_{\alpha_1} \otimes \ldots  \otimes  X_{\alpha_n} \rangle_{\cal C} = {\cal H}_{\G, S^2 - \{p_1, ..., p_n\}}(\alpha_1, ..., \alpha_n).
\ee
The weights $\alpha_k$ are assigned to the punctures $p_k$ on  the $n$-punctured sphere $S^2 - \{p_1, ..., p_n\}$. 
The space 
on the right 
 is  the cluster Hilbert space,  provided by the cluster  quantization  of the symplectic leave of the space 
${\rm Loc}_{\G, S^2 - \{p_1, ..., p_n\}}$,  parametrising  $\G-$local systems on $S^2 - \{p_1, ..., p_n\}$ with given 
semisimple parts $\alpha_k$ of  monodromies around the punctures $p_k$.  The two key properties of the multiplicity spaces are the following:

\begin{itemize}

\item The  braid group of $S^2 - \{p_1, ..., p_n\}$ projective action on the multiplicity spaces. 

\vskip 1mm\item The analog of the associativity  of the tensor product: 
\be \la{TQFTCL1}
 \begin{split}
   \langle V_{\alpha_1} \otimes V_{\alpha_2} \otimes  V_{\alpha_3}    \otimes  X_{\alpha_4} \rangle_{\cal C} = &\int_{\beta    \in \H/W(\R_{>0})}
   \langle V_{\alpha_1} \otimes V_{\alpha_2} \otimes    X_{\beta} \rangle_{\cal C} \otimes \langle V_{\beta^*} \otimes V_{\alpha_3} \otimes    X_{\alpha_4} \rangle_{\cal C} d\beta\\
   =& \int_{\gamma    \in \H/W(\R_{>0})}   
   \langle V_{\alpha_2} \otimes V_{\alpha_3} \otimes    X_{\gamma} \rangle_{\cal C} \otimes \langle V_{\gamma^*} \otimes V_{\alpha_4} \otimes    X_{\alpha_1} \rangle_{\cal C} d\gamma.\\
   \end{split}
   \ee
\end{itemize}

The first property follows from Theorem \ref{Th1.12}: the braid group of $S^2 - \{p_1, \ldots ,  p_n\}$ acts on the cluster quantization space on the right of (\ref{TQFTCL}). 
The second follows from the Modular Functor Conjecture \ref{MK}, 
applied to the two loops $\beta$ and $\gamma$ on  $S^2 - \{p_1, p_2, p_3, p_4\}$ on Figure \ref{pin16a}.

\begin{figure}
 \begin{center}
 \begin{tikzpicture}[scale=1.3]

\node at (-1.3, -0.8) {$\alpha_1$};
\node at (-1.3, 0.8) {$\alpha_2$};
\node at (1.3, 0.8) {$\alpha_3$};
\node at (1.3, -0.8) {$\alpha_4$};

 \begin{scope}[shift={(1.05,0.61)}]
 \draw  [blue][thick, rotate=30]  (0,0) ellipse (0.07 and 0.1);
 \end{scope}
  \begin{scope}[shift={(1.05,-0.61)}]
 \draw  [blue][thick, rotate=-30]  (0,0) ellipse (0.07 and 0.1);
 \end{scope}
  \begin{scope}[shift={(-1.05,-0.61)}]
 \draw  [blue][thick, rotate=30]  (0,0) ellipse (0.07 and 0.1);
 \end{scope}
  \begin{scope}[shift={(-1.05,0.61)}]
 \draw  [blue][thick, rotate=-30]  (0,0) ellipse (0.07 and 0.1);
 \end{scope}
 \node at (-0.2,0) {$\beta$};
\draw [blue][thick] (1, 0.7) arc (-60:-120:2);
\draw  [blue][thick] (1, -0.7) arc (60: 120:2);
\draw  [blue][thick] (1.13, 0.55) arc (150:210:1.1);
\draw  [blue][thick] (-1.13, -0.55) arc (-30:30:1.1);
\draw [blue][thick] (0, 0.432) arc (90:270:0.1 and 0.432);
\draw  [blue][thick, dashed] (0, 0.432) arc (90:-90:0.1 and 0.432);

\begin{scope}[shift={(5,0)}]
 \node at (-1.3, -0.8) {$\alpha_1$};
\node at (-1.3, 0.8) {$\alpha_2$};
\node at (1.3, 0.8) {$\alpha_3$};
\node at (1.3, -0.8) {$\alpha_4$};
 \begin{scope}[shift={(1.05,0.61)}]
 \draw  [blue][thick, rotate=30]  (0,0) ellipse (0.07 and 0.1);
 \end{scope}
  \begin{scope}[shift={(1.05,-0.61)}]
 \draw  [blue][thick, rotate=-30]  (0,0) ellipse (0.07 and 0.1);
 \end{scope}
  \begin{scope}[shift={(-1.05,-0.61)}]
 \draw  [blue][thick, rotate=30]  (0,0) ellipse (0.07 and 0.1);
 \end{scope}
  \begin{scope}[shift={(-1.05,0.61)}]
 \draw  [blue][thick, rotate=-30]  (0,0) ellipse (0.07 and 0.1);
 \end{scope}
 \node at (0,0.15) {$\gamma$};
\draw [blue][thick] (1, 0.7) arc (-60:-120:2);
\draw  [blue][thick] (1, -0.7) arc (60: 120:2);
\draw  [blue][thick] (1.13, 0.55) arc (150:210:1.1);
\draw  [blue][thick] (-1.13, -0.55) arc (-30:30:1.1);
\draw [blue][thick] (0.98, 0) arc (0:180:0.98 and 0.1);
\draw  [blue][thick, dashed] (0.98, 0) arc (0:-180:0.98 and 0.1);
\end{scope}
 \end{tikzpicture}
 \end{center}
\caption{Identity (\ref{TQFTCL1}) follows from the Modular Functor Conjecture  by cutting the  sphere with four punctures, decorated by the weights $\alpha_1, ..., \alpha_4$, by the loops $\beta$ and $\gamma$.}
\label{pin16a}
\end{figure}

 \vskip 2mm
 
 So to get a continuous analog of a TQFT given by multiplicity spaces with the above properties it is sufficient to consider just surfaces with punctures and 
 quantize the related moduli of $\G-$local systems. 
 
 However to go one step further, and construct 
 a continuous analog   of  the braided monoidal category ${\cal C}_{\mathfrak{g}}$ of the quantum group principal series representations    we need   cluster quantization of 
the  moduli spaces ${\mathscr P}_{\G, \bS}$ assigned to \underline{decorated} surfaces $\bS$.  
The basic decorated surface we need is 
the punctured disc  $\odot$ with two special  points, see Figure \ref{gs1}. The circle $S^1$ to which the category  ${\cal C}_{\mathfrak{g}}$  
is assigned is  the boundary of the punctured 
disc $\odot$ corresponding to the puncture, expanded as a hole. 
The   simple objects  -  the principal series representations $V_\alpha$ -  
 are   assigned to  the  generic symplectic leaves of the moduli space ${\mathscr L}_{\G, \odot}$. The latter are   obtained by fixing    the 
semi-simple part of the monodromy $\alpha\in \H/W(\R_{>0})$ of a $\G-$local system around the puncture. 
The space ${\mathscr L}_{\G, \odot}$ is essentially the  moduli space of Stokes data for the flat connections on $\C^\times$ with an irregular singularity - an order two pole - at $\infty$, and regular at $0$. 
To get  formula (\ref{TQFTCL}) we apply   Modular Functor Conjecture \ref{MK} to the loop $\gamma$ on the disc with $n-1$ punctures   and two 
special boundary points, which encircles the punctures. See Figure \ref{pin5a}  for the $n=3$ case. 
\vskip 2mm

In Section \ref{TQFT} we discuss the algebraic-geometric properties of the moduli spaces ${\mathscr P}_{\G, \bS}$ which underly the existence of  such extended TQFT's. 
\vskip 2mm

\subsubsection{The moduli spaces ${\mathscr L}_{\G, \odot_n}$.}  \la{2.5.1}
 \begin{figure}[ht]
\epsfxsize 200pt
\center{
\begin{tikzpicture}[scale=0.6]
\draw[->, >=stealth, thick] (3.5,0.3) -- (6, 0.3);
\draw[<-, >=stealth, thick] (3.5,-0.3) -- (6, -0.3);
\node at (4.75, 0.7) {cut};
\node at (4.75, -0.7) {glue};
\draw (0,0) ellipse (2.5cm and 1cm);
\draw[dashed] (0,0) ellipse (8mm and 1cm);
\node[red, label=above: {\small $s$}] at (0,1) {\small $\bullet$};
\node[red, label=below: {\small $t$}] at (0,-1) {\small $\bullet$};
\node[red] at (0,0) {\Large $\circ$};
\node[red] at (1.5, 0) {\Large $\circ$};
\node[red] at (-1.5, 0) {\Large $\circ$};
\draw (8,1) arc(90:270:1.2cm and 1cm);
\draw (11,-1) arc(-90:90:1.2cm and 1cm);
\draw[dashed] (8,1)--(8,-1);
\draw[dashed] (11,1)--(11,-1);
\draw[dashed] (9.5,0) ellipse (8mm and 1cm);
\node[red] at (8,1) {\small $\bullet$};
\node[red] at (8,-1) {\small $\bullet$};
\node[red] at (7.5,0) {\Large $\circ$};
\node[red] at (9.5,0) {\Large $\circ$};
\node[red] at (9.5,1) {\small $\bullet$};
\node[red] at (9.5,-1) {\small $\bullet$};
\node[red] at (11.5,0) {\Large $\circ$};
\node[red] at (11,1) {\small $\bullet$};
\node[red] at (11,-1) {\small $\bullet$};
\end{tikzpicture}
 }
\caption{Cutting the decorated surface $\odot \ast \odot \ast \odot$ by  two arcs into three decorated surfaces $\odot$.}
\label{pin16}
\end{figure}

Recall the decorated surface $\odot_n = \odot \ast \cdots \ast \odot$   given by a   disc with $n$ punctures and two ordered special  points $s, t$.\footnote{We use slightly different conventions: in the previous Section  the punctures where located on the axis through the special points, while now we picture them   ordered   from the left to the right.}   
 It can be    obtained by gluing of $n$ copies of the decorated surfaces $\odot$ along boundary intervals, see Figure \ref{pin16}.
Recall the outer  monodromy map 
\be \la{MUOUT}
 \mu_{\rm out}: {\rm Loc}_{\G, \odot_n} \lra \H.
\ee
  Its fiber over the unit $e \in \H$ is denoted by  
  $$
  {\mathscr L}_{\G, \odot_n}:= \mu_{\rm out}^{-1}(e)\subset {\rm Loc}_{\G, \odot_n}. 
    $$
The gluing of the moduli spaces ${\rm Loc}_{\G, \odot_n}$ results in the product of outer monodromies. So   
the gluing  provides a gluing map of the moduli spaces ${\mathscr L}_{\G, \odot_n}$:
\be \la{12.30.18.10}
{\mathscr L}_{\G, \odot_n} \ast {\mathscr L}_{\G, \odot_m} \lra {\mathscr L}_{\G, \odot_{n+m}}.
\ee
Denote by $\widetilde {\mathscr L}_{\G, \odot_n}$  the fiber over $e$ of the  map 
$\mu_{\rm out}: {\mathscr P}_{\G, \odot _n} \to\H$. The gluing   induces a map 
$$
\widetilde {\mathscr L}_{\G, \odot_n} \ast \widetilde {\mathscr L}_{\G, \odot_m} \lra \widetilde {\mathscr L}_{\G, \odot_{n+m}}.
$$
 
\subsubsection{Factorization property.} The key feature of  the  space ${\mathscr L}_{\G, \odot_n}$ is that the map (\ref{12.30.18.10}) is almost an isomorphism - it is a finite dominant map. 
Namely, given a collection of $(n-1)$ arcs $\beta$ connecting    special points $s, t$ and subdividing the surface into $n$   punctured discs, 
see Figure \ref{pin16}, there is  a gluing map: 
\be \la{CRPP}
\begin{split}
&F_{\beta}: {\mathscr L}_{\G, \odot} \times \ldots \times {\mathscr L}_{\G, \odot} \stackrel{}{\lra} {\mathscr L}_{\G, \odot_n}.\\
\end{split}
\ee

The space ${\mathscr L}_{\G, \odot_n}(\R_{>0})$ of the 
real positive points of  ${\mathscr L}_{\G, \odot_n}$ is well defined. Indeed, the  space ${\mathscr P}_{\G, \odot_n}$ has a   cluster Poisson structure by Theorem \ref{MTH}, 
so the manifold ${\mathscr P}_{\G, \odot_n}(\R_{>0})$ is well defined. The space ${\mathscr L}_{\G, \odot_n}(\R_{>0})$ is obtained from it by taking 
the quotient by the  action of the Weyl group $W^n$, which is positive by the part 1) of Theorem \ref{MTH}, and then taking the fiber of the outer monodromy map (\ref{MUOUT}). 

Let $S$ be the subgroup of   order two elements in the Cartan group $\H$. Then $|S|=2^r$, where $r={\rm dim}~\H$. 

\bl 1. The  map $F_\beta$ is a Galois cover with the Galois group $S^{n-1}$ over the generic point. 

\vskip 1mm2. It induces an isomorphism over the spaces of real positive points:
\be \la{MAPFR}
F^+_{\beta}: {\mathscr L}_{\G, \odot}(\R_{>0}) \times \ldots \times {\mathscr L}_{\G, \odot}(\R_{>0}) \stackrel{\sim}{\lra} {\mathscr L}_{\G, \odot_n}(\R_{>0}).\ee
\el 

\begin{proof}  1. Let us  restrict a framed $\G$-local system on $\odot_n$ to    punctured discs $\odot$ obtained by cutting along the arcs $\beta$. 
Assume that the pair of flags at each arc  $\beta_i$ is generic. Then we can take a pinning $p_i$ for the pair of flags sitting on the arc. 
Since the pinning on the left  boundary segment of $\odot_n$ is given, one can find inductively pinnings $p_1, \ldots , p_{n-1}$ so that  the   
  ``outer monodromy is $e$" condition holds on each of the left $n-1$ discs. Finally,  the 
  outer  monodromy around rightmost disc is   $e$ since it is so for the space $\odot_n$. For each arc $\beta_1, \ldots, \beta_{n-1}$ the group $S$ acts simply transitively on the set of such  pinnings $p_i$.  
  Indeed, the action of an element $h \in \H$ on pinnings amounts to the action of $h^2$ on the outer monodromy.  
  
  2. The positivity implies that the pair of flags at each arc  $\beta_i$ is generic. So the map (\ref{MAPFR}) is onto. The action of any non-trivial element of the group $S$ destroys positivity, so the map is injective. \end{proof} 

Therefore the space ${\mathscr L}_{\G, \odot_n}$ can be almost factorised, in many different ways,  into a product of  moduli spaces  ${\mathscr L}_{\G, \odot}$. Note that   the factorization property fails for the  moduli spaces ${\mathscr P}_{\G, \odot_n}$. Indeed,   if we replace the ${\mathscr L}$-spaces by the ${\mathscr P}$-spaces,  dimensions of  the spaces in (\ref{CRPP}) would be different. 

\subsubsection{Braid group action.} By   Theorem \ref{MTH}  the  braid group ${\rm Br}_n$, realized as the subgroup of the mapping 
 class group  of     $\odot_n$ moving the punctures,   acts as a cluster  modular group  on   ${\mathscr P}_{\G, \odot_n}$, and preserves  
 ${\mathscr L}_{\G, \odot_n}$. 
 
 \bt \la{T1.22} There is a unique birational action of the braid group ${\rm Br}_n$ on the space  ${\mathscr L}_{\G, \odot} \times \ldots \times {\mathscr L}_{\G, \odot} $ in (\ref{CRPP}), lifting the action of  
 on ${\mathscr L}_{\G, \odot_n}$   and preserving  ${\mathscr L}_{\G, \odot_n}(\R_{>0})$.  It  relates different projections $F_\beta$. 

Similar results hold for the space $\widetilde {\mathscr L}_{\G, \odot_n}$. In particular, for $n=2$ we get the   birational  automorphisms  \be
\begin{split}
 &{\cal R}:  {\mathscr L}_{\G, \odot} \times {\mathscr L}_{\G, \odot} \lra {\mathscr L}_{\G, \odot} \times {\mathscr L}_{\G, \odot},\\
 &\widetilde   {\cal R}:  \widetilde {\mathscr L}_{\G, \odot} \times \widetilde {\mathscr L}_{\G, \odot} \lra \widetilde {\mathscr L}_{\G, \odot} \times 
  \widetilde {\mathscr L}_{\G, \odot}.\\
  \end{split}
  \ee
   lifting the    mapping class group  element  
  interchanging the two punctures on $\odot_2$: 
\begin{figure}[ht]
\begin{center}
\begin{tikzpicture}[scale=0.8]
\draw (0,0) ellipse (1.5cm and .5cm);
\node[red] at (90:.5) {\tiny $\bullet$};
\node[red] at (-90:.5) {\tiny $\bullet$};
\node[blue] at (0:.5) {\tiny $\bullet$};
\node[blue] at (180:.5) {\tiny $\bullet$};
 \draw[blue, -latex] (180:.5) arc (150:30:.57);
  \draw[blue, directed] (0:.5) arc (-30:-150:.57);
\end{tikzpicture}
\end{center}
 \caption{The brading is provided by    the mapping class group element  interchanging the punctures.}
\label{pin25}
\end{figure}
 \et 
 The uniqueness  is clear from   isomorphism (\ref{MAPFR}). Theorem \ref{T1.22} will not be used in this paper.




\subsubsection{Monoidal structure.} Gluing  of   ${\mathscr P}$-spaces is compatible with the cluster Poisson amalgamation. A collection of  arcs $\beta$  cutting the decorated surface $\odot_n$  gives rise to a  
  $W^n$-equivariant map of $\ast$-algebras
 \be
\mathcal{A}_\hbar({\mathscr P}_{\G, \odot \ast ... \ast \odot})   \stackrel{}{\lra} \mathcal{A}_\hbar({\mathscr P}_{\G, \odot} \times \ldots \times {\mathscr P}_{\G, \odot}) = 
\mathcal{A}_\hbar({\mathscr P}_{\G, \odot} )\otimes \ldots \otimes \mathcal{A}_\hbar({\mathscr P}_{\G, \odot}).
  \ee 
Passing to the   outer monodromy =1 quotients,  we  arrive at the map of $\ast$-algebras induced by  map (\ref{CRPP}):
 \be \la{CRPP1}
F^*_{\beta}: \mathcal{A}_\hbar({\mathscr L}_{\G, \odot \ast ... \ast \odot})   \stackrel{}{\lra} \mathcal{A}_\hbar({\mathscr L}_{\G, \odot} \times \ldots \times {\mathscr L}_{\G, \odot}).
\ee
The quantization  delivers a map of  Hilbert spaces intertwining, just as in Conjecture \ref{MK},  the map   (\ref{CRPP1}):
\be \la{12.30.18.1}
F^\circ_{\beta}:  {\cal H}({\mathscr L}_{\G, \odot \ast ... \ast \odot}) \stackrel{}{\lra} {\cal H}({\mathscr L}_{\G, \odot}) \otimes \ldots \otimes {\cal H}({\mathscr L}_{\G, \odot}).
\ee 
Using the fact that the map (\ref{MAPFR}) is an isomorphism, one can see that the map (\ref{12.30.18.1}) is an isomorphism.  

\vskip 1mm
Next, there is a canonical projection, obtained by taking a simple loop $\gamma$ in   
$\odot \ast ... \ast \odot$ containing all the punctures, and restricting to the complement of the loop $\gamma$  which does not contain the punctures:
\be \la{12.30.18.4}
{\rm Res}_\gamma: {\mathscr L}_{\G, \odot \ast ... \ast \odot} \lra {\mathscr L}_{\G, \odot}.
\ee
It should give rise to   the dual map of $\ast-$algebras. 
In fact, we consider a subalgebra ${\mathscr A}'_\hbar({\mathscr L}_{\G, \odot}) \subset \mathcal{A}_\hbar({\mathscr L}_{\G, \odot})$, see (\ref{Oprime}), which should coincide with the whole algebra, and then   get a map of $\ast$-algebras
$$
{\rm Res}_\gamma^*: {\mathscr A}'_\hbar({\mathscr L}_{\G, \odot}) \lra \mathcal{A}_\hbar({\mathscr L}_{\G, \odot \ast ... \ast \odot}).
$$
It allows to view the Hilbert space ${\cal H}({\mathscr L}_{\G, \odot \ast ... \ast \odot})$ as an ${\mathscr A}'_\hbar({\mathscr L}_{\G, \odot})$-module. 
Using the factorization map (\ref{12.30.18.1}), we get an ${\mathscr A}'_\hbar({\mathscr L}_{\G, \odot})$-module structure on the tensor product 
$$
{\cal H}({\mathscr L}_{\G, \odot}) \otimes_\beta \ldots \otimes_\beta {\cal H}({\mathscr L}_{\G, \odot}).
$$
The notation $\otimes_\beta$  emphasizes that it depends on the isotopy class of the collection of arcs $\beta$.

Recall the decomposition of the $\mathcal{A}_\hbar({\mathscr L}_{\G, \odot})$-module according to the action of the center, providing    the principal series of representations  
 ${\cal H}({\mathscr L}_{\G, \odot})_\alpha$ of the $\ast$-algebra $\mathcal{A}_\hbar({\mathscr L}_{\G, \odot})$:
 \be \la{12.30.18.14}
{\cal H}({\mathscr L}_{\G, \odot}) = \int {\cal H}({\mathscr L}_{\G, \odot})_\alpha d\alpha, \qquad \alpha \in \H(\R_{>0})/W.
\ee
The  $\otimes_\beta-$product respects   decomposition (\ref{12.30.18.14}). 
So we get tensor products of elementary objects:
\be \la{12.30.18.114}
{\cal H}({\mathscr L}_{\G, \odot})_{\alpha_1} \otimes_\beta \ldots \otimes_\beta {\cal H}({\mathscr L}_{\G, \odot})_{\alpha_n}.
\ee
Incorporating the dependence on $\beta$, we get local systems of Hilbert spaces on   ${\cal M}_{0, n}$.

\subsubsection{Decomposition of the tensor product.} To decompose the tensor product of the principal series ${\mathscr A}'_\hbar({\mathscr L}_{\G, \odot})$-modules 
(\ref{12.30.18.114}), let us 
cut the decorated surface $\odot_n$ along the loop $\gamma$,  getting a sphere with $n$ punctures and a punctured disc $\odot$ with two special points, see Figure \ref{pin5a}:
  \be \la{ALPHACUT}
  \odot\ast \ldots \ast \odot = S^2-\{p_1, ..., p_{n}\} \cup_\gamma    \odot.
\ee
We use the notation
  $$
  S^2_n:= S^2 -\{p_1, ..., p_{n}\}.
  $$ 
 Then  modular functor Conjecture \ref{MK} applied to the loop $\gamma$ provides   an isomorphism 
\be \la{12.30.18.2}
{\rm MFK}_\gamma: {\cal H}({\mathscr L}_{\G, \odot_{n-1}}) \stackrel{\sim}{\lra} {\cal H}({\mathscr L}_{\G, \odot}) \otimes   {\cal H}({\rm Loc}_{\G, S^2_n}).
\ee
Combining factorization isomorphism (\ref{12.30.18.1}) with (\ref{12.30.18.2}), we arrive at the isomorphism 
\be \la{12.30.18.3}
{\rm MFK}_\gamma\circ {F^\circ_\alpha}^{-1}:  \underbrace{{\cal H}({\mathscr L}_{\G, \odot}) \otimes \ldots \otimes {\cal H}({\mathscr L}_{\G, \odot})}_{\text{$n-1$ factors}}
 \stackrel{\sim}{\lra} {\cal H}({\mathscr L}_{\G, \odot}) \otimes   {\cal H}({\rm Loc}_{\G, S^2_n}).
\ee 
Decomposing it  with respect to the characters of the center, we get:
\be \la{12.30.18.13}
\begin{split}
& {\cal H}({\rm Loc}_{\G, S_n^2}) = \int {\cal H}({\rm Loc}_{\G, S^2_n})_{\alpha_1, ..., \alpha_n} d\alpha_1 \ldots d\alpha_n.\\
 &  {\cal H}({\mathscr L}_{\G, \odot})_{\alpha_1} \otimes \ldots \otimes {\cal H}({\mathscr L}_{\G, \odot})_{\alpha_{n-1}} =\int 
 {\cal H}({\mathscr L}_{\G, \odot})_{\alpha_{n}}  \otimes   {\cal H}({\rm Loc}_{\G, S^2_n})_{\alpha_1, ...,  \alpha_{n}}d\alpha_{n}.\\
 \end{split}
\ee

 \begin{figure}[ht]
\epsfxsize 200pt
\center{
\begin{tikzpicture}[scale=1]
\draw[dashed, red] (0,0) ellipse (1cm and 1cm);
\draw (0,0) ellipse (0.4cm and 0.7cm);
\draw[thick] (0,0.4) circle (0.7mm);
\draw[thick] (0,-0.4) circle (0.7mm);
\node[red] at (0,1) {{\tiny $\bullet$}};
\node[red] at (0,-1) {{\tiny $\bullet$}};
\draw[dashed, red] (5,0) ellipse (1cm and 1cm);
\node[red] at (5,1) {{\tiny $\bullet$}};
\node[red] at (5,-1) {{\tiny $\bullet$}};
\draw[thick] (8.4,0.53) circle (0.7mm);
\draw[thick] (8.4,-0.53) circle (0.7mm);
\draw (8.4,0.46) arc(90:270:0.3cm and .46cm);
\draw[thick] (0,-0.4) circle (0.7mm);
\draw (5,-.6) arc(115:65:4cm and 4cm);
\draw (5,.6) arc(-115:-65:4cm and 4cm);
\draw (6.5,.23) arc(90:-90:0.12cm and .23cm);
\draw[dashed] (6.5,.23) arc(90:270:0.12cm and .23cm);
\node at (-0.6,0) {$\gamma$};
\node at (6.5,-.5) {$\gamma$};
\end{tikzpicture}
 }
\caption{Decomposing the tensor product of  the principal series.}
\label{pin5a}
\end{figure}

Let us summarize the structures  provided by the   Modular Functor Conjecture isomorphism: 

\begin{enumerate}

\item  The action of the $\ast$-algebra ${\mathscr A}'_\hbar({\mathscr L}_{\G, \odot})$ on the   factor 
${\cal S}({\rm Loc}_{\G, S^2_n})$ is trivial, and its action on ${\cal S}({\mathscr L}_{\G, \odot \ast ... \ast \odot})$ is provided by the map (\ref{12.30.18.4}). 

\vskip 1mm \item  The action of  ${\mathscr A}'_\hbar({\mathscr L}_{\G, \odot})$ on  ${\cal S}({\mathscr L}_{\G, \odot}) \otimes \ldots \otimes {\cal S}({\mathscr L}_{\G, \odot})$ induced by 
the action on ${\cal S}({\mathscr L}_{\G, \odot \ast ... \ast \odot})$  and  isomorphism (\ref{12.30.18.1}) is  the action on the tensor product of ${\mathscr A}'_\hbar({\mathscr L}_{\G, \odot})$-modules ${\cal S}({\mathscr L}_{\G, \odot})$. 
\end{enumerate}

This just means that  isomorphisms (\ref{12.30.18.3}) and (\ref{12.30.18.13}) provide a decomposition of the tensor product of   principal series representations of the $\ast$-algebras 
${\mathscr A}'_\hbar({\mathscr L}_{\G, \odot})$  into an integral of the  principal series  
representations   with the multiplicity spaces given by the Hilbert spaces ${\cal H}({\rm Loc}_{\G, S^2_n})_{\alpha_1, ..., \alpha_n}$: 
 
 \bt \la{MFCQI} Assume    Conjecture \ref{MK}. Then   the local system of Hilbert spaces ${\cal H}({\rm Loc}_{\G, S^2_n})_{_{\alpha_1, ..., \alpha_n}}$ on ${\cal M}_{0,n}$   is identified with the local system  of invariants of the tensor product of principal series of 
  representations ${\cal H}_{\alpha_1}, \ldots , {\cal H}_{\alpha_n}$ of the  $\ast$-algebra   $  {\mathscr A}'_\hbar({\mathscr L}_{\G, \odot})$:
  \be \la{TENV}
{\cal H}({\rm Loc}_{\G, S^2_n})_{\alpha_1, ..., \alpha_n}=   \Bigl({\cal H}_{\alpha_1}\otimes  \ldots \otimes {\cal H}_{\alpha_n}\Bigr)^{\mathcal{A}_\hbar({\mathscr L}_{\G, \odot})}.
\ee  
\et   
So the Hilbert space ${\cal H}({\rm Loc}_{\G, S})_{\underline{\alpha} }$ for a  punctured surface $S$   generalizes  the space of tensor product invariants  of the principal 
  series    representations of  the $\ast$-algebra  ${\mathscr A}'_\hbar({\mathscr L}_{\G, \odot})$. 
  
\subsubsection{Conclusion.} The geometric properties of   moduli spaces ${\mathscr P}_{\G, \bS}$ and     ${\mathscr L}_{\G, \odot}$  
       show that representations of the $\ast$-algebra  ${\mathscr A}'_\hbar({\mathscr L}_{\G, \odot})$ in Hilbert spaces ${\cal H}({\rm Loc}_{\G, S})_{\underline{\alpha} }$ 
       form a  continuous analog of a braided monoidal category. The monoidal structure is given by the gluing map for the  moduli spaces ${\mathscr L}_{\G, \odot}$, see Figure \ref{pin16}, and the 
     brading is given by the map ${\rm Br}$, see Figure \ref{pin25}. 
     The key facts underlying this  are 
     
\begin{itemize}

\item  Cluster Poisson nature of the moduli spaces ${\mathscr P}_{\G, \bS}$,  ${\mathscr L}_{\G, \odot}$ and  their   group of symmetries $\Gamma_{\G, \bS}$.  
     
\vskip 1mm\item The very existence and the cluster nature of the gluing map.  
     
\end{itemize}

 This braided monoidal ``category" describes the principal series  $\ast$-representations of the 
  modular double  $\mathcal{A}_\hbar(\frak{g})$ of the quantum group $\U_q(\frak{g})$, see  Section  \ref{SEC1.8}.

  \medskip
  
   \subsection{Principal series of $\ast$-representations of quantum groups}  \la{SEC1.8}
   
 \medskip 
  \subsubsection{The principal series of unitary  representations of the group $\G(\R)$.}  \la{SECT5.2.1}
  The  group $\G $ acts on the principal affine space ${\cal A} $ 
  from the left, and the Cartan group $\H $ acts from the right.  
  The arising  
  representation  of the group $\G(\R)$ in an appropriate space of   functions on ${\cal A}(\R)$ is decomposed according to  quasicharacters   $\lambda: \H(\R) \to \C^*$. 
 For each  $\lambda$ we get a representation  in the subspace   of the $\lambda$-homogeneous functions. They form   the principal series  representations. 
The  group $\G(\R)$ acts in the  Hilbert space  $L_2({\cal A}(\R);  \mu^{1/2})$ of half-densities on ${\cal A}(\R)$.\footnote{For any manifold $M$, the space of complex valued half-densities   with compact support is a pre Hilbert space with the scalar product $\langle f,g\rangle:= \int_M f\overline  g$. 
Indeed, $f\overline g$ is a density with compact support on $M$. Its completion is denoted by $L_2(M;  \mu^{1/2})$. The real points of the  space ${\cal A}=\G/\U$ form a manifold. So applying the above construction we get the Hilbert space  $L_2({\cal A}(\R);  \mu^{1/2})$.}   
  It is decomposed into the integral of unitary  representations $V_\lambda$:  
  $$
  L_2({\cal A}(\R);  \mu^{1/2}) = \int V_\lambda d\lambda.
    $$
Here $\lambda =   \rho+ \alpha$, where   
   $\alpha$ are unitary characters of the group $\H(\R)$.\footnote{The character $\rho$ is
  the eigenvalue of the group $\H(\R)$ acting on the square root of   the    volume form  on  ${\cal A}$.} 
    The representations $V_\lambda$ form the unitary principal series of representations of $\G(\R)$.  
There is a unitary action of the Weyl group  by the Gelfand-Graev \cite{GG73}  intertwiners
\be 
\begin{split}
 &  {\cal I}_w: L_2({\cal A}(\R);  \mu^{1/2}) \lra   L_2({\cal A}(\R);  \mu^{1/2}), \qquad w\in W,\\
 &{\cal I}_w: V_{\rho + \alpha } \stackrel{\sim}{\lra} V_{  \rho + w(\alpha) }.\\
\end{split}
\ee
So   representations $V_{\rho + w(\alpha)  } $ for different $w\in W$ are equivalent. The intertwiners are a key feature of the principal series of representations. 

Here is an algebraic counterpart of this picture. Denote by ${\rm Diff}_{\cal A}$ the ring of polynomial differential operators on the space ${\cal A}$.\footnote{ 
The   affine closure of   ${\cal A}$ is   singular. One uses Grothendieck's definition of differential operators.}  
The action of the group $\G$ on ${\cal A}$ gives rise to an injective map of the universal enveloping algebra $\U({\g})$ of the Lie algebra ${\g}$ of $\G$ 
$$
\kappa: \U({\g}) \hlra {\rm Diff}_{\cal A}.
$$
The action of $ \U({\g}) $ commutes with the right action of $\H$. Let  ${\rm Diff}_{\cal A}^\H$ be the subring of ${\rm Diff}_{\cal A}$ commuting with the   action of $\H$. 
It  contains the symmetric algebra $S(\mathfrak{h}) = \U(\mathfrak{h})$ of the Cartan Lie algebra $\mathfrak{h}$.

Recall the Weyl group  $\circ-$action    on  functions on $\mathfrak{h}$, given by  $(w\circ f)(x+ \rho):= f(w(x)+ \rho)$. 
Denote by ${\cal Z}_{\U({\g})}$ the center of  $\U({\g})$. Recall the Harich-Chandra  isomorphism  
\be \la{HC}
{\rm HC}: {\cal Z}_{\U({\g})}\stackrel{\sim}{\lra}  S(\mathfrak{h})^W.
\ee 
It is characterized by the property that    the action of an element $Z \in {\cal Z}_{\U({\g})}$ on ${\cal A}$ coincides with the action of the element 
${\rm HC}(Z)\in S(\mathfrak{h})^W$ on ${\cal A}$. 
There is an isomorphism, where the $\otimes$-product  is  for the  embedding  ${\cal Z}_{\U({\g})}\hra \U({\g})$ and   the Harich-Chandra map: 
\be \la{IUD}
 \U({\g}) \otimes_{{\cal Z}_{\U({\g})}} S(\mathfrak{h})  = {\rm Diff}_{\cal A}^\H.
\ee
So    we get a non-commutative Galois extension 
$\U({\g})\subset {\rm Diff}_{\cal A}^\H$ with the Galois group $W$. 
The action $\ast$ of the  group $W$ on  ${\rm Diff}_{\cal A}^\H$ is induced by the   intertwining operators ${\cal I}_w$:
$$
{\cal I}_w(Df) = (w\ast D){\cal I}_w(f), \qquad \forall D \in {\rm Diff}_{\cal A}.
$$ 
Under the isomorphism (\ref{IUD}) it corresponds to 
 the $\circ$-action  of $W$ on $S(\mathfrak{h})$. So  we have:
\be \la{IUD1}
\U(\g)    = {\rm Diff}_{\cal A}^{\H, W}.
\ee

We conclude that the {\it realization} of the principal series representations in functions on ${\cal A}$ reveals a much larger algebra ${\rm Diff}_{\cal A}$ acting on the principal series as a whole.

\vskip 2mm \paragraph{\bf Example:  $\G= {\rm SL}_2$.} Then ${\cal A} = {\Bbb A}^2 - \{0\}$ and  ${\rm Diff}_{\cal A}= \C[x_1, x_2, \partial_{x_1}, \partial_{x_2}]$. 
The group ${\mathbb G}_m$ acts by $(x_1, x_2) \lms (hx_1, hx_2)$. The unitary intertwiner acting on functions is given by 
the Fourier transform 
$$
f(x_1, x_2) \lms {\cal I}f(y_1, y_2):= \int f(x_1, x_2){\rm e}^{2\pi i(x_1y_2-x_2y_1)}dx_1dx_2.
$$
The  canonical half-form $|dx_1dx_2|^{1/2}$  allows to identify functions with half-forms. Note that ${\cal I}^2 = {\rm Id}$ by the classical Fourier transform inversion formula. So the intertwiner ${\cal I}$ 
provides the unitary action of the Weyl group of ${\rm SL}_2$ on the principal series.   
 
\vskip 3mm
 We wanted to have an analog of the principal series of unitary representations for the quantum group $\U_q(\mathfrak{g})$, which is not just a collection of representations, 
but captures all the extra structure the unitary principal series has, as described above. 

\subsubsection{What should we mean by an infinite dimensional  representation  of a quantum group?} 
  
  Let us look at the very notion of a group representation. Usually a representation space is constructed using some choices: 
  we have a collection of representations $\rho_c$ in vector spaces $V_c$, parametrized by $c \in {\cal C}$. 
  However usually the group acts transitively on the set ${\cal C}$ of  choices, and given two choices $c_1, c_2$, and any group element $g$ such that $g c_1=c_2$ gives rise to an equivalence 
  of representations  $\rho_{c_1}$ and $\rho_{c_2}$. Therefore we may forget about the choices involved in the construction. Here are two examples.
  
 \vskip 1mm  1. The  principal series of representations of $\G(\R)$ is constructed by decomposing    the induced representation ${\rm Ind}_{\U}^\G \C$ from a maximal unipotent subgroup 
  $\U \subset \G$. Since   any two of them are $\G(\R)$-conjugated, the resulting notion does not depend on the choice of $\U$. 
  
 \vskip 1mm  2. The  space of the   Weil representation of the group $\widetilde {\rm Sp}(V)$ is defined by using a decomposition of a  symplectic vector space $V$ into a 
  sum $V = L_1\oplus L_2$ of   Lagrangian subspaces. Since the group ${\rm Sp}(V)$ acts transitively on  such decompositions, the resulting representations are   equivalent. 
 
  \vskip 2mm
  
 If one wants to catch   representations of a semi-simple Lie group by looking at   representations of its Lie algebra, the only adequate way to do this is to study Harich-Chandra modules, that is representations of pairs $(\mathfrak{g}, K)$, where $K$ is a maximal compact subgroup of the group $\G$,   
 such that each representation of $K$ has finite multiplicity. By the  Casselman-Wallach 
   theorem \cite{C},  \cite{BK},  the category of smooth representations 
  of $\G$ is equivalent to the category of Harish-Chandra modules.  
 
 The category 
  of all representations of the Lie algebra $\mathfrak{g}$ is way too big, and is certainly not a reasonable object to study.   For example,   given an  infinite dimensional 
    representation of 
  $\G$ there are infinitely many inequivalent representations of the Lie algebra $\mathfrak{g}$ associated with it.

  \vskip 2mm
  
   Let us now turn to the very notion of the quantum group and its representations. 
   
   First of all, the quantum group $\U_q(\mathfrak{g})$ is   not a group.       
   So one must include a different type of  extra data   to capture  familiar properties of group representations, and  
  to guarantee that representations defined using different choices are equivalent,   providing such equivalences. Here is what we see.

  \begin{enumerate} 
  
  \item     Drinfeld \cite{D} and Jimbo  \cite{J} defined  $\U_q(\mathfrak{g})$    as a Hopf algebra  to capture
    tensor product  of representations.
              
 \vskip 1mm\item   Faddeev proposed  \cite{Fa1} that   infinite dimensional representations    should be  
   modules  over a modular double of the Hopf algebra $\U_q(\mathfrak{g})$, involving $q$ and $q^\vee$, and rigidifying representations. 
  
 Quantization of cluster Poisson varieties \cite{FG07} suggested\footnote{Indeed, the Hilbert space produced by the   quantization   is a module over the tensor product of the  $q$-deformed   algebra of functions on the  cluster Poisson variety, and the $q^\vee$-deformed algebra for   its cluster Langlands dual.} that   the    modular double should involve the  Langlands duality. 
  So the quantum group modular double   is a Hopf algebra
\be \la{ALGh}
  \mathcal{A}_\hbar(\mathfrak{g}):= \U_q(\mathfrak{g})\otimes_\C\U_{q^\vee}({\mathfrak{g}^\vee}).
\ee

\vskip 1mm \item   The Hopf algebra $\U_q(\mathfrak{g})$ was originally defined via  Chevalley generators. However   Lusztig's  action \cite{L} of the 
  braid group ${\Bbb B}_{\mathfrak{g}}$  by automorphisms of $\U_q(\mathfrak{g})$  does not respect the  generators.     
   A priori it is unclear why it should preserve   equivalence classes of   representations.   
   So a construction of a principal series of representation 
   of   $\mathcal{A}_\hbar(\mathfrak{g})$ is meaningless unless it    provides these equivalences. 
           
\vskip 1mm \item  The principal series of representations of a quantum group   
   should come with  an   action of the Weyl group $W$ by intertwiners, quantizing the Gelfand-Graev intertwiners.  
   
 \vskip 1mm  \item One should have {\it realizations} of both quantum principal series of representations and the space of invariants of their tensor products, which make transparent their relationship with the 
   classical realizations of these spaces. 
 
 \end{enumerate}
 
 Our key point is that the cluster nature and quantization of the  moduli spaces ${\mathscr P}_{\G, \bS}$ provides a principal series of representations of the modular double 
 (\ref{ALGh}) with all these features.

\subsubsection{The principal series of $\ast$-representations for   $\widetilde {\mathscr L}_{\G, \odot}$.} 
Recall 
 the   space  $\widetilde {\mathscr L}_{\G, \odot}$  from Section \ref{2.5.1}.   
Specifying  Theorem \ref{Th1.14} to 
the    $\bS=\odot$ case,   we recall    key   features of  the spaces  ${\mathscr P}_{\G, \odot}$,   $\widetilde {\mathscr L}_{\G, \odot}$ and  ${\mathscr L}_{\G, \odot}$: 
  \be
\la{twist.x.poisson}
\begin{gathered}
 \xymatrix{     
{\mathscr L}_{\G, \odot}  &\ar[l]_{\pi} \widetilde {\mathscr L}_{\G, \odot}  \ar[d]  \ar[r] &  {\mathscr P}_{\G, \odot} \ar[d]^{\mu_{\rm out}}  & \qquad \qquad \qquad \qquad \\
& e \ar[r]& \H &\qquad \qquad \qquad \qquad \\}
\end{gathered}
\begin{gathered}
 \xymatrix{   
         \H  &  \ar[l]_{\mu_p} {\mathscr L}_{\G, \odot} \ar[r]^{\rho_{s}} \ar[d]^{{\cal W}_{t}}_{{\cal W}_{s}} &\H\\
          &{\Bbb A}^{r} \oplus{\Bbb A}^{r}&}
\end{gathered}
\ee  
 
 \begin{itemize}

\vskip 1mm  \item The  map $\pi: \widetilde {\mathscr L}_{\G, \odot} \lra {\mathscr L}_{\G, \odot}$      is a Galois cover at the generic point with the Galois group $W$.
    
\vskip 1mm\item There are two projections of the space ${\mathscr L}_{\G, \odot}$ onto the Cartan group, given by the monodromy around the puncture $p$ and by the  $\H$-invariant (\ref{HSi}) of the two  decorated flags 
  near  the    point $s$:
 \be \la{PII}
\begin{split}
&\mu_p:  {\mathscr L}_{\G, \odot}\lra \H /W.\\
&\rho_{s}: {\mathscr L}_{\G, \odot  }\lra \H.\\
\end{split}
  \ee
\end{itemize}   

The cluster quantization, applied to the  space $\widetilde {\mathscr L}_{\G, \odot}$, provides  the principal series of representations of 
the     modular double $\ast-$algebra 
$$
{\cal A}_{\hbar}(\widetilde {\mathscr L}_{\G, \odot}):= {\cal O}_q(\widetilde {\mathscr L}_{\G, \odot})\otimes {\cal O}_{q^\vee}(\widetilde {\mathscr L}_{\G^\vee, \odot}).
$$
 By Theorem \ref{Th8.5A}   the Weyl group    acts by  
automorphisms of the  $\ast$-algebra ${\cal O}_q(\widetilde {\mathscr L}_{\G, \odot})$,  and hence the $\ast$-algebra ${\cal A}_{\hbar}(\widetilde {\mathscr L}_{\G, \odot})$ for each of the two $\ast$-structures on it.  
  Set
\be \la{QLSP}
\begin{split}
&{\cal O}_q({\mathscr L}_{\G, \odot}):={\cal O}_q(\widetilde {\mathscr L}_{\G, \odot})^W,\\
&\mathcal{A}_\hbar({\mathscr L}_{\G, \odot}):= {\cal A}_{\hbar}(\widetilde {\mathscr L}_{\G, \odot})^W.\\
\end{split}
\ee
It follows that each of the the algebras (\ref{QLSP}) inherit the $\ast$-algebra structure(s).

\bt \la{Th2.20}  Let $\G$ be a split semi-simple algebraic group over $\Q$. Then:  

\begin{enumerate}
\vskip 1mm \item   The   representation of the $\ast$-algebra ${\cal A}_{\hbar}(\widetilde {\mathscr L}_{\G, \odot})$ in the   Hilbert space ${\cal H}_{\G, \odot}$  
is decomposed   into an integral of $\ast$-representations, paramatrized by the points $\alpha \in \H(\R_{>0})$:
$$
{\cal H}_{\G, \odot} = \int_{\alpha \in \H(\R_{>0})} {\cal H}_{\G, \odot; \alpha}d\alpha.
$$

\vskip 1mm\item  The Weyl group $W$ acts   on  ${\cal H}_{\G, \odot}$ by unitary projective operators ${\rm I}_w$,   where $w\in W$, intertwining the represenation of the $\ast$-algebra 
${\cal A}_{\hbar}(\widetilde {\mathscr L}_{\G, \odot})$ on the subspace ${\cal S}_{\G, \odot} \subset {\cal H}_{\G, \odot}$:
\be \la{1.1.19.1}
\begin{split}
&{\cal I}_w: {\cal H}_{\G, \odot; \alpha} \stackrel{}{\lra} {\cal H}_{\G, \odot; w(\alpha)}.\\
&{\cal I}_w\circ a (f) = w(a)\circ {\cal I}_w (f), \qquad \forall a \in  {\cal A}_{\hbar}(\widetilde {\mathscr L}_{\G, \odot}), \qquad\forall f \in {\cal S}_{\G, \odot}.\\
\end{split}
\ee

\vskip 1mm\item The intertwiners ${\cal I}_w$  provide unitary equivalences of    represenations of the $\ast$-algebra $\mathcal{A}_\hbar({\mathscr L}_{\G, \odot})$:
 \be
\begin{split}
&{\cal I}_w: {\cal H}_{\G, \odot; \alpha} \stackrel{\sim}{\lra} {\cal H}_{\G, \odot; w(\alpha)}.\\
\end{split}
\ee

\vskip 1mm\item   The braid group $\B_\G$  acts  
 by automorphisms of  the algebras ${\cal O}_q(\widetilde {\mathscr L}_{\G, \odot})$ and ${\cal O}_q({\mathscr L}_{\G, \odot})$, as well as the  $\ast$-algebras $\mathcal{A}_\hbar(\widetilde {\mathscr L}_{\G, \odot})$ and $\mathcal{A}_\hbar({\mathscr L}_{\G, \odot})$.

 \vskip 1mm\item   The braid group $\B_\G$ acts on  ${\cal H}_{\G, \odot}$ by unitary projective operators 
 intertwining the  action of  ${\mathbb B}_\G$
     on the $\ast$-algebra ${\cal A}_{\hbar}(\widetilde {\mathscr L}_{\G, \odot})$. 
\end{enumerate}

   \et
    
 \begin{proof}

1) Follows from the fact that $\mu^*{\cal O}(\H)$ lies in the center of ${\cal O}_q(\widetilde {\mathscr L}_{\G, \odot})$.

2)  The  projection $\mu: \widetilde {\mathscr L}_{\G, \odot}  \lra \H$ given  by the monodromy around the puncture  intertwines the action of the   group $W$  
  on these spaces. 
 By Theorem \ref{THH6.4}, the  action of   $W$  on the space $\widetilde {\mathscr L}_{\G, \odot} $  is   cluster Poisson. 
So by  Theorem \ref{MTHCQ}, there is a unique unitary projective intertwiner  (\ref{1.1.19.1}).

3) Since by  definition the Weyl group acts trivially on  subalgebra (\ref{QLSP}), the claim follows. 

4) The braid group   acts by generalized cluster Poisson automorphisms of the space $\widetilde {\mathscr L}_{\G, \odot}$. Therefore    
it  acts by automorphisms of the algebra ${\cal O}_q(\widetilde {\mathscr L}_{\G, \odot})$. Since the actions of the braid and Weyl groups commute, 
the  braid group   acts by automorphisms of ${\cal O}_q({\mathscr L}_{\G, \odot})$. The rest is clear. 

 5) Follows from 4).\end{proof}

\subsubsection{Quantum groups via  quantum geometry of  spaces ${\mathscr L}_{\G, \odot}$.}   
 
Recall the $\ast$-algebra structure  on the algebra ${\cal O}_q({\mathscr L}_{\G, \odot})$, induced by the embedding (\ref{QLSP}).   
  For each of the two special points $s, t$, there are  regular functions parametrized by the simple positive roots: 
 $$
   {\rm W}_{s, i},  {\rm W}_{t, i}\in {\cal O}_q({\mathscr L}_{\G, \odot}), ~i \in {\rm I}.
  $$ 
\vskip 2mm
At this point it is handy to introduce a variant of the algebra of quantum functions on   ${\mathscr L}_{\G, \odot}$. Precisely,   there is a quantum deformation 
${\rm F}_q$ of the  isomorphism ${\rm F}$ in (\ref{31}), providing   $\ast-$subalgebras
\be \la{Oprime}
\begin{split}
&{\cal O}'_q({\rm  Loc}_{\G, \odot}):=  {\rm F}_q^*\Bigl({\cal O}_q(\B^-) \otimes  {\cal O}_q(\B^+)\Bigr)  \subset {\cal O}_q({\rm Loc}_{\G, \odot}),\\
&{\cal O}'_q({\mathscr L}_{\G, \odot}):=  \qquad{\cal O}'_q({\rm  Loc}_{\G, \odot}){/  \mu_{\rm out} =e}   \qquad\subset {\cal O}_q({\mathscr L}_{\G, \odot}),\\
&{\cal A}'_q({\mathscr L}_{\G, \odot}):=  \qquad {\cal O}'_q({\rm  Loc}_{\G, \odot}) \otimes_\C  {\cal O}'_{q^\vee}({\mathscr L}_{\G, \odot}).\\
\end{split}
\ee
By Theorem \ref{4}, the   $q \to 1$ specializations of these algebras coincide, e.g. ${\cal O}'({\mathscr L}_{\G, \odot}) = {\cal O}({\mathscr L}_{\G, \odot})$.

It is straightforward to see that the composition map $\circ$ in (\ref{circm}) gives rise to a map of algebras 
\be \la{coprm}
\circ^*: {\cal O}'_q({\mathscr L}_{\G, \odot}) \lra {\cal O}'_q({\mathscr L}_{\G, \odot}) \otimes {\cal O}'_q({\mathscr L}_{\G, \odot}).
\ee

\bl The algebra  ${\cal O}'_q({\mathscr L}_{\G, \odot})$ carries a  Hopf algebra structure with the coproduct map (\ref{coprm}), and the antipode provided by 
inverse map on ${\mathscr L}_{\G, \odot}$ from Section \ref{SEC1.7}.
\el

\begin{proof} The composition map (\ref{circm})  is   associative in the category of cluster Poisson varieties. The  
inverse map on ${\mathscr L}_{\G, \odot}$  is a cluster Poisson map. 
\end{proof}

Denote by ${\cal Z}_{\U_q({\mathfrak{g}})}$ the center of the  algebra $\U_q({\mathfrak{g}})$. 
The algebra  $\U_q({\mathfrak{g}})$ is equipped with a coproduct  
 $\Delta: \U_q({\mathfrak{g}}) \lra \U_q({\mathfrak{g}}) \otimes \U_q({\mathfrak{g}})$, an antipode and a counit, making it into a   Hopf algebra.

\bt \la{Th2.21}  Let $\G$ be an adjoint split semi-simple algebraic group over $\Q$. Then:  

\begin{enumerate}   
\vskip 1mm \item There is a canonical  map   of $\ast$-algebras, defined on the generators (\ref{RESCX}) as follows:  
\be \la{MKa}
\begin{split}
\kappa: ~&\U_q({\g}) \lra {\cal O}_q({\mathscr L}_{\G, \odot});\\
&{\bf E_i} \lms {\rm W}_{s,  i}, \qquad{\bf F_{i}} \lms {\rm W}_{t,  i^*}, \qquad{\bf K_i} \lms {\rm K}_{s, i}.\\
\end{split}
\ee
 
\vskip 1mm \item The coproduct $\Delta$   is induced by   geometric coproduct (\ref{coprm}) i.e. there is     a commutative diagram
 
  \begin{displaymath}
    \xymatrix{
        \U_q({\g})  \ar[r]^{\Delta\qquad } \ar[d]_{\kappa} &\U_q({\g}) \otimes \U_q({\g})       \ar[d]^{\kappa \otimes \kappa} \\
          {\cal O}'_q({\mathscr L}_{\G, \odot}) \ar[r]^{\circ^*  \qquad} &    {\cal O}'_q({\mathscr L}_{\G, \odot}) \otimes {\cal O}'_q({\mathscr L}_{\G, \odot})}
         \end{displaymath}
 
     \item The antipode and  counit in $\U_q({\g})$ are induced by  map (\ref{inv}) and the counit map on ${\cal O}_q({\mathscr L}_{\G, \odot}) $.
      
\vskip 1mm\item  Assuming that  $\G$ is simply-laced, the map $\kappa$ intertwines    Lusztig's  braid group   action   on $\U_q({\mathfrak{g}})$ with the      braid group action   on ${\cal O}_q({\mathscr L}_{\G, \odot})$, and it is injective. 
\end{enumerate}
\et

The braid group action was defined by Lusztig in \cite{L6} for simply-laced groups $\G$, and in \cite{L7} in general, see also  \cite{L}. 
Another approach is due to  Soibelman \cite{So}, see also \cite{LS}, \cite{LS1}. 
  We prove the comparison 4)  for simply laced groups just to simplify the exposition: it can be done in general. 
In our applications  we  use only the geometric cluster braid group action. 
 Theorem \ref{Th2.21}   proved   in Section \ref{SSECC8}.   
 
 \bt \la{UQ1} \cite{S22}  If $\mathfrak{g}$ is simply-laced, then the map $\kappa$ provides an isomorphism 
 \be
 \kappa: ~\U_q({\g}) \stackrel{=}{\lra} {\cal O}_q({\mathscr L}_{\G, \odot}).
 \ee
 \et

 \bcon \la{UQ} The map $\kappa$ is an isomorphism for any $\mathfrak{g}$.
  \econ
 
 The classical counterparts of Theorem \ref{UQ1} and Conjecture \ref{UQ} were proved in \cite{S20}.

\subsubsection{Geometric R-matrix.} Recall    the   mapping class group element  $\beta$ for the 
twice punctured disc with two special  points,  interchanging the two punctures, Figure \ref{pin25}. 
It  is  a cluster Poisson transformation by Theorem \ref{MTH}. 
So  it provides a quantum cluster transformation, the geometric analog of the R-matrix: 
\be
{\cal R}: {\cal O}_q({\mathscr L}_{\G, \odot}) \otimes {\cal O}_q({\mathscr L}_{\G, \odot}) \lra {\cal O}_q({\mathscr L}_{\G, \odot}) \otimes {\cal O}_q({\mathscr L}_{\G, \odot}).
\ee
The element $\beta$ coincides with the action of the braid group element $\mu(w_0) \in {\mathbb B}_\G$, provided by Theorem \ref{TTHH2x}. 
So, given a reduced word ${\bf i} = (i_1, \ldots , i_N)$ for $w_0$, we get a decomposition 
of  ${\cal R}$:
$$
{\cal R} = \mu(s_{i_1}) \circ \ldots \circ \mu(s_{i_N}).
$$
It is a geometric version of the factorisation of  
  the ${\cal R}-$matrix using the q-Weyl group     
 \cite{KR},  \cite{So1},  \cite{LS1}.

\subsubsection{Historical comments.} Infinite dimensional representations of $\U_q(sl_2)$ were studied in \cite{Sch}.
 
 The quantum modular double of $\U_q(sl_2)$ and its realization in a quantum torus algebra together with   
   series of its $\ast-$representations were introduced by L.D. Faddeev  \cite{Fa2}, see also 
   \cite{Fa1}.
   
   This construction was developed and studied further in \cite{PT1}, \cite{PT2},     \cite{BT},  \cite{KLS}. 
    
    In  \cite{GKL},  
       an embedding of  $\U_q(sl_n)$    to a quantum torus algebra was constructed. Its  generalisation for  $\U_q(\mathfrak{g})$ was proposed, without proof,   in \cite{GKL2}.   Another 
      embedding  to a quantum torus algebra for  $\U_q(sl_m)$ was constructed in \cite{FI}.

 Cluster Poisson varieties  were introduced and quantized in \cite{FG03b}, 
 and  used to quantize moduli spaces of ${\rm PGL_m}-$local systems \cite{FG03a}. Principal series of $\ast-$representations of
  quantized cluster Poisson varieties were constructed in \cite{FG07}. 
 The cluster  amalgamation technique was introduced    to    establish   cluster Poisson nature of Poisson Lie groups in \cite{FG05}. 
Based on this, V. Fock and the first author suggested that   $\U_q(\mathfrak{g})$  
is  the q-deformed algebra of $W-$invariant regular functions on the  reduced moduli space of $\G-$local systems on a punctured disc with two special points $\odot$. 
    
    Following this suggestion, A. Shapiro and G. Schrader \cite{SS1} found a beautiful realization of $\U_q(sl_m)$ in the quantum torus algebra provided by 
    the special cluster Poisson coordinate system \cite{FG03a}, assigned to a punctured disc with two special points.  
 

 I. Ip constructed in   \cite{II2} a quiver $\mathfrak{D}_\mathfrak{g}^{\bf i}$ for every reduced decomposition ${\bf i}$ of $w_0$, which is mutation equivalent to our quivers related to $\mathscr{P}_{\G, \odot}$. Theorem 4.14 of \cite{II2}  asserts that there exists an embedding $\xi_{\bf i}$ of the algebra $\mathcal{U}_q(\mathfrak{g})$ into a quotient of the quantum torus algebra $\mathcal{D}_{\mathfrak{g}}^{\bf i}$ associated to the quiver $\mathfrak{D}_\mathfrak{g}^{\bf i}$, and the embeddings for different reduced  decompositions 
 ${\bf i}, {\bf i}'$ are related by  quantum cluster ${\mathscr X}-$transformations. However, the proof of Theorem 4.14 is incomplete since both claims  - the construction of   $\xi_{\bf i}$ and its properties - crucially use  the assumption that the quantum cluster transformation from $\mathcal{D}_{\mathfrak{g}}^{\bf i}$ to   $\mathcal{D}_{\mathfrak{g}}^{{\bf i}'}$   does not depend on  the sequence of mutations from ${\bf i}$ to ${\bf i}'$. Unfortunately, this assumption has not been proved or given any references in \cite{II2}. The   proof of the similar crucial Theorem 5.1 in   \cite{II1} is incorrect.

\subsubsection{The principal series of  representations of the modular double of  $\U_q({\g})$.}  
Combining Theorems \ref{Th2.20} and \ref{Th2.21} we  
arrive at    the principal series ${\cal V}_{\lambda}$, $\lambda \in \H(\R_{>0})$  of representations of the modular double $\mathcal{A}_\hbar(\frak{g})$ of the quantum group  $\U_q({\frak{g}})$. Namely, we  take the map
\be \la{poiop}
\kappa\otimes   \kappa^\vee: \mathcal{A}_\hbar(\frak{g})=  \U_q({\g}) \otimes_\C \U_{q^\vee}(\mathfrak{g}^\vee) \hlra \mathcal{A}_\hbar({\mathscr L}_{\G, \odot})
\ee
and compose it with the principal series  of  representations of $\ast$-algebra $\mathcal{A}_\hbar({\mathscr L}_{\G, \odot})$.

 \bc
 For any $\lambda\in \H(\R_{>0})$,  and any $w \in W$, the principal series  representations ${\cal V}_{w(\lambda)}$   of the quantum group modular double $\ast$-algebra 
 $ \mathcal{A}_\hbar({\mathfrak{g}})$ 
 are canonically isomorphic. 
 \ec
 
 \begin{proof}  Follows from Theorem \ref{Th2.20}, part 3). 
   \end{proof} 
   
 At the moment,   the classical Gelfand-Naimark {realization} 
of the  principal representation of the group $\G(\R)$ in ${\rm L}_2(\G/\U(\R))$ on the principal affine space seems to be  unrelated 
to any  {realization} of its  quantum counterpart -  the representation of the subalgebra 
$$
\mathcal{A}_\hbar(\mathfrak{g}) \subset  \mathcal{A}_\hbar({\mathscr L}_{\G, \odot})
$$
in the cluster Schwartz space ${\cal S}({\mathscr L}_{\G, \odot})$ of 
 the principal   representation of the  $\ast$-algebra  $\mathcal{A}_\hbar({\mathscr L}_{\G, \odot})$. 
 
 In Section \ref{Section 5.4.5} we explain how to connect the two constructing
    a more geometric realization of the algebra $\mathcal{A}_\hbar({\mathscr L}_{\G, \odot})$ in the cluster Hilbert space related 
 to ${\rm L}_2(\G/\U(\R_{>0}))$. We elaborate this further in Section \ref{SEC9.4}. The main construction  uses crucially the cluster structure on  $\G$.

 \subsubsection{A quantum analog of the regular representation of the group $\G(\R)$.} \la{4.4.8}
   Denote by ${\rm Cl}_{2,2}$ a cylinder with an ordered pair of special points on each of the two of its boundary components. The outer monodromies around the boudary components 
   provide a map  
   \be
  \mu_{\rm out}: {\mathscr P}_{\G, {\rm Cl}_{2,2}}   \lra \H \times \H.
   \ee
   Denote by   ${\cal R}_{\G, {\rm Cl}_{2,2}}$ its fiber over the unit $e\in \H \times \H$.  Part 3) of Theorem \ref{UEAB} provides a  map 
\be \la{ewqwe}
   {\cal A}_\hbar(\mathfrak{g}) \otimes_\C {\cal A}_\hbar(\mathfrak{g}) \lra {\cal A}_\hbar({\cal R}_{\G, {\rm Cl}_{2,2}} ).
\ee
Cutting the cylinder ${\rm Cl}_{2,2}$ by a loop $\alpha$, see the Figure \ref{pin100}, we get a map of  moduli spaces
$$
{\rm Res}_\alpha: {\cal R}_{\G, {\rm Cl}_{2,2}} \lra {\mathscr L}_{\G, \odot}  \times {\mathscr L}_{\G, \odot}.
$$   
 \begin{figure}[ht]
 \begin{center}
 \begin{tikzpicture}[scale=0.4]
 \draw (0,0) ellipse (0.5 and 1.25);
\draw (0, -1.25) -- (5, -1.25);
\draw (5, -1.25) arc (-90:90:0.5 and 1.25);
\draw [dashed] (5, -1.25) arc (270:90:0.5 and 1.25);
\draw [blue] (2.5, -1.25) arc (-90:90:0.5 and 1.25);
\draw [blue, dashed] (2.5, -1.25) arc (270:90:0.5 and 1.25);
\draw (5, 1.25) -- (0, 1.25);  
\node [red] at (-0.5,0) {\small $\bullet$};
\node [red] at (0.5,0) {\small $\bullet$};
\node [red] at (5.5,0) {\small $\bullet$};
\node [red] at (4.5,0) {\small $\bullet$};
\node [blue] at (2.5,0) {\small $\alpha$};
 \draw (13,0) circle (1.25);
  \draw (17,0) circle (1.25);
  \node [red] at (11.75,0) {\small $\bullet$};
\node [red] at (14.25,0) {\small $\bullet$};
\node [red] at (15.75,0) {\small $\bullet$};
\node [red] at (18.25,0) {\small $\bullet$};
\node [blue, thick] at (13,0) {$\circ$};
\node [blue, thick] at (17,0) {$\circ$};
\draw [-latex, thick] (8,0) -- (10, 0);
\end{tikzpicture}
\end{center}
\caption{Cutting a cylinder with $2+2$ special points by a loop $\alpha$.}
\label{pin100}
\end{figure}
    Modular Functor Conjecture \ref{MK} predicts an isomorphism of Hilbert spaces, compatible with the action of the corresponding $\ast-$algebras on the Schwarz spaces:
   \be \la{GRDa}
 {\cal H}({\cal R}_{\G, {\rm Cl}_{2,2}}) = \int_{\lambda\in \H(\R_{>0})}  {\cal H}({\cal R}_{\G, \odot})_\lambda  \otimes {\cal H}({\cal R}_{\G, \odot})_{\lambda^{-1}}d\lambda.
 \ee        
 Precisely,  the action of the quantum group modular double $\mathcal{A}_\hbar(\mathfrak{g}) \otimes  \mathcal{A}_\hbar(\mathfrak{g}) $ on   ${\cal S}({\cal R}_{\G, {\rm Cl}_{2,2}})$ provided by 
 (\ref{ewqwe}) should be intertwined 
   with  its action on  ${\cal S}({\cal R}_{\G, \odot})_\lambda  \otimes {\cal S}({\cal R}_{\G, \odot})_{\lambda^{-1}}$ provided by (\ref{poiop}).
     
The decomposition (\ref{GRDa}) is a quantum analog of the decomposition  
 \be
 {\rm L}_2(\G(\R)) = \int_{\chi}  {\cal V}_\chi  \otimes {\cal V}_{\chi^{*}}d\chi 
  \ee
of the regular representation of  the group $\G(\R) \times \G(\R)$  via  unitary representations ${\cal V}_\chi$  of  $\G(\R)$.    
Yet to fully establish this analogy we have to relate the    regular representation of the group $\G(\R) \times \G(\R)$ in ${\rm L}_2(\G(\R))$ to 
  its  quantum counterpart -  the representation of the algebra 
$
\mathcal{A}_\hbar(\mathfrak{g}) \otimes  \mathcal{A}_\hbar(\mathfrak{g})
$ 
in the cluster Schwarz space ${\cal S}({\cal R}_{\G, {\rm Cl}_{2,2}})$, provided by the map (\ref{ewqwe}) followed by 
 the principal series representation of ${\cal A}_\hbar({\cal R}_{\G, {\rm Cl}_{2,2}})$. 
 
 In Section \ref{Section 5.4.5} we  connect the quantum and classical realizations by constructing
    a   realization of the algebra ${\cal A}_\hbar({\cal R}_{\G, {\rm Cl}_{2,2}})$ in the cluster Hilbert space related 
 to ${\rm L}_2(\G(\R_{>0}))$. 
We  use essentially the cluster structure on  $\G$.   The details are elaborated in Section \ref{SEC16.8}, 
using crucially the results of Section \ref{SECTa17.8}.

  \subsection{The  TQFT of    representations of   enhanced algebras ${\cal O}_q({\mathscr P}_{\G, \bS})$}  \la{5.2}

\medskip

In Section \ref{5.2} we explain how the category of  finite dimensional   
class one representations of the quantum universal enveloping algebra $\U_q(\mathfrak{g})$ could be realized as a genus zero sector of an  
algebraic-geometric TQFT. 

This  TQFT is given by a   class of  ${\cal O}_q({\mathscr P}_{\G, \bS})-$modules $V_{\G, \bS; \lambda}$,  {\it enhanced by the action of a  product  ${\rm H}^{\rm c_\bS}$ of the Cartan groups}. 
The vector spaces  $V_{\G, \bS; \lambda}$ are finite-dimensional if and only if the genus of $\bS$ is  $0$.
They   
should provide a traditional TQFT, which  is nothing else but the classical Wess-Zumino-Witten TQFT.    

\subsubsection{The enhanced quantum group.} Traditionally,  the quantum group   is understood as a Hopf algebra    $\U_q(\mathfrak{g})$.  
It is   determined by the   braided monoidal  category  of its finite dimensional representations. 

However this perspective has the following issue. 

There are   $2^r$ non-isomorphic one dimensional representations of $\U_q(\mathfrak{g})$, given by 
\be
{\bf E}_i, {\bf F}_i \lms 0, \quad {\bf K}_i \lms \pm 1. 
\ee

On the other hand, the category of all finite dimensional representations  of $\U_q(\mathfrak{g})$ contains the  subcategory $\U_q(\mathfrak{g})-{\rm mod}$ of the  class one  representations. 
They are deformations of  
finite dimensional representations of $\mathfrak{g}$. The  simple    class one representations are   parametrised by the dominant weights. 
Any irreducible  finite dimensional representation of $\U_q(\mathfrak{g})$ is a tensor product  of a class one representation and a  
  one dimensional representation. 

The true  category of finite dimensional representations of the quantum group  should be the category of class one representations, unless  $q$ is an even  root of unity.  
This can be achieved by modifying   the very notion of the quantum group: \vskip 1mm

{\it The quantum group is not just the Hopf algebra $\U_q(\mathfrak{g})$, but rather 
a Harish-Chandra like data} 
\be \la{dataxxc} 
\Bigl(\U_q(\mathfrak{g}), \rm H, \tau \Bigr).
\ee
\vskip 1mm

 Here  $\tau$ denotes an action of the   Cartan group $\rm H$   by automorphisms of the algebra $\U_q(\mathfrak{g})$ such  that, using the standard generators, normalized  as  in  Section \ref{SSEECC11.3}, we have:
\be \la{wtEFK}
\tau(h) ({\bf E}_i) = \alpha_i(h) {\bf E}_i, \quad \tau(h) ({\bf F}_i) = \alpha^{-1}_i(h) {\bf F}_i, \quad \tau(h) ({\bf K}_i )=  {\bf K}_i. 
\ee
We refer to such a data as an {\it enhanced quantum group}. Its representations are defined as follows.

  \bd \la{D5.10} A finite dimensional representation of    
  the enhanced quantum group  (\ref{dataxxc}) is  given by
  
  \begin{enumerate}
  
\vskip 1mm  \item  A finite dimensional representation $\rho_V$ of the algebra $\U_q(\mathfrak{g})$ in a vector space $V$, together with   
  
 \vskip 1mm \item An action $\rho_{\rm H}$ of the Cartan group ${\rm H}$ on $V$, compatible   as follows:
  
  \end{enumerate}

\begin{itemize}

\vskip 1mm\item One has 
 \be \la{qgd}
 \begin{split}
 & \rho_V\circ \tau(h)(A) =   \rho_{\rm H}(h) \rho_V(A)  \rho_{\rm H}(h)^{-1} \qquad \forall A \in \U_q(\mathfrak{g}).  \\
 \end{split}
 \ee

\vskip 1mm\item  For each simple 
  coroot  $\alpha^\vee_i:{\Bbb G}_m \hra {\rm H}$ we have 
\be
\rho_V({\bf K}_{i}) = \tau(\alpha_i^\vee(q)) \in {\rm Aut}(V).
\ee 
\end{itemize}
\ed
We call such representations {\it enhanced representations} of the quantum group $\U_q(\mathfrak{g})$. \vskip 1mm

In particular, take a one dimensional enhanced 
representation of $\U_q(\mathfrak{g})$. Then the Cartan group $\rm H$ acts by a charachter $\chi$,   the elements ${\bf E}_i, {\bf F}_i$ must act by zero since by (\ref{wtEFK}) 
they have non-trivial weights $\alpha^{\pm 1}_i$, and  elements  
${\bf K}_i$ act by multiplication by $\chi\circ \alpha^\vee_i(q)$. Since   
$$
0=[{\bf E}_i, {\bf F}_j] = \delta_{ij}\frac{{\bf K}^{-1}_i-{\bf K}_i}{q_i-q_i^{-1}},
$$
 the character $\chi$ must be trivial, unless $q$ is a  $2n-$th root of unity and $\chi(z) = z^n$.

\vskip 2mm

Let us show how the enhanced quantum group   (\ref{dataxxc})  is build in our quantized moduli spaces. Denote by $c_\bS$ the number of colored boundary intervals on $\bS$.  
We  assigned to a colored decorated surface $\bS$ the algebra ${\cal O}_q({\mathscr P}_{\G, \bS})$. Let us define the category 
of {\it enhanced representations of the algebra} ${\cal O}_q({\mathscr P}_{\G, \bS})$, denoted by 
\be \la{datas}
({\cal O}_q({\mathscr P}_{\G, \bS}), {\rm H}^{c_\bS}, \tau)-{\rm mod}. 
\ee

When  $\bS = \odot$,  reducing the moduli space ${\mathscr P}_{\G, \odot}$ to the one ${\mathscr L}_{\G, \odot}$, and restricting to  finite dimensional representations, 
we get the category of finite dimensional modules over the 
enhanced quantum group   (\ref{dataxxc}). \vskip 1mm


Let us implement this plan.

\subsubsection{Enhanced  representations of the algebra ${\cal O}_q({\mathscr P}_{\G, \bS})$.} 
 For each colored boundary interval ${\rm I}$ on $\bS$, there is an action (\ref{ALa})  of the 
 Cartan group $\rm H$      
  on the moduli space ${\mathscr P}_{\G, \bS}$.   It gives  rise  
to an action $\tau_{\rm I}$ of the group $\rm H$ by automorphisms of  ${\cal O}_q({\mathscr P}_{\G, \bS})$.    
So we get an action $\tau$ of the  algebraic torus 
$
{\rm H}^{c_\bS}
$ 
  by automorphisms of the algebra  ${\cal O}_q({\mathscr P}_{\G, \bS})$.  \vskip 1mm

Recall that for each special point $s$ on $\bS$ shared by two colored boundary intervals there are elements
\be
{\bf K}_{s, \alpha_i}^{\pm 1} \in {\cal O}_q({\mathscr P}_{\G, \bS}).
\ee

\bd  \la{def5.11} The enhanced representations (\ref{datas})  of the algebra 
$ {\cal O}_q({\mathscr P}_{\G, \bS})$ are  given by: 
 \begin{enumerate}
\vskip 1mm \item A    representation $\rho_V$ of the algebra ${\cal O}_q({\mathscr P}_{\G, \bS})$ in a  vector space $V$:
 \be
 \rho_V: {\cal O}_q({\mathscr P}_{\G, \bS}) \lra {\rm End}(V), 
  \ee 
  
\vskip 1mm  \item An   action  $\rho_{\rm H}$ of the algebraic  group $\rm H^{c_\bS}$ in the   space $V$, compatible with the representation $\rho_V$:
 \be
 \begin{split}
 &g_{\rm H}: \H \lra {\rm Aut}(V), \\
 & \rho_V(\tau_h(A)) =   \rho_{\rm H}(h) \rho_V(A)  \rho_{\rm H}(h^{-1}), 
 \\
 \end{split}
 \ee
\end{enumerate}
subject to the following condition:
\begin{itemize}
\vskip 1mm \item For each special point $s$ on $\bS$ shared by two   colored intervals ${\rm I}_s^-, {\rm I}_s^+$, where ${\rm I}_s^+$ follows after $s$, the action of the element 
${\bf K}_{s, \alpha_i}$ on $V$ is related to the action $\tau_{{\rm I}_s^r}$ of the Cartan group  ${\rm H}$   as follows:
\be \la{228a}
\rho_V({\bf K}_{s, \alpha_i}) = \rho_{\rm H}(\tau_{{\rm I}_s^+}(\alpha_i^\vee(q_i))).
\ee

\end{itemize}
 
\ed 

Our next goal is to discuss how the quantum algebra  ${\cal O}_q({\mathscr P}_{\G, \bS})$ behaves under cutting 
of the surface $\bS$. 
\vskip 1mm

\subsubsection{The monoidal category of  class one representations of $\U_q(\mathfrak{g})$ via the moduli spaces ${\mathscr L}_{\G, \odot}$} \la{5.3.3} 

The moduli space ${\mathscr L}_{\G, \odot}$ comes with   
  the canonical action of the Cartan  group $\rm H$ on the pinning $p$:
\be
\gamma: \H \times {\mathscr L}_{\G, \odot} \lra {\mathscr L}_{\G, \odot}.
\ee
It is promoted to the action $\gamma$ of the Cartan group $\H$  by automorphisms of the algebra ${\cal O}_q({\mathscr L}_{\G, \odot})$. \vskip 1mm

According  to Theorem  \ref{Th2.21} and Theorem  \ref{UQ1}/Conjecture \ref{UQ}, there is a canonical injective map
\be  \la{543212}
\kappa: \U_q(\mathfrak{g}) \stackrel{\sim}{\lra} {\cal O}_q({\mathscr L}_{\G, \odot}) 
\ee
which is   an isomorphism in the  simply laced case, and conjecturally     in general. 
It intertwines the actions of the Cartan group $\rm H$ on $\U_q(\mathfrak{g})$ and    ${\cal O}_q({\mathscr L}_{\G, \odot})$, and respects the coproduct.  
So it provides a tensor functor
\be \la{43212}
({\cal O}_q({\mathscr L}_{\G, \odot}), {\rm H})-{\rm mod} \stackrel{\sim}{\lra} (\U_q(\mathfrak{g}), {\rm H})-{\rm mod}
\ee
 which is   an equivalence of tensor categories  in the simply-laced case. Restricting to   finite dimensional representations, the monoidal category on the right is canonically equivalent to    the  
 monoidal category  of finite dimensional class one representations of 
  $\U_q(\mathfrak{g})$.

\subsubsection{The Wess-Zumino-Witten TQFT from   enhanced representations of algebras ${\cal O}_q({\mathscr P}_{\G, \bS})$}

 The center of the algebra ${\cal O}_q({{\cal P}_{\G, \bS}})$ for generic $q$ is identified with 
     the algebra ${\cal O}({\rm C}_{\G, \bS})$ of regular functions ${\cal O}({\rm C}_{\G, \bS})$ on the Casimir torus   ${\rm C}_{\G, \bS}$. 
     It always  
  contains   the product of the   Cartan groups ${\rm H}$ over the   punctures $p$ on $\bS$.    In the case when all boundary segments are colored,
   the Casimir torus ${\rm C}_{\G, \bS}$  is described in Theorem \ref{TH1.16a}.  \vskip 2mm

Let  $\alpha$ be a non-boundary loop  on   $\bS$. Let $\bS_\alpha$ be the surface obtained by cutting $\bS$ by the loop $\alpha$. 
One has 
\be
{\rm C}_{\G, \bS_\alpha} = {\rm C}_{\G, \bS} \times {\rm H}_{\alpha_+} \times  {\rm H}_{\alpha_-}.
\ee
Any dominant coweights  
$$
 \lambda \in X_*({\rm C}_{\G, \bS}), \qquad \mu  \in X_*({\rm H})  
$$
  give rise to the characters of the algebras ${\cal O}({\rm C}_{\G, \bS})$ and ${\cal O}({\rm C}_{\G, \bS_\alpha})$ described by the points  
\be \la{centralchar}
\lambda(q) \in {\rm C}_{\G, \bS}(\C), \qquad ( \lambda, \mu, \mu)(q) \in {\rm C}_{\G, \bS_\alpha}(\C).
\ee
\vskip 1mm
 
In particular, let $\bS$ be a surface without boundary with $n$ punctures.  Then ${\rm C}_{\G, \bS} = {\rm H}^n$ by Theorem \ref{TH1.16a}. So   the center of the algebra ${\cal O}_q({\cal P}_{\G, \bS})$ for generic $q$ 
 is canonically isomorphic to  ${\cal O}({\rm H}^n)$. Its  characters are   points of the complex torus ${\rm H}^n(\C)$.  
Let $g = g(\bS)$ be the genus of $\bS$. The mapping class group $\Gamma_{\bS} = \pi_1({\cal M}_{g,n})$ acts  by cluster transformations  on the algebras ${\cal O}_q({\mathscr P}_{\G, \bS})$ 
and hence   on their finite dimensional  representations. \vskip 2mm

Denote by $V_\lambda$ the finite dimensional class one   $\U_q(\mathfrak{g})-$module with the highest weight $\lambda$. 
Recall  the mulitplicity space   
\be \la{msp}
\Bigl(V_{\lambda_1}\otimes \ldots \otimes V_{\lambda_n}\Bigr)^{\U_q(\mathfrak{g})}.  
\ee
It has a natural action of the  braid group $\Gamma_{S^2_n}$. The Kazhdan-Lusztig theory \cite{KL} identifies  these $\Gamma_{S^2_n}-$modules    
 with the spaces of conformal blocks for the WZW theory.\vskip 2mm
 
 Given a simple non-trivial loop $\alpha$ on $\bS$, Conjecture \ref{QMONC} predicts  the quantum monodromy subalgebra
 $$
 {\rm M}^*_\alpha {\cal O}({\rm H}/W) \subset {\cal O}_q({\cal P}_{\G, \bS}).  
  $$
The algebraic Modular Functor Conjecture \ref{koskin} describes its centralizer in terms of the surface $\bS_\alpha$:   
  \be \la{MFKiso}
   {\cal O}_q({\cal P}_{\G, \bS})^{{\rm M}^*_\alpha {\cal O}({\rm H}/W)} =  \Bigl( {\cal O}_q({\cal P}_{\G, \bS_\alpha})/{\cal I}_\alpha\Bigr)^{W_\Delta}.
      \ee

      For the convenience of the reader, we recall Conjecture \ref{MFK2i} in a sharper form.
      
      \bcon   \la{MFK2}  $\bS $ be a sphere with $n$ punctures, and $q$ is not a root of unity. Then 
       
\begin{enumerate}

\vskip 1mm\item    The  
   category of enhanced finite dimensional    ${\cal O}_q({{\cal P}_{\G, \bS}})-$modules  is semi-simple. 
The 
irreducible modules ${\cal V}_{\G, \bS, \alpha}$ 
are parametrised by   integral coweights  at the punctures:
\be \la{insi}
\alpha   = (\alpha _1, \ldots , \alpha_n)  \in X_*({\rm H}^n).
\ee
The center  ${\cal O}({\rm H}^n)$   acts 
on  the module ${\cal V}_{\G, \bS, \alpha}$ by the character     
$\alpha(q) \in {\rm H}^n(\C)$.  

\noindent
The   $W^n-$action on the algebra ${\cal O}_q({{\cal P}_{\G, \bS}})$ provides  intertwiners 
\be
{\cal I}_w: {\cal V}_{\G, \bS, \alpha} \lra {\cal V}_{\G, \bS, w(\alpha)}, \qquad w \in W^n.
\ee 

\vskip 1mm\item Let $V_{\G, \bS, \alpha}$ be the ${\cal O}_q({{\rm Loc}_{\G, \bS}})-$module given by the restriction of   ${\cal V}_{\G, \bS, \alpha}$ 
to  the subalgebra ${\cal O}_q({{\rm Loc}_{\G, \bS}})$.  
 The simple  enhanced ${\cal O}_q({{\rm Loc}_{\G, \bS}})-$modules  
are the ones $V_{\G, \bS, \lambda}$,    parametrised by     dominant coweights:
\be \la{ins}
\lambda   = (\lambda _1, \ldots , \lambda_n)  \in X^+_*({\rm H}^n).
\ee

\vskip 1mm\item  Given a simple non-trivial loop $\alpha$ on $\bS$,  
 the restriction of    $V_{\G, \bS, \lambda}$     to the    quantum monodromy  centralizer 
 $$
{\cal O}_q({\rm Loc}_{\G, \bS})^{{\rm M}^*_\alpha {\cal O}({\rm H}/W)} \subset {\cal O}_q({\rm Loc}_{\G, \bS})
$$
is isomorphic, using  isomorphism (\ref{MFKiso}), to a direct sum  over the set of  dominant integral coweights  
\be \la{MFCPR}
V_{\G, \bS, \lambda} = \bigoplus_{\mu \in    X^+_*({\rm H})} V_{\G, \bS_\alpha, \lambda; \mu} 
\ee
 of the ${\cal O}_q({\rm Loc}_{\G, \bS_\alpha})-$modules  $V_{\bS_\alpha, \lambda; \mu}$    
   with the  central  
 character    
$$
(\lambda, \mu, \mu)(q) \in {\rm C}_{\G, \bS_\alpha}(\C)/(W^n\times W_\Delta).
$$

 \vskip 1mm\item  There is an isomorphism of $\Gamma_{\bS}-$modules 
\be \la{mspi}
V_{\G, \bS,   \lambda_1, ..., \lambda_n} = \Bigl(V_{\lambda_1}\otimes \ldots \otimes V_{\lambda_n}\Bigr)^{\U_q(\mathfrak{g})}. 
\ee 
\end{enumerate}
 \econ

The next Conjecture complements this picture,    extending it to the higher genus $g(\bS)>0$ case. 

\bcon   \la{5.13} 
For any hyperbolic   surface $\bS$  with $n$ punctures, there are enhanced     
  ${\cal O}_q({\rm Loc}_{\G, \bS})-$modules $V_{\G, \bS, \lambda}$, parametrized   by  $\lambda \in X_*({\rm C}_{\G, \bS})/W^n$ such that:

\begin{enumerate} \item 
They satisfy the Modular Functor Conjecture property (\ref{MFCPR}).

\vskip 1mm\item There are similar   ${\cal O}_q({\mathscr P}_{\G, \bS})-$modules $V_{\G, \bS, \alpha}$, where  $\alpha \in X_*({\rm C}_{\G, \bS})$, satisfying formula (\ref{MFCPR}).

\vskip 1mm\item When $\bS$ is a surface without boundary, but possibly with punctures, the representation $V_{\G, \bS, \lambda}$ is isomorphic, as a $\Gamma_{\bS}-$modules, to the 
space of conformal blocks for the WZW theory related to $(\G, \bS)$.

\end{enumerate}
\econ

Few comments are in order.

\begin{enumerate}

\item  The new feature of our construction   is that  the    spaces $ V_{\G, \bS, \{\lambda_p\}}$ of  WZW conformal blocks    
carry an action of   the algebra ${\cal O}_q({\rm Loc}_{\G, \bS})$.  
For the   multiplicity   space  (\ref{msp}) it is the algebra ${\cal O}_q({\rm Loc}_{\G, S^2_n})$.

\vskip 1mm\item If $g(\bS)>0$,    non-trivial      ${\cal O}_q({\rm Loc}_{\G, \bS})-$modules are infinite dimensional. This is clear from the property (\ref{MFCPR}) since the sum 
there will be infinite. 

\vskip 1mm\item  As   explained in Section \ref{5.3.3},  we recover the class one representations of $\U_q(\mathfrak{g})$ from the colored decorated surface $\odot$. The standard argument 
 using Modular Functor conjecture isomorphism (\ref{MFCPR}) reduces understanding of the spaces $V_{\G, \bS, \lambda}$ to the case when $\bS$ is a pair of pants. In this case   (\ref{msp}) follows from Conjecture \ref{MFK2} using   isomorphism   (\ref{MFCPR}) and Figure \ref{pin5a}.  
  \end{enumerate}

\subsubsection{The algebraic-geometric avatar   ${\cal R}_{\cal F}$ of a TQFT} 
Consider the category of colored decorated surfaces, with the morphisms  generated by   gluings of  pairs of colored boundary intervals: 
\be \la{Mor2}
\gamma: (\bS, {\rm I}_1, {\rm I}_2) \stackrel{}{\lra} {\bS'}:= {\bS}/({\rm I}_1\sim {\rm I}_2), 
\ee 
and  embeddings of colored decorated surfaces considered modulo isotopies:
\be \la{Mor1}
i:\bS \stackrel{}{\hra} \bS'.
\ee

Denote by ${\rm H}_{{\rm I}}$  the Cartan group assigned to a colored 
boundary interval $\rm I$.  There is the diagonal subgroup
\be \la{DS}
{\rm H}_\Delta \subset {\rm H}_{{\rm I}_1} \times {\rm H}_{{\rm I}_2}\subset {\rm H}^{c_\bS}.
\ee
A gluing map $\gamma$  in (\ref{Mor2}) induces an embedding
\be \la{NEM}
{\rm H}^{c_{\bS'}} \times {\rm H}_\Delta \hra {\rm H}^{c_\bS}.
\ee

\bl \la{mor3} The gluing map $\gamma$ in (\ref{Mor2}) induces an injection  from ${\cal O}_q({\mathscr P}_{\G, \bS'})$ to 
the ${{\rm H}_\Delta}-$invariants    on   ${\cal O}_q({\mathscr P}_{\G, \bS})$:
\be \la{240} 
\pi^*_{\bS \to \bS'}: {\cal O}_q({\mathscr P}_{\G, \bS'}) \stackrel{}{\hra} {\cal O}_q({\mathscr P}_{\G, \bS})^{{\rm H}_\Delta}. 
\ee 
\el

\begin{proof}  The gluing map is a regular map of moduli spaces. From the cluster point of view, it 
amounts to   the amalgamation of the corresponding cluster Poisson varieties. The latter is described by taking the invariants 
of the torus provided by the  frozen variables glued under the amalgamation. In our case this is the Cartan group ${\rm H}_\Delta$. The clusters after the amalgamation contain the ones before as a subset, since after the amalgamation we allow some frozen variables to mutate.   Therefore the gluing gives rise to  a map of algebras of cluster regular functions. It is evidently injective.  
\end{proof}

The map (\ref{240}) may not  be an isomorphism, since the algebra on the right have frozen variables inverted, while on the left this may not be the case. 

\bcon \la{QMONCa}   
 Given an embedding $i$ in (\ref{Mor1}), there exists a unique map of algebras 
\be
{\cal O}_q({\mathscr P}_{\G, \bS'}) \lra {\cal O}_q({\mathscr P}_{\G, \bS}) 
\ee
which 
 commutes with the quantum monodromy maps for all loops $\alpha$ on $\bS$, and  
   with the gluing maps. 
\econ

Our goal is to define a   functor ${\cal R}_{\cal F}$ on the category of colored decorated surfaces $\bS$, which assigns to an $\bS$ the category   from Definition \ref{def5.11}:
\be \la{RF}
\begin{split}
 &{\cal R}_{\cal F}: \left\{\mbox{\rm Colored Decorated Surfaces} \right\}\lra {\rm Cat},\\
&{\cal R}_{\cal F}(\bS) = ({\cal O}_q({\mathscr P}_{\G, \bS}), {\rm H}^{c_\bS}, \tau)-{\rm mod}.\\
\end{split}
\ee

\bt  The functor ${\cal R}_{\cal F}$  in (\ref{RF}) 
  is functorial for the automorphisms    and   gluings $\gamma$ in (\ref{Mor2}). 
Assuming Conjecture \ref{QMONC}, it is functorial for  embeddings $i$ in (\ref{Mor1}). 
\et

\begin{proof} 
1. {\it Functoriality for   automorphisms of    $\bS$}. The action of the mapping class group $\Gamma_\bS$   on   ${\mathscr P}_{\G, \bS}$ gives rise to its action  by 
automorphisms of the algebra ${\cal O}_q({\mathscr P}_{\G, \bS'})$ and the torus ${\rm H}^{c_{\bS'}}$. \vskip 1mm

2. {\it Functoriality for the gluing maps $\gamma$}. Given a gluing map (\ref{Mor2}), let us define the functor 
\be \la{224}
{\cal R}_{\cal F}(\bS \stackrel{\gamma}{\lra} {\bS'}): ({\cal O}_q({\mathscr P}_{\G, \bS}), {\rm H}^{c_{\bS}}, \tau_\bS)-{\rm mod} ~~ \lra ~~({\cal O}_q({\mathscr P}_{\G, {\bS'}}), {\rm H}^{c_{\bS'}}, \tau_{\bS'})^{{\rm H}_\Delta}-{\rm mod}
\ee
by taking the invariants of the  subgroup ${\rm H}_\Delta$:
 \be
 {\cal R}_{\cal F}(\bS \stackrel{\gamma}{\lra} {\bS'}): V_{\G, \bS} \lms  V_{\G, {\bS}}^{{\rm H}_\Delta}.
 \ee
The group ${\rm H}^{c_{\bS'}}$ acts on   $V_{\G, {\bS}}$ via    embedding (\ref{NEM}). 
So   it acts on the space of invariants $V_{\G, {\bS}}^{{\rm H}_\Delta}$. The compatibility condition (\ref{228a}) is evident.  Thanks to  Lemma \ref{mor3}, the algebra 
${\cal O}_q({\mathscr P}_{\G, \bS'})$ acts on the   ${\rm H}_\Delta-$invariants.  \vskip 1mm

3. {\it Functoriality for restrictions}. Given an embedding of colored decorated surface (\ref{Mor1}), we need  a functor
\be
{\cal R}_{\cal F}(\bS \stackrel{i}{\hra} \bS'): ({\cal O}_q({\mathscr P}_{\G, \bS}, \tau_{\rm H}), {\rm H}^{c_{\bS}}, \tau_\bS)-{\rm mod} ~~\lra ~~({\cal O}_q({\mathscr P}_{\G, {\bS'}}), {\rm H}^{c_{\bS'}}, \tau_{\bS'})-{\rm mod}.
\ee
Below we skip $\tau$ in the notation. For this we need to define 
  a map of the  enhanced quantum algebras:
\be \la{MENA}
i_q^*: ({\cal O}_q({\mathscr P}_{\G, \bS'}), {\rm H}^{c_{\bS'}}) \lra ({\cal O}_q({\mathscr P}_{\G, \bS}), {\rm H}^{c_{\bS}}).
\ee 
Its  existence   is exactly what   Conjecture \ref{QMONCa} predicts. 
\end{proof}

\medskip
\subsection{Canonical realisations and quasiclassical limits} \la{SC2.7}

\medskip

\subsubsection{Relating the quantum and classical representation theory.} 
There are several important vector spaces studied  in the  representation theory of Lie groups: 

\vskip 1mm 1. The principal series representation  ${\cal V}$ 
  of  the  group $\G(\R)$ is realized in   $ {\rm L_2}({\cal A}_\G(\R), \mu_{\cal A})$, where  ${\cal A}_\G= \G/\U$.    
 Its  $n$-th tensor power  ${\cal V}^{\otimes n}$ is realized  in 
  $  {\rm L_2}({\cal A}_\G^n(\R), \mu_{\cal A})$. 
  
 \vskip 1mm 2.The subspace of $\G(\R)$-invariants in the $n$-th tensor power  ${\cal V}^{\otimes n}$  is given by 
  \be \la{MSPX}
   {\rm L_2}({\cal A}_\G^n(\R), \mu_{{\cal A}_\G})^{\G(\R)} =  {\rm L_2}({\rm Conf_n}({\cal A}_\G), \mu_{\cal A}). 
\ee
It can   be identified with   
${\rm Hom}_{\G(\R)}({\cal V}, {\cal V}^{\otimes n-1})$. 
So we call  spaces (\ref{MSPX}) the {\it  multiplicity spaces}. 

 \vskip 1mm3. The regular representation of $\G \times\G$ in   functions  on $\G$. 
 
 \vskip 2mm
 We want  to find  explicit realizations of   Hilbert spaces   appearing   in the quantum representation theory: 
  
 \vskip 1mm 1.  Principal series  representations   of the modular double 
  $\mathcal{A}_\hbar(\mathfrak{g})$ of the 
  quantum group $\U_q(\mathfrak{g})$. 
  
  \vskip 1mm2. Multiplicity spaces =  spaces of invariants for   tensor powers of the principal series     of   
  $\mathcal{A}_\hbar(\mathfrak{g})$. 
  
 \vskip 1mm 3. The regular representation of $\mathcal{A}_\hbar(\mathfrak{g}) \times \mathcal{A}_\hbar(\mathfrak{g})^{\rm op}$ defined in Section \ref{SEC1.8}.
  
  \vskip 2mm
  
  First of all, let us recall   the definitions of these Hilbert spaces: 
  
  \vskip 1mm 1.   By   definition, the principal series    of  
  $\mathcal{A}_\hbar(\mathfrak{g})$ is realized in   the Hilbert space ${\cal H}({\mathscr L}_{\G, \odot})$.

 \vskip 1mm 2. Here the starting point is Theorem \ref{MFCQI}, which we recast as follows.  
  \bt  \la{MFCQIX}
Assuming Conjecture \ref{MK},  the multiplicity space for the $(n-1)$-st tensor power of the principal series $\ast-$representations of  
  $\mathcal{A}_\hbar(\mathfrak{g})$   is identified with the   Hilbert space ${\cal H}({\rm Loc}_{\G; S^2_n})$.
    
  The parameters $(\alpha_1, \ldots , \alpha_{n}) \in (\H/W)(\R_{>0})^{n}$ of   representations match the  symplectic leaves of the space ${\rm Loc}_{\G; S^2_n}(\R_{>0})$, 
  parametrized by the monodromies    around the punctures.  
  \et
  
\vskip 1mm 3. It is just the Hilbert space assigned to the cluster Poisson variety ${\cal R}_{\G, {\rm Cl}_{2,2}}$, see Section \ref{SEC1.8}. 
  
\vskip 2mm

  So the quest for   realizations of quantum  principal series representations,  their  multiplicity spaces, and the quantum regular representation 
   boils down to the question  of realization of the Hilbert space ${\cal H}({\mathscr L}_{\G, \odot})$,  ${\cal H}({\rm Loc}_{\G; S^2_n})$ and ${\cal H}({\cal R}_{\G, {\rm Cl}_{2,2}})$.

  A priori these spaces look very differently than their classical counterparts. 
  For example,  the space ${\cal H}({\rm Loc}_{\G; S^2_n})$ carries a braid group action, 
  while the space  (\ref{MSPX}) has just $\Z/n\Z$ action. So   to trace the connection between the two spaces we have to break the symmetry. 
    
   \vskip 2mm
   Note  that the Hilbert space  obtained by  quantization of a cluster Poisson variety ${\cal M}$ in general does not have a canonical realization as a Hilbert space 
   $L^2({\cal Y})$ for some natural space ${\cal Y}$. 
   
  \vskip 2mm
On the other hand, any cluster Poisson space  ${\mathscr X}$  gives rise to its {\it cluster symplectic double} ${\cal D}$ which  carry a 
 symplectic structure. Its algebra of regular functions   has a   $q-$deformation ${\cal O}_q({\cal D})$. For each cluster there is a  
 representation of the  modular double 
  $\mathcal{A}_\hbar({\cal D})$  in    ${\rm L_2}({\cal A}(\R_{>0}), \mu_{\cal A})$. The representations  for different clusters are related by unitary intertwiners\cite{FG07}. 
  Here ${\mathscr A}$ is a cluster $K_2$-variety dual to 
         ${\mathscr X}$, so that the pair $({\mathscr A}, {\mathscr X})$  is a cluster ensemble,  ${\mathscr A}(\R_{>0})$   is the set of its positive  points, and $\mu_{\mathscr A}$ the cluster  volume form.   
          We review  these construction   in Section \ref{Sec8a}.     
              \vskip 2mm 
               This suggests   the idea to  find   realizations of   Hilbert spaces ${\cal H}({\mathscr L}_{\G, \odot})$,  ${\cal H}({\rm Loc}_{\G; S^2_n})$ and ${\cal H}({\cal R}_{\G, {\rm Cl}_{2,2}})$ 
               using a symplectic double.                
               However, the  spaces ${\mathscr L}_{\G, \odot}$ and ${\rm Loc}_{\G; S^2_n}$ are Poisson, but not symplectic.     
                 
               To handle this issue, we develop in Section \ref{SEC9.1} a general notion of the   {\it double of a  cluster Poisson variety  with frozen variables}. 
              It is a Poisson space, obtained, roughly speaking, by doubling the unfrozen variables and keeping a single copy of the frozen ones. Its 
                center is generated by the 
               frozen variables.  If the set of  frozen variables is empty, we recover the symplectic double. 
                              
 \vskip 2mm
Let us now explain how we implement the these ideas.

\subsubsection{Main constructions.} Take a decorated surface $\bS$.  
Let us inflate each puncture on $\bS$ to a boundary component without special points. Denote by $\bS^\circ$ the same surface 
with the opposite orientation. Gluing  the inflated  surfaces $\bS$ and $\bS^\circ$ along the matching boundary components, we get a surface $\bS_{\cal D}$. It comes with punctures corresponding to the special points of $\bS$, and a collection of loops $\gamma = \{\alpha_i\}$ corresponding to the punctures on $\bS$.

 \begin{figure}[ht]
 \begin{center}
 \begin{tikzpicture}
 \draw[red] (0:1) --(60:1)--(120:1)--(180:1)--(240:1)--(300:1)--(360:1);
 \draw[red] (0,0) ellipse (0.3 and 0.15);
 \begin{scope}[shift={(4,0)}]
  \draw[blue] (0,0) ellipse (0.3 and 0.15);
 \foreach \count in {1,2,3,4,5,6}
{
\draw [blue, rotate=60*\count] (1, 0.1) arc (-90:-210:0.477);
\draw [blue, rotate=60*\count]  (1,0) ellipse (0.07 and 0.1);
}
 \end{scope}
 \end{tikzpicture}
 \end{center}
\caption{The decorated surface $\bS$  on the left is a hexagon with a hole. Its topological double  $\bS_{\cal D}$ is shown  on the right. It is a torus with six punctures, 
shown as small ovals on the right. The corners of the hexagon are its special points. They generate the six punctures on the topological double.}
\label{pin102}
\end{figure}

Let ${\cal D}_{\G, \bS}$ be the cluster symplectic double of the cluster Poisson variety ${\mathscr P}_{\G, \bS}$. 
   The corresponding quantum algebra  ${\cal O}_q({\cal D}_{\G, \bS})$  has  a  subalgebra  ${\bf B}$  
 generated by frozen variables at the boundary intervals of $\bS$. It  is 
  the center  of ${\cal O}_q({\cal D}_{\G,  \bS})^{\bf B }$ for generic $q$. Theorem \ref{11.26.18.1Xaaa} below is stated precisely in Section \ref{Sec8}. The simplest example  
  when $\G = {\rm PGL}_2$ and $\bS$ is a polygon is considered in Section \ref{Sec12.2n}. We discuss the general cluster set-up in Section \ref{SEC9.1}, and prove Theorem \ref{11.26.18.1Xaaa} 
   in Section \ref{SEC9.2}.   

  \bt  \la{11.26.18.1Xaaa}  
Assume that $\bS$ has no punctures. Then there are  injective  maps of   algebras:
\be \la{MAPZEaa}
\begin{split}
&\zeta_{\cal O}^*:     {\cal O}_q({\mathscr P}_{\G, \bS_{\cal D} } ) \hra {\cal O}_q({\cal D}_{\G,  \bS})^{\bf B },\\
&\zeta_{\mathscr A}^*:     \mathcal{A}_\hbar({\mathscr P}_{\G, \bS_{\cal D} } ) \hra \mathcal{A}_\hbar({\cal D}_{\G, \bS})^{\bf B }.\\ 
\end{split}
\ee
\et 
The first map (\ref{MAPZEaa}) induces the  second.   
 Cluster  varieties $({\mathscr A}_{\G,  \bS}, {\mathscr P}_{\G,  \bS})$  form a cluster ensemble by Theorem \ref{MTHa}. Since 
${\cal D}_{\G, \bS}$ is the cluster symplectic double of   ${\mathscr P}_{\G,  \bS}$, there  are cluster realizations of the $\ast-$algebra  $\mathcal{A}_\hbar({\cal D}_{\G, \bS})$ in ${\rm L_2}({\mathscr A}_{\G,  \bS}(\R_{>0}))$, depending on a choice of a cluster  on ${\mathscr A}_{\G, \bS}$. Combining  them with 
the embedding $\zeta_{\mathscr A}^*$ in (\ref{MAPZEaa}), we get  realizations of the $\ast-$algebra $\mathcal{A}_\hbar({\mathscr P}_{\G, \bS_{\cal D} })$. 

If the surface $\bS$ has punctures, Theorem \ref{11.26.18.1Xaaa} is valid, and established by combining the proof of Theorem \ref{11.26.18.1Xaaa} with  the main construction of \cite{FG14}. We leave this as an exercise for the reader. 

\subsubsection{Realizations of quantum multiplicity spaces.} Here is a sketch,   elaborated further in Section \ref{sec16.5:RMSP}. 

Let $\bS$ be  an   $n-$gon ${\P}_n$. Its topological double  
is an $n-$punctured  sphere $S_n^2$.  Let ${\cal D}_{\G, \P_n}$ be the cluster symplectic double of the cluster Poisson variety ${\mathscr P}_{\G, \P_n}$
   The algebra  ${\cal O}_q({\cal D}_{\G, \P_n})$  has  a  subalgebra  ${\bf B}$  
 generated by frozen variables.  The latter are attached to the sides of the polygon, and generate a subalgebra 
              ${\cal O}({\H'}^n)$, where ${\H'}$ is the Cartan group of the universal cover of $\G$.    Theorem \ref{11.26.18.1Xa}  provides 
  injective  maps of   algebras:
\be \la{MAPZE}
\begin{split}
&\zeta_{\cal O}^*:     {\cal O}_q({\mathscr X}_{\G, S^2_n} ) \hra {\cal O}_q({\cal D}_{\G, \P_n})^{\bf B }. \\
&\zeta_{\mathscr A}^*:     \mathcal{A}_\hbar({\mathscr X}_{\G, S^2_n} ) \hra \mathcal{A}_\hbar({\cal D}_{\G, \P_n})^{\bf B }.\\
\end{split}
\ee

Combining (\ref{MAPZE}) with the 
 representation of the $\ast$-algebra $\mathcal{A}_\hbar({\cal D}_{\G, \P_n})$ in  
${\rm L_2}({\mathscr A}_{  \G, \P_n}(\R_{>0}))$ we get a realization of the quantum multiplicity space. 
There is 
  an isomorphism of  cluster $K_2$-varieties:
$$
{\mathscr A}_{  \G, \P_n} =      {\rm Conf}_n({\cal A}_{  \G}).
$$
So  we realize the quantum multiplicity space in  functions on the same space as the classical one. 
 The map $\zeta_{\cal O}^*$ induces a canonical $W^n$-equivariant  injective  map  
$$
 {\rm Center}~{\cal O}_q({\mathscr X}_{  \G, S^2_n}) =  {\cal O}({\H^n}) \stackrel{ }{\hra}  {\bf B} \stackrel{\sim}{=}{\cal O}({{\H'}^n}).
$$

\subsubsection{Realization of  quantum principal series.} Let us give a sketch, which we  elaborate in Section \ref{SEC9.4}. 

 Let $\Delta$ be a triangle with two colored sides. 
The corresponding cluster Poisson moduli space   ${\mathscr P}_{\G, \Delta} $   
parametrizes   triples of flags    with a pair of pinnings, see also  Definition \ref{DEF13.12}.  
  We consider the double ${\cal D}_{\G, \Delta} $ of the space ${\mathscr P}_{\G, \Delta} $ related to the set of  frozen variables assigned to one of the colored sides.   
    Theorem \ref{Th9.13} relates it  to the  space  ${\cal R}_{\G, \odot}$. 
    
Since  by Theorem \ref{MTHa} $({\mathscr A}_{\G, \Delta},  {\mathscr P}_{\G, \Delta})$ form  a cluster ensemble, we get a realization of the $\ast$-algebra $\mathcal{A}_\hbar({\cal R}_{\G, \odot})$, and hence, via Theorem \ref{Th2.21}, of  $\mathcal{A}_\hbar(\frak{g})$, in 
$L_2({\mathscr A}_{\G, \Delta}(\R_{>0}), \mu_{\mathscr A})$.

 Finally, ${\mathscr A}_{\G, \Delta} $ is identified with the principal affine space ${\cal A}_\G$.   So we realize the quantum principal series  of  
 $\mathcal{A}_\hbar(\frak{g})$    on the same space as the principal series  for $\G(\R)$.

\subsubsection{Realization of  the regular representation of $\mathcal{A}_\hbar(\mathfrak{g})\times \mathcal{A}_\hbar(\mathfrak{g})$ in  ${\rm L}_2({\G}({\Bbb R}_{>0}))$.} 
\la{Section 5.4.5}
The cylinder ${\rm Cl}_{2,2}$ with $2+2$ special points is the topological double of the colored rectangle  
 $\square$ with $4$ special points given by the vertices, and the colored sides given by    the two opposite sides of the rectangle, 
see Section \ref{PLGG}.  

 The corresponding cluster Poisson variety 
 is   ${\mathscr P}_{\G, \square}$. The related cluster $K_2-$variety is ${\mathscr A}_{\G', \square}$. The pair $({\mathscr A}_{\G', \square}, {\mathscr L}_{\G, \square})$ is a cluster ensemble by Theorem \ref{MTHa}. Similarly to Theorem \ref{11.26.18.1Xaaa}, there is a  cluster realization of the $\ast-$algebra $\mathcal{A}_\hbar(\G, {\rm Cl}_{2,2})$  
 in ${\rm L}_2({\mathscr A}_{\G', \square}(\R_{>0}))$, provided by Theorem \ref{24.527.1013}, proved in Section \ref{SECTa17.8}. 
 
 \begin{figure}[ht]
 \begin{center}
 \begin{tikzpicture}[scale=0.7]
   \draw (0,0) -- (0,1);
   \draw (3,0) -- (3,1);
   \draw[dashed, -latex] (0,0)-- (3,0);
   \draw[dashed, -latex] (0,1) -- (3,1);
    \node[red]  at (0,0) {\tiny $\bullet$}; 
 \node[red]  at (0,1) {\tiny $\bullet$}; 
  \node[red]  at (3,0) {\tiny $\bullet$}; 
   \node[red]  at (3,1) {\tiny $\bullet$}; 
   \node at (4,.5) {$\ast$};
   \begin{scope}[shift={(5,0)}]
      \draw (0,0) -- (0,1);
   \draw (3,0) -- (3,1);
   \draw[dashed, -latex] (0,0)-- (3,0);
   \draw[dashed, -latex] (0,1) -- (3,1);
    \node[red]  at (0,0) {\tiny $\bullet$}; 
 \node[red]  at (0,1) {\tiny $\bullet$}; 
  \node[red]  at (3,0) {\tiny $\bullet$}; 
   \node[red]  at (3,1) {\tiny $\bullet$}; 
   \end{scope}
   \draw[-latex, thick] (9, .5) -- (10, .5);
   \begin{scope}[shift={(11,0.5)}, scale=.5]
    \draw (0,0) ellipse (0.5 and 1.25);
\draw (0, -1.25) -- (5, -1.25);
\draw (5, -1.25) arc (-90:90:0.5 and 1.25);
\draw [dashed] (5, -1.25) arc (270:90:0.5 and 1.25);
\draw (5, 1.25) -- (0, 1.25);  
\node [red] at (-0.5,0) {\tiny $\bullet$};
\node [red] at (0.5,0) {\tiny $\bullet$};
\node [red] at (5.5,0) {\tiny $\bullet$};
\node [red] at (4.5,0) {\tiny $\bullet$};
\end{scope}
 \end{tikzpicture}
 \end{center}
\caption{The topological double of a colored rectangle  $\square$ is a cylinder ${\rm Cl}_{2,2}$.}
\label{pin110}
\end{figure}

By Theorem \ref{Th4.3.3}, the cluster Poisson moduli space $  {\mathscr P}_{\G, \square} $ 
 is the  Poisson Lie group $\G$.  
 The finite cover $$
 {\mathscr A}_{\G', \square} \lra {\mathscr P}_{\G, \square}
 $$
   identifies the positive points: 
 $$
 {\mathscr A}_{\G', \square}(\R_{>0})  = 
 {\mathscr P}_{\G, \square}(\R_{>0}).
 $$
 So  the $\ast-$algebra $\mathcal{A}_\hbar(\G, {\rm Cl}_{2,2})$ has a cluster  realization in  ${\rm L}_2({\G}({\Bbb R}_{>0}))$. We elaborate all this in Section \ref{SEC16.8}.

\subsubsection{Quasiclassical limit in  the cluster set-up.}   Here is a   framework 
 for   quasiclassical limit.     
    
\vskip 1mmi) Any  vector  space ${\rm V}_{\rm quan}$   in the quantum representation theory appears as the Hilbert space of  
quantization of a cluster Poisson variety ${\cal M}$. 
The latter has a canonical projection $\pi: {\cal M} \lra \H_{\cal M}$ onto a split torus $\H_{\cal M}$. 
The   algebra $\pi^*{\cal O}(\H_{\cal M})$ is the center of the Poisson algebra ${\cal O}({\cal M})$ and its $q-$deformation  ${\cal O}_q({\cal M})$ for generic $q$. 

\vskip 1mmii) The classical vector space ${\rm V}_{\rm class}$  is realized in  a space of functions on a space ${\mathscr X}$. 

\vskip 2mm
The quasiclassical limit relating the quantum and the classical picture works as follows. 

\begin{itemize} \item The space ${\mathscr X}$ has a cluster Poisson variety structure, with a given set of frozen variables. 
\end{itemize}
The cluster Poisson variety  ${\mathscr X}$ gives rise to the quantum symplectic double    algebra 
  ${\cal O}_q({\cal D}_{\mathscr X})$. The   frozen variables on ${\mathscr X}$ generate a commutative subalgebra ${\bf B} \subset {\cal O}_q({\cal D}_{\mathscr X})$.  
   
 \begin{itemize} \item Let $ {\cal O}_q({\cal D}_{\mathscr X})^{\bf B} $ be  the centralizer  of ${\bf B} $. Then there are injective maps
\be \nonumber
\begin{split}
&i: {\cal O}_q({\cal M}) \stackrel{}{\hra}  {\cal O}_q({\cal D}_{\mathscr X})^{\bf B},  \\
& i: \pi^*{\cal O}(\H_{\cal M}) \hra {\bf B}.\\
\end{split}
\ee
The targets are finite dimensional over the source algebras. 
  
\vskip 1mm   \item Combining the embedding $i$ with the canonical realization of ${\cal O}_q({\cal D}_{\mathscr X})$   in  ${\rm L_2}({\mathscr A}(\R_{>0}), \mu_{\mathscr A})$, we get a 
    {\it canonical cluster realization 
   of ${\cal O}_q({\cal M})$  in  $L_2({\mathscr A}(\R_{>0}), \mu_{\mathscr A})$.}
     \end{itemize}

\vskip 2mm     
     
     In all   examples there is a finite cover ${\mathscr A} \lra {\mathscr X}$ inducing an isomorphism  ${\mathscr A}(\R_{>0}) = {\mathscr X}(\R_{>0})$. It  provides an identification 
          $L_2({\mathscr A}(\R_{>0}), \mu_{\mathscr A}) = L_2({\mathscr X}(\R_{>0}), \mu_{\mathscr X})$. We conclude that 
     \begin{itemize} \item There is 
   canonical cluster realization 
   of ${\cal O}_q({\cal M})$  in  ${\rm L_2}({\mathscr X}(\R_{>0}), \mu_{\mathscr X})$. Thus 
$
  {\rm V}_{\rm quan} \sim  {\rm V}_{\rm class}.
$
    \end{itemize}        
          
So we realize the quantum vector space    ${\rm V}_{\rm quan}$  in the space of functions on the   positive points    of the 
same space ${\mathscr X}$ where  the classical vector space    ${\rm V}_{\rm class}$  is realized. Yet there is no   canonical isomorphism  
 ${\rm V}_{\rm quan} =  {\rm V}_{\rm class} $ since the vector space ${\rm V}_{\rm quan}$ has a multitude of realisations parametrised by  clusters on 
 ${\mathscr X}$, and hence many isomorphisms 
 with the space    ${\rm V}_{\rm class}$.

\vskip 2mm

 In the traditional   set-up ${\mathscr X}$ is just a manifold, called  the configuration space, and ${\cal M}:= T^*{\mathscr X}$. It is quantized to the algebra differential operators on ${\mathscr X}$, acting   in the space of functions on ${\mathscr X}$.

 We stress that, unlike in the classical case, the cluster Poisson structure on ${\mathscr X}$ plays a key role. Furthermore, ${\cal O}_q({\cal M})$ is realised not in the functions on 
 ${\mathscr X}$, but in functions on the dual space  ${\mathscr A}(\R_{>0})$, and the realisation itself depends on the choice of a cluster.

\subsubsection{Conclusions.} There are many analogies between the classical  principal series representations of $\G(\R)$ and the quantum ones for the quantum group modular 
double $\ast-$algebra $\mathcal{A}_\hbar(\mathfrak{g})$: 
  
\begin{itemize}

\vskip 1mm \item The quantum principal series $\ast-$representations of  
$\mathcal{A}_\hbar(\mathfrak{g})$ are realized in   ${\rm L}_2({\mathscr A}(\R_{>0}), \mu_{\mathscr A})$.   The classical one are realised in ${\rm L}_2({\mathscr A}(\R), \mu_{\mathscr A})$ . 
 
\vskip 1mm\item The   Weyl group action   by  quantum unitary   intertwiners  is the analog of Gelfand-Graev unitary intertwiners in ${\rm L}_2({\mathscr A}(\R), \mu_{\mathscr A})$. 
So for any   $\lambda \in \H(\R_{>0})$, the quantum principal series  $\ast-$representations ${\cal V}_{w(\lambda)}$, $w \in W$,  of  
$\mathcal{A}_\hbar(\mathfrak{g})$, 
 are canonically isomorphic. 
\vskip 1mm\item  The commutative subalgera $\mu_p^*{\cal O}_\H$ is the analog of the  center   $\widehat {\cal Z}_{\U({\g})}$ of  $   \widehat {\U}({\g})$.\footnote{As 
 Section \ref{SEC12.4} shows,   $\mu_p^*{\cal O}_q(\H)$ might be  slightly different then $\widehat {\cal Z}_{\U({\g})}$.}

\vskip 1mm\item  The algebra $\mathcal{A}_\hbar({\cal D}_\Delta)$ is the analog of the algebra  of differential operators ${\rm Diff}_{\mathscr A}$. 
  The  actions of   $\mathcal{A}_\hbar({\cal D})$ in functions on 
${\mathscr A}(\R_{>0}) $  are the analog of the action of   ${\rm Diff}_{\mathscr A}$ on functions on ${\mathscr A}(\R)$. 

  The analogs of isomorphisms (\ref{IUD}) and (\ref{IUD1}) are the isomorphisms 
\be \nonumber
\begin{split}
&{\cal O}_q({\cal R}_{\G, \odot}) = {\cal O}_q({\cal D}_\Delta)^\H,\\
&{\cal O}_q({\mathscr L}_{\G, \odot}) =   {\cal O}_q({\cal D}_\Delta)^{\H, W}.\\
\end{split}
\ee

\vskip 1mm  \item Multiplicity spaces of the classical and quantum principal series representations can be realized in the  spaces of functions related to the same variety   ${\rm Conf}_n({\cal A})$.  
  
 \item The quantum regular representation of the $\ast-$algebra $\mathcal{A}_\hbar(\mathfrak{g}) \otimes \mathcal{A}_\hbar(\mathfrak{g})$ is realized in ${\rm L}_2(\G(\R_{>0}))$.

 \vskip 1mm \item There is a  unitary projective  action of the braid group ${\Bbb B}_{\mathfrak g}$  on each principal series  representation ${\cal V}_\lambda$, 
  quantizing Lusztig's   braid group action on $\U_q(\mathfrak{g})$. This proves that the braid group action on $\U_q(\mathfrak{g})$ preserves 
  the isomorphism class of representation.   
  \end{itemize}

Finally, we want to stress again the key new aspects of the quantum representation theory: 

\vskip 2mm
{\it Although   representations are realized 
in functions on classical spaces, the way we realize them depends crucially on a choice of a cluster coordinate system on those spaces. 

Therefore quantum representation theory can not be developed without establishing the cluster nature of the relevant classical spaces. 
The comparison between the realizations in different cluster coordinate system is the backbone of quantum representation theory.} 

    \medskip 

\subsection{Integrable systems on the moduli  spaces of {$\G$}-local systems on $S$} \la{SECT2.8a}

\medskip

\subsubsection{The integrable system for ${\rm SL_2}$.} Let $S$ be a genus $g>1$ surface with $n$ punctures. Pick a collection $\gamma$ 
of simple loops $\gamma= \{\alpha_1, ..., \alpha_{3g-3+n}\}$, cutting $S$ into pair of pants. 

If $\G = {\rm SL}_2$, then the functions ${\rm M}_{\alpha_i}$  on the space ${\rm Loc}_{{\rm SL_2}, S}$ given by the eigenvalues of the monodromies of ${\rm SL}_2-$local systems along the loops $\alpha_i$ commute, and 
generate an algebra of polynomials in $3g-3+n$ variables.  Note that 
$$
{\rm dim}~{\rm Loc}_{{\rm SL_2}, S} = -3 \chi(S) = 2\cdot (3g-3+n) +n.
$$
The  space ${\rm Loc}_{{\rm SL_2}, S}$ is   Poisson   with the center generated by the regular functions $C_1, ..., C_n$ given by the 
eigenvalues of the monodromies  around the punctures. 
The fibers   
 with the given  values  of the Casimir function $C_i$  are the generic symplectic leaves. 
The monodromies   ${\rm M}_{\alpha_i}$ are Hamiltonians of an integrable system at the generic fiber. 
This integrable system plays an important role in many  areas, e.g. in Quantum Field Theory \cite{AGT},  \cite{NS},  \cite{T10}, \cite{JN}. 
We suggest a generalization for a simply-laced $\G$, which we think should play a similar role. 

\begin{figure}[ht]
 \begin{center}
 \begin{tikzpicture}[scale=1.3]
  \begin{scope}[shift={(7,0)}]
  \draw[blue] (0,0) circle (0.5);
 \foreach \count in {1,2,3,4,5,6}
{
\draw [rotate=60*\count] (60:.68) -- (120: .68);
\draw [rotate=60*\count] (60:.68) -- (60: 1.3);
\draw [blue, rotate=60*\count] (0, 0.5) arc (-90:-270:0.06 and 0.089);
\draw [blue, dashed, rotate=60*\count] (0, 0.5) arc (-90:90:0.06 and 0.089);
\draw [blue, rotate=60*\count] (1, 0.1) arc (-90:-210:0.477);
\draw [blue, rotate=60*\count]  (1,0) ellipse (0.07 and 0.1);
}
\end{scope}
 \begin{scope}[shift={(3.5,0)}]
  \draw[red] (0,0) circle (0.5);
 \foreach \count in {1,2,3,4,5,6}
{
\draw [rotate=60*\count] (60:.68) -- (120: .68);
\draw [rotate=60*\count] (60:.68) -- (60: 1.3);
\draw [red, rotate=60*\count] (1, 0.1) arc (-90:-210:0.477);
}
 \end{scope}

 \foreach \count in {1,2,3,4,5,6}
{
\draw [rotate=60*\count] (60:.68) -- (120: .68);
\draw [rotate=60*\count] (60:.68) -- (60: 1.3);
\draw [dashed, red, rotate=60*\count] (0,0) -- (30:1.155) -- (90:1.155);
}
 \end{tikzpicture}
 \end{center}
\caption{A 3-valent ribbon graph $\Gamma$  (black, on the right) gives rise to a surface $\bS_{\cal D}$ with a pair of pants decomposition (blue, on the left), and  a decorated surface $\bS_\Gamma$ (red, in the middle), equivalent to a triangulated surface  at the right (shown by punctured red). The surface $\bS_{\cal D}$ is the double of $\bS_\Gamma$.}
\label{pin101}
\end{figure}

\subsubsection{An integrable system for any   group $\G$.}  Below $\G'$ denotes the universal cover of the adjoint group $\G$. 
Let  $\gamma$ be a pair of pants decomposition of a surface $S$, given by a collection of $3g-3$ simple loops $\{\alpha_i\}$. Then there is a moduli space 
${\mathscr P}_{{\G'}, S; \gamma}$  \cite[Section 2.2.1]{FG14} and a projection map
$$
{\mathscr P}_{{\G'}, S; \gamma}\lra {\mathscr P}_{{\G'}, S}.
$$
Namely, the stack ${\mathscr P}_{{\G'}, S; \gamma}$ parametrises $\G'-$local systems on $S$ equipped with a choice of an invariant flag 
 for each loop $\alpha_i$ of the collection $\gamma$ .  Forgetting the invariant flags at the loops, we get the projection onto ${\mathscr P}_{{\G'}, S}$. It is a 
 Galois cover over the generic point with the Galois group $W^{3g-3}$, where $W$ is the Weyl group.  
 
Let us  present the surface $S$ as a topological double of a decorated surface $\bS_\Gamma$, assigned to  a 3-valent ribbon graph $\Gamma$, possibly with legs.  
Recall that  a ribbon graph $\Gamma$ determines a decorated surface $\bS_\Gamma$, obtained by gluing the ribbons of $\Gamma$,  
  with a collection of face paths, shown   by solid red curves on Figure \ref{pin101},  \ref{pin101+}. 
  The boundary of the surface $\bS_\Gamma$ is a disjoint union of circles, given by the 
  boundaries of the holes, and intervals, given by the face arcs.  Shrinking  the holes of $\bS_\Gamma$  to punctures, and the  face arcs to special points, we get an equivalent surface     glued from the triangles $t_v$ dual to the links of the vertices $v$ of $\Gamma$, depicted
  as punctured red triangles. 

\begin{figure}
 \begin{center}
 \begin{tikzpicture}[scale=1.3]
 \node at (0,0) {$\bullet$};
 \node at (0,0.2) {$v$};
  \node at (1,0) {$\bullet$};
   \node at (1,0.2) {$w$};
 \draw[thick] (133:1)--(0,0)--(1,0)--+(47:1);
 \draw[thick] (-133:1)--(0,0);
 \draw[thick] (1,0)--+(-47:1);
 \node at (0.5, -0.6) {$\Gamma$};

\begin{scope}[shift={(4,0)}]
\draw[thick] (133:1)--(0,0)--(1,0)--+(47:1);
 \draw[thick] (-133:1)--(0,0);
 \draw[thick] (1,0)--+(-47:1);
\draw [red][thick] (1.5, 0.7) arc (-60:-120:2);
\draw  [red][thick] (1.5, -0.7) arc (60: 120:2);
\draw  [red][thick] (1.63, 0.55) arc (150:210:1.1);
\draw  [red][thick] (-0.63, -0.55) arc (-30:30:1.1);
\end{scope}

 \begin{scope}[shift={(8.5,0)}]
 \begin{scope}[shift={(1.05,0.61)}]
 \draw  [blue][thick, rotate=30]  (0,0) ellipse (0.07 and 0.1);
 \end{scope}
  \begin{scope}[shift={(1.05,-0.61)}]
 \draw  [blue][thick, rotate=-30]  (0,0) ellipse (0.07 and 0.1);
 \end{scope}
  \begin{scope}[shift={(-1.05,-0.61)}]
 \draw  [blue][thick, rotate=30]  (0,0) ellipse (0.07 and 0.1);
 \end{scope}
  \begin{scope}[shift={(-1.05,0.61)}]
 \draw  [blue][thick, rotate=-30]  (0,0) ellipse (0.07 and 0.1);
 \end{scope}
 \node at (-0.5,0) {$\mathscr{P}_v$};
  \node at (0.5,0) {$\mathscr{P}_w$};
\draw [blue][thick] (1, 0.7) arc (-60:-120:2);
\draw  [blue][thick] (1, -0.7) arc (60: 120:2);
\draw  [blue][thick] (1.13, 0.55) arc (150:210:1.1);
\draw  [blue][thick] (-1.13, -0.55) arc (-30:30:1.1);
\draw [blue][thick] (0, 0.432) arc (90:270:0.1 and 0.432);
\draw  [blue][thick, dashed] (0, 0.432) arc (90:-90:0.1 and 0.432);
 \end{scope}
 \end{tikzpicture}
 \end{center}
\caption{The 3-valent ribbon graph $\Gamma$  on the left gives rise to the decorated surface $\bS_\Gamma$, shown as    the red rectangle in the middle. The 
 topological double of the  $\bS_\Gamma$  is a sphere with four punctures $\bS_{\cal D}$ on the right. The pair $\Gamma \subset \bS_\Gamma$ in the middle determines    a pair of pants decomposition of the surface $\bS_{\cal D}= {\mathscr P}_v \cup {\mathscr P}_w$ on the right.}
\label{pin101+}
\end{figure}

For any surface $S$, we can find a ribbon graph $\Gamma$ such that the topological double of  $\bS_\Gamma$, obtained by gluing $\bS_\Gamma$ and its mirror $\bS_\Gamma^\circ$ along the 
matching boundary components,  is identified with the original surface $S$, see the  blue surfaces on Figures \ref{pin101},  \ref{pin101+}. The topological double  of $\bS_\Gamma$ comes with  
  a decomposition into pairs of pants ${\mathscr P}_v$, labeled by the  vertices $v$   of $\Gamma$, matching the  initial pair of pants decomposition $\gamma$ of $S$:
\be \nonumber
S = 
\bS_\Gamma \cup \bS_\Gamma^\circ, \qquad S= \cup_v {\mathscr P}_v.
\ee

Since the $\G'-$local systems parametrised by the space ${\mathscr P}_{{\G'}, S; \gamma}$ are equipped with the flags over each loop $\alpha_i$, invariant under the monodromy around the loop $\alpha_i$, the semisimple part of the monodromy around the loop $\alpha_i$ takes values in the Cartan group $\H$ of $\G'$. 
So it provides the $\H-$valued monodromy $\H_{\alpha_i}$.   Applying to it the characters of $\H$, given by the simple dominant  weights $\omega_j$, we arrive that the 
commuting Hamiltonians $$
H_{\alpha_i}^j: = \omega_j(\H_{\alpha_i}).
$$
Let us complement them  to a collection of commuting Hamiltonians, whose number is exactly  
 half of the dimension of the symplectic leaves of ${\mathscr P}_{{\G}, S; \gamma}$.

   For each  vertex $v$ of $\Gamma$,    pick a   cluster  coordinate system ${\bf c}_v = \{A^{{\bf c}_v}_1, \ldots , A^{{\bf c}_v}_{m+3r}\}$ on ${\rm Conf}_3({\cal A})$.
   Here $$
   m = {\rm dim} {\rm Conf}_3({\cal B}) = {\rm dim}\U - r, \ \ \ r= {\rm dim}\H = {\rm rk}\G.
   $$ 
   Corollary \ref{12.8} provides a set of Poisson commuting  
   rational functions $ \{B^{{\bf c}_v}_1, \ldots , B^{{\bf c}_v}_{m+3r}\}$ on ${\mathscr P}_{\G', v}$. These functions    include the   $3 r$ Hamiltonians 
given by the   monodromies around the boundary   of ${\mathscr P}_{\G, v}$. The rest of  the functions, 
assigned to  the unfrozen cluster coordinates on ${\rm Conf}_3({\cal A})$,  are denoted by 
$$
 \{B^{{\bf c}_v}_1, \ldots , B^{{\bf c}_v}_{m}\}.
 $$ 
 We pull them back   to  ${\mathscr P}_{{\G}, S; \gamma}$ via the restriction map ${\mathscr P}_{{\G}, S; \gamma} \to {\mathscr P}_{{\G},  v}$ and, abusing notation,   denote them  by $\{B^{{\bf c}_v}_1, \ldots , B^{{\bf c}_v}_{m}\}$. When $v$ runs through all vertices of $\Gamma$, these functions, complemented by the 
$ r \cdot (3g-3+n)$ Hamiltonians $H_{\alpha_i}^j$,  provide a collection of rational functions 
\be \la{HAML} 
 \{H_{\alpha_i}^j \} ~\cup   ~\{B^{{\bf c}_v}_1, \ldots , B^{{\bf c}_v}_{m}\}, \qquad \mbox{where} ~v \in \{\mbox{\rm vertices of $\Gamma$}\}.
\ee   
So all together, since the number of pair of pants is $2g-2+n$,  we got 
$$
r \cdot (3g-3+n) + m \cdot (2g-2+n) = \frac{1}{2} \cdot ({\rm dim}{\mathscr P}_{{\G}, S} - n r)
$$
commuting Hamiltonians, which is precisely the  half of the dimension of the generic symplectic leaf of  ${\mathscr P}_{{\G}, S}$. 

To check the above identity, note that ${\rm dim}\G=2m+3r$, and so 
$$
{\rm dim}{\mathscr P}_{{\G}, S} -nr = -{\rm dim}\G \cdot \chi(S) -nr= (2m+3r)(2g-2+n) -nr = 2\cdot ((2m+3r)(g-1) + (m+r)n).
$$ 
  \bt  \la{11.26.18.1Xa}  
Let $\Gamma$ be a trivalent ribbon graph. 
Then  the   functions (\ref{HAML}) Poisson commute   in ${\cal O}({\mathscr P}_{\G', S; \gamma} )$, 
providing a polarisation =  integrable system on the symplectic leaves of ${\mathscr P}_{\G', S; \gamma}$.
\et

 \begin{proof} As was mentioned above, for each pair of pants ${\mathscr P}_v$, Corollary \ref{12.8}   provides a set of Poisson commuting  
   rational functions $ \{B^{{\bf c}_v}_1, \ldots , B^{{\bf c}_v}_{m+3r}\}$ on ${\mathscr P}_{\G', v}$. 
   The functions $\{B^{{\bf c}_v}_1, \ldots , B^{{\bf c}_v}_{m}\}$   assigned to different pair of pants ${\mathscr P}_v$ and ${\mathscr P}_w$ Poisson commute:
   $$
   \{B^{{\bf c}_v}_s, B^{{\bf c}_w}_t\}=0    \ \ \ \forall \ 1 \leq s, t\leq m.
   $$
    \end{proof}

Unlike the ${\rm SL_2}-$case, different clusters ${\bf c}_v$ on ${\rm Conf}_3({\cal A})$   lead to different integrable systems.

\subsubsection{Towards a generalisation of the AGT conjecture to any $\G$.} The AGT conjecture 
 \cite{AGT}  relates   Nekrasov's partition functions \cite{N} of Gaiotto's  4d ${\cal N}=2$ superconformal field theory on $\R^4$ \cite{G} to conformal blocks in Liouville theory. 
  These 4d SCFT   should correspond to infra-red limit of the compactification of the hypothetical  six-dimensional superconformal (2,0) theory of type $A_1$
on a genus $g$ Riemann surface $\Sigma$ with $n$ punctures.  

Compactifying  the six-dimensional (2,0) theory of type $\G$ on $\Sigma$ one should get more general 
4d SCFT, referred to below as ${\rm SCFT_4}(\G; \Sigma)$. Recall that 
the pair of pants decomposition of $S$ assigned to $\Gamma$ determines a Deligne-Mumford boundary divisor ${\rm DM}_\Gamma \subset \overline {\cal M}_{g,n}$. 

\vskip 2mm
1.   We suggest that, given a pair of pants decomposition of $\Sigma$  related to a ribbon graph $\Gamma$, and  
\be \nonumber
\mbox{\it a collection of cluster coordinate systems $\{{\bf c}_v\}$ for  ${\rm Conf}_3({\cal A})$,  where $\{v\} = \{ \mbox{vertices of}$ $\Gamma \}$}
\ee
 one should be able to produce a Lagrangian 
description ${\rm SCFT}_4(\G; \Sigma, \Gamma, \{{\bf c}_v\})$  of the  ${\rm SCFT_4}(\G; \Sigma)$, when $\Sigma$ 
degenerates according to the Deligne-Mumford stratum ${\rm DM}_\Gamma$, related   
 a la Nekrasov-Shatashvili \cite{NS} to the integrable system   assigned to the triple $(\G; \Gamma, \{{\bf c}_v\})$. 

\vskip 2mm
2. As   explained in Section \ref{SC2.7}, a data $(\G; \Gamma, \{{\bf c}_v\})$ provides a realisation of the   space of conformal blocks 
for the pair $(\G, S)$ as   ${\rm L}_2({\mathscr A}_{\G, \bS_\Gamma})$. 
Another collection of clusters $\{{\bf c}'_v\}$, obtained from $\{{\bf c}_v\}$ by a cluster tranaformation,  amounts  
to a unitary  intertwiner  relating the two realizations:
\be \la{INN}
{\rm I}_{\{{\bf c}_v\} \to \{{\bf c}'_v\}}: {\rm L}_2({\mathscr A}_{\G, \bS_\Gamma}) \lra {\rm L}_2({\mathscr A}_{\G, \bS_\Gamma}).
\ee

\vskip 2mm
3. 
The Nekrasov  partition function for ${\rm SCFT}_4(\G; \Sigma, \Gamma, \{{\bf c}_v\})$ should  match  the expansion of conformal blocks of the $\G-$Toda theory, realised in 
${\rm L}_2({\mathscr A}_{\G,  \bS_\Gamma})$, near the divisor  ${\rm DM}_\Gamma$.

\vskip 2mm

4. Theories  ${\rm SCFT}_4(\G; \Sigma, \Gamma, \{{\bf c}_v\})$ and ${\rm SCFT}_4(\G; \Sigma, 
\Gamma, \{{\bf c}'_v\})$  should be  related by an  equivalence -  a generalising $S-$duality transformation - 
 assigned to the cluster transformation 
${\bf c}_v \to {\bf c}'_v$.  

By 3),  their partition functions correspond to conformal blocks. The latter,  realized 
in  ${\rm L}_2({\mathscr A}_{\G, \bS_\Gamma})$, should be related by the unitary intertwiner (\ref{INN}).

\medskip

\section{Representation theory and quantized  moduli spaces of local systems } \la{sec1.9}  
  
\medskip

 Section \ref{sec1.9} contains one of the most important applications of the constructions of the paper. However it is not directly used in the rest of the paper. Therefore, readers may choose to skip it initially and revisit it later as needed.

\subsubsection{The set-up.}

Recall  the Teichm\"uller space ${\cal T}_{g,n}$ parametrizing   data $(\Sigma, S;  \varphi)$, where $\Sigma$ is a genus $g$ Riemann   surface 
with $n$ punctures $p_1, ..., p_n$, $S$ is an oriented topological surface with punctures, and $$
\varphi: S \stackrel{\sim}{\lra} \Sigma
$$
is a   homeomorphism considered up to an isotopy. The Teichm\"uller space is contractable. 

The mapping class group $\Gamma_S= {\rm Diff}(S) / {\rm Diff}_0(S)$  acts on the   space ${\cal T}_{g,n}$. 

The quotient 
${\cal T}_{g,n}/\Gamma_S$ is the moduli space ${\cal M}_{g,n}$, so $\Gamma_S = \pi_1({\cal M}_{g,n})$. 
    There are $n+1$ line bundles    
    over the moduli space ${\cal M}_{g,n}$: the determinant line bundle $  {\rm det} ~\Omega^1_\Sigma$, and the cotangent line bundles  $T^*_{p_i}\Sigma$ at   the points  $p_i\in \Sigma$. 
    The product of the  punctured line bundles is a $(\C^*)^{n+1}-$bundle $ \widehat {\cal M}_{g,n}$ over ${\cal M}_{g,n}$:
\be \la{TORUSB}
  \widehat {\cal M}_{g,n}:= \left({\rm det} ~\Omega^1_\Sigma- \{\mbox{the zero section}\}\right) \times \prod_{i=1}^n \left(T^*_{p_i}\Sigma - \{0\}\right).
\ee  
The {\it extended mapping class 
group} defined by    $\widehat \Gamma_S:= \pi_1(   \widehat {\cal M}_{g,n})$  is therefore  a central extension:
\be \nonumber
0 \lra \Z^{n+1} \lra \widehat \Gamma_S \lra  \Gamma_{S}\lra 0.
\ee
The extension class is described by 
the  Chern classes    of these line bundles in $H^2({\cal M}_{g,n}) = H^2(\Gamma_S)$.

The {\it extended Teichm\"uller space $\widehat {\cal T}_{g,n}$} is the universal cover of    $\widehat {\cal M}_{g,n}$. So the group $\widehat \Gamma_S$ acts on it, and
$$
\widehat {\cal T}_{g,n}/ \widehat \Gamma_S =     \widehat {\cal M}_{g,n}.
$$

A connection on a vector bundle is {\it integrable} if it gives rise to the parallel transport along paths. 
 An {\it integrable flat connection} is   the same thing as a {\it local system} of  vector spaces. 
A connection on an infinite dimensional  vector bundle may not be integrable.

A local system on  $\widehat {\cal M}_{g,n}$ is the same thing as a 
$\widehat \Gamma_S$-equivariant  local system  on   $\widehat {\cal T}_{g,n}$.  

Given a point $\tau \in   \widehat {\cal M}_{g,n}$, 
there are canonical functors:  
\be \nonumber
\begin{split}
& \{\mbox{Representations of the group $\widehat \Gamma_S$ in   infinite dimensional  vector spaces}\} \stackrel{\sim}{\lra}  \\
&\mbox{\{$\widehat \Gamma_S$-equivariant local systems of infinite dimensional   vector bundles  
on   $\widehat {\cal T}_{g,n}$\}} \lra  \\
&\mbox{\{$\widehat \Gamma_S$-equivariant  infinite dimensional  vector bundles  
 with flat connection
on   $\widehat {\cal T}_{g,n}$\}}.\\\end{split}
\ee
The first    is an equivalence. 
 The second is  a fully faithful embedding. Its  image consists of integrable flat connections. 
 
\vskip 2mm
  Conjecture \ref{GENDKL}    relates  two kinds of data 
 of  very  different nature.  
 
 \vskip 2mm
 The first, which we call {\it de Rham data}, come from Representation Theory. 
It  describes a continuous analog of   fusion tensor categories for   
 certain continuous series of infinite dimensional   highest weight representations of the $W$-algebra   ${\cal W}_{\mathfrak{g}^\vee}$.  
 Precisely, a de Rham data is given by 
 
 \begin{itemize}
 
 \item  {\it De Rham   bundle},  given   by the   coinvariants\footnote{Given a Lie algebra $\mathfrak{g}$ acting on a vector space $V$, the   space of coinvariants 
 is  the quotient $V_{\mathfrak{g}}:= V / \mathfrak{g}V$. On the other hand, there is a space $(V^*)^{\mathfrak{g}} = (V_{\mathfrak{g}})^*$ of invariants in the   dual vector space $V^*$. If $V$ is finite dimensional, the two spaces are dual to each other. Otherwise  the  space of invariants $(V_{\mathfrak{g}})^*$ is much bigger  than the space of coinvariants $V_{\mathfrak{g}}$. }  
 of  the oscillatory  representations of the $W$-algebra   ${\cal W}_{\mathfrak{g}^\vee}$, 
 assigned to  punctures of a Riemann surface $\Sigma$. It is a bundle  of infinite dimensional vector spaces on the moduli space $\widehat {\cal M}_{g,n}$, with  discrete topology  
and  a flat non-integrable connection.
 \footnote{The dual to this space  is often called the space of {\it conformal blocks}. 
 For example, let $\G = {\rm SL}_2$, and  let $V:= {\cal V}_\Sigma$ be a tensor product of Verma modules 
  and $\mathfrak{g}:= {\rm Vir}_\Sigma$   the Lie algebra of meromorphic vector fields on $\Sigma$ with poles at the punctures. In \cite{FBZ}, or \cite[page 82]{VT}, 
   the space of conformal blocks is  defined as the space of invariants of ${\rm Vir}_\Sigma$ 
   acting on the dual vector space $({\cal V}_\Sigma)^*$. We avoid calling it the space of conformal blocks since by  Conjecture \ref{GENDKL} the conformal blocks   should lie in a different  space ${\cal S}^*({\rm Loc}_{\G, S})$.}  
  
 \end{itemize}


 \vskip 2mm
 The second, which we call {\it Betti data}, 
is provided by the quantization of the  moduli space ${\rm Loc}_{\G, S}$. 

\vskip 2mm
Let us elaborate now on the definitions of the de Rham and Betti data. 

\subsubsection{$W$-algebras.} Recall the $W$-algebra ${\cal W}_{\mathfrak{g}}$ related to a  
  Lie algebra ${\mathfrak{g}}$.\footnote{We consider only   $W$-algebras assigned to the principal nilpotent element of ${\mathfrak{g}}$.}   
 For  $\mathfrak{g} = {sl_2}$ it  is   the Virasoro Lie algebra. 
Recall the Langlands dual Lie algebra    ${\mathfrak{g}^\vee}$   for the    Lie algebra  ${\mathfrak{g}}$.    
  
The  W-algebra ${\cal W}_{\mathfrak{g}}$ contains  the Virasoro subalgebra. The central charge of a representation of the $W$-algebra is defined as the 
central charge of its Virasoro subalgebra.  Denote by  ${\rm rk}_{\mathfrak{g}}$  the rank of $\mathfrak{g}$,  by $\rho_{\mathfrak{g}}$ the half sum of the positive roots of $\mathfrak{g}$, and  by 
$h^\vee$ the dual Coxeter number.

Let $\widehat {\mathfrak{g}}$ be the Kac-Moody Lie algebra  assigned to a  Lie algebra $\mathfrak{g}$. Given a level $k$ representation $V$ of $\widehat {\mathfrak{g}}$, 
the quantum Drinfeld-Sokolov reduction, given by the semi-infinite cohomology $H^{\infty/2}_{\mathfrak{n}(t)}(V, \psi)$, is a representation of 
 the associated W-algebra ${\cal W}_{\mathfrak{g}}$. Its Virasoro subalgebra   acts  
 with the central charge
\be \la{CCH}
\begin{split}
c_\mathfrak{g}:= &{\rm rk}_{\mathfrak{g}}\Bigl(1-  h^\vee(h^\vee+1)\frac{(k+h^\vee-1)^2}{k+h^\vee}\Bigr)\\
= &{\rm rk}_{\mathfrak{g}}\left(1 +  h^\vee(h^\vee+1)Q^2\right).\\
\end{split}
\ee
 It is useful to express $Q$ via the Planck constant $\hbar$   as follows\footnote{We use $\beta$ rather then the traditional $b$.}
\be \la{huw2}
\begin{split} 
&Q^2 = \hbar + \hbar^{-1}+2,\\
&Q = \beta+ \beta^{-1}, \qquad{\rm Re}(\beta) \geq 0, \qquad \hbar = \beta^2.\\ 
\end{split}
\ee
Note that   the condition $Q^2  \geq 0$ just means that there are   two options for the Planck constant:   
\be \nonumber 
\begin{split} 
& \hbar \in \R_{>0} \qquad\mbox{or} \qquad|\hbar|  = 1.\\ 
\end{split}
\ee
The Planck constant $\hbar$ related to the level $k$ by 
\be  \la{hbark}
\hbar = -(k+h^\vee).
\ee
To check (\ref{CCH}) note in terms of the Planck constant, we have 
\be \nonumber
\begin{split}
&\frac{(k+h^\vee-1)^2}{k+h^\vee} =   
{k+h^\vee} + \frac{ 1}{k+h^\vee} -2 = -(\hbar +\hbar^{-1} + 2).\\
\end{split}
\ee

\subsubsection{Oscillatory series of   representations of  W-algebras.} The algebra ${\cal W}_{\mathfrak{g}}$ has a series 
 of unitary highest weight representations, which we call the {\it oscillatory series}. Their definition follows from the works of Feigin - Frenkel \cite{FF1}, \cite{FF2}. 
Oscillatory representations $V_\lambda$  of  ${\cal W}_{\mathfrak{g}}$ are  parametrized by    
  $\alpha \in   \mathfrak{h}^\ast$ and a non-negative real number $Q \in \R_{\geq 0}$, 
  with the weight $\lambda$ and the central charge $c$  given by\footnote{Beware of the clash of  notation: $\hbar$ is the Planck constant, $h^\vee$ is the dual Coxeter number,  $h$ is the highest weight of the Verma module for the Virasoro Lie algebra, and if this was not enough,
   $\mathfrak{h}$ is the Cartan Lie algebra.} 
\be \la{huw1}
\begin{split} 
& \lambda = Q\rho_{\mathfrak{g}} + i \alpha \in Q\rho_{\mathfrak{g}} + i \mathfrak{h}^\ast,\\
&c =
 {\rm rk}_{\mathfrak{g}}\left(1 +  h^\vee(h^\vee+1)Q^2\right).\\
\end{split}
\ee
The representations on the same orbit of the $\circ$-action of the Weyl group $W$ are equivalent:
\be \la{WWW}
 V_\lambda \stackrel{\sim}{=}  V_{w\circ  \lambda}\qquad \forall w \in W, \qquad w\circ \lambda:= w(\lambda - Q\rho_{\mathfrak{g}}) +Q\rho_{\mathfrak{g}}.
\ee

 Oscillatory representations form a subclass of  unitary highest weight representations. 
 
\subsubsection{Example.} For $\mathfrak{g} = {sl_2}$ these are the unitary representations 
 of the Virasoro Lie algebra  \cite[\S 3.4]{KRR}, parametrized by  the central charge $c$ and the highest weight $h$, such that  
\be
\begin{split}
 c\geq 1, \qquad h\geq (c-1)/24.
\end{split}
\ee

\begin{figure}[ht]
\begin{center}
\begin{tikzpicture}[scale=0.6]
\draw[dashed, directed] (0,0) --(8,0);
\draw[dashed] (0,0) -- (0,1);
\draw[directed] (0,1) -- (0,4);
\draw[fill=blue] (0,1) --(8,1) --(8,2)--(0,1);
\node at (-.3, -.3) {\small $0$};
\node at (-.3, 1) {\small $1$};
\node at (8.3, 0) {\small $h$};
\node at (0, 4.3) {\small $c$};
\node at (10, 2) {\small $h\geq (c-1)/24$};
\node at (8.8, 3) {\small $c\geq 1$};
\end{tikzpicture}
\end{center}
\caption{The oscillatory representations sector in the unitary rectangle $h\geq 0, c\geq 1$.}
\label{hsdiagram}
\end{figure}

Indeed, in this case $h^\vee=2$, so $c= 1+6Q^2$, and the   highest weight $h$  
 of the highest weight vector     
 of the oscillatory representation  $V_{\frac{1}{2}Q+i\alpha}$ of the Virasoro Lie algebra  is given by 
\be \la{hQ}
h = \frac{Q^2}{4}+ \alpha^2. 
\ee

We recall \cite[\S 12.5]{KRR} that 
  all unitary highest weight representations of the Virasoro Lie algebra are  the Verma modules with 
  $c \geq 1, h \geq 0$, and     the maximal quotients of the Verma modules  
  (the unitary  discrete series corresponding to the minimal models) with 
  \be \nonumber
  \begin{split}
&  c = 1-\frac{6}{(m+2)(m+3)}, \qquad  m=0,1,2, \ldots \\
&h = \frac{[(m+3)r-(m+2)s]^2-1}{4(m+2)(m+3)}, \qquad  r,s \in {\Bbb N}, \qquad 1 \leq s \leq r \leq m+1.\\
\end{split}
\ee

  \vskip 2mm
There are two kinds of data  on the space  $\widehat {\cal T}_{g,n}$, given by families of   $\widehat \Gamma_S$-equivariant  infinite dimensional  vector bundles with flat connections,   
parametrized by  the same space of parameters.

\subsubsection{De Rham data on  $\widehat {\cal T}_{g,n}$.}  Using the chiral algebra assigned to  ${\cal W}_{\mathfrak{g}}$, one can assign to a given    Riemann surface $\Sigma$  with punctures $p_1, ..., p_n$   
the  space of  coinvariants   for  any collection of  highest weight representations ${\rm V}_{\lambda_k}$ of  ${\cal W}_{\mathfrak{g}}$  attached to the punctures $p_k$,  denoted by
\be \la{COINV2}
  \left( {\rm V}_{\lambda_1} \otimes \ldots \otimes {\rm V}_{\lambda_n} \right)_{\Sigma, {\cal W}_{\mathfrak{g}}}, \qquad \lambda_k \in \mathfrak{h}_\C^\ast.
\ee 
It has a distinguished complex line  spanned     by the tensor product of 
   the highest weight vectors:
    \be \la{HWL}
   {\mathfrak{L}}_{\Sigma, \underline{\alpha}}= \langle v_{\lambda_1} \otimes \ldots \otimes v_{\lambda_m}\rangle_{\C}.
\ee
  The  coinvariants (\ref{COINV2})  give rise to 
  a $\widehat \Gamma_S-$equivariant   vector bundle  on the   Teichm\"uller space $\widehat {\cal T}_{g,n}$: 
\be \la{COINV3}
 \Delta_{{\cal W}_{\mathfrak{g}}} \left( {\cal V}_{\lambda_1} \otimes \ldots \otimes {\cal V}_{\lambda_n} \right) 
\ee  
We refer to it as the {\it de Rham bundle}.  
 The following  basic fact is well known, and can be found in \cite{FBZ}.
  \bt \la{FBZ} The  vector bundle of coinvariants (\ref{COINV3}) on the space $\widehat {\cal T}_{g,n}$  carries a $\widehat \Gamma_S-$equivariant  flat connection $\nabla_{\rm coinv}$.  
  \et
  
      We can restrict our attention to   oscillatory representations of the $W$-algebra ${\cal W}_{\mathfrak{g}^\vee}$, 
  assigning   to each puncture $p_k$   a highest weight of  the Langlands dual Lie algebra $\mathfrak{g}^\vee$ of the following shape:
  \be \la{lambdaalpha}
  \lambda_k = Q\rho   + \sqrt{-1}\alpha_k, \qquad \alpha_k \in   \mathfrak{h}_{\mathfrak{g}^\vee}^\ast, \qquad \rho:= \rho_{\mathfrak{g}^\vee} .
  \ee

\subsubsection{Betti data on  $\widehat {\cal T}_{g,n}$.}   Let $S$ be an oriented topological  surface   with $n$ punctures.   The quantization of the Betti moduli   space ${\rm Loc}_{\G, S}$, 
provided by  Theorem \ref{MTH}, delivers a cluster triple of spaces (\ref{Gtriple}),
as well as:

\vskip 1mm

1. A   representation of the  group $\widehat \Gamma_S$ in the Hilbert space ${\cal H}({\rm Loc}_{\G, S})$, preserving   ${\cal S}({\rm Loc}_{\G, S})$.
\vskip 1mm
 
2. A $\widehat \Gamma_S$-equivariant   representation  of the $\ast$-algebra $\mathcal{A}_\hbar({\rm Loc}_{\G, S})$ in the space ${\cal S}({\rm Loc}_{\G, S})$. 
\vskip 1mm

The 
Hilbert space ${\cal H}({\rm Loc}_{\G, S})$ is   decomposed  
into   an integral of representations parametrized by  the points  
  $\underline{\alpha} = (\alpha_1, ..., \alpha_n)$, where  $\alpha_{k} \in \H/W(\R_{>0})$, assigned to the punctures:
$$
{\cal H}({\rm Loc}_{\G, S}) = \int {\cal H}({\rm Loc}_{\G, S})_{\underline{\alpha}}d\underline{\alpha}, \qquad \underline{\alpha}    \in (\H/W)(\R_{>0})^{n}.
$$
So we get a triple of  representations of the group  $\widehat \Gamma_S$:  
\be \nonumber
  {\cal S}({\rm Loc}_{\G, S})_{\underline{\alpha}} \subset   {\cal H}({\rm Loc}_{\G, S})_{\underline{\alpha}} \subset   {\cal S}^*({\rm Loc}_{\G, S})_{\underline{\alpha}}. 
\ee
We  view them as $\widehat \Gamma_S-$equivariant  local systems of  infinite dimensional  vector spaces on the  
space  $\widehat {\cal T}_{g,n}$:
\be \nonumber
 \underline  {\cal S}({\rm Loc}_{\G, S})_{\underline{\alpha}} \subset  \underline  {\cal H}({\rm Loc}_{\G, S})_{\underline{\alpha}} \subset  \underline
  {\cal S}^*({\rm Loc}_{\G, S})_{\underline{\alpha}}. 
\ee
Abusing notation, we denote  the same way  infinite dimensional local systems on the moduli space $\widehat {\cal M}_{g,n}$. For example,  we  take  the product 
$\widehat {\cal T}_{g,n} \times {\cal S}({\rm Loc}_{\G, S})_{\underline{\alpha}}$,  viewed as a trivialized  flat bundle  over $\widehat {\cal T}_{g,n}$, and  consider its quotient under the action of the group $\widehat \Gamma_S$ provided by its representation in ${\cal S}({\rm Loc}_{\G, S})_{\underline{\alpha}}$:
$$
 \underline  {\cal S}({\rm Loc}_{\G, S})_{\underline{\alpha}} :=  \Bigl(\widehat {\cal T}_{g,n} \times  {\cal S}({\rm Loc}_{\G, S})_{\underline{\alpha}} \Bigr)/ \widehat\Gamma_S.
$$  
     \vskip 2mm
 Let us  relate  the labels of   de Rham and Betti  bundles.         
 Recall the isomorphisms 
\be \nonumber
\begin{split}
{\rm exp}:  ~&\mathfrak{h} \stackrel{\sim}{\lra}  \H(\R_{>0}), \\
&\mathfrak{h}  = \mathfrak{h}_{\mathfrak{g}^\vee}^\ast.\\
\end{split}
\ee
Combining them, we get a canonical isomorphism
\be \la{IDDE}
\mathfrak{h}_{\mathfrak{g}^\vee}^\ast = \H(\R_{>0}).
\ee
So we can match   unitary oscillatory
 representations  of the algebra 
${\cal W}_{\mathfrak{g}^\vee}$ to  the points of $\ \H(\R_{>0})$:
\be \la{ALA}
{\rm V}_{Q\rho+   \sqrt{-1}\alpha} \qquad \longleftrightarrow \qquad \alpha \in \H(\R_{>0}), \qquad \rho:= \rho_{{\mathfrak{g}}^\vee}.
\ee
 
Thanks to (\ref{WWW}),   de Rham and Betti  bundles are parametrized by the same space:     
$(\H/W)(\R_{>0})^n$.   The 
   homeomorphism $\varphi: S \to \Sigma$ defining a point of the Teichm\"uller space     identifies, via (\ref{IDDE}),  the  labels at the punctures on $S$ and $\Sigma$. 
The de Rham bundles depend on a central charge $c$, while the Betti  bundles depend on the Planck constant $\hbar$.

  \bcon \la{GENDKL}  Let $S$ be a genus $g$ surface with $n$ punctures. Then
  
 \begin{enumerate}
 
 \vskip 1mm\item    There is a non-degenerate  
   pairing of  vector bundles with flat connections on the   moduli space $\widehat {\cal M}_{g,n}$, 
        continuous on  $ \underline {\cal S}({\rm Loc}_{\G, S})_{\underline{\alpha}}$:
  \be \la{CONJIpair}
  {\cal C}_{\underline{\alpha}}:   \Delta_{{\cal W}_{\mathfrak{g}^\vee}} \left( {\cal V}_{\lambda_1} \otimes \ldots \otimes {\cal V}_{\lambda_n} \right) \bigotimes     
 \underline {\cal S}({\rm Loc}_{\G, S})_{\underline{\alpha}} \lra \C. 
\ee
Here $\lambda_i$ and $\alpha_i$ are related by (\ref{lambdaalpha}), and the central charge $c$  relates to the  parameter $Q$   by 
\be \la{huw1}
\begin{split} 
&c =
 {\rm rk}_{\mathfrak{g}}\left(1 +  h^\vee(h^\vee+1)Q^2\right).\\
\end{split}
 \ee
 
\vskip 1mm\item   The pairing (\ref{CONJIpair}) extends to  complex values of $\lambda_i$ and $\alpha_i$  related by (\ref{lambdaalpha}), and complex values of $c$ and $\hbar$ related by (\ref{huw1}), but may   become  degenerate for certain values of the parameters. 
  
\vskip 1mm\item  Since the  pairing   (\ref{CONJIpair}) is continuous on $\underline {\cal S}$, it  induces a map  of 
 vector bundles with flat connections:
\be \la{CONJI}
 {\cal C}_{\underline{\alpha}}:  \Delta_{{\cal W}_{\mathfrak{g}^\vee}}\left( {\cal V}_{\lambda_1} \otimes \ldots \otimes {\cal V}_{\lambda_n} \right)  \stackrel{}{\lra}    
\underline {\cal S}^*({\rm Loc}_{\G, S})_{\underline{\alpha}}.
\ee
The vectors in the image of the map (\ref{CONJI}) are the conformal blocks for the Toda theory. 
 \end{enumerate}

  \econ
  
In particular, the image of the  highest weight  vector in (\ref{HWL}) under the map (\ref{CONJI})  should provide a $\widehat \Gamma_S$-equivariant map to the vector space 
of distributions: 
\be \la{GEQM}
\begin{split}
&  {\rm C}_{\underline{\alpha}}:  \widehat {\cal T}_{g, n} \lra  {\cal S}^*({\rm Loc}_{\G, S})_{ \underline{\alpha} }.\\
&{\rm C}_{\underline{\alpha}}(\gamma (\tau)) =  \gamma \cdot  {\rm C}_{\underline{\alpha}}( \tau) \qquad \forall \gamma \in \widehat\Gamma_S, \qquad \forall \tau \in \widehat {\cal T}_{g,n}.\\
\end{split}
\ee  
The   element 
${\rm C}_{\underline{\alpha}}(\tau)$ 
  is the {conformal block}  for the primary fields in  the Toda theory.
 
 \vskip   2mm
 
  Let us elaborate  Conjecture  \ref{GENDKL},  starting with an example illustrating the nature of  coinvarints (\ref{COINV3}).
  
\subsubsection{Example.}    The sheaf of differential operators $ {\cal D}_X$   on a smooth variety   $X$ is 
    an infinite dimensional vector bundle with a flat connection. The connection is given by the   map 
    ${\cal T}_X \otimes {\cal D}_X \lra {\cal D}_X$, $\tau \otimes D\lms\tau D$, where ${\cal T}_X$ is the tangent sheaf.      
    The connection is non-integrable. It does not have a monodromy.

    The    coinvariants (\ref{COINV3})  for $\G = {\rm PGL_2}$ is  a twisted ${\cal D}-$module on  ${\cal M}_{g,n}$, twisted by the line bundle 
\be \la{chh}
 ({\rm det}~\Omega^1_{\Sigma})^c   \otimes  (T^*_{p_1}\Sigma)^{h_1} \otimes \ldots \otimes (T^*_{p_n}\Sigma)^{h_n}, \qquad [\Sigma] \in {\cal M}_{g,n}.
\ee
Here  $c$ is the central charge,  and $h_i$ the highest weight of the Verma module at the puncture $p_i$ on  $\Sigma$.  It is identified, as a vector bundle with a flat connection, 
     with the ${\cal D}-$module  
     of   twisted differential 
    operators on the line bundle (\ref{chh}). 
        
  \vskip 2mm
  
  Conjecture \ref{GENDKL} predicts a non-degenerate continuous pairing between two infinite dimensional 
  flat vector bundles  on $\widehat {\cal M}_{g,n}$. The first  carries a non-integrable 
  connection, and thus does not have a monodromy. The second has an interesting monodromy. How can this happen? Here is an explanation. 
  
    \subsubsection{Coinvariants and the cluster action of the mapping class group.} \la{6.0.8}
  Given a point $\tau \in \widehat {\cal T}_{g,n}$,   Conjecture \ref{GENDKL} predicts the existence of a canonical continuous pairing between the discreet space of coinvariants and the topological Frechet vector space    ${\cal S}({\rm Loc}_{\G, S; {\bf c}})$:  \be  \la{CONJIpairx}
  {\rm C}_{ {\bf c}, \tau}:  \left( {\rm V}_{\lambda_1} \otimes \ldots \otimes {\rm V}_{\lambda_n} \right)_{\Sigma_\tau, {\cal W}_{\mathfrak{g}^\vee}} \bigotimes     
 {\cal S}({\rm Loc}_{\G, S; {\bf c}}) \lra \C.
\ee 
As we discussed in Section \ref{1.1.2}, The space    ${\cal S}({\rm Loc}_{\G, S; {\bf c}})$ is 
 defined using a   combinatorial  data ${\bf c}$, which determines a cluster coordinate system on the space ${\mathscr P}_{\G,S}$. 
Changing the  data ${\bf c}$ to   ${\bf c}'$,   related by a cluster transformation ${\bf c} \to {\bf c}'$,  does not change the space of coinvariants, and amounts to the intertwiner 
$$
{\cal I}_{{\bf c} \to {\bf c}'}: {\cal S}({\rm Loc}_{\G, S; {\bf c}})  \lra {\cal S}({\rm Loc}_{\G, S; {\bf c}'}),
$$
which determines  
  the compatibility condition:
\be \la{PRNG}
{\rm C}_{ {\bf c'}, \tau} (v \otimes {\cal I}_{{\bf c} \to {\bf c}'}(s)) = {\rm C}_{ {\bf c}, \tau} (v \otimes s).
\ee
   
   Next, just as in Section \ref{1.0.7},         when $\tau$ varies, we have a family of 
 pairings  ${\rm C}_{{\bf c}, \tau}$, providing a family of maps from the \underline{constant}   space ${\cal S}({\rm Loc}_{\G, S; {\bf c}})$ to the variable dual space of invariants:
   \be \la{CONJIpairx2}
  {\rm C}'_{  {\bf c}, \tau}:       
 {\cal S}({\rm Loc}_{\G, S; {\bf c}})  \lra \left( {\rm V}^*_{\lambda_1} \otimes \ldots \otimes {\rm V}^*_{\lambda_n} \right)_{\Sigma_\tau, }^{{\cal W}_{\mathfrak{g}^\vee}}. 
\ee  
So any vector $s \in {\cal S}({\rm Loc}_{\G, S; {\bf c}})$ gives rise to a section ${\cal C}'_{  {\bf c}, \tau}(s)$ of the vector bundle of invariants  
 \be \la{COINV2a}
  \left( {\cal V}^*_{\lambda_1} \otimes \ldots \otimes {\cal V}^*_{\lambda_n} \right)^{{\cal W}_{\mathfrak{g}^\vee}}.
\ee
Denote by $\nabla_{\rm inv}$ the flat connection on this  bundle dual to the connection $\nabla_{\rm coinv}$ on the coinvariants from  Theorem \ref{FBZ}.     
Conjecture \ref{GENDKL}  implies that   the section ${\cal C}'_{{\bf c}, \tau}(s)$  
 is flat:
\be \la{Gamaeq}
\begin{split}
&\nabla_{\rm inv}{\cal C}'_{  {\bf c}, \tau}(s)=0, \qquad \forall s \in {\cal S}({\rm Loc}_{\G, S; {\bf c}}).\\
\end{split}
\ee
Although  the connection $\nabla_{\rm coinv}$  is not   integrable,   the map (\ref{GEQM})  is $\widehat \Gamma_S-$equivariant. 
 Indeed, 
the  group $\widehat \Gamma_S$ acts  on the  combinatorial data $\{{\bf c}\}$. 
  Let $\gamma$ be a loop on   $\widehat {\cal M}_{g,n}^*$, and $[\gamma]\in \widehat \Gamma_S$ its homotopy class. 
  Let $\widetilde \gamma: \tau \to \tau'$  be a lift of $\gamma$ to  $\widehat {\cal T}_{g,n}$.   Moving along the path 
  $\widetilde \gamma$ we change the  data ${\bf c}$ to  $[\gamma]({\bf c})$. So it amounts to a pairing
  \be \la{CONJIpairx3}
 {\cal C}_{  [\gamma]{\bf c}, \tau}:  \left( {\rm V}_{\lambda_1} \otimes \ldots \otimes {\rm V}_{\lambda_n} \right)_{\Sigma_\tau, {\cal W}_{\mathfrak{g}^\vee}}    \bigotimes        
 {\cal S}({\rm Loc}_{\G, S; [\gamma]{\bf c}}) \lra \C. 
\ee 
Thanks to (\ref{PRNG}), this is equivalent to the action of the intertwiner ${\rm I}_{{\bf c} \to [\gamma]{\bf c}}$.  This imples  (\ref{GEQM}).

\bc  Conjectures \ref{MK} $\&$ \ref{GENDKL} + Theorem \ref{MFCQI} imply that one should have  a non-degenerate pairing between the coinvariants of  oscillatory representations of the $W$-algebra ${\cal W}_{{\mathfrak{g}}^\vee}$ and   invariants 
 of the modular quantum group  $\mathcal{A}_\hbar({\mathfrak{g}})$, continuous on the second factor and invariant under the  braid group action:
 \be \la{CONJQG}
  \left( {\rm V}_{\lambda_1} \otimes \ldots \otimes {\rm V}_{\lambda_n} \right)_{\C\P^1-\{z_1, ..., z_n\}, {\cal W}_{{\mathfrak{g}^\vee}}} \bigotimes   
  \Bigl({\cal S}_{\alpha_1}\otimes  \ldots \otimes {\cal S}_{\alpha_n}\Bigr)^{\mathcal{A}_\hbar(\mathfrak{g})} \lra \C.
\ee
  \ec

\vskip 2mm
\subsubsection{The general irregular case.}

  Let $S$ be a genus $g$ surface with $n$ punctures. Let $\beta_p$ be an element of the cyclic closure of the positive braid semigroup ${\Bbb B}^+_{\mathfrak{g}}$ 
 describing the irregularity type of a  $\G-$bundle with meromorphic connection on a Riemann surface  of type $S$  at a puncture $p$, see Section \ref{5.3.5} and 
 \cite[Section 11.1]{GK}. 
Recall the Betti stack of framed Stokes data of type $\beta = \{\beta_p\}$:
 \be
 {\cal M}_{\B}(\G, S; \beta).
 \ee
 It  carries a cluster Poisson structure, equivariant under the action of the {\it wild mapping class group} $\Gamma_{\G, S; \beta}$. 
 Its Poisson center is described by $\Gamma_{\G, S; \beta}$-equivariant   projection to a Casimir torus ${\rm C}_{\G, S; \beta}$:
 \be
{\cal M}_{\B}(\G, S; \beta) \lra {\rm C}_{\G, S; \beta}.
  \ee
 A description of the   torus ${\rm C}_{\G, S; \beta}$ and this projection follows from the fact that   $ {\cal M}_{\B}(\G, S; \beta)$ has a cluster Poisson structure. 
 In particular, the  torus ${\rm C}_{\G, S; \beta}$ is a product over the punctures $p$ of $S$ of  tori ${\rm C}_{\G, S; \beta_p} $:
  \be
 {\rm C}_{\G, S; \beta} = \prod_p  {\rm C}_{\G, S; \beta_p}.
 \ee
Given an  $\alpha \in  {\rm C}_{\G, S; \beta}(\R_{>0})$, the quantized Betti moduli space ${\cal M}_{\B}(\G, S; \beta)$ provides  a cluster Schwarz   space 
 $$
  \underline {\cal S}({\cal M}_{\B}(\G, S; \beta))_{ {\alpha}}
  $$
     with a unitary representation of the extended wild mapping class group. 
 
 Take a $\lambda \in  {\rm C}_{\G, S; \beta}(\C)$ related to $\alpha$ as in (\ref{lambdaalpha}), with the central charge $c$  related to $\hbar$   by (\ref{huw1}). 
 Generalising Conjecture \ref{GENDKL}, we have   the following Conjecture. 
 
\bcon \la{GENDKL1}  

 i) There is a natural infinite dimensional flat bundle   
 ${\rm CB}_{\G, S; \lambda}$  of irregular conformal blocks of type $\beta_p$ at the punctures $p$ over a   
 complex torus bundle $\widehat {\cal P}_{\G, S; \beta}$ over  
 ${\cal P}_{\G, S, \beta}$, similar to (\ref{TORUSB}).

\vskip 1mm
 
  ii) There is a non-degenerate  
   pairing of  vector bundles with flat connections on  $\widehat {\cal P}_{\G, S; \beta}$, 
        continuous on  $ \underline {\cal S}({\cal M}_{\B}(\G, S; \beta))_{ {\alpha}}$:
  \be \la{CONJIpairI}
  {\cal C}_{\underline{\alpha}}:   {\rm CB}_{\G, S; \lambda} \bigotimes     
 \underline {\cal S}({\cal M}_{\B}(\G, S; \beta))_{ {\alpha}} \lra \C. 
\ee
  \econ

 A number of comments are in order.

    
\vskip 2mm \paragraph{\bf 1. }  
Conjecture \ref{GENDKL} relates the Toda theory to the quantized   moduli space ${\rm Loc}_{\G, S}$. To state it precisely it is 
    not enough  to have just a unitary representation of the group $\widehat \Gamma_\bS$ in a Hilbert space ${\cal H}({\rm Loc}_{\G, S})$;  we need   the whole quantization data. Indeed, we need 
      a representation of the algebra $\mathcal{A}_\hbar(\G, S)$ by unbounded operators in the space ${\cal H}({\rm Loc}_{\G, S})$   to define the  subspace ${\cal S}({\rm Loc}_{\G, S})$. And we must know that the action of the group   $\widehat \Gamma_\bS$ preserves this subspace.

\vskip 2mm \paragraph{\bf 2.} 
 To describe coinvariants of arbitrary  unitary oscillatory representations we need to quantize the moduli space $ {\rm Loc}_{\G, S}$ for  both $\hbar>0$ and $|\hbar|=1$, i.e. for any $Q>0$.
  For  example, for    $sl_2$
$$
c_{sl_2}=  1-  6\frac{(k+1)^2}{k+2} = 1+6Q^2 = 13+6(\hbar + \hbar^{-1}).
$$
So if   $\hbar>0$, then $c>  25$.  
If  $|\hbar|=1$, then $1 \leq c \leq 25$. 
So to get the oscillatory series of unitary Verma modules $M_{h, c}$ with $c \geq 1, h \geq \frac{c-1}{24}$ we need the  $|\hbar|=1$ case.

\vskip 2mm \paragraph{\bf 3.} 
The names de Rham and Betti   reflect the fact that de Rham data are defined by using   a Riemann surface, while  Betti data are
  defined 
 using just a topological surface.

Recall that a weight $n$ variation of Hodge structures ${\cal H}$ is Calabi-Yau if $h^{n,0}=1$. Then   $k$-th covariant derivatives of the line $F^{n} {\cal H}$ lie in $F^{n-k} {\cal H}$. 

The de Rham and Betti bundles with flat connections should be realizations of a variation of   {\it infinite dimensional Calabi-Yau motives}.  Indeed, the map (\ref{CONJI}) should be  the comparison map between  the de Rham  and  Betti realizations. 
The map (\ref{GEQM}) should describe the deepest $F^{\frac{\infty}{2}, 0}{\cal H}$ part of the Hodge filtration. The Hodge filtration should be recovered by covariant derivatives of   section (\ref{GEQM}). 
 The periods of this motive  should describe   Feynman integrals for the correlation functions in the Toda theory: 
 $$
 \int e^{-\frac{1}{\hbar}{\cal S}(\varphi)}{\cal V}_{z_1}(\varphi) \ldots {\cal V}_{z_n}(\varphi) {\cal D}\varphi.
 $$
 Here the space of fields $\{\varphi\}$ consists of functions on $\Sigma$ with values in the Cartan Lie algebra $\mathfrak{h}$. 
 Oscillatory representations of the W-algebra are the spaces of states in the Toda theory. Insertions ${\cal V}_{z_1}(\varphi) \ldots {\cal V}_{z_n}(\varphi)$ in the Feynman integral 
 correspond to the vectors in the space of coinvariants.   The highest weight vector of the representation ${\rm V}_\lambda$    corresponds to the insertion of 
 $e^{\alpha \cdot \varphi}$, where $\alpha$ and $\lambda$ are related by (\ref{lambdaalpha}).

\vskip 2mm \paragraph{\bf 4.} 
Let us explain  why  the $\widehat \Gamma_S$-equivariant map (\ref{GEQM}) predicted   by Conjecture \ref{GENDKL}      should exist.

  Given a Riemann surface  $\Sigma \in {\cal M}_{g, n}$,  a group $\G$, and $\underline{\alpha}$, there is a subspace of opers 
    $$
    {\rm Oper}_{\G, \Sigma; \underline{\alpha}}\subset {\rm Loc}^{\rm DR}_{\G, \Sigma; \underline{\alpha}}
    $$ 
    in the de Rham moduli space of holomorphic $\G$-connection on $\Sigma$, whose  singularities at the punctures are specified by the weights  $\underline{\alpha}$. 
    Given  a marking $\varphi: S \stackrel{\sim}{\lra} \Sigma$ describing a point $\tau \in {\cal T}_{g,n}$, the monodromies   
    of opers  provide a   subspace in the Betti moduli space:
\be \la{LAGa}
   \varphi^*{\rm Oper}_{\G, \Sigma; \underline{\alpha}} \subset {\rm Loc}^{\rm Betti}_{\G, S; \underline{\alpha}}.
\ee
 The   space   on the right is symplectic. The subspace of opers  is Lagrangian. 
   We expect that this Lagrangian   can be quantized, providing a line ${\rm L}_\tau\subset {\cal S}^*({\rm Loc}_{\G, S})_{\underline{\alpha}}$. 
   Then the image of the highest weight line ${\rm C}_{\underline{\alpha}}  ({\mathfrak{L}}_{\Sigma, \underline{\alpha}})$ in (\ref{HWL})  should be  
 this line, and  
     the projectivization of the canonical map   (\ref{GEQM}) is given by 
   $
    \tau \in \widehat {\cal T}_{g,n} \lms    {\rm L}_\tau.
  $

\vskip 2mm \paragraph{\bf 5.} 
A key    feature of  the map (\ref{CONJI}) is that the  quantization  spaces ${\cal H}({\rm Loc}_{\G, S}) \subset {\cal S}^*({\rm Loc}_{\G, S})$  are cluster vector spaces, and therefore 
\underline{do not} have a unique  realization. For example, ${\cal H}({\rm Loc}_{\G, S})$ 
is   presented by   an infinite collection of  Hilbert spaces related by   unitary projective transformations. Precisely, it is a functor 
\be \la{FUNCTa}
   {\cal H}_{\G, S}: {\rm M}_{\G, S} \lra  {\rm Hilb}
\ee
 from a combinatorial   {\it cluster modular groupoid} ${\rm M}_{\G, S}$  of the cluster Poisson space ${\mathscr P}_{\G, S}$  to the category of Hilbert spaces ${\rm Hilb}$ with  morphisms given by  unitary projective  maps. Indeed, to construct a   realization of the  space ${\cal H}({\rm Loc}_{\G, S})$ we pick a combinatial data ${\bf c}$ providing a cluster 
 Poisson coordinate system on the  space ${\mathscr P}_{\G, S}$, defining a point of the groupoid ${\rm M}_{\G, S}$.   
 
The mapping class group $\Gamma_S$ maps  to the  group $\pi_1({\rm M}_{\G, S})$. 
 This   way the    group $\Gamma_S$ acts by unitary projective transformations on the triple of  quantization spaces. 
   One   imagines that   functor (\ref{FUNCTa})   describes   a local system on a space by  using an infinite 
   collection of base points and canonical paths between them.

 \vskip 2mm \paragraph{\bf 6.} 
 The Hilbert spaces ${\cal H}({\rm Loc}_{\G, S})$ for $\hbar$ and ${\cal H}({\rm Loc}_{\G^\vee, S})$ for $\hbar^{-1}$ are canonically identified.  
  
  On the other hand, the involution $\hbar \to \hbar^{-1}$ amounts to the involution $(k+h^\vee) \to (k+h^\vee)^{-1}$ which appears in the Feigin-Frenkel duality 
for representations of   ${\cal W}_{\mathfrak{g}}$ and ${\cal W}_{\mathfrak{g}^\vee}$:
$$
\hbar \lra \hbar^{-1} \qquad \longleftrightarrow \qquad (k+h^\vee) \to (k+h^\vee)^{-1}. 
$$

{\it So the Langlands type duality built into the construction of    the spaces ${\cal H}({\rm Loc}_{\G, S})$ corresponds under   (\ref{CONJI}) 
 to the Feigin-Frenkel  duality for 
    representations of $W$-algebras.}

\vskip 2mm \paragraph{\bf 7.} 
The Hilbert spaces ${\cal H}({\rm Loc}_{\G, S; \hbar})$ for $\hbar$ and the dual space ${\cal H}({\rm Loc}_{\G, S^\circ; -\hbar})^*$ for $-\hbar$ can be    identified. Here $S^\circ$ is the surface $S$ with the opposite orientation. 

The involution  $\hbar \to -\hbar$ amounts to the involution $(k+h^\vee) \to -(k+h^\vee)$, see (\ref{hbark}), 
and  the involution 
\be \la{hbarFF}
c_\mathfrak{g} \longleftrightarrow 2{\rm rk}_{\mathfrak{g}}(1+2h^\vee(h^\vee+1)) - c_\mathfrak{g}. 
\ee
For example, involution  (\ref{hbarFF}) for $\frak{g}=sl_2$ amounts to   the Feigin-Fuchs  
involution\footnote{The Feigin-Fuchs duality is stated in  \cite{FFu}  as $(c, h')  \longleftrightarrow (26-c, -1-h')$. 
It uses  the generators  $L'_n = z^{n+1}d/dz$ of the Virasoro algebra. The generators   used in conformal field theory      
differ by the sign: $L_n =  - z^{n+1}d/dz$.  For example,  Virasoro  Verma modules  $M_{c, h}$  are unitarizable if   $c \geq 1, h \geq 0$. So $h'=-h$, and $h'\to-1-h'$ corresponds to $h \to 1-h$. }
\be\la{FFIN}
(c, h)  \longleftrightarrow (26-c, 1-h)
\ee 
which appears in the contravariant duality for representations of the Virasoro algebra \cite{FFu}.  
 Indeed, if $\frak{g}=sl_2$, then 
 $c = 13 + 6(\hbar+\hbar^{-1})$,   and   $\hbar \to -\hbar$  amounts to  $c \to 26-c$.

Involution (\ref{hbarFF}) should correspond to a contravariant duality on representations of  ${\cal W}_{\mathfrak{g}}$.  

\vskip 2mm

{\it So the  $\hbar \longleftrightarrow -\hbar$ duality   for the  spaces ${\cal H}({\rm Loc}_{\G, S})$ should correspond under   (\ref{CONJI}) 
 to the Feigin-Fuchs type  duality for 
    representations of $W$-algebras.}

\vskip 2mm \paragraph{\bf 8.} %
 Conjectural pairing (\ref{CONJQG}) is merely  a  combination of   isomorphism (\ref{TENV})  
 with    conjectural pairing (\ref{CONJIpair}) for $\C\P^1-\{z_1, ..., z_n\}$. 
  One might think about it as a conjectural   analog of the Kazhdan-Lusztig theory \cite{KL} for the principal series of representations of quantum group. 
  The Kazhdan-Lusztig theory  relates the modular tensor category of 
   finite dimensional  representations of ${\U}_q({\mathfrak{g}})$ to the fusion category of integrable representations of affine Kac Moody algebras at the   negative level.

\vskip 2mm \paragraph{\bf 9.} 
Does the space of coinvariants 
$\left( {V}_{\lambda_1} \otimes \ldots \otimes {V}_{\lambda_n} \right)_{{\cal W}_{\mathfrak{g}}^\vee}$ has a unitary structure for a unitary oscillatory representations $V_\lambda$?  
Note that   quotients of an infinite dimensional  
Hermitian vector space  do not necessarily inherit a unitary structure. And the target of the map (\ref{CONJI}) is not a Hilbert space.

    Conformal blocks, which a priori lie in the space ${\cal S}^*({\rm Loc}_{\G, S})$, 
  may not lie in the Hilbert space ${\cal H}({\rm Loc}_{\G, S})$. If they are, their norm is the correlation function. Otherwise 
it should be regularised. 
Assume that the map (\ref{CONJI}) 
 lands in the Hilbert space  ${\cal H}({\rm Loc}_{\G, S})_{\underline{\alpha}}$,  providing a   
       map of  $\widehat \Gamma_S$-equivariant   vector bundles on  $\widehat {\cal T}_{g,n}$:
  \be \nonumber
  {\rm C}_{\underline{\alpha}}:  \Delta_{{\cal W}_{\mathfrak{g}^\vee}} \left( {\cal V}_{\lambda_1} \otimes \ldots \otimes {\cal V}_{\lambda_n} \right) \stackrel{}{\lra}     
  {\cal H}({\rm Loc}_{\G, S})_{\underline{\alpha}}.
\ee
Then the   image of the  highest weight  line (\ref{HWL}) under  this map   gives a $\Gamma_S$-equivariant map
\be \la{GEQMa}
  {\cal C}_{\underline{\alpha}}:  \widehat {\cal T}_{g, n} \lra \C{\P}({\cal H}({\rm Loc}_{\G, S})_{ \underline{\alpha} }).
\ee

Recall that for any  complex vector space $V$ with a positive definite hermitian metric,  
 the projective space $\C\P(V)$ of lines in $V$ 
 carries  a canonical functorial Fubini-Studi K\"ahler form $\Omega_{\rm FS}$: for any subspace $V'\subset V$, 
 the  Kahler form on $\C\P(V')$  for the   induced metric on $V'$ coincides with the Kahler form induced by the embedding $\C\P(V')  \to \C\P(V)$. 
There is the tautological  Hermitian line bundle  
${\mathscr L} \to \C\P(V)$.

The pull back ${\cal C}^*_{\underline{\alpha}}({\mathscr L})$ 
  of the tautological  line bundle is a $\Gamma_S$-equivariant Hermitian line bundle on $\widehat {\cal T}_{g, n}$. 
  Its curvature    is  a Kahler form  on  $\widehat {\cal T}_{g, n}$ - the pull back   ${\cal C}^*_{\underline{\alpha}}(\Omega_{\rm FS})$  of the Fubini-Studi form.   
 One should have a canonical potential ${\rm C}_{g,n; \underline{\alpha}}$ for the Kahler form  ${\cal C}^*_{\underline{\alpha}}(\Omega_{\rm FS})$,   
given by the logarithm of the norm  of a   $\Gamma_S$-equivariant section ${s}$ of the line bundle ${\cal C}^*_{\underline{\alpha}}({\mathscr L})$   on the Teichm\"uller space:
\be \la{huw1ax}
\begin{split}
&{\rm C}_{g,n; \underline{\alpha}}:=\log ||{s}||, \qquad \overline \partial \partial \log ||{s}|| = {\cal C}^*_{\underline{\alpha}}(\Omega_{\rm FS}).\\
\end{split}
 \ee
The function ${\rm C}_{g,n; \underline{\alpha}}$ descends to   the  space ${\mathscr L}^*_{g, n} $. It is the correlation function for Toda theory.

\vskip 2mm \paragraph{\bf 10.} 
The key case of Conjecture \ref{GENDKL} is   when $S$  is a pair of pants ${\mathscr P}_3$, and  $\Sigma =  \C\P^1-\{z_1, z_2, z_3\}$. 

In this case Conjecture \ref{GENDKL}, combined with the realization construction from Section \ref{sec16.5:RMSP},   
 reduces to  the following conjecture  
relating the   coinvarinats of three oscillatory representations of the $W-$algebra ${\cal W}_{{\mathfrak{g}^\vee}}$ and the space of functions on   configurations   of triples of decorated flags. 
To state it,    recall canonical projection,   see Section \ref{SECC1.3}: 
   $$
 {\rm Conf}_3({\mathcal A}_\G ) \lra \H_\G^3.
  $$
 Passing to the positive  points, and taking the fiber over a point $(\lambda_1, \lambda_2, \lambda_3) \in \H_\G(\R_{>0})^3$, we get a subspace 
 ${\rm Conf}_3({\mathcal A}_\G )_{{\lambda_1}, {\lambda_2}, {\lambda_3}} \subset {\rm Conf}_3({\mathcal A}_\G(\R_{>0}))$.

\bcon \la{A3C} For any $({{\lambda_1}, {\lambda_2}, {\lambda_3}} ) \in H_\G^3(\R_{>0})$, there exists a natural non-degenerate pairing
 \be \la{NPzz}
  \left( {\rm V}_{\lambda_1} \otimes {\rm V}_{\lambda_2}  \otimes {\rm V}_{\lambda_3} \right)_{\C\P^1-\{z_1, z_2, z_3\}, {\cal W}_{{\mathfrak{g}^\vee}}} 
  \bigotimes {\cal S}({\rm Conf}_3({\mathcal A}_\G )_{{\lambda_1}, {\lambda_2}, {\lambda_3}}) \lra \C.
 \ee
 \econ
 
 If $\G = {\rm PGL}_2$, the space ${\rm Conf}_3({\mathcal A}_\G )_{{\lambda_1}, {\lambda_2}, {\lambda_3}}$ is just a point, and 
 the space of coinvariants for   the Virasoro  algebra $({\rm V}_{\lambda_1} \otimes {\rm V}_{\lambda_2} \otimes {\rm V}_{\lambda_3})_{\rm Vir}$  is one-dimensional. 
 So  Conjecture \ref{A3C} is valid. The  pairing (\ref{NPzz}),  suitably normalized,  is given by the  three-point correlation function for the Liouville model   \cite{DO}, \cite{ZZ}. 
 
 Taking a decomposition  of the surface $S$ into pair of pants, and combining Conjecture \ref{A3C} with the Modular Functor Conjecture \ref{MK}, we arrive at a conjectural presentation 
 of the Hilbert space ${\cal H}_{\G, S}$ as an integral of   Hilbert spaces    assigned to triangles related to the pair of pants  decomposition.

\vskip 2mm \paragraph{\bf 11.} 
If $\G = {\rm PGL_2}$, a pair of pants decomposition $\alpha$ of $S$ gives rise to a Hilbert space ${\rm L}_2^\sigma:= {\rm L}^\sigma_2(\R_+^{3g-3+n})$. The   coordinates 
in $\R_+^{3g-3+n}$ match the logarithms of  lengths of the geodecics on $\Sigma$ realizing $\alpha$. It is identified by  \ref{MK} with 
${\cal H}({{\rm PGL_2}, S})$.   The unitary action of an element  $\gamma \in \Gamma_S$ should correspond to Teschner's unitary operators 
$ {\rm L}^\sigma_2   \to {\rm L}^{\gamma(\sigma)}_2 $  \cite{T}. On the other hand,  gluing  conformal blocks according to the decomposition 
$\alpha$ \cite[5.2.1]{T10} provides 
a collection of invariants parametrised by the same   $\R_+^{3g-3+n}$. One can view this  \cite{T10}   as a strong indication of the existence of a map  
\be \nonumber
 {\cal H}( {\rm PGL_2}, S)  \lra \left( {\rm V}^*_{\lambda_1} \otimes \ldots \otimes {\rm V}^*_{\lambda_n} \right)_{\Sigma}^{{\rm Vir}}. 
\ee  
This is  stronger   than (\ref{CONJIpairx2}), where we conjecture only a  similar map of the   subspace ${\cal S}( {\rm PGL_2}, S)$. 



\medskip

\section{Cluster Poisson structure of the  space ${\mathscr P}_{\rm PGL_m, \bS}$ via amalgamation} \la{SSECC3}

\medskip

\subsection{Cluster Poisson amalgamation for the space ${\mathscr P}_{\rm PGL_2, \bS}$} \la{SEC5.100}

\medskip

Let $V_2$ be a two dimensional vector space. Recall the cross-ratio (\ref{CRRA}): 
\[
r(x_1, x_2, x_3, x_4):= 
\frac{\omega (x'_1,x'_2)\omega(x'_3,x'_4)}{\omega(x'_1,x'_4)\omega(x'_2,x'_3)}.
\]
 
\begin{figure}[ht]
\epsfxsize 200pt
\center{
\begin{tikzpicture}[scale=0.5]
\foreach \count in {1,2,3}
{
\node at (90-120*\count:2.3) {${x}_\count$};
\draw [dashed, red, arrows={-Latex[length=2mm 3 0]}]  (90-120*\count:2cm) -- (-30-120*\count:2cm);
}
\node at (30:1.4) {$r$};
\node at (150:1.4) {$q$};
\node at (270:1.4) {$p$};
\begin{scope}[shift={(8,0)}]
\foreach \count in {1,2,3}
{
\draw [dashed, red, arrows={-Latex[length=2mm 3 0]}]  (90-120*\count:2cm) -- (-30-120*\count:2cm);
}
\node at (30:1.4) {${\rm X}_r$};
\node at (150:1.4) {${\rm X}_q$};
\node at (270:1.4) {${\rm X}_p$};
\end{scope}
\end{tikzpicture}
 }
\caption{A triangle with pinnings   encoding a 
configuration $(x_1, x_2, x_3, p,q,r)$ of 6 points on 
${\mathbb P}^1$,  described by coordinates $X_p, X_q, X_r$. A pinning is shown by an oriented edge with vertices   labelled by points, 
and a middle pinning point. }
\label{fga8}
\end{figure}
Recall, see Section \ref{3.4.2},  that a pinning $\pi$ for ${\rm PGL_2}$  can be described by a triple of distinct points on ${\mathbb P}^1$:
$$
\pi = (x_1, x_2; p).
$$ 
The points $x_i$ are the flags.   
The dual pinning $\pi^*$ can be described as a triple 
$$
\pi^*:=(x_2, x_1; p^*),
$$
 where    $p^*$ is the unique point on  ${\mathbb P}^1$ such that 
\be \la{12.21.08.1x}
r(x_1,p, x_2,p^*)=1.
\ee

Our elementary object is a triple of flags with pinnings,  shown on the left of 
Figure \ref{fga8}. It  encodes   a 
configuration of $6$ points $(x_1, x_2, x_3, p,q,r)$ on 
${\mathbb P}^1$. We assign the coordinates $X_r, X_p, X_q$  to the 
 pinnings:
\be  \la{12.21.08.4}
X_r:= r(x_1, x_2, x_3, r), \quad X_p:=  r(x_2, x_3, x_1, p), \quad X_q :=
r(x_3, x_1, x_2, q).
\ee

\begin{figure}[ht]
\epsfxsize 200pt
\center{
\begin{tikzpicture}[scale=1.0]
\draw (1,-1) -- (0,0) -- (1,1);    
\draw (2,-1) -- (3,0) -- (2,1); 
\draw[dashed, red, arrows={-Latex[length=2mm 3 0]}] (1,1) -- (1,-1);
\draw[dashed, red, arrows={-Latex[length=2mm 3 0]}] (2,-1) -- (2,1);
\node[label=left:$x_2$] at (0,0) {};
\node[label=above:$x_3$] at (1,1) {};
\node[label=right:$p$] at (0.8,-0.2) {};
\node[label=below:$x_1$] at (1,-1) {};
\node[label=right:$y_2$] at (3,0) {};
\node[label=above:$y_3$] at (2,1) {};
\node[label=left:$q$] at (2.2,0.2) {};
\node[label=below:$y_1$] at (2,-1) {};
\begin{scope}[shift={(5,0)}]
\draw (1,-1) -- (0,0) -- (1,1);    
\draw (2,-1) -- (3,0) -- (2,1); 
\draw[dashed, red, arrows={-Latex[length=2mm 3 0]}] (1,-1) -- (1,1);
\draw[dashed, red, arrows={-Latex[length=2mm 3 0]}] (2,-1) -- (2,1);
\node[label=left:$x_2$] at (0,0) {};
\node[label=above:$x_3$] at (1,1) {};
\node[label=right:$p^*$] at (0.8,0.2) {};
\node[label=below:$x_1$] at (1,-1) {};
\node[label=right:$y_2$] at (3,0) {};
\node[label=above:$y_3$] at (2,1) {};
\node[label=left:$q$] at (2.2,0.2) {};
\node[label=below:$y_1$] at (2,-1) {};
\end{scope}
\begin{scope}[shift={(6,0)}]
\draw (5,-1) -- (4,0) -- (5,1) -- (6,0) -- (5,-1);    
\draw[thick] (5,-1) -- (5,1); 
\node[label=left:$x_2$] at (4,0) {};
\node[label=above:$x_3$] at (5,1) {};
\node[label=right:$x_4$] at (6,0) {};
\node[label=below:$x_1$] at (5,-1) {};
\end{scope}
\end{tikzpicture}
 }
\caption{Gluing two triangles with pinnings. First, we replace the pinning $(x_3, x_1; p)$ by the dual pinning $(x_1, x_3; p^*)$, and then glue the 
two triangles, decorated by the points,   matching the pinnings $(x_1, x_3; p^*)$ and $(y_1, y_3; q)$.}
\label{fga7}
\end{figure}
Consider configurations of 
points $(x_1, x_2, x_3, p)$ and $(y_1, y_2, y_3, q)$ on the two copies ${\mathbb P}_x^1$ and   ${\mathbb P}_y^1$ of the projective line, 
see Figure \ref{fga7}. The points $(x_1, x_2, x_3)$ 
are assigned to the vertices of a triangle $\tau_x$, and the point $p$ determines  
a pinning  $(x_3, x_1; p)$. Similarly the points $(y_1, y_2, y_3)$ 
are assigned to the vertices of another triangle $\tau_y$, and the point $q$ determines 
a pinning  $(y_1, y_3; q)$. Let $p^*$ be the unique point on  ${\mathbb P}_x^1$ such that 
(\ref{12.21.08.1x}) holds. 
Then there is a unique isomorphism $\varphi: {\mathbb P}_y^1\lra {\mathbb P}_x^1$ such that 
$$
\varphi(y_1) = x_1, \ \ \varphi(y_3) = x_3, \ \ \varphi(q) = p^*.
$$
We set   $x_4:= \varphi(y_4)\in {\mathbb P}^1_x$. So the configuration of points $(x_1, x_2, x_3, x_4)$  
is obtained by gluing  configurations $(x_1, x_2, x_3)$ and $(y_1, y_2, y_3)$ matching the pinnings 
$(x_1, x_3; p^*)$ and $(y_1, y_3; q)$. We assign the points  $(x_1, x_2, x_3, x_4)$ to the vertices of rectangle $r$ obtained by gluing   triangles $\tau_x$ and $\tau_y$, so that the edge $E$, obtained by gluing the edges of the triangles, carries the points $x_1, x_3$. 
We assign to the edge $E$ the coordinate
$$
X_{E}:= r(x_1, x_2, x_3, x_4).
$$

\bl \la{12.21.08.3}
The $X$-coordinate $X_E$  
is the product of the ones assigned to the 
pinnings of the  triangles: 
$$
X_{E} = X_pX_q.
$$
\el

\begin{proof} 
 One has to prove that 
$
r(x_1, x_2, x_3, p)r(x_3, x_4, x_1, q) = r(x_1, x_2, x_3, x_4). 
$ 
This follows from the identity
\be
\begin{split}
&r(x_1, x_2, x_3 ,p) r(x_3, x_4, x_1, q) r(x_1, p, x_3, q) =  r(x_1, x_2, x_3 ,x_4),  \\
\end{split}
\ee
which is valid   for any points $x_1, ..., x_4, p, q$, and easy to check. \end{proof}
  
Given a decorated surface $\bS$, 
let us define a coordinate system on the moduli space 
${\mathscr P}_{\rm PGL_2, \bS}$.    
Pick  an ideal triangulation 
$T$ of $\bS$.  
Each boundary interval of the triangulation carries a pinning.  
Let us assign to  {every} edge $E$ of   $T$ a rational function  $X^T_E$ 
on the moduli space ${\mathscr P}_{\rm PGL_2, \bS}$. If $E$ is an internal edge,  go to the universal 
cover of $\bS$, lift the triangulation to a triangulation $T'$ of the cover, 
take an edge $E'$ projecting to $E$, and a rectangle $r_{E'}$ it determines there. 
The vertices of the rectangle are decorated by a configuration of flags 
$(x_1, x_2, x_3, x_4)$ provided by the framing on $\bS$, 
so that $(x_1x_3)$ decorates $E'$. Set 
$X_E:= r(x_1, x_2, x_3, x_4)$. 

Now let $E$ be a boundary edge with a pinning $p_E$.
We go  to the universal cover and take a lift $E'$ of the edge. 
It  is contained in the unique ideal triangle $t_{E'}$, and there is a pinning 
assigned to the edge $E'$, provided by the pinning on $\bS$.  We assign 
to the pair $(t_{E'}, p_{E'})$ 
the corresponding coordinate  (\ref{12.21.08.4}), and declare it 
the coordinate $X^T_{E}$. 
  If the rectangle $r_{E'}$ or triangle $t_{E'}$ 
projects isomorphically onto 
the surface, it is not necessary to go to the universal cover to define the coordinate -- 
we can  use the ideal rectangle/triangle on $\bS$. 

\vskip 1mm
\bt \la{12.21.8.10}
Let $T$ be  an ideal triangulation of  $\bS$.

\begin{enumerate}

\vskip 1mm \item    
Rational functions $\{X^T_E\}$ 
are amalgamations 
of the ones assigned to triangles of $T$.

\vskip 1mm\item  Rational functions $\{X^T_E\}$ 
provide a rational coordinate system on the space  ${\mathscr P}_{{\rm PGL_2}, \bS}$.

\vskip 1mm \item  Rational coordinate systems assigned to different ideal triangulations 
of $\bS$ provide a cluster Poisson variety structure  
on the moduli space  ${\mathscr P}_{{\rm PGL_2}, \bS}$. 

\vskip 1mm\item Forgetting   pinnings, we recover the cluster Poisson structure of the 
  space  ${\mathscr X}_{{\rm PGL_2}, \bS}$ defined in \cite{FG03a}.
\end{enumerate} 
\et

\begin{proof} Claim 1) follows 
immediately from Lemma \ref{12.21.08.3}. Claims 2), 3)  are general properties
 of cluster amalgamation.  Claim 4) follows from the very definitions. \end{proof}
 
  For example, for the torus with a hole and two marked points on the boundary, 
see Figure \ref{fga9}, 
a coordinate system is obtained by amalgamating $8$ triangles with pinnings.  

\begin{figure}[ht]
\epsfxsize 200pt
\center{
\begin{tikzpicture}[scale=1.0]
\draw (0, 2) -- (0,0) -- (2,0) -- (2,2)--(0,2) -- (1, 1.7);    
\draw (1,0.3) -- (2,2) -- (1, 1.7) -- (0,0) -- (1,0.3) -- (2,0); 
\draw [dashed, red, arrows={-Stealth[length=1.5mm 3 0]}] (1,0.3) to [out=110, in=250] (1,1.7);
\draw [dashed, red, arrows={-Stealth[length=1.5mm 3 0]}] (1,1.7) to [out=-70, in=70] (1,0.3);
\begin{scope}[shift={(6,1)}, scale=0.4]
\draw [thick]  (-2.6,0) circle (1.5mm);
\draw (-3.5,0) .. controls (-3.5,2) and (0, 2.5) .. (3, 2);
\draw[yscale=-1] (-3.5,0) .. controls (-3.5,2) and (0, 2.5) .. (3, 2);
\draw (-2,0.2) .. controls (-1.5,-0.3) and (-1, -0.5) .. (0, -0.5) .. controls (1,-0.5) and (1.5, -0.3) .. (2,0.2);
\draw (-1.75,0) .. controls (-1.5,0.3) and (-1, 0.5) .. (0, 0.5) .. controls (1,0.5) and (1.5, 0.3) .. (1.75,0);
\draw [dashed, red, arrows={-Stealth[length=1.5mm 3 0]}] (3,-2) to [out=120, in=240] (3,2);
\draw [dashed, red, arrows={-Stealth[length=1.5mm 3 0]}] (3,2) to [out=-60, in=60] (3,-2);
\end{scope}
\end{tikzpicture}
 }
\caption{An ideal triangulation, with pinnings,  
of a punctured torus with   two special   points.}
\label{fga9}
\end{figure}

\medskip 

\subsection{Cluster Poisson amalgamation for ${\mathscr P}_{\rm PGL_m, \bS}$} \la{SSECC3.2}

\medskip

We define a cluster Poisson structure on the space ${\mathscr P}_{\rm PGL_m, \bS}$ by amalgamating the spaces ${\mathscr P}_{{\rm PGL_m}, t}$ assigned to the triangles of an ideal triangulation $T$ of $\bS$.

\subsubsection{Pinnings for ${\rm PGL_m}$.} 
A generic pair of flags $(F_1, F_2)$  is equivalent to a decomposition 
of the vector space $V_m$  into a sum of 
$1-$dimensional subspaces: 
$V_m = E_1 \oplus \ldots \oplus E_m$ 
such that 
\be \nonumber
\begin{split}
&F_1 = \{E_1, \quad  E_1\oplus E_2, \quad  \ldots, \quad E_1 \oplus E_2
\oplus \ldots  \oplus E_m  
 \}, \\
&F_2 = \{E_m, \quad  E_m\oplus E_{m-1}, \quad  \ldots, \quad E_m
\oplus E_{m-1} \oplus ... \oplus E_1  
 \}.  \\
\end{split}
\ee
A pinning   is described by a triple 
$(F_1, F_{2}; p)$
where $F_1$ and $F_2$ are generic flags in ${\mathbb P}(V_m)$, 
and $p$ is a $1$-dimensional subspace in   generic position to the flags. 
Indeed, the subspace $p$ determines 
a basis $(e_1, ..., e_m)$ of $V_m$ 
defined up to a common non-zero scalar such that each $e_i\in E_i$ and
$e_1+ ... + e_m \in p$.  We   denote the pinning by $p: F_1\stackrel{}{\to}F_2$. This definition depends on a choice of orientation. Changing the orientation, we get its opposite pinning
$p^\ast: F_2\stackrel{}{\to}F_1$, where $p^\ast$ is the $1$-dimensional subspace containing $\sum_{i=1}^m (-1)^i e_i$.

\subsubsection{Cluster Poisson structure of the   space ${\mathscr P}_{{\rm PGL_m}, t}$.} \label{SSECC3.2.2}
A generic point of the   space ${\mathscr P}_{{\rm PGL_m}, t}$ can be described 
by a configuration   
\be \la{1.4.09.1}
(F_1, F_2, F_3, p, q, r) \quad \mbox{modulo the action of ${\rm PGL_m}$}. 
\ee
Here $F_i$ are flags in $ {\mathbb P}(V_m)$, and 
$(F_1, F_{2}; p)$ etc. are pinnings, see the left of Figure \ref{fga12}.  

An essential part of the cluster Poisson structure 
was defined in \cite{FG03a}. We complement it by 
adding pinnings to the sides of the triangle, and describing the coresponding 
frozen coordinates. 

\begin{figure}[ht]
\epsfxsize 200pt
\center{
\begin{tikzpicture}[scale=0.6]
\foreach \count in {1,2,3}
{
\node at (90-120*\count:2.3) {${\rm F}_\count$};
\draw [thick, arrows={-Stealth[length=2mm 3 0]}]  (90-120*\count:2cm) -- (-30-120*\count:2cm);
}
\node at (30:1.4) {$r$};
\node at (150:1.4) {$q$};
\node at (270:1.4) {$p$};
\begin{scope}[shift={(6,-1)}]
\foreach \n in {1,2,3}
{\node[circle,draw,minimum size=5pt,inner sep=0pt] (\n) at (60:\n) {};
\node[circle,draw,minimum size=5pt,inner sep=0pt] (b\n) at (0:\n) {};
\node[circle,draw,minimum size=5pt,inner sep=0pt] (c\n) at ([shift={(4,0)}]120:\n) {};
}
\node[circle,fill,draw,minimum size=5pt,inner sep=0pt] (A) at ([shift={(1,0)}]60:1) {};
\node[circle,fill,draw,minimum size=5pt,inner sep=0pt] (B) at ([shift={(1,0)}]60:2) {};
\node[circle,fill,draw,minimum size=5pt,inner sep=0pt] (C) at ([shift={(2,0)}]60:1) {};
 \draw[directed,dashed] (1) -- (2); 
  \draw[directed,dashed] (2) -- (3); 
   \draw[directed,dashed] (b3) -- (b2); 
  \draw[directed,dashed] (b2) -- (b1); 
 \draw[directed,dashed] (c3) -- (c2); 
  \draw[directed,dashed] (c2) -- (c1); 
\foreach \from/\to in {c1/C,C/A, A/1, c2/B,B/2,c3/3,1/b1, 2/A, A/b2, 3/B,B/C,C/b3, b3/c1, C/c2, b2/C, b1/A, A/B, B/c3}
                 \draw[directed, thick] (\from) -- (\to);    
\end{scope}
\end{tikzpicture}
 }
\caption{Triples of flags with pinnings, and the 
Poisson tensor $\varepsilon$ ($\G={\rm PGL_4}$).}
\label{fga12}
\end{figure}

Let ${\rm I}^{(m)}$ 
be the set of all non-negative integral solutions of the equation
$$
a_1+a_2+a_3 = m, \qquad a_i \in \Z_{\geq 0},
$$ 
excluding the solutions with two $0$'s. 
It is parametrised by vertices of the $m$-triangulation of a triangle 
shown  on Figure \ref{fga12}. 
The solutions including $0$ form the frozen  
subset ${\rm F}^{(m)}$ of ${\rm I}^{(m)}$. 
The corresponding vertices are shown by little circles. 
The non-zero solutions are parametrised by the 
boldface points. 
Let us assign to each such solution $(a_1, a_2, a_3)$ a coordinate 
$X_{a_1, a_2, a_3}$. 

Recall the triple ratio of three flags in ${\mathbb P}^2 = {\mathbb P}(V_3)$. 
Given vectors $a,b,c\in V_3$, let 
$\Delta(a, b, c)$ be value of a  volume in $V_3$ on 
$a\wedge b\wedge c$. 
We describe a flag in ${\mathbb P}^2$ by a pair of vectors $(a,b)$ in $V_3$, so that the flag is $\langle a \rangle \subset \langle a, b\rangle$. 
The triple ratio of the flags $A, B, C$ described by 
three pairs of vectors $(a_1, a_2), (b_1, b_2), (c_1, c_2)$ in $V_3$ 
is given by
$$
r_3(A, B, C):= \frac{\Delta(a_1, a_2, b_1)\Delta(b_1, b_2, c_1)\Delta(c_1, c_2, a_1)}
{\Delta(a_1, a_2, c_1)\Delta(b_1, b_2, a_1)\Delta(c_1, c_2, b_1)}.
$$ 

We start from a non-zero solution. It provides us with  a subspace 
\be \la{1.09.3.1}
V_{a_1, a_2, a_3}:= F_1^{(a_1-1)} \oplus F_2^{(a_2-1)} \oplus F_3^{(a_3-1)}.  
\ee
There is a $3-$dimensional quotient
 $
V_m/V_{a_1, a_2, a_3}. 
 $
The flag 
$F_i$ gives rise to a flag $\overline F_i$ in the quotient\footnote{We 
skip the triple $(a_1, a_2, a_3)$ from the notation of this flag.}, 
obtained by projecting the flag
$$
F_i^{(a_1)}/F_i^{(a_1-1)}  \subset
F_i^{(a_1+1)}/F_i^{(a_1-1)}.   
$$
The value of the coordinate $X_{a_1, a_2, a_3}$ on 
configuration (\ref{1.4.09.1}) is the triple ratio 
of flags $\overline F_i$:
$$
X_{a_1, a_2, a_3}:= r_3(\overline F_1, \overline F_2,\overline F_3). 
$$
Take a solution $(a_1, 0, a_3)$. 
It corresponds to a point on the side of the triangle 
with the flags $F_3, F_1$. 
Since  $F_2^{(-1)}=0$,  subspace (\ref{1.09.3.1})  
reduces to a subspace
$
V_{a_1, 0, a_3}:= F_1^{(a_1-1)} \oplus F_3^{(a_3-1)}. 
$
So we get a $2-$dimensional quotient 
$V_m/V_{a_1, 0, a_3}$.  
The flags $F_1, F_2, F_3$ provide three $1-$dimensional subspaces 
in this quotient, given by 
$$
\overline F_1:= F_1^{(a_1)}, \quad \overline F_2:= F_2^{(1)}, \quad 
  \overline F_3:= F_3^{(a_3)}.
$$
There is the fourth $1-$dimensional subspace $\overline R$ given by the projection 
of the pinning $1-$dimensional subspace $r$ to  the quotient. 
The value of the coordinate $X_{a_1, 0, a_3}$ on a 
configuration (\ref{1.4.09.1}) is the cross ratio 
of these $1-$dimensional subspaces: 
\be \la{1.4.09.2}
X_{a_1, 0, a_3}:= r(\overline F_1, \overline F_2, \overline F_3, \overline R). 
\ee
The other frozen coordinates are defined similarly. 

We define the Poisson tensor $\varepsilon$   by the rule shown on Figure \ref{fga12}. 
Here two vertices $ij$ are connected by a solid arrow going from $i$ to $j$ 
if and only if 
$\varepsilon_{ij}=1$, and by dotted arrow   from $i$ to $j$ if and only if 
$\varepsilon_{ij}=1/2$. The function $\varepsilon_{ij}$ is skewsymmetric, 
and takes values $0, \pm 1/2, \pm 1$. So it is determined by the picture. 
The multipliers are equal to $1$. This defines a quiver.

\subsubsection{Cluster structure of the moduli space ${\mathscr P}_{{\rm PGL_m}, \bS}$.}  
 \label{SSECC3.2.3}

Cluster amalgamation of these quivers and the corresponding cluster Poisson varieties 
over triangles $t$ of an ideal triangulation of $\bS$ 
provides  a cluster Poisson structure on the moduli space 
${\mathscr P}_{{\rm PGL_m}, \bS}$. Note that after the amalgamation 
the value of the $\varepsilon$-function between the vertices at  
the same edge of the ideal triangulation is zero. 

Forgetting the   pinnings  and the 
frozen coordinates, we get the    quiver
describing the cluster Poisson structure of the moduli space 
${\mathscr X}_{{\rm PGL_m}, \bS}$ in \cite{FG03a}.  
By the very definition, the cluster Poisson coordinates 
assigned to the internal vertices of the triangles are 
the   functions $X_{a_1, a_2, a_3}$ from 
 \cite{FG03a}. 

Finally, 
recall that the coordinate assigned to a quadruple of flags $(F_1, F_2, F_3, F_4)$ 
and a pair of integers $(a_1, a_3)$ such that  $a_1+a_3=m$ is the cross ratio
$$
X_{a_1, a_3}:= 
r(\overline F^{(a_1)}_1, 
\overline F^{(1)}_2, \overline F^{(a_3)}_3, \overline F^{(1)}_4).
$$
Here $\overline F^{(a)}_k$ is the one dimensional subspace in 
the two dimensional quotient $V_2/(F^{(a_1-1)}_1\oplus F^{(a_3-1)}_3)$ 
obtained by projection of the subspace $F^{(a)}_k$. It is easy to show that  
\[r(\overline{F}^{(a_1)}_1, \overline{R^\ast},{ \overline F^{(a_3)}_3}, {\overline R})=1.\] 
Therefore Lemma \ref{12.21.08.3} implies that  the   cluster coordinates 
assigned after the amalgamation 
to the vertices at the edges of the ideal triangulation 
are the coordinates from \cite{FG03a}.  
This boils down to 
$X_{a_1, a_3} = X_{a_1, 0, a_3}X_{a_1, a_3, 0}$ since these 
coordinates are products of the cluster coordinates 
(\ref{1.4.09.2}). 
Summarising, we get

\bt \la{12.21.8.102}
Given an ideal triangulation $T$ of  $\bS$, the 
 amalgamation of cluster varieties  ${\mathscr P}_{{\rm PGL_m}, t}$ over the triangles $t$ 
of  $T$ 
provides a cluster Poisson  
structure on the moduli space  ${\mathscr P}_{{\rm PGL_m}, \bS}$. 
Forgetting  the  pinnings, we recover the
 cluster Poisson structure on ${\mathscr X}_{{\rm PGL_m}, \bS}$  defined 
in {\cite{FG03a}. }
\et


\medskip

\section{Cluster  coordinates on the moduli spaces ${\mathscr A}_{\G, \bS}$ and ${\mathscr P}_{\G', \bS}$} \la{Clcoord}
 
\medskip

In Section \ref{Clcoord}, it is convenient to assume first that all boundary intervals of $\bS$ are colored. 
The  passage to the colored decorated surfaces $\bS$, where only some of the boundary intervals are colored, is straightforward: one simply removes the frozen vertices and the frozen coordinates corresponding to the non-colored sides.

\subsection{Decompositions for a pair of decorated flags} \la{SSec2.3X}
\medskip

 Let $w\in W$. 
Consider the subset of the set   of vertices of the Dynkin diagram given by 
\be \la{simple.index.positive.r.hh}
{\rm I}(w):=\left\{j\in {\rm I} ~ \middle |~ w\cdot \alpha_j^\vee  \mbox{ is negative}\right\}. 
\ee
It gives rise to a subgroup  of the Cartan subgroup ${\rm H}$ denoted by 
\be
\la{H.5.26.16.hh}
{\rm H}(w)= \{\prod_{i\in {\rm I}(w)} \alpha_i^\vee(b_i)~\small | ~ b_i \in \mathbb{G}_m\}.
\ee
Note that the following statements are equivalent: 
\[ {\rm I}(w)={\rm I} \qquad \Longleftrightarrow\qquad  {\rm H}(w)={\rm H} \qquad \Longleftrightarrow\qquad  w=w_0. \]

Let  ${\bf i}=(i_1,\ldots, i_n)$ be a reduced word of $w$. We get a sequence of distinct coroots
\be
\la{reduced.word.positive.roots}
\beta_{k}^{\bf i}:= s_{i_n}\ldots s_{i_{k+1}}\cdot  \alpha_{i_k}^\vee, \hskip 7mm k \in \{1,\ldots, n\}.
\ee
They are precisely the positive coroots $\alpha^\vee$ such that $w\cdot \alpha^\vee$ are negative. 
Take their sum 
\be  \la{rhow}
\rho_w:= \sum_k \beta_{k}^{\bf i}.
\ee
Note that $w\cdot \alpha^\vee$ is negative if and only if $w_0w\cdot \alpha^\vee$ is positive. Therefore the subset  ${\rm I}(w_0w)$ of ${\rm I}$ is the complement of ${\rm I}(w)$. So one has 
$
\rho_w + \rho_{w_0w} =\rho_{w_0}.
$

\bl We have
\be
\rho_w(-1) = \overline{w^{-1}}\cdot \overline{w}.
\ee
\el
\begin{remark} For the longest element $w_0$, we have $\rho_{w_0}(-1)= s_{\G}$.
\end{remark}
\begin{proof} 
By definition $\overline{w^{-1}} =\overline{s}_{i_n}\ldots  \overline{s}_{i_1}$ and ${\overline{w}}= {\overline{s}}_{i_1}\ldots  {\overline{s}}_{i_n}$. 
Note that  $ {\overline{s}}_{i}^2=\alpha_{i}^\vee(-1)$. Therefore
\[
\overline{w^{-1}}\cdot \overline{w} = \overline{s}_{i_n}\ldots {\overline{s}}_{i_2} \alpha_{i_1}^\vee(-1) {\overline{s}}_{i_2}\ldots  {\overline{s}}_{i_n}  = \beta_{1}^{\bf i}(-1) {\overline{s}}_{i_n}\ldots {\overline{s}}_{i_2} {\overline{s}}_{i_2}\ldots  {\overline{s}}_{i_n}. 
\]
Repeating the same process for $i_k$, we prove the Lemma.
\end{proof}

\blc
\la{decomp.partial.dec.flag.15.02}
Let  $(\A_l, \A_r)$ be a pair of decorated flags with
\be
w(\A_l, \A_r)=w, \hskip 7mm h(\A_r, \A_l)\in {\rm H}(w).
\ee
There exists a unique chain of decorated flags
\begin{center}
\begin{tikzpicture}
\node [circle,draw=red,fill=red,minimum size=3pt,inner sep=0pt,label=below:{\small {$\A_l={\A}_0\qquad ~$}}] (a) at (0,0) {};
 \node [circle,draw=red,fill=white,minimum size=3pt,inner sep=0pt,label=below:{\small $~{\A}_1$}]  (b) at (1,0){};
 \node [circle,draw=red,fill=white,minimum size=3pt,inner sep=0pt,label=below:{\small $~{\A}_2$}] (c) at (2,0) {};
 \node [circle,draw=red,fill=white,minimum size=3pt,inner sep=0pt,label=below:{\small $~{\A}_3$}] (d) at (3,0) {};
  \node [circle,draw=red,fill=white,minimum size=3pt,inner sep=0pt,label=below:{\small $~{\A}_{4}$}] (e) at (4,0) {};
  \node [circle,draw=red,fill=red,minimum size=3pt,inner sep=0pt,label=below:{\small {$\qquad \qquad ~{\A}_{n}=\A_r$}}] (f) at (5.5,0) {};
  \node[blue] at (0.5, 0.2) {{\small $s_{i_1}$}};
   \node[blue] at (1.5, 0.2) {{\small $s_{i_2}$}};
    \node[blue] at (2.5, 0.2) {{\small $s_{i_3}$}};
     \node[blue] at (3.5, 0.2) {{\small $s_{i_4}$}};
     \node[blue] at (4.75, 0.2) {$\cdots$};
      \node at (4.75, -0.45) {$\cdots$};
 \foreach \from/\to in {a/b, b/c, c/d, e/f, d/e}
                   \draw[thick] (\from) -- (\to); 
 \end{tikzpicture}
\end{center}
such that
\begin{itemize}
\vskip 1mm\item  $w(\A_{k}, \A_{k-1})= s_{i_k}$ for each $k\in \{1,\ldots, n\}$.
\vskip 1mm\item If  $\beta_{{\bf i},k}$ is simple, then $h(\A_{k}, \A_{k-1}) \in \alpha_{i_k}^\vee({\Bbb G}_m)$. Otherwise $h(\A_{k}, \A_{k-1})=1$.
\end{itemize}
\elc

\begin{proof} 
By Corollary \ref{USFX}, there is a unique chain of flags
\be
\la{decompos.flag.chain.0}
\{\B_l=\B_0, \, \B_1,\, \B_2, \ldots,\, \B_n=\B_r\} 
\ee
such that $w(\B_{k}, \B_{k-1})=s_{i_k}$ for each  $k \in \{1,\ldots, n\}$. 
By definition,
\[
h(\A_r, \A_l)=\prod_{i \in {\rm I}(w)} \alpha^\vee_i (b_i) \in {\rm H}(w).
\]
Let $\A_n := \A_r$. The decorations for the other $\B_k$'s are recursively defined such that  
\be
\la{special.hk}
h(\A_{k}, \A_{k-1})=    \left\{  \begin{array}{ll}  
      \alpha_{i_k}^\vee (b_i)& \mbox{if  $ \beta_{k}^{\bf i}=\alpha_i^\vee$ is simple},  \\
     1 & \mbox{otherwise}. \\
   \end{array}\right.
\ee
Therefore
\be
\la{special.hk1}
h_k := s_{i_n}\ldots s_{i_{k+1}}\big( h(\A_{k}, \A_{k-1})\big) = \left\{  \begin{array}{ll} 
      \alpha_i^\vee(b_i) & \mbox{if  $\beta_{k}^{\bf i}=\alpha_i^\vee$ is simple},  \\
      1 & \mbox{otherwise}. \\
   \end{array}\right.
\ee 
 By \eqref{chain.formula.for.cartan.5.8}, we have
$
h(\A_n, \A_0) =h_n\ldots h_1= h(\A_r, \A_l). 
$ 
Therefore $\A_0= \A_l$. The uniqueness of the chain follows by construction.
\end{proof}

\bl \la{basic.decomp.chain.w0}
Suppose that $(\A_l, \A_r)$ is  generic. Let ${\bf i}=(i_1,\ldots, i_m)$ be a reduced word of $w_0$. Set
\[v_k:= s_{i_k}\cdots s_{i_1}, \hskip14mm u_k:= s_{i_{k+1}}\cdots s_{i_{n}}.
\]
The decoration of the intermediate $\A_k$ in the chain   is chosen such that 
\[
h(\A_l, \A_k)\in {\rm H}(v_k)\cdot \rho_{v_k}(-1), \hskip 7mm
h(\A_r, \A_k) \in {\rm H}(u_k).
\]
\el
\begin{proof}  
The second condition follows  from  construction \eqref{special.hk}. Consider the opposite reduced word ${\bf i}^{\circ}:= (i_m, \ldots, i_1)$. Note that $\beta_{k}^{{\bf i}^{\circ}}$ is simple if and only if $\beta_{m+1-k}^{\bf i}$ is simple.
The reversed chain  satisfies a similar property with respect to ${\bf i}^{\circ}$. If $h(\A_k, \A_{k-1})=\alpha_{i_k}^\vee(b_k)$
then $
h(\A_{k-1}, \A_k)= \alpha_{i_k}^\vee(-b_k) 
$.  So we get the first condition with the extra factor $\rho_{v_k}(-1)$  on the right.
\end{proof}

\medskip

\subsection{Cluster $K_2$-coordinates  on the moduli space ${\mathscr A}_{\G, \bS}$} \la{SeCC2.4A}

\medskip

In Section \ref{SeCC2.4A}, $\G$ is a simply connected group and $\G'$ is the adjoint group.

Let  ${\bf i}=(i_1,\ldots, i_m)$ be a reduced word of $w_0$ and let ${\bf i}^\ast=(i_1^\ast, \ldots, i_m^\ast)$. By Lemma \ref{decomp.partial.dec.flag.15.02}, for a generic pair $(\A_2, \A_3)$ there is a unique chain  $\{\A_2 =\A_2^0,\,\A_2^1,\, \ldots, \A_2^m =\A_{3}\}$  with respect to  ${\bf i}^\ast$.

For every fundamental weight $\Lambda_i$,  there is a unique $\G$-invariant regular function $\Delta_i$ on $\mathcal{A}^2$  such that $\Delta_i(h\cdot [\U], [\U^-])=\Lambda_i(h)$ for all $h\in {\rm H}$.

\bd \la{DEF9.1} Pick  a reduced word ${\bf i}=(i_1,\ldots, i_m)$ of $w_0$, and a vertex $\A_1$ of the triangle. 
The associated collection of $K_2$-cluster $\{A_k\}$ on  ${\rm Conf}_3^{\times}({\cal A})$ consists of 

\begin{itemize}

\item the frozen cluster coordinates  given by  $\Delta_i(\A_1, \A_{2}), ~\Delta_i({\A_1, \A_3}),~\Delta_i(\A_3, \A_2)$ for $i\in {\rm I}$; and

\item the non-frozen cluster coordinates given by
\[
A_p=\Delta_{i_p}(\A_1, \A_{2}^p) ,
\]
where $p$ runs through the indices $1, ..., m$ such that $i_p$ is not the rightmost simple reflection $i$ appearing in ${\bf i}$ for all $i\in {\rm I}$, and $\A_2^p$ is a decorated flag in Figure \ref{2017.5.9.6.24h*}.
\end{itemize}
\ed

\begin{figure}[ht]
\epsfxsize 200pt
\center{
\begin{tikzpicture}[scale=1.1]
\node [circle,draw=red,fill=red,minimum size=3pt,inner sep=0pt, label=above:{\small ${\A}_1$}] (a) at (0,0) {};
\node [circle,draw=red,fill=red,minimum size=3pt,inner sep=0pt, label=below:{\small ${{\A}_2=\A_2^0\qquad \qquad}$}] (b) at (-1.8,-3) {};
 \node [circle,draw=red,fill=white, minimum size=3pt,inner sep=0pt, label=below:{\small $~{\A}_2^2$}]  (c) at (-0.36,-3){};
  \node [circle,draw=red,fill=white, minimum size=3pt,inner sep=0pt, label=below:{\small $~{\A}_2^1$}]  (g) at (-1.08,-3){};
 \node [circle,draw=red,fill=white, minimum size=3pt,inner sep=0pt, label=below:{\small $~{\A}_2^{3}$}] (h) at (0.36,-3) {};
 \node [circle,draw=red,fill=red,minimum size=3pt,inner sep=0pt, label=below:{\small ${\qquad \qquad \qquad{\A}_2^m=\A_3}$}] (e) at (1.8,-3) {};
 \node[blue] at (1, -2.75) {$\cdots$};
  \node at (1.1, -3.4) {$\cdots$};
 \foreach \from/\to in {a/b, b/g, g/c, h/e, a/e, c/h}
                  \draw[thick] (\from) -- (\to);
\foreach \from/\to in {a/c, a/g, a/h}                   
              \draw[thick] (\from) -- (\to);                      
 \node[blue] at (-1.44,-2.8) {{\small $s_{i_1^\ast}$}};
  \node[blue] at (-0.72,-2.8) {{\small ${s}_{i_2^\ast}$}}; 
   \node[blue] at (0,-2.8) {{\small $s_{i_3^\ast}$}};                
 \end{tikzpicture}
 }
 \caption{Introducing a cluster structure on ${\rm Conf}_3^{\times}({\cal A})$.}
 \label{2017.5.9.6.24h*}
 \end{figure}

In Section \ref{SECT4} we give another definition of functions $\{A_k\}$   on  ${\rm Conf}_3({\cal A})$, clarifying their nature. 
Note that the group $\H^3$ acts on ${\rm Conf}_3({\cal A})$ by rescaling decorated flags. There is a  decomposition 
$$
{\cal O}({\rm Conf}_3({\cal A})) = \bigoplus_{\lambda, \mu, \nu}(V_\lambda \otimes V_\mu \otimes V_\nu)^\G.
$$
Here $\lambda, \mu, \nu$ are dominant weights,   $V_\lambda$ is the irreducible finite dimensional   $\G$-module with the highest weight $\lambda$, and 
$(V_\lambda \otimes V_\mu \otimes V_\nu)^\G$  the eigenspace of $\H^3$ with the eigencharacter $(\lambda, \mu, \nu)$.
 
We interpret    functions $A_k$  as specific vectors in certain triple tensor products invariants with  
  the {\it multiplicity one  property}: 
$$
A_k \in (V_\lambda \otimes V_\mu \otimes V_\nu)^\G, \qquad{\rm dim}(V_\lambda \otimes V_\mu \otimes V_\nu)^\G=1.
$$
 The collection of triples $(\lambda, \mu, \nu)$ which appear in a given cluster coordinate system 
 is determined by the reduced decomposition ${\bf i}$. Proposition \ref{equivalence.cluster.a.ba} tells that both definitions are equivalent. 

\vskip 2mm

Next, take an ideal triangulation ${\cal T}$ of a decorated surface $\bS$, and cut the surface into triangles. For each of the triangles $t$, 
 pick a vertex $v_t$ and a reduced word ${\bf i}_t$ of $w_0$. 
Then Definition \ref{DEF9.1} provides a   collection of functions  $\{A_{t, k}\}$ on the moduli space ${\mathscr A}_{ \G, t}$ assigned to each   triangle $t$. 

There is a    {\it restriction map}, given by the restriction of a decorated $\G$-local system on $\bS$ to  $t$:
$$
{\rm Res}_{t}: {\mathscr A}_{\G, \bS }   \lra {\mathscr A}_{\G, t}.
 $$

\bd \la{DEF9.2} Take an ideal triangulation ${\cal T}$ of a decorated surface $\bS$,  and pick for each of the triangles $t$ of   ${\cal T}$   a vertex $v_t$, 
and a reduced word ${\bf i}_t$ of $w_0$. We assign to any such a data $({\cal T}, \{{\bf i}_t\}, \{{v}_t\})$
a collection of  regular functions  on  the space $ {\mathscr A}_{\G, \bS}$ is given by the functions 
\be \la{RDT}
\{{\rm Res}^*_{t} A_{t, k}\}.
\ee
Here, for a given triangle $t$, the functions $\{A_{t, k}\}$ are the ones $\{A_k\}$ from Definition \ref{DEF9.1}. 

  The frozen variables are  the functions assigned to the boundary intervals of $\bS$. 
 
\ed
Note that the functions  assigned to any edge $E$ of the triangulation ${\cal T}$ 
depend only on the two decorated flags at the vertices of the edge. They do not depend on the rest of the data. 

In Section \ref{sec.3.1} we assign to   a triangle $t$ with a vertex $v$, and a reduced decomposition ${\bf i}$ of $w_0$ 
 a  quiver ${\bf Q}_{t, v, {\bf i}}$. 
Its frozen vertices are assigned to the sides of the triangle $t$.  The frozen variables at each side are parametrized by the set ${\rm I}$ of  vertices of the Dynkin diagram. 
So, given an ideal triangulation ${\cal T}$ of $\bS$, one can amalgamate  quivers ${\bf Q}_{t, v, {\bf i}}$ assigned to the triangles $t$ of the 
triangulation ${\cal T}$, getting a quiver 
\be \la{RDT1}
{\bf Q}_{ {\cal T}, \{{\bf i}_t\}, \{{v}_t\}} := \ast_{t} {\bf Q}_{t, v, {\bf i}}.
\ee

\bt \la{MTHAs} There is a unique cluster $K_2$-structure on the   space ${\mathscr A}_{ \G, \bS}$   such that 
for each data $({\cal T}, \{{\bf i}_t\}, \{{v}_t\})$,   functions   (\ref{RDT}) form a cluster 
$K_2$-coordinate system associated with   quiver (\ref{RDT1}).
\et

 The proof of Theorem \ref{MTHAs} is given in Sections \ref{pcs.sec}-\ref{sec.conf3}. It consists of proving  two major claims:
  \begin{itemize}
  
\vskip 1mm  \item  
  
  {\it Cyclic invariance}. We  prove in Theorem \ref{main.result2.cluster.cyclic} 
  that changing the  vertex $v_t$ of the triangle $t$, or changing the reduced decomposition ${\bf i}_t$ of $w_0$,  amounts to a cluster   $K_2$-transformation.  
  
\vskip 1mm  \item {\it Indepence on flips.} Given a quadrangle $Q$, let us cut it by a diagonal  into    triangles $t_1$ and $t_2$. We prove in Section \ref{Sect8.5} that 
  the   cluster $K_2$-structure on the  space  ${\mathscr A}_{\G, Q}$  
 given by the amalgamation of   ${\mathscr A}_{\G, t_1}$ and  ${\mathscr A}_{\G, t_2}$ does not depend on the choice of the diagonal.
   \end{itemize}

\begin{figure}[ht]
\epsfxsize 200pt
\center{
\begin{tikzpicture}[scale=0.7]
\node [circle,draw=red,fill=red,minimum size=3pt,inner sep=0pt] (a) at (90:2) {};
\node [circle,draw=red,fill=red,minimum size=3pt,inner sep=0pt] (b) at (210:2) {};
 \node [circle,draw=red,fill=red,minimum size=3pt,inner sep=0pt] (c) at (330:2) {};
 \node [] at (90:2.4) {{\small $\B_1$}};
\node [] at (210:2.4) {{\small $\B_2$}};
 \node [] at (330:2.4) {{\small $\B_3$}};
  \node [blue] (a1) at (105:2.3) {{\small $\A_1$}};
   \node [blue] (a2) at (75:2.3) {{\small $\A_1'$}};
     \node [blue] (b1) at (225:2.3) {{\small $\A_2$}};
   \node [blue] (b2) at (195:2.3) {{\small $\A_2'$}};
     \node [blue] (c1) at (345:2.3) {{\small $\A_3$}};
   \node [blue] (c2) at (315:2.3) {{\small $\A_3'$}};
 \foreach \from/\to in {a/b, b/c, c/a}
                  \draw[thick] (\from) -- (\to);
  \foreach \from/\to in {a1/b2, b1/c2, c1/a2}
                  \draw[blue, <->] (\from) -- (\to);                 
 \end{tikzpicture}
 }
 \caption{The space $\mathscr{P}_{\G, t}$; here $p_{12}=(\A_1, \A_2')$ is a pinning over  $(\B_1, \B_2)$ and so on.}
 \label{2018.9.23.15.10ss}
 \end{figure}

\subsection{Cluster Poisson coordinates on the moduli space  $\mathscr{P}_{\G, \bS}$} \la{CSEC4.3}

\medskip

In Section \ref{CSEC4.3}, about the space $\mathscr{P}_{\G, \bS}$,  we assume that $\G$ is the adjoint group. 

\subsubsection{Cluster Poisson coordinates on $\mathscr{P}_{\G, t}$.}
  Let $(\B_1, \B_2)$ be a generic pair of flags of $\G$. Recall that a pinning for  $(\B_1, \B_2)$ is a pair of decorated flags $(\A_1, \A_2) \in (\G/\U)^2$ over $(\B_1, \B_2)$ such that $h(\A_1, \A_2)=1$.

The space $\mathscr{P}_{\G, t}$ parametrizes $\G$-orbits of data
$(\B_1, \B_2, \B_3; p_{12}, p_{23}, p_{31})$ where $(\B_i, \B_{i+1})$ are flags of generic position and $p_{i, i+1}$ is a pinning for $(\B_i, \B_{i+1})$ for each $i\in \Z/3\Z$, see 
Figure \ref{2018.9.23.15.10ss}.

 Let ${\bf i}$ be a reduced word of $w_0$.     
Take the decomposition of $(\B_2, \B_3)$ with respect to ${\bf i}^\ast$:
\be
\la{2018.9.22.17.41ssX}
\xymatrix{\B_2=\B_2^0 &\ar[l]_{\small \qquad s_{i_1^\ast}} \B_2^1&\ar[l]_{\small s_{i_2^\ast}} \B_2^2&\ar[l]_{\small ~s_{i_3^\ast}} \cdots&\ar[l]_{\small ~~~s_{i_m^\ast}\qquad }\B_2^m=\B_3.}
\ee
Suppose that the pairs $(\B_1, \B_2^j)$ are  generic  for $j=0,\ldots, m$. 

\begin{example}
A flag for ${\rm PGL}_3$ can be presented by a pair $(p, L)$ in the projective space $\mathbb{P}^2$, where $p$ is a point on a line $L$. As illustrated by the following figure, let $(p_1, L_1)$ and $(p_2, L_2)$ be a generic pair of flags. 
\begin{center}
\begin{tikzpicture}
\draw[thick] (-1,-0.5) -- (2,1);
\draw[thick] (3, -.5) --(0,1);
\node[label=above: {\small $p_1\in L_1$}] at (-.4,-.2) {$\bullet$};
\node[label=above: {\small $p_2\in L_2$}] at (2.4,-.2) {$\bullet$};
\draw[thick, dashed] (-1, -.2) -- (3, -.2);
\end{tikzpicture}
\end{center}
 Let $p$ be the intersection of $L_1$ and $L_2$. With respect to the word $(1,2,1)$, we have the decomposition:
\[
\xymatrix{(p_1, L_1) &\ar[l]_{\small s_{1}} (p, L_1) & \ar[l]_{\small  s_{2}} (p, L_2) & \ar[l]_{\small  s_{1}} (p_2, L_2).}
\]
Let $L$ be the line connecting $p_1$ and $p_2$. With respect to the word $(2,1,2)$, we have the decomposition:
\[
\xymatrix{(p_1, L_1) &\ar[l]_{\small  s_{2}} (p_1, L) & \ar[l]_{\small  s_{1}} (p_2, L) & \ar[l]_{\small  s_{2}} (p_2, L_2).}
\]
\end{example}

Let $i\in {\rm I}$ be a vertex of the Dynkin diagram. Let $\B_{2, i}^j$ be the unique flag such that 
\[
w(\B_{2}^j, \B_{2, i}^j)=w_0s_i, \hskip 10mm w(\B_{2, i}^j, \B_1)=s_i.
\]
We prove in Lemma \ref{6.30.17.33.hh} that  $\B_{2,i}^{j-1}=\B_{2,i}^{j}$ whenever $i_j\neq i$. 
Let us take all the distinct flags among $\B_{2,i}^j$ and relabel them by $\B_{i\choose k},~ 0\leq k\leq n_i$ from the left to the right as shown below:

 
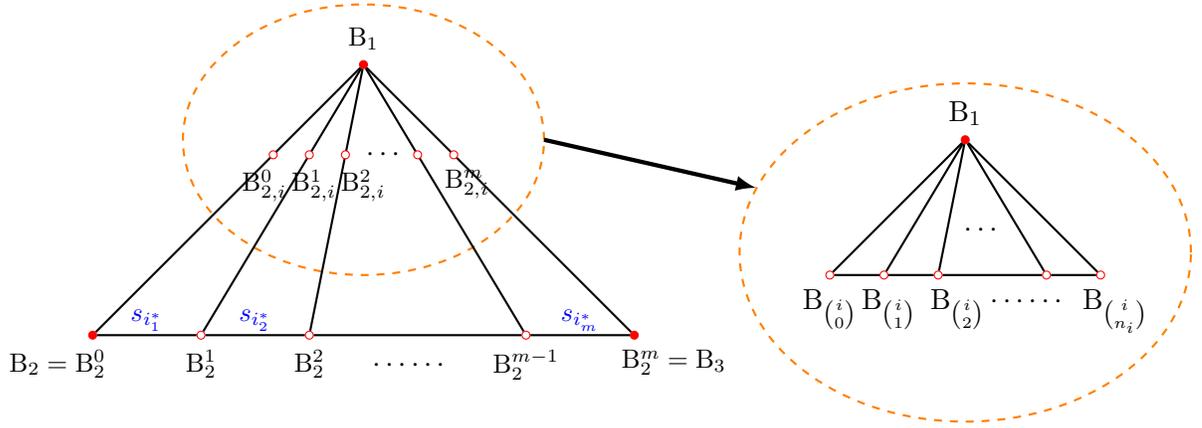
\begin{center}
\begin{tikzpicture}[scale=2]
   \draw[dashed, orange, thick] (0, -0.5) ellipse (1.2cm and 0.9cm);   
\node [circle,draw=red,fill=red,minimum size=3pt,inner sep=0pt, label=above:{\small ${\B}_1$}] (a) at (0,0) {};
\node [circle,draw=red,fill=red,minimum size=3pt,inner sep=0pt, label=below:{\small ${{\B}_2=\B_2^0\qquad \qquad }$}] (b) at (-1.8,-1.8) {};
 \node [circle,draw=red,fill=white, minimum size=3pt,inner sep=0pt, label=below:{\small ${\B}_2^2$}]  (c) at (-0.36,-1.8){};
  \node [circle,draw=red,fill=white, minimum size=3pt,inner sep=0pt, label=below:{\small ${\B}_2^1$}]  (g) at (-1.08,-1.8){};
 \node [circle,draw=red,fill=white, minimum size=3pt,inner sep=0pt, label=below:{\small ${\B}_2^{m-1}$}] (h) at (1.08,-1.8) {};
 \node [circle,draw=red,fill=red,minimum size=3pt,inner sep=0pt, label=below:{\small ${\qquad \qquad ~{\B}_2^m=\B_3}$}] (e) at (1.8,-1.8) {};
 \node at (.1, -.6) {$~\cdots$};
 \node at (.27, -2) {$~\cdots\cdots$};
 \foreach \from/\to in {a/b, b/g, g/c, h/e, a/e, c/h}
                  \draw[thick] (\from) -- (\to);
\foreach \from/\to in {a/c, a/g, a/h}                   
              \draw[thick] (\from) -- (\to);                      
 \node[blue] at (-1.44,-1.7) {{\small $s_{i_1^\ast}$}};
  \node[blue] at (-0.72,-1.7) {{\small ${s}_{i_2^\ast}$}}; 
   \node[blue] at (1.44,-1.7) {{\small $s_{i_m^\ast}$}};  
 \node [circle,draw=red,fill=white,minimum size=3pt,inner sep=0pt, label=below:{\small ${\B_{2,i}^0\qquad}$}]  at (-.6,-.6) {};
 \node [circle,draw=red,fill=white, minimum size=3pt,inner sep=0pt, label=below:{\small ${\qquad {\B}_{2,i}^2}$}]   at (-0.12,-.6){};
  \node [circle,draw=red,fill=white, minimum size=3pt,inner sep=0pt, label=below:{\small ${~{\B}_{2,i}^1}$}]   at (-.36,-.6){};
 \node [circle,draw=red,fill=white, minimum size=3pt,inner sep=0pt]  at (.36,-.6) {};
 \node [circle,draw=red,fill=white,minimum size=3pt,inner sep=0pt, label=below:{\small ${\qquad{\B}_{2,i}^m}$}]  at (.6,-.6) {}; 
 \draw [ultra thick, -latex] (1.2,-0.5) -- (2.61,-.82);     
\begin{scope}[shift={(4,-.5)}, scale=1.5]
 \node [circle,draw=red,fill=red,minimum size=3pt,inner sep=0pt, label=above:${\B}_1$] (z) at (0,0) {};
  \node [circle,draw=red,fill=white,minimum size=3pt,inner sep=0pt, label=below:${\B_{i \choose 0}}$] (y0) at (-.6,-.6) {};
 \node [circle,draw=red,fill=white, minimum size=3pt,inner sep=0pt, label=below:${\qquad {\B}_{i \choose 2}}$] (y2)  at (-0.12,-.6){};
  \node [circle,draw=red,fill=white, minimum size=3pt,inner sep=0pt, label=below:${~{\B}_{i \choose 1}}$]  (y1) at (-.36,-.6){};
 \node [circle,draw=red,fill=white, minimum size=3pt,inner sep=0pt]  (y) at (.36,-.6) {};
 \node [circle,draw=red,fill=white,minimum size=3pt,inner sep=0pt, label=below:${\qquad{\B}_{i \choose n_i}}$]  (yn) at (.6,-.6) {};    
  \foreach \from/\to in {z/y1, z/y, z/y2, z/y0, z/yn, y0/y1,y1/y2, y2/y, y/yn}
                  \draw[thick] (\from) -- (\to);  
 \node at (.08, -.4) {$\cdots$};
 \node at (.25, -.72) {$~\cdots\cdots$};               
   \draw[dashed, orange, thick] (0, -0.5) ellipse (1cm and 0.75cm);               
 \end{scope}
 \end{tikzpicture}
 \end{center}
 Recall the pair of decorated flags $\A_1, \A_1'$    provided by the pinnings sharing the vertex $\B_1$. Set:
 \be
 \la{defxikX}
 X_{i \choose k}= \left\{    \begin{array}{ll} 
      \displaystyle{  \mathcal{W}_i\left(\A_1, \B_{i \choose 0}, \B_{i \choose 1}\right)}& \hskip 5mm \mbox{if } k=0, \\
      & \\
      \displaystyle{r^+\left(\B_1, \B_{i \choose k-1}, \B_{i \choose k}, \B_{i \choose k+1}\right)} &  \hskip 5mm \mbox{if } 0<k<n_i, \\
     &  \\
       \displaystyle{\mathcal{W}_i\left(\A_1', \B_{i \choose n_i-1}, \B_{i \choose n_i}\right)}^{-1}  &  \hskip 5mm \mbox{if } k=n_i. \\
    \end{array} \right.
 \ee
Here $\mathcal{W}_i$ denotes the potential function and $r^+$ is the cross ratio. It is easy to see that 
\be \la{211a}
\displaystyle{r^+\left(\B_1, \B_{i \choose k-1}, \B_{i \choose k}, \B_{i \choose k+1}\right) = \frac{\mathcal{W}_i\left(\A_1, \B_{i \choose k}, \B_{i \choose k+1}\right)} {\mathcal{W}_i\left(\A_1, \B_{i \choose k-1}, \B_{i \choose k}\right)}}.
\ee

Next, using the pinning $p_{23}=(\A_2, \A_3')$ over $(\B_2, \B_3)$,   sequence \eqref{2018.9.22.17.41ssX} has a decomposition
\[
\xymatrix{\A_2=\A_2^0 &\ar[l]_{\small \qquad s_{i_1^\ast}}\A_2^1&\ar[l]_{\small s_{i_2^\ast}}\A_2^2&\ar[l]_{\small s_{i_3^\ast}}\cdots &\ar[l]_{\small s_{i_m^\ast}\ \  } \A_2^m=\A_3'}
\]
such that the $h$-distance between any neighbored decorated flag is $1$. 
Let ${\rm A}_1$ be an arbitrary decorated flag over $\B_1$. We define the {\it primary coordinates}:
\be\la{PC}
P_{{\bf i}, k} = \frac{\Delta_{i_k}\left({\rm A}_1, \A_2^{k}\right)}{\Delta_{i_k}\left({\rm A}_1, \A_2^{k-1}\right) }, \hskip 15mm \forall ~ k=1,\ldots, m.
\ee
They are naturally attached to the edges of the bottom side of the right triangle   on the Figure above. They do not depend on the decoration of $\B_1$ chosen. 

Recall the positive coroot
$\beta_{k}^{\bf i}$ in \eqref{reduced.word.positive.roots}. 
If $\beta_{k}^{\bf i}=\alpha_i^\vee$ is a simple positive coroot, then we set 
\be
\la{bottomXY}
X_{i \choose -\infty} = P_{{\bf i}, k}. 
\ee

\bd \la{DEFXC} The  functions $X_{i \choose k}$ in (\ref{defxikX}) and    $X_{i \choose -\infty}$ in (\ref{bottomXY}) are  the cluster Poisson coordinates on   $\mathscr{P}_{\G', t}$ assigned to a reduced word ${\bf i}$ for $w_0$ and the vertex $\B_1$ of the triangle $t$.
\ed

\bl \la{L5.10} i) The potential $\mathcal{W}_i(\A_1, \B_2, \B_3)$ is  a sum of cluster monomials: 
\be \la{FFPW}
\mathcal{W}_i(\A_1, \B_2, \B_3)=\sum_{0\leq k < n_i} X_{i \choose 0}X_{i \choose 1}\cdots X_{i \choose k}.
\ee
Meanwhile,
\be  
\la{hdist.pspace}
\alpha_i\left(h({\A_1', \A_1})\right)=X_{i \choose 0}X_{i \choose 1}\cdots X_{i \choose n_{i}}.
\ee

ii) Let ${\bf i}$ be a reduced word of $w_0$ starting with $i_1=i$. One has
 \be
 \la{W.dist2}
 \begin{split}
& \mathcal{W}_{i^\ast}(\A_2, \B_3, \B_1) = P_{{\bf i}, 1},\\
& {\alpha}_{i^\ast}\left(h(\A_2', \A_2)\right)= P_{{\bf i}, 1} X_{i \choose 0}.\\
 \end{split}
 \ee
 
 iii) Let ${\bf i}$ be a reduced word of $w_0$ ending with $i_m=i$. One has
 \be
  \la{W.dist3}
 \begin{split}
& \mathcal{W}_{i^\ast}(\A_3, \B_1, \B_2) = X_{i \choose n_i},\\
& {\alpha}_{i^\ast}\left(h(\A_3', \A_3)\right)=  X_{i \choose n_i} P_{{\bf i},m}.\\
 \end{split}
 \ee
 
\el
\begin{proof}
i) Note that 
\[
\mathcal{W}_i(\A_1, \B_2, \B_3)=\sum_{0\leq k <n_i} \mathcal{W}_i\left(\A_1, \B_{i \choose k}, \B_{i \choose k+1}\right). 
\]
For $0\leq k< n_i$, the first two lines in (\ref{defxikX}) together with (\ref{211a}) imply that  
\[
X_{i \choose 0}X_{i \choose 1}\cdots X_{i \choose k} = \mathcal{W}_i\left(\A_1, \B_{i \choose k}, \B_{i \choose k+1}\right).
\]
Using (\ref{defxikX}) and (\ref{211a}), we get 
\[
X_{i \choose 0}X_{i \choose 1}\cdots X_{i \choose n_{i}} = \frac{ \displaystyle{\mathcal{W}_i\left(\A_1, \B_{i \choose n_i-1}, \B_{i \choose n_i}\right)}} { \displaystyle{\mathcal{W}_i\left(\A_1', \B_{i \choose n_i-1}, \B_{i \choose n_i}\right)}}=\alpha_i\left(h({\A_1', \A_1})\right)= \alpha_i\left(h(\A_1', \A_2')\right).
\]  
ii) The first formula of \eqref{W.dist2} is proved via \eqref{leftpotential.1}. Using \eqref{bfz.lusztig}, we get
\[
P_{{\bf i}, 1} X_{i \choose 0} = \alpha_i\left(h(\A_1, \A_2)^{-1}\right)= \alpha_{i^\ast}(h(\A_2', \A_2)).
\]
iii) It follows by using the same arguments as in ii).
\end{proof}

\subsubsection{Cluster Poisson coordinates on $\mathscr{P}_{\G, \bS}$.} Pick a datum   $\left({\cal T}, \{{\bf i}_t\}, \{{v}_t\}\right)$, where ${\cal T}$ is an ideal triangulation of $\bS$, ${\bf i}_t$ is a 
reduced word for $w_0$ assigned to a triangle $t$ of ${\cal T}$, and $v_t$ is a vertex of the triangle $t$.  Let us amalgamate the cluster Poisson coordinate systems assigned in 
Definition \ref{DEFXC} to the spaces ${\mathscr P}_{\G, t}$, when $t$ runs through the triangles of ${\cal T}$, following the   amalgamation procedure \cite{FG05}, recalled in 
Section \ref{sec2}. The resulting coordinates on   $\mathscr{P}_{\G', \bS}$ are described as follows. 

The ideal triangulation ${\cal T}$ provides  the gluing  map associated with   the triangles $t$ of ${\cal T}$. 

It is obtained by taking first 
the product of the moduli spaces ${\mathscr P}_{\G, t}$ over all triangles $t$ of the triangulation ${\cal T}$;  taking then the quotient by the   action of the product of the copies of the Cartan group $\H$ parametrised by the internal edges of ${\cal T}$; and finally using    an open embedding to $\mathscr{P}_{\G', \bS}$: 
\be \la{REST}
{\gamma}_{\cal T}: \prod_{t \in {\cal T}} \mathscr{P}_{\G', t} \lra \frac{\prod_{t \in {\cal T}} \mathscr{P}_{\G', t}}{\H^{\{\mbox{\rm internal edges of ${\cal T}$}\}} } \hra \mathscr{P}_{\G', \bS}.
\ee
So to define a collection of rational functions on the space $\mathscr{P}_{\G', {\cal T}}$ we consider the following collection of rational functions on  
$\prod_{t \in {\cal T}} \mathscr{P}_{\G', t}$,  invariant under the action of the group 
$\H^{\{\mbox{\rm internal edges of ${\cal T}$}\}}$:

\begin{enumerate}

\vskip 1mm\item  For each   triangle $t$ of   ${\cal T}$, we take all unfrozen cluster Poisson coordinates $X^{t, {\rm nf}}_j$ on  $\mathscr{P}_{\G', t}$. 
They are  preserved by the action of the group $\H^3$. So we get a collection of rational functions on the space 
$\mathscr{P}_{\G', {\cal T}}$, denoted by 
\be  \la{e3}
X^{t, {\rm nf}}_j.
\ee

\vskip 1mm\item  Recall that we assign to each side of a triangle $t$ a set of functions parametrized by the set of simple positive roots ${\rm I}$. These $3{\rm rk}(\G)$ functions    are the 
frozen coordinates on  $\mathscr{P}_{\G', t}$. 

Given an internal edge $e$ of   ${\cal T}$ shared by    triangles $t_1$ and $t_2$, 
  take the frozen coordinates $X^{t_1, e}_i$,   $i \in {\rm I}$, assigned to the edge $e$ of   $t_1$, and similar coordinates $X^{t_2, e}_i$,   and define a collection of coordinates assigned to the edge $e$ by taking their products:
\be \la{e1}
X^{e}_i:=  X^{t_1, e}_i \cdot  X^{t_2, e}_i, \qquad e = t_1 \cap t_2: \mbox{internal edges of ${\cal T}$}, \qquad i \in {\rm I}.
\ee

\vskip 1mm\item  If  $f$ is a boundary interval, there is a single triangle $t$ containing $f$. Take the  frozen coordinates assigned to the edge $f$ of $t$. 
This way we get all frozen coordinates on $\mathscr{P}_{\G', \bS}$:
\be \la{e2}
X^{f}_i:=  X^{t, f}_i, \qquad f:  \mbox{boundary intervals on $\bS$}, \qquad i \in {\rm I}.
\ee

\end{enumerate}

\bd \la{e4} The collection of functions (\ref{e3}), (\ref{e1}) and (\ref{e2}), when $t$ runs through all triangles of ${\cal T}$, form  the set of cluster Poisson coordinates asigned to the datum 
$\left({\cal T}, \{{\bf i}_t\}, \{{v}_t\}\right)$. 
\ed

 \bt \la{MTHPs} There is a unique cluster Poisson structure on the space ${\mathscr P}_{ \G, \bS}$   such that 
for each data $\left({\cal T}, \{{\bf i}_t\}, \{{v}_t\}\right)$,  the  functions  from Definition  \ref{e4} form a cluster 
Poisson coordinate system associated with   quiver (\ref{RDT1}).
\et

The proof of Theorem \ref{MTHPs}  is given in Section  \ref{SEC9.2x}.

 \vskip 2mm \paragraph{\bf Examples.}      Let $(V_2, \langle \ast, \ast\rangle )$ be a 2-dimensional symplectic vector space. A pinning for ${\rm PGL}_2$ can also be viewed as a pair of vectors $l, l'$ with $\langle l,l'\rangle=1$, modulo the $\mathbb{Z}/2$-action
 \[
 (l,l') \longmapsto (-l, -l').
 \]
 By definition, the moduli space ${\mathscr P}_{{\rm PGL_2}, t}$  for a triangle $t$ parametrizes the triples of pinnings $\left( (l_1, l_2'), (l_2, l_3'), (l_3, l_1') \right)$ such that $(l_i, l_{i'})=0$ for $i\in \Z/3\Z$, modulo the action of ${\rm PGL_2}$.

 \begin{center}
\begin{tikzpicture}[scale=0.6]
  \node  (a1) at (105:2.3) {$l_1$};
   \node  (a2) at (75:2.3) {$l_1'$};
     \node  (b1) at (225:2.3) {$l_2$};
   \node  (b2) at (195:2.3) {$l_2'$};
     \node  (c1) at (345:2.3) {$l_3$};
   \node  (c2) at (315:2.3) {$l_3'$};
                  \draw (90:2) -- (210:2) -- (330:2) -- (90:2);
  \foreach \from/\to in {a1/b2, b1/c2, c1/a2}
                  \draw[<->, >=stealth] (\from) -- (\to);                 
 \end{tikzpicture}
 \end{center}

  The cluster Poisson coordinate $X_{13}$ assigned to the side $(13)$ coincides with  the potential ${\cal W}_3$ at the vertex, defined  using the vector $l_3$, as the picture above formula (\ref{EXFO}) illustrates: 
\be \la{EXFO}
 X_{13} = r(l_1, l_2, l_3, l_1'+l_3)= \frac{  \langle l_1, l_2\rangle  \langle l_3, l'_1+l_3\rangle }{ \langle l_2, l_3\rangle \langle l_1, l'_1+l_3\rangle   } = \frac{ \langle l_1, l_2\rangle    }{ \langle l_2, l_3\rangle \langle l_1, l_3\rangle }  = {\cal W}_{3}.
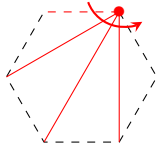
\ee
 Note that this expression is invariant under the maps $(l_i, l'_{i+1}) \lms (-l_i, -l'_{i+1})$.

2. Let $\bS = {\rm P}_n$ be a polygon with $n$ vertices. The moduli space ${\mathscr P}_{ {\rm PGL_2},  {\rm P}_n}$ carries a cluster Poisson variety structure 
of finite type $A_{n}$,  with the  
frozen variables assigned to the sides, and unfrozen ones to the edges of a  triangulation ${\cal T}$ of the polygon. 
Pick a vertex $v$ of the polygon. 
Take a counterclockwise 
oriented arc $\alpha_v$ near the vertex $v$, and denote by $E_1, ..., E_m$   
the  sides and diagonals of the triangulation ${\cal T}$   intersecting the arc, 
counted in the order the oriented arc intersects them. Then it is easy to see that the potential ${\cal W}_v$ at the vertex $v$ is given by 
\be \la{FPWa}
{\cal W}_v:= X_{E_1}  + X_{E_1} X_{E_2} + \ldots + X_{E_1} X_{E_2}  \ldots X_{E_{m-1}}.
\ee
 
 \begin{center}
\begin{tikzpicture}
  \draw[dashed] (60:1) -- (0:1) -- (-60:1) -- (-120:1)--(-180:1) -- (-240:1);
 \draw[red, dashed] (60:1) -- (120:1);     
 \node[red] at (60:1) {{\small $\bullet$}};
 \draw[red] (60:1) -- (180:1);   
  \draw[red] (60:1) -- (240:1);
   \draw[red] (60:1) -- (300:1);   
   \draw[red, -stealth,thick] (85:1) arc(200:300:5mm);    
 \end{tikzpicture}
 \end{center}
 
3. Formula (\ref{FPWa}) allows to  interpret formula (\ref{FFPW}) for  the potential geometrically.

\subsubsection{Cluster   structures of the moduli spaces ${\mathscr P}_{\G, t}$  and ${\mathscr A}_{\G, t}$ by amalgamation.}  
 The cluster structures of the moduli spaces ${\mathscr P}_{\G, \bS}$ is obtained by   the  amalgamation of the moduli spaces ${\mathscr P}_{\G, t}$   assigned to triangles $t$ of an ideal triangulation of a decorated surface $\bS$.  
 
 The cluster Poisson structure on the  space ${\mathscr P}_{t}$ is   defined by the amalgamation of elementary 
cluster Poisson varieties 
${\mathscr P}_{s}$ and $\overline {\mathscr P}_{s}$, assigned to the generators $s$ of the Weyl group $W$. The space ${\mathscr P}_{s}$ was introduced in \cite{FG05}.  Let  $r:= {\rm rk}(G)$ be the rank of $\G$. Then we have   
 $$
{\dim}{\mathscr P}_{s} = r+1, \qquad  {\dim}\overline{\mathscr P}_{s} = r+2, \qquad   {\mathscr P}_{s}  \subset \overline{\mathscr P}_{s}.
$$

\bd \la{DEFP} Let $s$ be a generator of the Weyl group. The moduli space ${\mathscr P}_{s}$ parametrizes the orbits of the group $\G$ acting on the following data:
 $
(\B_0, \B_1, \B_2; p_{10}, p_{12}).
 $
 
Here $(\B_0, \B_1, \B_2)$ is a triple of flags such that the pairs $(\B_0, \B_1)$ and $(\B_0, \B_2)$ are generic, $p_{01}$ and $p_{02}$ are pinnings 
 for them, and the pair 
 $(\B_1,  \B_2)$ is in the relative position $s$. \ed

  A triple of flags  $(\B_0, \B_1, \B_2)$ as in Definition \ref{DEFP} determines uniquely a pair of flags 
 $(\B_1^s, \B_2^s)$ such that $\B_i^s$ is in the relative position $s$ to   $\B_0$, and $w_0s$ to $\B_i$, $i=1,2$. 
 The flags in   the relative position $s$ to  $\B_0$ form a projective line, denoted by  ${\rm P}_{\B_0, s}^1$.

 \bd  The   space $\overline {\mathscr P}_{s}$ parametrizes   $\G$-orbits on the space of  data
  $(\B_0, \B_1, \B_2; p_{10}, p_{12}, p_s)$, where   $(\B_0, \B_1, \B_2; p_{10}, p_{12})$ is as in Definition \ref{DEFP}, and $p_s$ is an {\em $s$-pinning for the triple $(\B_0, \B_1, \B_2)$}, defined as an extra flag ${\rm F}$ in the relative position $s$ to $\B_0$.
 \ed
 
There is a   function $P_s$ on the space $\overline {\mathscr P}_{s}$:    the cross-ratio   of points 
 $(\B_0, \B_1^s, \B_2^s, F)$ on  ${\rm P}_{\B_0, s}^1$:
 $$
 P_s:= r(\B_0, \B_1^s,  {\rm F}, \B_2^s,).
 $$
  
We picture ${\mathscr P}_{s}$ and 
 $\overline {\mathscr P}_{s}$ by ``wedges'' with two sides and a short base. 
On   pictures, the amalgamation amounts to gluing 
the ``wedges'' along their   sides, see  
  Figure \ref{fga31}.

\begin{figure}[ht]
\epsfxsize 200pt
\center{
\begin{tikzpicture}[scale=0.5]
\tikzset{>=open triangle 45}
\foreach \count in {1,2,3}
{
\draw [dashed, ->]  (30-120*\count:2cm) -- (-90-120*\count:2cm);
}
\draw [dashed]  (-90:2cm) -- (-0.57,1);
\draw [dashed]  (-90:2cm) -- (0.57,1);
\begin{scope}[shift={(-7,0)}]
\draw [dashed, ->]  (0,-2) -- (-0.57,1);
\draw [dashed, ->]  (0,-2) -- (0.57,1);
\draw [dashed] (-0.57,1) -- (0.57,1);
\draw [dashed, ->]  (0.3,-2) -- (0.87,1);
\draw [dashed, ->]  (0.3,-2) -- (2.01,1);
\draw [dashed, ->]  (0.87,1) -- (2.01,1);
\draw [xscale=-1, dashed, ->]  (0.3,-2) -- (0.87,1);
\draw [xscale=-1, dashed, ->]  (0.3,-2) -- (2.01,1);
\draw [xscale=-1, dashed, <-]  (0.87,1) -- (2.01,1);
\end{scope}
\begin{scope}[shift={(-13,0)}]
\draw [dashed, ->]  (0,-2) -- (-0.57,1);
\draw [dashed, ->]  (0,-2) -- (0.57,1);
\draw [dashed, ->] (-0.57,1) -- (0.57,1);
\end{scope}
\begin{scope}[shift={(-15,0)}]
\draw [dashed, ->]  (0,-2) -- (-0.57,1);
\draw [dashed, ->]  (0,-2) -- (0.57,1);
\draw [dashed] (-0.57,1) -- (0.57,1);
\node [label=above: {\small $\B_1$}] at (-0.57,1) {};
\node [label=above: {\small $\B_2$}] at (0.57,1) {};
\node [label=below: {\small $\B_0$}] at (0,-2) {};
\end{scope}
\end{tikzpicture}
 }

\caption{Cluster varieties ${\mathscr P}_{s}$ and 
$\overline {\mathscr P}_{s}$ (left); and their amalgamation   into 
${\mathscr P}_{G, t}$ (right).}
\label{fga31}
\end{figure} 
\noindent
  Given a reduced word ${\bf i}$,  $w_0 = s_{i_1} ... s_{i_n}$, we 
  amalgamate   elementary cluster varieties ${\mathscr P}_{s_{i_k}}$ or 
$\overline {\mathscr P}_{s_{i_k}}$:
\be \la{AMALGP}
\stackrel{(\overline{\qquad})}{ {\mathscr P}}_{s_{i_1}}\ast \ldots \ast 
\stackrel{(\overline{\qquad})}{ {\mathscr P}}_{s_{i_n}}. 
\ee
To specify which of ${\mathscr P}_{s_{i_k}}$ or 
$\overline {\mathscr P}_{s_{i_k}}$ is used,  
consider the sequence of coroots $\beta_{{\bf i},k}$ in (\ref{reduced.word.positive.roots}) for   ${\bf i}$. 
Each positive coroot appears  just once. 
Let $\{j_1 <  \ldots <j_r\}$ 
be the sequence of $i_k$'s such that 
$\beta^{\bf i}_{k}:= s_{i_n} ... s_{i_{k+1}}\alpha^\vee_{i_k}$ is a simple positive coroot. 
We get an ordered set of simple 
positive coroots 
\be \la{mea11aa}
\{\beta_{j_1}, \beta_{j_2}, \ldots , \beta_{j_r}\}. 
\ee
We use the  
variety  
$\overline {\mathscr P}_{s_{i_k}}$ if and only if $i_k \in \{j_1, ..., j_r\}$.  
 By the amalgamation, the 
space (\ref{AMALGP}) inherits a cluster Poisson structure. To obtain 
the   cluster Poisson structure  on the space  ${\mathscr P}_{\G, t}$ we amalgamate (\ref{AMALGP})  with 
the Cartan group $\H$ equipped with a cluster 
structure determined by sequence (\ref{mea11aa}).   
This does not change the variety.   
We prove that the obtained cluster Poisson structure does not depend 
on the vertex of the triangle   used to cut the triangle $t$ into 
the   wedges.

 \vskip 3mm
    The cluster   $K_2$-structure of the space ${\mathscr A}_{\G, \bS}$ is obtained using the same amalgamation  pattern.


\medskip

\section{Tensor invariants and cluster coordinates on the  space ${\rm Conf}_3({\cal A})$} \la{SECT4}
\medskip
 
For ${\rm G}={\rm SL}_{n}$, there is a special cluster $K_2$-coordinate system on the  space ${\rm Conf}_3({\cal A})$ given by a collection of
  functions $\Delta_{a,b,c}$, $a+b+c=n$, \cite{FG03a}. 
  
  In Section \ref{SECT4} we define a supply of functions $\Delta_{\lambda, \mu, \nu}$ and their specializations 
 $\Gamma_\lambda$,   generalizing 
 functions $\Delta_{a,b,c}$ to any  $\G$. We show in Proposition \ref{equivalence.cluster.a.ba} that they contain  cluster coordinate systems on the space ${\rm Conf}_3({\cal A})$ 
 defined in Section \ref{SeCC2.4A}.

\subsection{Tensor invariants and  coordinate rings of configuration spaces} \la{CURSEC5.1}

\medskip

We denote by  ${\rm Conf}_m^{\times}({\cal A})$ the open subspace of  ${\rm Conf}_m({\cal A})$ parametrizing decorated flags 
$(\A_1, \ldots, \A_m)$ such that the underlying flags $(\B_i, \B_{i+1})$ are in generic position, $i \in \Z/m\Z$.

\subsubsection{Representations and the coordinate ring of ${\cal A}$.}
\medskip

The Peter-Weyl Theorem implies that 
$
{\cal O}({\cal A}) = \bigoplus_{\lambda} V_{\lambda},
$ 
where the  sum is over all finite dimensional  irreducible representations $V_\lambda$ of $\G$. Precisely, $V_{\lambda}$ consists of  functions $F \in {\cal O}({\cal A})$ such that
$$
F(\A\cdot h) = \lambda(h)F(\A), \hskip 7mm \forall h\in {\rm H}, ~~\forall \A\in {\cal A}.
$$
The representation $V_{\lambda}$ admits a weight decomposition
\[
V_{\lambda}=\bigoplus_{\mu} V_{\lambda}(\mu), \hskip 10mm \mbox{where } V_\lambda(\mu):=\{F\in V_{\lambda} ~|~ F(h \cdot \A) = \mu(h)F(\A)\}.
\]

\bl \la{nontrivial.highest.weight}
Let $F\in V_{\lambda}(\mu)$ be a nonzero function. Then 
$F([{\rm U}])\not =  0$ if and only if $\mu = \lambda$.
\el 
\begin{proof}
Note that
$ [\U] = h^{-1}\cdot [\U]\cdot h
$ for all $h\in {\rm H}$. Therefore 
\[
F ([\U]) = F(h^{-1} \cdot [\U]\cdot h)= F([\U])\cdot (\lambda-\mu)(h).
\]
If $\mu \neq \lambda$, then $F([\U])=0$. 
Suppose that $\mu=\lambda$ and $F([\U])=0$. Note that $\dim V_{\lambda}(\lambda)=1$.  So $f([\U])=0$ for all $f \in V_{\lambda}$. The group  $\G$ acts transitively on ${\cal A}$. 
So $F=0$.  Contradiction. 
\end{proof}

\bd
\la{highest.weight.vector.tt}
The function $F_{\lambda}$  is the unique   function in $V_{\lambda}(\lambda)$ such that
$
F_{\lambda}([\U])=1. 
$ 
\ed

The function $F_{\lambda}$  exists by Lemma \ref{nontrivial.highest.weight}.  We abbreviate $F_{\Lambda_i}=F_i$ for $i\in {\rm I}$.

\bl For any dominant weights $\lambda, \mu$, we have
$
F_{\lambda}F_{\mu} = F_{\lambda+\mu}.
$
\el
\begin{proof} Clearly $F_{\lambda}F_{\mu}\in V_{\lambda+\mu}(\lambda+\mu)$, and 
$F_{\lambda}([\U])F_{\mu}([\U]) =1$.
\end{proof} 

The Weyl group $W$ acts on  ${\cal O}({\cal A})$ by  
$w^* F(\A):= F(\overline{w}^{-1}\cdot \A)$. The next Lemma is well known.

\bl 
\la{w.invairnat.high.w}
Let $\lambda$ be a dominant weight. Then $w^*F_\lambda= F_{\lambda}$ for any $w\in W$ such that $w\cdot \lambda =\lambda$.
\el

Let $\mu$ be a weight. Then there is a unique dominant weight $\lambda$ such that    $\mu=w\lambda$ for some    $w\in W$. Set
\[
F_{\mu}(\A):= F_{\lambda}(\overline{w}^{-1}\cdot \A) \in V_{\lambda} (w \lambda).
\]
  By Lemma \ref{w.invairnat.high.w}, the function $F_\mu$ is well defined: it   does not depend on $w$ chosen.   The functions  $F_\mu$ are essentially the {\it generalized minors} of Fomin-Zelevinsky.

\subsubsection{Double tensor invariants and ${\rm Conf}_2({\cal A})$.}
Recall that $\lambda^\ast = -w_0(\lambda)$.  We have 
\[
{\mathcal{O}\big({\rm Conf}_2({\cal A})\big)= {\cal O}({{\cal A}^2})^{\G}{=} \bigoplus_{\lambda}} \big(V_{\lambda}\bigotimes V_{\lambda^*}\big)^{\G}.
\]
Note that  $\dim \big(V_{\lambda}\bigotimes V_{\lambda^*}\big)^{\G} =1$.

\bl Let $\Delta\in \big(V_{\lambda}\bigotimes V_{\lambda^*}\big)^{\G}$ be a nontrivial function. Then 
$
\Delta([{\rm U}], [\U^{-}])\neq 0.
$
\el
\begin{proof} Suppose $\Delta([{\rm U}], [\U^{-}])=0$. Then for any generic pair $(\A_1, \A_2)$, we have
\be
\la{basic.invariants.pair.decorated}
\Delta(\A_1, \A_2) = \Delta \big(h(\A_1, \A_2) \cdot [\U], [\U^{-}]\big)= \lambda\big(h(\A_1, \A_2)\big) \Delta([\U], [\U^{-}]) =0.
\ee
So the function $\Delta$ is zero on a  Zariski dense subset of ${\cal A}^2$. Thus it is zero. Contradiction.
\end{proof}
\bd Let $\Delta_{\lambda}\in \big(V_{\lambda}\bigotimes V_{\lambda^*}\big)^{\G}$ be the   function such that
$
\Delta_{\lambda}([\U], [\U^{-}])= 1.
$ 
\ed
We abbreviate $\Delta_{\Lambda_i}=\Delta_i$ for $i\in {\rm I}$.

\bl
For any dominant weights $\lambda$, $\mu$, we have
$
\Delta_{\lambda} \Delta_{\mu} = \Delta_{\lambda+\mu}.
$

If $(\A_1, \A_2)\in {\rm Conf}^\times_2({\cal A})$, then 
$
\Delta_{\lambda}(\A_1, \A_2) = \lambda\big(h(\A_1, \A_2)\big).
$
\el
\begin{proof} It is a direct consequence of the definition of $\Delta_\lambda$ and \eqref{basic.invariants.pair.decorated}.
\end{proof}
\bl \la{delta=f}
$
\Delta_{\lambda}(\A, [\U^-]) = F_{\lambda}(\A).
$
\el
\begin{proof} Clearly the function $\Delta_{\lambda}(\ast, [\U^-])\in V_{\lambda}$. Note that
\begin{align}
\Delta_{\lambda}(h \cdot \A, [\U^-]) &= \Delta_{\lambda}(\A, h^{-1}\cdot [\U^-])= \Delta_{\lambda}(\A, [\U^-]\cdot w_0(h^{-1})) \nonumber \\ 
&= \Delta_{\lambda}(\A, [\U^-]) \lambda^*(w_0(h^{-1})) 
= \Delta_{\lambda}(\A, [\U^-])  \lambda(h). \nonumber
\end{align}
Therefore $\Delta_{\lambda}(\ast, [\U^-])\in V_{\lambda}(\lambda)$. The Lemma follows from  
$
\Delta_{\lambda}([\U], [\U^-]) = F_{\lambda}([\U])=1.
$
\end{proof}

\bl 
\la{generic.5.22.19.23}
Suppose that $w(\A_1, \A_2)=uw_0$. Then
\[
\Delta_{\lambda}(\A_1, \A_2)=   \left\{\begin{array}{ll} 
      \lambda\left(h(\A_1, \A_2)\right) & \mbox{if $u\cdot \lambda=\lambda$}, \\
      0 & \mbox{otherwise}. \\
   \end{array}\right.
\]
In particular $\Delta_i(\A_1, \A_2)\neq 0$ for all $i\in {\rm I}$ if and only if $(\A_1, \A_2)$ is of generic position.
\el
\begin{proof} Let $(\A_1, \A_2) = (\overline{u}^{-1}\cdot [\U]\cdot h, [\U^-])$. 
By Lemma \ref{delta=f}, we have
\[
\Delta_{\lambda}(\A_1, \A_2) = F_\lambda (\overline{u}^{-1}\cdot [\U]\cdot h)= F_\lambda (\overline{u}^{-1}\cdot [\U]) \lambda(h).\]
If $u\lambda =\lambda$, then 
$F_{\lambda}(\overline{u}^{-1}\cdot [\U])= F_{\lambda}([\U])=1$. 
Otherwise, 
$$
F_\lambda (\overline{u}^{-1} [\U]\cdot h) = F_\lambda(u^{-1}(h) \overline{u}^{-1} [\U]) = (u\lambda)(h) \cdot F_\lambda(\overline{u}^{-1} [\U]).
$$ 
So $\Delta_{\lambda}(\A_1, \A_2) =0$.
\end{proof}

\bl 
\label{prop.11.12.13} 
For each pair $(\A_1, \A_2)\in {\rm Conf}_2({\cal A})$, we have
$
\Delta_{\lambda}(\A_1, \A_2) = \Delta_{\lambda^*} (\A_2\cdot s_\G, \A_1).
$
\el

\begin{proof} 
By definition $(\A_1, \A_2)\lms \Delta_{\lambda^*} (\A_2\cdot s_\G, \A_1)$ is a regular function in $(V_\lambda\bigotimes V_{\lambda^\ast})^\G$. Meanwhile
\[
\Delta_{\lambda^*} ([\U^-]\cdot s_\G, [\U]) = \Delta_{\lambda^*}(s_\G^2\cdot [\U], \overline{w}_0 \cdot [\U]) = \Delta_{\lambda^*}( [\U],[\U^-])=1.
\]\end{proof}

\begin{remark} Lemma \ref{prop.11.12.13} implies that the linear basis $\{\Delta_{\lambda}\}$ of ${\cal O}\big({\rm Conf}_2({\cal A})\big)$ is preserved by the automorphism of 
${\rm Conf}_2({\cal A})$ given by $(\A_1, \A_2)\lms (\A_2\cdot s_{\G}, \A_1).$
\end{remark}

\subsubsection{Triple tensor invariants and ${\rm Conf}_3({\cal A})$.} There is a  decomposition
$$
{\cal O}({\rm Conf}_3({\cal A}))= \bigoplus_{\lambda, \mu, \nu} (V_{\lambda}\otimes V_{\mu}\otimes V_{\nu})^\G.
$$
Here  $(V_{\lambda}\otimes V_{\mu}\otimes V_{\nu})^\G$ consists of $\G$-invariant functions $\Delta$ of ${\cal A}^3$ such that
\be
\la{prop.1.5hhh}
\Delta (\A_1\cdot h_1, \A_2 \cdot h_2, \A_3\cdot h_3) =\Delta(\A_1, \A_2, \A_3) \lambda(h_1)\mu(h_2)\nu(h_3).
\ee

\bd A triple of dominant weights $(\lambda, \mu, \nu)$ is  {admissible} if there is an   $w\in W$ so that \be
\la{basic.triple.condition}
w \lambda = \nu^\ast - \mu, ~~~~\mbox{where $\nu^\ast := -w_0(\nu)$}.
\ee\ed

\bl \la{LLL1} For any admissible triple  $(\lambda, \mu, \nu)$   of dominant weights  we have 
$$
\dim (V_{\lambda}\otimes V_{\mu}\otimes V_{\nu})^\G =1.
$$
\el

\begin{proof} 
Note that
$
\dim (V_{\lambda}\otimes V_{\mu}\otimes V_{\nu})^\G \leq \dim V_{\lambda}(\nu^\ast - \mu) = \dim V_{\lambda} (w\cdot \lambda)= 1.
$

Due to the PRV conjecture proved by Kumar \cite{Ku}, 
$
\dim (V_{\lambda}\otimes V_{\mu}\otimes V_{\nu})^\G \geq 1.
$
\end{proof}

\medskip

\subsection{The functions $\Delta_{\lambda, \mu, \nu}$}

\medskip

\blc  Let  $(\lambda, \mu, \nu)$ be an admissible triple. Then there is a unique  function 
\be
\la{prop.1.4.hhh}
\Delta_{\lambda, \mu, \nu}\in (V_{\lambda}\otimes V_{\mu}\otimes V_{\nu})^\G, ~~ \mbox{normalized by} ~~\Delta_{\lambda, \mu, \nu}(\A, [{\rm U}], [{\rm U}^-]) = F_{w\lambda}(\A)\cdot \mu(s_{\G}).
\ee
\elc

\begin{proof}
Let $\delta_{\lambda, \mu, \nu}\in (V_{\lambda}\otimes V_{\mu}\otimes V_{\nu})^\G$.  
Define a function $f$ on ${\cal A}$ by setting 
\[
f(\A):=  \delta_{\lambda, \mu, \nu}( \A, [{\rm U}], [{\rm U}^-]).
\]

\bl 
If $\delta_{\lambda, \mu, \nu} \not =0$, then the function $f$ is a nontrivial function in $V_{\lambda}(w\lambda)$.
\el
\begin{proof}
Clearly $f\in V_{\lambda}$. We have $f \in V_{\lambda}(w\lambda)$ since
\begin{align}
&f(h\cdot \A)=   \delta_{\lambda, \mu, \nu}(h\cdot \A, [{\rm U}], [{\rm U}^-])  =  \Delta_{\lambda, \mu, \nu}(\A, h^{-1}\cdot [{\rm U}], h^{-1}\cdot [{\rm U}^-])=\nonumber\\
& \delta_{\lambda, \mu, \nu}(\A,  [{\rm U}]\cdot h^{-1},  [{\rm U}^-]\cdot  w_0(h^{-1}) )
 = \delta_{\lambda, \mu, \nu}(\A,  [{\rm U}],  [{\rm U}^-] )  \cdot  \big(\nu^*-\mu \big) (h) 
 = f(\A) \cdot (w\lambda)(h). \nonumber
\end{align}
If $f=0$ then for any $(\A_1, \A_2, \A_3)$ in the $\G$-orbit of $(\A, [\U], [\U^-]\cdot h )$, we have
\[
\delta_{\lambda, \mu, \nu}(\A_1, \A_2, \A_3) = f(\A) \nu(h) =0.
\]
Such $(\A_1, \A_2, \A_3)$ are  Zariski  dense in ${\cal A}^3$. So $\delta_{\lambda, \mu, \nu}=0$.   Contradiction. 
\end{proof}
Since $\dim (V_{\lambda}\otimes V_{\mu}\otimes V_{\nu})^\G =1$, the normalization condition single out a unique function. \end{proof}


\subsubsection{An alternative description of the function $\Delta_{\lambda, \mu, \nu}$.}

Let $(\A_1, \A_2, \A_3)$ be a triple with generic $(\A_2, \A_3)$. For any  $w\in W$,  
 by Lemma \ref{lem.part.flag.decom}, there is a unique flag $\B_4$ such that
\be
\la{intermediate.flag.ss}
 w(\B_2, \B_4) = w w_0, \hskip 7mm w(\B_3, \B_4)= w^*. 
\ee
 Let $(\lambda, \mu, \nu)$ be an admissible triple, and $w\in W$ is as in \eqref{basic.triple.condition}.  
 Choose a decorated flag $\A_4$ over  $\B_4$, and set
\be
\la{2017.4.14.23.hs}
A_{\lambda, \mu, \nu} :=  \Delta_{\lambda}(\A_1, \A_4)\cdot {\mu}\left(h(\A_2, \A_4)\cdot \rho_{w_0w^{-1}}(-1)\right)\cdot  {\nu}\left(h(\A_3, \A_4)\right).
\ee
\bl
The function \eqref{2017.4.14.23.hs} is independent of the decoration of $\B_4$ chosen.
\el
\begin{proof} Set
$h_{ij} = h(\A_i, \A_j). $
Let us rescale $\A_4\lms \A_4 \cdot h$. Then 
\begin{align}
 \Delta_{\lambda}(\A_1, \A_4) &\lms  \Delta_{\lambda}(\A_1, \A_4) \cdot \lambda(h^*) , \nonumber\\
 {\mu}(h_{24}) &\lms  \mu(h_{24})\cdot (w^{-1} \mu)(h^*), \nonumber\\
  {\nu}(h_{34}) &\lms \nu(h_{34})\cdot  (w^{-1}w_0 \nu)(h^*). \nonumber
 \end{align}
 Note that
 \[
\lambda + w^{-1}\mu + w^{-1}w_0\nu =\lambda+ w^{-1}(\mu -\nu^\ast )\stackrel{(\ref{basic.triple.condition})}{=} 0.
 \]
 The product \eqref{2017.4.14.23.hs} remains intact.
\end{proof}

\bt
The  restriction of the  function $\Delta_{\lambda, \mu, \nu}$ to ${\rm Conf}_3^{\times}({\cal A})$ is equal to   $A_{\lambda, \mu, \nu}$.
\et
\begin{proof} The function $\Delta_{\lambda, \mu, \nu}$ restricted on ${\rm Conf}_3^{\times}({\cal A})$ is characterized by \eqref{prop.1.5hhh} and \eqref{prop.1.4.hhh}.
Clearly $A_{\lambda, \mu, \nu}$ satisfies \eqref{prop.1.5hhh}.  Suppose that $(\A_1, \A_2, \A_3)=(\A, [\U], [\U^{-}])$. Choose $\A_4 = \overline{w}\cdot [\U^{-}]$. Then $h_{24}\rho_{w_0w^{-1}}(-1)=s_\G$ and $h_{34}=1$. Therefore 
\[
A_{\lambda, \mu, \nu }\big(\A, [\U], [\U^{-}]\big) =\Delta_{\lambda}( \A,   \overline{w}\cdot [\U^{-}] )\cdot \mu(s_{\G})= \Delta_{\lambda}(  \overline{w}^{-1}\cdot\A,   [\U^{-}] ) \cdot \mu(s_{\G})= F_{w\lambda}(\A)\cdot \mu(s_{\G}).
\]
So $A_{\lambda, \mu, \nu}$ satisfies \eqref{prop.1.4.hhh}.
\end{proof}

\subsection{The functions $\Gamma_\lambda$ and cluster coordinates on the space ${\rm Conf}_3({\cal A})$} \la{CSEC12.3}

\medskip \paragraph{The function $\Gamma_\lambda$.}
The fundamental weights  $\Lambda_1, \ldots, \Lambda_r$ form a basis of the weight lattice.  Let $[n]_+=\max\{n,0\}$. Every weight $\lambda = \sum_{i\in {\rm I}} n_i \Lambda_i$ admits a decomposition 
 \be 
 \la{triple.dominant.weights}
\lambda=\lambda_+ - \lambda_-, \hskip 7mm \mbox{where~}
 \lambda_+ = \sum_{i\in {\rm I}} [n_i]_+ \Lambda_i, \hskip 5mm \lambda_-= \sum_{i\in {\rm I}} [-n_i]_+ \Lambda_i.
 \ee
 The $W$-orbit of the weight $\lambda$  contains a unique dominant weight $ \overline{\lambda} = w^{-1}{\lambda}$, where $w\in W$.  
 The  triple $(\overline{\lambda},\lambda_-, \lambda_+^*)$  satisfies  condition \eqref{basic.triple.condition}, and thus is admissible. So we arrive at

 \blc 
\la{function.gamma.alpha}
Every weight $\lambda$ naturally gives rise to a  function
 \be
 \Gamma_{\lambda}:= \Delta_{\overline{\lambda},\lambda_-, \lambda_+^*} \in {\cal O}\big({\rm Conf}_3({\cal A})\big).
 \ee
\elc

\bl 
\la{decomp.alpha.plus.minus}
Let ${\lambda}=w \overline{\lambda}$ be on the $W$-orbit of a dominant weight $\bar{\lambda}$. We have
\[
\lambda_+ = \sum_{i\in {\rm I}(w_0w^{-1})} \langle \alpha_i^\vee, \lambda\rangle \Lambda_i , \hskip 7mm  \lambda_- = - \sum_{i  \in {\rm I}(w^{-1})} \langle \alpha_i^\vee, \lambda\rangle \Lambda_i.
\]
\el

\begin{proof} Recall that  the subset ${\rm I}(w_0w^{-1})$ of ${\rm I}$ is the complement of ${\rm I}(w^{-1})$. 
Meanwhile,
\[
\lambda= \sum_{i\in {\rm I}}  \langle \alpha_i^\vee, \lambda\rangle \Lambda_i = \sum_{i\in I} \langle w^{-1}\cdot \alpha_i^\vee, \bar{\lambda}\rangle \Lambda_i.
\]
By definition, the coroot $w^{-1}\cdot \alpha_i^\vee$ is negative if and only if $i \in {\rm I}(w^{-1})$.
\end{proof}

Let $(\A_1, \A_2, \A_3)\in {\rm Conf}_3^{\times}({\cal A})$. By Lemma \ref{basic.decomp.chain.w0}, there is a decoration of   $\B_4$ in  \eqref{intermediate.flag.ss} such that 
\[
h(\A_3, \A_4) \in {\rm H}(w^{*-1}), \hskip 7mm h(\A_2, \A_4) \cdot \rho_{w_0w^{-1}}(-1) \in {\rm H}(w_0w^{-1}).
\]
\bp 
\la{equivalence.cluster.a.ba}
$
\Gamma_{\lambda} (\A, \A_l, \A_r) = \Delta_{\overline{\lambda}}(\A, \A_k).
$
\ep
\begin{proof} By Lemma \ref{decomp.alpha.plus.minus}, we have
$\lambda_+^* \left(h(\A_3, \A_4) \right) =1$ and $\lambda_-\left(h(\A_2, \A_4) \cdot \rho_{w_0w^{-1}}(-1)\right) =1$. 
\end{proof}

 Every reduced word ${\bf i}=(i_1,\ldots, i_m)$ of $w_0$ gives rise to a sequence of
{\it chamber weights}:
\be
\la{triple.chamber.weights}
\gamma_k^{\bf i}:= s_{i_m}\cdots s_{i_{k+1}}\cdot \Lambda_{i_k}.
\ee
It gives rise to a sequence of  regular functions $\Gamma_{{\bf i},k} : = \Gamma_{\gamma_k^{\bf i}}$ on ${\rm Conf}_3({\cal A})$.

If $(\A_l, \A_r)$ is in generic position, by Lemma \ref{basic.decomp.chain.w0}, we get a   sequence of  decorated flags
\[
\xymatrix{\A_l=\A_0 \ar[r]^{\small ~~~~s_{i_1^\ast}}&\A_1\ar[r]^{\small s_{i_2^\ast}}&\A_2\ar[r]^{\small s_{i_3^\ast}}&\cdots\ar[r]^{\small s_{i_m^\ast}~~~~}&\A_m=\A_r}
\]
By Proposition \ref{equivalence.cluster.a.ba}, we have
\be
\la{equivalence.cluster.a.ba.eq}
\Gamma_{{\bf i},k} (\A, \A_l, \A_r) = \Delta_{i_k}(\A, \A_k).
\ee

 \begin{example}
Let ${\rm G}={\rm SL}_{n+1}$. We consider the reduced word
$
{\bf i}=(1, 2, 1, 3,2,1..., n, ..., 1).
$ 
Then the functions $\Gamma_{{\bf i},{k}}$ coincide with the canonical functions $\Delta_{a,b,c}$ in \cite{FG03a}.
\end{example}

\subsubsection{Cluster $K_2$-coordinate systems on ${\rm Conf}_3(\mathcal{A})$, revisited.}  Let $\omega_{{\bf i}, k}=s_{i_m}\cdots s_{i_{k+1}}$. Define
\[
\begin{split}
&V_{{\bf i}, i}=\left \{ (\Lambda_i, [\omega_{{\bf i}, k} \cdot \Lambda_{i}]_-,  [\omega_{{\bf i}, k} \cdot \Lambda_{i}]_+^\ast)~|~ k=0,\ldots,m\right\}, ~~~~i\in {\rm I},\\
&V'=\left\{ (0, \Lambda_i, \Lambda_{i^*})~|~ i \in {\rm I}\right\}.\\
\end{split}
\]
Observe that the quiver assigned in Section \ref{sec.3.1} to a reduced decomposition ${\bf i}$ of $w_0$ has the  vertices  labelled by the set
$V:= V' \cup \cup_{i\in {\rm I}} V_{{\bf i}, i}$. 
Note that  each $v:=(\lambda, \mu, \nu)\in V$ is admissible. 
So we get a collection of regular functions
\[
A_v=\Delta_{\lambda, \mu, \nu} \in {\cal O}({\rm Conf}_3(\mathcal{A})), ~~~~v \in V.
\]
 They  coincide with 
the cluster $K_2$-coordinate system on ${\rm Conf}_3(\mathcal{A})$ defined in Definition \ref{DEF9.1}.


\medskip

\section{Partial configuration spaces}
\la{pcs.sec}

\medskip

In this Section, $\G'$ is simply connected, and $\G$ is its adjoint group.

\medskip 

\subsection{Main result}

\medskip

Recall the subgroup ${\rm H}(w)$ in \eqref{H.5.26.16.hh} and $\rho_w(-1)={\overline{w^{-1}}}{\overline{w}}$.

\bd \la{def10.1} Let $(u, v)\in W\times W$. The partial configuration space ${\rm Conf}_u^v({\cal B})$ parametrizes $\G$-orbits of quadruples $(\B_l,\B_r,  \B^l, \B^r)$ of flags such that\footnote{Note that $w(\A_l, \A_r)=u^\ast$ is equivalent to $w_-(\A_l, \A_r)=u$. }
\[
w(\B_l, \B^l) = w(\B_r, \B^r)= w_0,  \hskip 10mm w(\B_l, \B_r)= u^*,  \hskip 5mm w(\B^l, \B^r)= v.
\]
A configuration in ${\rm Conf}_u^{v}({\cal B})$ is illustrated on Figure \ref{partial.confbruhatB}:
\begin{figure}[ht]
\begin{center}
\epsfxsize100pt
\begin{tikzpicture}[scale=0.6]
\node [circle,draw=red,fill=red,minimum size=3pt,inner sep=0pt,label=below:{\small ${\B}_l$}]  (a) at (-1,-3){};
\node [circle,draw=red,fill=red,minimum size=3pt,inner sep=0pt,label=above:{\small ${\B}^l$}] (b) at (-1,0) {};
 \node [circle,draw=red,fill=red,minimum size=3pt,inner sep=0pt,label=above:{\small ${\B}^r$}] (c) at (0.8,0) {};
  \node [circle,draw=red,fill=red,minimum size=3pt,inner sep=0pt,label=below:{\small $~~{\B}_r$}] (d) at (0.8,-3) {};
 \foreach \from/\to in {b/a, c/d}
                   \draw[thick] (\from) -- (\to);
   \draw[draw, directed, thick] (b) -- (c);
\draw[draw, directed, thick] (a) -- (d);    
 \node[blue] at (-0.1,0.33) {$v$};
 \node[blue] at (-0.1,-2.5) {$u^\ast$};
  \node[blue] at (-1.55,-1.5) {$w_0$};
 \node[blue] at (1.35,-1.5) {$w_0$};
 \end{tikzpicture}
 \caption{The partial configuration space ${\rm Conf}_u^v({\cal B})$.}
 \label{partial.confbruhatB}
\end{center}
\end{figure}

 The space ${\rm Conf}_u^{v}({\cal A})$ parametrizes $\G'$-orbits of quadruples $(\A_l,\A_r,  \A^l, \A^r)$ of decorated flags with quadruples of underlying flags in $ {\rm Conf}_u^{v}({\cal B})$ and 
\[
h(\A_r, \A_l)\in {\rm H}(u^\ast),  \hskip 15mm h(\A^r, \A^l)\in {\rm H}(v)\cdot \rho_v(-1).
\]
A configuration in ${\rm Conf}_u^{v}({\cal A})$ is illustrated on Figure \ref{partial.confbruhat}: 
\begin{figure}[ht]
\begin{center}
\epsfxsize100pt
\begin{tikzpicture}[scale=0.6]
\node [circle,draw=red,fill=red,minimum size=3pt,inner sep=0pt,label=below:{\small ${\A}_l$}]  (a) at (-1,-3){};
\node [circle,draw=red,fill=red,minimum size=3pt,inner sep=0pt,label=above:{\small ${\A}^l$}] (b) at (-1,0) {};
 \node [circle,draw=red,fill=red,minimum size=3pt,inner sep=0pt,label=above:{\small ${\A}^r$}] (c) at (0.8,0) {};
  \node [circle,draw=red,fill=red,minimum size=3pt,inner sep=0pt,label=below:{\small $~~{\A}_r$}] (d) at (0.8,-3) {};
 \foreach \from/\to in {b/a, c/d}
                   \draw[thick] (\from) -- (\to);
   \draw[draw, directed, thick] (b) -- (c);
\draw[draw, directed, thick] (a) -- (d);    
 \node[blue] at (-0.1,0.33) {$v$};
 \node[blue] at (-0.1,-2.5) {$u^\ast$};
  \node[blue] at (-1.55,-1.5) {$w_0$};
 \node[blue] at (1.35,-1.5) {$w_0$};
 \end{tikzpicture}               
\caption{The partial configuration space ${\rm Conf}_u^v({\cal A})$.}
\label{partial.confbruhat}
\end{center}
\end{figure}

\ed

\bt \la{Th2.3}
\la{main.result1.a.structure}
\begin{enumerate}
\item Every reduced word ${\bf i}$ of $(u, v)$ gives rise to  a (weighted) quiver
\[
 {\bf J}({\bf i})=\left\{  {J}({\bf i})\supset  {J}_{\rm uf}({\bf i}),  {\varepsilon}({\bf i})=(\varepsilon_{ij}), d({\bf i})=\{d_i\}\right\}.
\]
Here ${J}({\bf i})$ is a set of vertices, ${J}_{\rm uf}({\bf i})$ is a subset of unfrozen vertices,  $ \varepsilon({\bf i})$ is the exchange matrix, and $d({\bf i})$ is the set of multipliers such that $\widehat{\varepsilon}_{ij}=\varepsilon_{ij}d_j^{-1}$ is skewsymmerric.  Quivers defined using different reduced words are related by sequences of quiver mutations and permutations of vertices.

\vskip 1mm\item  Every ${\bf i}$  gives rise to regular maps
\[ {\rm T}^{\mathscr A}_{\bf i}:=( {\Bbb G}_m)^{\# {J}({\bf i})} ~ \xhookrightarrow{~~c_{{\mathscr A}, {\bf i}}~~}
 ~{\rm Conf}_{u}^{v}({\cal A}) \xrightarrow{~~\pi_{\bf i}~~} {\Bbb A}^{\# {J}({\bf i})}.
 \]
Here $\pi^{\bf i}$ is  presented by a coordinate system  $\{A_{{\bf i}, j} \}_{ j\in {J}({\bf i})}$, and 
 $c_{{\mathscr A}, {\bf i}}$ is an open embedding such that $\pi^{\bf i}\circ c_{{\mathscr A}, {\bf i}}$ is the standard injection. 
 Transition maps between  coordinate systems  using different reduced words are cluster $K_2-$transformations encoded by the same sequences of cluster mutations and permutations as in 1).

\vskip 1mm\item Every ${\bf i}$ gives rise to an open algebraic torus embedding
\[
 {\rm T}^{\mathscr X}_{\bf i}:=( {\Bbb G}_m)^{\# {J}_{\rm uf}({\bf i})}~ \xhookrightarrow{~~c_{{\mathscr X}, {\bf i}}~~} ~{\rm Conf}_u^{v}({\cal B}).
\]
The inverse map  $c_{\mathscr{X}, {\bf i}}^{-1}$ is presented by a birational  coordinate system $\left\{X_{{\bf i}, j} \right\}_{j\in {J}_{\rm uf}({\bf i}) }$.
Transition maps between coordinate systems using different reduced words are cluster Poisson transformations encoded by the the same sequence as in 1).

\vskip 1mm\item The following diagram commutes 
\be
\la{projection.2017.3.27.12.00}
\begin{gathered}
 \xymatrix{
        {\rm T}^{\mathscr A}_{\bf i} \ar[d]_{p_{\bf i}}\ar[r]^{c_{{\mathscr A},{\bf i}}~~~} &{\rm Conf}^{v}_u({\cal A})\ar[d]^{p} \\
        {\rm T}^{\mathscr X}_{\bf i} \ar[r]_{c_{{\mathscr X}, {\bf i}}~~~} &{\rm Conf}_u^{v}({\cal B})
 }
\end{gathered}
\ee
where  $p$ is the natural projection defined by forgetting decorations, and $p_{\bf i}$ is defined by 
\be
\la{2017.4.19.14.39tt}
p_{\bf i}^* X_{{\bf i}, j} = \prod_{k\in {J}({\bf i})} A_{{\bf i}, k}^{\varepsilon_{jk}}, \hskip 1cm \forall j \in {J}_{\rm uf}({\bf i}).
\ee
\end{enumerate}
\et
In Section \ref{partial.conf.part3}, we investigate a moduli space $\mathcal{P}_{u}^v$ which is naturally isomorphic to a double Bruhat cell in $\G$. Cluster Poisson structures of double Bruhat cells and their underlying seeds ${\bf J}({\bf i})$ were constructed in \cite{FG05}. Part 3 of Theorem \ref{Th2.3}  can be deduced from there as well. 
We prove Theorem \ref{main.result1.a.structure} in the rest of  Section \ref{pcs.sec}.

\subsection{Construction of underlying quivers} 
\la{sec.3.1} 
\medskip

The Weyl group $W\times W$ of $\G' \times \G'$ is generated by the simple reflections
\[
s_k, ~~~ k\in {\rm I}\cup{\overline{\rm I}}=\{1, \ldots, r\}\cup\{\overline{1}, \ldots, {\overline{r}}\},
\]
where the generators of the second $W$ are indexed by ${\overline{\rm  I}}=\{\overline{1}, \ldots, \overline{r}\}$. A reduced word ${\bf i}$ for $(u, v)\in W\times W$ will be represented by a sequence of indices $(i_1,\ldots, i_m)$ from ${\rm I} \cup \overline{\rm I}$.

Let $(\ast ,\ast)$ be the canonical pairing between the root lattice and coroot lattice of $\G'$. 
Recall the Cartan matrix  
$C_{ij}=(\alpha_i, \alpha_j^\vee)$
and the multipliers $d_j ={\langle\alpha_j^\vee, \alpha_j^\vee\rangle} \in \{1, 2,  3\}$ such that $\widehat{C}_{ij}=C_{ij}d_j^{-1}$ is symmetric.

\subsubsection{Elementary quivers.} \la{11.2.1}
 
The weighted quiver ${\bf J}(i)$ consists of an underlying set
\[
{J}(i):= \left({\rm I}-\{i\}\right) \cup \{i_l\}\cup \{i_r\} \cup\{i_e\}.
\]
There is a {\it decoration} map $\pi: {J}(i)\rightarrow {\rm I}$ which sends $i_l, ~ i_r$ and $i_e$ to $i$, and is the identity map on ${\rm I}-\{i\}$. The multipliers on ${J}(i)$ is defined by pulling back the multipliers on ${\rm I}$. 
The skew-symmetrizable matrix $\varepsilon(i)$ is indexed by ${J}(i)\times {J}(i)$,  defined as follows
\be
\label{exchange.matrix.ele.quiver}
\varepsilon(i)_{i_l,j}=\frac{- C_{ij}}{2}, ~~~\varepsilon(i)_{i_r,j}=\frac{C_{ij}}{2}, ~~~~\varepsilon(i)_{i_r,i_l}=\varepsilon(i)_{i_l,i_e} =\varepsilon(i)_{i_e, i_r}=1; \hskip 7mm \varepsilon(i)_{jk}=0 ~~\mbox{if } i\notin\{j,k\}.
\ee

\vskip 1mmThe quiver ${\bf J}(\overline{i})$ is obtained from ${\bf J}(i)$
by setting $\varepsilon(\overline{i})= - \varepsilon({i})$.

\vskip 1mmThe sub-quiver ${\bf J}_\circ(i)$ is obtained by deleting the extra element $i_e$ from ${\bf J}(i)$. 

\vskip 1mm The sub-quiver  ${\bf J}_\circ(\overline{i})$ is obtained by deleting the extra element $i_e$ from ${\bf J}(\overline{i})$. 

\vskip 1mm An elementary quiver can be presented by a directed graph with vertices labelled by the underlying set $J$ and arrows encoding the  exchange matrices $\varepsilon=(\varepsilon_{jk})$, where 
\[
\varepsilon_{jk}= \# \{\mbox{arrows from $j$ to $k$}\} - \# \{\mbox{arrows from $k$ to $j$}\}.
\] 
Note that  we employ here the notation $\# \{\mbox{arrows from $a$ to $b$}\}$ for the {\it total weight} of the arrows from $a$ to $b$, which is a half-integer.  
On the pictures the arrows from $a$ to $b$ are either dashed, and then counted with the weight $\frac{1}{2}$, 
or solid, and then counted with the weight $1$. 

For non simply laced cases, we will employ special arrows as shown on Example \ref{typeB3quiver}.

\begin{example} The quivers ${\bf J}_\circ(2), {\bf J}(3), {\bf J}(\overline{2})$ for the type $A_4$  are decribed by the following graphs. 
The green dot corresponds to $i_e$. Elements of ${J}(i)$ are placed at levels labelled by ${\rm I}$, with the extra $i_e$ placed at yet a different level. Elements at each level are ordered from the left to the right.
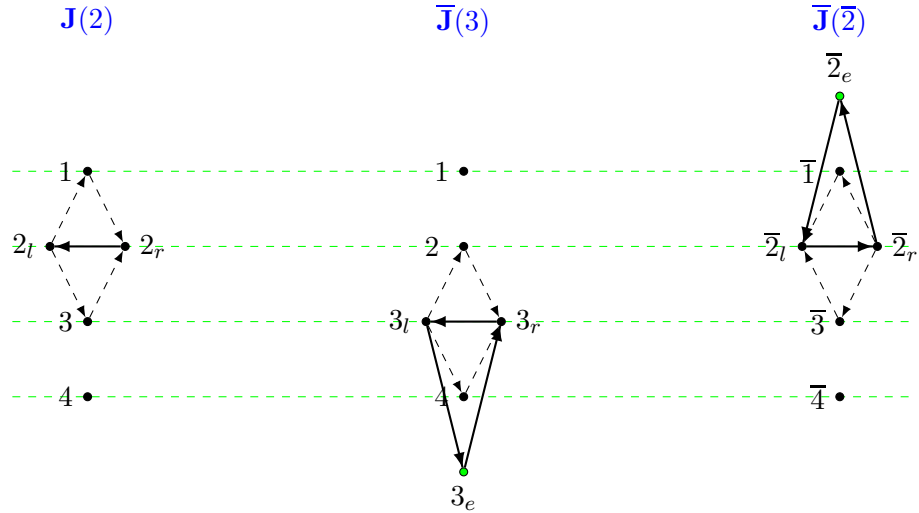
\begin{center} 
 \begin{tikzpicture}[scale=1.0]
                   \draw[green,dashed] (-6,0) -- (6,0);     
  \draw[green,dashed] (-6,-1) -- (6,-1);
   \draw[green,dashed] (-6,-2) -- (6,-2);
    \draw[green,dashed] (-6,-3) -- (6,-3);
    
 \node [blue] at (-5,2) {${\bf J}_\circ(2)$};
 \node [blue] at (0,2) {${\bf J}(3)$};
  \node [blue] at (5,2) {${\bf J}(\overline{2})$};
\node [circle,draw,fill,minimum size=3pt,inner sep=0pt,label=left:$1$] (f) at (0,0) {};
 \node [circle,draw,fill,minimum size=3pt,inner sep=0pt,label=left:$3_l$]  (b) at (-0.5,-2){};
 \node [circle,draw,fill,minimum size=3pt,inner sep=0pt,label=right:$3_r$] (c) at (0.5,-2) {};
 \node [circle,draw,fill,minimum size=3pt,inner sep=0pt,label=left:$2~$] (d) at (0,-1) {};
  \node [circle,draw=black,fill=green,minimum size=3pt,inner sep=0pt,label=below:$3_e$] (e) at (0,-4) {};
  \node [circle,draw,fill,minimum size=3pt,inner sep=0pt,label=left:$4$] (a) at (0,-3) {};
\foreach \from/\to in {c/b, b/e,e/c}
                   \draw[directed, thick] (\from) -- (\to);
 \foreach \from/\to in {b/a, a/c,b/d,d/c}
                   \draw[directed, dashed] (\from) -- (\to);

\node [circle,draw,fill,minimum size=3pt,inner sep=0pt,label=left:$1$] (a) at (-5,0) {};
 \node [circle,draw,fill,minimum size=3pt,inner sep=0pt,label=left:$2_l$]  (b) at (-5.5,-1){};
 \node [circle,draw,fill,minimum size=3pt,inner sep=0pt,label=right:$2_r$] (c) at (-4.5,-1) {};
 \node [circle,draw,fill,minimum size=3pt,inner sep=0pt,label=left:$3$] (d) at (-5,-2) {};
  \node at (-5,-4.5) {};
   \node [circle,draw,fill,minimum size=3pt,inner sep=0pt,label=left:$4$] (f) at (-5,-3) {};
 \foreach \from/\to in {b/a, a/c,b/d,d/c}
                   \draw[directed, dashed] (\from) -- (\to);
\foreach \from/\to in {c/b}
                   \draw[directed, thick] (\from) -- (\to);
 
\node [circle,draw,fill,minimum size=3pt,inner sep=0pt,label=left:$\overline{1}~$] (A) at (5,0) {};
 \node [circle,draw,fill,minimum size=3pt,inner sep=0pt,label=left:$\overline{2}_l$]  (B) at (4.5,-1){};
 \node [circle,draw,fill,minimum size=3pt,inner sep=0pt,label=right:$\overline{2}_r$] (C) at (5.5,-1) {};
 \node [circle,draw,fill,minimum size=3pt,inner sep=0pt,label=left:$\overline{3}$] (D) at (5,-2) {};
  \node [circle,draw=black,fill=green,minimum size=3pt,inner sep=0pt,label=above:$\overline{2}_e$] (E) at (5,1) {};
    \node [circle,draw,fill,minimum size=3pt,inner sep=0pt,label=left:$\overline{4}$] (F) at (5,-3) {};
\foreach \from/\to in {B/C, E/B,C/E}
                   \draw[directed, thick] (\from) -- (\to);
 \foreach \from/\to in {A/B, C/A,D/B,C/D}
                   \draw[directed, dashed] (\from) -- (\to);              
 \end{tikzpicture}
                   
\end{center}
\end{example}

\vskip 5mm

\begin{example}
\la{typeB3quiver}
The quivers ${\bf J}(1), {\bf J}(2)$ for type $B_3$ is described by the following graph:

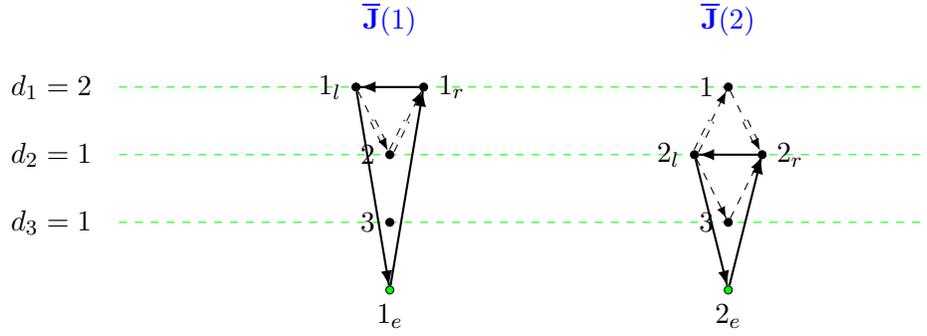
\begin{center}
\begin{tikzpicture}[scale=0.9]   
\node at (-7,0) {$d_1=2$};
\node at (-7,-1) {$d_2=1$};
\node at (-7,-2) {$d_3=1$};
         \draw[green,dashed] (-6,0) -- (6,0);     
  \draw[green,dashed] (-6,-1) -- (6,-1);
   \draw[green,dashed] (-6,-2) -- (6,-2);    
 \node [blue] at (-2,1) {${\bf J}(1)$};        
  \node [blue] at (3,1) {${\bf J}(2)$};      
\node [circle,draw,fill,minimum size=3pt,inner sep=0pt, label=left:${1}_l$] (a0) at (-2.5,0) {};
\node [circle,draw,fill,minimum size=3pt,inner sep=0pt,  label=right:${1}_r$]  (b0) at (-1.5,0){};
\node [circle,draw,fill,minimum size=3pt,inner sep=0pt,  label=left:$2$] (c0) at (-2,-1) {};
\node [circle,draw,fill,minimum size=3pt,inner sep=0pt,  label=left:$3$] (d0) at (-2,-2) {};
 \node [circle,draw=black,fill=green,minimum size=3pt,inner sep=0pt, label=below: $1_e$] (e0) at (-2,-3) {};
\foreach \from/\to in {b0/a0, a0/e0, e0/b0}
                   \draw[directed, thick] (\from) -- (\to);
                   \draw[--z>] (a0) to (c0);
                   \draw[z-->] (c0) to (b0);     
                   
 \node [circle,draw,fill,minimum size=3pt,inner sep=0pt, label=left:${2}_l$] (a1) at (2.5,-1) {};
\node [circle,draw,fill,minimum size=3pt,inner sep=0pt,  label=right:${2}_r$]  (b1) at (3.5,-1){};
\node [circle,draw,fill,minimum size=3pt,inner sep=0pt,  label=left:$1$] (c1) at (3,0) {};
\node [circle,draw,fill,minimum size=3pt,inner sep=0pt,  label=left:$3$] (d1) at (3,-2) {};
 \node [circle,draw=black,fill=green,minimum size=3pt,inner sep=0pt, label=below: $2_e$] (e1) at (3,-3) {};
\foreach \from/\to in {b1/a1, a1/e1, e1/b1}
                   \draw[directed, thick] (\from) -- (\to);
                   \draw[--z>] (c1) to (b1);
                   \draw[z-->] (a1) to (c1);     
                       \draw[directed, dashed] (a1) to (d1);
                   \draw[directed, dashed] (d1) to (b1);                              
               
 \end{tikzpicture}
 \end{center}

\end{example}


\subsubsection{The quivers ${\bf K}({\bf i})$ and ${\bf H}({\bf i})$.} 
Let $(u,v)\in W\times W$. For each reduced word ${\bf i}=(i_1, \ldots, i_m)$ of $(u,v)$. 
we get a chain of  distinct positive roots and coroots of $\G' \times \G'$
\be
\la{sequence.beta.906}
\alpha_k^{\bf i}:=s_{i_m}\ldots s_{i_{k+1}}\cdot \alpha_{i_{k}},\hskip 7mm
\beta_k^{\bf i}:=s_{i_m}\ldots s_{i_{k+1}}\cdot \alpha_{i_{k}}^\vee, \hskip 7mm k\in\{1,..., m\}.
\ee
The quiver ${\bf K}({\bf i})$ consists of $m$ many frozen vertices, whose exchange matrix is 
\[
\varepsilon_{jk}= \left\{ \begin{array}{ll} \frac{{\rm sgn}(k-j)}{2} {(\alpha_j^{\bf i}, \beta_k^{\bf i})} &\mbox{if } i_j, i_k\in {\rm I}\\
\frac{{\rm sgn}(j-k)}{2} {(\alpha_j^{\bf i}, \beta_k^{\bf i})}  &\mbox{if } i_j, i_k\in \overline{\rm I}\\
0, &\mbox{otherwise}\\
\end{array}
\right.
\]
The multipliers of the $k$th vertex are $d_k :={\langle \alpha_{i_k}^\vee, \alpha_{i_k}^\vee \rangle}$.

The quiver ${\bf H}({\bf i})$ is a full subquiver of ${\bf K}({\bf i})$ consisting of vertices $k$ such that $\beta_k^{\bf i}$ are simple. 
Let 
\[{\rm I}(u, v)={\rm I}(u)\cup {\rm I}(v) = \{j_1,..., j_s\}\cup \{l_1, \ldots, l_t\}. 
\] 
be the set  of elements $j\in {\rm I}$ and $l\in \overline{{\rm I}}$ such that the coroots $u\cdot \alpha_j^\vee$ and $v\cdot \alpha_{l}^\vee$ are negative.  
They are precisely the simple coroots appearing in the chain $\beta_k^{\bf i}$.
Therefore, the vertices of ${\bf H}({\bf i})$ are naturally parametrized by $I(u,v)$.

\subsubsection{Amalgamation.}
Recall the amalgamation of quivers introduced in \cite[Section 2.2]{FG05}.
\bd
\la{quiver.amla.fg4.sec22}
The quiver ${\bf Q}({\bf i})$ is the amalgamation of quivers ${\bf J}(i)$ and ${\bf K}({\bf i})$ determined by the reduced word ${\bf i}$:
\[
{\bf Q}({\bf i}):= {\bf J}(i_1)\ast\ldots\ast {\bf J}(i_m)\ast {\bf K}({\bf i})\]
Precisely, for every $i\in {\rm I}$ and for every $j=1,\ldots, m-1$, the right element of ${\bf J}(i_j)$  at level $i$ is glued with the  left element of ${\bf J}(i_{j+1})$ at level $i$.
The extra vertex of each ${\bf J}(i_k)$ is glued with the $k$th vertex of ${\bf K}({\bf i})$. 
\vskip 1mm

The quiver ${\bf J}({\bf i})$ is a full subquiver of ${\bf Q}({\bf i})$ by deleting all the ``non-simple" vertices of ${\bf K}({\bf i})$.
\vskip 1mm
The quiver ${\bf J}_\circ({\bf i})$ is a full subquiver of ${\bf J}({\bf i})$ by deleting all the vertices of ${\bf H}({\bf i})$.

\vskip 1mm

The quiver ${\bf J}_{\rm uf}({\bf i})$ is a subquiver of ${\bf J}({\bf i})$ by deleting the leftmost and rightmost vetices at every level $i \in {\rm I}$. It is the unfrozen part of ${\bf Q}({\bf i})$, ${\bf J}({\bf i})$, and ${\bf J}_\circ({\bf i})$.
\ed

\begin{example} 
Consider the root system $A_4$. Let ${\bf i}=(2, \overline{2}, \overline{3}, 3, \overline{2})$. The sequence $\beta_k^{\bf i}$ is
\[
\alpha_2^\vee+\alpha_3^\vee, ~~~~~ \alpha_{\bar 2}^\vee,~~~~~ \alpha_{\bar 2}^\vee+\alpha_{\bar 3}^\vee,~~~~~\alpha_3^\vee, ~~~~~\alpha_{\bar 3}^\vee.
\]
The quiver ${\bf J}({\bf i})$ is illustrated by the following Figure, where the  hallow dots are frozen.
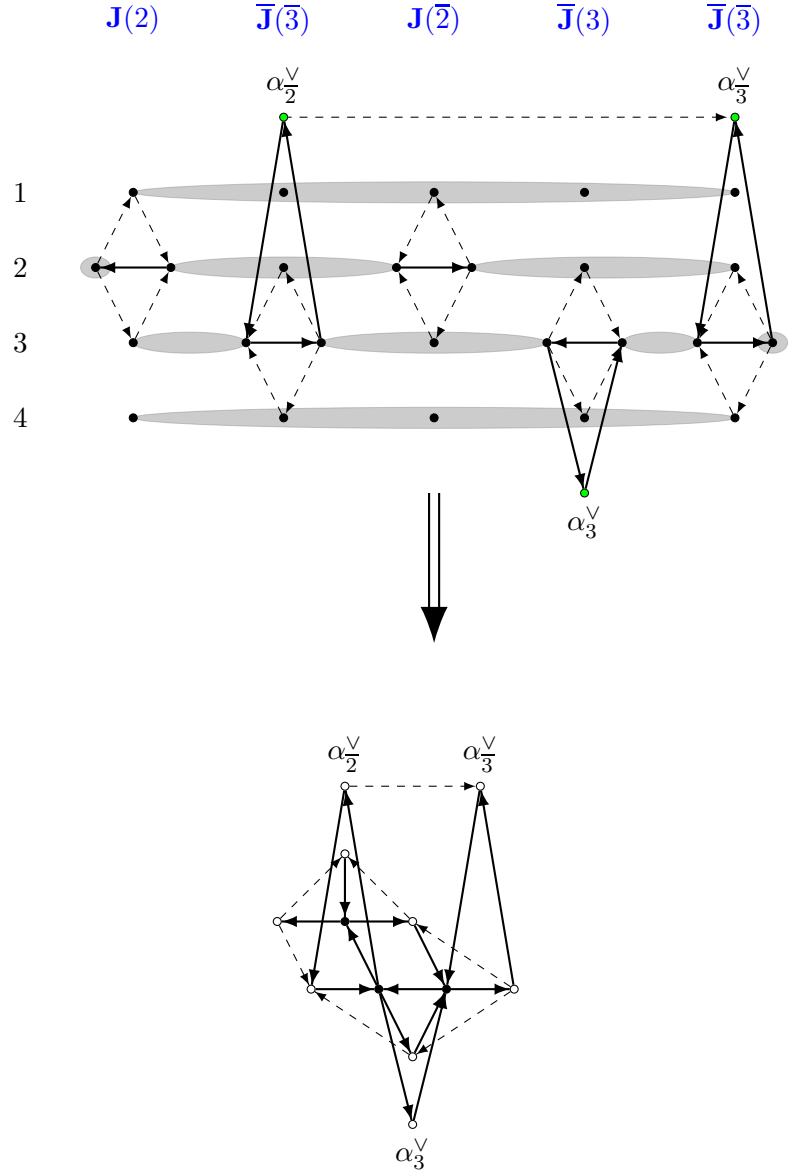
\begin{center}
\begin{tikzpicture}[scale=1]
 \draw[gray, fill=gray, opacity=0.4, very thin] (4,0) ellipse (40mm and 1.4mm);  
  \draw[gray, fill=gray, opacity=0.4,  very thin] (-0.5,-1) ellipse (2mm and 1.4mm); 
 \draw[gray, fill=gray, opacity=0.4, very thin] (2,-1) ellipse (15mm and 1.4mm);   
 \draw[gray, fill=gray, opacity=0.4,  very thin] (6.25,-1) ellipse (17.5mm and 1.4mm); 
 \draw[gray, fill=gray, opacity=0.4,  very thin] (0.75,-2) ellipse (7.5mm and 1.4mm);  
 \draw[gray, fill=gray, opacity=0.4,  very thin] (4,-2) ellipse (15mm and 1.4mm);    
   \draw[gray, fill=gray, opacity=0.4,  very thin] (8.5,-2) ellipse (2mm and 1.4mm);  
 \draw[gray, fill=gray, opacity=0.4,  very thin] (7,-2) ellipse (5mm and 1.4mm);                 
 \draw[gray, fill=gray, opacity=0.4, very thin] (4,-3) ellipse (40mm and 1.4mm);
  \draw[dashed,directed] (2, 1) -- (7.9, 1); 
\node[blue] at (0,2.3) {${\bf J}_\circ(2)$};
\node[blue] at (2,2.3) {${\bf J}({\overline{3}})$};
\node[blue] at (4,2.3) {${\bf J}({\overline{2}})$};
\node[blue] at (6,2.3) {${\bf J}({{3}})$};
\node[blue] at (8,2.3) {${\bf J}({\overline{3}})$};
\node at (-1.5, 0) {$1$};
\node at (-1.5, -1) {$2$};
\node at (-1.5, -2) {$3$};
\node at (-1.5, -3) {$4$};
%
\node [circle,draw,fill,minimum size=3pt,inner sep=0pt] (a) at (0,0) {};
 \node [circle,draw,fill,minimum size=3pt,inner sep=0pt]  (b) at (-0.5,-1){};
 \node [circle,draw,fill,minimum size=3pt,inner sep=0pt] (c) at (0.5,-1) {};
 \node [circle,draw,fill,minimum size=3pt,inner sep=0pt] (d) at (0,-2) {};
  \node [circle,draw,fill,minimum size=3pt,inner sep=0pt] (f) at (0,-3) {};
\foreach \from/\to in {c/b}
                   \draw[directed, thick] (\from) -- (\to);
 \foreach \from/\to in {b/a, a/c,b/d,d/c}
                   \draw[directed,  dashed] (\from) -- (\to);
\node [circle,draw,fill,minimum size=3pt,inner sep=0pt] (A) at (2,-3) {};
 \node [circle,draw,fill,minimum size=3pt,inner sep=0pt]  (B2) at (1.5,-2){};
 \node [circle,draw,fill,minimum size=3pt,inner sep=0pt] (B1) at (2.5,-2) {};
 \node [circle,draw,fill,minimum size=3pt,inner sep=0pt] (C) at (2,-1) {};
    \node [circle,draw,fill,minimum size=3pt,inner sep=0pt] (D) at (2,0) {};
 \node [circle,draw,fill=green,minimum size=3pt,inner sep=0pt, label=above: $\alpha_{\overline{2}}^\vee$] (E) at (2,1) {};
\foreach \to/\from in {B1/B2, B2/E, E/B1}
                   \draw[directed, thick] (\from) -- (\to);
 \foreach \to/\from in {A/B1, C/B1,B2/A,B2/C}
                   \draw[directed, dashed] (\from) -- (\to);
\node [circle,draw,fill,minimum size=3pt,inner sep=0pt] (A) at (4,-3) {};
 \node [circle,draw,fill,minimum size=3pt,inner sep=0pt]  (B) at (4,-2){};
 \node [circle,draw,fill,minimum size=3pt,inner sep=0pt] (C1) at (3.5,-1) {};
 \node [circle,draw,fill,minimum size=3pt,inner sep=0pt] (C2) at (4.5,-1) {};
    \node [circle,draw,fill,minimum size=3pt,inner sep=0pt] (D) at (4,0) {};
\foreach \from/\to in {C1/C2}
                   \draw[directed, thick] (\from) -- (\to);
 \foreach \from/\to in {B/C1, D/C1,C2/B,C2/D}
                   \draw[directed,  dashed] (\from) -- (\to);          
\node [circle,draw,fill,minimum size=3pt,inner sep=0pt] (A) at (6,0) {};
\node [circle,draw,fill,minimum size=3pt,inner sep=0pt]  (B) at (6,-1){};
\node [circle,draw,fill,minimum size=3pt,inner sep=0pt] (C1) at (5.5,-2) {};
\node [circle,draw,fill,minimum size=3pt,inner sep=0pt] (C2) at (6.5,-2) {};
\node [circle,draw,fill,minimum size=3pt,inner sep=0pt] (D) at (6,-3) {};
\node [circle,draw,fill=green,minimum size=3pt,inner sep=0pt, label=below: $\alpha_3^\vee$] (E) at (6,-4) {};
\foreach \from/\to in {C2/C1, C1/E, E/C2}
                   \draw[directed, thick] (\from) -- (\to);
 \foreach \from/\to in {B/C2, D/C2,C1/B,C1/D}
                   \draw[directed,  dashed] (\from) -- (\to);          
\node [circle,draw,fill,minimum size=3pt,inner sep=0pt] (A) at (8,-3) {};
\node [circle,draw,fill,minimum size=3pt,inner sep=0pt]  (B1) at (7.5,-2){};
\node [circle,draw,fill,minimum size=3pt,inner sep=0pt] (B2) at (8.5,-2) {};
\node [circle,draw,fill,minimum size=3pt,inner sep=0pt] (C) at (8,-1) {};
\node [circle,draw,fill,minimum size=3pt,inner sep=0pt] (D) at (8,0) {};
\node [circle,draw,fill=green,minimum size=3pt,inner sep=0pt, label=above: $\alpha_{\overline{3}}^\vee$] (E) at (8,1) {};
\foreach \to/\from in {B2/B1, B1/E, E/B2}
                   \draw[directed, thick] (\from) -- (\to);
 \foreach \to/\from in {C/B2, A/B2,B1/A,B1/C}
                   \draw[directed, dashed] (\from) -- (\to);          
%
 \draw[thick,double distance=3pt, arrows = {- Latex[length=0pt 3 0]}] (4,-4) -- (4,-6);     
  \end{tikzpicture}
  \vskip 1cm
   \begin{tikzpicture}[scale=0.9]    
\node [circle,draw,minimum size=3pt,inner sep=0pt] (A1) at (16.5,0) {};
 \node [circle,draw,minimum size=3pt,inner sep=0pt]  (B1) at (15.5,-1){};
 \node [circle,draw,fill,minimum size=3pt,inner sep=0pt] (B2) at (16.5,-1) {};
 \node [circle,draw,minimum size=3pt,inner sep=0pt] (B3) at (17.5,-1) {};
 \node [circle,draw,minimum size=3pt,inner sep=0pt] (C1) at (16,-2) {};
 \node [circle,draw,fill,minimum size=3pt,inner sep=0pt] (C2) at (17,-2) {};
 \node [circle,draw,fill,minimum size=3pt,inner sep=0pt] (C3) at (18,-2) {};
  \node [circle,draw,minimum size=3pt,inner sep=0pt] (C4) at (19,-2) {};
  \node [circle,draw,minimum size=3pt,inner sep=0pt] (D1) at (17.5,-3) {};
   \node [circle,draw,minimum size=3pt,inner sep=0pt, label=below: $\alpha_3^\vee$] (E2) at (17.5,-4) {};
     \node [circle,draw,minimum size=3pt,inner sep=0pt, label=above: $\alpha_{\overline{2}}^\vee$] (F1) at (16.5,1) {};
  \node [circle,draw,minimum size=3pt,inner sep=0pt, label=above: $\alpha_{\overline{3}}^\vee$] (F2) at (18.5,1) {};
\foreach \from/\to in {B2/B1, B2/B3, C1/C2,C3/C2,C3/C4,C2/B2,B3/C3,A1/B2,C2/D1,D1/C3}
                 \draw[directed, thick] (\from) -- (\to);
\foreach \from/\to in {C2/F1,F1/C1, C4/F2, F2/C3,C2/E2,E2/C3}
                 \draw[directed, thick] (\from) -- (\to);                 
\foreach \from/\to in {B1/A1,B1/C1,D1/C1, B3/A1, C4/B3, C4/D1,F1/F2}
               \draw[directed, dashed] (\from) -- (\to);  
 \end{tikzpicture}
                   
\end{center}
\end{example}

\subsubsection{An alternative notation.} \la{10.2.4} Let us reserve an alternative notation for the parametrization of the  vertices of ${\bf J}({\bf i})$. For $i\in {\rm I}$, denote by $n_{i}({\bf i})$ the number of occurrences of $i$ and $\bar{i}$ in ${\bf i}$. The vertices of ${\bf J}({\bf i})$ on the level $i$ from the left to the right can be indexed by
\[
V_i({\bf i}):=\left \{{i \choose k} ~\middle |~ 0\leq k\leq n_i({\bf i})\right \}.
\]
The vertices of ${\bf H}({\bf i})$ are indexed by 
\[
J_\infty := \left\{{i \choose -\infty} ~\middle |~ i \in {\rm I}(u)\right\} \, \bigsqcup \, \left\{{i \choose +\infty} ~\middle |~ \overline{i} \in {\rm I}(v)\right \}.
\]
Therefore  the underlying set of ${\bf J}(i)$ is indexed by
\be
\la{alternative.notation.J0}
{J}({\bf i}) :=\left \{{i \choose k} ~\middle |~ i\in {\rm I}, ~~ 0\leq k\leq n_i({\bf i})\right \} \bigsqcup J_\infty,
\ee
where the non-frozen part is the subset
\be
\la{alternative.notation.J}
{J}_{\rm uf}({\bf i}):=\left \{{i \choose k} ~\middle |~ i\in {\rm I}, ~~ 0< k< n_i({\bf i})\right\}.
\ee
For $1\leq k \leq n_i({\bf i})$, we will identify the symbol ${i \choose k}$ with the  position of  the $k^{th}$ $i$ or $\overline{i}$ appearing in the word ${\bf i}$. For example, if ${\bf i}=(2, \overline{2}, \overline{3}, 3, \overline{2})$, then 
\be
\la{convention.notation}
{2 \choose 1}=1,\hskip 3mm {2 \choose 2}=2, \hskip 3mm {2 \choose 3}=5, \hskip 3mm {3 \choose 1}=3, \hskip 3mm {3 \choose 2}=4. 
\ee

\subsection{Proof of Theorem \ref{Th2.3}, part 1}
\la{thm.3.6.st}
\medskip

It suffices to check the following four basic cases (plus an easy lemma indicating possible locations of simple positive coroots).
\begin{itemize}
\item[] \hskip -1cm{\bf Case 1}: $({i}, {\overline{i}}) \sim  ({\overline{i}}, {i})$.
\end{itemize}

The ${\bf J}(i, \overline{i})$ has one unfrozen vertex at level $i$, as shown on the following figure.
\vskip 5mm
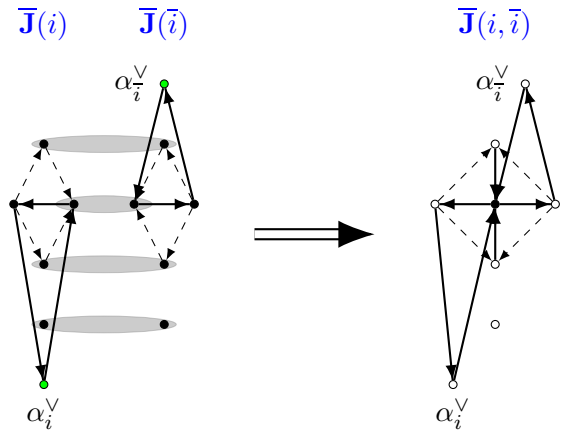
\begin{center}
\begin{tikzpicture}[scale=0.7]

 \draw[gray, fill=gray, opacity=0.4, very thin] (1,0) ellipse (12mm and 1.4mm);  
    \draw[gray, fill=gray, opacity=0.4, very thin] (1,-1) ellipse (8mm and 1.4mm);   
    \draw[gray, fill=gray, opacity=0.4,  very thin] (1,-2) ellipse (12mm and 1.4mm);                
    \draw[gray, fill=gray, opacity=0.4, very thin] (1,-3) ellipse (12mm and 1.4mm);

\node[blue] at (0,2) {${\bf J}(i)$};
\node[blue] at (2,2) {${\bf J}({\overline{i}})$};
\node[blue] at (7.5,2) {${{\bf J}}(i, \overline{i})$};

\node [circle,draw,fill,minimum size=3pt,inner sep=0pt] (a) at (0,0) {};
 \node [circle,draw,fill,minimum size=3pt,inner sep=0pt]  (b) at (-0.5,-1){};
 \node [circle,draw,fill,minimum size=3pt,inner sep=0pt] (c) at (0.5,-1) {};
 \node [circle,draw,fill,minimum size=3pt,inner sep=0pt] (d) at (0,-2) {};
  \node [circle,draw,fill=green,minimum size=3pt,inner sep=0pt, label=below: $\alpha_i^\vee$] (e) at (0,-4) {};
  \node [circle,draw,fill,minimum size=3pt,inner sep=0pt] (f) at (0,-3) {};
\foreach \from/\to in {c/b, b/e,e/c}
                   \draw[directed, thick] (\from) -- (\to);
 \foreach \from/\to in {b/a, a/c,b/d,d/c}
                   \draw[directed, dashed] (\from) -- (\to);
\node [circle,draw,fill,minimum size=3pt,inner sep=0pt] (A) at (2,0) {};
 \node [circle,draw,fill,minimum size=3pt,inner sep=0pt]  (B) at (1.5,-1){};
 \node [circle,draw,fill,minimum size=3pt,inner sep=0pt] (C) at (2.5,-1) {};
 \node [circle,draw,fill,minimum size=3pt,inner sep=0pt] (D) at (2,-2) {};
  \node [circle,draw,fill=green,minimum size=3pt,inner sep=0pt, label=left: $\alpha_{\overline{i}}^\vee$ ] (E) at (2,1) {};
    \node [circle,draw,fill,minimum size=3pt,inner sep=0pt] (F) at (2,-3) {};
\foreach \from/\to in {B/C, E/B,C/E}
                   \draw[directed, thick] (\from) -- (\to);
 \foreach \from/\to in {A/B, C/A,D/B,C/D}
                   \draw[directed, dashed] (\from) -- (\to);
    \draw[thick,double distance=3pt, arrows = {- Latex[length=0pt 3 0]}] (3.5,-1.5) -- (5.5,-1.5);       
    \node [circle,draw,minimum size=3pt,inner sep=0pt] (A1) at (7.5,0) {};
 \node [circle,draw,fill,minimum size=3pt,inner sep=0pt]  (B1) at (7.5,-1){};
 \node [circle,draw,minimum size=3pt,inner sep=0pt] (C1) at (8.5,-1) {};
  \node [circle,draw, minimum size=3pt,inner sep=0pt] (G1) at (6.5,-1) {};
 \node [circle,draw,minimum size=3pt,inner sep=0pt] (D1) at (7.5,-2) {};
  \node [circle,draw,minimum size=3pt,inner sep=0pt, label=left: $\alpha_{\overline{i}}^\vee$] (E1) at (8,1) {};
    \node [circle,draw,minimum size=3pt,inner sep=0pt] (F1) at (7.5,-3) {};
      \node [circle,draw,minimum size=3pt,inner sep=0pt, label=below: $\alpha_i^\vee$] (H1) at (6.8,-4) {};
\foreach \from/\to in {B1/C1, B1/G1, E1/B1,C1/E1, H1/B1, G1/H1}
                   \draw[directed, thick] (\from) -- (\to);
\foreach \from/\to in {A1/B1,D1/B1}
                   \draw[directed, thick] (\from) -- (\to);
 \foreach \from/\to in {C1/A1,G1/A1, G1/D1, C1/D1}
                   \draw[directed, dashed] (\from) -- (\to);  
 \end{tikzpicture}
 \end{center}
Mutating at this unfrozen vertex gives rise to the quiver ${\bf J}(\overline{i}, i)$, as shown below.
\vskip 5mm
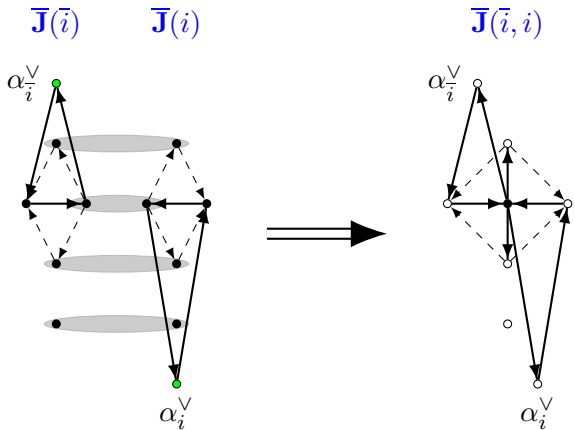
\begin{center}
\begin{tikzpicture}[scale=0.8]
 \draw[gray, fill=gray, opacity=0.4, very thin] (1,0) ellipse (12mm and 1.4mm);  
    \draw[gray, fill=gray, opacity=0.4, very thin] (1,-1) ellipse (8mm and 1.4mm);   
    \draw[gray, fill=gray, opacity=0.4,  very thin] (1,-2) ellipse (12mm and 1.4mm);                
    \draw[gray, fill=gray, opacity=0.4, very thin] (1,-3) ellipse (12mm and 1.4mm);
    \node[blue] at (2,2) {${\bf J}(i)$};
\node[blue] at (0,2) {${\bf J}({\overline{i}})$};
\node[blue] at (7.5,2) {${{\bf J}}(\overline{i},i)$};

\node [circle,draw,fill,minimum size=3pt,inner sep=0pt] (a) at (2,0) {};
 \node [circle,draw,fill,minimum size=3pt,inner sep=0pt]  (b) at (1.5,-1){};
 \node [circle,draw,fill,minimum size=3pt,inner sep=0pt] (c) at (2.5,-1) {};
 \node [circle,draw,fill,minimum size=3pt,inner sep=0pt] (d) at (2,-2) {};
  \node [circle,draw,fill=green,minimum size=3pt,inner sep=0pt, label=below: $\alpha_i^\vee$] (e) at (2,-4) {};
  \node [circle,draw,fill,minimum size=3pt,inner sep=0pt] (f) at (2,-3) {};
\foreach \from/\to in {c/b, b/e,e/c}
                   \draw[directed, thick] (\from) -- (\to);
 \foreach \from/\to in {b/a, a/c,b/d,d/c}
                   \draw[directed,  dashed] (\from) -- (\to);
\node [circle,draw,fill,minimum size=3pt,inner sep=0pt] (A) at (0,0) {};
 \node [circle,draw,fill,minimum size=3pt,inner sep=0pt]  (B) at (-0.5,-1){};
 \node [circle,draw,fill,minimum size=3pt,inner sep=0pt] (C) at (0.5,-1) {};
 \node [circle,draw,fill,minimum size=3pt,inner sep=0pt] (D) at (0,-2) {};
  \node [circle,draw,fill=green,minimum size=3pt,inner sep=0pt, label=left: $\alpha_{\overline{i}}^\vee$] (E) at (0,1) {};
    \node [circle,draw,fill,minimum size=3pt,inner sep=0pt] (F) at (0,-3) {};
\foreach \from/\to in {B/C, E/B,C/E}
                   \draw[directed, thick] (\from) -- (\to);
 \foreach \from/\to in {A/B, C/A,D/B,C/D}
                   \draw[directed,  dashed] (\from) -- (\to);
    \draw[thick,double distance=3pt, arrows = {- Latex[length=0pt 3 0]}] (3.5,-1.5) -- (5.5,-1.5);       
    \node [circle,draw,minimum size=3pt,inner sep=0pt] (A1) at (7.5,0) {};
 \node [circle,draw,fill,minimum size=3pt,inner sep=0pt]  (B1) at (7.5,-1){};
 \node [circle,draw,minimum size=3pt,inner sep=0pt] (C1) at (8.5,-1) {};
  \node [circle,draw,minimum size=3pt,inner sep=0pt] (G1) at (6.5,-1) {};
 \node [circle,draw,minimum size=3pt,inner sep=0pt] (D1) at (7.5,-2) {};
  \node [circle,draw,minimum size=3pt,inner sep=0pt, label=left: $\alpha_{\overline{i}}^\vee$] (E1) at (7,1) {};
    \node [circle,draw,minimum size=3pt,inner sep=0pt] (F1) at (7.5,-3) {};
      \node [circle,draw,minimum size=3pt,inner sep=0pt, label=below: $\alpha_i^\vee$] (H1) at (8,-4) {};
\foreach \to/\from in {B1/C1, B1/G1, E1/B1,G1/E1, H1/B1, C1/H1}
                   \draw[directed, thick] (\from) -- (\to);
\foreach \to/\from in {A1/B1,D1/B1}
                   \draw[directed, thick] (\from) -- (\to);
 \foreach \to/\from in {C1/A1,G1/A1, G1/D1, C1/D1}
                   \draw[directed, dashed] (\from) -- (\to);           
 \end{tikzpicture}
                   
\end{center}


\begin{itemize}
\item[] \hskip -1cm  {\bf Case 2:}  $(i, j, i)\sim (j, i, j)$, where $C_{ji}=C_{ij}=-1$. 
\end{itemize}

The quiver ${\bf J}(i, j, i)$ has an unfrozen vertex at level $i$, as shown below.
 \vskip 5mm
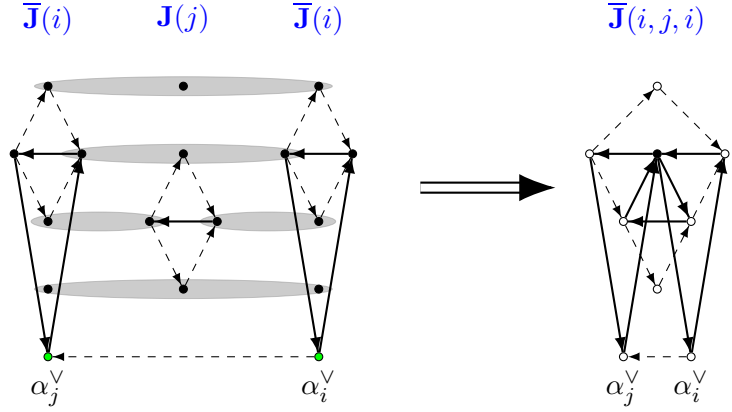
\begin{center}
\begin{tikzpicture}[scale=0.9]
 \draw[gray, fill=gray, opacity=0.4, very thin] (2,0) ellipse (22mm and 1.4mm);  
 \draw[gray, fill=gray, opacity=0.4, very thin] (2,-1) ellipse (18mm and 1.4mm);   
 \draw[gray, fill=gray, opacity=0.4,  very thin] (0.75,-2) ellipse (10mm and 1.4mm);  
 \draw[gray, fill=gray, opacity=0.4,  very thin] (3.25,-2) ellipse (10mm and 1.4mm);                
 \draw[gray, fill=gray, opacity=0.4, very thin] (2,-3) ellipse (22mm and 1.4mm);
  
\node[blue] at (0,1) {${\bf J}(i)$};
\node[blue] at (2,1) {${\bf J}_\circ(j)$};
\node[blue] at (4,1) {${\bf J}(i)$};
\node[blue] at (9,1) {${{\bf J}}(i,j,i)$};

\node [circle,draw,fill,minimum size=3pt,inner sep=0pt] (a) at (0,0) {};
\node [circle,draw,fill,minimum size=3pt,inner sep=0pt]  (b) at (-0.5,-1){};
\node [circle,draw,fill,minimum size=3pt,inner sep=0pt] (c) at (0.5,-1) {};
\node [circle,draw,fill,minimum size=3pt,inner sep=0pt] (d) at (0,-2) {};
\node [circle,draw,fill=green,minimum size=3pt,inner sep=0pt, label=below: $\alpha_j^\vee$] (e) at (0,-4) {};
\node [circle,draw,fill,minimum size=3pt,inner sep=0pt] (f) at (0,-3) {};
\foreach \from/\to in {c/b, b/e,e/c}
                   \draw[directed, thick] (\from) -- (\to);
 \foreach \from/\to in {b/a, a/c,b/d,d/c}
                   \draw[directed,  dashed] (\from) -- (\to);
\node [circle,draw,fill,minimum size=3pt,inner sep=0pt] (A) at (2,-3) {};
 \node [circle,draw,fill,minimum size=3pt,inner sep=0pt]  (C) at (1.5,-2){};
 \node [circle,draw,fill,minimum size=3pt,inner sep=0pt] (B) at (2.5,-2) {};
 \node [circle,draw,fill,minimum size=3pt,inner sep=0pt] (D) at (2,-1) {};
    \node [circle,draw,fill,minimum size=3pt,inner sep=0pt] (F) at (2,0) {};
\foreach \from/\to in {B/C}
                   \draw[directed, thick] (\from) -- (\to);
 \foreach \from/\to in {A/B, C/A,D/B,C/D}
                   \draw[directed,  dashed] (\from) -- (\to);
\node [circle,draw,fill,minimum size=3pt,inner sep=0pt] (a1) at (4,0) {};
\node [circle,draw,fill,minimum size=3pt,inner sep=0pt]  (b1) at (3.5,-1){};
\node [circle,draw,fill,minimum size=3pt,inner sep=0pt] (c1) at (4.5,-1) {};
\node [circle,draw,fill,minimum size=3pt,inner sep=0pt] (d1) at (4,-2) {};
\node [circle,draw,fill=green,minimum size=3pt,inner sep=0pt,label=below: $\alpha_i^\vee$] (e1) at (4,-4) {};
\node [circle,draw,fill,minimum size=3pt,inner sep=0pt] (f1) at (4,-3) {};
\foreach \from/\to in {c1/b1, b1/e1,e1/c1}
                   \draw[directed, thick] (\from) -- (\to);
 \foreach \from/\to in {b1/a1, a1/c1,b1/d1,d1/c1, e1/e}
                   \draw[directed,  dashed] (\from) -- (\to);                
    \draw[thick,double distance=3pt, arrows = {- Latex[length=0pt 3 0]}] (5.5,-1.5) -- (7.5,-1.5);       
    \node [circle,draw,minimum size=3pt,inner sep=0pt] (A1) at (9,0) {};
 \node [circle,draw,fill,minimum size=3pt,inner sep=0pt]  (B2) at (9,-1){};
 \node [circle,draw,minimum size=3pt,inner sep=0pt] (B3) at (10,-1) {};
  \node [circle,draw,minimum size=3pt,inner sep=0pt] (B1) at (8,-1) {};
 \node [circle,draw,minimum size=3pt,inner sep=0pt] (C1) at (8.5,-2) {};
 \node [circle,draw,minimum size=3pt,inner sep=0pt] (C2) at (9.5,-2) {};
  \node [circle,draw,minimum size=3pt,inner sep=0pt, label=below: $\alpha_i^\vee$] (E2) at (9.5,-4) {};
    \node [circle,draw,minimum size=3pt,inner sep=0pt] (D1) at (9,-3) {};
      \node [circle,draw,minimum size=3pt,inner sep=0pt, label=below: $\alpha_j^\vee$] (E1) at (8.5,-4) {};
\foreach \from/\to in {B3/B2, B2/B1, C2/C1,B1/E1,E1/B2, B2/E2, E2/B3}
                   \draw[directed, thick] (\from) -- (\to);
\foreach \from/\to in {C1/B2,B2/C2}
                   \draw[directed, thick] (\from) -- (\to);
 \foreach \from/\to in {B1/A1,A1/B3,B1/C1, C2/B3, C1/D1, D1/C2, E2/E1}
                   \draw[directed,  dashed] (\from) -- (\to);  
 \end{tikzpicture}
\end{center}
Mutating at this unfrozen vertex, and re-placing it at level $j$, we get the quiver ${\bf J}(j,i,j)$ as follows
 \vskip 5mm

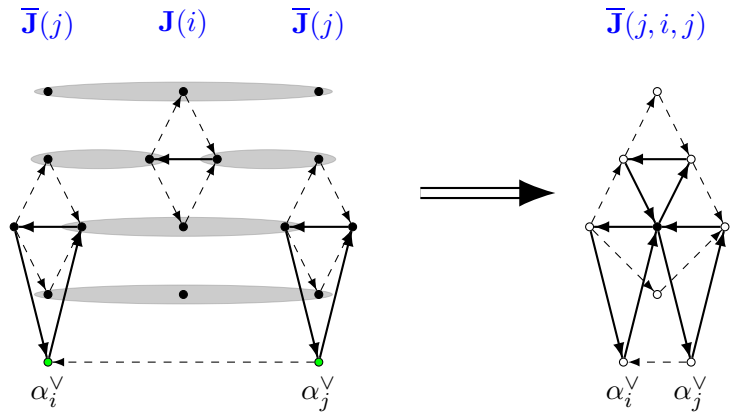
\begin{center}
\begin{tikzpicture}[scale=0.9]
 \draw[gray, fill=gray, opacity=0.4, very thin] (2,0) ellipse (22mm and 1.4mm);  
    \draw[gray, fill=gray, opacity=0.4, very thin] (0.75,-1) ellipse (10mm and 1.4mm); 
     \draw[gray, fill=gray, opacity=0.4, very thin] (3.25,-1) ellipse (10mm and 1.4mm);   
    \draw[gray, fill=gray, opacity=0.4,  very thin] (2,-2) ellipse (18mm and 1.4mm);                
    \draw[gray, fill=gray, opacity=0.4, very thin] (2,-3) ellipse (22mm and 1.4mm);

\node[blue] at (2,1) {${\bf J}_\circ(i)$};
\node[blue] at (0,1) {${{\bf J}}(j)$};
\node[blue] at (4,1) {${{\bf J}}(j)$};
\node[blue] at (9,1) {${{\bf J}}(j,i,j)$};

\node [circle,draw,fill,minimum size=3pt,inner sep=0pt] (a) at (2,0) {};
\node [circle,draw,fill,minimum size=3pt,inner sep=0pt]  (b) at (1.5,-1){};
\node [circle,draw,fill,minimum size=3pt,inner sep=0pt] (c) at (2.5,-1) {};
\node [circle,draw,fill,minimum size=3pt,inner sep=0pt] (d) at (2,-2) {};
\node [circle,draw,fill,minimum size=3pt,inner sep=0pt] (f) at (2,-3) {};
\foreach \from/\to in {c/b}
                   \draw[directed, thick] (\from) -- (\to);
 \foreach \from/\to in {b/a, a/c,b/d,d/c}
                   \draw[directed,  dashed] (\from) -- (\to);
\node [circle,draw,fill,minimum size=3pt,inner sep=0pt] (A) at (0,-3) {};
 \node [circle,draw,fill,minimum size=3pt,inner sep=0pt]  (B) at (0.5,-2){};
 \node [circle,draw,fill,minimum size=3pt,inner sep=0pt] (C) at (0-.5,-2) {};
 \node [circle,draw,fill,minimum size=3pt,inner sep=0pt] (D) at (0,-1) {};
  \node [circle,draw,fill=green,minimum size=3pt,inner sep=0pt, label=below: $\alpha_i^\vee$] (E) at (0,-4) {};
    \node [circle,draw,fill,minimum size=3pt,inner sep=0pt] (F) at (0,0) {};
\foreach \from/\to in {B/C, E/B,C/E}
                   \draw[directed, thick] (\from) -- (\to);
 \foreach \from/\to in {A/B, C/A,D/B,C/D}
                   \draw[directed,  dashed] (\from) -- (\to);
                   \node [circle,draw,fill,minimum size=3pt,inner sep=0pt] (A1) at (4,-3) {};
 \node [circle,draw,fill,minimum size=3pt,inner sep=0pt]  (B1) at (4.5,-2){};
 \node [circle,draw,fill,minimum size=3pt,inner sep=0pt] (C1) at (3.5,-2) {};
 \node [circle,draw,fill,minimum size=3pt,inner sep=0pt] (D1) at (4,-1) {};
  \node [circle,draw,fill=green,minimum size=3pt,inner sep=0pt, label=below: $\alpha_j^\vee$] (E1) at (4,-4) {};
    \node [circle,draw,fill,minimum size=3pt,inner sep=0pt] (F1) at (4,0) {};
\foreach \from/\to in {B1/C1, E1/B1,C1/E1}
                   \draw[directed, thick] (\from) -- (\to);
 \foreach \from/\to in {A1/B1, C1/A1,D1/B1,C1/D1, E1/E}
                   \draw[directed,  dashed] (\from) -- (\to);
    \draw[thick,double distance=3pt, arrows = {- Latex[length=0pt 3 0]}] (5.5,-1.5) -- (7.5,-1.5);       
    \node [circle,draw,minimum size=3pt,inner sep=0pt] (A1) at (9,0) {};
 \node [circle,draw,fill,minimum size=3pt,inner sep=0pt]  (C2) at (9,-2){};
 \node [circle,draw, minimum size=3pt,inner sep=0pt] (C3) at (10,-2) {};
  \node [circle,draw,minimum size=3pt,inner sep=0pt] (C1) at (8,-2) {};
 \node [circle,draw,minimum size=3pt,inner sep=0pt] (B1) at (8.5,-1) {};
  \node [circle,draw,minimum size=3pt,inner sep=0pt] (B2) at (9.5,-1) {};
  \node [circle,draw,minimum size=3pt,inner sep=0pt, label=below: $\alpha_i^\vee$] (D1) at (8.5,-4) {};
    \node [circle,draw,minimum size=3pt,inner sep=0pt] (E1) at (9,-3) {};
      \node [circle,draw,minimum size=3pt,inner sep=0pt, label=below: $\alpha_j^\vee$] (D2) at (9.5,-4) {};
\foreach \from/\to in {B2/B1, C2/C1, C3/C2, C1/D1, D1/C2, C2/D2, D2/C3}
                   \draw[directed, thick] (\from) -- (\to);
\foreach \from/\to in {B1/C2,C2/B2}
                   \draw[directed, thick] (\from) -- (\to);
 \foreach \from/\to in {C1/E1, E1/C3, C1/B1, B2/C3, B1/A1, A1/B2, D2/D1}
                   \draw[directed, dashed] (\from) -- (\to);           
 \end{tikzpicture}
                   
\end{center}

\begin{itemize}
\item[] \hskip -1cm  {\bf Case 3:}  $(i, j, i, j)\sim (j, i, j, i)$, where $C_{ij}=-1$ and $C_{ji}=-2$. 
\end{itemize}

The quiver ${\bf J}(i,j,i,j)$ has two unfrozen vertices, placed at the levels $i$ and  $j$ respectively. 
Below we show how to assemble the quiver from the elementary ones:

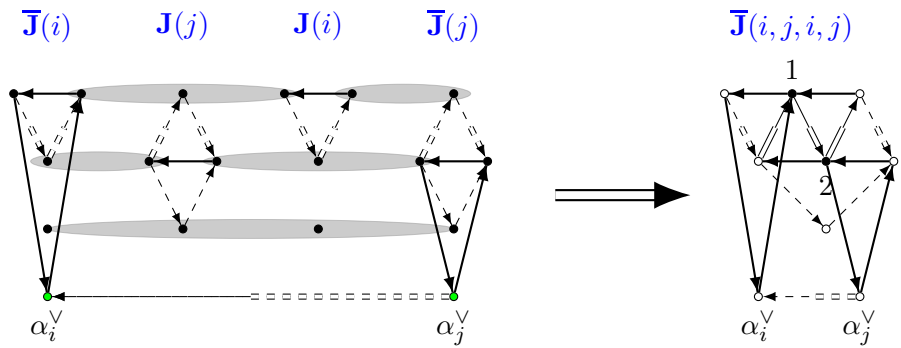
\begin{center}
\begin{tikzpicture}[scale=0.9]
 \draw[gray, fill=gray, opacity=0.4, very thin] (0,0) ellipse (17mm and 1.4mm);  
  \draw[gray, fill=gray, opacity=0.4, very thin] (3.25,0) ellipse (10mm and 1.4mm);  
    \draw[gray, fill=gray, opacity=0.4, very thin] (-1.25,-1) ellipse (10mm and 1.4mm); 
     \draw[gray, fill=gray, opacity=0.4, very thin] (2,-1) ellipse (17mm and 1.4mm);   
    \draw[gray, fill=gray, opacity=0.4,  very thin] (1,-2) ellipse (30mm and 1.4mm);                
    \node[blue] at (-2,1) {${\bf J}(i)$};
    \node[blue] at (2,1) {${\bf J}_\circ(i)$};
\node[blue] at (0,1) {${{\bf J}}_\circ(j)$};
\node[blue] at (4,1) {${{\bf J}}(j)$};
\node[blue] at (9,1) {${{\bf J}}(i,j,i,j)$};

\node [circle,draw,fill,minimum size=3pt,inner sep=0pt] (a0) at (-2.5,0) {};
\node [circle,draw,fill,minimum size=3pt,inner sep=0pt]  (b0) at (-1.5,0){};
\node [circle,draw,fill,minimum size=3pt,inner sep=0pt] (c0) at (-2,-1) {};
\node [circle,draw,fill,minimum size=3pt,inner sep=0pt] (d0) at (-2,-2) {};
 \node [circle,draw,fill=green,minimum size=3pt,inner sep=0pt, label=below: $\alpha_i^\vee$] (e0) at (-2,-3) {};
\foreach \from/\to in {b0/a0, a0/e0, e0/b0}
                   \draw[directed, thick] (\from) -- (\to);

\node [circle,draw,fill,minimum size=3pt,inner sep=0pt] (a) at (0,0) {};
 \node [circle,draw,fill,minimum size=3pt,inner sep=0pt] (b) at (-.5,-1){};
 \node [circle,draw,fill,minimum size=3pt,inner sep=0pt] (c) at (.5,-1) {};
 \node [circle,draw,fill,minimum size=3pt,inner sep=0pt] (d) at (0,-2) {};
\foreach \from/\to in {c/b}
                   \draw[directed, thick] (\from) -- (\to); 
 \foreach \from/\to in {b/d,d/c}
                   \draw[directed, dashed] (\from) -- (\to);

\node [circle,draw,fill,minimum size=3pt,inner sep=0pt] (a2) at (1.5,0) {};
\node [circle,draw,fill,minimum size=3pt,inner sep=0pt]  (b2) at (2.5,0){};
\node [circle,draw,fill,minimum size=3pt,inner sep=0pt] (c2) at (2,-1) {};
\node [circle,draw,fill,minimum size=3pt,inner sep=0pt] (d2) at (2,-2) {};
\foreach \from/\to in {b2/a2}
                   \draw[directed, thick] (\from) -- (\to);

 \node [circle,draw,fill,minimum size=3pt,inner sep=0pt] (A1) at (4,-2) {};
 \node [circle,draw,fill,minimum size=3pt,inner sep=0pt]  (B1) at (4.5,-1){};
 \node [circle,draw,fill,minimum size=3pt,inner sep=0pt] (C1) at (3.5,-1) {};
 \node [circle,draw,fill,minimum size=3pt,inner sep=0pt] (D1) at (4,0) {};
  \node [circle,draw,fill=green,minimum size=3pt,inner sep=0pt, label=below: $\alpha_j^\vee$] (E1) at (4,-3) {};
\foreach \from/\to in {B1/C1, E1/B1,C1/E1}
                   \draw[directed, thick] (\from) -- (\to);
 \foreach \from/\to in {A1/B1, C1/A1, E1/e0}
                   \draw[directed, dashed] (\from) -- (\to);                

 \foreach \from/\to in {a0/c0, a/c, a2/c2,D1/B1}
                   \draw[--z>] (\from) to (\to);
                
  \foreach \from/\to in {c0/b0, b/a, c2/b2,C1/D1, E1/e0}
                   \draw[z-->] (\from) to (\to);

\draw[thick,double distance=3pt, arrows = {- Latex[length=0pt 3 0]}] (5.5,-1.5) -- (7.5,-1.5);       

    \node [circle,draw,minimum size=3pt,inner sep=0pt] (A1) at (8,0) {};
     \node [circle,draw,fill,minimum size=3pt,inner sep=0pt, label=above: $1$] (A2) at (9,0) {};
      \node [circle,draw,minimum size=3pt,inner sep=0pt] (A3) at (10,0) {};
      
       \node [circle,draw,minimum size=3pt,inner sep=0pt] (B1) at (8.5,-1) {};
       \node [circle,draw,fill,minimum size=3pt,inner sep=0pt, label=below: $2$] (B2) at (9.5,-1) {};
       \node [circle,draw,minimum size=3pt,inner sep=0pt] (B3) at (10.5,-1) {};
 
  \node [circle,draw,minimum size=3pt,inner sep=0pt] (C1) at (9.5,-2) {};

  \node [circle,draw,minimum size=3pt,inner sep=0pt, label=below: $\alpha_i^\vee$] (D1) at (8.5,-3) {};
      \node [circle,draw,minimum size=3pt,inner sep=0pt, label=below: $\alpha_j^\vee$] (D2) at (10,-3) {};
\foreach \from/\to in {A3/A2, A2/A1, B3/B2, B2/B1}
                   \draw[directed, thick] (\from) -- (\to);
\foreach \from/\to in {A2/B2}
                   \draw[-z>] (\from) -- (\to);
\foreach \from/\to in {B1/A2,B2/A3}
                   \draw[z->] (\from) -- (\to);                   
\foreach \from/\to in {A1/B1,A3/B3}
                   \draw[--z>] (\from) -- (\to);           
               \draw[z-->] (D2)--(D1);            
 \foreach \from/\to in {B1/C1, C1/B3}
                  \draw[directed, dashed] (\from) -- (\to);           
 \foreach \from/\to in {A1/D1, D1/A2, B2/D2, D2/B3}
                  \draw[directed, thick] (\from) -- (\to);                            
 \end{tikzpicture}
 \end{center}                  
Let us apply the mutations $(2,1,2)$ to the above quiver. The following Figure demonstrates that the 
resulted quiver coincides with the quiver ${\bf J}(j,i,j,i)$:
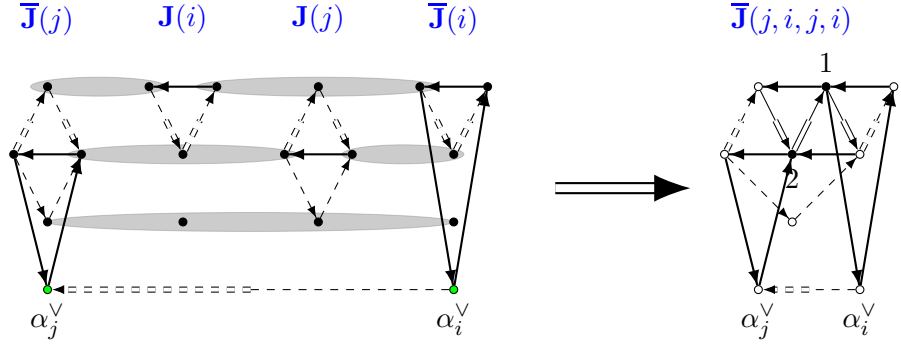
\begin{center}
\begin{tikzpicture}[scale=0.9]
 \draw[gray, fill=gray, opacity=0.4, very thin] (-1.25,0) ellipse (10mm and 1.4mm);  
 \draw[gray, fill=gray, opacity=0.4, very thin] (2,0) ellipse (18mm and 1.4mm);  
 \draw[gray, fill=gray, opacity=0.4, very thin] (0,-1) ellipse (17mm and 1.4mm); 
 \draw[gray, fill=gray, opacity=0.4, very thin] (3.25,-1) ellipse (9mm and 1.4mm);   
 \draw[gray, fill=gray, opacity=0.4,  very thin] (1,-2) ellipse (30mm and 1.4mm);                
    
\node[blue] at (-2,1) {${\bf J}(j)$};
\node[blue] at (2,1) {${\bf J}_\circ(j)$};
\node[blue] at (0,1) {${{\bf J}}_\circ(i)$};
\node[blue] at (4,1) {${{\bf J}}(i)$};
\node[blue] at (9,1) {${{\bf J}}(j,i,j,i)$};

\node [circle,draw,fill,minimum size=3pt,inner sep=0pt] (A1) at (-2,0) {};
\node [circle,draw,fill,minimum size=3pt,inner sep=0pt] (A2) at (-.5,0){};
\node [circle,draw,fill,minimum size=3pt,inner sep=0pt] (A3) at (.5,0) {};
\node [circle,draw,fill,minimum size=3pt,inner sep=0pt] (A4) at (2,0) {};
\node [circle,draw,fill,minimum size=3pt,inner sep=0pt] (A5) at (3.5,0) {};
\node [circle,draw,fill,minimum size=3pt,inner sep=0pt] (A6) at (4.5,0) {};

\node [circle,draw,fill,minimum size=3pt,inner sep=0pt] (B1) at (-2.5,-1) {};
\node [circle,draw,fill,minimum size=3pt,inner sep=0pt] (B2) at (-1.5,-1){};
\node [circle,draw,fill,minimum size=3pt,inner sep=0pt] (B3) at (0,-1) {};
\node [circle,draw,fill,minimum size=3pt,inner sep=0pt] (B4) at (1.5,-1) {};
\node [circle,draw,fill,minimum size=3pt,inner sep=0pt] (B5) at (2.5,-1) {};
\node [circle,draw,fill,minimum size=3pt,inner sep=0pt] (B6) at (4,-1) {};

\node [circle,draw,fill,minimum size=3pt,inner sep=0pt] (C1) at (-2,-2) {};
\node [circle,draw,fill,minimum size=3pt,inner sep=0pt] (C2) at (0,-2){};
\node [circle,draw,fill,minimum size=3pt,inner sep=0pt] (C3) at (2,-2) {};
\node [circle,draw,fill,minimum size=3pt,inner sep=0pt] (C4) at (4,-2) {};

\node [circle,draw,minimum size=3pt,inner sep=0pt] (A7) at (8.5,0) {};
\node [circle,draw,fill,minimum size=3pt,inner sep=0pt, label=above: $1$] (A8) at (9.5,0) {};
\node [circle,draw,minimum size=3pt,inner sep=0pt] (A9) at (10.5,0) {};
      
\node [circle,draw,minimum size=3pt,inner sep=0pt] (B7) at (8,-1) {};
\node [circle,draw,fill,minimum size=3pt,inner sep=0pt, label=below: $2$] (B8) at (9,-1) {};
\node [circle,draw,minimum size=3pt,inner sep=0pt] (B9) at (10,-1) {};

\node [circle,draw,minimum size=3pt,inner sep=0pt] (C5) at (9,-2) {}; 
\node [circle,draw,minimum size=3pt,inner sep=0pt, label=below: $\alpha_j^\vee$] (D3) at (8.5,-3) {}; 
\node [circle,draw,minimum size=3pt,inner sep=0pt, label=below: $\alpha_i^\vee$] (D4) at (10,-3) {}; 

\node [circle,draw,fill=green,minimum size=3pt,inner sep=0pt, label=below: $\alpha_j^\vee$] (D1) at (-2,-3) {};
\node [circle,draw,fill=green,minimum size=3pt,inner sep=0pt, label=below: $\alpha_i^\vee$] (D2) at (4,-3) {};

\foreach \from/\to in {A6/A5, A3/A2, B5/B4, B2/B1, B1/D1, D1/B2, A5/D2, D2/A6, A9/A8, A8/A7, B9/B8, B8/B7, B7/D3, D3/B8, A8/D4, D4/A9}
                  \draw[directed, thick] (\from) -- (\to); 
\foreach \from/\to in {B1/C1, C1/B2, B4/C3, C3/B5, B7/C5, C5/B9}
                  \draw[directed, dashed] (\from) -- (\to);                   
\foreach \from/\to in {A1/B2, A2/B3, A4/B5, A5/B6, D2/D1,D4/D3}
                  \draw[--z>] (\from) -- (\to); 
\foreach \from/\to in {B1/A1, B3/A3, B4/A4, B6/A6,B7/A7,B9/A9}
                  \draw[z-->] (\from) -- (\to);                   
\foreach \from/\to in {A7/B8, A8/B9}
                  \draw[-z>] (\from) -- (\to); 
\foreach \from/\to in {B8/A8}
                  \draw[z->] (\from) -- (\to);     
\draw[thick,double distance=3pt, arrows = {- Latex[length=0pt 3 0]}] (5.5,-1.5) -- (7.5,-1.5);                                                                     
\end{tikzpicture}
\end{center}

\begin{itemize}
\item[] \hskip -1cm  {\bf Case 4:}  $(i, j, i, j, i, j)\sim (j, i, j, i, j, i)$, where $C_{ij}=-1$ and $C_{ji}=-3$.
\end{itemize}

The quiver ${\bf J}(i,j,i,j,i,j)$ has four unfrozen vertices, placed at the levels $i$ and  $j$. 

\begin{center}
\begin{tikzpicture}[scale=1.0]         
    
\node[blue] at (1.5,1) {${{\bf J}}(i,j,i,j,i,j)$};

\node [circle,draw,minimum size=3pt,inner sep=0pt, label=above:$i_l$] (A1) at (0,0.1) {};
\node [circle,draw,fill,minimum size=3pt,inner sep=0pt, label=above:$1$] (A2) at (1,0.1) {};
\node [circle,draw,fill,minimum size=3pt,inner sep=0pt, label=above:$2$] (A3) at (2,0.1) {};
\node [circle,draw,minimum size=3pt,inner sep=0pt, label=above:$i_r$] (A4) at (3,0.1) {};
      
\node [circle,draw,minimum size=3pt,inner sep=0pt, label=below:$j_l$] (B1) at (0.5,-1) {};
\node [circle,draw,fill,minimum size=3pt,inner sep=0pt, label=below:$3$] (B2) at (1.5,-1) {};
\node [circle,draw,fill,minimum size=3pt,inner sep=0pt, label=below:$4$] (B3) at (2.5,-1) {};
\node [circle,draw,minimum size=3pt,inner sep=0pt, label=below:$j_r$] (B4) at (3.5,-1) {};

\node [circle,draw,minimum size=3pt,inner sep=0pt, label=below: $\alpha_i^\vee$] (D1) at (0.5,-2.1) {}; 
\node [circle,draw,minimum size=3pt,inner sep=0pt, label=below: $\alpha_j^\vee$] (D2) at (3,-2.1) {};

\foreach \from/\to in {A4/A3, A3/A2, A2/A1, B4/B3, B3/B2, B2/B1,A1/D1, D1/A2, B3/D2, D2/B4}
                  \draw[directed, thick] (\from) -- (\to);                 
\foreach \from/\to in {D2/D1}
                  \draw[x-->] (\from) -- (\to); 
\foreach \from/\to in {A1/B1,A4/B4}
                  \draw[--x>] (\from) -- (\to);                   
\foreach \from/\to in {A2/B2, A3/B3}
                  \draw[-x>] (\from) -- (\to); 
\foreach \from/\to in {B1/A2,B2/A3,B3/A4}
                  \draw[x->] (\from) -- (\to);                                                                     
\end{tikzpicture}
\end{center}
Following \cite[p.23]{FG05}, we apply the sequence of mutations $(4, 3, 2, 1, 4, 2, 4, 3, 1, 4)$ to the above quiver. A computer check shows that the resulted quiver coincides with the quiver ${\bf J}(j,i,j,i,j,i)$.

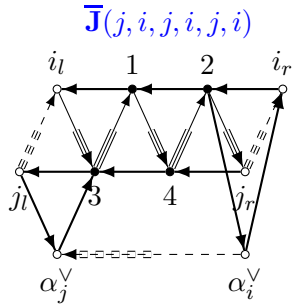
\begin{center}
\begin{tikzpicture}[scale=1.0]         
    
\node[blue] at (2.5,1) {${{\bf J}}(j,i,j,i,j,i)$};

\node [circle,draw,minimum size=3pt,inner sep=0pt, label=above:$i_l$] (A1) at (1,0.1) {};
\node [circle,draw,fill,minimum size=3pt,inner sep=0pt, label=above:$1$] (A2) at (2,0.1) {};
\node [circle,draw,fill, minimum size=3pt,inner sep=0pt, label=above:$2$] (A3) at (3,0.1) {};
\node [circle,draw,minimum size=3pt,inner sep=0pt, label=above:$i_r$] (A4) at (4,0.1) {};

\node [circle,draw,minimum size=3pt,inner sep=0pt, label=below:$j_l$] (B1) at (0.5,-1) {};
\node [circle,draw,fill,minimum size=3pt,inner sep=0pt, label=below:$3$] (B2) at (1.5,-1) {};
\node [circle,draw,fill,minimum size=3pt,inner sep=0pt, label=below:$4$] (B3) at (2.5,-1) {};
\node [circle,draw,minimum size=3pt,inner sep=0pt, label=below:$j_r$] (B4) at (3.5,-1) {};

\node [circle,draw,minimum size=3pt,inner sep=0pt, label=below: $\alpha_j^\vee$] (D1) at (1,-2.1) {}; 
\node [circle,draw,minimum size=3pt,inner sep=0pt, label=below: $\alpha_i^\vee$] (D2) at (3.5,-2.1) {};

\foreach \from/\to in {A4/A3, A3/A2, A2/A1, B4/B3, B3/B2, B2/B1,B1/D1, D1/B2, A3/D2, D2/A4}
                  \draw[directed, thick] (\from) -- (\to);                 
\foreach \from/\to in {D2/D1}
                  \draw[--x>] (\from) -- (\to); 
\foreach \from/\to in {B1/A1,B4/A4}
                  \draw[x-->] (\from) -- (\to);                   
\foreach \from/\to in {A1/B2, A2/B3, A3/B4}
                  \draw[-x>] (\from) -- (\to); 
\foreach \from/\to in {B2/A2,B3/A3}
                  \draw[x->] (\from) -- (\to);                                                                     
\end{tikzpicture}
\end{center}


\subsection{Triples of decorated flags and partial potentials}

\label{triple2020.1.28}
\medskip 
 
 Recall the configuration space ${\rm Conf}^*({\cal A}, {\cal B}, {\cal B})$ in (\ref{ABA}).  
 Let $w\in W$. Define the subspace ${\rm Conf}^w({\cal A}, {\cal B}, {\cal B})$ of ${\rm Conf}^*({\cal A}, {\cal B}, {\cal B})$ by requiring that the $w_-$-position of the  pair $(\B_1, \B_2)$ is  $w$. 
So 
\be
{\rm Conf}^*({\cal A}, {\cal B}, {\cal B}):= \bigsqcup_{w\in W} {\rm Conf}^w({\cal A}, {\cal B}, {\cal B}).
\ee
Since the space of frames is a $\G$-torsor, 
  every  $(\A,\B_1, \B_2)\in {\rm Conf}^*({\cal A}, {\cal B}, {\cal B})$ has a unique representative which identifies $\big(\A, \B_1\big)$ with the frame $([{\rm U}], \B^{-})$. There is a unique $u\in \U$ such that 
\[
(\A,\B_1, \B_2)=([\U], \B^-, u\cdot \B^-) \in {\rm Conf}^*({\cal A}, {\cal B}, {\cal B}).
\]
The element $u$ is called the {\it angle invariant} of $(\A, \B_1, \B_2)$: 
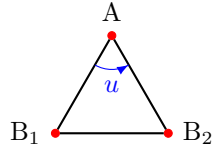
\begin{figure}[h]
\epsfxsize 100pt
\center{
\begin{tikzpicture}[scale=0.5]
\node [circle,draw=red,fill=red,minimum size=3pt,inner sep=0pt,label=above:{\small ${\A}$}] (a) at (0,0) {};
 \node [circle,draw=red,fill=red,minimum size=3pt,inner sep=0pt,label=left:{\small ${\B}_1$}]  (b) at (-120:3){};
 \node [circle,draw=red,fill=red,minimum size=3pt,inner sep=0pt,label=right:{\small ${\B}_2$}] (c) at (-60:3) {};
 \foreach \from/\to in {a/b, a/c, b/c}
                   \draw[thick] (\from) -- (\to);
  \draw[blue, directed] (-120:.9) arc (-120:-60:0.9);                  
  \node[blue] at (0,-1.36) {{\small $u$}};
 \end{tikzpicture}}
\label{angle.map.figure.22}
\end{figure}

In this way we get a natural isomorphism  
\be
\la{angle.map}
{\bf an}: {\rm Conf}^\ast ({\cal A}, {\cal B},{\cal B})  \stackrel{\sim}{\lra} \U , \hskip 7mm (\A, \B_1, \B_2) = ([\U], \B^-, u \cdot \B^-) \lms u.
\ee
 Let 
$
\U_w := \U \cap \B^- w \B^-
$. Since
$
w_-(\B^-, u \cdot \B^-) =w$ if and only if   $u \in \U_w$, we get  
\bl \la{angle.map.partial.corr.h}
Under the isomorphism  \eqref{angle.map}, we have
$
{\rm Conf}^w ({\cal A}, {\cal B},{\cal B})\stackrel{\sim}{=}\U_{w}. 
$
\el

\subsubsection{Partial potential of ${\rm Conf}^*({\cal A}, {\cal B}, {\cal B})$.} 
Fix an $i\in {\rm I}$. Let $(\A, \B_1, \B_2)\in {\rm Conf}^*({\cal A}, {\cal B}, {\cal B})$. As shown on Figure \ref{proj.map.figure.22.1}, there is a unique flag $\B_1'$ such that 
$w(\B_1', \A) = s_i$ and $w(\B_1, \B_1') =w_0 s_i$. Imposing the same condition to the pair $(\A, \B_2)$, we get another flag $\B_2'$. 

The triple $(\A, \B_1', \B_2')$ is canonically identified with a point in the  space $ {\rm Conf}^*({\cal A}, {\cal B}, {\cal B})_{{\rm SL}_2}$ for the group ${\rm SL}_2$. Therefore we get a projection
\be
p_i : {\rm Conf}^*({\cal A}, {\cal B}, {\cal B}) \lra  {\rm Conf}^*({\cal A}, {\cal B}, {\cal B})_{{\rm SL}_2}.
\ee
Composing $p_i$ with   potential   \eqref{sl2.potential.6.3}, we get a {\it partial potential} $
{\cal W}_i := {\cal W}\circ p_i 
$ on  ${\rm Conf}^*({\cal A}, {\cal B}, {\cal B}) $.

\begin{figure}[h]
\epsfxsize 100pt
\center{
\begin{tikzpicture}[scale=0.65]
\node [circle,draw=red,fill=red,minimum size=3pt,inner sep=0pt,label=above:{\small ${\A}$}] (a) at (0,0) {};
 \node [circle,draw=red,fill=red,minimum size=3pt,inner sep=0pt,label=left:{\small ${\B}_1$}]  (b) at (-1.5,-3){};
 \node [circle,draw=red,fill=red,minimum size=3pt,inner sep=0pt,label=right:{\small ${\B}_2$}] (c) at (1.5,-3) {};
  \node [circle,draw=red,fill=white,minimum size=3pt,inner sep=0pt,label=left:{\small ${\B}_1'$}]  (d) at (-0.5,-1){};
   \node [circle,draw=red,fill=white,minimum size=3pt,inner sep=0pt,label=right:{\small ${\B}_2'$}]  (e) at (0.5,-1){};
    \foreach \from/\to in {b/d, c/e}
                   \draw[thick, directed] (\from) -- (\to);
               \draw[thick] (b) -- (c);   
               \draw[thick] (d) -- (a);
               \draw[thick] (e) -- (a); 
               \node [blue] at (-0.45, -0.5) {{\small $s_i~~$}};
                \node [blue] at (0.45, -0.5) {{\small $~~s_i$}};
                 \node [blue] at (-1.45, -2) {{\small $w_0s_i~~$}};
                   \node [blue] at (1.45, -2) {{\small $~~w_0s_i$}};
          \draw[dotted, orange, thick] (0, -0.5) circle (1.2cm);         
 \end{tikzpicture}}
\caption{The projection map $p_i$.}
\label{proj.map.figure.22.1}
\end{figure}
\bl 
\la{6.30.17.45.hh}
The partial potential ${\cal W}_i$ is additive, i.e.,
\[
{\cal W}_i(\A, \B_1, \B_2)+ {\cal W}_i(\A, \B_2, \B_3) = {\cal W}_i(\A, \B_1, \B_3).
\]
\el
\begin{proof} It is a direct consequence of Lemma \ref{6.30.17.43.hh}.
\end{proof}

\bl
\la{6.30.17.33.hh}
Let $(\A, \B_1, \B_2)\in {\rm Conf}^{s_j}({\cal A}, {\cal B}, {\cal B})$. If $j\neq i$, then 
$
{\cal W}_i(\A, \B_1, \B_2) =0.
$
\el
\begin{proof}
Let $\B_1', \B_2'$ be defined as above.  By definition, it  suffices to show that $\B_1'=\B_2'$. 
Note that $w(\B_2, \B_1)= s_{j^\ast}$ and $w(\B_1, \B_1')= w_0 s_i$.  Since $j\neq i$, we have
$
l(s_{j^\ast} w_0 s_i) = l(w_0s_i) -1. 
$ 
So 
$w(\B_2, \B_1')$ is either $s_{j^*} w_0 s_i$ or $w_0 s_i$. 
However, if $w(\B_2, \B_1')=s_{j^*} w_0 s_i$, then $w(\B_2, \A)=s_{j^*} w_0$, which contradicts the assumption that $(\A, \B_2)$ is a frame. 
Thus 
$w(\B_2, \B_1')=w_0s_i$ and $w(\B_1', \A)=s_i$.  By uniqueness of the choices of $\B_1'$ and $\B_2'$, we get $\B_1'=\B_2'$. \end{proof}

\bp  The partial potential ${\cal W}_i$  is   the composition of the map \eqref{angle.map}
and  $\chi_i$ in \eqref{additive.isomorphism.5.21}:  
\[
{\cal W}_i = {\chi_i} \circ {\bf an}.
\]
\ep

\begin{proof}  
Let $(\A_0, \B, \B')\in {\rm Conf}^*({\cal A}, {\cal B}, {\cal B})$. There is a sequence
$
\B=\B_1,\B_2, \ldots, \B_m=\B'
$
such that $(\A_0, \B_k, \B_{k+1})\in {\rm Conf}^{s_{i_k}}({\cal A}, {\cal B}, {\cal B})$ for $k=1, ..., m-1$. 
By Lemma \ref{6.30.17.45.hh}, 
\[
{\cal W}_i(\A_0, \B, \B') =\sum_{k=1}^{m-1} {\cal W}_i(\A_0, \B_k, \B_{k+1}).
\]
Note that ${\bf an}(\A_0, \B, \B')=\prod_{k=1}^{m-1} {\bf an}(\A_0, \B_k, \B_{k+1})$. Therefore 
\be
\begin{split}
&{\chi_i} \big( {\bf an}(\A_0, \B, \B') \big)=
\sum_{k=1}^{m-1} \chi_i\big({\bf an}(\A_0, \B_k, \B_{k+1})\big).
\nonumber \\
\end{split}
\ee
So it suffices to show that for any $k$ we have 
 ${\cal W}_i(\A_0, \B_k, \B_{k+1})=\chi_i\big({\bf an}(\A_0, \B_k, \B_{k+1})\big)$. 
 
By Lemma \ref{angle.map.partial.corr.h}, we have
$
{\rm Conf}^{s_{i_k}}({\cal A}, {\cal B}, {\cal B}) \stackrel{\sim}{=} {\rm U}_{s_{i_k}}=\{x_{i_k}(a) ~|~ a\neq 0\}.
$

Thus $(\A_0,\B_k, \B_{k+1})=([\U], \B^-, x_{i_k}(a_k)\cdot \B^-)\in {\rm Conf}^{s_{i_k}}({\cal A}, {\cal B}, {\cal B}) $ for some $a_k\neq 0$.

If $i_k=i$, then 
\[
{\cal W}_i(\A_0, \B_{k}, \B_{k+1}) =a_k = \chi_i ({{\bf an}(\A_0, \B_{k}, \B_{k+1})}).
\]
Otherwise, by Lemma \ref{6.30.17.33.hh},  ${\cal W}_i(\A_0, \B_{k}, \B_{k+1}) =0= \chi_i \big({{\bf an}(\A_0, \B_{k}, \B_{k+1})}\big).$
\end{proof}

\subsection{Elementary configuration spaces}
\medskip

\subsubsection{The elementary configuration space ${\cal A}({{i}})$.} Let $i\in {\rm I}$. The  space ${\cal A}(i)$ is defined as the space  ${\rm Conf}_{s_i}^e({\cal A})$. It parametrizes $\G$-orbits of triples $(\A, \A_l, \A_r)$ such that
\be
w(\A, \A_l)= w(\A, \A_r)=w_0,\hskip 7mm w(\A_r, \A_l)=s_{i^\ast}, \hskip 4mm h(\A_r, \A_l) \in {\rm H}(s_{i^\ast}).
\ee
Imposing an extra condition $h(\A_r, \A_l)=1$, we get a subspace ${\cal A}_\circ(i)$ of ${\cal A}(i)$.
Below the  left figure describes a triple in  ${\cal A}(i)$, and the right figure describes a triple in  ${\cal A}_\circ(i)$.

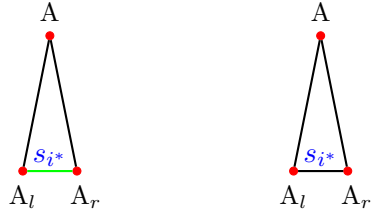
\begin{center}
\begin{tikzpicture}[scale=0.6]
\node [circle,draw=red,fill=red,minimum size=3pt,inner sep=0pt,label=above:{\small ${\A}$}] (a) at (0,0) {};
 \node [circle,draw=red,fill=red,minimum size=3pt,inner sep=0pt,label=below:{\small ${\A}_l$}]  (b) at (-0.6,-3){};
 \node [circle,draw=red,fill=red,minimum size=3pt,inner sep=0pt,label=below:{\small $~~{\A}_r$}] (c) at (0.6,-3) {};
 \foreach \from/\to in {b/a, c/a}
                   \draw[thick] (\from) -- (\to);
  \draw[green, thick] (c) -- (b);                 
 \node[blue] at (0,-2.7) {${ \small s_{i^\ast}}$};
 
 \node [circle,draw=red,fill=red,minimum size=3pt,inner sep=0pt,label=above:{\small ${\A}$}] (a1) at (6,0) {};
 \node [circle,draw=red,fill=red,minimum size=3pt,inner sep=0pt,label=below:{\small ${\A}_l$}]  (b1) at (5.4,-3){};
 \node [circle,draw=red,fill=red,minimum size=3pt,inner sep=0pt,label=below:{\small $~~{\A}_r$}] (c1) at (6.6,-3) {};
 \foreach \from/\to in {b1/a1, c1/a1}
                   \draw[thick] (\from) -- (\to);
  \draw[thick] (c1) -- (b1);                 
 \node[blue] at (6,-2.7) {${\small s_{i^\ast}}$};
 
  \end{tikzpicture}
 \end{center}

\bd
We define a set of functions on ${\cal A}(i)$ parametrized by ${J}(i)$:
\be
\la{ele.cor.cl.1}
\forall (\A, \A_l, \A_r)\in {\cal A}(i), \hskip 10mm
{\displaystyle{A_j := \left\{    \begin{array}{ll} 
     \Delta_j(\A, \A_l) \hskip 7mm & \mbox{if } j\in {\rm I}-\{i\}\\
       \Delta_{i}(\A, \A_l) & \mbox{if } j=i_l\\
       \Delta_{i}(\A, \A_r) & \mbox{if } j=i_r\\
      \Lambda_{i^\ast}(h(\A_r, \A_l)) & \mbox{if }j=i_e. \\
   \end{array}\right.}}
\ee
Recall the matrix $\varepsilon(i)$ of ${\bf J}(i)$ from \eqref{exchange.matrix.ele.quiver}. 
We define  a {\it canonical}   2-form of ${\cal A}(i)$ as follows
\be
{\Omega}(i) =  \sum_{j, k \in {J}(i)} \widetilde \varepsilon(i)_{jk}  {\rm d} \log (A_j) \wedge   {\rm d} \log (A_k), \hskip 10mm \mbox{where } \widetilde \varepsilon(i)_{jk}= d_j {\varepsilon(i)}_{jk} 
\ee
\ed

Let $(\A, \A_l, \A_r)\in {\cal A}(i)$. 
There is an angle invariant $u\in {\rm U}$ such that
\be
\la{lemma6.3.gs1}
(\A, \B_l, \B_r) = ([\U], \B^-, u\cdot \B^-) \in {\rm Conf}^*({\cal A}, {\cal B}, {\cal B}).
\ee
Note that $w_-(\B^-, u \cdot \B^-)=s_i$ if and only if $u\in \U_{s_{i}}$. Therefore $u=x_{i}(b)$ for some $b \neq 0$.
The potential of ${\cal A}(i)$ at the top vertex  is a function
\be
\la{lusztig.local.data.1}
 {\cal A}(i)\lra {\Bbb G}_m, ~~~~(\A, \A_l, \A_r)\lms b:={\cal W}_{i}(\A, \B_l, \B_r).
\ee
The potential $b$ is related to Lusztig's positive atlas for $\U$.  The following Lemma generalizes a local version of Theorem 2.19 in \cite{FZ98}, relating generalized minors and Lusztig's parametrization. 
\bl \la{L7.12}
\la{bfz.lusztig.lemma}
The potential $b$ in \eqref{lusztig.local.data.1} is 
\be 
\la{bfz.lusztig}
b= \frac{A_{i_e}}{A_{i_l}A_{i_r}} \prod_{j\in {\rm I}-\{i\}}A_j^{-C_{ij}}.\ee
\el
\begin{proof} Suppose that
\[
(\A, \A_l, \A_r) = ([\U], [\U^{-}], x_{i^*}(b)\cdot [\U^{-}]).
\]
Then $A_j=1$ if $j\neq i_e$. By Part 3 of Lemma \ref{idenity.basic.invariants.tt}. we have
$
A_{i_e} = b.
$
Identity \eqref{bfz.lusztig} holds. 

It is enough to show that the right hand of \eqref{bfz.lusztig} is invariant under rescaling
\[
(\A, \A_l, \A_r) \lms (\A, \,\A_l\cdot h_1, \, \A_r\cdot h_2).
\]
Let $h^*=w_0(h^{-1})$. By \eqref{ele.cor.cl.1}, we get $A_j\lms A_j t_j$, where the rescaling factors are
\[
t_j = \left\{    \begin{array}{ll} 
       \Lambda_{j}(h_1^*) \hskip 7mm & \mbox{if } j\in {\rm I}-\{i\}\\
         \Lambda_{i}(h_1^*) & \mbox{if } j=i_l\\
       \Lambda_{i}(h_2^*)& \mbox{if } j=i_r\\
      \Lambda_{i}\big(h_2^* s_{i}(h_1^{*})^{-1}\big) & \mbox{if }j=i_e. \\
   \end{array}\right.
\]
Note that
\[
s_{i}(\Lambda_{i})  =\Lambda_{i}-\alpha_{i}= \Lambda_{i} - \sum_{j\in {\rm I}} C_{ji}\Lambda_j = - \Lambda_i - \sum_{j\in {\rm I}-\{i\}} C_{ij}\Lambda_j. 
\]
Hence
\[
t_{i_e} = \Lambda_i(h_2^*) \cdot  s_i(\Lambda_i)(h_1^*)^{-1} = t_{i_l} t_{i_r} \prod_{j\in {\rm I}-\{i\}} t_j ^{C_{ij}}.
\]
Thus the right hand side of \eqref{bfz.lusztig} keeps intact.
\end{proof}

A result similar to Lemma \ref{L7.12} was independently proved in \cite[Theorem 5.39]{IIO}.


\bl 
\la{basic.formula.16.2.1.18.3.26}
Let $(\A, \A_l, \A_r)\in  {\cal A}(i)$.
If $j \neq i$, then $\Delta_j(\A, \A_l)= \Delta_j(\A, \A_r).$
\el
\begin{proof}
By \eqref{lemma6.3.gs1}, on the $\G$-orbit of $(\A_l, \A_r, \A')$, there exists a triple
\be
\la{lemma6.3.gs12}
\left([\U],\, [\U^-]\cdot h_1, \, x_i(b)[\U^-]\cdot h_2\right).
\ee
By the definition of $\Delta_j$, we get
\be
\la{basic.formula.16.2.1.18.3.2611}
\Delta_{j}(\A,\A_l) =\Lambda_{j^*}(h_1), ~~~~~ \Delta_{j}(\A, \A_r)=\Lambda_{j^*}(h_2).
\ee
It suffices to show that $\Lambda_{j^*}(h_2h_1^{-1})=1$ for $j\neq i$, or equivalently $h_2h_1^{-1}\in {\rm H}(s_{i^*})$.

Note that
$h s_{k} (h^{-1}) = \alpha_{k}^{\vee}(\alpha_{k}(h)).$
Therefore
\begin{align}
h(\A_r, \A_l)&= h\big(x_i(b)[\U^-]\cdot h_2, [\U^-]\cdot h_1\big)  \nonumber\\
&= h\big(x_i(b) [\U^-], [\U^-]\big) \cdot h_2 s_{i^*}(h_1^{-1})\nonumber\\
&= \alpha_{i^*}^\vee(b) \cdot  h_1 s_{i^*} (h_1^{-1}) \cdot h_2 h_1^{-1}. \nonumber \\
 &= \alpha_{i^*}^\vee\big(b\alpha_{i*}(h^{-1})\big) \cdot h_2 h_1^{-1}. \nonumber 
\end{align}
The condition $h(\A_r, \A_l)\in {\rm H}(s_{i^*})$ implies that 
$
h_2 h_1^{-1} \in {\rm H}(s_{i^*}).
$
\end{proof}

\bp
\la{coordinate.ele.conf}
The functions $\{A_j\}_{j \in{J}(i)}$ give rise to an isomorphism  
\[
{\cal A}(i) \stackrel{\sim}{\lra} ({\Bbb G}_m)^{r+2}.
\]
\ep
By definition, ${\cal A}_\circ(i)$ is a sub-torus of ${\cal A}(i)$ determined  by requiring $A_{i_e}=1$.

\begin{proof}
By \eqref{lemma6.3.gs12}, a configuration  $(\A, \A_l, \A_r)\in {\cal A}(i)$ is determined by its invariants $h_{1}, h_{2}, $ and $x_{i}(b)$. By \eqref{basic.formula.16.2.1.18.3.2611} and  Lemma \ref{basic.formula.16.2.1.18.3.26}, we have
\[
h_{1}^*=\alpha_1^\vee(A_1)\ldots\alpha_{i}^\vee(A_{i_l}) \ldots \alpha_r^\vee(A_r), \hskip 7mm h_{2}^* =  \alpha_1^\vee(A_1)\ldots\alpha_{i}^\vee(A_{i_r}) \ldots \alpha_r^\vee(A_r).
\]
By \eqref{bfz.lusztig},  $x_{i}(b)$ by is recovered by the $A$'s. Therefore a configuration $(\A_l, \A_r, \A')$ is uniquely determined by the set $\{A_j\}_{j\in {J}(i)}$, and vice versa. 
\end{proof}

\subsubsection{The elementary configuration space ${\cal A}({\overline{i}})$.} We define ${\cal A}({\overline{i}}):={\rm Conf}_{e}^{s_{{i}}}({\cal A})$. It parametrises $\G$-orbits of triples $(\A, \A^l, \A^r)\in {\cal A}^3$ such that
\[
w(\A, \A^l)= w(\A, \A^l)=w_0, ~~~~ w(\A^l, \A^r)= s_{{i}}, ~~ h(\A^r, \A^l)\in {\rm H}_{i}.
\]
It has a coordinate system parametrized by the set ${J}(\overline{i})$:
\be
\la{ele.cor.cl.2}
\forall (\A, \A^l, \A^r)\in {\cal A}(\overline{i}), \hskip 10mm
A_{{J}} := \left\{    \begin{array}{ll} 
       \Delta_{j}( \A^l, \A) \hskip 7mm & \mbox{if } {J}\in \overline{{\rm I}}-\{\overline{i}\}\\
       \Delta_{i}( \A^l, \A) & \mbox{if } \overline{j}=\overline{i}_l\\
       \Delta_{i}( \A^r, \A) & \mbox{if } \overline{j}=\overline{i}_r\\
      \Lambda_{i}(h(\A^l, \A^r)) & \mbox{if }\overline{j}=\overline{i}_e. \\
   \end{array}\right.
\ee

The space ${\cal A}(\overline{i})$ carries a 2-form $\Omega(\overline{i})$ encoding the matrix of ${\bf J}(\overline{i})$, and it has a subspace  ${\cal A}_\circ(\overline{i})$ defined by the equation $A_{\overline{i}_e}=1$. The partial potential of ${\cal A}(\overline{i})$ at the bottom vertex is
\be
\la{bfz.lusztig2t}
\overline{b}:= {\cal W}_{i^*}(\A, \B^r, \B^l) = \frac{A_{\overline{i}_e}}{A_{\overline{i}_l}A_{\overline{i}_r}}\prod_{\overline{j}\in \overline{{\rm I}}-\{\overline{i}\}} A_{\overline{j}}^{-C_{ij}}.
\ee
\bl Let $\overline{b}$ be the potential as above. Let  $h_l:=h(\A^l, \A), ~h_r:=h(\A_r, \A).$ Then
\be
(\A, \A^l, \A^r)= ([\U^-],~ h_l \cdot [\U], ~y_i(\overline{b})h_r\cdot [\U]).
\ee
\el
\subsubsection{The reflections.}
Given an element $i \in {\rm I}$, the reflection $r$  is an isomorphism 
\be
\begin{split}
&r:~{\cal A}(i) \stackrel{\sim}{\lra} {\cal A}(\overline{i}), \\
&(\A', \A_l, \A_r) \lms (\A, \A^l, \A^r), \\
\end{split}
\ee
where $\A=\A_r$, $ \A^r=\A'$, 
and $\A^l$ is the unique decorated flag such that
\[
w(\A_l, \A^l) = w_0s_i,   \hskip 7mm w(\A^l,  \A') = s_{i}, \hskip 7mm h(\A^l, \A')= w_0\big(h(\A_r, \A_l)^{-1}\big).
\]
By Lemma \ref{lem.part.flag.decom}, we have $w(\A^l, \A)=w_0$. Therefore  $r$ is well defined. By definition, the image of $\mathcal{A}_\circ(i) \subset  {\mathcal{A}}(i)$ under $r$ is $\mathcal{A}_\circ(\overline{i})$. 
The reflection $r$ is illustrated by the following figure. 
\begin{center}
\begin{tikzpicture}[scale=0.9]
 \path [fill=gray!30] (0,0) -- (-1,-0.5) -- (0.6,-3) -- cycle;  
 \node [circle,draw=red,fill=red,minimum size=3pt,inner sep=0pt,label=above:{\small ${\A}'$}] (a) at (0,0) {};
 \node [circle,draw=red,fill=red,minimum size=3pt,inner sep=0pt,label=below:{\small ${\A}_l$}]  (b) at (-0.6,-3){};
 \node [circle,draw=red,fill=red,minimum size=3pt,inner sep=0pt,label=below:{\small ${\A}_r$}] (c) at (0.6,-3) {};
 \node [circle,draw=red,fill=white,minimum size=3pt,inner sep=0pt,label=left:{\small ${\A}^l$}] (d) at (-1,-0.5) {};
 \node[blue] at (0,-2.8) {\small $s_{i^*}$};
 \node[blue] at (-0.5,-0) {\small $s_{i}$};
 \node[blue] at (-1.3,-2) {\small $w_0s_i$}; 
  
 \foreach \from/\to in {a/b, a/c}
           \draw[thick] (\from) -- (\to);
  \draw[thick, green] (b) -- (c);
  \draw[thick, green] (a) -- (d);                        
  \draw[thick, directed] (b) -- (d);    
  \draw[draw, dashed] (c) -- (d);              
  \draw[thick,->] (2.5,-1.5) -- (4.5,-1.5);        
  \node [] at (3.5, -1.2) {$r$};
\begin{scope}[shift={(7,-3)}]

\path [fill=gray!30] (0,0) -- (-0.6,3) -- (0.6,3) -- cycle;  
\node [blue] at (0,3.2) {\small $s_{i}$};
\node [label=below:{\small ${\A} \mbox{=}{\A}_r$}] at (0.5,0) {};
\node [label=above:{\small ${\A}^r \mbox{=}{\A}'$}] at (1,2.9) {};
\node [label=above:{\small ${\A}^l$}] at (-0.8,2.9) {};
\node [circle,draw=red,fill=red,minimum size=3pt,inner sep=0pt] (a) at (0,0) {};
\node [circle,draw=red,fill=red,minimum size=3pt,inner sep=0pt]  (b) at (-0.6,3){};
\node [circle,draw=red,fill=red,minimum size=3pt,inner sep=0pt] (c) at (0.6,3) {};
\draw[thick, green] (b) -- (c);
\foreach \from/\to in {a/b, a/c}
                   \draw[thick] (\from) -- (\to);
 \end{scope}                  
 \end{tikzpicture}
\end{center}

\bl 
\la{identity.reflect.f}
In terms of the coordinates \eqref{ele.cor.cl.1}$\&$\eqref{ele.cor.cl.2}, the reflection $r$ is expressed as
\be
\la{identity.reflect.f.1}
 r^*A_{\overline{j}} =\left\{ \begin{array}{ll}     
 {\displaystyle A_{i_l}^{-1}A_{i_e}^2  \prod_{k\in {\rm I}-\{i\}} A_k^{-C_{ik}}} & \mbox{if } \overline{j}={\overline{i}_l}\\
        A_{j} & \mbox{if } \overline{j}\neq {\overline{i}_l}.\\
   \end{array} \right.
\ee
The potential $\overline{b}$ of ${\cal A}(\overline{i})$ at the vertex $\A$ is transformed by the map $r$ as follows: 
\be
\la{identity.reflect.f.3} r^*\overline{b} =  \frac{A_{i_l}}{A_{i_r}A_{i_e}}.
\ee
\el
\begin{proof}
The second formula of  \eqref{identity.reflect.f.1}  is due to the definition of $r$ and Lemma \ref{basic.formula.16.2.1.18.3.26}.
By \eqref{chain.formula.for.cartan.5.8},
\[
h(\A', \A_l) {=} h(\A', \A^l) \cdot s_i \big(h(\A^l, \A_l)\big). 
\]
Note that $h(\A', \A^l)\in \alpha_{i}^\vee({\Bbb G}_m)$. Therefore
\be
\la{2016.5.8.11.48.1}
h(\A^l, \A_l) = s_i\big(h(\A', \A_l) h(\A', \A^l)^{-1}\big) = s_i(h(\A', \A_l)) h(\A', \A^l).
\ee
Meanwhile
\be
\la{2016.5.8.11.48.1.tt}
(w_0s_{i^*}) \big(h(\A_l, \A_r)\big) = w_0(h(\A_l, \A_r)^{-1}) = h(\A', \A^l).
\ee
Therefore
\be
\la{2016.5.8.11.48.2}
h(\A^l, \A_r) \stackrel{\eqref{chain.formula.for.cartan.5.8}}{=} h(\A^l, \A_l) \cdot (w_0s_{i^*}) \big(h(\A_l, \A_r)\big) \stackrel{\eqref{2016.5.8.11.48.1}\&\eqref{2016.5.8.11.48.1.tt}}{=} h(\A', \A^l)^2 \cdot s_i\big(h(\A', \A_l)\big). 
\ee
The first formula of  \eqref{identity.reflect.f.1} follows since
\[
r^*A_{{\overline{i}_l}} =  \Lambda_i(h(\A^l, \A_r)) 
=\Lambda_i\Big(h(\A', \A^l)^{2} \cdot s_i\big(h(\A', \A_l)\big)\Big)= A_{i_e}^2 A_{i_l}^{-1} \prod_{k\in {\rm I}-\{i\}} A_k^{-C_{ik}}.
\]
Formula \eqref{identity.reflect.f.3} follows by plugging \eqref{identity.reflect.f.1}  into \eqref{bfz.lusztig}. 
\end{proof}

\bp 
\la{2017.2.24.21.20h}
The reflection $r$ preserves the   canonical 2-form:
$
r^*\Omega(\overline{i}) = \Omega({i}).
$
\ep

\begin{proof}
Set $a_j:= {\rm d} \log(A_j)$. 
Recall the potential $b$ of ${\cal A}(i)$, and $\overline{b}$ of ${\cal A}(\overline{i})$. We set
\[ p := {\rm d} \log b = -a_{i_l}-a_{i_r}+a_{i_e}-\sum_{j\in {\rm I}-\{i\}}  C_{ij} {a_j}.
\hskip 1cm
 \overline{p} :=  {\rm d} \log \overline{b} = -a_{\overline{i}_l}-a_{\overline{i}_r}+a_{\overline{i}_e}-\sum_{\overline{j} \in \overline{{\rm I}}-\{\overline{i}\}}  C_{ij} {a_{\overline{j}}}.\]
A direct calculation shows that
\[
\Omega(i) 
=(a_{i_l}- a_{i_r}) \wedge (p+ a_{i_e}), \hskip 1cm 
\Omega(\overline{i}^\ast) =(\overline{p}+ a_{\overline{i}_e}) \wedge (a_{\overline{i}_l}- a_{\overline{i}_r}). 
\]
By Lemma \ref{identity.reflect.f}, we have
\[
r^*(\overline{p}) = a_{i_l}-a_{i_r}-a_{i_e}, \hskip 7mm r^*(a_{\overline{i}_e})=a_{i_e},  \hskip 7mm r^*(a_{\overline{i}_l}) =p+  a_{i_r}+a_{i_e} , \hskip 7mm  r^*(a_{\overline{i}_r}) = a_{i_r}.
\]
Therefore
$r^*(\overline{p}+a_{\overline{i}_e}) = a_{i_l}-a_{i_r}$ and $r^*(a_{\overline{i}_l}-a_{\overline{i}_r}) = p_{i_l}+a_{i_e}.
$ 
Hence
$r^*\Omega({\overline{i}}) =\Omega(i).$
\end{proof}


\subsection{Proof of Theorem \ref{Th2.3}, part 2}

\medskip

Let us fix a reduced word ${\bf i}=(i_1, \ldots, i_m)$ of $(u, v)$. We consider the space
\be
\la{2017.3.5.16.26}
{\cal A}(i_1)\times \cdots \times {\cal A}(i_m).
\ee

\bd The space ${\cal A}({\bf i})$ is the subspace of \eqref{2017.3.5.16.26} such that
\begin{itemize}
\item for $k=1,\ldots, m-1$, the right pair of decorated flags in $ {\cal A}_{(\circ)}(i_k)$ and the left pair in $ {\cal A}_{(\circ)}(i_{k+1})$ belong to the same $\G$-orbit.
\item we take the subspace $\mathcal{A}_\circ(i_k) \subset  {\mathcal{A}}(i_k)$ in every component of the product unless the coroot $\beta_{k}^{\bf i} $ is simple.
\end{itemize}
\ed

\bl We have 
\[
{\cal A}({\bf i})\stackrel{\sim}{=}{\rm T}_{\bf i}^{\mathscr A}:= ({\Bbb G}_m)^{\# {J}({\bf i})}.
\]
\el
\begin{proof} By Proposition \ref{coordinate.ele.conf}, every elementary configuration space is isomorphic to a torus. Thus \eqref{2017.3.5.16.26} is a torus. By definition, ${\cal A}({\bf i})$ is a sub-torus of $\eqref{2017.3.5.16.26}$,  consistent with the amalgamation of elementary quivers.
\end{proof}

Gluing  adjacent pairs of decorated flags for all $k\in\{1,\ldots, m-1\}$ we get a regular map 
\[
c_{{\mathscr A},{\bf i}}: ~ {\cal A}({\bf i}) \lra {\rm Conf}_u^v({\cal A}).
\]
\bl The map $c_{{\mathscr A},{\bf i}}$ is an open embedding.
\el
\begin{proof} We prove the Lemma by defining the inverse birational map $c_{{\mathscr A},{\bf i}}^{-1}$. Let 
\[
(\A_l,\A_r,  \A^l, \A^r) \in {\rm Conf}_u^v({\cal A}).
\]
Picking the ${\rm I}$-part of indices in ${\bf i}$, we get a reduced word ${\bf i}_u=(j_1, \ldots, j_p)$ of $u$. 
By Lemma \ref{decomp.partial.dec.flag.15.02}, the word ${\bf i}_u^\ast$ uniquely determines a decomposition of the bottom pair $(\A_l, \A_r)$:
\be
\la{bochain11st}
\A_l= \A_0\stackrel{j_1^*}{\lra} \A_1\stackrel{j_2^*}{\lra}\ldots \stackrel{j_p^*}{\lra}\A_p = \A_r.
\ee
Picking the $\overline{I}$-part, we get a reduced word ${\bf i}_v= (k_1, \ldots, k_q)$ of $v$.
For $(\A^l, \A^r)$, note the appearance of $\rho_v(-1)$ in their $h$-condition. We replace the second condition of Lemma \ref{decomp.partial.dec.flag.15.02} by
\begin{itemize}
\item If $\beta_{s}^{{\bf i}_v}$ is simple, then $h(\A^{s-1}, \A^{s}) \in \alpha_{k_s}^\vee({\Bbb G}_m)$, otherwise\footnote{It is equivalent to the condition $h(\B^{s}, \B^{s-1})=\alpha_{k_s}^\vee(-1)$.}
 $h(\B^{s-1}, \B^{s})=1$. 
\end{itemize}
Similarly, it gives rise to a sequence of intermediate decorated flags
\be
\la{tochain12st}
\A^l= \A^0\stackrel{k_1}{\lra} \A^1\stackrel{k_2}{\lra}\cdots \stackrel{k_q}{\lra}\A^q = \A^r.
\ee 
Using the chains \eqref{bochain11st} and \eqref{tochain12st}, the ${\bf i}$ gives rise to a decomposition of ${\rm Conf}_u^v({\cal A})$ into elementary configuration spaces, as explained by the following example. 

\begin{example} Consider the root system of type $A_4$. Let ${\bf i}=(2, \overline{2}, \overline{3}, 3, \overline{2})$ be a reduced word for $(u,v)$.  As  illustrated by Figure \ref{coord.decom.12.00h}, the ${\bf i}$ provides a decomposition of ${\rm Conf}_u^v({\cal A})$ into elementary configuration spaces. 

\begin{figure}[h]
\epsfxsize 100pt
\center{
\begin{tikzpicture}[scale=0.95]
\node [circle,draw=red,fill=red,minimum size=3pt,inner sep=0pt] (a) at (-1.8,0) {};
\node [circle,draw=red,fill=white,minimum size=3pt,inner sep=0pt,label=above:{\small ${\A}^1$}] (e) at (-0.6,0) {};
\node [circle,draw=red,fill=white,minimum size=3pt,inner sep=0pt,label=above:{\small ${\A}^2$}] (f) at (0.6,0) {};
\node [circle,draw=red,fill=red,minimum size=3pt,inner sep=0pt] (g) at (1.8,0) {};
 \node [circle,draw=red,fill=red,minimum size=3pt,inner sep=0pt]  (b) at (-1.2,-3){};
 \node [circle,draw=red,fill=red,minimum size=3pt,inner sep=0pt] (c) at (1.2,-3) {};
 \node [circle,draw=red,fill=white,minimum size=3pt,inner sep=0pt,label=below:{\small ${\A}_1$}] (d) at (0,-3) {};
 \node at (-2.1, 0.4) {\small {${\A}^l=\A^0~$}};
  \node at (2.28, 0.4) {\small $~{\A}^3=\A^r$};
   \node at (-1.5, -3.4) {\small ${\A}_l=\A_0~$};
  \node at (1.68, -3.4) {\small $~{\A}_2=\A_r$};
 \foreach \from/\to in {a/b, b/d,  g/c, e/f}
                   \draw[thick] (\from) -- (\to);
 \foreach \from/\to in {a/e, f/g, d/c}                   
                   \draw[thick, green] (\from) -- (\to);
 \foreach \from/\to in {a/d, e/d, f/d, f/c}                   
                   \draw[dashed] (\from) -- (\to);                 
 \node[blue] at (-0.7,-2.8) {\small $s_{2^\ast}$};
  \node[blue] at (0.7,-2.8) {\small $s_{3^\ast}$};
  \node[blue] at (-1.2,0.2) {\small $s_{{2}}$};
   \node[blue] at (0,0.2) {\small $s_{{3}}$};
   \node[blue] at (1.2,0.2) {\small $s_{{2}}$};  
\draw[thick, ->] (3, -1.5) --(4.5, -1.5);     

\draw[thick] (6, -3) -- (7.2,-3) -- (5.4,0)--cycle;

\draw[thick]  (5.7,0.1) -- (7.5, -2.9) -- (6.9,0.1);
\draw[thick, green]  (5.7,0.1) -- (6.9,0.1);

\draw[thick] (7.8, -3) -- (8.4,0) -- (7.2,0)--cycle;

\draw[thick] (8.1, -3) -- (8.7,0) -- (9.3,-3);
\draw[thick, green] (8.1, -3) -- (9.3,-3);

\draw[thick] (9,0) -- (9.6, -3) -- (10.2,0);
\draw[thick, green] (10.2,0) -- (9,0);

 \end{tikzpicture}
        }
\caption{The decomposition of ${\rm Conf}_u^v({\cal A})$ into elementary configuration spaces.}
\label{coord.decom.12.00h}
\end{figure}

\end{example}
In general, we obtain a rational decomposition map
\[
d^{\bf i}=(d^{\bf i}_1, \ldots, d^{\bf i}_m): ~ ~~{\rm Conf}_u^v({\cal A})\lra {\cal A}(i_1)\times \ldots \times {\cal A}(i_m). \]
Here $d^{\bf i}$ is well-defined if  all the pairs of intermediate decorated flags connecting by dashed lines in Figure  \ref{coord.decom.12.00h} are of generic position. In this case, the image of $d^{\bf i}$ is  ${\cal A}({\bf i})$.
\end{proof}
\begin{remark} 
Each ${\cal A}(i)$ has a local coordinate system indexed by ${J}(i)$. Using the map $\pi^{\bf i}$, we get a local coordinate system of ${\rm Conf}_u^v({\cal A})$ indexed by ${J}({\bf i})$. For example, each $k\in\{1,\ldots, m\}$ corresponds to a dashed line and therefore a pair of intermediate decorated flags $(\A^{t_k},\A_{b_k})$ as in Figure \ref{coord.decom.12.00h}. We a get a regular function \[A_k:=\Delta_{i_k}(\A^{t_k},\A_{b_k}).\] These coordinates give rise to a regular map
\be
\pi^{\bf i}: {\rm Conf}_u^v(\mathcal{A})\lra {\Bbb A}^{\# {J}({\bf i})}. 
\ee

The exchange matrix $\varepsilon({\bf i})_{ij}$ of ${\bf J}({\bf i})$ determines a 2-form on ${\rm Conf}_u^v({\cal A})$:
\be
\la{canonical.2.form.u.v}
\Omega_{u,v}=\Omega({\bf i}):=  \sum_{i, j \in {\bf J}({\bf i}) } \widetilde \varepsilon({\bf i})_{ij}  {\rm d} \log (A_i) \wedge d\log (A_j).
\ee
 Theorem \ref{main.result1.a.structure} tells that the 2-form $\Omega_{u,v}$ is independent of  ${\bf i}$. 
By definition, we have a decomposition
\be
\la{canonical.2.form.u.v.deco}
\Omega_{u,v} := \Omega_{\bf i}^{\bf J} + \Omega_{\bf i}^{\bf H}, \hskip 1cm \mbox{where } \Omega_{\bf i}^{\bf J}:=\sum_{k=1}^{m} (d^{\bf i}_k)^*\Omega(i_k),
\ee
and $\Omega_{\bf i}^{\bf H}$ corresponds to the quiver ${\bf H}({\bf i})$.

\end{remark}

\paragraph{\bf Proof of Theorem \ref{main.result1.a.structure},  part 2.}
It remains to show that the transition map of the coordinate systems given by two reduced words ${\bf i}$ and ${\bf i}'$ are cluster $K_2$-transformations. 
It  reduces to   the following basic cases.

\begin{itemize}
\item[] \hskip -1cm {\bf Case 1:}  $(i,\overline{i}) \sim (\overline{i}, i)$.
\end{itemize}
 
Let $(\A_l, \A_r, \A^l, \A^r)\in 
{\rm Conf}_{s_i} ^{{s}_{{i}}}({\cal A}).$ As shown on Figure \ref{coord.decom.17.15.p1h}, it has 2 different decompositions. 
\begin{figure}[ht]
\epsfxsize 90pt
\center{
\begin{tikzpicture}
\node [circle,draw=red,fill=red,minimum size=3pt,inner sep=0pt,label=above:{\small ${\A}^l$}] (a) at (0,0) {};
\node [circle,draw=red,fill=red,minimum size=3pt,inner sep=0pt,label=above:{\small ${\A}^r$}] (b) at (1.2,0) {};
 \node [circle,draw=red,fill=red,minimum size=3pt,inner sep=0pt,label=below:{\small ${\A}_l$}]  (d) at (0,-3){};
 \node [circle,draw=red,fill=red,minimum size=3pt,inner sep=0pt,label=below:{\small ${\A}_r$}] (c) at (1.2,-3) {};
 \foreach \from/\to in {b/c, d/a}
                  \draw[thick] (\from) -- (\to);
  \foreach \from/\to in {a/b, d/c}
                  \draw[thick,green] (\from) -- (\to);                 
\foreach \from/\to in {a/c}                   
              \draw[dashed] (\from) -- (\to);
     \draw[blue, thick, ->] (0,-0.5) arc (-90:0:0.5);                
 \node[blue] at (0.6,-2.8) {\small $s_{i^\ast}$};
  \node[blue] at (0.6,0.2) {\small ${s}_{{i}}$}; 
    \node [circle,draw=red,fill=white,minimum size=3pt,inner sep=0pt,label=right:{\small ${\B}_2$}] at (0.4,-1) {};  
      \node [circle,draw=red,fill=white,minimum size=3pt,inner sep=0pt,label=left:{\small ${\B}_1$}] at (0,-1) {};  
  \draw[thick,<->] (2.5,-1.5) -- (4.5,-1.5);            

\begin{scope}[shift={(5.8,0)}]
\node [circle,draw=red,fill=red,minimum size=3pt,inner sep=0pt,label=above:{\small ${\A}^l$}] (a) at (0,0) {};
\node [circle,draw=red,fill=red,minimum size=3pt,inner sep=0pt,label=above:{\small ${\A}^r$}] (b) at (1.2,0) {};
 \node [circle,draw=red,fill=red,minimum size=3pt,inner sep=0pt,label=below:{\small ${\A}_l$}]  (d) at (0,-3){};
 \node [circle,draw=red,fill=red,minimum size=3pt,inner sep=0pt,label=below:{\small ${\A}_r$}] (c) at (1.2,-3) {};
 \foreach \from/\to in {b/c, d/a}
                  \draw[thick] (\from) -- (\to);
  \foreach \from/\to in {a/b, d/c}
                  \draw[thick,green] (\from) -- (\to);  
\foreach \from/\to in {d/b}                   
              \draw[dashed] (\from) -- (\to);
                   
 \node[blue] at (0.6,-2.8) {\small $s_{i^\ast}$};
  \node[blue] at (0.6,0.2) {\small ${s}_{{i}}$}; 
  \node [circle,draw=red,fill=white,minimum size=3pt,inner sep=0pt,label=left:{\small ${\B}_1$}] at (0,-1) {};   
   \draw[blue, thick, ->] (0,-.5) arc (-90:0:0.5);      
 \end{scope}   
 \end{tikzpicture}

        }
\caption{Decompositions of ${\rm Conf}_{s_i}^{s_i}({\cal A})$.}
\label{coord.decom.17.15.p1h}
\end{figure}

For $j\in {\rm I}-\{i\}$, we have
\[
A_j := \Delta_j(\A^l, \A_l)=\Delta_j(\A^l, \A_r)=\Delta_j(\A^r, \A_l) =\Delta_j(\A^r, \A_r).
\]
For $i$, we have
\[
A_{i_l}= \Delta_i (\A^l, \A_l), ~~~~~ Y= \Delta_i(\A^l, \A_r), ~~~~~Z= \Delta_i(\A^r, \A_l),~~~~  A_{i_r} = \Delta_i (\A^r, \A_r), 
\]
\[
A_{i_e}= \Lambda_{i^\ast}(h(\A_r, \A_l)), ~~~~~~ A_{\overline{i}_e}= \Lambda_{i}(h(\A^l, \A^r)).
\]
Here  $Y,~ Z$ correspond to the  non-frozen vertices of the  quivers in case 1 of Section \ref{thm.3.6.st}. The following Lemma shows that the transition map from $(i, \overline{i})$ to $( \overline{i}, i)$ is a cluster $K_2$-mutation at the unfrozen vertex.

\bl \la{CL5.3} We have
\be
\la{cluster.partial.identity.may2}
YZ = A_{i_l} A_{i_r} + A_{i_e}A_{\overline{i}_e} \prod_{j\in {\rm I}-\{i\}} A_j^{-C_{ij}}.
\ee
\el

\begin{proof}
Set $\B_l=\pi(\A_l)$ and so on. As shown on Figure \ref{coord.decom.17.15.p1h}, we 
take $\B_1, \B_2\in {\cal B}$ such that
\[
w(\B_l, \B_1)= w_0s_{i}, ~~~~~w(\B_1, \B^l)= s_{i}; ~~~ w(\B_r, \B_2)= w_0s_{i}, ~~~~~w(\B_2, \B^l)= s_{i}.
\]
Recall the potential of ${\rm Conf}({\cal A}, {\cal B}, {\cal B})_{{\rm SL}_2}$.
By the additivity of potential, we get
\[
 {\cal W}(\A^l, \B_1, \B_2)+{\cal W}  (\A^l, \B_2, \B^r)={\cal W}(\A^l, \B_1, \B^r). 
\]
By \eqref{bfz.lusztig} and \eqref{identity.reflect.f.3}, we get 
\[
 {\cal W}  (\A^l, \B_1, \B_2) =\frac{A_{{i}_e}}{ A_{i_l}Y} \prod_{j\in {\rm I}-\{i\}} A_j^{-C_{ij}}, \hskip 9mm
 {\cal W}(\A^l, \B_2, \B^r)=\frac{A_{i_r}}{YA_{\overline{i}_e}},
 \]
 \[
 {\cal W}(\A^l, \B_1, \B^r)= \frac{Z}{A_{i_l}A_{\overline{i}_e}}.
\]
Therefore
\[
\frac{Z}{A_{i_l}A_{\overline{i}_e}} =\frac{A_{i_r}}{YA_{\overline{i}_e}}+\frac{A_{{i}_e}}{ A_{i_l}Y} \prod_{j\in {\rm I}-\{i\}} A_j^{-C_{ji}}.
\]
Multiplying both sides by $YA_{i_l} A_{\overline{i}_e}$, we get \eqref{cluster.partial.identity.may2}.
\end{proof}

Let us consider the situation where $A_{i_e}=A_{\overline{i}_e'}=1$. We define the local primary coordinates
\[
P_1= \frac{Y}{A_{i_l}}, \hskip 7mm P_2= \frac{Y}{A_{i_r}}; \hskip 10mm
P_1'= \frac{A_{i_l}}{Z}, \hskip 7mm P_2'= \frac{A_{i_r}}{Z}
\]
and the cluster Poisson coordinate
\[
X= A_{i_l}^{-1}A_{i_r}^{-1}\prod_{j\in {\rm I}-\{i\}}A_j^{C_{ij}}
\]
The follow result will be useful in the future.

\begin{corollary} 
\label{local.primary}
We have
\[
P_1'=P_2(1+X^{-1})^{-1}, \hskip 7mm P_2'=P_1(1+X^{-1})^{-1}
\]
\end{corollary}
\begin{proof} It follows by plugging the Formula \eqref{cluster.partial.identity.may2} directly.
\end{proof}

\begin{itemize}
\item[] \hskip -1cm  {\bf Case 2:}  $(i, j, i)\sim (j, i, j)$, where $C_{ji}=C_{ij}=-1$. 
\end{itemize}

The space ${\rm Conf}_{s_is_js_i}^e({\cal A})$ has two decompositions. 

\begin{figure}[ht]
\epsfxsize 100pt
\center{
\begin{tikzpicture}
\node [circle,draw=red,fill=red,minimum size=3pt,inner sep=0pt,label=above:{\small ${\A}'$}] (a) at (0,0) {};
\node [circle,draw=red,fill=red,minimum size=3pt,inner sep=0pt,label=below:{\small ${\A}_l$}] (b) at (-1.8,-3) {};
 \node [circle,draw=red,fill=white, minimum size=3pt,inner sep=0pt,label=below:{\small ${\A}_2$}]  (c) at (-0.6,-3){};
 \node [circle,draw=red,fill=white, minimum size=3pt,inner sep=0pt,label=below:{\small ${\A}_3$}] (d) at (0.6,-3) {};
 \node [circle,draw=red,fill=red,minimum size=3pt,inner sep=0pt,label=below:{\small ${\A}_r$}] (e) at (1.8,-3) {};
 \foreach \from/\to in {a/b, d/c,  e/a}
                  \draw[thick] (\from) -- (\to);
 \foreach \from/\to in {b/c, d/e}
                  \draw[thick, green] (\from) -- (\to);                  
\foreach \from/\to in {a/c,a/d}                   
              \draw[dashed] (\from) -- (\to);
                   
 \node[blue] at (-1.2,-2.8) {\small $s_{i^\ast}$};
  \node[blue] at (0,-2.82) {\small ${s}_{j^\ast}$}; 
   \node[blue] at (1.2,-2.8) {\small $s_{i^\ast}$};
  \draw[thick,<->] (2.5,-1.5) -- (4.5,-1.5);            

\begin{scope}[shift={(7,0)}]
\node [circle,draw=red,fill=red,minimum size=3pt,inner sep=0pt,label=above:{\small ${\A}'$}] (a) at (0,0) {};
\node [circle,draw=red,fill=red,minimum size=3pt,inner sep=0pt,label=below:{\small ${\A}_l$}] (b) at (-1.8,-3) {};
 \node [circle,draw=red,fill=white, minimum size=3pt,inner sep=0pt,label=below:{\small ${\A}_4$}]  (c) at (-0.6,-3){};
 \node [circle,draw=red,fill=white, minimum size=3pt,inner sep=0pt,label=below:{\small ${\A}_5$}] (d) at (0.6,-3) {};
 \node [circle,draw=red,fill=red,minimum size=3pt,inner sep=0pt,label=below:{\small ${\A}_r$}] (e) at (1.8,-3) {};
 \foreach \from/\to in {a/b, d/c,  e/a}
                  \draw[thick] (\from) -- (\to);
 \foreach \from/\to in {b/c, d/e}
                  \draw[thick, green] (\from) -- (\to);  
\foreach \from/\to in {a/c,a/d}                   
              \draw[dashed] (\from) -- (\to);
                   
 \node[blue] at (-1.2,-2.82) {\small $s_{j^\ast}$};
  \node[blue] at (0,-2.8) {\small ${s}_{i^\ast}$}; 
   \node[blue] at (1.2,-2.82) {\small $s_{j^\ast}$};      
\end{scope}   
 \end{tikzpicture}
        }
\caption{Two different decompositions.}
\label{coord.decom.17.15.p2h}
\end{figure}

Consider functions
$$
\begin{array}{ll}
A_{i_l}:= \Delta_i(\A', \A_l), ~~~~~~~~  &A_{j_l}:= \Delta_j(\A', \A_j), \\
A_{i_r}:= \Delta_{i}(\A', \A_r),~~~ &A_{j_r} := \Delta_{j}(\A', \A_r), \\
A_{i_e}: =\Lambda_{i^\ast} (h(\A_r, \A_l)), ~~~~~~~~~~ &A_{j_e} :=\Lambda_{j^\ast} (h(\A_r, \A_l))=\Lambda_{i^\ast}(h(\A_2, \A_l)),\\
Y:=\Delta_{i}(\A', \A_2)=\Delta_{i}(\A', \A_3),  ~~~~~~ &Z:=\Delta_{j}(\A', \A_4)=\Delta_{j}(\A', \A_5).  \\
\end{array}
$$
\bl We have
\be
\la{cluster.partial.identity.may1}
YZ =  A_{i_l} A_{i_e} A_{j_r} +  A_{j_l} A_{j_e} A_{i_r}.
\ee
\el
\begin{proof} 
Recall the partial potential ${\cal W}_{i}$. We have\footnote{It is equivalent to the fact that if $x_{i}(a)x_j(b)x_i(c)= x_{j}(a')x_i(b') x_{j}(c')$, then $a+c =b'$.}
\[
{\cal W}_{i}(\A', \B_l, \B_r) = {\cal W}_{i}(\A', \B_l, \B_2)+ {\cal W}_{i}(\A', \B_3, \B_r)=  {\cal W}_{i}(\A', \B_4, \B_5).
\]
For $k\in {\rm I}-\{i, j\}$, we set $A_k= \Delta_k(\A', \A_l)$. Set 
$P := \prod_{k\in {\rm I}-\{i, j\}} A_k^{-C_{ki}}.
$ 
By \eqref{bfz.lusztig}, we get
\[
{\cal W}_{i}(\A', \B_l, \B_2) = \frac{A_{j_e}}{ A_{i_l}Y} A_{j_l} P,                 
\hskip 7mm 
{\cal W}_{i}(\A', \B_3, \B_r)= \frac{A_{i_e}}{ Y A_{i_r}} A_{j_r} P, 
\hskip 7mm
{\cal W}_{i}(\A', \B_4, \B_5) =  \frac{1}{A_{i_l}A_{i_r}} Z P.
\]
Therefore
\[ 
\frac{A_{j_e}}{ A_{i_l}Y} A_{j_l}P +  \frac{A_{i_e}}{ YA_{i_r}} A_{j_r}P  =  \frac{1}{A_{i_l}A_{i_r}} ZP .
\]
Multiplying both sides by $A_{i_l} Y A_{i_r}P^{-1}$, we get \eqref{cluster.partial.identity.may1}.
\end{proof}

\begin{itemize}
\item[] \hskip -1cm  {\bf Case 3:}  $(i, j, i, j)\sim (j, i, j, i)$, where  $C_{ij}=-1$ and $C_{ji}=-2$.
\end{itemize}
 
Let $u=s_is_js_is_j$. We have two different decompositions for ${\rm Conf}_u^e(\mathcal{A})$.

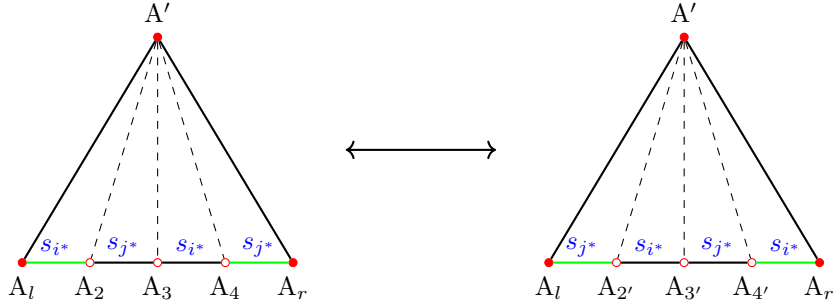
\begin{center}
\begin{tikzpicture}
\node [circle,draw=red,fill=red,minimum size=3pt,inner sep=0pt,label=above:{\small ${\A}'$}] (a) at (0,0) {};
\node [circle,draw=red,fill=red,minimum size=3pt,inner sep=0pt,label=below:{\small${\A}_l$}] (b) at (-1.8,-3) {};
 \node [circle,draw=red,fill=white, minimum size=3pt,inner sep=0pt,label=below:{\small ${\A}_2$}]  (c) at (-0.9,-3){};
 \node [circle,draw=red,fill=white, minimum size=3pt,inner sep=0pt,label=below:{\small ${\A}_3$}] (d) at (0,-3) {};
  \node [circle,draw=red,fill=white, minimum size=3pt,inner sep=0pt,label=below:{\small ${\A}_4$}] (e) at (0.9,-3) {};
 \node [circle,draw=red,fill=red,minimum size=3pt,inner sep=0pt,label=below:{\small ${\A}_r$}] (f) at (1.8,-3) {};
 \foreach \from/\to in {a/b,  d/c, d/e,  f/a}
                  \draw[thick] (\from) -- (\to);
 \foreach \from/\to in {b/c, e/f}
                  \draw[thick, green] (\from) -- (\to);   
\foreach \from/\to in {a/c,a/d,a/e}                   
              \draw[dashed] (\from) -- (\to);
                   
 \node[blue] at (-1.35,-2.8) {\small $s_{i^\ast}$};
  \node[blue] at (-.45,-2.8) {\small ${s}_{j^\ast}$}; 
   \node[blue] at (.45,-2.8) {\small $s_{i^\ast}$};
   \node[blue] at (1.35,-2.8) {\small $s_{j^\ast}$};
  \draw[thick,<->] (2.5,-1.5) -- (4.5,-1.5);            

 \begin{scope}[shift={(7,0)}]
\node [circle,draw=red,fill=red,minimum size=3pt,inner sep=0pt,label=above:{\small ${\A}'$}] (a) at (0,0) {};
\node [circle,draw=red,fill=red,minimum size=3pt,inner sep=0pt,label=below:{\small ${\A}_l$}] (b) at (-1.8,-3) {};
 \node [circle,draw=red,fill=white, minimum size=3pt,inner sep=0pt,label=below:{\small ${\A}_{2'}$}]  (c) at (-0.9,-3){};
 \node [circle,draw=red,fill=white, minimum size=3pt,inner sep=0pt,label=below:{\small ${\A}_{3'}$}] (d) at (0,-3) {};
  \node [circle,draw=red,fill=white, minimum size=3pt,inner sep=0pt,label=below:{\small ${\A}_{4'}$}] (e) at (0.9,-3) {};
 \node [circle,draw=red,fill=red,minimum size=3pt,inner sep=0pt,label=below:{\small ${\A}_r$}] (f) at (1.8,-3) {};
 \foreach \from/\to in {a/b,  d/c, d/e,  f/a}
                  \draw[thick] (\from) -- (\to);
 \foreach \from/\to in {b/c, e/f}
                  \draw[thick, green] (\from) -- (\to);                  
\foreach \from/\to in {a/c,a/d,a/e}                   
              \draw[dashed] (\from) -- (\to);
                   
 \node[blue] at (-1.35,-2.8) {\small $s_{j^\ast}$};
  \node[blue] at (-.45,-2.8) {\small ${s}_{i^\ast}$}; 
   \node[blue] at (.45,-2.8) {\small $s_{j^\ast}$};
   \node[blue] at (1.35,-2.8) {\small $s_{i^\ast}$};  
\end{scope}   
 \end{tikzpicture}
\end{center}

Let us set 
$$
\begin{array}{ll}
A_{i_l}:= \Delta_i(\A', \A_l), ~~~~~~~~  &A_{j_l}:= \Delta_j(\A', \A_j), \\
A_{i_r}:= \Delta_{i}(\A', \A_r),~~~ &A_{j_r} := \Delta_{j}(\A', \A_r), \\
A_{i_e}: =\Lambda_{i^\ast} (h(\A_r, \A_l)), ~~~~~~~~~~ &A_{j_e} :=\Lambda_{j^\ast} (h(\A_r, \A_l)),\\
Y_i:=\Delta_{i}(\A', \A_2)=\Delta_{i}(\A', \A_3),  ~~~~~~ &Y_j:=\Delta_{j}(\A', \A_3)=\Delta_{j}(\A', \A_4),  \\
Z_i:=\Delta_{i}(\A', \A_{3'})=\Delta_i(\A', \A_{4'}),       &Z_j:= \Delta_{j}(\A', \A_{2'})= \Delta_{j}(\A', \A_{3'}).\\
\end{array}
$$

\begin{lemma}We have
\be
Z_iY_iY_j= A_{i_e}A_{j_r}Y_i^2+ A_{i_r}A_{j_e}\left(A_{i_e}A_{i_r}A_{j_l}+A_{i_l} Y_j\right),
\ee
\be
Z_j Y_i^2Y_j = A_{i_e}^2 A_{j_l}A_{j_r}Y_i^2+ \left(A_{i_e}A_{i_r}A_{j_l}+A_{i_l} Y_j\right)^2A_{j_e}.
\ee
\end{lemma}
\begin{remark}
The Lemma refines of \cite[Thm.1.16]{FZ98} by taking the extra frozen variables $A_{i_e}$ and $A_{j_e}$ into consideration. A direct calculation shows that it is compatible with the mutation sequence $\mu_j \circ \mu_i\circ \mu_j$ considered in the proof of Theorem \ref{thm.3.6.st}, case 3.
\end{remark}
\begin{proof} For $k\in {\rm I} -\{i,j\}$, we set $A_k =\Delta_k(\A', \A_l)$. Set
\[
P_i :=\prod_{k\in {\rm I}-\{i,j\}} A_{k}^{-C_{ik}}, \hskip 10mm P_j :=\prod_{k\in {\rm I}-\{i,j\}} A_{k}^{-C_{jk}}.
\]
By \eqref{bfz.lusztig},  the partial potentials corresponding to the first decomposition are
\be
\label{pot11015}
b_1={\displaystyle \frac{A_{i_e}}{A_{i_l}Y_i}} A_{j_l}P_i,  \hskip 7mm 
b_2= {\displaystyle \frac{1}{A_{j_l}Y_j}}Y_i^2P_j, \hskip 7mm 
b_3=\frac{1}{Y_iA_{i_r}}Y_j P_i,  \hskip 7mm 
b_4=\frac{A_{j_e}}{Y_{j}A_{j_r}}A_{i_r}^2 P_j. 
\ee
The partial potentials corresponding to the second decompositions are
\be
\label{pot11016}
b_1'={\displaystyle \frac{A_{j_e}}{A_{j_l}Z_j}} A_{i_l}^2P_j,  \hskip 7mm 
b_2'= {\displaystyle \frac{1}{A_{i_l}Z_i}}Z_jP_i, \hskip 7mm 
b_3'=\frac{1}{Z_jA_{j_r}}Z_i^2 P_j,  \hskip 7mm 
b_4'=\frac{A_{i_e}}{Z_{i}A_{i_r}}A_{j_r} P_i. 
\ee
Recall the formula
\be
\label{decmo s1s2s1s2}
x_{i}(b_1)x_j(b_2)x_i(b_3)x_j(b_4)=
x_j(b_2b_3^2b_4q^{-1})x_i(qp^{-1})x_j(p^2q^{-1})x_i(b_1b_2b_3p^{-1}),
\ee
where
\be
\label{pot11018}
p=b_1b_2+(b_1+b_3)b_4, ~~~q=b_1^2b_2+(b_1+b_3)^2b_4.
\ee
Note that the coordinates on the right hand side of \eqref{decmo s1s2s1s2} are $b_1', b_2', b_3', b_4'$ respectively. Hence
\be
\label{pot11019}
b_2'b_3'= p, ~~~~~~ (b_2')^2 b_3'= q.
\ee
Plugging \eqref{pot11015}, \eqref{pot11016}, and \eqref{pot11018} into \eqref{pot11019}, we prove the Lemma.
\end{proof}

\begin{itemize}
\item[] \hskip -1cm  {\bf Case 4:}  $(i, j, i, j, i, j)\sim (j, i, j, i, j, i)$, where $C_{ji}=-3$ and $C_{ij}=-1$. 
\end{itemize}
 
Let $u=s_is_js_is_js_is_j$. We have two different decompositions for ${\rm Conf}_u^e(\mathcal{A})$.

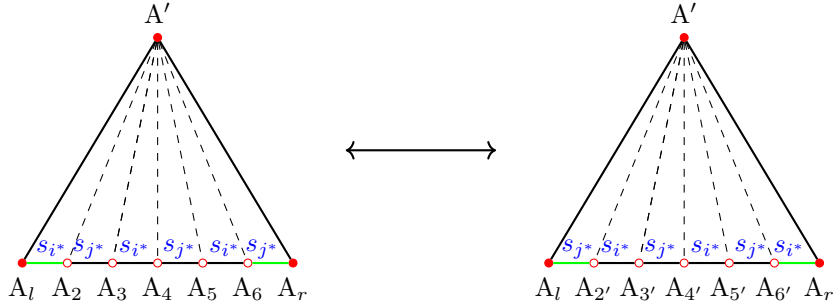
\begin{center}
\begin{tikzpicture}
\node [circle,draw=red,fill=red,minimum size=3pt,inner sep=0pt,label=above:{\small ${\A}'$}] (a) at (0,0) {};
\node [circle,draw=red,fill=red,minimum size=3pt,inner sep=0pt,label=below:{\small ${\A}_l$}] (b) at (-1.8,-3) {};
 \node [circle,draw=red,fill=white, minimum size=3pt,inner sep=0pt,label=below:{\small ${\A}_2$}]  (c) at (-1.2,-3){};
 \node [circle,draw=red,fill=white, minimum size=3pt,inner sep=0pt,label=below:{\small ${\A}_3$}] (d) at (-0.6,-3) {};
  \node [circle,draw=red,fill=white, minimum size=3pt,inner sep=0pt,label=below:{\small ${\A}_4$}] (e) at (0,-3) {};
  \node [circle,draw=red,fill=white, minimum size=3pt,inner sep=0pt,label=below:{\small ${\A}_5$}] (f) at (0.6,-3) {};
  \node [circle,draw=red,fill=white, minimum size=3pt,inner sep=0pt,label=below:{\small ${\A}_6$}] (g) at (1.2,-3) {};
 \node [circle,draw=red,fill=red,minimum size=3pt,inner sep=0pt,label=below:{\small ${\A}_r$}] (h) at (1.8,-3) {};
 \foreach \from/\to in {a/b,  c/d, d/e,  e/f, f/g, h/a}
                  \draw[thick] (\from) -- (\to);
 \foreach \from/\to in {b/c, g/h}
                  \draw[thick, green] (\from) -- (\to);   
\foreach \from/\to in {a/c,a/d,a/e,a/d, a/f, a/g}                   
              \draw[dashed] (\from) -- (\to);
                   
 \node[blue] at (-1.4,-2.8) {\small ${s_{i^\ast}}$};
  \node[blue] at (-.9,-2.82) {\small ${s}_{j^\ast}$}; 
   \node[blue] at (-.3,-2.8) {\small $s_{i^\ast}$};
   \node[blue] at (.3,-2.82) {\small $s_{j^\ast}$};
     \node[blue] at (.9,-2.8) {\small $s_{i^\ast}$};
   \node[blue] at (1.4,-2.82) {\small $s_{j^\ast}$};
  \draw[thick,<->] (2.5,-1.5) -- (4.5,-1.5);            

\begin{scope}[shift={(7,0)}]
\node [circle,draw=red,fill=red,minimum size=3pt,inner sep=0pt,label=above:{\small ${\A}'$}] (a) at (0,0) {};
\node [circle,draw=red,fill=red,minimum size=3pt,inner sep=0pt,label=below:{\small ${\A}_l$}] (b) at (-1.8,-3) {};
 \node [circle,draw=red,fill=white, minimum size=3pt,inner sep=0pt,label=below:{\small ${\A}_{2'}$}]  (c) at (-1.2,-3){};
 \node [circle,draw=red,fill=white, minimum size=3pt,inner sep=0pt,label=below:{\small ${\A}_{3'}$}] (d) at (-0.6,-3) {};
  \node [circle,draw=red,fill=white, minimum size=3pt,inner sep=0pt,label=below:{\small ${\A}_{4'}$}] (e) at (0,-3) {};
  \node [circle,draw=red,fill=white, minimum size=3pt,inner sep=0pt,label=below:{\small ${\A}_{5'}$}] (f) at (0.6,-3) {};
  \node [circle,draw=red,fill=white, minimum size=3pt,inner sep=0pt,label=below:{\small ${\A}_{6'}$}] (g) at (1.2,-3) {};
 \node [circle,draw=red,fill=red,minimum size=3pt,inner sep=0pt,label=below:{\small ${\A}_r$}] (h) at (1.8,-3) {};
 \foreach \from/\to in {a/b,  c/d, d/e,  e/f, f/g, h/a}
                  \draw[thick] (\from) -- (\to);
 \foreach \from/\to in {b/c, g/h}
                  \draw[thick, green] (\from) -- (\to);   
\foreach \from/\to in {a/c,a/d,a/e,a/d, a/f, a/g}                   
              \draw[dashed] (\from) -- (\to);
                   
 \node[blue] at (-1.4,-2.82) {\small ${s_{j^\ast}}$};
  \node[blue] at (-.9,-2.8) {\small ${s}_{i^\ast}$}; 
   \node[blue] at (-.3,-2.82) {\small $s_{j^\ast}$};
   \node[blue] at (.3,-2.8) {\small $s_{i^\ast}$};
     \node[blue] at (.9,-2.82) {\small $s_{j^\ast}$};
   \node[blue] at (1.4,-2.8) {\small $s_{i^\ast}$};
\end{scope}   
 \end{tikzpicture}
\end{center}

With respect to the first decomposition ${\bf i}=(i,j,i,j,i,j)$, we obtain a coordinate system $K_{\bf i}$ of ${\rm Conf}_u^v(\mathcal{A})$, which consists of 10 coordinates:
\be
\la{sysg2cor1}
\begin{array}{ll}
A_{i_l}:= \Delta_i(\A', \A_l), ~~~~~~~~  &A_{j_l}:= \Delta_j(\A', \A_j), \\
A_{i_r}:= \Delta_{i}(\A', \A_r),~~~ &A_{j_r} := \Delta_{j}(\A', \A_r), \\
A_{i_e}: =\Lambda_{i^\ast} (h(\A_r, \A_l)), ~~~~~~~~~~ &A_{j_e} :=\Lambda_{j^\ast} (h(\A_r, \A_l)),\\
Y_1:=\Delta_{i}(\A', \A_2)=\Delta_{i}(\A', \A_3),  ~~~~~~ &Y_3:=\Delta_{j}(\A', \A_3)=\Delta_{j}(\A', \A_4),  \\
Y_2:=\Delta_{i}(\A', \A_{4})=\Delta_i(\A', \A_{5}),       &Y_4:= \Delta_{j}(\A', \A_{5})= \Delta_{j}(\A', \A_{6}).\\
\end{array}
\ee
The coordinates are parametrized by the 10 vertices of the first quiver ${\bf J}({\bf i})$ shown in the proof of Theorem \ref{thm.3.6.st}, case 4: the coordinate $Y_i$ is assigned to the unfrozen vertex $i$, etc. 

With respect to the second decomposition ${\bf i}'=(j,i,j,i,j,i)$, we obtain a coordinate system $K_{{\bf i}'}$, which consists of the first 6 functions on \eqref{sysg2cor1} and 4 new unfrozen functions 
\be
\la{sysg2cor2}
\begin{array}{ll}
Z_1:=\Delta_{i}(\A', \A_{3'})=\Delta_{i}(\A', \A_{4'}),  ~~~~~~ &Z_3:=\Delta_{j}(\A', \A_{2'})=\Delta_{j}(\A', \A_{3'}),  \\
Z_2:=\Delta_{i}(\A', \A_{5'})=\Delta_i(\A', \A_{6'}),       &Z_4:= \Delta_{j}(\A', \A_{4'})= \Delta_{j}(\A', \A_{5'}).\\
\end{array}
\ee
The coordinates are parametrized by  the 10 vertices of the second quiver ${\bf J}({\bf i}')$ shown in the proof of Theorem \ref{thm.3.6.st}, case 4.

\begin{lemma} The transition from the coordinate system $K_{\bf i}$ to $K_{{\bf i}'}$ coincides with the cluster $K_2$-transformations under the sequence of mutations $(4, 3, 2, 1, 4, 2, 4, 3, 1, 4)$.
\end{lemma}
\begin{proof} Let us set $\A_l:= \A_1$ and $\A_r:= \A_7$. For $k\in \{1,\ldots, 6\}$,  consider the partial potential
\[
b_k:= \mathcal{W}_{i_k}(\A', \B_k, \B_{k+1}), \hskip 7mm \mbox{where } (i_1,\ldots,i_6)=(i,j,i,j,i,j).
\]
Putting them together, we obtain a unipotent element
\[
u=x_i(b_1)x_j(b_2)x_i(b_3)x_j(b_4)x_i(b_5)x_j(b_6).
\]
Since $x_i(b)= H_i(b){\bf E}_i H_i(b^{-1})$, 
 the element $u$ can be expressed using the coordinates in \cite{FG05}:
\be \la{FG3EFH}
 u:= H_i(X_{i_l})H_j(X_{j_l}){\bf E}_i H_i(X_1) {\bf E}_jH_j(X_3) {\bf E}_iH_i(X_2) {\bf E}_j H_j(X_4) {\bf E}_i {\bf E}_jH_i(X_{i_r}) H_j(X_{j_r})
\ee
where
\be
\la{sysg2poissoncor1}
\begin{array}{ll}
X_{i_l}:= b_1, ~~~~~~~~  &X_{j_l}:= b_2, \\
X_1:=b_3/b_1,  ~~~~~~ &X_3:=b_4/b_2,  \\
X_2:=b_5/b_3,       &X_4:= b_6/b_4,\\
X_{i_r}:= 1/b_5,~~~ &X_{j_r} := 1/b_6. \\
\end{array}
\ee
Similarly, with respect to the second decomposition $(i_1',\ldots, i_6')=(j,i,j,i,j,i)$ and the potential $b_k':= \mathcal{W}_{i_k'}(\A', \B_{k'}, \B_{k+1'})$, we obtain another decomposition
\[
 u:= H_i(X_{i_l}')H_j(X_{j_l}'){\bf E}_j H_j(X_3') {\bf E}_iH_i(X_1') {\bf E}_jH_j(X_4') {\bf E}_i H_i(X_2') {\bf E}_j {\bf E}_iH_i(X_{i_r}') H_j(X_{j_r}')
\]
where
\be
\la{sysg2poissoncor2}
\begin{array}{ll}
X_{i_l}':= b_2', ~~~~~~~~  &X_{j_l}':= b_1', \\
X_1':=b_4'/b_2',  ~~~~~~ &X_3':=b_3'/b_1',  \\
X_2':=b_6'/b_4',       &X_4':= b_5'/b_3',\\
X_{i_r}':= 1/b_6',~~~ &X_{j_r}' := 1/b_5'. \\
\end{array}
\ee
By \cite[Thm 3.5 (5)]{FG05}, the transition map from \eqref{sysg2poissoncor1} to \eqref{sysg2poissoncor2} coincides with the cluster Poisson transformation under the sequence $(4, 3, 2, 1, 4, 2, 4, 3, 1, 4)$. 
Note that the first 6 coordinates  in \eqref{sysg2poissoncor1} recovers the potential $b_1, \ldots, b_6$. Together with four of the  frozen coordinates, they form a new mixed coordinate system 
\[
P_{\bf i}:=\{A_{i_l}, A_{j_l}, A_{i_r}, A_{j_r}, X_{i_l}, X_{j_l}, X_1, X_2, X_3, X_4\}.
\]
Note that $P_{\bf i}$ is related to $K_{\bf i}$ by an invertible linear base change $p_{\bf i}$.

Similarly, for the  second decomposition ${\bf i}'$ we get a new mixed coordinate system
\[
P_{{\bf i}'}:=\{A_{i_l}, A_{j_l}, A_{i_r}, A_{j_r}, X_{i_l}', X_{j_l}', X_1', X_2', X_3', X_4'\}.
\]
It is related to $K_{{\bf i}'}$ by an invertible  linear base change $p_{{\bf i}'}$. 
Putting the systems together, we get a diagram
\begin{center}
\begin{tikzcd}
  K_{\bf i} \arrow[r, "\mu^{\mathscr{A}}" ] \arrow[d, "p_{\bf i}"]
    & K_{{\bf i}'} \arrow[d, "p_{{\bf i}'}" ] \\
  P_{\bf i} \arrow[r, "\mu^{\mathscr{X}}"]
&P_{{\bf i}'} \end{tikzcd}
\end{center}
Here $p_{\bf i}, p_{{\bf i}'}$ are essentially the $p$ maps relating $X-$and $A-$coordinates. By  \cite[Thm 3.5 (5)]{FG05}, $\mu^{\mathscr{X}}$ is the cluster Poisson transformation under the sequence $(4, 3, 2, 1, 4, 2, 4, 3, 1, 4)$. To make the above diagram commute, $\mu^{\mathscr{A}}$ must be the cluster $K_2$ transformation under the  same sequence. 
\end{proof}

\subsection{Proof of Theorem \ref{Th2.3}, part 3}
\label{partial.conf.part3}
\medskip

Consider the {\it double Bruhat cell}
\[
\G_{u,v}:= \B^- u \B^{-} \cap \B v \B \subset \G.
\]
Below we rapidly recall the cluster Poisson structure of $\G_{u,v}$  introduced in \cite{FG05}. 

The standard pinning $([\U], [\U^-])$ determines a datum  $({\rm H}, {\rm B}, x_i,  y_i; i\in {\rm I})$ of $\G$. We set
\[
{E}_i := x_i(1), \hskip 12mm {E}_{\overline{i}}:= y_i(1), ~~~~i\in {\rm I}.
\]
Let $H_1,\ldots, H_r$ be the fundamental coweights of $\G$. Abusing notation, we also write $H_{\overline{i}}=H_i$. We have
\be
\la{2017.5.5.18.54}
x_i(a)= H_i(a){E}_i H_i(a^{-1}), \hskip 8mm y_i(b)= H_i(b^{-1}){E}_{\overline{i}} H_i(b). 
\ee
For $i\neq j$, we have 
$$
{E}_i H_j(a) = H_j(a) {E}_i, \hskip 8mm {E}_{\overline{i}} H_{j}(b) = H_j(b) {E}_{\overline{i}}.
$$

Every reduced word  ${\bf i}=(i_1,\ldots, i_m)$  of $(u,v)$ gives rise to a  quiver ${\bf J}({\bf i})$ in Definition \ref{quiver.amla.fg4.sec22} with cluster Poisson coordinates $\{X_j\}$.  There is an open embedding map
\begin{align}
ev_{\bf i}: ~~ {\Bbb G}_m^{r+m} &\lra  \G_{u,v}, \nonumber \\
 \{X_j\} &\lms g= \left(\prod_{i\in {\rm I}} H_i(X_{i \choose 0})\right)\cdot E_{i_1} H_{i_1}(X_1)\cdots E_{i_m} H_{i_m}(X_m).  \la{g.decomposition}
\end{align}

For any reduced words ${\bf i}$, ${\bf i}'$ of $(u,v)$,  the transition map $ev_{{\bf i}'}^{-1}\circ ev_{\bf i}$ is a cluster Poisson transformation, which is deduced from the following  identities ({\it cf.} \cite[Proposition 3.6]{FG05}):
\[
{E}_{\overline{i}} H_{i}(X) {E}_i = \Big(\prod_{j\neq i} H_j(1+X)^{-C_{ji}}\Big) H_i(1+X^{-1})^{-1} {E}_i H_i(X^{-1}){E}_{\overline{i}} H_i(1+X^{-1})^{-1}.
\]
If $C_{ij}=C_{ji}=-1$, then 
\[
{E}_i H_i(X) { E}_j {E}_i = H_i(1+X) H_{j}(1+X^{-1})^{-1} {E}_j H_j(X^{-1}) {E}_i {E}_j  H_j(1+X) H_{i}(1+X^{-1})^{-1}. 
\]
See similar identities for non-simply laced cases in {\it loc.cit.} The quotient space ${\rm H}\backslash {\rm G}_{u,v}/ {\rm H}$ inherits a cluster Poisson structure from ${\rm G}_{u,v}$ by deleting all the frozen vertices. 

\vskip 2mm

Let $\mathcal{P}_{u}^v$ be  the moduli space  parametrizing $\G$-orbits in ${\rm Conf}_u^v(\mathcal{B})$ plus a pinning $p_l$ over its  left pair $(\B^l, \B_l)$ and a pinning $p_r$ over the its right pair $(\B^r, \B_r)$.
\bl 
\la{2017.4.18.21.18hh}
There is a natural isomorphism
\[
\G'_{u,v} \stackrel{=}{\lra} \mathcal{P}_{u}^v, \hskip 9mm g \lms \left(p_l= ([\U], [\U^-]), ~p_r= (g\cdot [\U], g \cdot [\U^-])\right)
\]
\el 
\begin{proof} For each $\G'$-orbit in $\mathcal{P}_{u}^v$, there is a unique representative such that the left pinning $p_l=([\U], [\U^-])$. For such a representative, there is a unique element $g\in \G'$ taking $p_l$ to $p_r$.
The Lemma follows from the fact that  $w_-([\U^-], g \cdot [\U^-])=u$ if and only if $g\in \B^- u \B^-$, and $w([\U], g\cdot [\U])=v$ if and only if $g\in \B v\B$.
\end{proof}

There is a natural projection  $\mathcal{P}_{u}^v \lra {\rm Conf}_u^v(\mathcal{B})$  forgetting pinnings. The following diagram commutes
\be
\begin{gathered}
 \xymatrix{
        \G'_{u, v}\ar[r]^{=} \ar[d]& \mathcal{P}_u^v \ar[d] \\
        {\rm H}\backslash {\rm G}'_{u,v}/ {\rm H} \ar[r]^{=} & {\rm Conf}_u^v(\mathcal{B})
 }
\end{gathered}
\ee
The isomorphisms in diagram  induce cluster Poisson structures on  $\mathcal{P}_u^v$ and  ${\rm Conf}_u^v({\cal B})$,   concluding the proof of Theorem \ref{Th2.3}, part 3.

\subsection{Proof of Theorem \ref{Th2.3}, part 4}
\medskip

Let ${\bf i}$ and ${\bf i}'$ be reduced words of $(u, v)$. By  part 2 of Theorem \ref{Th2.3}, the transition map $\alpha_{{\bf i}, {\bf i}'}:= c^{-1}_{\mathscr{A},{\bf i}'}\circ c_{\mathscr{A}, {\bf i}}$ is a cluster $K_2$-transformation. By part 3, the transition map $\chi_{{\bf i}, {\bf i}'}:= c^{-1}_{\mathscr{X},{\bf i}'}\circ c_{\mathscr{X}, {\bf i}}$ is a cluster Poisson transformation encoded by the same sequence of cluster mutations and permutations as $\alpha_{{\bf i}, {\bf i}'}$.
As a general property of cluster ensemble, the following diagram commutes:
\be
\begin{gathered}
 \xymatrix{
        {\rm T}^{\mathscr A}_{\bf i} \ar[d]_{p_{\bf i}}\ar[r]^{\alpha_{{\bf i}, {\bf i}'}} & {\rm T}^{\mathscr A}_{{\bf i}'} \ar[d]^{{p_{\bf i}'}} \\
        {\rm T}^{\mathscr X}_{\bf i} \ar[r]_{\chi_{ {\bf i}, {\bf i}'}} & {\rm T}^{\mathscr X}_{{\bf i}'}
 }
\end{gathered}
\ee
Therefore it suffices to prove part 4 of the Theorem for one specific reduced word ${\bf i}$. 

Let ${\bf i}=(i_1, \ldots, i_p, \overline{i}_{p+1}, \ldots, \overline{i}_n)$ be a reduced word of $(u, v)$ such that
\[
i_1, \ldots, i_p \in {\rm I}, \hskip 10mm \overline{i}_{p+1}, \ldots, \overline{i}_{n}\in \overline{\rm I}.
\]

Take  $(\A_l, \A_r, \A^l, \A^r)\in {\rm Conf}_u^v(\mathcal A)$ such that every pair of decorated flags connected by a dashed line in the following decomposition is of generic position: 
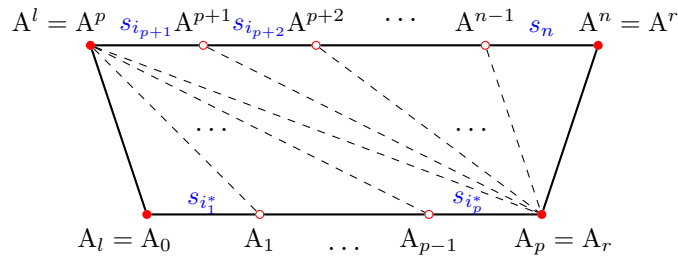
\begin{center}
\begin{tikzpicture}[scale=0.75]
\node [circle,draw=red,fill=red,minimum size=3pt,inner sep=0pt, label=above:{\small ${{\A}^l=\A^p}\qquad$}] (a1) at (-2,0) {};
\node [circle,draw=red,fill=white,minimum size=3pt,inner sep=0pt,label=above:{\small ${\A}^{p+1}$}] (a2) at (0,0) {};
\node [circle,draw=red,fill=white,minimum size=3pt,inner sep=0pt,label=above:{\small ${\A}^{p+2}$}] (a3) at (2,0) {};
\node [circle,draw=red,fill=white,minimum size=3pt,inner sep=0pt,label=above:{\small ${\A}^{n-1}$}] (a4) at (5,0) {};
\node [circle,draw=red,fill=red,minimum size=3pt,inner sep=0pt, label=above:{\small $\qquad{{\A}^{n}=\A^r}$}] (a5) at (7,0) {};
 \node [circle,draw=red,fill=red,minimum size=3pt,inner sep=0pt, label=below:{\small ${{\A}_{l}=\A_0}\qquad$}]  (b1) at (-1,-3){};
 \node [circle,draw=red,fill=white,minimum size=3pt,inner sep=0pt, label=below:{\small $\A_1$}] (b2) at (1,-3) {};
 \node [circle,draw=red,fill=white,minimum size=3pt,inner sep=0pt,label=below:{\small ${\A}_{p-1}$}] (b3) at (4,-3) {};
  \node [circle,draw=red,fill=red,minimum size=3pt,inner sep=0pt,label=below:{\small ${\qquad{\A}_{p}=\A_r}$}] (b4) at (6,-3) {};
 \foreach \from/\to in {a1/b2, a1/b3, a1/b4, a2/b4, a3/b4, a4/b4}
                   \draw[dashed] (\from) -- (\to);
  \foreach \from/\to in {a1/b1, a1/a2, a2/a3, a3/a4, a4/a5, b1/b2, b2/b3, b3/b4, a5/b4}
                   \draw[thick] (\from) -- (\to);                            
 \node[blue] at (0,-2.8) {\small $s_{i_1^\ast}$};
 \node[] at (.15, -1.5) {$\ldots$};
 \node[] at (4.75, -1.5) {$\ldots$};
  \node[blue] at (4.7,-2.8) {\small $s_{i_p^\ast}$};
  \node[blue] at (-1,0.3) {\small $s_{i_{p+1}}$};
    \node[blue] at (1,0.3) {\small $s_{i_{p+2}}$};
 \node at (2.5, -3.6) {$\ldots$};
  \node at (3.5, .45) {$\ldots$};
   \node[blue] at (6,0.3) {\small $s_{{n}}$};  
 \end{tikzpicture}
 \end{center}

Let us consider the potential function for each angle
\be
b_s=\left\{ \begin{array}{ll}
 {\cal W}_{i_s}(\A^l, \B_{s-1}, \B_{s}), ~~~~  &\mbox{if } 1 \leq s \leq p;  \\
 {\cal W}_{i_{s}^*}(\A_r, \B^{s}, \B^{s-1}),~~~~~& \mbox{if } p < s \leq n.
 \end{array}\right.
 \ee
 Let $pr$ be the natural projection from  $\G'$ to $\G$. Set
 \[
 h= pr(h(\A^l, \A_r)^{-1})= H_{1}(c_1)\ldots H_{r}(c_r)\in \G.
 \]
Define
 \[
 g= x_{i_1}(b_1)\ldots x_{i_p}(b_p)\cdot h \cdot y_{i_{p+1}}(b_{p+1})\ldots y_{i_{n}}(b_{n}) \in \G.
 \]
 By \eqref{2017.5.5.18.54}, we have
 \[
 g=  \prod_{s=1}^p H_{j_s}(b_s) {E}_{j_s} H_{j_s}(b_s^{-1})   \cdot \prod_{i\in {\rm I}} H_{i}(c_i) \cdot  \prod_{s=p+1}^n H_{i_s}(b_{s}^{-1}) {E}_{\overline{i}_{s}} H_{i_{s}}(b_{s}).
 \]
By the construction of $g$, we have
\[
(\B^l, \B_l, \B^r, \B_r) = (\B, \B^-, g\cdot \B, g\cdot \B^-).
\]   

Let $i\in {\rm I}$. 
Recall the alternative notations \eqref{alternative.notation.J}. For each ${i \choose l}\in J_{\rm uf}({\bf i})$, we get
\be
\la{Xv}
p^*X_{i \choose l}   =\left\{  \begin{array}{ll}
    {\displaystyle { b_{i \choose {l+1}}}{  b_{i \choose l}}^{-1}} \hskip 7mm &\mbox{if }  {i \choose l+1} \leq p, \\
 \\  {\displaystyle   c_i b_{i \choose l}^{-1} b_{i \choose {l+1}}^{-1}} & \mbox{if } {i \choose l}\leq p < {i \choose {l+1}}, \\
 \\  {\displaystyle  b_{i \choose {l}} b_{i \choose {l+1}}^{-1}} & \mbox{if } p<{i \choose {l}}.\\
   \end{array}\right.
\ee
Note that $b_s$ are Laurent monomials of $A$-variables of ${\rm Conf}_u^v({\cal A})$ given by \eqref{bfz.lusztig}$\&$\eqref{bfz.lusztig2t}, and
\[
c_i = \alpha_i (h(\A^l, \A_r)^{-1}) =\prod _{j\in I} \Lambda_j(h(\A^l, \A_r))^{-C_{ij}}.
\] 
Plugging them back to \eqref{Xv}, a careful comparison shows that $p^*X_{i \choose l}$ is presented by the a Laurent monomial of $A$-variables as in \eqref{2017.4.19.14.39tt}, which concludes the proof of Theorem \ref{Th2.3}, part 4.

\begin{example} Let $\G={\rm PGL}_3$ and $u=v=w_0$. Consider the word ${\bf i}=(2,1,2, \overline{1}, \overline{2}, \overline{1})$. Its quiver and cluster variables are illustrated below. 
\begin{center}
\begin{tikzpicture}[scale=1.3]
\draw [gray=20!] (-.1, 3.1) -- (3.1, 3.1) -- (3.1, -.1) -- (-.1, 3.1) -- (-.1, -.1) -- (3.1, -.1);
\node (A1) at (1, 3) {$A_1$}; 
\node (A2) at (2,3) {$A_2$};
\node (A3) at (0, 2) {$A_3$}; 
\node (A4) at (1,2) {$A_4$};
\node (A5) at (2, 2) {$A_5$}; 
\node (A6) at (3,2) {$A_6$};
\node (A7) at (0, 1) {$A_7$}; 
\node (A8) at (1,1) {$A_8$};
\node (A9) at (2, 1) {$A_9$}; 
\node (A10) at (3,1) {$A_{10}$};
\node (A11) at (1, 0) {$A_{11}$}; 
\node (A12) at (2,0) {$A_{12}$};
\foreach \from/\to in {A1/A4, A4/A5, A5/A1, A2/A5, A5/A6, A6/A2, A5/A9, A9/A10, A10/A5, A4/A3, A3/A8, A8/A4, A8/A7, A7/A11, A11/A8, A9/A8, A8/A12, A12/A9}
\draw[directed] (\from) --(\to);
\foreach \from/\to in {A1/A2, A6/A10, A7/A3, A12/A11}
\draw[dashed, directed] (\from)--(\to);

\begin{scope}[xshift=6cm]
\draw [gray=20!] (-.1, 3.1) -- (3.1, 3.1) -- (3.1, -.1) -- (-.1, 3.1) -- (-.1, -.1) -- (3.1, -.1);
\node (X1) at (1, 3) {}; 
\node (X2) at (2,3) {};
\node (X3) at (0, 2) {}; 
\node (X4) at (1,2) {$X_{1 \choose 1}$};
\node (X5) at (2, 2) {$X_{1 \choose 2}$}; 
\node (X6) at (3,2) {};
\node (X7) at (0, 1) {}; 
\node (X8) at (1,1) {$X_{2 \choose 1}$};
\node (X9) at (2, 1) {$X_{2 \choose 2}$}; 
\node (X10) at (3,1) {};
\node (X11) at (1, 0) {}; 
\node (X12) at (2,0) {};
\foreach \from/\to in {X1/X4, X4/X5, X5/X1, X2/X5, X5/X6, X6/X2, X5/X9, X9/X10, X10/X5, X4/X3, X3/X8, X8/X4, X8/X7, X7/X11, X11/X8, X9/X8, X8/X12, X12/X9}
\draw[directed] (\from) --(\to);
\foreach \from/\to in {X1/X2, X6/X10, X7/X3, X12/X11}
\draw[dashed, directed] (\from)--(\to);
\end{scope}
\end{tikzpicture}
\end{center}
By Lemma \ref{bfz.lusztig.lemma}, we get 
\[
b_{1 \choose 1}=\frac{A_8}{A_3A_4}, \qquad b_{1 \choose 2}=\frac{A_1A_9}{A_4A_5}, \qquad  b_{1 \choose 3}=\frac{A_2A_{10}}{A_5A_6},
\]
\[
b_{2 \choose 1}=\frac{A_{11}A_3}{A_7A_8}, \qquad  b_{2 \choose 2}=\frac{A_{12}A_4}{A_8A_9}, \qquad b_{2 \choose 3}=\frac{A_5}{A_9A_{10}}. 
\]
Note that
$\Lambda_1(h(\A^l, \A_r))=A_4$ and  $\Lambda_2(h(\A^l, \A_r))=A_9$.
Therefore
\[
c_1=\frac{A_9}{A_4^2}, \qquad c_2=\frac{A_4}{A_9^2}.
\]
We get 
\[
p^\ast X_{1 \choose 1}= c_1 b_{1 \choose 1}^{-1} b_{1\choose 2}^{-1}=\frac{A_3A_5}{A_1A_8}, \qquad p^\ast X_{1 \choose 2}= b_{1 \choose 2}b_{1\choose 3}^{-1}=\frac{A_1A_6A_9}{A_2A_4A_{10}},
\]
\[
p^\ast X_{2 \choose 1} = b_{2 \choose 2}b_{2 \choose 1}^{-1} = \frac{A_4A_7A_{12}}{A_3A_9A_{11}}, \qquad 
p^\ast X_{2 \choose 2} = c_2 b_{2 \choose 2}^{-1} b_{2 \choose 3}^{-1} =\frac{A_8A_{10}}{A_5A_{12}}.
\]

\end{example}

 \medskip
 
 \section{Cluster $K_2$-structure of the moduli spaces ${\rm Conf}_3^\times({\cal A})$ and ${\mathscr A}_{\G, \bS}$}
\la{sec.conf3}

\medskip

In Section \ref{sec.conf3} we again assume first that all boundary intervals of $\bS$ are colored. The passage to the case when only some of them are colored is trivial: we just drop the 
frozen coordinates and the frozen vertices assigned to the non-colored intervals.


The space ${\rm Conf}_3^\times({\cal A})$ parametrizes $\G'$-orbits of triples $(\A_1, \A_2,  \A_3)\in {\cal A}^3$ where each pair $(\A_i, \A_{j})$  is in generic position.
It inherits a cluster $K_2$-structure from ${\rm Conf}^e_{w_0}({\cal A})$ by the isomorphism
\be
\la{2017.2.26.10.09hh}
\begin{aligned}
\iota: ~{\rm Conf}_3^\times({\cal A}) &\stackrel{\sim}{\lra} {\rm Conf}^e_{w_0}({\cal A}), \\
~~~~~(\A_1, \A_2, \A_3) & \lms (\A_l, \A_r, \A^l, \A^r)=(\A_2, \A_3, \A_1, \A_1).
\end{aligned}
\ee

Elaborating this, we assign $\A_1, \A_2, \A_3$ to vertices of a triangle in counterclockwise order. According to Lemma \ref{decomp.partial.dec.flag.15.02}, every reduced word ${\bf i}$ of $w_0$ gives rise to a further decomposition of the bottom pair $(\A_2, \A_3)$, as illustrated on the left of Figure \ref{2017.2.23.9.46h}.  Using the machinery of amalgamation \cite{FG05}, we introduce a cluster structure to ${\rm Conf}_3^\times({\cal A})$, as constructed in Section \ref{pcs.sec}. In particular, we obtain an open torus embedding for each reduced word ${\bf i}$ of $w_0$
\be
\la{cai}
c_{{\mathscr A}, {\bf i}}: {\rm T}_{\bf i}^{\mathscr A}:=({\Bbb G}_m)^{l(w_0)+2r}~ \hookrightarrow{~} {\rm Conf}_3^{\times}({\cal A}).
\ee
Theorem \ref{Th2.3} asserts that the cluster structure is independent of ${\bf i}$ chosen.
\vskip 2mm

There is an involution of $\G'$
\be
\la{involution.ast.tt}
\ast: \G' \lra \G', ~~~~{E}_{i}\lms{E}_{i^\ast}, ~~~~{E}_{\overline{i}}\lms {E}_{\overline{i}^\ast}~~~~ h \lms w_0(h^{-1}).
\ee
It induces an involution of ${\cal A}$ still denoted by $\ast$. 
Define the twisted map
\be
\la{twisted.map.a}
\eta:  {\rm Conf}_3^{\times}({\cal A})\stackrel{\sim}{\lra} {\rm Conf}_3^\times({\cal A}), \hskip 7mm (\A_1, \A_2, \A_3)\lms (\ast \A_3\cdot s_\G, \ast \A_1, \ast \A_2).
\ee
\bt
\la{main.result2.cluster.cyclic}
The map $\eta$ is a cluster $K_2$-automorphism of $ {\rm Conf}_3^{\times}({\cal A})$. 
\et
We prove Theorem \ref{main.result2.cluster.cyclic} by presenting an explicit sequence of cluster mutations and permutations that gives rise to $\eta$.
The proof  occupies the rest of this section. The most difficult part of the proof is the one concerning the frozen variables.

\subsection{The twisted map of ${\rm Conf}_3^\times({\cal B})$}
\la{sec.conf3B}
\medskip

The space
$
{\rm Conf}_{3}^\times({\cal B})
$
consists of $\G$-orbits of triples $(\B_1, \B_2, \B_3)\in {\cal B}^3$ where each pair of flags is in generic position.
It  admits a cluster Poisson structure induced by the isomorphism 
\be
\la{2017.3.19.10.09hh}
\begin{aligned}
\iota: ~{\rm Conf}_3^\times({\cal B}) &\lra {\rm Conf}^e_{w_0}({\cal B}), \\
~~~~~(\B_1, \B_2, \B_3)&\lms (\B_l, \B_r, \B^l, \B^r)=(\B_2, \B_3, \B_1, \B_1).
\end{aligned}
\ee
Let us consider the {\it twist map}
\be
\la{eta.dt.conf3b}
{\eta}:  {\rm Conf}_3^{\times}({\cal B})\stackrel{\sim}{\lra} {\rm Conf}_3^\times({\cal B}), \hskip 7mm (\B_1, \B_2, \B_3)\lms (\ast \B_3, \ast \B_1, \ast \B_2).
\ee
\bt 
\la{b.cluster.poisson.1541}
The map $\eta$ is a cluster Poisson automorphism of ${\rm Conf}_3^\times({\cal B})$. 
\et

\begin{proof}
Let ${\bf i}=(i_1, \ldots, i_m)$  be a reduced word of $w_0$.  The pairs $(\B_1, \B_2)$ and  $(\B_2, \B_3)$ admit decompositions
\be
\begin{gathered}
\B_1=\B_0'' \xrightarrow{~s_{i_1}~} \B_1''\xrightarrow{~s_{i_2}~}\cdots \xrightarrow{~s_{i_m}~} \B_m''=\B_2 , \\
\B_2=\B_0' \xrightarrow{~s_{i_1^\ast}~} \B_1'\xrightarrow{~s_{i_2^\ast}~}\cdots \xrightarrow{~s_{i_m^\ast}~} \B_m'=\B_3. 
\end{gathered}
\ee
The decompositions are illustrated by the following figure.
\begin{center}
\begin{tikzpicture}[scale=1.5]
\node [circle,draw=red,fill=red,minimum size=3pt,inner sep=0pt, label=above:{\small ${\B}_1$}] (a) at (0,0) {};
\node [circle,draw=red,fill=red,minimum size=3pt,inner sep=0pt, label=below:{\small ${\B}_2$}] (b) at (-1.8,-3) {};
 \node [circle,draw=red,fill=white, minimum size=3pt,inner sep=0pt, label=below:{\small ${\B}_2'$}]  (c) at (-0.36,-3){};
  \node [circle,draw=red,fill=white, minimum size=3pt,inner sep=0pt, label=below:{\small ${\B}_1'$}]  (g) at (-1.08,-3){};
 \node [circle,draw=red,fill=white, minimum size=3pt,inner sep=0pt, label=below:{\small ${\B}_{m-1}'$}] (h) at (1.08,-3) {};
 \node [circle,draw=red,fill=red,minimum size=3pt,inner sep=0pt, label=below:{\small ${\B}_3$}] (e) at (1.8,-3) {};
  \node [circle,draw=red,fill=white, minimum size=3pt,inner sep=0pt, label=left:{\small ${\B}_1''$}]  (c1) at (-0.36,-0.6){};
  \node [circle,draw=red,fill=white, minimum size=3pt,inner sep=0pt, label=left:{\small ${\B}_2''$}]  (g1) at (-0.72,-1.2){};
 \node [circle,draw=red,fill=white, minimum size=3pt,inner sep=0pt, label=left:{\small ${\B}_{m-1}''$}] (h1) at (-1.44,-2.4) {};
 \node [label=below:$\cdots$] at (0.36, -3) {};
\node [rotate=59.03] at (-1.38, -1.8) {$\cdots$};
 \foreach \from/\to in {b/g, g/c, h/e, a/e, a/c1, c1/g1, h1/b}
                  \draw[thick] (\from) -- (\to);
\foreach \from/\to in {c/h, g1/h1}
                  \draw[thick] (\from) -- (\to);
\foreach \from/\to in {c1/g, g1/c, h1/h}                   
              \draw[dashed] (\from) -- (\to);            
 \node[blue] at (-1.44,-3.14) {\small $s_{i_1^\ast}$};
  \node[blue] at (-0.72,-3.14) {\small ${s}_{i_2^\ast}$}; 
   \node[blue] at (1.44,-3.14) {\small $s_{i_m^\ast}$};          
    \node[blue] at (-0.38,-0.3) {\small $s_{i_1}$};
  \node[blue] at (-0.76,-0.9) {\small ${s}_{i_2}$}; 
   \node[blue] at (-1.82,-2.7) {\small $s_{i_m}$};        
    \end{tikzpicture}
 \end{center}
Set
\be
\la{2017.2.24.14.43sh}
v_k:= s_{i_k}s_{i_{k-1}}\ldots s_{i_1}, \hskip 7mm u_k :=s_{i_{k+1}}\ldots s_{i_m}.
\ee
There is a sequence of isomorphisms - rotations of dashed segments in the above figure:
\be
\la{comp.pois}
\begin{split}
&{\bf R}: ~{\rm Conf}^e_{w_0}({\cal B}) \xrightarrow{~r_1~} {\rm Conf}_{u_1}^{v_1}({\cal B})\xrightarrow{~r_2~} {\rm Conf}_{u_2}^{v_2}({\cal B})\xrightarrow{~r_3~}\cdots\xrightarrow{~r_m~}{\rm Conf}_{e}^{w_0}({\cal B}), \\
 &(\B_2,\B_3, \B_1, \B_1)\stackrel{r_1}{\lms}(\B_1', \B_3, \B_1'', \B_1)\stackrel{r_2}{\lms}(\B_2', \B_3, \B_2'', \B_1) \stackrel{r_3}{\lms} \cdots \stackrel{r_m}{\lms} (\B_3, \B_3, \B_2, \B_1)\\
\end{split}
 \ee
where the isomorphism $r_k: {\rm Conf}_{u_{k-1}}^{v_{k-1}}({\cal B}) \lra {\rm Conf}_{u_k}^{v_k}({\cal B})$ is shown on Figure \ref{2017.3.23.15.17sh}.
\begin{figure}[ht]
\epsfxsize150pt
\center{
 \begin{tikzpicture}[scale=0.8]
  \path [fill=gray!30] (-8,-3) --(-7,0)-- (-7,-3) -- (-8,-3);         
   \path [fill=gray!30] (-1,0) --(0,0)-- (0,-3) -- (-1,0);         
 \node [circle,draw=red,fill=red,minimum size=3pt,inner sep=0pt,label=below:{\small ${\B}_{k-1}'$}] (a) at (-8,-3) {};
\node [circle,draw=red,fill=red,minimum size=3pt,inner sep=0pt,label=above:{\small ${\B}_{k-1}''$}] (e) at (-7,0) {};
\node [circle,draw=red,fill=red,minimum size=3pt,inner sep=0pt,label=above:{\small ${\B}_1$}] (b) at (-5.5,0) {};
 \node [circle,draw=red,fill=white,minimum size=3pt,inner sep=0pt,label=below:{\small ${\B}_k'$}]  (d) at (-7,-3){};
 \node [circle,draw=red,fill=red,minimum size=3pt,inner sep=0pt,label=below:{\small ${\B}_3$}] (c) at (-5.5,-3) {};
 \foreach \from/\to in {a/e, e/b, b/c, d/c, d/a}
                  \draw[thick] (\from) -- (\to);
\foreach \from/\to in {e/d}
                  \draw[dashed] (\from) -- (\to);                     
  \draw[thick,->] (-4,-1.5) -- (-2,-1.5);   
               \node at (-3,-1.3) {$r_k$};         
\node [circle,draw=red,fill=red,minimum size=3pt,inner sep=0pt,label=above:{\small ${\B}_k''$}] (a) at (-1,0) {};
\node [circle,draw=red,fill=white,minimum size=3pt,inner sep=0pt,label=above:{\small ${\B}_{k-1}''$}] (e) at (0,0) {};
\node [circle,draw=red,fill=red,minimum size=3pt,inner sep=0pt,label=above:{\small ${\B}_1$}] (b) at (1.5,0) {};
 \node [circle,draw=red,fill=red,minimum size=3pt,inner sep=0pt,label=below:{\small ${\B}_k'$}]  (d) at (0,-3){};
 \node [circle,draw=red,fill=red,minimum size=3pt,inner sep=0pt,label=below:{\small ${\B}_3$}] (c) at (1.5,-3) {};
 \foreach \from/\to in {a/e, e/b, b/c, d/c, d/a}
                  \draw[thick] (\from) -- (\to);
\foreach \from/\to in {e/d}
                  \draw[dashed] (\from) -- (\to);                               
 \end{tikzpicture}
  }
\caption{The reflection map $r_k$.}
\label{2017.3.23.15.17sh}
\end{figure}

Let us consider the following reduced words of $(u_k, v_k)$
\be
\la{switch.seed.kk}
{\bf i}_k:= (\overline{i}_{k}, i_{k+1}, \ldots, i_{m},~\overline{i}_{k-1}, \ldots, \overline{i}_1),\hskip 10mm
{\bf i}_k':= (i_{k+1}, \ldots, i_{m},~\overline{i}_k, \ldots, \overline{i}_1).
\ee
Note that ${\bf i}_{k}$ is obtained from ${\bf i}_{k-1}'$ by changing the first index $i_k$ to $\overline{i}_k$. The unfrozen subquiver ${\bf J}_{\rm uf}({\bf i}_{k-1}')$   coincides with ${\bf J}_{\rm uf}({\bf i}_{k})$. In terms of  decomposition \eqref{g.decomposition}, the isomorphism $r_k$ is equivalent to a change of the first ${E}_{i_k}$ to ${E}_{\overline{i}_k}$, which keeps the unfrozen $X$-coordinates intact. So we get a   commutative diagram
\be
\la{twist.x.poisson}
\begin{gathered}
 \xymatrix{
        {\rm T}^{\mathscr X}_{{\bf i}_{k-1}'} \ar[d]_{id}\ar[r]^{c_{{\mathscr X},{\bf i}_{k-1}'}~~~~} &{\rm Conf}_{u_{k-1}}^{v_{k-1}}({\cal B}) \ar[d]_{r_k} \\
        {\rm T}^{\mathscr X}_{{\bf i}_k} \ar[r]_{c_{{\mathscr X}, {\bf i}_k}~~~} &{\rm Conf}_{u_{k}}^{v_{k}}({\cal B})
 }
\end{gathered}
\ee
In other words,  $r_k$ identifies the cluster Poisson structures on both spaces. 

For ${k}\in \{1,2,\ldots, m\}$, denote by $t_k$ the number of occurrences of $i={i_k}$ in the sub word
$( i_{k+1}, i_{k+2}, \ldots,  i_m)$.  
Recall the notations  in \eqref{alternative.notation.J}. By the proof of Theorem \ref{thm.3.6.st}, the sequence 
\be
L_k({\bf i}):=\mu_{i_k \choose t_k}\circ \ldots \circ \mu_{i_k\choose 2} \circ  \mu_{i_k \choose 1}
\ee
shifts the first $\overline{i}_k$ to the right, and takes the quiver ${\bf J}_{\rm uf}({\bf i}_k)$ to ${\bf J}_{\rm uf}({\bf i}_{k}')$. By  Part 3 of Theorem \ref{Th2.3}, the transition map $c_{{\mathscr X}, {\bf i}_{k}'}^{-1}\circ c_{{\mathscr X}, {\bf i}_{k}}$ is the cluster Poisson transformation presented by $L_k({\bf i})$:
\be
\la{twist.x.poisson.lk}
\begin{gathered}
 \xymatrix{
        {\rm T}^{\mathscr X}_{{\bf i}_{k}} \ar[d]_{L_{k}({\bf i})}\ar[r]^{c_{{\mathscr X},{\bf i}_{k}}~~~~} &{\rm Conf}_{u_{k}}^{v_{k}}({\cal B})  \\
        {\rm T}^{\mathscr X}_{{\bf i}_k'} \ar[ur]_{c_{{\mathscr X}, {\bf i}_k'}~~~} &
 }
\end{gathered}
\ee
Set $\overline{\bf i}={\bf i}_m= (\overline{i}_m, \ldots, \overline{i}_1)$. 
Combining \eqref{twist.x.poisson} and \eqref{twist.x.poisson.lk} for all $k$, we get 
\be
\begin{gathered}
 \xymatrix{
        {\rm T}^{\mathscr X}_{{\bf i}} \ar[d]_{L_m({\bf i})\circ \cdots \circ L_1({\bf i})}\ar[r]^{c_{{\mathscr X},{\bf i}}~~~~} &{\rm Conf}_{w_0}^{e}({\cal B}) \ar[d]^{\bf R} \\
        {\rm T}^{\mathscr X}_{\overline{\bf i}} \ar[r]_{c_{{\mathscr X}, \overline{\bf i}}~~~} &{\rm Conf}_e^{w_0}({\cal B})
 }
\end{gathered}
\ee

Lastly, let us consider the following isomorphism
\[
{\bf S}: {\rm Conf}_e^{w_0}({\cal B}) \lra {\rm Conf}_{w_0}^e({\cal B}), \hskip 7mm (\B_3, \B_3, \B_2, \B_1) \lra (\ast \B_1, \ast \B_2, \ast \B_3, \ast \B_3).
\]
\begin{figure}[ht]
\epsfxsize 200pt
\center{
\begin{tikzpicture}[scale=1.1]
\node [circle,draw=red,fill=red,minimum size=3pt,inner sep=0pt, label=above:{\small $\ast{\B}_3$}] (a) at (0,0) {};
\node [circle,draw=red,fill=red,minimum size=3pt,inner sep=0pt, label=below:{\small $\ast{\B}_1$}] (b) at (-1.8,-3) {};
 \node [circle,draw=red,fill=white, minimum size=3pt,inner sep=0pt, label=below:{\small $\ast{\B}_2''$}]  (c) at (-0.36,-3){};
  \node [circle,draw=red,fill=white, minimum size=3pt,inner sep=0pt, label=below:{\small $\ast{\B}_1''$}]  (g) at (-1.08,-3){};
 \node [circle,draw=red,fill=white, minimum size=3pt,inner sep=0pt, label=below:{\small $\ast{\B}_{m-1}''$}] (h) at (1.08,-3) {};
 \node [circle,draw=red,fill=red,minimum size=3pt,inner sep=0pt, label=below:{\small $\ast{\B}_2$}] (e) at (1.8,-3) {};
 \foreach \from/\to in {a/b, b/g, g/c, h/e, a/e}
                  \draw[thick] (\from) -- (\to);
\foreach \from/\to in {c/h}
                  \draw[thick] (\from) -- (\to);
\foreach \from/\to in {a/c, a/g, a/h}                   
              \draw[thick] (\from) -- (\to);   
  \node [label=below: $\dotsi$] at (0.36, -3) {};                     
 \node[blue] at (-1.41,-2.8) {\small $s_{i_1^\ast}$};
  \node[blue] at (-0.72,-2.8) {\small ${s}_{i_2^\ast}$}; 
   \node[blue] at (1.44,-2.8) {\small $s_{i_m^\ast}$};          
    \draw[thick,->] (-4,-1.5) -- (-2,-1.5);   
               \node at (-3,-1.3) {${\bf S}$};   
\node [circle,draw=red,fill=red,minimum size=3pt,inner sep=0pt, label=below:{\small ${\B}_3$}] (a) at (-6,-3) {};
\node [circle,draw=red,fill=red,minimum size=3pt,inner sep=0pt, label=above:{\small ${\B}_{1}$}] (b) at (-4.2,0) {};
 \node [circle,draw=red,fill=white, minimum size=3pt,inner sep=0pt, label=above:{\small ${\B}_{2}''$}]  (c) at (-5.64,0){};
  \node [circle,draw=red,fill=white, minimum size=3pt,inner sep=0pt, label=above:{\small ${\B}_{1}''$}]  (g) at (-4.92,0){};
 \node [circle,draw=red,fill=white, minimum size=3pt,inner sep=0pt, label=above:{\small ${\B}_{m-1}''$}] (h) at (-7.08,0) {};
 \node [circle,draw=red,fill=red,minimum size=3pt,inner sep=0pt,  label=above:{\small ${\B}_{2}$}] (e) at (-7.8,0) {};
 \foreach \from/\to in {a/b, b/g, g/c, h/e, a/e}
                  \draw[thick] (\from) -- (\to);
\foreach \from/\to in {c/h}
                  \draw[thick] (\from) -- (\to);
\foreach \from/\to in {a/c, a/g, a/h}                   
              \draw[thick] (\from) -- (\to);               
 \node[blue] at (-4.56,-0.2) {$s_{i_1}$};
  \node[blue] at (-5.28,-0.2) {${s}_{i_{2}}$}; 
   \node[blue] at (-7.41,-0.2) {$s_{i_m}$};     
     \node [label=above: $\dotsi$] at (-6.36, 0) {};           
 \end{tikzpicture}
 }
 \caption{The isomorphism ${\bf S}$.}
 \label{2017.2.26.6.24sh}
 \end{figure}

Recall the notation $n_i({\bf i})$ from Section \ref{10.2.4}. The quiver ${\bf J}_{\rm uf}(\overline{\bf i})$ is isomorphic to ${\bf J}_{\rm uf}({\bf i})$ by the following permutation of vertices
\[
\sigma: J_{\rm uf}(\overline{\bf i}) \stackrel{\sim}{\lra} J_{\rm uf}({\bf i}), \hskip 7mm 
\sigma{i \choose k}={i \choose n_{i}({\bf i})-k}. 
\]
As shown on Figure \ref{2017.2.26.6.24sh}, the following diagram commutes
\bg
 \xymatrix{
      {\rm T}^{\mathscr X}_{\overline{\bf i}} \ar[d]_{\sigma} \ar[r]^{c_{{\mathscr X}, \overline{\bf i}}~~~} &{\rm Conf}_e^{w_0}({\cal B}) \ar[d]^{\bf S}  \\
        {\rm T}^{\mathscr X}_{{\bf i}} \ar[r]^{c_{{\mathscr X},{\bf i}}~~~~} &{\rm Conf}_{w_0}^e({\cal B}). 
 }
\eg
Let $P({\bf i})$ be the cluster Poisson transformation presented by the sequence\footnote{In fact, it is a maximal green sequence in the sense of Keller.} 
\be
\la{cluster.seq.eta}
\sigma \circ L_n({\bf i})\circ \ldots \circ L_2({\bf i})\circ L_1({\bf i}).
\ee
The isomorphism $\iota$ identifies ${\rm Conf}_3^\times({\cal B})$ with ${\rm Conf}_{w_0}^{e}({\cal B})$. By the definition of the twist map,  
\[\eta= \iota^{-1}\circ {\bf S}\circ {\bf R}\circ \iota.
\] 
Therefore the following diagram commutes
\be
\la{twist.x.poisson}
\begin{gathered}
 \xymatrix{
        {\rm T}^{\mathscr X}_{\bf i} \ar[d]_{P({\bf i})}\ar[r]^{c_{{\mathscr X},{\bf i}}~~~} &{\rm Conf}_3^\times({\cal B}) \ar[d]^{\eta} \\
        {\rm T}^{\mathscr X}_{\bf i} \ar[r]_{c_{{\mathscr X}, {\bf i}}~~~} &{\rm Conf}_3^\times({\cal B})
 }
\end{gathered}
\ee
which shows that $\eta$ is a cluster Poisson automorphism.
\end{proof}

The sequence of mutations    in the proof of Theorem \ref{b.cluster.poisson.1541} was inspired by the proof of  \cite[Prop. 3.2]{L4}. 

\begin{remark} 
\la{dt.cluster.eta.double.bruhat}
The map $\eta$ is the cluster Donaldson-Thomas transformation of ${\rm Conf}_3^\times({\cal B})$.\footnote{This was proved for the  case when $\G={\rm PGL}_m$  in \cite{GS16}. For general $\G$, it was proved in \cite{W}.}
In particular, let us tropicalize the diagram \eqref{twist.x.poisson}, getting
\be
\la{tropical.p.x.1132}
\begin{gathered}
 \xymatrix{
       \Z^{ m-r} \ar[d]_{P({\bf i})^t}\ar[r]^{c_{{\mathscr X},{\bf i}}^t~~~~~~~} &{\rm Conf}^{\times}_3({\cal B})(\Z^t)\ar[d]^{\eta^t} \\
        \Z^{ m-r}  \ar[r]_{c_{{\mathscr X}, {\bf i}}^t~~~~~~} &{\rm Conf}_3^{\times}({\cal B})(\Z^t)
 }
\end{gathered}
\ee
The left map $P({\bf i})^t$  is piecewise linear  that sends each positive unit vector $e_i$ to $-e_i$.
\end{remark}

\subsection{Cyclic invariance of the 2-form $\Omega$}
\la{sec.4.5}
Recall the 2-form $\Omega_{u,v}$ in \eqref{canonical.2.form.u.v}.
Let $\Omega:= \iota^*(\Omega_{w_0, e})$ be the induced 2-form of ${\rm Conf}_3^\times({\cal A})$.
\bp 
\la{2017.2.23.13.49h}
The twist  map $\eta$ of ${\rm Conf}_3^\times(\mathcal{A})$ preserves the canonical 2-form:
$
\eta^\ast \Omega = \Omega.
$
\ep
The proof of Proposition \ref{2017.2.23.13.49h} requires a little preparation.

\paragraph{\bf The isomorphism ${\bf R}$.} 
As shown on Figure \ref{2017.2.23.9.46h}, we consider the isomorphism
\[
{\bf R}: {\rm Conf}_{w_0}^{e}({\cal A})\stackrel{\sim}{\lra}{\rm Conf}_e^{w_0} ({\cal A}),\hskip 7mm  (\A_2, \A_3, \A_1, \A_1)\lra (\A_3, \A_3, \A_2\cdot h, \A_1)
\]
where $h\in {\rm H}$ is chosen such that
\[
h(\A_2\cdot h, \A_1)= w_0\big(h(\A_3, \A_2)^{-1}\big)=h(\A_2, \A_3)s_\G.
\]
Therefore
$
h =h(\A_2, \A_3)h(\A_2,\A_1\cdot s_\G)^{-1}.
$
\begin{figure}[ht]
\epsfxsize 200pt
\center{
\begin{tikzpicture}[scale=1.2]
\node [circle,draw=red,fill=red,minimum size=3pt,inner sep=0pt, label=above:{\small ${\A}_1$}] (a) at (0,0) {};
\node [circle,draw=red,fill=red,minimum size=3pt,inner sep=0pt, label=below:{\small ${\A}_2$}] (b) at (-1.8,-3) {};
 \node [circle,draw=red,fill=white, minimum size=3pt,inner sep=0pt, label=below:{\small ${\A}_2'$}]  (c) at (-0.36,-3){};
  \node [circle,draw=red,fill=white, minimum size=3pt,inner sep=0pt, label=below:{\small ${\A}_1'$}]  (g) at (-1.08,-3){};
 \node [circle,draw=red,fill=white, minimum size=3pt,inner sep=0pt, label=below:{\small ${\A}_{m-1}'$}] (h) at (1.08,-3) {};
 \node [circle,draw=red,fill=red,minimum size=3pt,inner sep=0pt, label=below:{\small ${\A}_3$}] (e) at (1.8,-3) {};
 \foreach \from/\to in {a/b, b/g, g/c, h/e, a/e, c/h}
                  \draw[thick] (\from) -- (\to);
\foreach \from/\to in {a/c, a/g, a/h}                   
              \draw[thick] (\from) -- (\to);            
 \node[blue] at (-1.44,-2.8) {\small $s_{i_1^\ast}$};
  \node[blue] at (-0.72,-2.8) {\small ${s}_{i_2^\ast}$}; 
   \node[blue] at (1.44,-2.8) {\small $s_{i_{m}^\ast}$};     
   \node at (.36, -3.4) {$\cdots$};     
    \draw[thick,->] (2,-1.5) -- (4,-1.5);   
               \node at (3,-1.3) {${\bf R}$};   
\node [circle,draw=red,fill=red,minimum size=3pt,inner sep=0pt, label=below:{\small ${\A}_3$}] (a) at (6,-3) {};
\node [circle,draw=red,fill=red,minimum size=3pt,inner sep=0pt, label=above:{\small ${\A}_{1}$}] (b) at (7.8,0) {};
 \node [circle,draw=red,fill=white, minimum size=3pt,inner sep=0pt, label=above:{\small ${\A}_{2}''$}]  (c) at (6.36,0){};
  \node [circle,draw=red,fill=white, minimum size=3pt,inner sep=0pt, label=above:{\small ${\A}_{1}''$}]  (g) at (7.08,0){};
 \node [circle,draw=red,fill=white, minimum size=3pt,inner sep=0pt, label=above:{\small ${\A}_{m-1}''$}] (h) at (4.92,0) {};
 \node [circle,draw=red,fill=red,minimum size=3pt,inner sep=0pt,  label=above:{\small ${\A}_{2}\cdot h~~~$}] (e) at (4.2,0) {};
 \foreach \from/\to in {a/b, b/g, g/c, h/e, a/e, c/h}
                  \draw[thick] (\from) -- (\to);
\foreach \from/\to in {a/c, a/g, a/h}                   
              \draw[thick] (\from) -- (\to);               
 \node[blue] at (7.44,-0.2) {\small $s_{i_1}$};
  \node[blue] at (6.72,-0.2) {\small ${s}_{i_{2}}$}; 
    \node at (5.6, 0.3) {$\cdots$}; 
   \node[blue] at (4.56,-0.2) {\small $s_{i_m}$};        
 \end{tikzpicture}
 }
 \caption{The isomorphism ${\bf R}$.}
 \label{2017.2.23.9.46h}
 \end{figure}
\bl 
\la{2017.2.26.10.09h}
The isomorphism ${\bf R}$ preserves the 2-forms:
$
{\bf R}^*(\Omega_{e}^{w_0}) =\Omega_{w_0}^{e}.
$
\el
\begin{proof} 
Fix a reduced word ${\bf i}=(i_1, \ldots, i_m)$  of $w_0$. By Lemma \ref{decomp.partial.dec.flag.15.02}, the pair $(\A_2, \A_3)$ admits a decomposition
\[
\A_2=\A_0' \xrightarrow{~s_{i_1^\ast}~}\A_1'\xrightarrow{~s_{i_2^\ast}~}\cdots\xrightarrow{~s_{i_m^\ast}~} \A_n'=\A_3\]
We consider the decomposition of $(\A_1, \A_2\cdot h)$
\[
\A_1=\A_0'' \xrightarrow{~s_{i_1}~}\A_1''\xrightarrow{~s_{i_2}~}\cdots\xrightarrow{~s_{i_m}~} \A_m''=\A_2\cdot h
\]
such that 
\be
\la{2017.2.24.21.24h}
h(\A_k'', \A_{k-1}'')=w_0\Big(h(\A_k', \A_{k-1}')^{-1}\Big). 
\ee
Recall $(u_k, v_k)$ in \eqref{2017.2.24.14.43sh}.
The map ${\bf R}$ can be decomposed into a sequence of isomorphisms
\bg
\la{seq.of.refs.a}
{\bf R}:~ {\rm Conf}_{w_0}^{ e}({\cal A}) \stackrel{r_1}{\lra}{\rm Conf}_{u_1}^{ v_1}({\cal A})\stackrel{r_2}{\lra} \cdots\stackrel{r_m}{\lra}{\rm Conf}_{e}^{ w_0}({\cal A}), \\
 (\A_2,\A_3, \A_1, \A_1)\stackrel{r_1}{\lms}(\A_1', \A_3, \A_1'', \A_1)\stackrel{r_2}{\lms}\cdots \stackrel{r_m}{\lms} (\A_3, \A_3, \A_2\cdot h, \A_1).
\eg
where the map $r_k$ is illustrated by Figure \ref{2017.3.23.15.17h}.
\begin{figure}[ht]
\epsfxsize150pt
\center{
 \begin{tikzpicture}[scale=0.8]
 \path [fill=gray!30] (-8,-3) -- (-7,0) -- (-7,-3) -- (-8,-3);
  \path [fill=gray!30] (-1,0) -- (0,0) -- (0,-3) -- (-1,0);
 \node [circle,draw=red,fill=red,minimum size=3pt,inner sep=0pt,label=below:{\small ${\A}_{k-1}'$}] (a) at (-8,-3) {};
\node [circle,draw=red,fill=red,minimum size=3pt,inner sep=0pt,label=above:{\small ${\A}_{k-1}''$}] (e) at (-7,0) {};
\node [circle,draw=red,fill=red,minimum size=3pt,inner sep=0pt,label=above:{\small ${\A}_1$}] (b) at (-5.5,0) {};
 \node [circle,draw=red,fill=white,minimum size=3pt,inner sep=0pt,label=below:{\small ${\A}_k'$}]  (d) at (-7,-3){};
 \node [circle,draw=red,fill=red,minimum size=3pt,inner sep=0pt,label=below:{\small ${\A}_3$}] (c) at (-5.5,-3) {};
 \foreach \from/\to in {a/e, e/b, b/c, d/c, d/a}
                  \draw[thick] (\from) -- (\to);
\foreach \from/\to in {e/d}
                  \draw[dashed, thick] (\from) -- (\to);                     
  \draw[thick,->] (-4,-1.5) -- (-2,-1.5);   
               \node at (-3,-1.3) {$r_k$};         
\node [circle,draw=red,fill=red,minimum size=3pt,inner sep=0pt,label=above:{\small ${\A}_k''$}] (a) at (-1,0) {};
\node [circle,draw=red,fill=white,minimum size=3pt,inner sep=0pt,label=above:{\small ${\A}_{k-1}''$}] (e) at (0,0) {};
\node [circle,draw=red,fill=red,minimum size=3pt,inner sep=0pt,label=above:{\small ${\A}_1$}] (b) at (1.5,0) {};
 \node [circle,draw=red,fill=red,minimum size=3pt,inner sep=0pt,label=below:{\small ${\A}_k'$}]  (d) at (0,-3){};
 \node [circle,draw=red,fill=red,minimum size=3pt,inner sep=0pt,label=below:{\small ${\A}_3$}] (c) at (1.5,-3) {};
 \foreach \from/\to in {a/e, e/b, b/c, d/c, d/a}
                  \draw[thick] (\from) -- (\to);
\foreach \from/\to in {e/d}
                  \draw[dashed, thick] (\from) -- (\to);                  
 \end{tikzpicture}
  }
\caption{The reflection map $r_k$.}
\label{2017.3.23.15.17h}
\end{figure}

By \eqref{canonical.2.form.u.v.deco}, the 2-form $\Omega_{u_k,v_k}$ of ${\rm Conf}_{u_k}^{v_k}({\cal A})$ can be decomposed into 2 parts
\[\Omega_{u_k,v_k} =\Omega_{{\bf i}_k}^{\bf J}+\Omega_{{\bf i}_k}^{\bf H} = \Omega_{{\bf i}_k'}^{\bf J}+\Omega_{{\bf i}_k'}^{\bf H}.\]
Note that $\Omega_{{\bf i}_k}^{\bf H}=\Omega_{{\bf i}_k'}^{\bf H}$ and therefore $\Omega_{{\bf i}_k}^{\bf J}= \Omega_{{\bf i}_k'}^{\bf J}$.  Condition \eqref{2017.2.24.21.24h} implies that
\[
{\bf R}^*(\Omega_{{\bf i}_m}^{\bf H})= \Omega_{{\bf i}_0}^{\bf H}.
\]
Note that the left triangle in the map $r_k$ is the reflection map of elementary configuration spaces and the right rectangle keeps intact. 
By Proposition \ref{2017.2.24.21.20h}, we have
\[
r_k^*(\Omega_{{\bf i}_k}^{\bf J})= \Omega_{{\bf i}_{k-1}'}^{\bf J}= \Omega_{{\bf i}_{k-1}}^{\bf J}.
\]
Therefore ${\bf R}^*(\Omega_{{\bf i}_m}^{\bf J})= \Omega_{{\bf i}_0}^{\bf J}$, which concludes the proof. 
\end{proof}

\subsubsection{The rescaling map $m$.}
Let $i\in {\rm I}$. There is a $\mathbb{G}_m$-action on ${\rm Conf}_3^\times({\cal A})$
\[
m_i:\,{\rm Conf}_3^\times({\cal A}) \times \mathbb{G}_m \lra {\rm Conf}_3^\times({\cal A}),\hskip 1cm (\A_1, \A_2, \A_3)\times t_i \longrightarrow (\A_1, \A_2, \A_3\cdot \alpha_{i}^\vee(t_i)).
\]
Recall the frozen variables of  ${\rm Conf}_3^\times({\cal A})$ associated to the pairs $(\A_1, \A_3)$ and $(\A_2, \A_3)$:
\[
A_{j}^b:= \Delta_{j}(\A_2, \A_3), \hskip 10mm
A_{j}^r:= \Delta_j(\A_3\cdot s_\G, \A_1).
\]
\bl  
\la{2017.3.23.14.44h}
We have: 
\[
m_i^*\Omega=\Omega+\sum_{j\in {\rm I}} \widetilde{C}_{ji}  {\rm d} \log (A_{j^*}^b / A_{j}^r) \wedge  {\rm d} \log (t_i), \hskip 7mm \mbox{where } \widetilde{C}_{ji}=d_j C_{ji}.
\]
\el
\begin{proof}
Let ${\bf i}$ be a reduced word of $w_0$ ending with $i$. 
Under the action $m_i$, we get
\[
A_{i^\ast}^b \lms A_{i^\ast}^b  t_i, \hskip 7mm A_{i}^r\lms A_{i}^r t_i.
\]
The rest coordinates in $ {\bf J}({\bf i})$ remains intact. The Lemma follows by an explicit calculation.
\end{proof} 
\bp 
\la{2017.2.26.20.49}
The 2-form $\Omega$ is invariant under the automorphism
\[
m:~{\rm Conf}_3^{\times}({\cal A}) \stackrel{\sim}{\lra} {\rm Conf}_3^{\times}({\cal A}), \hskip 7mm (\A_1, \A_2, \A_3) \lms (\A_1, \A_2, \A_3\cdot h),
\]
where $h:=  h(\A_3, \A_1)h(\A_3, \A_2)^{-1}$.
\ep
\begin{proof} Set $h=\prod_{i\in I} \alpha_i^\vee(t_i)$, where
$$
 t_i= \Lambda_{i}(h)=\Lambda_i\big( h(\A_3, \A_1)\big)\cdot \Lambda_{i}\big(h(\A_3, \A_2)^{-1}\big)=  A_{i}^r/ A_{i^\ast}^b.
$$
Therefore
$
 {\rm d} \log (t_i)
= {\rm d} \log (A_{i}^r / A_{i^*}^b). 
$ 
By Lemma \ref{2017.3.23.14.44h},  we get
\[
m^*\Omega = \Omega -  \sum_{i, j \in {\rm I}} \widetilde{C}_{ji}\, {\rm d} \log (t_j) \wedge     {\rm d} \log ( t_i)= \Omega.
\]
Note that the matrix $\widetilde{C}_{ji}$ is symmetric. \end{proof}

\paragraph{\bf Proof of Proposition \ref{2017.2.23.13.49h}.}
Let us consider the rotation map
\[
{\bf S}: ~{\rm Conf}_e^{w_0}({\cal A})\lra {\rm Conf}_{w_0}^{e}({\cal A}), \hskip 1cm (\A_3, \A_3, \A_2, \A_1) \lms (\ast\A_1, \ast\A_2, \ast\A_3 \cdot s_{\G}, \ast\A_3 \cdot s_{\G}).
\]
\begin{figure}[ht]
\epsfxsize 200pt
\center{
\begin{tikzpicture}[scale=1.2]
\node [circle,draw=red,fill=red,minimum size=3pt,inner sep=0pt, label=above:{\small $\ast{\A}_3\cdot s_\G$}] (a) at (0,0) {};
\node [circle,draw=red,fill=red,minimum size=3pt,inner sep=0pt, label=below:{\small $\ast{\A}_1$}] (b) at (-1.8,-3) {};
 \node [circle,draw=red,fill=white, minimum size=3pt,inner sep=0pt, label=below:{\small $\ast{\A}_2'$}]  (c) at (-0.36,-3){};
  \node [circle,draw=red,fill=white, minimum size=3pt,inner sep=0pt, label=below:{\small $\ast{\A}_1'$}]  (g) at (-1.08,-3){};
 \node [circle,draw=red,fill=white, minimum size=3pt,inner sep=0pt, label=below:{\small $\ast{\A}_{m-1}'$}] (h) at (1.08,-3) {};
 \node [circle,draw=red,fill=red,minimum size=3pt,inner sep=0pt, label=below:{\small $\ast{\A}_2$}] (e) at (1.8,-3) {};
 \foreach \from/\to in {a/b, b/g, g/c, h/e, a/e, c/h}
                  \draw[thick] (\from) -- (\to);
\foreach \from/\to in {a/c, a/g, a/h}                   
              \draw[thick] (\from) -- (\to);            
 \node[blue] at (-1.44,-2.8) {\small $s_{i_1^\ast}$};
  \node[blue] at (-0.72,-2.8) {\small ${s}_{i_2^\ast}$}; 
   \node[blue] at (1.44,-2.8) {\small $s_{i_m^\ast}$};       
\node at (0.36, -3.3) {$\cdots$};      
    \draw[thick,->] (-4,-1.5) -- (-2,-1.5);   
               \node at (-3,-1.3) {${\bf S}$};   
\node [circle,draw=red,fill=red,minimum size=3pt,inner sep=0pt, label=below:{\small ${\A}_3$}] (a) at (-6,-3) {};
\node [circle,draw=red,fill=red,minimum size=3pt,inner sep=0pt, label=above:{\small ${\A}_{1}$}] (b) at (-4.2,0) {};
 \node [circle,draw=red,fill=white, minimum size=3pt,inner sep=0pt, label=above:{\small ${\A}_{2}'$}]  (c) at (-5.64,0){};
  \node [circle,draw=red,fill=white, minimum size=3pt,inner sep=0pt, label=above:{\small ${\A}_{1}'$}]  (g) at (-4.92,0){};
 \node [circle,draw=red,fill=white, minimum size=3pt,inner sep=0pt, label=above:{\small ${\A}_{m-1}'$}] (h) at (-7.08,0) {};
 \node [circle,draw=red,fill=red,minimum size=3pt,inner sep=0pt,  label=above:{\small ${\A}_{2}$}] (e) at (-7.8,0) {};
 \foreach \from/\to in {a/b, b/g, g/c, h/e, a/e, c/h}
                  \draw[thick] (\from) -- (\to);
\foreach \from/\to in {a/c, a/g, a/h}                   
              \draw[thick] (\from) -- (\to);               
 \node[blue] at (-4.56,-0.2) {\small $s_{i_1}$};
  \node[blue] at (-5.28,-0.2) {\small ${s}_{i_{2}}$}; 
   \node[blue] at (-7.44,-0.2) {\small $s_{i_m}$};      
\node at (-6.2, .3) {$\cdots$};       
 \end{tikzpicture}
 }
 \caption{The isomorphism ${\bf S}$.}
 \label{2017.2.26.6.24h}
 \end{figure}
 \bl 
 \label{2017.2.26.20.49h}
 The map ${\bf S}$ preserves the 2-form:
 \[
{\bf S}^*(\Omega_{w_0}^{ e}) =\Omega_{e}^{ w_0}.
 \]
 \el
 
 \begin{proof}
 Let $\{\A_1=\A_0', \A_1',\ldots, \A_n'=\A_2\}$ be a decomposition of the pair $(\A_1, \A_2)$, which gives rise to a quiver of  ${\rm Conf}_e^{w_0}({\cal A})$. 
 We rotate the triangle by $180^\circ$ and applying $\ast$ action to each decorated flag, obtaining a quiver of ${\rm Conf}_{w_0}^{e}({\cal A})$. 
 By Proposition \ref{prop.11.12.13}, we have 
 \[
 \Delta_{i}(\ast \A_3\cdot s_\G, \ast \A_k') = \Delta_i(\A_k',\A_3).
 \]
Hence the cluster variables are permuted correspondingly. Thus the $2$-forms are ${\bf S}$-invariant.
\end{proof}

Note that the twisted map  ${\eta}= m\circ \iota^{-1}\circ {\bf S}\circ {\bf R}\circ \iota$. 
Combining Lemma \ref{2017.2.26.10.09h}, Proposition \ref{2017.2.26.20.49}, Lemma \ref{2017.2.26.20.49h}, and \eqref{2017.2.26.10.09hh}, 
we finish the proof of Proposition \ref{2017.2.23.13.49h}.

\subsection{Frozen vertices of the quiver $ {\bf J}({\bf i})$.}
Let ${\bf q}=(\varepsilon_{ij})$ be a weighted quiver with neither loops nor 2 cycles. Deleting all of its frozen vertices, we obtain a full subquiver ${\bf q}_0$. 
\bl[{\cite[Proposition 2.10]{GS16}}]
\la{frozen.vs.xtrop}
Every frozen vertex of  ${\bf q}$ canonically gives rise to an integral tropical point of the underlying cluster Poisson variety ${\mathscr X}_{|{\bf q}_0|}$.
\el 
\begin{proof}
Let us parametrize the non-frozen vertices in ${\bf q}$    by $\{1,\ldots, m\}$. A frozen vertex $j$ gives rise to a vector
$
f_j:= (\varepsilon_{1j}, \ldots, \varepsilon_{mj}) \in \Z^m,
$  where $\varepsilon_{ij}$'s are entries of the exchange matrix associated with ${\bf q}$.
For the mutated quiver $\mu_k({\bf q})$, it becomes
\[
f_j':=  (\varepsilon_{1j}', \ldots, \varepsilon_{mj}') \in \Z^m, \hskip 9mm \mbox{where } \varepsilon_{ij}':=\left\{   \begin{array}{ll} 
-\varepsilon_{ij} & \mbox{if } i=k, \\
\varepsilon_{ij} -\varepsilon_{ik}\min\{0, -{\rm sgn}(\varepsilon_{ik})\varepsilon_{kj}\} & \mbox{if } i\neq k.
   \end{array}\right.
\]
Note that the above formula coincides with the tropicalization of the cluster Poisson mutation $\mu_k$. Hence we get an integral tropical point
$$
{\bf f}_j:=c_{{\mathscr X}, {\bf q}_0}^t(f_j) =  c_{{\mathscr X}, \mu_k({\bf q}_0)}^t(f_j')  \in \mathscr{X}_{|{\bf q}_0|}(\Z^t). 
$$
\end{proof}
Let $(\A_1, \A_2, \A_3)\in {\rm Conf}_3^{\times}({\cal A})$. The $3r$ frozen variables of ${\rm Conf}_3^{\times}({\cal A})$ are
\be
\la{frozen.var.ed}
A_i^l:=\Delta_i(\A_1, \A_2), \hskip 7mm A_i^b:=\Delta_i(\A_2, \A_3), \hskip 7mm A_i^r:= \Delta_i(\A_3\cdot s_\G, \A_1), \hskip 6mm i\in {\rm I}.
\ee
which corresponds to the $r$ frozen vertices placed on each side of the triangle respectively. By Lemma \ref{frozen.vs.xtrop}, we obtain $3r$ tropical points denoted by
\[
{\bf l}_i, ~ {\bf b}_i, ~{\bf r}_i \in {\rm Conf}_3^\times({\cal B})(\Z^t), \hskip 10mm i\in {\rm I}.\]

The tropicalization of the twisted map $\eta$ is an bijection
\[
\eta^t: {\rm Conf}_3^{\times}({\cal B})(\Z^t)\xrightarrow{~~\sim~~}{\rm Conf}_3^{\times}({\cal B})(\Z^t).
\]
\bl 
\la{2017.4.2.14.14hh}
The $\eta^t$ permutes the frozen vertices in the following way
\be
\la{2017.4.2.14.13hh}
\eta^t({\bf l}_i) = {\bf b}_{i^\ast}, \hskip 7mm \eta^t({\bf b}_i) = {\bf r}_{i^*},  \hskip 7mm \eta^t({\bf r}_i)={\bf l}_{{i}^*}. 
\ee
\el
\begin{proof}
Let ${\bf i}=(i_1,\ldots, i_m)$ be a reduced word of $w_0$. It gives rise to a bijection
\[
c_{{\mathscr X},{\bf i}}^t:  ~~\Z^{m-r}\xrightarrow{~\sim~} {\rm Conf}_3^{\times}({\cal B})(\Z^t).
\]

Suppose that $i_1=i$. Recall the quiver $ {\bf J}({\bf i})$ in Definition \ref{quiver.amla.fg4.sec22}. The frozen vertex corresponding to ${\bf l}_i$ is placed at the most left on level $i$. So   the preimage of ${\bf l}_i$ under $c_{{\mathscr X},{\bf i}}^t$  is a positive unit vector
\[
(c_{{\mathscr X},{\bf i}}^t)^{-1} ({\bf l}_i) = e_{i \choose 1}.
\]
Similarly, we have 
$
(c_{{\mathscr X},{\bf i}}^t)^{-1} ({\bf b}_{i^*}) =- e_{i \choose 1}.
$
By Remark \ref{dt.cluster.eta.double.bruhat},  we get $\eta^t({\bf l}_i) = {\bf b}_{i^\ast}$.

Suppose that $i_m = i$. Similarly, we get
\[
(c_{{\mathscr X},{\bf i}}^t)^{-1} ({\bf b}_i) = e_{i \choose n_i -1}, \hskip 7mm  (c_{{\mathscr X},{\bf i}}^t)^{-1} ({\bf r}_{i^*}) = -e_{i \choose n_i -1}.
\]
Therefore $\eta^t({\bf b}_i) = {\bf r}_{i^*}$.

The map $\eta^3$ is the involution $(\B_1, \B_2, \B_3)\lms (\ast\B_1, \ast \B_2, \ast \B_3)$.
Let ${\bf i}^\ast = (i_1^\ast, \ldots, i_m^\ast)$. Let $I$ be a bijection from ${\bf J}_\circ({\bf i})$ to ${\bf J}_\circ({{\bf i}^\ast})$ by relabeling level $i$ by $i^*$. The following diagram commutes
\be
\begin{gathered}
 \xymatrix{
        {\rm T}^{\mathscr X}_{\bf i} \ar[d]_{I}\ar[r]^{c_{{\mathscr X},{\bf i}}~~~} &{\rm Conf}_{3}^\times({\cal B})\ar[d]^{\eta^3} \\
        {\rm T}^{\mathscr X}_{{\bf i}^\ast} \ar[r]_{c_{{\mathscr X}, {\bf i}^*}~~~} &{\rm Conf}_3^{\times}({\cal B})
 }
\end{gathered}
\ee
Therefore $(\eta^t)^3 ({\bf l}_i) = {\bf l}_{i^*}$. 
Note that ${\bf r}_i =\eta^t({\bf b}_{i^*})= (\eta^t)^2({\bf l}_i)$.
Therefore
$
\eta^t({\bf r}_i) = (\eta^t)^3 ({\bf l}_i) = {\bf l}_{{i}^*}.
$
\end{proof}
Let $\tau$ be the permutation of frozen vertices of $ {\bf J}({\bf i})$ following \eqref{2017.4.2.14.13hh}. The following result is a direct consequence of Theorem \ref{main.result2.cluster.cyclic} and  Lemma \ref{2017.4.2.14.14hh}.
\bp The following sequence preserves the quiver $ {\varepsilon}({\bf i})$ of $ {\bf J}({\bf i})$ up to arrows among frozen vertices:\be
\la{2017.4.2.14.36hh}
\tau\circ \sigma \circ L_{m}({\bf i})\circ \ldots \circ L_2({\bf i})\circ L_1({\bf i}).
\ee
\ep
\begin{remark} In the framework of cluster algebras \cite{FZI}, the arrows between the frozen vertices were not considered. Note that the quiver is encoded by the canonical 2-form $\Omega$ of ${\rm Conf}_3^\times({\cal A})$.  Therefore Proposition \ref{2017.2.23.13.49h}   implies that \eqref{2017.4.2.14.36hh} preserves the arrows among frozen vertices as well. 
\end{remark}

\subsection{Proof of Theorem \ref{main.result2.cluster.cyclic}}
 Let $K({\bf i}):  {\rm T}^{\mathscr A}_{\bf i}\lra  {\rm T}^{\mathscr A}_{\bf i} $ be the cluster $K_2$-transformation presented by the sequence \eqref{2017.4.2.14.36hh}. Conjugating with \eqref{cai}, we obtain an automorphism of ${\rm Conf}_3^{\times}({\cal A})$ denoted by
\[
D:= c_{{\mathscr A}, {\bf i}} \circ K({\bf i})\circ c_{{\mathscr A}, {\bf i}}^{-1}.
\] 
\bl 
\label{well-define.Dmap}
The automorphism $D$ is independent of ${\bf i}$ chosen.
\el
\begin{proof}
Let ${\bf i}, {\bf i}'$ be reduced words of $w_0$. For the cluster Poisson variety ${\rm Conf}_3^{\times}({\cal B})$, by \eqref{twist.x.poisson}, the following diagram commutes
\be
\begin{gathered}
 \xymatrix{
        {\rm T}^{\mathscr X}_{\bf i} \ar[d]_{P({\bf i})}\ar[r]^{\chi_{{\bf i}',{\bf i}}} &{\rm T}^{\mathscr X}_{{\bf i}'} \ar[d]^{P({\bf i}')} \\
        {\rm T}^{\mathscr X}_{\bf i} \ar[r]_{\chi_{{\bf i}', {\bf i}}} &{\rm T}^{\mathscr X}_{{\bf i}'}
 }
\end{gathered}
\ee
where $\chi_{{\bf i}',{\bf i}}$ is the transition map  $c_{{\mathscr X}, {\bf i}'}^{-1}\circ c_{{\mathscr X},{\bf i}}$.

For ${\rm Conf}_3^{\times}({\cal A})$,  the cluster $K_2$-transformations $\alpha_{{\bf i}, {\bf i}'}:= c_{{\mathscr A}, {\bf i}'}^{-1}\circ c_{{\mathscr A}, {\bf i}}$, $K({\bf i})$, $K({\bf i}')$ are presented by the same sequences of cluster mutations as $\chi_{{\bf i}, {\bf i}'}$, $P({\bf i})$, $P({\bf i}')$ respectively. As a general result of cluster algebras,  if two sequences of quiver mutations give rise to the same cluster Poisson transformation, then they also give rise to the same cluster $K_2$-transformation. Meanwhile, the exchange graph of a cluster algebra only depends on its unfrozen part. 
Hence the follow diagram commutes
\be
\begin{gathered}
 \xymatrix{
        {\rm T}^{\mathscr A}_{\bf i} \ar[d]_{K({\bf i})}\ar[r]^{\alpha_{{\bf i}',{\bf i}}} &{\rm T}^{\mathscr A}_{{\bf i}'} \ar[d]^{K({\bf i}')} \\
        {\rm T}^{\mathscr A}_{\bf i} \ar[r]_{\alpha_{{\bf i}', {\bf i}}} &{\rm T}^{\mathscr A}_{{\bf i}'}
 }
\end{gathered}
\ee
The Lemma is proved.
\end{proof}

Recall the potential functions assigned to the left bottom angle
\[
W^L_i: {\rm Conf}_3^{\times}({\cal A})\lra{\Bbb A},
\hskip 7mm
(\A_1, \A_2, \A_3)\lms {\cal W}_i(\A_2, \B_3, \B_1).
\]
\bl
\la{birational.sio.hh.4.2}
The following map is a birational isomorphism:
\[
(p, e, W^L): ~{\rm Conf}_3^{\times}({\cal A})\lra {\rm Conf}_3^{\times}({\cal B}) \times {\rm H}^2 \times {\Bbb A}^r
\]
\[
(\A_1, \A_2, \A_3) \lms (\B_1, \B_2, \B_3) \times \big(h(\A_2, \A_1), h(\A_2, \A_3)\big) \times ({W}_1^L, \ldots, {W}_r^L).
\]
\el
\begin{proof}
A configuration $(\B_1, \B_2, \B_3)\in {\rm Conf}_3^{\times}({\cal B})$ with partial potentials $(W_1^L, \ldots, W_r^{L})$ uniquely determines a configuration $(\B_1, \A_2, \B_3)\in \G \backslash ({\cal B}\times {\cal A}\times {\cal B})$. The Cartan elements $h_1, h_3$ determine the decorations of $\B_1$ and $\B_3$. 
\end{proof}

\paragraph{\bf Proof of Theorem \ref{main.result2.cluster.cyclic}.} To distinguish the twisted maps \eqref{twisted.map.a} and \eqref{eta.dt.conf3b}, we denote the former by $\eta_A$, and the latter by $\eta_X$.
We prove Theorem \ref{main.result2.cluster.cyclic} by showing that  $D=\eta_A.$ 
By Lemma \ref{birational.sio.hh.4.2}, it is equivalent to prove that
\be \la{TROIZA}
(p, e, W^L)\circ \eta_A = (p, e, W^L) \circ D.
\ee

i) By \eqref{projection.2017.3.27.12.00}, we have
$
 c_{{\mathscr X},{\bf i}}\circ p_{\bf i} = p\circ c_{{\mathscr A},{\bf i}}. 
$ 
Since $K({\bf i})$ and $P({\bf i})$ are presented by same sequence of mutations, we get
$ p_{\bf i}\circ K({\bf i})= P({\bf i})\circ p_{\bf i}. 
$ 
Therefore
$$
p\circ D = p \circ c_{{\mathscr A},{\bf i}}\circ K({\bf i})\circ c_{{\mathscr A},{\bf i}}^{-1}=  c_{{\mathscr X},{\bf i}}\circ P({\bf i})\circ c_{{\mathscr X},{\bf i}}^{-1}\circ p =  \eta_X \circ p=p \circ \eta_A .
$$

ii) By the construction of $K({\bf i})$, it follows easily that
$
e\circ \eta_A= e\circ D.
$

iii) To finish the proof of (\ref{TROIZA}), it remains to show that 
\be
 \la{potential.lusztig.alg.4.2}
W^L_i\circ \eta_A =W^L_i \circ D.
\ee
By Lemma \ref{well-define.Dmap}, we may fix a reduced word ${\bf i}$ of $w_0$ starting with $i$. 
Consider the following associated quivers with cluster $K_2$ and Poisson variables
\[
 {\bf J}({\bf i})=\big( {J}({\bf i}),  {\varepsilon}({\bf i}), \{A_{{\bf i},j}\}\big), \hskip 10mm {\bf J}_{\rm uf}({\bf i})=\big({J}_{\rm uf}({\bf i}), {\varepsilon}_{\rm uf}({\bf i}), \{X_{{\bf i},j}\}\big).
\]
The set 
$V:=\left\{{i \choose 1}, \ldots , {i \choose n_i-1}\right\}
$ of non-frozen vertices on level $i$  
provides a type $A$ quiver
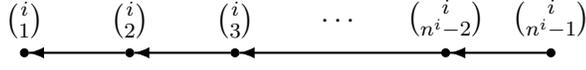
\begin{center}
\begin{tikzpicture}[scale=1.4]
\node [circle,draw,fill,minimum size=3pt,inner sep=0pt, label=above:$i \choose 1$] (a) at (0,0) {};
\node [circle,draw,fill,minimum size=3pt,inner sep=0pt, label=above:$i \choose 2$] (b) at (1,0) {};
 \node [circle,draw,fill, minimum size=3pt,inner sep=0pt, label=above:$i \choose 3$]  (c) at (2,0){};
  \node [circle,draw,fill, minimum size=3pt,inner sep=0pt, label=above:$i \choose n_i-2$]  (d) at (4,0){};
 \node [circle,draw,fill,minimum size=3pt,inner sep=0pt, label=above:$i \choose n_i-1$] (e) at (5,0) {};
 \node at (3, 0.3) {$\cdots$};
 \foreach \from/\to in {b/a, c/b, e/d}
                  \draw[directed, thick] (\from) -- (\to);
\foreach \from/\to in {d/c}                   
              \draw[directed, thick] (\from) -- (\to);            
     \end{tikzpicture}
    \end{center}
To simplify the notation, we write
\be
A_k := A_{{\bf i}, {i \choose k}}, \hskip 10mm X_k := X_{{\bf i}, {i \choose k}}=\prod_{j\in  {J}({\bf i})} A_{{\bf i}, j}^{\varepsilon_{j, {i \choose k}}}, \hskip 7mm X_k' :=\prod_{j\in  {J}({\bf i})-V} A_{{\bf i}, j}^{\varepsilon_{j, {i \choose k}}}
\ee
Recall the frozen variables $A_i^l$ and $A^b_{i^\ast}$ in \eqref{frozen.var.ed}. 
Take $\B', \B''$ as shown on Figure \ref{2017.4.2.9.46h}.  By \eqref{identity.reflect.f.3},  
\be
\la{leftpotential.1}
W_{i^*}^L= {\cal W}_{i^*}(\A_2, \B_3, \B_1) = {\cal W}(\A_2, \B', \B'')= \frac{A_{1}}{A_{i}^l A^b_{i^*}},
\ee
 By the definition of $D$, we have
\be \la{brutal.calculation.f.g.c}
\begin{split}
&D^* A_i^l = A_{i^*}^r= A_{{\bf i},{i \choose n_i}} := A_{n_i},  ~~~~~D^* A_{i^*}^b = A_i^l = A_{{\bf i}, {i \choose 0}}:=A_0,\\
&D^* A_{1}  = L_1({\bf i})^* A_{n_i-1}:= \big(\mu_{i \choose n_i-1}\circ \cdots \circ \mu_{i \choose 1}\big)^* A_{n_i -1}.\\
\end{split}
\ee
\begin{figure}[ht]
\epsfxsize 200pt
\center{
\begin{tikzpicture}[scale=0.9]
\node [circle,draw=red,fill=red,minimum size=3pt,inner sep=0pt, label=above:{\small ${\A}_1$}] (a) at (0,0) {};
\node [circle,draw=red,fill=red,minimum size=3pt,inner sep=0pt, label=below:{\small ${\A}_2$}] (b) at (-1.8,-3) {};
 \node [circle,draw=red,fill=white, minimum size=3pt,inner sep=0pt, label=left:{\small ${\B}''$}]  (c) at (-1.2,-2){};
  \node [circle,draw=red,fill=white, minimum size=3pt,inner sep=0pt, label=below:{\small ${\B}'$}]  (g) at (-0.88,-3){};
 \node [circle,draw=red,fill=red,minimum size=3pt,inner sep=0pt, label=below:{\small ${\A}_3$}] (e) at (1.8,-3) {};
 \foreach \from/\to in {a/c, b/c, g/b, a/e, g/e}
                  \draw[thick] (\from) -- (\to);
\foreach \from/\to in {a/g}                   
              \draw[dashed] (\from) -- (\to);            
 \node[blue] at (-1.34,-3.2) {\small $s_{i^\ast}$};
  \node[blue] at (-1.7,-2.5) {\small ${s}_{i^\ast}$};    
     \end{tikzpicture}
    }
 \caption{The partial potential $W_{i^*}^L$.}
 \label{2017.4.2.9.46h}
 \end{figure}
We compute the second line in \eqref{brutal.calculation.f.g.c} by using $F$-polynomials and $g$-vectors of   \cite{FZIV}.   
Set 
 \[
\prod_{j} A_j^{b_j} \oplus \prod_{j} A_j^{c_j} = \prod_{j} A_j^{\min\{b_j, c_j\}}.
 \]
For a Laurent polynomial $F$ in variables $A_j$, we define $F^{\rm Trop}$ by replacing each $+$ by $\oplus$. 
By \cite[Corollary 6.3]{FZIV},   the last formula in \eqref{brutal.calculation.f.g.c} can be presented by 
\be
\la{L.ex}
L_1({\bf i})^* A_{ n_i -1} =\frac{F(X_{1}, \ldots X_{ n_i-1})}{F(X_{1}', \ldots X_{n_i-1}')^{\rm Trop}} A_{1}^{g_1} \ldots A_{n_i-1}^{g_{n_i-1}}.
\ee
Following the algorithm of \cite[Proposition 5.1]{FZIV},  the above $F$-polynomial is 
$$
F(X_1, \ldots, X_{n_i-1}) = 1+ X_{1} + X_{1} X_{ 2}+ \ldots+ X_{1}\cdots X_{n_i-1}.
$$
Following \cite[Proposition 6.6]{FZIV}, the above $g$-vector is
\be
\la{g.vector.expression}
(g_1,\ldots, g_{n_i-1})  = (-1, 0, \ldots, 0) \in \Z^{n_i-1}.
\ee
Recall the potential
$
a_k := {\cal W}_i(\A_1, \B_{{i \choose k}-1}', \B_{i \choose k}'). 
$ 
By \eqref{Xv}, we have
$
X_{k} = a_{k+1} / a_k.
$ 
Therefore
\be
\la{F.poly.ex}
{F(X_{1}, \ldots X_{n_i-1})} = \frac{a_1+a_2+\ldots +a_{n_i}}{a_1}=\frac{{\cal W}_{i}(\A_1, \B_2, \B_3)}{a_1}.
\ee
Meanwhile
\be
X_{k}' =  \left\{ \begin{array}{ll} {\displaystyle a_2 A_2 / a_1},~~~~ &\mbox{if } k=1,\\
{\displaystyle a_{k+1}A_{k+1}  / (a_kA_{k-1})},~~~~ &\mbox{if } 1<k<n_i-1,\\
{\displaystyle a_{n_i} / (a_{n_{i}-1}A_{n_i-2})}  &\mbox{if }k=n_i-1.  \end{array} \right.
\ee
According to the explicit expression \eqref{bfz.lusztig} of $a_k$, we get
 \bg
 \la{F.trop.ex}
 F(X_{1}', \ldots X_{n_i-1}')^{\rm Trop} =\frac{a_1A_1\oplus a_2A_1A_2\oplus \ldots \oplus a_{n_i-1}A_{n_i-1}A_{n_i}\oplus a_{n_i}A_{n_i-1}}{a_1A_1}\\
 = (A_0 A_{n_i} a_1A_1)^{-1}.
 \eg
Plugging \eqref{g.vector.expression}, \eqref{F.poly.ex} and \eqref{F.trop.ex} into \eqref{L.ex}, we get
$
L_1({\bf i})^* A_{ n_i-1} = {\cal W}_{i}(\A_1, \B_2, \B_3) A_{0} A_{n_i}.
$
Therefore
$
D^* W_{i^*}^L = {\cal W}_{i}(\A_1, \B_2, \B_3) = \eta^*_A  W_{i^*}^L .
$

\subsection{Cluster $K_2$-structure of the moduli space ${\mathscr A}_{\G, \bS}$} \la{Sect8.5}
We have defined a cluster $K_2$-structure on the space 
${\cal A}_{\G, t}{=}{\rm Conf}_3^\times(\mathcal{A})$, and proved that it does not depend on the reduced decomposition ${\bf i}$ of $w_0$ and a choice of a vertex of the triangle $t$. It remains to prove that it 
does not change under a flip of triangulation. Let ${\bf q}$ be a quadrilateral where we perform a flip of a diagonal $E$. There are two cluster coordinate systems on ${\cal A}_{\G, \bf q}$ using the two amalgamation patterns corresponding to the two diagonals of ${\bf q}$. 
For example, see    the left and  the right pictures on Figure \ref{pin18} for $\G={\rm SL}_3$. 
Each of them corresponds to a particular reduced decomposition of $(w_0, w_0)$. By Theorem  \ref{Th2.3}, there is a cluster $K_2$-transformation relating them. 

 This concludes the proof of the existence of a cluster $K_2$-structure on the moduli space ${\mathscr A}_{\G, \bS}$, equivariant under the action of the group $\Gamma_\bS$, claimed in 
 Theorem \ref{MTHa}.
 
 \begin{figure}[ht]
 \begin{center}
 \begin{tikzpicture}[scale=0.5]
\draw (0,0)--(3,0)--(3,2)--(0,2)--(0,0);
\coordinate (a) at (0,0);
\coordinate (b) at (1,0);
\coordinate (c) at (2,0);
\coordinate (d) at (3,0); 
\coordinate (A) at (0,2);
\coordinate (B) at (1,2);
\coordinate (C) at (2,2);
\coordinate (D) at (3,2);
\draw[-latex] (3.5,1)--(4.5,1);
\foreach \from/\to in {B/d, C/d, A/b, A/c}
\draw (\from)--(\to);
\draw[ultra thick, red] (A)--(d);
\begin{scope}[shift={(5,0)}]
\draw[-latex] (3.5,1)--(4.5,1);
\draw (0,0)--(3,0)--(3,2)--(0,2)--(0,0);
\coordinate (a) at (0,0);
\coordinate (b) at (1,0);
\coordinate (c) at (2,0);
\coordinate (d) at (3,0); 
\coordinate (A) at (0,2);
\coordinate (B) at (1,2);
\coordinate (C) at (2,2);
\coordinate (D) at (3,2);
\foreach \from/\to in {A/b, B/c, C/d}
\draw (\from)--(\to);
\draw[ultra thick] (A)--(c);
\draw[ultra thick] (B)--(d);
\end{scope}
\begin{scope}[shift={(10,0)}]
\draw[-latex] (3.5,1)--(4.5,1);
\draw (0,0)--(3,0)--(3,2)--(0,2)--(0,0);
\coordinate (a) at (0,0);
\coordinate (b) at (1,0);
\coordinate (c) at (2,0);
\coordinate (d) at (3,0); 
\coordinate (A) at (0,2);
\coordinate (B) at (1,2);
\coordinate (C) at (2,2);
\coordinate (D) at (3,2);
\foreach \from/\to in {B/b, B/c, C/c}
\draw (\from)--(\to);
\draw[ultra thick] (A)--(b);
\draw[ultra thick] (C)--(d);
\end{scope}
\begin{scope}[shift={(15,0)}]
\draw[-latex] (3.5,1)--(4.5,1);
\draw (0,0)--(3,0)--(3,2)--(0,2)--(0,0);
\coordinate (a) at (0,0);
\coordinate (b) at (1,0);
\coordinate (c) at (2,0);
\coordinate (d) at (3,0); 
\coordinate (A) at (0,2);
\coordinate (B) at (1,2);
\coordinate (C) at (2,2);
\coordinate (D) at (3,2);
\foreach \from/\to in {B/a, C/b, D/c}
\draw (\from)--(\to);
\draw[ultra thick] (B)--(b);
\draw[ultra thick] (C)--(c);
\end{scope}
\begin{scope}[shift={(20,0)}]
\draw[-latex] (3.5,1)--(4.5,1);
\draw (0,0)--(3,0)--(3,2)--(0,2)--(0,0);
\coordinate (a) at (0,0);
\coordinate (b) at (1,0);
\coordinate (c) at (2,0);
\coordinate (d) at (3,0); 
\coordinate (A) at (0,2);
\coordinate (B) at (1,2);
\coordinate (C) at (2,2);
\coordinate (D) at (3,2);
\foreach \from/\to in {B/a, C/a, D/b, D/c}
\draw (\from)--(\to);
\draw[ultra thick] (C)--(b);
\end{scope}
\begin{scope}[shift={(25,0)}]
\draw (0,0)--(3,0)--(3,2)--(0,2)--(0,0);
\coordinate (a) at (0,0);
\coordinate (b) at (1,0);
\coordinate (c) at (2,0);
\coordinate (d) at (3,0); 
\coordinate (A) at (0,2);
\coordinate (B) at (1,2);
\coordinate (C) at (2,2);
\coordinate (D) at (3,2);
\foreach \from/\to in {B/a, C/a, D/b, D/c}
\draw (\from)--(\to);
\draw[ultra thick, red] (D)--(a);
\end{scope}
\end{tikzpicture}
\end{center}
\caption{Sequence of moves proving independence of the cutting diagonal.}
\label{pin18}
\end{figure}

\section{Cluster Poisson structure of  moduli spaces ${\mathscr P}_{\G, t}$ and ${\mathscr P}_{\G, \bS}$} \la{CSEC7}

\medskip

In Section \ref{CSEC7} we assume first that all boundary intervals of $\bS$ are colored. 
Then the  passage to the case where only some of the boundary intervals are colored is trivial: one removes the frozen vertices and the frozen coordinates corresponding to the non-colored sides.  

In Section \ref{CSEC7}  we continue to use the notation $\G$ for the adjoint group, and $\G'$ for the simply-connected group. 

 \medskip
 
 \subsection{Cluster Poisson coordinates on ${\mathscr P}_{\G, \bS}$} \la{sec5..1}

\medskip

Let $\G$ be an adjoint group.  Recall the moduli space ${\mathscr P}_{\G, t}$,  depicted   on Figure \ref{2018.9.23.15.10ss}. 
 Let ${\bf i}$ be a reduced word of $w_0$. Pick a vertex $\B_1$ of the triangle $t$.  Recall the quiver $ {\bf J}({\bf i})$ in Definition \ref{quiver.amla.fg4.sec22} and  \eqref{alternative.notation.J0}. In Section \ref{CSEC4.3}  we   introduced cluster Poisson coordinates on  ${\mathscr P}_{\G, t}$, see Definition \ref{DEFXC}. In Section \ref{sec5..1} we observe that they are 
  indexed by the vertices of ${\bf J}({\bf i})$, and establish their  properties.

\subsubsection{Subquivers labelled by $i\in {\rm I}$.}
Let $i\in {\rm I}$ be a vertex of the Dynkin diagram.  The  sub-quiver of ${\bf J}({\bf i})$ formed by the vertices of level $i$ is a type $A$ quiver
 \begin{center} 
\begin{tikzpicture}[scale=1.4]
\node [circle,draw,fill=white,minimum size=5pt,inner sep=0pt, label=above:$i \choose 0$] (a) at (0,0) {};
\node [circle,draw,fill,minimum size=5pt,inner sep=0pt, label=above:$i \choose 1$] (b) at (1,0) {};
 \node [circle,draw,fill, minimum size=5pt,inner sep=0pt, label=above:$i \choose 2$]  (c) at (2,0){};
  \node [circle,draw,fill, minimum size=5pt,inner sep=0pt, label=above:$i \choose n^i-1$]  (d) at (4,0){};
 \node [circle,draw,fill=white,minimum size=5pt,inner sep=0pt, label=above:$i \choose n^i$] (e) at (5,0) {};
 \foreach \from/\to in {b/a, c/b, e/d, d/c}
                  \draw[directed, thick] (\from) -- (\to);                 
              \draw[dotted, thick] (3.2,0.3) -- (2.8,0.3);            
     \end{tikzpicture}
    \end{center} 
 The rest of the vertices are parametrized by ${i \choose -\infty}$, placed on the side $\{2,3\}$.

Let us  summarize our definition of the cluster Poisson coordinates on the space  ${\mathscr P}_{\G', t}$:

\bd The cluster Poisson coordinates on the space  ${\mathscr P}_{\G, t}$ assigned to a reduced word ${\bf i}$ for $w_0$ and the vertex $\B_1$ of the triangle $t$ are given by

\begin{itemize}

\item   For each $ i \in {\rm I}$,  functions $X_{i \choose k}$ in (\ref{defxikX}), assigned to   the type $A$ subquiver of  ${\bf J}({\bf i})$ above. 

\item Functions $X_{i \choose -\infty}$ in (\ref{bottomXY}), assigned to the subquiver ${\bf H}({\bf i})$ of ${\bf J}({\bf i})$.
\end{itemize} 

\ed

Recall the space ${\mathscr P}_{u}^v$ in Section \ref{partial.conf.part3}. There is a natural map from ${\mathscr P}_{\G, t}$ onto ${\mathscr P}_{w_0}^e$ by forgetting the pinning $p_{23}$. The coordinates $X_{i \choose k}$ in (\ref{defxikX}) coincide with the cluster Poisson coordinates of ${\mathscr P}_{w_0}^e$. In particular, the cluster Poisson coordinates for unfrozen vertices are independent of the pinnings. The frozen vertices are placed on the sides of $t$. Their cluster Poisson coordinates only depend on the pinning at the related side.  

\subsubsection{Action of the Cartan group.}  
For $\{a, b\}= \{1, 2\}, \{2, 3\}$ or $\{3, 1\}$, let $i_{ab}$ be the frozen vertex on the side $\{a, b\}$ assigned to $i\in {\rm I}$:
\[
i_{12}:={i \choose 0}; \hskip 5mm i_{23}:={i \choose - \infty};  \hskip 5mm  i_{31}={i^\ast \choose n_{i^\ast}}.
\]
The Cartan group ${\rm H}$ acts on the space of pinnings by 
\be
\la{cartan.pinning.action.hhh}
 (\A, \A') \longmapsto \left(\A\cdot h, \A' \cdot w_0(h)\right), ~~~~ \forall h \in {\rm H}.
\ee
This action gives rise to an action $\tau$ of ${\rm H}^3$ on ${\mathscr P}_{\G, t}$ by altering pinnings $p_{ab}$ by $h_{ab}\in {\rm H}$. 
\begin{lemma} 
\label{rescale.hh}
Under the action of $\tau$, we have
\[
\tau^* X_{i_{ab}} = \alpha_i(h_{ab}^{-1}) \cdot X_{i_{ab}}.
\]
\end{lemma}
\begin{proof} Let us prove the case when $\{a,b\}=\{2,3\}$. The other two cases follow easily from their definitions. 
 For ${\bf i}=(i_1,\ldots,i_m)$, let us define
\[
\gamma_1 = s_{i_m}\ldots s_{i_{k}}\cdot \Lambda_{i_k}, \hskip 14mm 
\gamma_2 = s_{i_m}\ldots s_{i_{k+1}}\cdot \Lambda_{i_k}. 
\]
Then $\alpha_i =\Lambda_i-s_i \cdot \Lambda_i$ is a simple positive root  for all $i\in {\rm I}$. Therefore 
$
\gamma_2-\gamma_1=s_{i_m}\cdots s_{i_{k+1}}\cdot \alpha_{i_k} := \alpha^{\bf i}_k.
$
Recall the function $\Gamma_{\lambda}$ in Definition \ref{function.gamma.alpha}. 
We have 
\[
\frac{\Gamma_{\lambda}\Bigl(\A_1, \A_2\cdot h, \A_3'\cdot w_0(h)\Bigr)}{\Gamma_{\lambda}\left(\A_1, \A_2, \A_3'\right)}=\left( [\lambda]_-- [\lambda]_+\right) (h)= \lambda(h^{-1}).
\]
By Proposition \ref{equivalence.cluster.a.ba}, we get 
\be
\la{cartan.action.pinning}
\frac{\tau^*P_{{\bf i},k}}{P_{{\bf i},k}} = \frac{\Gamma_{\gamma_2}\left(\A_1, \A_2\cdot h, \A_3'\cdot w_0(h)\right)}{\Gamma_{\gamma_1}\left(\A_1, \A_2\cdot h, \A_3'\cdot w_0(h)\right)}\cdot \frac{\Gamma_{\gamma_1}\left(\A_1, \A_2, \A_3'\right)}{\Gamma_{\gamma_2}\left(\A_1, \A_2, \A_3'\right)}=\left(\gamma_1-\gamma_2\right)(h_{23})= \alpha_k^{\bf i}(h_{23}^{-1}).
\ee
The Lemma is a special case when $\alpha_k^{\bf i}$ is a simple positive root.
\end{proof}

\subsubsection{Amalgamation of ${\mathscr P}-$spaces.} Let $q$ be a quadrilateral and ${\mathscr P}_{\G, q}$ be the space associated to $q$. 
There is an amalgamation map
\[
\gamma: ~ {\mathscr P}_{\G, t} \times {\mathscr P}_{\G,t}\longrightarrow {\mathscr P}_{\G, q}.
\]
The map glues a side $e_1$ of the first triangle and a side $e_2$ of the second triangle by identifying the pinnings $p_{e_1}$ and $p_{e_2}$.  
The frozen vertices associated to the side $e_1$, respectively $e_2$, are labelled by $i_{e_1}$, respectively $i_{e_2}$, for $i\in {\rm I}$.
The space ${\mathscr P}_{\G,q}$ inherits a cluster Poisson structure from ${\mathscr P}_{\G, t} \times {\mathscr P}_{\G,t}$ by amalgamating the vertices $i_{e_1}$ with $i^\ast_{e_2}$. Let $e$ be the diagonal of $q$ corresponding to $e_1$ and $e_2$.  There are $r$ many variables associated to the diagonal $e$ given by
\be
\la{x.variable.i_e}
X_{i_e}:= X_{i_{e_1}} X_{i_{e_2}^\ast}.
\ee
\begin{lemma} The function $X_{i_e}$ is a well-defined function of ${\mathscr P}_{\G,q}$.
\end{lemma}
\begin{proof}
For a generic point $x \in {\mathscr P}_{\G, q}$, the fiber $\gamma^{-1}(x)$ is an ${\rm H}$-torsor. Recall that the action \eqref{cartan.pinning.action.hhh} on a pinning depending on the orientation chosen: the rescaling of $\A$ by $h$ is equivalent to the rescaling of $\A'$ by $w_0(h)$.  Therefore if the pinning $p_{e_1}$ is rescaled by $h$, then correspondingly $p_{e_2}$ should be rescaled by $w_0(h)$. By Lemma \ref{rescale.hh}, the variables become
\be
X_{i_{e_1}} \lms X_{i_{e_1}} \alpha_i(h^{-1}),\hskip 7mm X_{i_{e_2}^\ast} \lms X_{i_{e_2}^\ast} \alpha_{i^\ast}(w_0(h^{-1})) = X_{i_{e_2}^\ast} \alpha_{i}(h)
\ee
Therefore $X_{i_e}\lms X_{i_e}$. In other words, $X_{i_e}$ descends to a function of ${\mathscr P}_{\G,q}$.
\end{proof}

More generally, given   a decorated surface $\bS$  and   a triangulation ${\rm T}$  of $\bS$, there is an amalgamation map 
\[
\gamma:~ \prod_{t \in {\rm T}} {\mathscr P}_{\G, t} \longrightarrow {\mathscr P}_{\G, \bS}
\]
For each triangle $t\in {\rm T}$,   pick an angle $v_t$ and a reduced word ${\bf i}$. They give rise to a cluster Poisson chart on each ${\mathscr P}_{\G',t}$. By amalgamation, the data $({\rm T}, v_t, {\bf i}_t)$ corresponds to a cluster Poisson chart of  ${\mathscr P}_{\G, \bS}$, so that
 \vskip 1mm
 \begin{enumerate}
\item To the inner part of each triangle $t\in {\rm T}$ is a associated $l(w_0)-r$ many variables, given by the unfrozen variables of corresponding $\mathscr{X}_{\G,t}= {\mathscr P}_{\G',t}/ H^3$.

 \vskip 1mm
\item To each inner edge $e$ of $T$ is associated $r$ many unfrozen variables, defined by \eqref{x.variable.i_e}.
 \vskip 1mm
 \item To each boundary edge $e$ of $T$ is associated $r$ many frozen variables, determined by the pinning $p_e$ and the space $\mathscr{X}_{\G', t}$ where $t \supset e$. 
\end{enumerate}
 \vskip 1mm
 Theorem \ref{MTHPs} asserts that the cluster Poisson structure on   ${\mathscr P}_{\G, \bS}$ is independent of the data $({\rm T}, \{{\bf i}_t\}, \{v_t\})$. We shall prove Theorem \ref{MTHPs} in Section \ref{SEC9.2x}.
\medskip
\subsection{Projection from $\mathscr{A}_{\G', \bS}$ to ${\mathscr P}_{\G, \bS}$}
\la{SEC9.2x}
\medskip

Let us consider  a  regular projection  
\be 
\la{2018.9.22.18.42hhX}
\begin{split}
\pi_3: ~& {\rm Conf}_3^\times(\mathcal{A}_{\G}) \lra {\mathscr P}_{\G, t},\\
&(\A_1, \A_2, \A_3) \lms(\B_1, \B_2, \B_3, p_1, p_2, p_3).\\
\end{split}
\ee
Here  $p_i=(\A_i, \A_{i+1}')$, where $\A_{i+1}'$ is the unique flag above $\B_{i+1}$ such that $(\A_i, \A_{i+1}')$ is a pinning. Defining $\pi_3$ we   broke the symmetry, and used the orientation of the triangle $t$. 

Recall the coordinates $\{A_k\}$ of ${\rm Conf}_3^\times(\mathcal{A}_{\G})$ labelled by the same quiver $\left(  {J}({\bf i}), \varepsilon({\bf i})=\varepsilon_{ij}  \right)$. The frozen variable associated to $i\in {\rm I}$ on the side $\{a, b\}$  is denoted by  $A_{i_{ab}}$. Recall the Cartan matrix $C_{ij}$. Let
\be \la{CORT}
\mathbb{A}_{i_{ab}}:=\alpha_i\Bigl(h(\A_a, \A_b)^{-1}\Bigr) = \prod_{j\in{\rm I}} A_{j_{ab}}^{-C_{ij}}.
\ee

\bp \la{11.24.18.1} For the projection \eqref{2018.9.22.18.42hhX}, we have 
\be \la{CORT1X}
\pi_3^*X_k =  \left\{    \begin{array}{ll} 
      {\prod_{j\in  {J}({\bf i})} A_j^{ \varepsilon_{kj}}} & \mbox{if $k$ is unfrozen}, \\
      &\\
 \mathbb{A}_{k}^{1/2} \cdot {\prod_{j\in  {J}({\bf i})} A_j^{ \varepsilon_{kj}}} & \mbox{if $k$ is frozen}. \\
    \end{array} \right.
\ee
\ep

\paragraph{\bf Remark.} Although   $\mathbb{A}_{k}^{1/2}$ has half integral exponents, $\pi_3^*X_k$ is a monomial in $A_j$, even if $k$ is frozen. 
Indeed, in (\ref{CORT1X})  the half integral exponents of $\mathbb{A}_{k}^{1/2}$ are compensated by the half integral values of the exchange matrix $\varepsilon_{ij}$. 

\begin{proof} 
If $k$ is unfrozen, then this is a special case of Theorem \ref{Th2.3}, part 4. 

If $k=i_{12}={i \choose 0}$, then it follows from Lemma \ref{bfz.lusztig.lemma}.

If $k=i_{31}={i^\ast \choose n_{i^\ast}}$, then by the last two identities in \eqref{hdist.pspace}, we have
\[
X_{i_{31}}=\Bigl( \mathcal{W}_{i^\ast}(\A_1, \B_{i^\ast \choose n_{i^\ast} -1},\, \B_{i^\ast \choose n_{i^\ast}})\cdot \mathbb{A}_{i_{31}}\Bigr)^{-1} .
\]
Then it again follows from Lemma \ref{bfz.lusztig.lemma}. 

Now we prove the case when $k=i_{23}$. Suppose that $\alpha_i^\vee = \beta_{l}^{\bf i}$ in \eqref{reduced.word.positive.roots}. Define
\[
\gamma_1 = s_{i_n}\ldots s_{i_{l}}\cdot \Lambda_{i_l}, \hskip 14mm 
\gamma_2 = s_{i_n}\ldots s_{i_{l+1}}\cdot \Lambda_{i_l}. 
\]
Note that $\gamma_{2}-\gamma_1= \alpha_i$, and 
\[
 \langle \alpha_i^\vee, \gamma_2 \rangle = \langle s_{i_n}\cdots s_{i_{l+1}} \cdot \alpha_{i_l}^\vee, s_{i_n}\ldots s_{i_{l+1}}\cdot \Lambda_{i_l} \rangle =1. 
\]
Let  ${\rm I}_1\subset {\rm I}$ be the subset parametrising simple positive coroots that appears before $\alpha_{i}^\vee$ in the sequence  $\beta_{l}^{\bf i}$ and let  
${\rm I}_2= {\rm I}_1\sqcup \{ i \}$.
By Lemma \ref{decomp.alpha.plus.minus}, we have
\begin{align}
[\gamma_2]_+ -[\gamma_1]_+ &= \sum_{j\in {\rm I}_2} \langle \alpha_j^\vee, \gamma_2 \rangle \Lambda_j - \sum_{j\in {\rm I}_1} \langle \alpha_j^\vee, \gamma_1 \rangle \Lambda_j  = \langle \alpha_i^\vee, \gamma_2 \rangle \Lambda_i + \sum_{j\in I _1} \langle \alpha_j^\vee,  \gamma_2-\gamma_1 \rangle \Lambda_j \nonumber\\
&= \Lambda_i + \sum_{j \in I_1} C_{ij} \Lambda_j. \nonumber
\end{align}
 Recall the pinning $p_{23}=(\A_2, \A_3')$. Let $\A_3= \A_3' \cdot h$. By \eqref{frozen.var.ed}, we have 
$
A_{j_{23}} =\Delta_j(\A_2, \A_3)= \Lambda_{j^\ast} (h).
$
By \eqref{bottomXY} and Proposition \ref{equivalence.cluster.a.ba}, we get 
\begin{align}
X_{i_{23}} &= \frac{\Gamma_{\gamma_2}\left(\A_1, \A_2, \A_3\cdot h^{-1}\right)} {\Gamma_{\gamma_1}\left(\A_1, \A_2, \A_3\cdot h^{-1}\right)} = \frac{\Gamma_{\gamma_2}\left(\A_1, \A_2, \A_3\right)}{\Gamma_{\gamma_1}\left(\A_1, \A_2, \A_3\right)}\cdot \left([\gamma_2]_+^\ast - [\gamma_1]_+^\ast\right)(h^{-1}) \nonumber\\
&=  \frac{\Gamma_{\gamma_2}\left(\A_1, \A_2, \A_3\right)}{\Gamma_{\gamma_1}\left(\A_1, \A_2, \A_3\right)} \cdot   A_{i_{23}}^{-1} \prod_{j\in {\rm I}_1} A_{j_{23}}^{-C_{ij}}. \nonumber
\end{align}
Note that $\Gamma_1$ and $\Gamma_2$ in the last line are cluster $K_2-$variables associated to the unfrozen vertices connecting $i_{23}$, and the rest are frozen variables on the side $\{2,3\}$.
It concludes the proof of the Proposition.
\end{proof}

Similarly, there is a natural map
\be
\la{CORT1X.general.nvadfo}
\pi_{\bS}:~ \mathscr{A}_{\G',\bS} \lra {\mathscr P}_{\G, \bS}
\ee
such that the pinning $p_{ab}$ associated to each boundary edge $ab$ of $\bS$ is determined by the first decorated flag $\A_a$. Recall the variable $\mathbb{A}_{i_{ab}}$ in \eqref{CORT}. 
Recall the data $({\rm T}, v_t, {\bf i}_t)$  inducing cluster structures on $\mathscr{A}_{\G',\bS}$ and ${\mathscr P}_{\G, \bS}$.
\bt \la{11.24.18.1} For the projection $\pi_{\bS}$, we have 
\be \la{CORT1X.general}
\pi_\bS^*X_k =  \left\{    \begin{array}{ll} 
      {\prod_{j} A_j^{ \varepsilon_{kj}}} & \mbox{if $k$ is unfrozen}, \\
      &\\
 \mathbb{A}_{k}^{1/2} \cdot {\prod_{j} A_j^{ \varepsilon_{kj}}} & \mbox{if $k=i_{ab}$ is frozen}. \\
    \end{array} \right.
\ee
\et

\begin{proof}
If $k$ is frozen or $k$ is placed in the inner part of a triangle $t\in {\rm T}$, then it is a direct consequence of Proposition \ref{11.24.18.1}. If $k$ is placed on an inner edge $e$ of ${\rm T}$, then it suffices to consider the quadrilateral $q$ containing $e$ as a diagonal. Locally, we have a map
\[
\pi_4: ~ \mathscr{A}_{\G', q}={\rm Conf}_4(\mathcal{A}) \lra {\mathscr P}_{\G, q},~~~~ (\A_1, \A_2, \A_3, \A_4) \lms (\B_1, \ldots, \B_4; p_1, \ldots, p_4)
\]
where $p_i=(\A_i, \B_{i+1})$. Let us consider the following sequence of maps
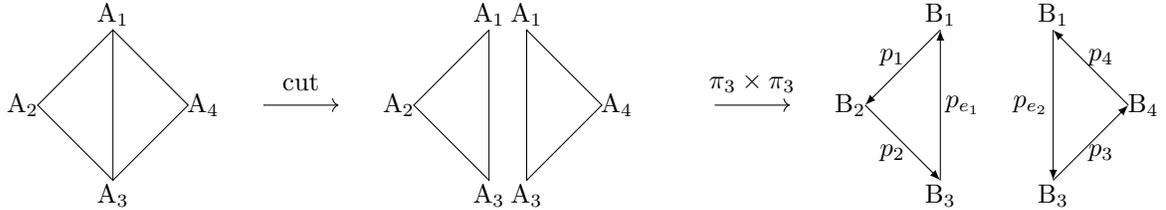
\begin{center}
\begin{tikzpicture}
\foreach \count in {1, 2, 3, 4}
{
\node at (90*\count:1.2) {{\small $\A_{\count}$}};
}
\draw  (90:1) -- (180:1) -- (270:1) -- (0:1) -- (90:1) -- (270:1);
\draw[->] (2,0) -- (3,0);
\node at (2.5,0.3) {\small {\rm cut}};
\draw[->] (8,0) -- (9,0);
\node at (8.5,0.3) {\small {$\pi_3\times \pi_3$}};
\draw[shift={(5,0)}]  (90:1) -- (180:1) -- (270:1) -- (90:1);
\foreach \count in {1, 2, 3}
{
\node[shift={(5,0)}] at (90*\count:1.2) {{\small $\A_{\count}$}};
}
\draw[shift={(5.5,0)}]  (270:1) -- (0:1) -- (90:1) -- (270:1);
\foreach \count in {1, 3, 4}
{
\node[shift={(5.5,0)}] at (90*\count:1.2) {{\small $\A_{\count}$}};
}
\draw[shift={(11,0)}, -latex]  (90:1) -- (180:1);
\draw[shift={(11,0)}, -latex]  (180:1) -- (270:1);
\draw[shift={(11,0)}, -latex]  (270:1) -- (90:1);
\foreach \count in {1, 2, 3}
{
\node[shift={(11,0)}] at (90*\count:1.2) {{\small $\B_{\count}$}};
}
\node[shift={(11,0)}] at (90*1+45:.9) {{\small $p_1$}};
\node[shift={(11,0)}] at (90*2+45:.9) {{\small $p_2$}};
\draw[shift={(12.5,0)},-latex]  (270:1) -- (0:1) ;
\draw[shift={(12.5,0)}, -latex]  (0:1) -- (90:1) ;
\draw[shift={(12.5,0)}, -latex]  (90:1) -- (270:1);
\node[shift={(12.5,0)}] at (90*3+45:.9) {{\small $p_3$}};
\node[shift={(12.5,0)}] at (90*4+45:.9) {{\small $p_4$}};
\node[shift={(11,0)}] at (.3,0) {{\small $p_{e_1}$}};
\node[shift={(11,0)}] at (1.2,0) {{\small $p_{e_2}$}};
\foreach \count in {1, 3, 4}
{
\node[shift={(12.5,0)}] at (90*\count:1.2) {{\small $\B_{\count}$}};
}
\end{tikzpicture}
\end{center}
Here $p_{e_1}=(\A_3, \B_1)$ and $p_{e_2}=(\A_1, \B_3)$. Let $h:= h(\A_3, \A_1)$. Note that to obtain the image of $\pi_4$, one has to alter the pinning $p_{e_1}$ by $h^{-1}$ first and then identify it with $p_{e_2}$. In other words,
\[
\pi_4= {\rm glue}\circ \tau_{h^{-1}}\circ (\pi_3\times \pi_3) \circ {\rm cut},
\]
where the map $\tau_{h^{-1}}$ alters the pinning $p_{e_1}$ by $h^{-1}$.
Combining Lemma \ref{rescale.hh} with \eqref{CORT1X}, we get
\[
\pi_4^\ast X_{i_e} = \alpha_i(h^{-1})(\pi_3\times \pi_3)^\ast (X_{i_{e_1}} X_{i_{e_2}^\ast}) = \prod_{j} A_j^{ \varepsilon_{i_e, j}}.
\]
\end{proof}

\paragraph{\bf Proof of Theorem \ref{MTHPs}.} For $\bS=t$, it suffices to prove the cyclic invariance for the cluster structure on ${\mathscr P}_{\G,t}$. By Theorem \ref{11.24.18.1}, we reduce it to the cyclic invariance of the cluster $K_2$ structure on ${\rm Conf}_3^\times(\mathcal{A})$. When $\bS=q$ is a quadrilateral, we   need to prove the flip invariance for  ${\mathscr P}_{\G,q}$. Using   Theorem \ref{11.24.18.1} again, we reduce it to the flip invariance of the cluster $K_2-$structure on ${\rm Conf}_4^\times(\mathcal{A})$. For general surfaces $\bS$, it follows from the amalgamation property of the space ${\mathscr P}_{\G, \bS}$.

\medskip \subsection{Cluster Poisson structure revisited} \la{SEC9.3x}

\medskip
In Section \ref{SEC9.3x} we introduce a different collection of coordinate systems on ${\mathscr P}_{\G, t}$. One of their benefits is that the Poisson brackets between the new coordinates 
in the cluster Poisson structure on ${\mathscr P}_{\G, t}$ 
are very simple, see Theorem \ref{TH9.5}.  

Recall that ${\mathscr P}_{\G, t}$ parametrizes configurations $(\B_1, \B_2, \B_3; p_{12}, p_{12}, p_{31})$ as  on Figure \ref{2018.9.23.15.10ss}. 
Let ${\bf i}=(i_1, \ldots, i_m)$ be a reduced word of $w_0$. It determines   primary coordinates $P_{{\bf i}, k}$ in \eqref{PC}. In this Section, we will fix an ${\bf i}$ and write $P_k$ instead of $P_{{\bf i}, k}$.
In addition to them, we consider   
\[
T_i:=\alpha_i(h(\A_1', \A_1)),\hskip 5mm  L_{i}= \alpha_i(h(\A_2', \A_2)), \hskip 5mm R_{i}=\alpha_i(h(\A_3', \A_3)), \hskip 6mm \forall i \in {\rm I}.
\]

\begin{lemma} \la{LEM9.4}
The following functions  form a coordinate system on ${\mathscr P}_{\G, t}$:
\be \la{FPKK}
\left\{P_{1}, \ldots, P_{ m}; L_1, \ldots, L_r; R_1, \ldots, R_r \right\}.
\ee
\end{lemma}

\begin{proof} 
Let us take $(\B_1, \A_2)=(\B, [\U^{-}])$. Recall the homomorphisms $\gamma_i$ in \eqref{5.21.2016.ss}. Define
\[
\phi_i(a) =\gamma_i    \begin{pmatrix} 
      a^{-1} & 1 \\
      0 & a \\
   \end{pmatrix}.
\]
Then we get
$
\A_{2}^{k}= \phi_{i_1}(P_{1})\ldots \phi_{i_k}(P_{k}) \cdot [\U^{-}]. $ 
In particular $\A_2^m=\A_3'$. The coordinates $L_1, \ldots, L_r$ recover the pinning $(\A_1, \A_2')$, and similarly $R_1, \ldots, R_r$ recover the pinning  $(\A_3, \A_1')$. Putting together, the coordinates \eqref{FPKK} determine an open embedding from $\mathbb{G}_m^{m+2r}$ into ${\mathscr P}_{\G,t}$, whose images are such that each $\A_2^k$ is of generic position to $\B_1$.
\end{proof}

\begin{lemma} 
The cluster Poisson coordinate $X_{i \choose k}$ is a monomial of primary coordinates:
\be
\label{X-P.transion}
X_{i \choose k} = \left\{ \begin{array}{ll} 
{ L_{i^\ast} P_{i \choose 1}^{-1} \prod_{s<{i \choose 1}} P_s^{-C_{i,i_s}},} & \mbox{if } k=0, \\
& \\
{ P_{i \choose k}^{-1} P_{i \choose k+1}^{-1} \prod_{{i \choose k}<s<{i \choose k+1}} P_s^{-C_{i,i_s}}} & \mbox{if } 0<k<n_i, \\
& \\
{ R_{i^\ast} P_{i \choose n_i}^{-1} \prod_{{i \choose n_i}<s} P_s^{-C_{i,i_s}}} & \mbox{if } k=n_i. \\
\end{array} \right.
\ee
\end{lemma}
\begin{proof}
Let us take $(\A_1, \A_2')=([{\rm U}], [{\rm U}^{-}])$. Let $g\in {\rm G}'$ be such that
$
(\A'_1, \A_3)= g\cdot ([{\rm U}], [{\rm U}^{-}]).
$

By the definition of cluster Poisson coordinates, we have
\be
\la{g.decomp.x.cord}
g=\Bigl(\prod_{i \in {\rm I}} H_i(X_{i \choose 0})\Bigr) {\bf E}_{i_1} H_{i_1}(X_1) \cdots {\bf E}_{i_m}H_{i_m}(X_m). 
\ee
Using the primary coordinates, we have
\be
\la{g.decomp.p.cord}
\begin{split}
&g= \Bigl(\prod_{i \in {\rm I}} H_i(L_{i^\ast})\Bigr) \phi_{i_1}(P_1) \cdots \phi_{i_m}(P_m) \cdot \prod_{i \in {\rm I}} H_i(R_{i^\ast}).\\
&\phi_i(P) = H_{i}(P^{-1}) {\bf E}_i H_i(P^{-1}) \prod_{j\neq i} H_j(P)^{-C_{ji}}.\\
\end{split}
\ee
Plugging it into \eqref{g.decomp.p.cord}, then comparing with \eqref{g.decomp.x.cord}, we prove the Lemma.
\end{proof}

Conversely, every primary coordinate $P_k$ can be expressed as a Laurent monomial of cluster Poisson variables. Moreover, $P_k$ does not depend on pinnings of sides $\{12\}$ and $\{13\}$. Therefore  
\be
\la{p-x.trans}
P_k = M_k \cdot \prod_{j \in {\rm I}} X_{j \choose -\infty}^{d_{k,j}},
\ee
where $M_s$ is a Laurent monomial of unfrozen cluster Poisson variables. Recall  the positive roots 
\[
\alpha_k^{\bf i}= s_{i_m}\ldots s_{i_{k+1}} \cdot \alpha_{i_k}.
\] 
\begin{lemma} The exponents $d_{k, j}$ in \eqref{p-x.trans} are such that
\be
\la{p-x.trans2}
\alpha_k^{\bf i} = \sum_{j \in {\rm I}}d_{k,j} \alpha_{j}.
\ee
\end{lemma}
\begin{proof} Recall the Cartan group action on the pinnings. 
By \eqref{cartan.action.pinning}, we have
$
\tau^\ast P_k = P_k \alpha_{k}^{\bf i}(h^{-1}). 
$ 
Meanwhile,
\[
\tau^\ast M_k =M_k, \hskip 5mm \tau^\ast X_{j \choose -\infty} = X_{j \choose -\infty} \alpha_{j}(h^{-1}).
\]
Comparing both sides of  \eqref{p-x.trans}, we prove the Lemma.
\end{proof}

Recall the Cartan matrix $C_{ij}$ with the multipliers $d_i=\langle\alpha_i^\vee, \alpha_i^\vee\rangle$. We consider a symmetric bilinear form on the root lattice such that 
\[
\left< \alpha_i, \alpha_j \right>= C_{ij} d_j^{-1}:=\widehat{C}_{ij}.
\] 


\bt \la{TH9.5} For the cluster Poisson structure on the space ${\mathscr P}_{\G', t}$, we have
\begin{align}
&\left\{ \log L_i, \log L_j\right\} = \left\{ \log R_i, \log R_j\right\}=0, \nonumber\\
&\left\{ \log L_i, \log R_{j^\ast}\right\} = \left< \alpha_i, \alpha_j\right>, \nonumber\\
&\left\{ \log P_{s}, \log P_{t} \right\} = {\rm sgn}(t-s) \cdot \langle \alpha_s^{\bf i}, \alpha_t^{\bf i} \rangle, \nonumber\\
&\left\{ \log L_{j}, \log P_{ s} \right\} = \left\{ \log P_s, \log R_{j^\ast}\right\} = \langle \alpha_j, \alpha_s^{\bf i}\rangle. \nonumber
\end{align}
\et

\begin{proof} i) Formula (\ref{hdist.pspace}) tells  us that
\be \la{F2.12}
T_i =\alpha_i(h({\A_1', \A_1}))=X_{i \choose 0}X_{i \choose 1}\cdots X_{i \choose n_{i}} .
\ee
Using this, and  looking at the quiver ${\bf Q}_{\alpha_i, \alpha_j}$ on Figure \ref{pin30}, we see that  $\{T_i, T_j\}=0$ 
for any $i, j \in {\rm I}$. Since the cluster Poisson bracket on the space ${\mathscr P}_{\G', t}$ is invariant under the rotation given by the cyclic shift of the vertices, we have 
$
\{L_i,L_j\}=\{R_i, R_j\}=0.
$
This  proves the two formulas in the first line of the claim. 

ii) The Poisson structure does not depend on the reduced word ${\bf i}$, so we assume that ${\bf i}$ starts with $i_1=j$. We have
\[
X_{j \choose 0}= \mathcal{W}_j(\A_1, \B_2, \B_2^1), \hskip 7mm X_{j^\ast \choose -\infty}=P_{1}= \frac{\Delta_{j}(\A_1, \A_2^1)}{\Delta_j(\A_1, \A_2)}.
\]
By \eqref{bfz.lusztig}, we get
\be
\la{leftkcoor}
X_{j \choose 0} X_{j^\ast \choose -\infty} = \frac{1}{\alpha_j(h(\A_1, \A_2))}=L_{j^\ast}.
\ee 
Note that $X_{j^\ast \choose -\infty}$ commutes with every $T_i$. Therefore
\[
\{\log T_i, \log L_{j^\ast}\}=\{\log T_i, \log X_{j \choose 0}\} = \widehat{C}_{ij},
\]
where the second identity follows from the pattern of the quiver ${\bf Q}_{\alpha_i, \alpha_j}$ again. By the cyclic invariance of the Poisson bracket, we get
\[
\left\{ \log L_i, \log R_{j^\ast}\right\}=\left \{\log T_i, \log L_{j^\ast}\right\}= \widehat{C}_{ij}.
\]

iii) Forgetting the pinnings $p_{12}$ and $p_{13}$ in ${\mathscr P}_{\G, t}$, we get a quotient space denoted by ${\mathscr P}_{\G, t}^{(1)}$. Note that ${\mathscr P}_{\G, t}^{(1)}$ inherits a cluster Poisson structure from  ${\mathscr P}_{\G, t}$ by taking all the unfrozen coordinates and the frozen coordinates associated to the side $\{23\}$. Following the proof of Lemma \ref{LEM9.4}, we conclude that $\{P_1, \ldots, P_m\}$ forms  a coordinate system on  ${\mathscr P}_{\G, t}^{(1)}$. 
Let us introduce a Poisson bracket $\{*,*\}_1$ by setting
\[
\left\{ \log P_s, \log P_t \right\}_1 = {\rm sgn}(t-s) \cdot \langle \alpha_s^{\bf i}, \alpha_{t}^{\bf i}\rangle. 
\]
We shall prove that $\{*,*\}_1$ coincides with the cluster Poisson structure on ${\mathscr P}_{\G, t}^{(1)}$.

\begin{lemma} 
\label{poisson.primary}
Let $p={i \choose k}$ be an unfrozen  vertex and let $q={i \choose k+1}$. Let $j \in \{1,\ldots, m\}$. We have
\[
\left\{ \log P_j, \log X_p \right\}_1 = \left\{ \begin{array}{ll} 2 d_i^{-1} &\mbox{if } j=p, \\
 - 2d_i^{-1} &\mbox{if }j =q,\\
0 &\mbox{otherwise}. \end{array}\right.
\]
\end{lemma}
\begin{proof} Let $p< t \leq m$. It follows easily by induction on $t$ that 
\be
\label{alphasumrel}
\alpha_{p}^{\bf i} -   s_{i_m}\ldots s_{i_t} ( \alpha_{i} )+ \sum_{ p < s < t} C_{i, i_s} \alpha_{s}^{\bf i} = 0
\ee
In particular, by setting $t=q$, we get
\be
\alpha_p^{\bf i} + \alpha_q^{\bf i} + \sum_{p<s< q} C_{i, i_s}\alpha_s^{\bf i}=0.
\ee
By \eqref{X-P.transion}, we have
\be
\log X_p = \log P_p + \log P_q + \sum_{p<s<q} C_{i, i_s}\log P_s.
\ee
Therefore we have 
\[
\left\{ \log P_j ,~ \log X_p \right\}_1 = \left\{\log P_j,~ \log P_p + \log P_q + \sum_{p<s<q} C_{i, i_s}\log P_s, \right\}_1.
\]
We split the proof of the Lemma   into the following cases.
\begin{itemize}
\vskip 1mm\item {\bf Case 1, $j<p$}. By the definition, we have
\[
\left\{ \log P_j,  \log X_p  \right\}_1 = -  \left< \alpha_j^{\bf i}, \alpha_p^{\bf i} + \alpha_q^{\bf i} + \sum_{p<s< q} C_{i, i_s}\alpha_s^{\bf i} \right> =0. \nonumber
 \]
\vskip 1mm \item {\bf Case 2, $j=p$}. Then 
 \[
 \left\{\log P_p,  \log X_p  \right\}_1 = - \left< \alpha_j^{\bf i},  \alpha_q^{\bf i} + \sum_{p<s< q} C_{i, i_s}\alpha_s^{\bf i} \right>=  \left< \alpha_p^{\bf i}, \alpha_p^{\bf i}\right>= 2 d_i^{-1}.
 \]
\vskip 1mm \item {\bf Case 3, $p<j<q$.} By \eqref{alphasumrel}, we get
 \begin{align}
 & \left\{ \log P_j,  \log X_p  \right\}_1 \nonumber\\
= & - \left\{\log P_j,~ \log P_p + \sum_{p<s<j} C_{i, i_s}\log P_s\right\}_1 -  \left\{\log P_j, ~\log P_q + \sum_{j<s<q} C_{i, i_s}\log P_s \right\}_1 \nonumber\\
=  & - \left< \alpha_j^{\bf i},~ s_{i_m}\ldots s_{i_{j}}(\alpha_i)\right> - \left< \alpha_j^{\bf i}, s_{i_m}\ldots s_{i_{j}}(\alpha_i)+ C_{i, i_j}\alpha_j^{\bf i} \right> \nonumber\\
=& -2\left< \alpha_j^{\bf i},~ s_{i_m}\ldots s_{i_{j}}(\alpha_i)\right> -C_{i, i_j}\left< \alpha_j^{\bf i}, ~\alpha_j^{\bf i} \right> \nonumber\\
 = &  2 \left<  \alpha_{i_j}, \alpha_i\right> - C_{i, i_j}\left< \alpha_{i_j}, \alpha_{i_j}\right>=0. \nonumber
 \end{align}
 \item {\bf Case 4, $j=q$}. By a similar argument to Case 2, we get 
$
 \left\{\log P_q,  \log X_p \right\}_1=- 2 d_i^{-1}.
$
\vskip 1mm \item {\bf Case 5, $q<j$.} By a similar argument to Case 1, we get 
$
 \left\{  \log P_j , \log X_p\right\}_1=0.
$
\end{itemize} 
\end{proof}

Recall the subspace ${\rm Conf}_3(\mathcal{A})_{e, e}\subset {\rm Conf}_3(\mathcal{A})$   parametrizing $\G'$-orbits of triples $(\A_1, \A_2, \A_3)$ with 
\[
h(\A_1,\A_2)=h(\A_2, \A_3)=1.
\] 
The space ${\rm Conf}_3(\mathcal{A})_{e, e}$ inherits a   cluster $K_2-$structure from ${\rm Conf}_3(\mathcal{A})$ with frozen variables associated to the sides $\{12\}$ and $\{23\}$ being 1. Furthermore, there is a natural isomorphism
\[
{\rm Conf}_3(\mathcal{A}_{e,e})\lra {\mathscr P}_{\G',t}^{(1)}, \hskip 7mm (\A_1, \A_2, \A_3) \lms (\B_1, \B_2, \B_3; p_{23}), \hskip 3mm \mbox{where } p_{23}=(\A_2, \A_3).
\]
Let us identify the above two spaces.
For $s={i \choose k}$, the cluster $K_2$ coordinate on ${\rm Conf}_3(\mathcal{A})_{e, e}$ is 
\be \la{CA}
A_s:={\Delta_{i} (\A_1, \A_2^k)} = \prod_{1\leq j \leq k} P_{i \choose j}.
\ee

\begin{lemma} Let $t$ be a non-frozen vertex and $s={i \choose k}$. We have 
\be \la{Th9.8}
\left\{ \log A_s, \log X_t \right\}_1 = 2 d_i^{-1}\delta_{st}.
\ee
\end{lemma}
\begin{proof} It follows from Lemma \ref{poisson.primary} and (\ref{CA}).
\end{proof}

The third part of Theorem \ref{TH9.5} is equivalent to the following result.
\bt
The Poisson structure $\{~,~\}_1$ on ${\mathscr P}_{\G,t}^{(1)}$ coincides with the cluster Poisson structure:
\[
\left\{ \log X_s, \log X_t \right\}_1 = 2 \widehat {\varepsilon}_{st}.
\]
\et
\begin{proof} If $s, t$ are both frozen, note the frozen coordinates $X_s$ associated to the side $\{23\}$ are precisely the coordinates $P_k$ with $\alpha_k^{\bf i}$ being simple roots. Hence the cluster Poisson bracket on them matches $\{*,*\}_1$. 

If one of $s,t$ is frozen, then this is a special case of Lemma \ref{poisson.primary}.

If $s,t$ are both unfrozen, it follows by substituting  $X_s = \prod_t A_t^{\varepsilon_{st}}$ to \eqref{Th9.8}.
\end{proof}

iv) Note that $T_i$ communicates with every unfrozen variable $X_s$, and 
\[
\left\{ \log T_{i}, \log X_{j \choose 0} \right\} = \widehat{C}_{ij}
\]
Recall the  cyclic rotation $r$ on ${\mathscr P}_{\G',t}$ and the frozen variables $X_{i \choose -\infty}$ associated to the side $\{23\}$. Since $r$ is a cluster automorphism, we have 
\[
r^\ast L_i = T_i, ~~~~~~~~~ r^\ast X_{i \choose -\infty} = X_{i \choose 0} \cdot Q_i
\]
where $Q_i$ are rational functions of unfrozen variables. Note that $r$ preserves the Poisson structure.  Therefore
\be
\left\{ \log L_{i}, \log X_{j \choose -\infty} \right\} = \left\{ \log T_{i}, \log X_{j \choose 0}\right\}+ \left\{ \log T_i, \log Q_j\right\}= \widehat{C}_{ij}
\ee
By \eqref{p-x.trans} and \eqref{p-x.trans2}, we get
\[\left\{ \log L_{i}, \log P_s \right\} = \left\{ \log L_{i}, \log M_s\right\} +\sum_{j\in {\rm I}} d_{s,j} \left\{ \log L_{i},  \log X_{j \choose -\infty}\right\} = \langle \alpha_i, \alpha_s^{\bf i} \rangle.
\]
This proves the first formula in the formula line of the claim of Theorem \ref{TH9.5}. The second formula follows by the same argument.
\end{proof}

The following result is a complement to Lemma \ref{poisson.primary}
\begin{lemma} 
\label{poisson.primary.comple}
 Let $j \in \{1,\ldots, m\}$. We have
\[
\left\{ \log P_j, \log X_{i \choose 0} \right\} = \left\{ \begin{array}{ll} -2 d_i^{-1} &\mbox{if } j={i \choose 1}, \\
0 &\mbox{otherwise}. \end{array}\right.
\]
\[
\left\{ \log P_j, \log X_{i \choose n_i} \right\} = \left\{ \begin{array}{ll} 2 d_i^{-1} &\mbox{if } j={i \choose n_i}, \\
0 &\mbox{otherwise}. \end{array}\right.
\]
\end{lemma}
\begin{proof}
To prove the first formula, we will use the identity
\[
-\alpha_{i^\ast} -   s_{i_m}\ldots s_{i_t} ( \alpha_{i} )+ \sum_{  s < t} C_{i, i_s} \alpha_{s}^{\bf i} = 0
\]
The rest of the proof goes through the same line as the proof of Lemma \ref{poisson.primary}. The second formula follows by a  similar argument.
\end{proof}

Consider the primary coordinates of the moduli space ${\mathscr P}_{\G, q}$, where $q$ is a quadrilateral. Every reduced word ${\bf i}=(i_1, \ldots, i_{2m})$ of $(w_0, w_0)$ provides a decomposition of ${\mathscr P}_{\G, q}$. Here is an illustrating example.
\begin{example} Let $\G={\rm PGL}_3$. The word ${\bf i}=(1, \overline{1}, 2,1, \overline{2}, \overline{1})$ gives rise to a decomposition of ${\mathscr P}_{{\rm PGL}_3, q}$. Here the top pinning is $p_{41}=(\A_4, \A_1')$, and the bottom pinning is $p_{23}=(\A_2, \A_3')$.
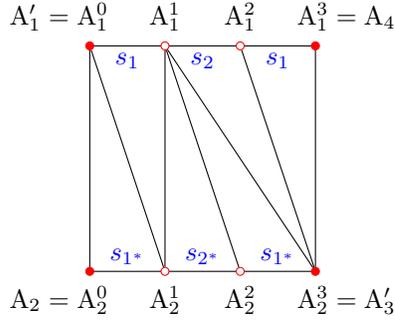
\begin{center}
\begin{tikzpicture}
\draw (1.5, 1.5) -- (-1.5, 1.5) -- (-1.5,-1.5) -- (1.5,-1.5) -- (1.5, 1.5);
\node[circle,draw=red,fill=red,minimum size=3pt,inner sep=0pt, label=above: {\small $\A_1'=\A_1^0~~~~~~~$}] (a1) at (-1.5, 1.5) {};
\node [circle,draw=red,fill=white,minimum size=3pt,inner sep=0pt, label=above: {\small $\A_1^1$}]  (a2) at (-.5, 1.5) {};
\node [circle,draw=red,fill=white,minimum size=3pt,inner sep=0pt, label=above: {\small $\A_1^2$}]  (a3) at (.5, 1.5) {};
\node[circle,draw=red,fill=red,minimum size=3pt,inner sep=0pt, label=above: {\small $~~~~~~~\A_1^3=\A_4$}]  (a4) at (1.5, 1.5) {};
\node[circle,draw=red,fill=red,minimum size=3pt,inner sep=0pt, label=below: {\small $\A_2=\A_2^0~~~~~~~$}]  (b1) at (-1.5, -1.5) {};
\node [circle,draw=red,fill=white,minimum size=3pt,inner sep=0pt, label=below: {\small $\A_2^1$}]  (b2) at (-.5, -1.5) {};
\node [circle,draw=red,fill=white,minimum size=3pt,inner sep=0pt, label=below: {\small $\A_2^2$}]  (b3) at (.5, -1.5) {};
\node[circle,draw=red,fill=red,minimum size=3pt,inner sep=0pt, label=below: {\small $~~~~~~~\A_2^3=\A_3'$}]  (b4) at (1.5, -1.5) {};
\foreach \from/\to in {a1/b2, a2/b2, a2/b3, a2/b4, a3/b4}
{\draw (\from) --(\to);}
\node [blue] at (-1, 1.3) {{\small $s_1$}};
\node [blue] at (0, 1.3) {{\small $s_2$}};
\node [blue] at (1, 1.3) {{\small $s_1$}};
\node [blue] at (-1, -1.3) {{\small $s_{1^\ast}$}};
\node [blue] at (0, -1.3) {{\small $s_{2^\ast}$}};
\node [blue] at (1, -1.3) {{\small $s_{1^\ast}$}};
\end{tikzpicture}
\end{center}
The space are decomposed into the following two types of elementary triangles
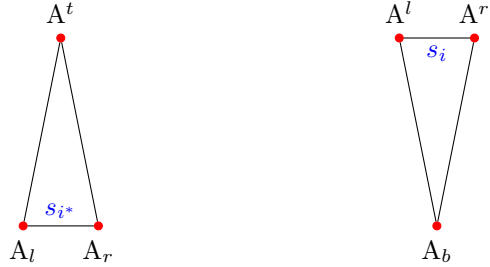
\begin{center}
\begin{tikzpicture}
\node[circle,draw=red,fill=red,minimum size=3pt,inner sep=0pt, label=above: {\small $\A^t$}] (a1) at (0, 2.5) {};
\node [circle,draw=red,fill=red,minimum size=3pt,inner sep=0pt, label=below: {\small $\A_l$}]  (a2) at (-.5, 0) {};
\node [circle,draw=red,fill=red,minimum size=3pt,inner sep=0pt, label=below: {\small $\A_r$}]  (a3) at (.5, 0) {};
\node [blue] at (0, .2) {{\small $s_{i^\ast}$}};
\draw (a1)--(a2)--(a3)--(a1);
\node[circle,draw=red,fill=red,minimum size=3pt,inner sep=0pt, label=below: {\small $\A_b$}] (b1) at (5, 0) {};
\node [circle,draw=red,fill=red,minimum size=3pt,inner sep=0pt, label=above: {\small $\A^l$}]  (b2) at (4.5, 2.5) {};
\node [circle,draw=red,fill=red,minimum size=3pt,inner sep=0pt, label=above: {\small $\A^r$}]  (b3) at (5.5, 2.5) {};
\node [blue] at (5, 2.3) {{\small $s_{i}$}};
\draw (b1)--(b2)--(b3)--(b1);
\end{tikzpicture}
\end{center}
The primary coordinates associated to the elementary triangles are 
\[
P= \frac{\Delta_i(\A^t, \A_r)}{\Delta_{i}(\A^t, \A_l)}, \hskip 1cm P'= \frac{\Delta_i(\A^l, \A_b)}{\Delta_{i}(\A^r, \A_b)}.
\]
Therefore we denote  the primary coordinates of ${\mathscr P}_{\G',q}$ from left to right are $P_{ 1}, ..., P_{ 2m}$.
\end{example}

Let ${\bf i}=(i_1, \ldots, i_{2m})$ be a reduced word of $(w_0, w_0)$. Recall the quiver ${\bf Q}({\bf i})$ in Definition \ref{quiver.amla.fg4.sec22}. 
Associated to each vertex of ${\bf Q}({\bf i})$ is a rational function of ${\mathscr P}_{\G', q}$: 
\begin{itemize}
\vskip 1mm\item the $s$th vertex of ${\bf K}({\bf i})$ is associated the primary coordinate $P_{s}$,
\vskip 1mm\item the variable ${i \choose k}$ is  associated the cluster Poisson variable $X_{{i \choose k}}$.
\end{itemize}
Recall the exchange matrix $\widehat{\varepsilon}_{jk}= {\varepsilon}_{jk}\cdot d_k^{-1}$.
\begin{proposition} \la{9.15}
\label{switchi,ibar.}
Let $X_{ j}$ be the   variable $P_{s}$ or  $X_{ {i \choose k}}$ associated to the vertex $j$ of ${\bf Q}({\bf i})$, we have 
\[
\left\{ \log X_{ j}, \log X_{k} \right\}= 2 \widehat{\varepsilon}_{jk}.
\]
\end{proposition}
\begin{remark} The functions $X_{j}$ do not form a coordinate system of ${\mathscr P}_{\G,q}$ since they are not algebraically independent. However  Proposition \ref{9.15} explicitly describes the Poisson structure of ${\mathscr P}_{\G, q}$. 
\end{remark}
\begin{proof} If ${\bf i}$ is a reduced word such that $i_1, \ldots, i_m \in {\rm I}$ and $i_{m+1}, \ldots, i_{2m}\in \overline{\rm I}$, then the coordinates are obtained by amalgamation of two copies of ${\mathscr P}_{\G', t}$. The Proposition follows from Theorem \ref{TH9.5}  $\&$ Lemmas \ref{poisson.primary},   \ref{poisson.primary.comple}.

Now let ${\bf i}=(\ldots, i_{s}, i_{s+1}, \ldots)=(\ldots, i, \overline{i}, \ldots)$ and let ${\bf i}'$ be the reduced word obtained by switching $i \leftrightarrow \overline{i}$. Following the proof of case 1 in Section \ref{thm.3.6.st}, we see that the quiver ${\bf Q}({\bf i}')$ is obtained from ${\bf Q}({\bf i})$ by a mutation at its $s$-th unfrozen vertex and then switching the $s$-th and $(s+1)$-th extra frozen vertices. Let $\{X_j\}$ and $\{X_j'\}$ be the functions associated to the vertices of ${\bf Q}({\bf i})$ and ${\bf Q}({\bf i}')$. It suffices to show that change of variables from $\{X_j\}$ to $\{X_j'\}$ satisfies the cluster Poisson mutation formula given by changing the   quiver ${\bf Q}({\bf i})$    to ${\bf Q}({\bf i}')$.

If $X_j'= X_{i \choose k}'$, then by definition they are cluster Poisson coordinates, and therefore satisfy the corresponding cluster Poisson mutation formula.

If $X_j'= P_k'$, where $k\notin \{s, s+1\}$, then by definition we have $P_k'=P_k$, and the formula follows.

By Corollary \ref{local.primary}, we get 
$
P_s'= P_{s+1}(1+X_s^{-1})^{-1}$ and $P_{s+1}'= P_{s}(1+X_s^{-1})^{-1},
$ 
which also satisfies the corresponding cluster mutation formula.
\end{proof}

\subsection{Proof of Theorem \ref{4.20.24.1}}\la{SECT13.4}

\begin{figure}[ht]
\begin{tikzpicture}
\draw[dashed, thick] (0,0)--(2,0)--+(60:2)--(2.25, 3.15);
\draw[dashed, thick] (2,0)+(60:2)--(4.5,1.732);
\draw[dashed, thick] (2,0)--+(-60:1.5);
\draw[dashed, thick] (0,0)--+ (120:1.5);
\draw[dashed, thick] (0,0)--+(-120:1.5);
\draw[-latex, thick] (60:1.5)--(-60:1.5);
\draw[-latex, thick] (-60:1.5)--(-1.5,0);
\draw[-latex, thick] (-1.5,0)--(60:1.5);
\draw[-latex, thick] (2,0)+(120:1.5)--+(1.5,0);
\draw[-latex, thick] (2,0)+(1.5,0)--+(-120:1.5);
\draw[-latex, thick] (2,0)+(-120:1.5)--+(120:1.5);
\draw[-latex, thick] (3,1.732)+(60:1.5)--+(-60:1.5);
\draw[-latex, thick] (3,1.732)+(-60:1.5)--+(-1.5,0);
\draw[-latex, thick] (3,1.732)+(-1.5,0)--+(60:1.5);
\draw[latex-, thick] (3.5,1.732)+(60:1.5)--+(-60:1.5);
\node (a) at (1.14, 0.5) {};
\node (b) at (2.48, 1.35) {};
\node (c) at (4.35,2.432) {};
\node (d) at (2.33, 1.1) {};
\draw[-latex, thick, red] (a) edge [bend right] (b);
\draw[-latex, thick, blue] (d) edge [bend left] (c);
\node[red, thick] at (1.8, 0.3) {$u$};
\node[blue, thick] at (3, 2.2) {$v$};

\draw[thick, latex-] (8,1) -- + (-120:1.8);
\draw[thick, latex-] (8,1) -- + (-60:1.8);
\node[red, thick] at (8, 0) {$u$};
\node (e) at (7.5,0.5) {};
\node (f) at (8.5, 0.5) {};
\draw[-latex, thick, red] (e) edge [bend right] (f);

\draw[thick, -latex] (12,-0.55) -- + (60:1.8);
\draw[thick, -latex] (12,-0.55) -- + (120:1.8);
\node[blue, thick] at (12, 0.5) {$v$};
\node (g) at (11.5,0) {};
\node (h) at (12.5, 0) {};
\draw[-latex, thick, blue] (g) edge [bend left] (h);

\end{tikzpicture}
\caption{Every loop along the punctured graph $\Gamma$ on $\bS$ is a composition of elementary arcs. Two of them are shown on the picture. Each elementary arc meets the sides of the  triangles at exactly three points: the origin, the end, and the internal intersection point.}
\label{13.4}
\end{figure}

\begin{proof} Take an ideal triangulation ${\cal T}$ of $\bS$. Consider the dual graph $\Gamma$, shown by the punctured lines on Figure \ref{13.4}. For each vertex $x$ of  $\Gamma$ consider a tiny bit smaller clockwise 
oriented triangle ${\rm T}_x$, obtained by shrinking the original triangle of the triangulation ${\cal T}$.  Its
 sides are perpendicular to the edges of $\Gamma$ sharing the vertex $x$, shown by green on Figure \ref{13.4}. 
So for each edge $E =(xy)$ of  $\Gamma$ there are two oriented sides of the  triangles ${\rm T}_x$ and  ${\rm T}_y$ intersecting the edge.  

A $\G-$local system ${\cal L}$ parametrised by a point of the moduli space $\mathscr{P}_{\G, \bS}$ can be described by a collection of pinnings 
assigned to the oriented sides of the  triangles ${\rm T}_x$ as follows. Take a pinning $p$ assigned to a side $s$ of a triangle ${\rm T}_x$. 
For each of the other two sides of the triangle ${\rm T}_x$, take  anther side intersecting the same edge of $\Gamma$. Let $s'$ be one of them.  
Let $q$ be the pinning at the side $s'$.  
Then the parallel transport from $s$ to $s'$   is described by the unique element $g_{s\to s'}\in \G$ such that $q = g_{s \to s'} p$. Given an oriented path on the graph $\Gamma$,  take the product 
of these elements of $\G$ along the path. 

The oriented sides $s$ and $s'$, viewed as the original sides of the triangulation ${\cal T}$, share a common vertex of ${\cal T}$. They are oriented either from  the common vertex, or to the common vertex, see Figure \ref{13.4}.  
In the first case, the element $g_{s \to s'}$ lies  in the Borel subgroup $\B$.  Otherwise it lies in the opposite Borel subgroup $\B_-$. 
 Any loop $l$ on  $\Gamma$ is a concatenation of elementary arcs of two types, $u$ and $v$, 
  each  connecting a pair of sides $s$ and $s'$ as above. It is a $u$-arc if  $s$ and $s'$ are oriented towards the common vertex, and a $v$-arc otherwise. 
 We conclude that: 
{\it The parallel transport  in   ${\cal L}$ of the pinning $p$ to the pinning $q$  along a path $l$ is   a described by the product $t_{p\to q}$ of the elements $g_{s\to s'}$. 
  Each of them lies either in $\B$ or in $\B_-$.}
 
For example, if $l= u_1 u_2 v_3 u_4 v_5$, the parallel transport  along the path $l$ is described by the element
$$
t_{p\to q}  = \U_{1}\U_{2}{\rm V}_{3}\U_{4}{\rm V}_{5}, \ \ \ \ \ \ \U_* \in \B, \ \ {\rm V}_* \in \B_-.
$$

 To prove Theorem \ref{4.20.24.1}, it remains to note that,  by the definition of the cluster Poisson coordinates, the elements $\U_*$ and ${\rm V}_*$ are products 
of   elementary unipotent elements ${\bf E}_i$, ${\bf F}_i$ and elements of the Cartan subgroup  $H_j(X_k)$ corresponding to the cluster Poisson coordinates $X_k$, see (\ref{FG3EFH}) and  \cite{FG05}. The latter act
 in Lusztig's canonical basis  by the diagonal matrices. The former act by unipotent matrices with the entries $1$ or $0$. Theorem \ref{4.20.24.1} follows.  \end{proof}


\medskip

 \section{Weyl,  braid and {${\rm Out}(\G)$} group 
actions  on   spaces   ${\mathscr P}_{\G,  \bS}$ and ${\mathscr A}_{\G,  \bS}$} \la{SEC5}

\medskip
\subsection{Preliminaries} \la{SEC12.1b}
\medskip

 \subsubsection{Shapes of some quivers.} We start from  general remarks on quivers defined by the amalgamation   in Section \ref{sec.3.1}, which  
 play key role in the    proofs in Sections \ref{Sec5.2} $\&$ \ref{SSEECC11.3}. 

A word  ${\bf i}$ representing a braid group element determines a quiver ${\bf Q}$, 
  equipped with a {\it level map} 
\be \la{LVLM}
l_{{\bf Q}}: {\rm V}_{{\bf Q}} \lra {\rm I} \cup \bullet.
\ee
Here  $\bullet$ is an extra element,  
so that $l_{{\bf Q}}^{-1}(\bullet)$ is the sub-quiver of the frozen elements ``on the bottom". 

1. For any simple positive root $\alpha$ 
there is a quiver ${\bf Q}_{\alpha}:= l_{{\bf Q}}^{-1}(\alpha) \subset {\bf Q}$ of the vertices on the level $\alpha$. 
It is a type $A$ quiver with the vertices   oriented from the right to the left. 

In particular, the quiver ${\bf Q}_{\alpha_i}$ for a reduced 
word ${\bf i}$ for $w_0$ is a quiver of type $A_{n_i+1}$:

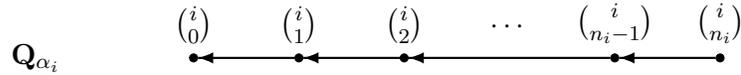
\begin{center} 
\begin{tikzpicture}[scale=1.4]
\node []  at (-1.5,0) {${\bf Q}_{\alpha_i}$};
\node [circle,draw,fill,minimum size=3pt,inner sep=0pt, label=above:${i \choose 0}$] (a) at (0,0) {};
\node [circle,draw,fill,minimum size=3pt,inner sep=0pt, label=above:${i \choose 1}$] (b) at (1,0) {};
 \node [circle,draw,fill, minimum size=3pt,inner sep=0pt, label=above:${i \choose 2}$]  (c) at (2,0){};
  \node [circle,draw,fill, minimum size=3pt,inner sep=0pt, label=above:${i \choose n_i-1}$]  (d) at (4,0){};
 \node [circle,draw,fill,minimum size=3pt,inner sep=0pt, label=above:${i \choose n_i}$] (e) at (5,0) {};
 \foreach \from/\to in {b/a, c/b, e/d, d/c}
                  \draw[directed, thick] (\from) -- (\to);
\node at (3,0.3) {$\cdots$};    
     \end{tikzpicture}
    \end{center}

2. For any  pair of simple positive roots $\alpha, \beta$ there is a sub-quiver 
${\bf Q}_{\alpha, \beta} := l_{{\bf Q}}^{-1}(\alpha, \beta)$. 
It has the following shape, where the weight of a solid  arrow is 
$ C_{\alpha \beta}$, and the punctured one  $ C_{\alpha \beta}/2$:

\begin{figure}[ht]
\centerline{
\begin{tikzpicture}[scale=0.8]
\node []  at (-3,-0.65) {${\bf Q}_{\alpha\beta}$};
\node []  at (7.3,0) {$\alpha$};
\node []  at (7.3, -1.3) {$\beta$};
\node[circle,draw, fill, minimum size=3pt,inner sep=0pt] (a1) at (0, 0) {};
\node[circle,draw, fill, minimum size=3pt,inner sep=0pt] (a2) at (1,0) {};
\node[circle,draw, fill, minimum size=3pt,inner sep=0pt] (a3) at (2, 0) {};
\node[circle,draw, fill, minimum size=3pt,inner sep=0pt] (a4) at (3, 0) {};
\node[circle,draw, fill, minimum size=3pt,inner sep=0pt] (a5) at (4, 0) {};
\node[circle,draw, fill, minimum size=3pt,inner sep=0pt] (a6) at (5,0) {};
\node[circle,draw, fill, minimum size=3pt,inner sep=0pt] (a7) at (6,0) {};
\node[circle,draw, fill, minimum size=3pt,inner sep=0pt] (b1) at (0.5, -1.3) {};
\node[circle,draw, fill, minimum size=3pt,inner sep=0pt] (b2) at (1.5,-1.3) {};
\node[circle,draw, fill, minimum size=3pt,inner sep=0pt] (b3) at (2.5, -1.3) {};
\node[circle,draw, fill, minimum size=3pt,inner sep=0pt] (b4) at (3.5, -1.3) {};
\node[circle,draw, fill, minimum size=3pt,inner sep=0pt] (b5) at (4.5, -1.3) {};
\node[circle,draw, fill, minimum size=3pt,inner sep=0pt] (b6) at (5.5,-1.3) {};
 \foreach \from/\to in {a7/a6, a6/a5, a5/a4, a4/a3, a3/a2, a2/a1, b6/b5, b5/b4, b4/b3, b3/b2, b2/b1, b1/a3, a3/b3, b3/a4, a4/b5, b5/a6, a6/b6}
                  \draw[directed] (\from) -- (\to); 
\foreach \from/\to in {a1/b1, b6/a7}
                  \draw[directed, dashed] (\from) -- (\to);                                       
 \draw[dotted, blue] (-1,0) -- (7,0);      
  \draw[dotted, blue] (-1,-1.3) -- (7,-1.3);         
     \end{tikzpicture}}
\caption{The shape of the quiver ${\bf Q}_{\alpha, \beta}$.}
\label{pin30}
\end{figure}

\subsubsection{Notation.}  \la{SECT14.1.2}
Let $\Delta$ be a root system with 
    the set of  positive simple roots ${\rm I}$ and the Cartan matirx $C_{ij}$. 
Let ${\mathbb B}_\Delta $ (respectively ${\mathbb B}_\Delta^+$) 
be the braid group (respectively semigroup) of $\Delta$. 
It is generated, as a group (respectively as a semigroup), 
 by the elements $s_i, i\in {\rm I}$, subject to the relations 
\begin{equation} \label{1}
\begin{array}{lcl}
s_is_j = s_js_i &\mbox{if}& C_{i j}=C_{ji}=0,\\
s_is_js_i = s_js_is_j &\mbox{if}& C_{i j}=C_{ji}=-1,\\
s_is_js_is_j = s_js_is_js_i &\mbox{if}& C_{ij}=2C_{ji}=-2,\\
s_is_js_is_js_is_j = s_js_is_js_is_js_i &\mbox{if}& C_{ij}=3C_{ji}=-3.
\end{array}
\end{equation}
 The 
Weyl group $W$ is the quotient of ${\mathbb B}_\Delta $ by the relations $s_i^2=1$. Denote by  $l(*)$   the length function on $W$. There is a   set theoretic section $\mu: W \to {\mathbb B}_\Delta^+$, 
such that 
\be \la{Smu}
\mu(ss') = \mu(s)\mu(s') \qquad \mbox{if $l(ss') = l(s)+l(s')$},
\ee
Abusing notation, we denote by $s_i$ the generators in $W$ and the elements $\mu(s_i) \in {\mathbb B}_\Delta^+$. 

 We   use the notation ${\mathbb B}_\Delta$  or ${\mathbb B}_{\mathfrak g}$ for the braid group assigned to  
the root system $\Delta$ of a semi-simple Lie algebra ${\mathfrak g}$, and the notation ${\mathbb B}$ when the root system is unambiguous.

If 
  $\Delta$ is irreducible of rank $>1$, 
 the  element 
 $
Z = \mu(w_0)^2 
\in {\mathbb B}^+
 $ 
generates the center of  the braid group ${\mathbb B}$, and of the braid semigroup. 
  For any presentation $\widetilde C$ of the Coxeter element, 
 $
Z = \mu(\widetilde C)^h. 
 $
The quotient  $\overline {\mathbb B}:= {\mathbb B}  / \langle Z\rangle$ is  called the {\it reduced braid group}. 

 \medskip
 
 \subsection{The group ${\rm Out}(\G)$ action on the moduli spaces ${\mathscr P}_{\G, \bS}$ and  ${\mathscr A}_{\G, \bS}$ is cluster} \la{Sec5.2*}
 
\medskip 

The group of outer automorphisms ${\rm Out}(\G)$ of the group $\G$ acts naturally on each of the moduli spaces ${\mathscr P}_{\G, \bS}$ and  ${\mathscr A}_{\G, \bS}$. 
 Indeed,  the group of inner automorphims of $\G$ acts trivially. These spaces are equipped with cluster structures of the relevant flavor.

 \bt \la{OCL} The  ${\rm Out}(\G)-$action on  the  cluster varieties  ${\mathscr P}_{\G, \bS}$ and  ${\mathscr A}_{\G, \bS}$  is cluster. 
 \et 

\begin{proof} It is sufficient to prove the claim for  the ${\mathscr P}$ and ${\mathscr A}$ spaces. 
Next, it is sufficient to prove that the action of ${\rm Out}(\G)$ is cluster for the spaces ${\mathscr P}_{\G, t}$ and ${\mathscr A}_{\G, t}$ related to a triangle $t$. 
 Indeed, the spaces ${\mathscr P}_{\G, \bS}$ and ${\mathscr A}_{\G, \bS}$ are obtained by amalgamation of the spaces related to  triangles. The action of ${\rm Out}(\G)$ 
 commutes with the amalgamation. Take a reduced decomposition ${\bf i}= (i_1, ..., i_m)$ of $w_0$. It gives rise to a cluster coordinate system ${\cal C}^*_{{\bf i}, t}$, where $\ast = {\mathscr P}$ or 
 $\ast =   {\mathscr A}$, 
  on each of the spaces ${\mathscr P}_{\G, t}$ and ${\mathscr A}_{\G, t}$, obtained by amalgamation of elementary cluster varieties ${\cal C}^*_{{\bf i}_k}$, followed up by amalgamation with 
  the Cartan group.  An element  $g \in {\rm Out}(\G)$ transform it to a cluster coordinate system ${\cal C}^*_{g({\bf i}), t}$  related to the reduced decomposition $g({\bf i})$ of $w_0$.\footnote{To check that $g({\bf i})$ is  a reduced decomposition of $w_0$ note that the length of the element of $W$ assigned to the reduced decomposition $g({\bf i})$ is maximal. Indeed, if it is smaller, the same will be true for the element related to the reduced decomposition $g^{-1}g({\bf i}) = {\bf i}$.   And $w_0$ is the unique element of maximal length in $W$.}
  Since any two reduced decompositions of $w_0$ are related by  braid relations (\ref{1}), and  any braid relation is a cluster transformation by Section \ref{pcs.sec}, it follows that the cluster coordinate systems 
  ${\cal C}^*_{g({\bf i}), t}$ and ${\cal C}^*_{{\bf i}, t}$ are related by a cluster transformations. The last step - the amalgamation with the Cartan group -  is respected by     
  these cluster transformations.  
  \end{proof}

  \paragraph{\bf Remarks.} 1. Theorem \ref{OCL} for the groups of type ${\rm A}_m$ has been proved in \cite[Section 9]{GS16}. 
  That proof was much more complicated,  taking a number   pages to accomplish, 
  since the   amalgamation construction of the cluster atlases was not available in full generality even for the moduli spaces 
  ${\mathscr X}_{{\rm PGL_m}, \bS}$,  ${\mathscr A}_{{\rm SL_m}, \bS}$, and so one had to stick to the special coordinate system introduced in \cite{FG03a}. 
  
  2. According to \cite[Section 9.3]{GS16}, the tropicalization of the canonical involution $\ast \in {\rm Out}(\G)$ acting on the space   ${\mathscr A}_{{\rm SL_m}, t}$ 
  coincides with the Schutzenberger involution. Therefore Theorem \ref{OCL} provides a cluster definition of the $\G-$analog of the Schutzenberger involution for any $\G$. 
  
  \medskip
  
   \subsection{Cluster nature of the Weyl group action on  spaces  ${\mathscr P}_{\G, \bS}$ and ${\mathscr A}_{\G, \bS}$} \la{Sec5.2}

\medskip

  \paragraph{\bf Cluster varieties from the  {cyclic envelope} of the braid semigroup 
${\mathbb B}^+$.} Denote by $[{\mathbb B}^+]$ the coinvariants of the cyclic shift $s_{j_1} \ldots s_{j_{n-1}}  s_{j_{n}}  \lms  s_{j_{n}} s_{j_1} \ldots s_{ j_{n-1}}$  on the braid semigroup  ${\mathbb B}^+$,  called   the {\it cyclic envelope} of  ${\mathbb B}^+$. 
A  cyclic word   
$ {\bf b}= s_{j_1} \ldots s_{j_n}$  
projects to an element   $[{\bf b}]\in [{\mathbb B}^+]$.

\bl \la{6.11}  Any element $[{\bf b}]\in [{\mathbb B}^+]$ gives rise to  cluster Poisson variety ${\mathscr X}_{\G, [{\bf b}]}$. 
\el

\begin{proof}    
We assign to a  word   ${\bf b} =   s_{j_1} \ldots s_{j_n}$  the amalgamation of  elementary cluster Poisson tori $ {\mathscr X}_{s_{i_1}}\ast \ldots \ast {\mathscr X}_{s_{i_n}}$, followed by the 
  amalgamation of the left frozen side    of   ${\mathscr X}_{s_{i_1} }$ with the right one of ${\mathscr X}_{s_{i_n} }$.   As in Figure \ref{X-space for b}, each elementary triangle gives rise to a moduli space $\mathscr{X}_{i_k}$, whose corresponding quiver is illustrated in Section \ref{sec.3.1}.
 The obtained cluster Poisson variety ${\mathscr X}_{\G, [{\bf b}]}$  is  invariant under the cyclic shift.    
Braid relation (\ref{1})  can be applied 
to any segment of a cyclic word in 
${\mathbb B}^+$. By 
\cite{FG05}, it gives rise to a cluster Poisson transformation.  So cluster Poisson variety ${\mathscr X}_{\G, [{\bf b}]}$ depend  only on  $[{\bf b}]\in [{\mathbb B}^+]$. \end{proof} 

\begin{figure}[ht]
\begin{tikzpicture}
\draw[thick] (0,0) circle (1.5cm); 
\node at (0,0) {$\bullet$};
\draw[thick] (0,0) -- (90:1.5);
\draw[thick] (0,0) -- (130:1.5);
\draw[thick] (0,0) -- (170:1.5);
\draw[thick] (0,0) -- (210:1.5);
\draw[thick] (0,0) -- (250:1.5);
\draw[thick] (0,0) -- (290:1.5);
\draw[thick] (0,0) -- (330:1.5);
\draw[thick] (0,0) -- (50:1.5);
\node at (15:1) {$\cdot$};
\node at (10:1) {$\cdot$};
\node at (5:1) {$\cdot$};
\node[blue] at (110: 1.7) {$s_{i_1}$};
\node[blue] at (150: 1.7) {$s_{i_2}$};
\node[blue] at (190: 1.7) {$s_{i_3}$};
\node[blue] at (230: 1.7) {$s_{i_4}$};
\node[blue] at (270: 1.7) {$s_{i_5}$};
\node[blue] at (310: 1.7) {$s_{i_6}$};
\node[blue] at (70: 1.7) {$s_{i_n}$};
\end{tikzpicture}
\caption{The moduli space ${\mathscr X}_{\G, [{\bf b}]}$.}
\label{X-space for b}
\end{figure}

 Generic points of ${\mathscr X}_{\G, [{\bf b}]}$ parametrise $\G$-local systems on a punctured disc 
with   special points $x_1, ..., x_n$  and the following    data: 

\vskip 2mm
1.  For each $k \in \Z/n$, a flat section $\beta_k$ of the flag local system ${\mathscr L}_{\cal B}$  
near the point $x_k$, so that 
$\beta_k$ and $\beta_{k+1}$ are in the relative position $s_{j_k}$.

3. A flat section $\beta$ of ${\mathscr L}_{\cal B}$ near the puncture. 

 \vskip 2mm
 Denote by ${\bf Q}$ the quiver obtained by the cyclic amalgamation of the elementary quivers.

 Recall the level map   (\ref{LVLM}). Let    ${\bf b}_i:= l_{{\bf Q}}^{-1}(i)$. It is a cyclic set, which  inherits 
 from the quiver ${\bf Q}$  a structure of a sub-quiver,  denoted by ${\bf q}_{i}$.  
 It is identified with an oriented  polygon.

\bl \la{6.11X}  
If    $|{\bf b}_i|>0$ for all $i\in {\rm I}$,    the   group $W$  acts   by birational automorphisms of ${\mathscr X}_{\G, [{\bf b}]}$.
\el 

\begin{proof}  Since   $|{\bf b}_i|>0$ for all $i\in {\rm I}$,  the monodromy of a  generic  $\G$-local system from ${\mathscr X}_{\G, [{\bf b}]}$ is regular.  
So  monodromy invariant    flags near the puncture    
  form a $W$-torsor. 
The   group $W$ acts by altering the invariant flag near the puncture, keeping the rest  intact.   
\end{proof}

There are canonical projections of cluster Poisson tori, and related cluster Poisson varieties:
\begin{equation}\label{eq78z}
\begin{split}
&\pi_i: {\mathscr X}_{{\rm G} , {\bf b}} \lra {\mathscr X}_{{\rm PGL}_2, {\bf b}_i}, \qquad \qquad\pi_i: {\mathscr X}_{{\rm G} , [{\bf b}]} \lra {\mathscr X}_{{\rm PGL}_2, [{\bf b}_i]}.\\
\end{split}
\end{equation}

\subsubsection{Cluster transformations $S_{ i}$.}
 A cyclic word  ${\bf b}$  is called {\it admissible} if  $|{\bf b}_i|>1$ for each $i \in {\rm I}$.  
 Let us order  vertices of the quiver  ${\bf q}_i$,  getting  a sequence of vertices $v_{1},   \ldots v_{m}$ of the 
  quiver ${\bf Q}$.  
 
\bd \la{DEF10.4}
Given an admissible cyclic word ${\bf b}$,     the cluster transformation   $S_{  i}$ is    the following composition of mutations and  a symmetry $\pi_{m-1, m} $  
exchanging $v_m \leftrightarrow v_{m-1}$:
\be \la{SEQM}
S_{  i}:=  \mu_{v_{1}}\circ   \ldots  \circ   \mu_{v_{m-1}} \circ  \pi_{m-1, m} \circ \mu_{v_{m-1}}  \circ    \ldots \circ \mu_{v_{1}}.
 \ee
 \ed  
 
 This definition  makes sense only if the cyclic word ${\bf b}$ is  admissible. 
 
 Denote by ${\rm D^*_{k}}$ a punctured disc   
  with $k$ special     points. 
 
\bt  \la{THH6.4} i) Let ${\bf b}$ be an admissible cyclic word. Then the assignment  
\be\la{ASSI}
s_i \lms S_i, \qquad\forall i\in {\rm I},
\ee
 gives rise to an action of the 
Weyl group $W$    by   cluster Poisson transformations of ${\mathscr X}_{\G, [{\bf b}]}$.

This action 
  coincides with the geometric action $\tau$ of $W$ on ${\mathscr X}_{\G, [{\bf b}]}$ from Lemma \ref{6.11X}.

\vskip 1mmii) The map (\ref{ASSI}) provides an action of the Weyl group by cluster Poisson transformations of the space ${\mathscr P}_{\G, {\rm D}^*_{k}}$, 
    $k>1$.
 
 \vskip 1mmiii)  Let $\bS$ be a decorated surface  with a puncture $p$.  Assume that  $\bS\not = {\rm D}^*_{1}$ if $\G$ of type ${\rm A_1, A_2}$,  and if $\bS$ has no boundary, the number of punctures is greater than $1$.     

Then the Weyl group  acts on the moduli space ${\mathscr P}_{\G, \bS}$  
   by cluster Poisson transformations. 
\et

Theorem \ref{THH6.4} for   moduli spaces ${\mathscr X}_{\rm PGL_m}$ has been proved in \cite[Sections 7-8]{GS16}. 

The proof of Theorem \ref{THH6.4} goes until Section \ref{13.3.4}.
\begin{proof} 

i) Let us show that the cluster transformation $S_i$ coincides with the transformation $\tau_{s_i}$  defining the geometric action of the generator $s_i$ on ${\mathscr X}_{\G, [{\bf b}]}$.
This would imply  that  elements $S_i$ satisfy the Weyl group relations.

\bl \la{L10.16} The cluster transformation $S_i$ preserves the quiver  ${\bf Q}$.  
\el

\begin{proof} The pair of quivers ${\bf q}_{i}\subset {\bf Q}$  satisfies  condition (176) of \cite[Definition 7.4]{GS16}:  
we have 
\be \la{176}
\sum_{i \in {\rm V}_{{\bf q}_{i}}}  \varepsilon_{ij}=0, \qquad \forall j \in {\rm V}_{\bf Q} - {\rm V}_{{\bf q}_{i}}.
\ee
Indeed, it is clear that the cyclic closure of the quiver ${\bf Q}_{\alpha, \beta}$ from Section \ref{SEC12.1} satisfies this condition. Therefore \cite[Theorem 7.7]{GS16} implies the claim. \end{proof}

\subsubsection{Cluster Poisson coordinates on the   torus ${\mathscr X}_{{\rm G} , {\bf b}}$.} Denote by $X_v$ the cluster Poisson coordinate on ${\mathscr X}_{{\rm PGL}_2, {\bf b}_i}$ assigned to the vertex $v$ of the polygon ${\bf q}_i$. 
We denote by $\widetilde v$ a vertex  of the quiver ${\bf Q}$ provided by the vertex $v$ of the quiver ${\bf q}_i$, and by 
$X_{\widetilde v}$ the corresponding cluster Poisson coordinate on ${\mathscr X}_{\G, {\bf b}}$. Let us recall the construction of these coordinates.   

For each special point $x_k$, consider an edge $E_k$ connecting it with the puncture,  as illustrated by the radiuses in Figure \ref{X-space for b}. Flat sections $\beta, \beta_k$, restricted to   $E_k$, define a generic pair of flags  
$(\B, \B_k)$. Then there is  a unique flag $\B^{(i)}_{k}$ in the   relative position $s_{j_k}$ to  $\B$ and $s_{j_k}w_0$ to $\B_{k}$.  
We   transport  $\B^{(i)}_{k}$ along the edge to the fiber at the point $x_k$. We get a collection of flags   
at the  special points, parametrized by    the cyclic set ${\bf b}_i$: 
\be \la{FLL2}
B^{(i)}_1, ... , B^{(i)}_{|{\bf b}_i|}.
\ee 
The coordinate $X_{\widetilde v}$ assigned to the element $v \in {\bf b}_i$ is given by 
$
X_{\widetilde v}:= r^+(\B, \B^{(i)}_{v-1}, \B^{(i)}_v, \B^{(i)}_{v+1}).
$
This definition make sense if $|{\bf b}_i|>1$: by $\B^{(i)}_{v-1}$ (respectively $\B^{(i)}_{v+1}$) we mean the flat section at the slot $v-1$ (respectively   $v+1$), 
parallel transported to $v$.

\subsubsection{A geometric version of projection (\ref{eq78z}).} 
A Borel subgroup $\B\subset \G$ and an element $i \in {\rm I}$ provide a parabolic subgroup ${\rm P}_i$ containing $\B$. Let $\U_i$ be the unipotent radical of $\rm P_i$. 
The quotient ${\rm P}_i/\U_i$ is a reductive group, which has a natural projection onto  the group  ${\rm PGL_2}$. 

So a $\G$-local system ${\cal L}$ on a punctured disc and an  invariant flag $\B$ near the puncture   
determine a  ${\rm PGL_2}$-local system ${\cal L}^{(i)}\subset {\cal L}$ on the punctured disc. A cyclic sequence of flags  $\B_1, ..., \B_n$    
in the fibers   ${\cal L}_{x_1}, ..., {\cal L}_{x_n}$  gives rise to a similar cyclic sequence of flags (\ref{FLL2}) in the ${\rm PGL_2}$-local system ${\cal L}^{(i)}$. 
So we get  a point of the moduli space ${\mathscr X}_{{\rm PGL_2} , [{\bf b}_i]}$. Summarising, we get  a canonical projection  
 \begin{equation}\label{eq78aa}
\begin{split}
&\pi_i: {\mathscr X}_{{\rm G} , [{\bf b}]} \lra {\mathscr X}_{{\rm PGL}_2, [{\bf b}_i]}.\\
\end{split}
\ee
By  the construction,  the   coordinate   $ X_{\widetilde v}$, $v \in {\bf b}_i$,  on  ${\mathscr X}_{{\rm G} , {\bf b}}$ is  the   pull back    of the  one  
   $X_v$    on    ${\mathscr X}_{{\rm G} , {\bf b}_i}$:
\be \la{PPPp}
 X_{\widetilde v}  =  \pi_i^*(X_v), \qquad \qquad\mbox{where $i \in {\rm I}$, and $v$ is a vertex of the quiver ${\bf q}_i$.}
\ee

Denote by $S_{{\bf q}_i}$ the cluster transformation (\ref{SEQM}) for the quiver ${\bf q}_i$. Evidently,   $S_{i}^*X_{\widetilde v} = \pi_i^*(S_{{\bf q}_i}^*X_v)$. We proved in \cite[Theorem 8.2]{GS16}  that    $S_{{\bf q}_i} = \tau_{s_i}$. 
 Thus there is a commutative diagram
 \begin{equation}\label{eq78ab}
\begin{gathered}
\xymatrix{
{\mathscr X}_{{\rm G} , [{\bf b}]} \ar[r]^{s_{i}} \ar[d]_{\pi_i}& {\mathscr X}_{{\rm G} , [{\bf b}]} \ar[d]^{\pi_i } &\\
{\mathscr X}_{{\rm PGL}_2, [{\bf b}_i]} \ar[r]_{s_i} & {\mathscr X}_{{\rm PGL}_2, [{\bf b}_i ]} &}
\end{gathered}
\end{equation}
 
 Let $i, j \in {\rm I}$. Set ${\bf b}_{ij}:= l_{\bf Q}^{-1}(\{i,j\})$. We consider several cases. 
 
 1) Assume   that $C_{ij}=C_{ji}=1$. Then there is a projection 
 $$
 \pi_{ij}: {\mathscr X}_{{\rm G} , [{\bf b}]}  \lra {\mathscr X}_{{\rm PGL_3}, [{\bf b}_{ij}]}.
 $$ 
Thanks to \cite[Theorem 8.2]{GS16}, the cluster Poisson transformation $S_i$ on ${\mathscr X}_{{\rm PGL_3}, [{\bf b}_{ij}]}$ coincides 
 with the geometric action $\tau_{s_i}$ of the Weyl group generator $s_i$ on ${\mathscr X}_{{\rm PGL_3}, [{\bf b}_{ij}]}$.  
 By the construction,  $\pi^*_{ij} \circ \tau_{s_i}$ coincides with the geometric action of $s_i$ on  ${\mathscr X}_{{\rm G} , [{\bf b}]}$.  
 By the definition, the map $\pi^*_{ij}$ intertwines the cluster transformations $S_i$ on ${\mathscr X}_{{\rm PGL_3}, [{\bf b}_{ij}]}$ and${\mathscr X}_{\G, [{\bf b}]}$. 
 The claim follows from these remarks. 
 
 2) If $C_{ij}=C_{ji}=0$, the claim follows by similar arguments using ${\rm PGL_2} \times  {\rm PGL_2} $.
 
3)  In the non-simply laced case we use the  cluster unfolding, see  \cite[Section 3.6]{FG05}, or Section 6 of the version 3 of \cite{FG03b},   to reduce the claim to the simply-laced case.

Let 
 $\sigma: {\rm A}_3 = \{1 - 2 - 3\}\to {\rm B}_2$  be a   folding of  Dynkin diagrams. Set  $a:= \sigma (1) = \sigma (3)$ and $b := \sigma(2)$. 
 One can choose a quiver ${\bf Q}_{A_3}$ for the  group of type $A_3$  of a shape ${\bf Q}^{1}_{A_3} - {\bf Q}^{2}_{A_3} -{\bf Q}^{3}_{A_3}$, with no arrows between   
  ${\bf Q}^{1}_{A_3}$ and   ${\bf Q}^{3}_{A_3}$, a quiver ${\bf Q}_{B_2}$ for the  group of type $B_2$  of a shape ${\bf Q}^{a}_{B_2} - {\bf Q}^{b}_{B_2}$, and a folding map 
  ${\bf Q}_{A_3} \to{\bf Q}_{B_2}$ which induces isomorphisms  
  ${\bf Q}^{1}_{A_3}  \stackrel{\sim}{=} {\bf Q}^{a}_{B_2}$ and ${\bf Q}^{3}_{A_3}   \stackrel{\sim}{=} {\bf Q}^{b}_{B_2}$.
    
    Then the cluster transformation $S_a$ 
 is a product of two  cluster transformations  
 $S_1 = \mu_1 \circ ... \circ \mu_k$ and $S_3 = \nu_1 \circ ... \circ \nu_k$, such that the mutations $\mu_i$ and $\nu_j$ commute. Therefore we have  
  $S_a =( \mu_1  \nu_1) \circ  (\mu_2  \nu_2)  \circ ... \circ (\mu_k  \nu_k)$, where $\mu_i\nu_i$ act as a mutation in the folded quiver.

ii)  The  moduli space ${\mathscr P}_{\G, {\rm D^*_{k}}}$ has a structure of cluster 
Poisson  variety ${\mathscr P}_{\G, [\mu(w_0)^{k}]}$. The cyclic word $\mu(w_0)^{k}$ is admissible if $k>1$. So the claim for unfrozen cluster Poisson variables follows from i). 
  The analog of Lemma \ref{L10.16} for frozen vertices holds since, by construction, each frozen vertex satisfies condition (\ref{176}). 
The claim for the frozen variables follows by the reduction to the ${\rm A}_1$ case from a similar result  for $\G = {\rm PGL_2}$ with frozen variables  \cite[Section 7]{GS16}.

 iii) Take an ideal triangulation of $\bS$ with at least two edges at the puncture $p$. Then $\bS$ can be glued from a punctured disc ${\rm D^*_{k}}$, with $k>1$ special  points, 
and a decorated surface $\bS'$. So the claim follows from ii). 
\end{proof} 

\subsubsection{The Weyl group action on the space  ${\mathscr A}_{\G, \bS}$.}  \la{13.3.4} As  in  \cite[Section 6]{GS16}, every puncture $p$ of $\bS$ gives rise to a geometric Weyl group action on the space  ${\mathscr A}_{\G, \bS}$.  The key ingredients for  this action are the partial potential functions $\mathcal{W}_{p, i}$, $i\in {\rm I}$,   
see  Section \ref{triple2020.1.28}. 
Lemma \ref{L7.12} explicitly presents the partial potential of an elementary configuration in terms of cluster coordinates. 
By Lemma \ref{6.30.17.45.hh}, $\mathcal{W}_{p, i}$ is the summation of the $i$th partial potentials of elementary configurations centered at $p$.
In this way, we obtain an explicit expression of $\mathcal{W}_{p, i}$  in terms of cluster coordinates.

The pair of spaces $({\mathscr A}_{\G, \bS}, {\mathscr P}_{\G, \bS})$ forms a cluster ensemble. For every puncture $p$, the  cluster  transformation from Definition \ref{DEF10.4} describing 
 the action of the generator $s_i \in W$ on the space   ${\mathscr P}_{\G, \bS}$ defines  a cluster transformation of the space ${\mathscr A}_{\G, \bS}$. 
 We claim that the latter coincides with the aforementioned geometric Weyl action on  ${\mathscr A}_{\G, \bS}$.
 Since the cluster transformation affects only    the level $i$ coordinates, the claim reduces to ${\rm SL_2}$, which has been done in \cite{GS16}. Namely,  we use 
 \cite[Theorem 7.7]{GS16} and  repeat the  three lines at the  end of  \cite[Section 8.1]{GS16}.
 
\medskip

\subsection{The braid group action  on the  moduli space  ${\mathscr X}_{\G,  \bS}$} \la{SEC5.1}
 
\medskip

In Sections \ref{SEC5.1}-\ref{sect10.4}, given a boundary component $\pi$ of $\bS$, we introduce actions of the group ${\mathbb B}_{\mathfrak g}^{(\pi)}$, defined in (\ref{BR*}), on the moduli spaces 
${\mathscr X}_{\G, \bS}$ , ${\mathscr P}_{\G, \bS}$ and ${\mathscr A}_{\G, \bS}$, and prove that it is given by cluster Poisson transformations for the space ${\mathscr X}_{\G, \bS}$, and quasi-cluster Poisson 
transformations for ${\mathscr P}_{\G, \bS}$. 
On  pictures we present the   case when   the number $d_\pi$  of special points on $\pi$  is even, i.e. ${\mathbb B}_{\mathfrak g}^{(\pi)} = {\mathbb B}_{\mathfrak g}$.

\subsubsection{Standard collections of flags.} Let   
 $(\B^+, \B^-)$ be a generic pair of Borel subgroups   in $\G$.
The set of Borel subgroups containing the Cartan subgroup  $\H:= \B^+\cap  \B^-$   is said to be
the {\it standard collection} of Borel subgroups 
for  $(\B^+, \B^-)$. 
The Weyl group $W=N({\rm H})/{\rm H}$ acts simply transitively 
on the standard collection. 
So one can label these Borel subgroups by   elements of  $W$,
 so that 
$\B^+$ is labelled by   $e$, and $\B^-$   by $w_0$. 
The group $\G$ acts by conjugation 
on the set of standard collections, with the stabiliser    $\H$.

\bl \la{12.10.08.1} 
Any collection of Borel subgroups 
$\{\B_s\}$, where $s \in W$, such that 
$\B_{s}$ and $\B_{t}$  are in the incidence relation $s^{-1}t$, 
is the standard collection for the pair $(\B_e, \B_{w_0})$.  
\el

\subsubsection{Braid group action on the space ${\mathscr X}_{\G, \bS}$.} 
Let    $S^1$ be an oriented circle which  
  carries   $d$   special points $x_i$, 
ordered cyclically by the circle orientation.  
  A {\it  framing} on a $\G$-local system ${\cal L}$ on $S^1$ is   
a collection $\{\beta_1, \ldots , \beta_{d}\}$ of 
flat sections  $\beta_i$      
of the flag bundle ${\cal L}_{\cal B}$ near $x_i$. 
  The moduli space ${\mathscr X}_{\G, S^1; d}$ parametrizes  pairs $({\cal L}, \beta)$, where ${\cal L}$ is a 
  $\G$-local system on $S^1$ with a framing $\beta$.

Choose   a special point $x_1$. 
Let us define an action $T_{x_1}$ of the  
 braid group ${\mathbb B}_\G$  by birational automorphisms 
of the space ${\mathscr X}_{\G, S^1; d}$. 
Take a reduced word ${\bf i} = (i_1, ..., i_N)$ for $w_0$:
\be \la{AST1} 
  w_0 = s_{i_1} s_{i_2} \ldots s_{i_N}.
 \ee

Consider  $(\beta_1, \beta_2)$ as a pair of flat sections on the arc $(x_1, x_2)$. 

If $d=1$, then $x_2=x_1$, while $\beta_2$ is by definition the flat section $\beta_1$ transported around the circle. Equivalently, we cut the circle $S^1$ at the point $x_1$, getting two flat sections at the ends. 

\begin{figure}[ht]
\epsfxsize 200pt
\center{
\begin{tikzpicture}[scale=0.8]
\draw [red, dashed]  (0,0) circle (20mm);
\node [red]  at (0,2) {\Large $\bullet$};
\node [red]  at (0,-2) {\Large $\bullet$};
\foreach \angle/\count in {330/4,150/2,210/3,30/1}
{
\node[thick] at (\angle:2cm) {$\circ$};
}
\foreach \angle/\count in {60/1,180/2,300/3}
{
\node at (\angle:1.8cm) {\small $s_2$};
\node at (\angle+60:1.8cm) {\small $s_1$};
}
\foreach \count in {6,5,3,2}
{
\node at (30+60*\count:2.6cm) {\small $\B_\count$};
}
\node at (90:2.6cm) {\small $\beta_1=\B_1$};
\node at (270:2.6cm) {\small $\beta_2=\B_4$};
    \end{tikzpicture}
 }
\caption{ A cyclic collection of 6 flags for   ${\rm PGL_3}$ for a framing $(\beta_1, \beta_2)$ and   $w_0=s_1s_2s_1$.} 
\label{pin1}
\end{figure}

The pair $(\beta_1, \beta_2)$ and  reduced word (\ref{AST1})   determine a standard collection 
of flags 
$$
 \B_1, \B_2, \ldots , \B_{N+1}, \qquad \beta_1 = \B_1, \beta_2=\B_{N+1}.
$$
We think about them as of   flat sections of ${\cal L}_{\cal B}$ located at $N+1$ points on the interval $[x_1, x_2]$, so that the flag $\B_1$ sits at $x_1$, the flag $\B_{n+1}$ 
sits at $x_2$ and other flags   between them. The points subdivide the interval $[x_1, x_2]$ into $N$ arcs. Each arc carries a generator $s_i$ of the  Weyl group. It tells the relative position of the two flags at the  arc ends.
If $d=1$, we are done. 
If $d>1$, we 
consider the {\it companion reduced word ${\bf i}^* = (i_1^*, ..., i_N^*)$} for $w_0$:
\be \la{AST}
  w_0 = s_{i^*_1} s_{i^*_2} \ldots s_{i^*_N}.
\ee
Just as above, the pair  of flat sections $(\beta_2, \beta_3)$  on the interval $[x_2, x_3]$,   together with the   reduced decomposition (\ref{AST}) of $w_0$, determines a standard collection 
of flags 
$$
\B_{N+1}, \B_{N+2}, \ldots , \B_{2N+1}, \qquad \beta_2 = \B_{N+1}, \beta_3 = \B_{2N+1}.
$$
We think about them as of the flat sections of ${\cal L}_{\cal B}$ near $N+1$ points on the interval $[x_2, x_3]$. 
Repeating this procedure $d$ times, we get $d N$ flat sections of ${\cal L}_{\cal B}$ cyclically located  on the  circle, see Figure \ref{pin1}. 
  Next, we need the following well known  Lemma.

 \bl \la{CENTRALZ} For any  reduced word $s_{i_1} ... s_{i_m}$,  there exists a reduced word of $w_0$  starting from $s_{i_1} ... s_{i_m}$. 
\el

\begin{proof}   Let $P$ be the set of positive roots and $-P$   the set of negative roots. For any $w\in W$, its length
$
l(w)= \# P \cap w(-P).
$ 
In particular, $w_0(-P)=P$ implies $l(w_0)= \# P.$ Meanwhile
\[
l(w^{-1}w_0)= \# P \cap w^{-1}w_0(-P)= \# P \cap w^{-1}(P)= \# w(P)\cap P.
\]
Therefore
$
l(w)+l(w^{-1}w_0)= \# P = l(w_0).
$
Take a reduced word of $w^{-1}w_0$. Then  the product of 
the  reduced decompositions of $w$ and $w^{-1}w_0$ is a reduced word of $w_0$.
\end{proof}

\subsubsection{The case when   $d$ is even.} Then, due to the   choice of a companion reduced word for $ w_0$, any pair of flags separated by   $N$ arcs 
is generic.  
Pick a generator $s_i \in  {\mathbb B}_{\mathfrak g}$, and a reduced word ${\bf i}$ of   $w_0$  starting from  $s_i$, see  Lemma \ref{CENTRALZ}. 
Let us define  an automorphism $T_{x_1}({s_i})$ of the moduli space 
${\mathscr X}_{\G,  S^1}$. It acts    by altering the decorations, leaving the local system intact:
\be
\begin{split}
&T_{x_1}({s_i}): ({\cal L},  \beta) \lra ({\cal L},  \beta').\\
\end{split}
\ee
The new  framings $\{ \beta'_1, \ldots ,   \beta'_{d}\}$  are obtained from the  ones $\{  \beta_1, \ldots ,   \beta_{d}\}$ by  moving  special points one step along the circle orientation,    see Figure \ref{pin2}. Special points are shown by  red   points.

\begin{figure}[ht]
\epsfxsize 200pt
\center{
\begin{tikzpicture}[scale=0.8]
\draw [red, dashed]  (0,0) circle (20mm);
\node [red]  at (0,2) {\Large $\bullet$};
\node [red]  at (0,-2) {\Large $\bullet$};
\foreach \angle/\count in {30/1,150/2,210/3,330/4}
{
\node[thick] at (\angle:2cm) {$\circ$};
}
\foreach \angle/\count in {60/1,180/2,300/3}
{
\node at (\angle:1.8cm) {\small $s_2$};
\node at (\angle+60:1.8cm) {\small $s_1$};
}
\foreach \count in {1,2,3,4,5,6}
{
\node at (30+60*\count:2.6cm) {\small $\B_\count$};
}
\draw [-latex, ultra thick] (3.6,0) -- (5.4,0);
\node at (4.5, 0.7) {${\rm T}(s_1)$};
\begin{scope}[shift={(9,0)}]
\draw [red, dashed]  (0,0) circle (20mm);
\node [red]  at (150:2cm) {\Large $\bullet$};
\node [red]  at (330:2cm) {\Large $\bullet$};
\foreach \angle/\count in {30/1,90/2,210/3,270/4}
{
\node[thick] at (\angle:2cm) {$\circ$};
}
\foreach \angle/\count in {60/1,180/2,300/3}
{
\node at (\angle:1.8cm) {\small $s_2$};
\node at (\angle+60:1.8cm) {\small $s_1$};
}
\foreach \count in {1,2,3,4,5,6}
{
\node at (30+60*\count:2.6cm) {\small $\B_\count$};
}
\end{scope}
    \end{tikzpicture}
 }
\caption{A braid group generator $s_1 \in {\mathbb B}_{\rm PGL_3}$ acts by altering the framings. The 
new framing  $(\B_2, \B_5)$ is obtained 
by moving the special   points    one step counterclockwise.}
\label{pin2}
\end{figure}
 More generally, let $w=s_{i_1} ... s_{i_m}$ be a reduced word  
for an element $w\in W$. We define  a transformation
$T_{x_1}({\mu(w)})$   by taking a reduced word for  $w_0$ which starts from $w=s_{i_1} ... s_{i_m}$, and moving   special points $m$ steps 
following  the circle orientation. By  the very construction 
we have 
\be   \la{12.10.08.35}
T_{x_1}({\mu(w)}) = T_{x_1}(s_{i_1}) ...T_{x_1}(s_{i_m}).
\ee
 
\bl \la{PROP5.5}
The transformations $T_{x_1}({s_i})$ satisfy the braid group relations. 
\el

 \begin{proof}  Given a braid group relation (\ref{1}),   by Lemma \ref{CENTRALZ} there exists a reduced word of $w_0$ which starts from either of the two  sides of (\ref{1}):
  $
 w_0 = s_1s_2s_1 \ldots w' = s_2s_1s_2 \ldots w'. 
  $ 
 By the very definition,  $T_{x_1}({s_1s_2s_1 \ldots})=T_{x_1}({s_2s_1s_2 \ldots})$.  
 So the claim    
 follows   from    (\ref{12.10.08.35}). \end{proof}

\subsubsection{The case when   $d$ is odd.} Let us define an action of the group $\B^*_{\mathfrak g}$ on 
${\mathscr X}_{\G, S^1; d}$ in this case.  
We have constructed above an action of $\B^*_{\mathfrak g}$ on the space 
 ${\mathscr X}_{\G, S^1; d}$ for even $d$. If $d$ is odd, ${\mathscr X}_{\G, S^1; d}$ can be recovered as the space of 
 objects in  ${\mathscr X}_{\G, S^1; 2d}$ equivariant under the action of the element $\mu(w_0)$:
 \be \la{eq}
 {\mathscr X}_{\G, S^1; d}= {\mathscr X}_{\G, S^1; 2d}^{\mu(w_0)}.
  \ee  
  Then the subgroup $\B^*_{\mathfrak g}$ acts on it, since it consists of  the elements  commuting with  $\mu(w_0)$. This   approach  makes clear where does the group 
  $\B^*_{\mathfrak g}$ came from.
      
\subsubsection{Another approach  for the subgroup $\widetilde \B^\ast_{\mathfrak g}$.} Recall the involution $w \lms w^*:= \mu(w_0)^{-1} w \mu(w_0)$.  
For any element $w\in W$ such that $w^*=w$,   
the transformation
$T_{x_1}({\mu(w)})$ is defined just as before. Note that due to the  condition $w = w^*$, and the   choice of a companion reduced word for $ w_0$, any pair of flags separated by   $N$ arcs 
is generic. Indeed,  the new collection of flags corresponds to a reduced decomposition of  $\mu(w)^{-1}\mu(w_0)^d\mu(w) = \mu(w_0)^d$. 
So we get an action of the group $\widetilde \B^\ast_{\mathfrak g}$.\footnote{There are  the  following elements  in   $\widetilde \B^{\ast}_{\mathfrak g}$, which probably generate it:
\be \la{1**}
s_is_{i^*} ~\mbox{if}~  C_{i, i^*}=0, \qquad s_is_{i^*}s_i ~\mbox{if}~  C_{i, i^*}= C_{i*, i}=-1, \qquad(s_is_{i^*})^2~\mbox{if} ~ C_{i, i^*} =-2, \qquad(s_is_{i^*})^3~\mbox{if} ~ C_{i, i^*} =-3.
\ee}

\vskip 2mm
 Let $\bS$ be a decorated surface.        
 Each boundary component   is identified with an oriented circle $S^1$. So all constructions above can be applied.

\bt \la{TTHH1}
  Choose a   special point 
$x_1$ at a boundary component $\pi$ of  $\bS$. Then 

\vskip 1mma) There is   birational  action $T_{x_1}$ of the  
 braid group ${\mathbb B}^{(\pi)}_{\mathfrak g}$  
on   ${\mathscr X}_{\G, \bS}$.  
One has 
 $T_{x_2}(w) = T_{x_1}(w^*)$.

\vskip 1mmb)   The  element $\mu(w_0)\in {\mathbb B}^*_{\mathfrak g}$ acts   as the element of   $\Gamma_\bS$ 
 shifting  by one   special points on  $\pi$.

\vskip 1mmc)  The braid subgroup $\widetilde {\mathbb B}^{(\pi)}_{\mathfrak g}$ acts by   automorphisms of the cluster Poisson structure on ${\mathscr X}_{\G, \bS}$.

 The group $\Gamma'_{\G, \bS}$   maps to the cluster modular group of   ${\mathscr X}_{\G, \bS}$. 
\et

\begin{proof} a) The first claim follows from Lemma \ref{PROP5.5} and its analog for odd $d$,  proven by an equivariance argument. 
 The second   is clear from the definition of the companion reduced word (\ref{AST}). 

 \vskip 2mm
 
 b) Moving the special point $x_1$ by one step   $N$ times amounts to moving it to the next point $x_2$.   
 
\vskip 2mm

c)  The braid group action related to 
a boundary component is  determined by the restriction of the 
framed local system to the boundary.  Therefore the actions of the braid groups at different boundary components commute, 
and   compatible with the action of the   group   $\Gamma_\bS$. 

It remains to show that the group $\widetilde {\mathbb B}^\pi_{\mathfrak g}$ maps to the cluster modular group of ${\mathscr X}_{\G, \bS}$.  

Assume   that $d$ is even. Write   the reduced decomposition ${\bf i}$  as  $w_0 = s_iw'$. Let ${\bf i}_+$ be the reduced decomposition   $w_0 = w's_i^*$. The original collection of intermediate flags 
is described by  a sequence of reduced decompositions ${\bf i}, {\bf i}^*,   \ldots, {\bf i}, {\bf i}^*$.  After the cyclic shift by one step, the resulting sequence of 
 intermediate flags is described by  a similar sequence of reduced decompositions ${\bf i_+}, {\bf i_+^*},   \ldots, {\bf i_+}, {\bf i^*_+}$.

\bl \la{III} There exists a sequence of braid relations transforming ${\bf i}{\bf i}^*$ to ${\bf i}_+{\bf i}_+^*$.  
\el

\begin{proof}  
  Write 
$
  \mu(w_0)^2  = s_iz'$. 
 The words $s_iz'$ and $z's_i$ define the same element $\mu(w_0)^2$ of the 
braid semigroup ${\mathbb B}^+$. Indeed, since $ \mu(w_0)^2$ is in the center,  $z's_i = s_i^{-1}(s_iz')s_i = s_iz'$. 
Thus one can find a sequence of 
elementary transformations (\ref{1})
 transforming  the first word to the second. 
Their composition 
transforms the seed defined by $s_iz'$ to the one defined by 
$z's_i $. So it provides a cluster transformation  
$ 
{\mathscr X}_{\mathbf i(Z's_i)} \lra {\mathscr X}_{\mathbf i(s_iz')}. 
$ 
 \end{proof}

Using a sequence of braid relations, we  
 transform   ${\bf i}_+ {\bf i}^*_+\lra {\bf i}{\bf i}$, and perform the corresponding sequence of mutations. The resulting
  quiver, describing the shifted collection of flags, is isomorphic to the original one. 
  
  Note that a braid realations in ${\bf i}$ give rise to 
 cluster Poisson transformations of the space 
${\mathscr X}_{\G,\bS}$. Indeed,  
 take  an ideal triangulation ${\cal T}$ of $\bS$ such that one side of a  triangle $t$ of ${\cal T}$ is given by an interval on the   boundary circle   of $\bS$ 
containing the segment of ${\bf i}$. So the part c) for even $d$ follows. 
 
The case of odd $d$ is  similar.  All  we need is an  analog of Lemma \ref{III} based on  the fact that if $w=w^*$, where $w\in W$, then writing 
a reduced decomposition $w_0 = ww'$ we have $w_0 = w'w$. 
\end{proof}

\medskip
\subsection{The braid group actions on the spaces ${\mathscr P}_{\G, \bS}$ and ${\mathscr A}_{\G, \bS}$}   \la{SECT9.4}

\medskip

In Section \ref{SECT9.4} we handle the cases of even and odd $d$  as before: we first run the construction for even $d$, and then reduce the case of an odd $d$ to it by either 
an equivariance argument, or directly. 

\subsubsection{Braid group action on the space ${\mathscr P}_{\G, \bS}$.}    As in Figure \ref{gs1}, a pinning $p_m: \beta_m \to \beta_{m+1}$ can be defined by a pair $(\alpha^r_{m+1}, \alpha^l_{m})$, where $\alpha^r_{m+1}$ is a  flat section 
of the decorated flag bundle ${\cal L}_{\mathscr A}$ over the space of positively oriented  tangent vectors to   $S^1$  
 near   $x_{m+1}$, projecting to $\beta_{m+1}$,   pictured on the right of $x_m$, and   $\alpha^l_{m}$ is a flat section near $x_{m}$ projecting to $\beta_{m}$, 
  pictured on the left of $x_{m}$, such that
 $$
 h(\alpha^r_{m+1}, \alpha^l_m) =1. 
 $$
 
  The moduli space ${\mathscr P}_{\G, S^1}$ parametrizes   triples $({\cal L}, \beta, p)$, where ${\cal L}$ is a 
  $\G$-local system on $S^1$ with a framing $\beta$, and $p$ is a collection of pinnings $p_m: \beta_m\to   \beta_{m+1}$, $m\in \Z/d\Z$,  
between the    framings at the  special points $(x_m, x_{m+1})$. 
 So pinnings provide a collection of   decorated flat sections 
\be \la{ORDEC}
\{\alpha^l_1, \alpha^r_1, \ldots ,  \alpha^l_{2d}, \alpha^r_{2d}\}.
\ee

\begin{figure}[ht]
\epsfxsize 200pt
\center{
\begin{tikzpicture}[scale=0.7]
\draw [red, dashed]  (0,0) circle (20mm);
\node [red, label=above:$\beta_1$] (a) at (0,2) {\Large $\bullet$};
\node [red, label=below:$\beta_2$] (b) at (0,-2) {\Large $\bullet$};
\node [thick] (c) at (-1.72,1) {$\circ$};
\node [thick] (d) at (1.72,1) {$\circ$};
\node [thick] (e) at (1.72,-1) {$\circ$};
\node [thick] (f) at (-1.72,-1) {$\circ$};
\draw[latex-] ([shift=(80:26mm)]0,0) arc (80:-80:26mm);
\draw[latex-] ([shift=(260:26mm)]0,0) arc (260:100:26mm);
\node (g) at (3,0) {\small $1$};
\node (h) at (-3,0) {\small $1$};   
\node (i) at (0.9,1.54) {\small $s_2$};
\node (j) at (1.8,0) {\small $s_1$};  
\node (k) at (0.9,-1.54) {\small $s_2$};
\node (l) at (-0.9,-1.54) {\small $s_1$};  
\node (m) at (-1.8,0) {\small $s_2$};
\node (n) at (-0.9,1.54) {\small $s_1$};  
\draw [latex-latex, ultra thick]  (5.4,0)-- (3.6,0);
\draw [red, dashed]  (9,0) circle (20mm);
\node [red] (a1) at (9,2) {\Large $\bullet$};
\node [red] (b1) at (9,-2) {\Large $\bullet$};
\node [thick, label=left:\textcolor{red}{\small $\A_2$}] (c1) at (7.28,1) {$\circ$};
\node [thick, label=right:\textcolor{red}{\small $\A_6$}] (d1) at (10.72,1) {$\circ$};
\node [thick, label=right:\textcolor{red}{\small $\A_5$}] (e1) at (10.72,-1) {$\circ$};
\node [thick, label=left:\textcolor{red}{\small $\A_3$}] (f1) at (7.28,-1) {$\circ$};
\node [red] (g1) at (12,0) {\small $1$};
\node [red] (h1) at (6,0) {\small $1$};   
\node (i1) at (9.9,1.54) {\small $s_2$};
\node (j1) at (10.8,0) {\small $s_1$};  
\node (k1) at (9.9,-1.54) {\small $s_2$};
\node (l1) at (8.1,-1.54) {\small $s_1$};  
\node (m1) at (7.2,0) {\small $s_2$};
\node (n1) at (8.1,1.54) {\small $s_1$};  
\node [red] at (6.9,2.2) {\small $1$};
\node [red]  at (6.9,-2.2) {\small $1$};   
\node [red] at (11.1,-2.2) {\small $1$};
\node [red]  at (11.1,2.2) {\small $1$};
\node at (9,3) {\small $\B_1$};  
\node at (9,-3) {\small $\B_4$};  
\node[red] (o) at (9.7, 2.54) {\small $\A_1^r$};  
\node[red] (p) at (9.7,-2.54) {\small $\A_4^l$};
\node[red] (q) at (8.3,-2.54) {\small $\A_4^r$};  
\node[red] (r) at (8.3,2.54) {\small $\A_1^l$};  
\draw[red] ([shift=(35:25mm)]9,0) arc (35:65:25mm);
\draw[red] ([shift=(-15:25mm)]9,0) arc (-15:15:25mm);
\draw[red] ([shift=(-65:25mm)]9,0) arc (-65:-35:25mm);
\draw[red] ([shift=(35:-25mm)]9,0) arc (35:65:-25mm);
\draw[red] ([shift=(-15:-25mm)]9,0) arc (-15:15:-25mm);
\draw[red] ([shift=(-65:-25mm)]9,0) arc (-65:-35:-25mm);
\begin{scope}[shift={(9,0)}]
\foreach \angle/\count in {75/1,285/2,255/3,105/4}
{
\draw[red] (\angle:1.8cm) -- (\angle:2.2cm);
}
\end{scope}
    \end{tikzpicture}
 }
\caption{Pinnings  $\longleftrightarrow$  collections of decorated flags, for    ${\rm PGL_3}$, and $w_0=s_1s_2s_1$.}
\label{pin4}
\end{figure}

Given   a pinning $p_{1}: \beta_1 \to  \beta_2$  and a   reduced decomposition ${\bf i}$  of $w_0$,   there is a unique collection of   decorated flags, as illustrated on the right of Figure \ref{pin4}:   
\be \la{INTF}
  \A_1 \stackrel{s_{i_1}}{\lra}   \A_2\stackrel{s_{i_2}} {\lra}  \ldots     \stackrel{s_{i_N}}{\lra} \A_{N+1}, \qquad \beta_1 = \B_1, \beta_2 = \B_{N+1}, \qquad h(\A_m, \A_{m+1})=1. 
\ee
The flags $(\B_m, \B_{m+1})$ are at the position $s_{i_m}$ and  the $h-$distance between them  is   $1$.

  Recall   the companion reduced decomposition ${\bf i}^*$   in  (\ref{AST}). Consider  the following data: 

\begin{itemize}

\vskip 1mm\item a cyclic collection of $d$ flat sections $\{\beta_1, \ldots , \beta_{d}\}$  of ${\cal L}_{\cal B}$   on $S^1$, referred to  as flags; 

\vskip 1mm\item  pinnings $p_m: \beta_m\to   \beta_{m+1}$ between them,  $m \in \Z/d\Z$;

\vskip 1mm\item 
a sequence  ${\bf i}, {\bf i}^*, {\bf i}, {\bf i}^*, \ldots, {\bf i}, {\bf i}^*$ of     reduced   words for $w_0$, assigned to pinnings.
 
 \end{itemize}
 
 It is equivalent to the following data, see Figure \ref{pin4}:   
 
 \begin{itemize}

\vskip 1mm\item $d$ sequences of decorated flags as in (\ref{INTF}), each containing $N+1$ decorated flags.  

\vskip 1mm\item The last decorated flag of each sequence and the first decorated flag of the next sequence  have the same underlying flag.  

\vskip 1mm\item The relative position of   decorated flags in each sequence is described by the corresponding  reduced word for the collection ${\bf i}, {\bf i}^*, {\bf i}, {\bf i}^*, \ldots, {\bf i}, {\bf i}^*$. 
 \end{itemize}
 
 Indeed, given the first type of data we proceed as in (\ref{INTF}). Given the second one, we assign to each sequence of decorated flags the first and the last ones in the sequence. 
 
Now let us define a shift of each   pinning $p_{m}: \beta_m \to  \beta_{m+1}$ by one step to the right. 

The pinnings $p_m: \beta_m \to \beta_{m+1}$ and $p_{m+1}:  \beta_{m+1} \to \beta_{m+2}$     
give rise  a collection of  
 flags
 $$
 \B_1, \ldots, \B_{N+1}, \ldots , \B_{2N+1}, 
 $$
 with  $\beta_m= \B_1,  \beta_{m+1} = \B_{N+1}, \beta_{m+2} = \B_{2N+1}$. Take 
 the  first   $N+2$ flags:  
  \be \la{INTFa}
  \B_1 \stackrel{s_{i_1}}{\lra}   \B_2\stackrel{s_{i_2}} {\lra}  \ldots     \stackrel{s_{i_N}}{\lra} \B_{N+1}\stackrel{s^*_{i_1}}{\lra} \B_{N+2}.
  \ee
  The pinning $p_{m}: \beta_m \to  \beta_{m+1}$ determines  decorations of the first $N+1$ flags:
   \be \la{INTFa.1}
  \A_1 \stackrel{s_{i_1}}{\lra}   \A_2\stackrel{s_{i_2}} {\lra}  \ldots     \stackrel{s_{i_N}}{\lra} \A_{N+1}.
  \ee  
  There exists a unique decoration $\A_{N+2}$ of the   flag $\B_{N+2}$ such that $h(\A_{N+1}, \A_{N+2})=1$. 
  One easily checks that one has 
  $h(\A_2, \A_{N+2})  =1$.

  So we get a new pinning
  $p_m':  \A_2 \lra \A_{N+2}$. 
  Therefore,  given a  collection of pinnings $(p_1, ..., p_{d})$, and    
   applying this construction for each $m \in \Z/d\Z$, we get   a new collection of pinnings $(p'_1 , ... , p'_{d})$. 
    It depends only on the first generator $s_{i_1}$ in the reduced word   ${\bf i}$. Indeed, given a generic pair of flags $(\B_1, \B_2)$, the first flag 
    of the intermediate collection of flags related to any reduced word for $w_0$ depends only on the first generator in the reduced word.

  \bd Suppose that the reduced word  ${\bf i}$ starts with $s_i$. Then the  transformation $T(s_i)$ maps pinnings $(p_1, ..., p_{d})$ to the pinnings $(p'_1 , ... , p'_{d})$. 
  \ed
  
\begin{figure}[ht]
\epsfxsize 200pt
\center{
\begin{tikzpicture}[scale=0.7]
\foreach \angle/\count in {30/1,150/2,210/3,330/4}
{\node[thick] at (\angle:2cm) {$\circ$};}
\foreach \angle/\count in {60/1,180/2,300/3}
{\node at (\angle:1.8cm) {\small $s_2$};
\node at (\angle+60:1.8cm) {\small $s_1$};}
\foreach \count in {5,6}
{\node at (30+60*\count:2.6cm) {\small $\B_\count$};}
\foreach \count in {2,3}
{\node[red] at (30+60*\count:2.6cm) {\small $\A_\count$};}
\begin{scope}[shift={(9,0)}]
\node[red ] at (160:2.6) {\small $\A_2$};
\node[red ] at (210:2.6) {\small $\A_3$};
\node[red ] at (270:2.6) {\small $\A_4^r$};
\node[red ] at (320:2.6) {\small $\A_5^r$};
\foreach \count in {1,6}
{\node at (30+60*\count:2.6cm) {\small $\B_\count$};}
\node[red] at (150:2cm) {\Large $\bullet$};
\node[red] at (-30:2cm) {\Large $\bullet$};
\foreach \angle/\count in {30/1,90/2,210/3,270/4}
{\node[thick] at (\angle:2cm) {$\circ$};}
\foreach \angle/\count in {60/1,180/2,300/3}
{\node at (\angle:1.8cm) {\small $s_2$};
\node at (\angle+60:1.8cm) {\small $s_1$};}
\end{scope}
\begin{scope}[shift={(-9,0)}]
\draw [red, dashed]  (9,0) circle (20mm);
\node [red] (a1) at (9,2) {\Large $\bullet$};
\node [red] (b1) at (9,-2) {\Large $\bullet$};
\node [red] (g1) at (6,0) {\small $1$}; 
\node [red] at (6.9,-2.2) {\small $1$};
\node [red]  at (6.9,2.2) {\small $1$};
\node (o) at (9.7, 2.54) {\small $\A_1^r$};  
\node (p) at (9.7,-2.54) {\small $\A_4^l$};
\node[red] (q) at (8.3,-2.54) {\small $\A_4^r$};  
\node[red] (r) at (8.3,2.54) {\small $\A_1^l$};  
\draw[red] ([shift=(220:25mm)]9,0) arc (220:245:25mm);
\draw[red] ([shift=(165:25mm)]9,0) arc (165:195:25mm);
\draw[red] ([shift=(115:25mm)]9,0) arc (115:145:25mm);
\begin{scope}[shift={(9,0)}]
\foreach \angle/\count in {75/1,285/2}
{
\draw[red] (\angle+180:1.8cm) -- (\angle+180:2.2cm);
\draw[black] (\angle:1.8cm) -- (\angle:2.2cm);
}
\end{scope} 
\end{scope}
\draw [ultra thick, -latex] (3.6,0) -- (5.4,0);
\node at (4.5, 0.7) {${\rm T}_x(s_1)$};
\draw [red, dashed]  (9,0) circle (20mm);
\node [red] (g1) at (6,0) {\small $1$};
\node [red]  at (7.5,-2.65) {\small $1$};   
\node [red] at (10.5,-2.7) {\small $1$};
\draw[red] ([shift=(-75:25mm)]9,0) arc (-75:-50:25mm);
\draw[red] ([shift=(-140:25mm)]9,0) arc (-140:-100:25mm);
\draw[red] ([shift=(-10:-25mm)]9,0) arc (-10:20:-25mm);
\begin{scope}[shift={(9,0)}]
\foreach \angle/\count in {-40/1,160/2}
{
\draw[red] (\angle:1.8cm) -- (\angle:2.2cm);
}
\end{scope}
    \end{tikzpicture}
 }
\caption{The action of a  generator $s_1 \in {\mathbb B}_{\rm sl_3}$ on  a  pinning $p_{1}: \beta_1 \to \beta_2$.  The 
new pinning   is obtained   by moving the pair of special   points    one step to the right.}
\label{pin3}
\end{figure}

Pick a special point $x\in S^1$,  a generator $s_i \in  {\mathbb B}_{\mathfrak g}$, and a reduced word for   $w_0$  starting from  $s_i$.  
We get an automorphism $T_{x}({s_i})$ of the  space 
${\mathscr P}_{\G,  S^1}$,  lifting the one of the space ${\mathscr X}_{\G,  S^1}$.
\be
\begin{split}
&T_{x}({s_i}): ({\cal L}, \beta; p) \lra ({\cal L},  \beta'; p').\\
\end{split}
\ee

Recall the discrete group $\Gamma_{\G, \bS}$ in (\ref{GGBG}) acting on the   space ${\mathscr P}_{\G, \bS}$, and its subgroup $\Gamma'_{\G, \bS}$.

\bt
 \la{TTHH2x}
  Pick   a   special point $x_1$    at a boundary component $\pi$ of  $\bS$.  

\begin{enumerate}

\vskip 1mm\item  There is an action    $T_{x_1}$ of the  
   braid group ${\mathbb B}^{(\pi)}_{\mathfrak{g}}$  
on   ${\mathscr P}_{\G, \bS}$.

\vskip 1mm\item    The   subgroup $\widetilde {\mathbb B}^{(\pi)}_{\mathfrak g}\subset {\mathbb B}^{(\pi)}_{{\mathfrak g}}$ acts by quasi-cluster Poisson transformations. 

\vskip 1mm\item The  group $\gamma'_{\G, \bS}$  
acts  by quasi-cluster Poisson automorphisms of    the space ${\mathscr P}_{\G, \bS}$. 
\end{enumerate}
\et 

\begin{proof}   Parts 1) $\&$ 3) are proved just as in Theorem \ref{TTHH1}. 
Part 2) is proved in Section \ref{sect10.4}.
\end{proof}  

\subsubsection{Braid group action on  ${\mathscr A}_{\G, \bS}$.}  Let ${\mathscr A}_{\Z}$ be the set of sequences $\{\A_k\}_{k\in \Z}$ of decorated flags such that every adjacent pair $(\A_k, \A_{k+1})$ is generic. 
 Let $i\in {\rm I}$. For each $k\in \mathbb{Z}$,   pick the unique $\A_k'$ so that \vskip 1mm

\vskip 1mm1. if $k$ is odd, then $w(\A_k, \A_k')= s_i,~w(\A_k', \A_{k+1})=s_i w_0, ~h(\A_k, \A_{k}')=1$,

\vskip 1mm2.  if $k$ is even, then $w(\A_k, \A_k')= s_{i^\ast},~w(\A_k', \A_{k+1})=s_{i^\ast} w_0, ~h(\A_k, \A_{k}')=1$. \vskip 1mm

Then the pair $(\A_k', \A_{k+1}')$ is generic.  So we get a transformation $T(s_i)$ of the set ${\mathscr A}_{\Z}$  by setting
 \[
 T(s_i): ~ \{\A_k\} \lms \{\A_k'\}.
 \]
 Just as in Lemma \ref{PROP5.5},   transformations $T(s_i)$ generate an action of the braid group $\B_{\mathfrak g}$ on ${\mathscr A}_\Z$. 
  
  \vskip 2mm
  
  A sequence $\{\A_k\}\in {\mathscr A}_\Z$ is said to be $(d, g)$-periodic if $\A_{k+d}= g\cdot \A_k$ for every $k\in \Z$. Let us introduce two more transformations on ${\mathscr A}_{\mathbb{Z}}$: \vskip 1mm

\vskip 1mm3.  the transformation $S_d$ shifting $\{ \A_k\}$ to $\{ \A_{k+d}\}$;

\vskip 1mm 4.  the transformation $M_g$ taking $\{ \A_k\}$ to $\{g\cdot \A_k\}$.
\vskip 1mm

A sequence $\A_\bullet:=\{ \A_k\}$ is $(d, g)$-periodic if and only if $S_d(\A_\bullet)= M_g(\A_\bullet)$.

  \bl 
  \la{asobcosp}
  Let $b$ belong to the centralizer of $w_0^d$ in ${\Bbb B}_{\mathfrak g}$. A sequence $\A_\bullet :=\{\A_k\}\in {\mathscr A}_\Z$ is $(d, g)$-periodic if and only if  $T(b)(\A_\bullet) $ is $(d, g)$-periodic.
  \el
  \begin{proof} 

By definition, the transformations $T, S, M$ satisfy the following relations
\[
T(b)\circ S_1 = S_1 \circ T(w_0b w_0^{-1}); \hskip 7mm T(b)\circ M_g = M_g \circ T(b).
\]  
If particular, if $b$ belongs to the centralizer of $w_0^d$ in ${\Bbb B}_{\mathfrak g}$, then we have
$
T(b) \circ S_d = S_d \circ T(b).
$
 
Hence $S_d(\A_\bullet)= M_g (\A_\bullet)$ if and only if $S_d(T(b)(\A_\bullet))= M_g(T(b)(\A_\bullet))$.  
  \end{proof}
  
 Recall  the moduli space ${\mathscr A}_{\G, \bS}$ parametrizing  pairs $(\mathcal{L}, \{\alpha\})$, where ${\cal L}$ is a twisted unipotent $\G$-local system, and $\{\alpha\}$ are flat sections of the associated bundle ${\cal L}_{\mathscr A}:={\cal L}\times_{\G} {\mathscr A}$ near special points and punctures of $\bS$. Let $\pi$ be a boundary component of $\bS$ carrying  special points $s_1, s_2,\ldots, s_d$. By the parallel transport of flat sections $\alpha_{s_k}$ to the point $s_1$, we get a $(d, g)$-periodic sequence in ${\mathscr A}_{\Bbb Z}$, where $g$ is the monodromy around $\pi$. Conversely, every $(d, g)$-periodic sequence sitting on the fiber of ${\cal L}_{\mathscr A}$ over $s_1$ provides the data of flat sections on the special points $s_1, \ldots, s_d$.  
Using Lemma \ref{asobcosp}, we apply the transformation $T(b)$ to the obtained $(d, g)$-periodic sequence by parallel transport to the fiber over $s_1$, and then transport them back to the other special points. It defines the action of ${\Bbb B}_{\mathfrak g}^\pi$ on $\mathcal{A}_{\G, \bS}$.
 
\medskip

\subsection{Quasi-cluster nature of the braid group actions} \la{sect10.4}

 \medskip
 
\subsubsection{Punctured disk case.} We shall begin with the case when $\bS$ is a punctured disk with $d$   special points. Let ${\bf i}=(i_1,\ldots, i_m)$ be a reduced word for $w_0$. We write ${\bf i}^\ast=(i_1^\ast, \ldots, i_m^\ast)$ and ${\bf i}_+=(i_2, \ldots, i_{m}, i_1^\ast)$. 
Take the central triangulation of $S$ as on Figure \ref{pspace.central.tri}, where $p_{k, k+1}=(\A_{k}^l, \A_{k+1}^r)$ are pinnings over $(\B_k, \B_{k+1})$.
Let  $h_k= h(\A_k^r, \A_k^l)$, or equivalently $\A_k^r=\A_k^l\cdot h_k$.  We set $\forall j \in {\rm I}$:
\[
K_{k,j}=\left\{ \begin{array}{ll} \alpha_j(h_k) &\mbox{if $k$ is odd},\\ \alpha_{j^\ast}(h_k) & \mbox{if $k$ is even}. \end{array} \right.  
\]
 \begin{figure}
 \begin{center}
\begin{tikzpicture}
\draw (0,0) circle [radius=1.5];
\draw (0,0) circle [radius=0.2];
\draw (0:0.2) -- (0:1.5);
\draw (60:0.2) -- (60:1.5);
\draw (120:0.2) -- (120:1.5);
\draw (180:0.2) -- (180:1.5);
\draw (240:0.2) -- (240:1.5);
\draw (300:0.2) -- (300:1.5);
\draw[dotted, thick] (75:1) arc (75:105:1cm);
\node[red] at (0, 0) {$\B$};
\node at (0:1.5) {$\bullet$};
\node at (60:1.5) {$\bullet$};
\node at (120:1.5) {$\bullet$};
\node at (180:1.5) {$\bullet$};
\node at (240:1.5) {$\bullet$};
\node at (300:1.5) {$\bullet$};
\node at (0:1.8) {$\B_4$};
\node at (60:1.8) {$\B_5$};
\node at (120:1.8) {$\B_{2d}$};
\node at (180:1.8) {$\B_1$};
\node at (240:1.8) {$\B_2$};
\node at (300:1.8) {$\B_3$};
\node[blue] at (30:1.8) {{\footnotesize $ p_{45}$}};
\node[blue] at (150:1.8) {{\footnotesize $ p_{2d,1}$}};
\node[blue] at (210:1.8) {{\footnotesize $ p_{12}$}};
\node[blue] at (270:1.8) {{\footnotesize $ p_{23}$}};
\node[blue] at (330:1.8) {{\footnotesize $ p_{34}$}};
\node at (30:0.5) {{\footnotesize $ {\bf i}^\ast$}};
\node at (150:.4) {{\footnotesize $  {\bf i}^\ast$}};
\node at (210:.4) {{\footnotesize $ {\bf i}$}};
\node at (270:.4) {{\footnotesize $  {\bf i}^\ast$}};
\node at (330:.4) {{\footnotesize $  {\bf i}$}};
\end{tikzpicture}
\end{center}
\caption{The space $\mathscr{P}_{\G, \bS}$ with central triangulation.}
\la{pspace.central.tri}
\end{figure}

We pick the central angle in each triangle, and use ${\bf i}$ and ${\bf i}^\ast$ alternatively.
Recall the primary coordinates $\{P_{k,1}, \ldots, P_{k, m}\}$ for each $(\B, \B_k, \B_{k+1}; p_{k, k+1})$. 
By Lemma \ref{LEM9.4}, the following functions form  a rational coordinate system of $\mathscr{P}_{\G, S}$:
\be \la{PRCO}
\p_{\bf i}=\left\{ K_{k,1}, \ldots, K_{k, r}; P_{k,1}, \ldots, P_{k, m} \right\}_{1\leq k \leq d}.
\ee

Recall the following coefficients in Theorem  \ref{TH9.5}: 
\[
b_{st}:=\langle \alpha_s^{\bf i}, \alpha_t^{\bf i}\rangle,  \hskip 7mm c_{jt}:=\langle \alpha_j,  \alpha_t^{\bf i}\rangle,  \hskip 7mm d_{ij}:=\langle \alpha_i,  \alpha_j\rangle:
\]
\bp With respect to the cluster Poisson bracket on $\mathscr{P}_{\G, \bS}$, we have
\begin{align} \left\{ \log P_{k, s}, \log P_{l, t}\right\}  &=  \left\{ \begin{array}{ll} {\rm sgn}(t-s)\cdot b_{st} & \mbox{if } k=l\\
0, & \mbox{otherwise}\\ \end{array} \right. \la{PP.relation} \\
\left\{ \log K_{k, j}, \log P_{l, t}\right\} &= \left\{ \begin{array}{ll} c_{jt} & \mbox{if } k=l\\
-c_{jt} & \mbox{if } k=l+1 \\
0, & \mbox{otherwise}\\ \end{array} \right. \la{KP.relation}
\end{align}
 If $d=1$, then all the coordinates $K_{k, i}$ commute. If $d\geq 2$ then
\be
\la{KK.relation}
\left\{ \log K_{k, i}, \log K_{l, j}\right\}= \left\{ \begin{array}{ll} d_{ij} & \mbox{if } l=k+1\\
-d_{ij} & \mbox{if } l=k-1 \\
0, & \mbox{otherwise}\\ \end{array} \right. 
\ee
\ep
\begin{proof} The space $\mathscr{P}_{\G, \bS}$ is obtained by gluing $d$ copies of $\mathscr{P}_{\G,t}$:
\[
g: ~ \left(\mathscr{P}_{\G,t}\right)^{d} \longrightarrow \mathscr{P}_{\G,S}.
\]
Furthermore, the gluing map $g$ preserves the Poisson structures on both spaces. 

Let $\{P_{k,1}, \ldots, P_{k,m}; L_{k,1}, \ldots, L_{k, r}; R_{k, 1}, \ldots R_{k, r}\}$ be the coordinates of the $k^{th}$ copy of $\mathscr{P}_{\G,t}$. By definition,  
 \be
 \la{k.pinning;cord}
 K_{k,j}=\left\{ \begin{array}{ll} R_{k-1,j}L_{k,j} &\mbox{if $k$ is odd},\\ R_{k-1,j^\ast}L_{k,j^\ast} & \mbox{if $k$ is even}. \end{array} \right.
 \ee
 The rest follows directly from Theorem  \ref{TH9.5}.
\end{proof}

Let $i_1=i^\ast$.  Let $\p_{{\bf i}_+}=\{K_{k,j}'; P_{k, s}'\}$ be the   coordinates (\ref{PRCO}) on $\mathscr{P}_{\G, S}$ associated  to the reduced word ${\bf i}_+$. Recall the automorphism $T(s_i)$ of the space $\mathscr{P}_{\G, S}$. 
\begin{lemma} \la{L2.8.19.1} We have
\begin{align}
T(s_{i})^\ast P_{k,s}'&= \left \{\begin{array}{ll} P_{k, s+1}, & \mbox{if } 1\leq s<m;\\ P_{k+1, 1} K_{k+1, i}^{-1} & \mbox{if }s=m.\end{array}\right.\nonumber\\
T(s_{i})^\ast K_{k,j}'&= K_{k,j}K_{k, i}^{-C_{ji}}. \nonumber
\end{align}
\end{lemma}
\begin{proof} We shall prove the first set of formulas for $k=1$ and the second set for $k=2$. The proof goes through the same line for general $k$. Take the following decompositions of  $p_{1,2}$ and $p_{2,3}$:
\[
p_{1,2}:\hskip 3mm \A_1^l=\A_{1,0}\stackrel{s_{i_1^\ast}}{\lra}\A_{1,1}\stackrel{s_{i_2^\ast}}{\lra}\cdots\stackrel{s_{i_m^\ast}}{\lra}\A_{1,m}=\A_{2}^r.
\]
\[
p_{2,3}:\hskip 3mm \A_2^l=\A_{2,0}\stackrel{s_{i_1}}{\lra}\A_{2,1}\stackrel{s_{i_2}}{\lra}\cdots\stackrel{s_{i_m}}{\lra}\A_{2,m}=\A_{3}^r.
\]
The map $T(s_i)$ shifts the pinning $p_{1,2}$ to
\[
p_{1,2}': \hskip 3mm \A_{1,1}\stackrel{s_{i_2^\ast}}{\lra}\A_{1,2}\stackrel{s_{i_3^\ast}}{\lra}\cdots \stackrel{s_{i_m^\ast}}{\lra} \A_{1,m}\stackrel{s_{i_1}}{\lra} \A_{2,1}'.
\]
In particular, the last pair are
\be
\la{shift.cartanh2}
(\A_{1,m}, \A_{2,1}') = (\A_2^r, \A_{2,1}') = (\A_2^l \cdot h_2, \A_{2,1}\cdot s_{i^\ast}(h_2))
\ee
By the definition of the primary coordinates,  
$
T(s_{i})^\ast P_{1,s}'= P_{1, s+1},~ \forall s=1,\ldots, m-1. 
$ 

For $s=m$, by \eqref{shift.cartanh2} and \eqref{cartan.action.pinning}, we get
\[
T(s_i)^\ast P_{1,m}' = P_{2,1} \alpha_{i^\ast}(h_2^{-1}) = P_{2,1} K_{2,i}^{-1}.
\]
From \eqref{shift.cartanh2} we also get $h_2'= s_{i^\ast}(h_2)$. Therefore
\[
T(s_i)^\ast K_{2, j}' = \alpha_{j^\ast} \left(h_2'\right) = (s_{i^\ast}\cdot \alpha_{j^\ast})(h_2)= K_{2,j} K_{2,i}^{-C_{ji}}.
\]
\end{proof}

\begin{proposition} 
\la{poisson.pre.braid.map.st}
The automorphism $T(s_i)$ preserves the cluster Poisson bracket of $\mathscr{P}_{\G, \bS}$.
\end{proposition}
\begin{proof}  We prove the Proposition by direct calculation. Note that
\[
\alpha_s^{{\bf i}_+}= s_i \cdot \alpha_{s+1}^{\bf i}, \qquad\forall 1\leq s<m;\hskip 5mm \alpha_{m}^{{\bf i}_+}= \alpha_i.
\]
Let us set
\[
b_{st}':=\langle \alpha_s^{{\bf i}_+}, \alpha_t^{{\bf i}_+}\rangle, \hskip 5mm c_{jt}':= \langle \alpha_j, \alpha_t^{{\bf i}_+}\rangle. 
\]

For the $K$$-K$ relations \eqref{KK.relation}, if $d\geq 2$, then 
\[
\left\{T(s_i)^* \log K_{k,j}', ~ T(s_i)^*\log K_{k+1, l}' \right\} = \langle s_i \cdot \alpha_j, ~ s_i \cdot \alpha_l\rangle = d_{jl}. 
\]
The other cases follow similarly.

For the $K$$-P$ relations \eqref{KP.relation}, we prove the $k=l$ case.  The other cases are similar. If $1\leq s<m$, then
\[
\left\{T(s_i)^* \log K_{k,j}', ~ T(s_i)^*\log P_{k, s}' \right\} = \left\{ \log K_{k,j}K_{k,i}^{-C_{ji}}, ~ \log P_{k, s+1} \right\}=  \langle s_i \cdot \alpha_j, ~ \alpha_{s+1}^{\bf i}\rangle = c_{js}'. 
\]
If $s=m$, then
\[
\left\{T(s_i)^* \log K_{k,j}', ~ T(s_i)^*\log P_{k, m}' \right\} = \left\{ \log K_{k,j}K_{k,i}^{-C_{ji}}, ~ \log P_{k+1, 1}K_{k+1, i}^{-1} \right\}=  \langle s_i \cdot \alpha_j, ~ -\alpha_i\rangle = c_{jm}'. 
\]

For the $P$$-P$ relations \eqref{PP.relation}, if $1\leq s<m$, then
\begin{align}
&\left\{T(s_i)^* \log P_{k,s}', ~ T(s_i)^*\log P_{k, m}' \right\} =\left\{ \log P_{k,s+1}, ~ \log P_{k+1, 1} K_{k+1, i}^{-1} \right\}\nonumber\\
= &\left\{  \log P_{k,s+1}, ~- \log  K_{k+1, i} \right\} = \langle \alpha_{s+1}^{\bf i}, - \alpha_i\rangle = \langle \alpha_{s}^{{\bf i}_+}, \alpha_m^{{\bf i}_+}\rangle = b_{sm}'. \nonumber
\end{align}
Similarly, we have the next two identities:
\be
\begin{split}
&\left\{T(s_i)^* \log P_{k,m}', ~ T(s_i)^*\log P_{k+1, s}' \right\} =\left\{ \log P_{k+1, 1} K_{k+1, i}^{-1} , \log P_{k+1, s+1} \right\}=  \langle \alpha_{1}^{{\bf i}} -\alpha_i,  \alpha_{s+1}^{\bf i}\rangle = 0.\\
&\left\{T(s_i)^* \log P_{k,m}', ~ T(s_i)^*\log P_{k+1, m}' \right\} =\left\{ \log P_{k+1, 1} K_{k+1, i}^{-1} , \log P_{k+2, 1} K_{k+2, i}^{-1}\right\}\nonumber\\
=&  \left\{ \log P_{k+1, 1} K_{k+1, i}^{-1} , - \log K_{k+2, i}\right\}= \langle \alpha_{1}^{{\bf i}} -\alpha_i,  - \alpha_{i}\rangle  = 0.
\end{split}
\ee
The other cases follow similarly.
\end{proof}

Let $\chi_{\bf i}=\{X_v\}$ be the cluster Poisson chart associated to the triangulation and reduced words as on Figure \ref{pspace.central.tri}. Equivalently, it corresponds to an open torus embedding still denoted by 
$
\chi_{\bf i}:~  \mathbb{G}_m^N \hookrightarrow \mathscr{P}_{\G, \bS}.
$
\bp 
\la{cluster.monomial.p.x.trans}
The transition map between the charts $\chi_{\bf i}$ and $\p_{\bf i}$ is given by  Laurent monomials, i.e. every $P_{k,s}$ and $K_{k,j}$ in $\p_{\bf i}$ can be expressed as a monomial in  $X_v$ in $\chi_{\bf i}$, and vice versa.
\ep

\begin{proof} 
As shown on \eqref{p-x.trans}, the $P$-coordinates can be expressed as Laurent monomials of cluster Poisson variables. 
For the $K$-coordinates, we need the following easy lemma 
\bl
Let $M$ be a Laurent monomial in a cluster chart $\chi = \{X_v\}$   of a cluster Poisson variety. If $\{M, X_v\}=0$   for every unfrozen $X_v$  in $\chi$, then $M$ is a Laurent monomial in every cluster Poisson chart of the variety.
\el
\begin{proof} We use the notations in Section \ref{sec.quasicl}. Let $M=\prod_{i} X_i^{n_i}$. We have
\[
\left\{M, X_v\right\} =\left( \sum_{i} n_i \widehat{\varepsilon}_{iv}\right) M X_v =0.
\]
As a consequence, we see that 
\[\sum_{i} n_i {\varepsilon_{iv}} = \sum_{i} n_i {\varepsilon_{iv}} d_v =0.
\]
Now use the mutation (\ref{pois.clust.mut}), we see that $M$ remains a Laurent monomial in its neighbored cluster chart. Repeating the same procedure recursively, we prove that $M$ is a Laurent monomial in every cluster Poisson chart of the variety.
\end{proof}
Since the $K$-coordinates commute with every unfrozen cluster Poisson variables,  they are Laurent monomials in every cluster chart. 
Since $\p_{\bf i}$ is a birational isomorphism,   the map $\p_{\bf i}\circ \chi_{\bf i}^{-1}: \mathbb{G}_m^N \lra \mathbb{G}_m^N$ is an isomorphism.  Therefore  conversely
 the coordinates $X_v$ in $\chi_{\bf i}$ are  Laurent monomials in   $\p_{\bf i}-$coordinates.
\end{proof}

The composition $T(s_i)\circ \chi_{\bf i}$ provides a new chart of $\mathscr{P}_{\G, S}$. By definition, the unfrozen part of $T(s_i)\circ \chi_{\bf i}$ coincides with the unfrozen part of $\chi_{{\bf i}_+}$ up to permutation. Combining Lemma \ref{L2.8.19.1} and Proposition \ref{cluster.monomial.p.x.trans}, the frozen variables of $\chi_{{\bf i}_+}$ can be presented as Laurent monomials of variables in $T(s_i)\circ \chi_{\bf i}$, and vice versa. Therefore $T(s_i)\circ \chi_{\bf i}$ is a quasi-cluster Poisson chart. Together with Lemma \ref{poisson.pre.braid.map.st}, we conclude that $T(s_i)$ is a quasi-cluster automorphism.

\subsubsection{Annulus case.} Let $\bS$ be an annulus with one special point on one side and $d$  special points on the other. We obtain  $\bS$ by gluing a quadrilateral together with $d-1$ triangles. So the cluster Poisson structure on $\mathscr{P}_{\G, \bS}$ is obtained by amalgamation  of $\mathscr{P}_{\G, q}$ and $d-1$ copies of $\mathscr{P}_{\G,t}$. 
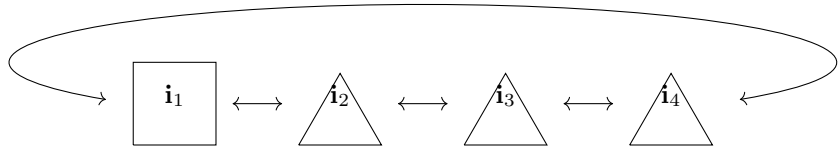
 \begin{center}
 \begin{tikzpicture}[scale=1.1]
   \node at (0.5,0.6) {\small ${\bf i}_{1}$};
\draw[<->] (3.5,1) [partial ellipse=-40:220:5cm and 0.7cm];
 \draw (0,0) -- (1,0) -- (1,1) -- (0,1) -- (0,0);
 \foreach \count in {2, 3, 4}
{ \begin{scope}[shift={(2*\count-2,0)}]
 \draw (0,0)--(0:1)--(60:1)--(0,0);
  \node at (0.5,0.6) {\small ${\bf i}_{\count}$};
 \draw[<->] (-.2, .5)--(-.8, . 5);
 \end{scope}
 }
 \end{tikzpicture}
 \end{center}
 
We make use of the reduced word  ${\bf i}_1$ of $(w_0, w_0)$ to give a decomposition of  $\mathscr{P}_{\G, q}$ , and the reduced words ${\bf i}_2, {\bf i}_3, \ldots, {\bf i}_{d}$ repeatedly to give decompositions of  the rest $d-1$ many $\mathscr{P}_{\G,t}$. Let us consider the amalgamation of  the quivers ${\bf Q}({\bf i}_1)$, ..., ${\bf Q}({\bf i}_{d})$:
\be
\la{quiver.annulaus,case}
{\bf Q}({\bf i}_1)\ast\ldots \ast {\bf Q}({\bf i}_{d}).
\ee
Here for each $j\in \Z/d \Z$ and for each $i\in {\rm I}$, the rightmost vertex of ${\bf Q}({\bf i}_j)$ at level $i$ is glued with the leftmost vertex of ${\bf Q}({\bf i}_j)$ at level $j$. The vertices of \eqref{quiver.annulaus,case} parametrize  functions of $\mathscr{P}_{\G, \bS}$. Precisely:
\begin{itemize}
\vskip 1mm\item the unfrozen part corresponds to unfrozen cluster Poisson variables $X_1, \ldots, X_{(d+1)m}$ of $\mathscr{P}_{\G, \bS}$;
\vskip 1mm\item the frozen part assigned to the top side of $q$ parametrizes primary coordinates $P_{0, 1}, \ldots, P_{0,m}$;
\vskip 1mm\item the frozen part assigned to the bottom sides of $q$ and $t$ parametrizes primary coordinates $P_{k, 1}, \ldots, P_{k, m}$ for $k=1,\ldots, d$.
\end{itemize}
By definition, the extra cluster Poisson variables are given by $P_{k, s}$ when  $\alpha_{s}^{{\bf i}_k}$ is a simple positive root. 
By  Proposition \ref{switchi,ibar.}, the Poisson structure on $\mathscr{P}_{\G, \bS}$ is encoded by the quiver \eqref{quiver.annulaus,case}.

Consider a  reduced word ${\bf i}=(i_1, \ldots, i_m)$ of $w_0$ starting with $i_1=i$. 

Below we assume that $d$ is even, and consider the following two cases.
 \begin{itemize}
\vskip 1mm\item[1.] ${\bf i}_1=({\bf i}, \overline{\bf i})$, and $({\bf i}_2, {\bf i}_3, \ldots, {\bf i}_{2d}) =({\bf i}^\ast, {\bf i}, {\bf i}^\ast, \ldots, {\bf i}^\ast)$. The corresponding quiver \eqref{quiver.annulaus,case} is denoted by ${\bf Q}_1$, whose vertices parametrize  unfrozen cluster variables and primary coordinates 
\[ X_1, \ldots, X_{(d+1)m}, P_{0, 1}, \ldots, P_{0,m}, \ldots, P_{d,1 }, \ldots, P_{d, m}.
\]

\vskip 1mm\item[2.]  ${\bf i}_1=(i_2, \ldots, i_m, \overline{\bf i}, i_1^\ast)$, and  $({\bf i}_2, {\bf i}_3, \ldots, {\bf i}_{d}) =({\bf i}_+^\ast, {\bf i}_+, {\bf i}_+^\ast, \ldots, {\bf i}_+^\ast)$ - recall that $d$ is even. 
 Let us denote the quiver \eqref{quiver.annulaus,case} by ${\bf Q}_2$ and the functions by
\[
X_1', \ldots, X_{(d+1)m}', P_{0, 1}', \ldots, P_{0,m}', \ldots, P_{d,1 }', \ldots, P_{d, m}'.
\]
\end{itemize}

Recall  the automorphism $T(s_i)$ of $\mathscr{P}_{\G, \bS}$. As in \eqref{k.pinning;cord}, we consider the functions $K_{k, j}$ assigned to $d$ many special points on the outer side of $\bS$.
\begin{lemma} For $j \in \Z/  (d+1)m \Z$, we have
\be
\la{unfrozen.permutation1}
 T(s_i)^\ast X_j'= X_{j+1}.
\ee
For $k=1,\ldots, d$, we have
\be 
\la{frozen.permutation2}
T(s_{i})^\ast P_{k,s}'= \left \{\begin{array}{ll} P_{k, s+1}, & \mbox{if } 1\leq s<m;\\ P_{k+1, 1} K_{k+1, i}^{-1} & \mbox{if }s=m.\end{array}\right.
\ee
For $k=0$, we have
\be
\la{frozen.permutation3}
T(s_i)^\ast P_{0, s}'= P_{0,s}.
\ee
\end{lemma}
\begin{proof} It follows by the same argument as in the proof of Lemma \ref{L2.8.19.1}.
\end{proof}
By the same argument in the proof of Proposition \ref{poisson.pre.braid.map.st}, the automorphism $T(s_i)$ preserves the cluster Poisson bracket of $\mathscr{P}_{\bS}$.
By \eqref{unfrozen.permutation1}, under $T(s_i)$, the unfrozen cluster Poisson variables of ${\bf Q}_2$ coincides with unfrozen variables of ${\bf Q}_1$ up to permutation. By \eqref{frozen.permutation2} - \eqref{frozen.permutation3}, the frozen cluster Poisson variables of ${\bf Q}_2$ are expressed as monomials in terms of cluster Poisson variables of ${\bf Q}_1$. So $T(s_i)$ is a quasi cluster automorphism. 

\vskip 2mm
The case when a boundary component has an odd number of special points is treated similarly.

\subsubsection{The total number of punctures and boundary components of  $\bS$ is greater than $2$.} In this case 
 one can  find a loop  $\alpha$ separating $\bS$ into an annulus with boundary components 
$\alpha$ and $\beta$, each containing 
  special points, and the rest.   The action of $T(s_i)$ on the component $\beta$, see the picture,   only changes the coordinates assigned to the annulus. 
  Thus we can reduce it to the annuls case. 
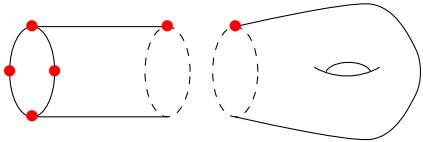
\begin{center}
\begin{tikzpicture}[scale=0.3]
\draw (0,0) ellipse (1cm and 2cm);
\draw[dashed] (6,0) ellipse (1cm and 2cm);
\draw (0,2) -- (6,2);
\draw (0,-2) -- (6,-2);
\draw[dashed] (9,0) ellipse (1cm and 2cm);
\draw plot [smooth] coordinates {(9,2) (15,3) (17,1) (17,-1) (15,-3)  (9,-2)};
\draw (15,0) arc (10:170:1cm and .5cm);
\draw (15.4,0.2) arc (-60:-120:2.9cm and 2.9cm); 
\node[red] at (0,2) {$\bullet$};
\node[red] at (6,2) {$\bullet$};
\node[red] at (9,2) {$\bullet$};
\node[red] at (0,-2) {$\bullet$};
\node[red] at (1,0) {$\bullet$};
\node[red] at (-1,0) {$\bullet$};
\end{tikzpicture}
\end{center}

\subsubsection{The quasi-cluster nature of the braid group action on the space ${\mathscr A}_{\G, \bS}$.} It is proved similarly to the case of the space ${\mathscr P}_{\G, \bS}$ 
considered above. The quasi-cluster transformations we need are obtained by using the same sequence of quiver mutations 
 for the spaces ${\mathscr A}_{\G, \bS}$ and ${\mathscr P}_{\G, \bS}$, followed by a monomial transformation on the space ${\mathscr P}_{\G, \bS}$, and the dual one for the space ${\mathscr A}_{\G, \bS}$. 

\medskip

\subsection{Equivariance of quantum cluster structures} \la{Sec8.4}

\medskip

Theorem \ref{Th8.5}   generalizes  \cite[Theorem 6.1]{FG07}, which settled the case when $\G$ is of type $A_n$ for the group $\Gamma_\bS$. Recall the discrete group $\Gamma_{\G, \bS}$, 
see (\ref{GGBG}), and its subgroup $\gamma'_{\G, \bS}$.

   \bt  \la{Th8.5} Let $\G$ be an arbitrary split semi-simple algebraic group over $\Q$. Then:

\vskip 1mmi) The   group $\gamma'_{\G, \bS}$ of $\bS$ acts by automorphisms of  algebras ${\cal O}_q({\mathscr P}_{\G, \bS})$,  ${\cal O}_q({\cal D}_{\G, \bS})$, given by quasi-cluster transformations, 
which are cluster on the  subgroup $\Gamma_\bS \times W^n$. 

\vskip 1mmii) The group $\gamma'_{\G, \bS}$ acts by unitary projective automorphisms of the principal series of representations of the $\ast$-algebra $\mathcal{A}_\hbar({\mathscr P}_{\G, \bS})$.

\et

Before we start the proof, let us explain the  core problem. We have   
the geometric action of the   group $\Gamma'_{\G, \bS}$, and proved that  it acts by   cluster transformations. 
 So the action of each  of its elements  can be lifted to a quantum cluster transformation. However we have to prove that different 
  cluster transformations providing the same geometric transformation 
 lift to the same quantum cluster transformation. Here is an example.

Each flip of an edge $E$ of an ideal  triangulation ${\rm T}$ of $\bS$ gives rise to a   quantum cluster transformation ${\bf c}^{\mathscr X}_q(E)$ of the quantized space 
${\mathscr P}_{\G, \bS}$ and a   quantum cluster transformation ${\bf c}^{\cal D}_q(E)$ of the  quantized symplectic double 
${\cal D}_{\G, \bS}$. The modular groupoid ${\rm M}_\bS$ of the decorated surface $\bS$ is generated by  flips. The only relations are the pentagon and 
quadrangle relations:

\vskip 1mm{\it The pentagon relation.} Given a pentagon ${\rm P}_5$ of an ideal triangulation ${\rm T}$, there is a sequence of five flips of the diagonals of the pentagon. 
It results in the cluster transformations 
\be
\begin{split}
&{\bf c}^{\ast}_q({\rm P}_5):= {\bf c}^{\ast}_q(E_5)\circ {\bf c}^{\ast}_q(E_4)\circ {\bf c}^{\ast}_q(E_3)\circ {\bf c}^{\ast}_q(E_2)\circ {\bf c}^{\ast}_q(E_1),
 \qquad\mbox{where $\ast$ is ${\mathscr X}$,  ${\mathscr A}$ or ${\cal D}$}.\\
\end{split}
\ee 
The pentagon relation says that this cluster transformation is trivial. 

\vskip 1mm{\it The quadrangle relation.} The flips at the non-adjacent edges commute.\footnote{Flips at two non-adjacent edges 
generate a quadrangle in the modular complex, hence the name.}

\vskip 1mmThe mapping class group $\Gamma_\bS$ is the fundamental group of the modular groupoid ${\rm M}_\bS$ of $\bS$. So we have to prove the quantum pentagon and quandrangle relations.

\begin{proof} i) $\implies$ ii).  Any quasi-cluster transformation ${\bf c}$ gives rise to an intertwiner ${\bf K}_{\bf c^\circ}$. 
Indeed, a quasi-cluster transformation is a composition of a cluster transformation followed by a Poisson monomial transformation. 
So one uses the main results of the quantization of cluster Poisson varieties \cite{FG03b} for $\hbar>0$ and Section \ref{SEC8.4} for $|\hbar| =1$  for the cluster part, and  Theorem \ref{Weil} for the monomial part. 
This way we define the action of the elements of the group $\Gamma_{\G, \bS}$ by unitary projective operators. 
It remains to prove that the relations between the generators 
produce operators $ \lambda{\rm Id}$,   $|\lambda|=1$. 

By \cite[Theorem 5.4]{FG07}, given an arbitrary quiver cluster transformation  ${\bf c}: {\bf q} \to {\bf q}$, if the quantum cluster transformations ${\bf c}^{\mathscr X}_q$ and 
${\bf c}^{\cal D}_q$ are the identity maps, then   ${\bf K}_{\bf c^\circ} = \lambda{\rm Id}$,   $|\lambda|=1$. 
This settles the   relations for   $\Gamma_\bS \times W^n$. 
The case of quasi-cluster transformations follows by   a trivial modification of the proof of  \cite[Theorem 5.4]{FG07}. Now it remains to use i).

 i) Let ${\bf q}$ be a  quiver.  Let ${\mathscr A}$ be the cluster $K_2$-variety, ${\mathscr X}$ the cluster Poisson variety,   ${\mathscr X}^\vee$  the Langlands  dual one, and 
  ${\cal D}$ the symplecitc double  assigned to ${\bf q}$. 
Let $\widetilde {\bf q}$ be a quiver containing  ${\bf q}$. It gives rise to a cluster $K_2$-variety   $\widetilde {\mathscr A}$ and cluster Poisson variety $\widetilde {\mathscr X}$. 
Recall the notation $V_{\bf q}$ 
for the set of vertices of  ${\bf q}$. 
  A  quasi-cluster transformation ${\bf c}$ induces a classical quasi-cluster transformation ${\bf c}^\ast$, where 
$\ast = {\mathscr A}, {\mathscr X}, {\mathscr X}^\vee, {\cal D}$ etc., 
and its quantum counterparts   
${\bf c}_q^{\mathscr X}, {\bf c}_q^{\cal D}$ etc.

 There  is a canonical projection $r: \widetilde {\mathscr X}\lra {\mathscr X}$, such that $r^*X_i = X_i$ for any cluster Poisson coordinate $X_i$. Therefore
   there is a canonical map of cluster varieties, given as a composition: 
\be \la{TCT}
\widetilde  {\mathscr A} \stackrel{p}{\lra} \widetilde {\mathscr X} \stackrel{r}{\lra} {\mathscr X} , \qquad X_i \lms \prod_{j \in V_{\widetilde {\bf q}}}A_j^{\varepsilon_{ij}}, 
~\qquad~i \in  V_{\bf q}.
\ee

\bp  [{\cite[Proposition 5.21]{FG07}}] \label{P5.21} Let $\widetilde {\bf q}$ be a 
quiver containing a quiver ${\bf q}$. Let 
 ${\bf c}: {\bf q} \to{\bf q}$ be a quiver cluster transformation. It induces a 
quiver cluster  transformation 
$\widetilde {\bf c}: \widetilde{\bf q} \to \widetilde{\bf q}$. 

Assume that the  cluster transformation  $\widetilde {\bf c}^{{\mathscr A}}$  is   trivial,   i.e. $( \widetilde {\bf c}^{{\mathscr A}})^*A_v = A_v$ for any cluster coordinate $A_v$ on $\widetilde {\mathscr A}$, 
and the  map (\ref{TCT}) is surjective. Then   
  the quantum cluster transformations ${\bf c}^{\mathscr X}_q$  and 
${\bf c}^{\cal D}_q$ of the $q$-deformed  
spaces ${\mathscr X}_{q}$ and ${\cal D}_{q}$ 
 are trivial. 
\ep

So  we need   a criteria of triviality of the cluster transformation $\widetilde {\bf c}^{\mathscr A}$ for sufficiently large $\widetilde {\mathscr A}$. 

It is sufficient to consider the double $\widehat {\bf q}$ of the quiver ${\bf q}$ with $V_{\widehat {\bf q}}:= V_{\bf q} \cup V'_{\bf q}$, where 
$V'_{\bf q}$ is a copy of the set $V_{\bf q}$,    the quiver $\widehat {\bf q}$ has no arrows between the vertices of $V'_{\bf q}$, and 
the only arrows between $V_{\bf q}$ and $V'_{\bf q}$ are single arrows from $v$ to $v'$ for each $v \in  V_{\bf q}$.

Let $v \in {\rm V}_{\bf q}$. 
Denote by $l_{\bf q}^v$ a point of ${\mathscr X}^\vee(\Z^t)$ which in the tropical cluster Poisson coordinate system related to the quiver ${\bf q}$  has just one non-zero coordinate - 
the one corresponding to the vertex $v \in {\rm V}_{\bf q}$, which is equal to $1$. 

\bp \la{P10.18} Let ${\bf c}: {\bf q}\to {\bf q}$ be a quiver cluster transformation such that for the  tropicalised cluster Poisson transformation  ${\bf c}^{{\mathscr X}^\vee}_t$ of the Langlands dual cluster Poisson variety ${\mathscr X}^\vee$ we have 
\be \la{TCPT}
{\bf c}^{{\mathscr X}^\vee}_t(l_{\bf q}^v) = l_{\bf q}^v, ~\qquad~\forall v \in {\rm V}_{\bf q}.
\ee
Then for any quiver $\widetilde {\bf q}$ containing  ${\bf q}$,  the  induced cluster transformation  $\widetilde {\bf c}^{{\mathscr A}}$  is trivial. 
\ep

\begin{proof} According to a formal variant of the Duality Conjectures \cite{FG03b}, proved in \cite{GHKK}, there is  a natural map to the completion $\widehat {\cal O}({\mathscr X})$ of the space ${\cal O}({\mathscr X})$:
\be 
\begin{split}
&{\rm I}_{\mathscr X}: {\mathscr X}^\vee(\Z^t) \lra \widehat {\cal O}({\mathscr X}),\\
&{\rm I}_{\mathscr X}(l_{\bf q}^v) = A_v, \qquad \forall {\bf q}, \forall v \in {\rm V}_{\bf q}.\\
\end{split}
\ee
Therefore (\ref{TCPT}) implies that $(\widetilde {\bf c}^{\mathscr A})^*A_v=A_v$ for all $v \in {\rm V}_{\bf q}$. Since  frozen $A$-variables do not mutate, 
 $(\widetilde {\bf c}^{\mathscr A})^*A_w=A_w$ for all $w \in {\rm V}_{\widetilde {\bf q}} - {\rm V}_{\bf q}$. 
\end{proof}

Given a puncture on $\bS$, the group  $W$ acts   geometrically  on the   spaces ${\mathscr A}_{\G, \bS}$, ${\mathscr P}_{\G, \bS}$, and ${\mathscr P}_{\G^\vee, \bS}$. 
These actions are cluster:  we defined a cluster transformation ${\bf c}^*(s_i)$ corresponding to each generator $s_i$ of $W$.  
Note that the Weyl groups of $\G$ and $\G^\vee$ coincide, and 
there is  an isomorphism of cluster Poisson varieties ${\mathscr P}_{\G^\vee, \bS} = {\mathscr P}^\vee_{\G, \bS}$. 
To prove that  the quantum cluster transformations ${\bf c}_q^{\mathscr X}(s_i)$ and ${\bf c}_q^{\cal D}(s_i)$  
 give rise to an action of   $W$ we must show that any product of the generators   $s_{i_1} s_{i_2} ... s_{i_n}$  
representing the unit in $W$ acts by the trivial quantum cluster transformation. By Proposition \ref{P10.18}, it is sufficient to prove that 
the tropicalised cluster Poisson transformation ${\bf c}^{{\mathscr X}^\vee_t}(s_{i_1} ) ... {\bf c}^{{\mathscr X}^\vee_t}(s_{i_n} ) $ is trivial. Since  it describes the tropicalization of  
the geometric action of  $W$ on the  space  ${\mathscr P}_{\G^\vee, \bS}$, the claim follows.  For the   group $\Gamma_\bS$,  we  start with the action of the modular 
groupoid ${\rm M}_\bS$ by cluster transformations, and the proof goes the same way. 

The arguments for the braid group ${\mathbb B}_{\mathfrak g}^\pi$ are literally the same, with the only proviso that we have to add monomial transformations of the frozen variables, which causes no problem. 
  \end{proof}


\medskip\section{Donaldson Thomas Transformation}
\label{DT=cluster}
 
 \medskip
 
 In Section \ref{DT=cluster} we assume first that all boundary intervals of $\bS$ are colored. The passage to the general case is trivial, just as before. \vskip 1mm

In this Section, we investigate the cluster and geometric nature of the Donaldson-Thomas transformations on the space $\mathscr{P}_{\G, \bS}$. We will develop a diagrammatic notation to illustrate the tropical points of   $\mathscr{P}_{\G, \bS}$ in Section \ref{diagramm.notations}.  We present a simple proof of  Theorem \ref{CDT} based on the diagrammatic notations. 

\subsection{Notations}
\label{diagramm.notations}
\medskip

Recall that each underlying quiver of $\mathscr{P}_{\G, t}$ has $3r$ frozen vertices. Every oriented side of $t$ corresponds to $r$ frozen vertices parametrized by the fundamental weights $\Lambda_1, \ldots, \Lambda_r$. When we change a side's orientation, its labeled weights change from $\Lambda$ to $-w_0(\Lambda)$.

Let us pick an angle $v$ of $t$ and a reduced word ${\bf i}=(i_1, \ldots, i_m)$ of $w_0$. It gives rise to  a quiver  via amalgamation. As in Figure \ref{2017.5.9.6.24h*}, the triangle $t$ is decomposed into $m$ many small triangles, whose sides correspond to   $m+1$ zig-zag paths denoted by $P_0, P_1, \ldots, P_m$ from $v$ to its opposite side. Each $P_j$ passes through $r$ many vertices parametrized by  $\Lambda_1, \ldots, \Lambda_r$. In particular, the paths $P_0$ and $P_m$ are assigned to two sides of $t$.

\begin{example}
The following quiver of $\mathscr{P}_{{\rm PGL}_4, t}$  corresponds to the word
$(1,2,3,1,2,1)$. 
We omit the frozen vertices at the bottom  side since they are irrelevant. There are seven zigzag paths illustrated by the red curves, each passing through three vertices labeled by $\Lambda_1, \Lambda_2, \Lambda_3$.
\[
\begin{tikzpicture}
\node (a1) at (0,0) {$\circ$};
\node (a2) at (1,0) {$\bullet$};
\node (a3) at (2,0) {$\bullet$};
\node (a4) at (3,0) {$\circ$};
\node (b1) at (.5, 1) {$\circ$};
\node (b2) at (1.5, 1) {$\bullet$};
\node (b3) at (2.5, 1) {$\circ$};
\node (c1) at (1, 2) {$\circ$};
\node (c2) at (2,2) {$\circ$};
\foreach \from/\to in {a2/a1, a3/a2, a4/a3, b2/b1, b3/b2, c2/c1, b1/a2, a2/b2, b2/a3, a3/b3, c1/b2, b2/c2}
\draw[directed] (\from) -- (\to);
\foreach \from/\to in {a1/b1, b1/c1, c2/b3, b3/a4}
\draw[dashed, directed] (\from) -- (\to);
\draw[red, -latex] (0.8,2) -- (-.2, 0);
\draw[red, -latex] (0.9,2) -- (.4, 1) -- (0.9,0);
\draw[red, -latex] (0.9,2) -- (1.4, 1) --(0.9,0);
\draw[red, -latex] (2.1,2) -- (1.1,0);
\draw[red, -latex] (2.2,2) -- (1.7,1) --(2.2,0);
\draw[red, -latex] (2.2,2) -- (2.7,1) --(2.2,0);
\draw[red, -latex] (2.3,2) -- (3.3,0);
\end{tikzpicture}
\]
These zigzag paths have been illustrated in Figure 28 of \cite[\S 3.1]{GS13}.
\end{example}

Now consider the space $\mathscr{P}_{\G, \bS}$. Let us fix a decorated ideal triangulation $\widehat{\mathcal{T}}:=(\mathcal{T},  \{v_t\}, \{{\bf i}_t\})$, where $\mathcal{T}$ is an ideal triangulation of $\bS$, $v_t$ is a vertex of each triangle $t$ of $\mathcal{T}$, and ${\bf i}_t$ is a reduced word of $w_0$ assigned to $t$. Correspondingly, we obtain a quiver $Q_{\widehat{\mathcal{T}}}$.

Let $\gamma=n_1\Lambda_1+\ldots + n_r \Lambda_r$. Pick a triangle $t$ of $\mathcal{T}$. For $k=0,\ldots, m$, there is a tropical point $l(t, {\bf i}, \gamma, k)$ in $\mathscr{P}_{\G, \bS}(\mathbb{Z}^t)$, whose coordinates associated to the vertices in $P_k$ are $(n_1, \ldots, n_r)$, and the rest are zero.
We shall illustrate the point $l(t, {\bf i}, \gamma, k)$ as the left graph in the following Figure in (O1). It is easy to see that the definition of $l(t, {\bf i}, \gamma, k)$ only depends on the choice of $t$ and ${\bf i}$. Changing the rest triangles of $\mathcal{T}$ via flip of diagonals will not change the coordinates of $l(t, {\bf i}, \gamma, k)$.

Similarly, if the paths are placed on the sides of $t$, then  we can illustrate its corresponding tropical points as the right graph in the figure in (O1). These tropical points do not depend on the choice of the angles $v$ and reduced words ${\bf i}$.
\medskip
\subsection{Three basic operations}
\medskip

Below $\gamma$ is a dominant weight. Let ${\bf i}=(i_1, \ldots,i_m)$ be a reduced word of $w_0$ and let
\[
w=s_{i_1} \ldots s_{i_k}, \hskip 7mm w'= s_{i_{k+1}^\ast}\cdots s_{i_m^\ast}.
\]
Note that $w'w=w_0$. As in \eqref{triple.dominant.weights}, we  set
\[
w(\gamma) =\gamma_+ - \gamma_-
\]
There are three basic operations summarized by the following figures. Note that we only draw the local pictures where the tropical points have been changed.

\begin{itemize}
\item[(O1)] The action of $w$ at a puncture $p$. 
\[
\begin{tikzpicture}[scale=2]
\draw (0,0) -- (-60:1) -- (-120:1) -- (0,0) ;
\node at (-0.02,-.15) {${\bf i}$}; 
\draw[red, -latex] (0,0) -- (-75: 0.9);
\node[red] at (0,-.45)  {${\tiny \gamma}$};
\node[label=below: $k$] at (-75:.9) {};
\node[label=above: $p$] at (0,0) {};
\draw[thick, -stealth] (1.2, -.4) -- (1.8, -.4);
\node at (1.5, -.3) {$w ~@ ~p$};
\begin{scope}[xshift=3cm]
\draw (0,0) -- (-60:1) ;
\node[label=above: $p$] at (0,0) {};
\node[red] at (-.5, -.45) {$\tiny w (\gamma)$};
\node[red] at (0, -1) {$\tiny  \gamma_{-}$};
\draw[red, -latex] (0,0) -- (-120: 1);
\draw[red, -latex] (-60:1) -- (-120: 1);
\end{scope}
\end{tikzpicture}
\]
The operation of $w'$ at a puncture $p$.
\[
\begin{tikzpicture}[scale=2]
\draw (0,0) -- (-120:1) ;
\node[label=above: $p$] at (0,0) {};
\node[red] at (.6, -.45) {$\tiny w(\gamma)$};
\draw[red, -latex] (0,0) -- (-60: 1);
\draw (-60:1) -- (-120: 1);
\draw[thick, -stealth] (1.2, -.4) -- (1.8, -.4);
\node at (1.5, -.3) {$w' ~@ ~p$};
\begin{scope}[xshift=3cm]
\draw (0,0) -- (-60:1);
\draw (-120:1) -- (0,0) ;
\draw[red, -latex] (-120:1) -- (-60: 1);
\draw[red, -latex] (0,0) -- (-75: 0.9);
\node at (-.02,-.13) {${\bf i}^\ast$}; 
\node[red] at (0, -1) {$\tiny  \gamma_{+}$};
\node[red] at (0,-.45)  {${\tiny w_0(\gamma})$};
\node[label=below: $k$] at (-75:.9) {};
\node[label=above: $p$] at (0,0) {};
\end{scope}
\end{tikzpicture}
\]

\item[(O2)] A flip of diagonal.
\[
\begin{tikzpicture}[scale=2]
\draw (0,0) -- (1,0) -- (1,1) -- (0,1) -- (0,0) ;
\draw[red, -latex] (0,1) -- (1,0);
\node[red] at (0.6,0.6)  {${\tiny \gamma}$};
\draw[thick, stealth-stealth] (1.7, .5) -- (2.3, 0.5);
\node at (2, .6) {Flip};
\begin{scope}[xshift=3cm]
\draw (0,0) -- (0,1); 
\draw (1,1) -- (1,0);
\draw[red, -latex] (1,1) -- (0,0);
\draw[red, latex-] (1,1) -- (0,1);
\draw[red, latex-] (1,0) -- (0,0);
\node[red] at (0.5, -.1) {$\tiny  \gamma$};
\node[red] at (0.5, .6) {$\tiny  w_0(\gamma)$};
\node[red] at (0.5, 1.1) {$\tiny  \gamma$};
\node at (1.5, 0.5) {$=$};
\end{scope}
\begin{scope}[xshift=5cm]
\draw (0,0) -- (0,1); 
\draw (1,1) -- (1,0);
\draw[red, latex-] (1,1) -- (0,0);
\draw[red, latex-] (1,1) -- (0,1);
\draw[red, latex-] (1,0) -- (0,0);
\node[red] at (0.5, -.1) {$\tiny  \gamma$};
\node[red] at (0.5, .6) {$\tiny  -\gamma$};
\node[red] at (0.5, 1.1) {$\tiny  \gamma$};
\end{scope}
\end{tikzpicture}
\]

\item[(O3)] The $\ast$-involution.
\[
\begin{tikzpicture}[scale=2]
\draw (0,0) -- (-60:1) -- (-120:1) -- (0,0) ;
\draw[red, -latex] (0,0) -- (-75: 0.9);
\node at (-.02,-.13) {${\bf i}^\ast$}; 
\node[red] at (0,-.45)  {${\tiny w_0(\gamma)}$};
\node[label=below: $k$] at (-75:.9) {};
\node[label=above: $p_1$] at (0,0) {};
\node at (-125:1.1) {$p_2$};
\node at (-55:1.1) {$p_3$};
\draw[thick, stealth-stealth] (-1, -.4) -- (-2, -.4);
\node at (-1.5, -.3) {$\ast$-involution};
\begin{scope}[xshift=-3cm]
\draw (0,0) -- (-60:1) -- (-120:1) -- (0,0) ;
\draw[red, -latex] (0,0) -- (-75: 0.9);
\node at (-.02,-.12) {${\bf i}$}; 
\node[red] at (0,-.45)  {${\tiny -\gamma}$};
\node[label=below: $k$] at (-75:.9) {};
\node[label=above: $p_1$] at (0,0) {};
\node at (-125:1.1) {$p_2$};
\node at (-55:1.1) {$p_3$};
\end{scope}
\end{tikzpicture}
\]
\end{itemize}

Here the proof of (O3) is clear. 

Now we prove (O1). Let $i_k=i$. Let $\omega$ be a weight such that 
$\langle \alpha_i^\vee,  \omega \rangle \geq 0$.  As in \eqref{reduced.word.positive.roots}, the word ${\bf i}$ gives rise to a sequence  $\{\beta_k^{\bf i}\}$ of positive coroots. We  consider the action of $s_i$ on the puncture $p$. 

If $\beta_k^{\bf i}$ is not simple, using Theorem 7.7 and Lemma 10.10 of \cite{GS16},
the tropical point $l(t, {\bf i}, \alpha, k)$ under the action of $s_i$ is 
\be
\la{2022.5.12.si.trop}
s_i\cdot l(t, {\bf i}, \omega, k) = l\left(t, {\bf i}, s_i(\omega), k-1\right).
\ee
The action is illustrated by the following graph.
\[
\begin{tikzpicture}[scale=2.5]
\draw (0,0) -- (-60:1) -- (-120:1) -- (0,0) ;
\node at (-0.02,-.15) {${\bf i}$}; 
\draw[red, -latex] (0,0) -- (-75: 0.9);
\node[red] at (0,-.45)  {${\tiny \omega}$};
\node[label=below: $k$] at (-75:.9) {};
\node[label=above: $p$] at (0,0) {};
\draw[thick, -stealth] (1.2, -.4) -- (1.8, -.4);
\node at (1.5, -.3) {$s_i ~@ ~p$};
\begin{scope}[xshift=3cm]
\draw (0,0) -- (-60:1) ;
\draw (0,0) -- (-60:1) -- (-120:1) -- (0,0) ;
\node at (-0.02,-.15) {${\bf i}$}; 
\draw[red, -latex] (0,0) -- (-85: 0.88);
\node[red] at (-.1,-.45)  {${\tiny s_i(\omega)}$};
\node[label=below: $k-1$] at (-85:.88) {};
\node[label=above: $p$] at (0,0) {};
\end{scope}
\end{tikzpicture}
\]

If $\beta_k^{\bf i}=\alpha_j^\vee$ is simple, then there is an extra frozen vertex on the bottom side of the triangle. Besides the change of tropical points illustrated by the above graph, the value of the extra frozen vertex becomes
\[
\langle \alpha_{i}^\vee, \omega \rangle =\langle w(\alpha_{i}^\vee), w(\omega)\rangle = \langle \alpha_{j^\ast}^\vee, w(\omega) \rangle.
\]

Let $\gamma$ be a dominant weight as in (O1). Let $s_i, u \in W$ be such that $l(s_iu)=l(u)+1$. Note that $u^{-1}(\alpha_i^\vee)$ is a positive coroot. It follows that
\[
\langle \alpha_i^\vee, u(\gamma) \rangle  = \langle u^{-1}(\alpha_i^\vee), \gamma \rangle\geq 0. 
\]
Now we may apply the action of $w=s_{i_1}\ldots s_{i_k}$  on the puncture $p$ recursively, obtaining
\[
w \cdot l(t, {\bf i}, \gamma, k) = l(t, {\bf i}, w(\gamma), 0).
\]
By Lemma \ref{decomp.alpha.plus.minus}, the extra frozen variables on the bottom side becomes $\gamma_-$. Therefore, we prove the first operation of (O1). The proof for the second operation of (O2) goes along the same line. 

\medskip 

Now we prove (O2). As illustrated by Figure \ref{pin18}, the flip of a diagonal of the quadrilateral can be factored into a sequence of elementary flips as in Figure \ref{coord.decom.17.15.p1h}.  We shall use a reduced word ${\bf i}=(i_1, \ldots, i_m)$ for the top side of the quadrilateral, and the word ${\bf i}'=(i_m^\ast, i_{m-1}^\ast, \ldots, i_1^\ast)$ for the bottom side. 
Let $\omega$ be a weight such that $\langle \alpha_i^\vee, \omega\rangle \geq 0$. By a direct calculation, under an elementary flip, the tropical point becomes 
\[
\begin{tikzpicture}[scale=2.5]
\draw (0.6,0) -- (1,0) -- (1,1) -- (0.6,1) -- (0.6,0) ;
\draw[red, -latex] (0.6,1) -- (1,0);
\node[red] at (0.7,0.5)  {${\tiny \omega}$};
\draw[thick, -stealth] (1.7, .5) -- (2.3, 0.5);
\node at (0.8, -.1) {$\tiny  s_{i^\ast}$};
\node at (0.8, 1.1) {$\tiny  s_{i}$};
\node at (2, .6) {Flip};
\begin{scope}[xshift=2.4cm]
\draw (0.6,0) -- (1,0) -- (1,1) -- (0.6,1) -- (0.6,0) ;
\draw[red, latex-] (0.6,0) -- (1,1);
\node at (0.8, -.1) {$\tiny  s_{i^\ast}$};
\node[red] at (0.7, .6) {$\tiny  s_i(\omega)$};
\node at (0.8, 1.1) {$\tiny  s_{i}$};
\end{scope}
\end{tikzpicture}
\]
If there are extra frozen variables associated with the top and the bottom sides, then the tropical value of the corresponding frozen variable becomes $\langle \alpha_i^\vee, \omega\rangle.$ Following the same procedure as in the proof of (O1), we prove (O2). 

\medskip
\subsection{Proof of DT conjecture}
\medskip

Following Section 10 of \cite{GS16},  it suffices to check the following four cases. 

{\bf Case 1: three punctures}. The following procedure is similar to Figure 50 on page 106 of \cite{GS16}
\[
\begin{tikzpicture}[scale=2.2]
\draw (0,0) -- (-60:1) -- (-120:1) -- (0,0) ;
\node at (-.02,-.12) {${\bf i}$}; 
\draw[red, -latex] (0,0) -- (-75: 0.9);
\node[red] at (0,-.45)  {${\tiny \gamma}$};
\node[label=below: $k$] at (-75:.9) {};
\node[label=above: $p_1$] at (0,0) {};
\node at (-125:1.1) {$p_2$};
\node at (-55:1.1) {$p_3$};
\draw[thick, -stealth] (1.2, -.4) -- (1.8, -.4);
\node at (1.5, -.3) {$w~@ ~p_1$};
\begin{scope}[xshift=3cm]
\draw (0,0) -- (-60:1) ;
\node[label=above: $p_1$] at (0,0) {};
\node[red] at (-.6, -.45) {$\tiny w (\gamma)$};
\node[red] at (0, -1) {$\tiny  \gamma_{-}$};
\draw[red, -latex] (0,0) -- (-120: 1);
\draw[red, -latex] (-60:1) -- (-120: 1);
\draw[thick, -stealth] (1.2, -.4) -- (1.8, -.4);
\node at (1.5, -.3) {$w_0 ~@ ~p_3$};
\node at (-125:1.1) {$p_2$};
\node at (-55:1.1) {$p_3$};
\end{scope}

\begin{scope}[yshift=-45]
\draw  (-60:1) -- (-120:1);
\draw[red, -latex] (0,0) -- (-120: 1);
\draw[red, -latex] (-60:1) -- (0: 0);
\node[red] at (-0.5,-.45)  {${\tiny \gamma_{+}}$};
\node[red] at (0.55,-.45)  {${\tiny w_0(\gamma_{-})}$};
\node[label=above: $p_1$] at (0,0) {};
\node at (-125:1.1) {$p_2$};
\node at (-55:1.1) {$p_3$};
\draw[thick, -stealth] (1.2, -.4) -- (1.8, -.4);
\node at (1.5, -.3) {$w_0 ~@ ~p_2$};
\begin{scope}[xshift=3cm]
\draw (0,0) -- (-120:1) ;
\node[label=above: $p_1$] at (0,0) {};
\node[red] at (.6, -.45) {$\tiny w(\gamma)$};
\node[red] at (0, -1) {$\tiny  w_0(\gamma_{+})$};
\draw[red, -latex] (0,0) -- (-60: 1);
\draw[red, -latex] (-60:1) -- (-120: 1);
\draw[thick, -stealth] (1.2, -.4) -- (1.8, -.4);
\node at (1.5, -.3) {$w' ~@ ~p_1$};
\node at (-125:1.1) {$p_2$};
\node at (-55:1.1) {$p_3$};
\end{scope}
\end{scope}

\begin{scope}[yshift=-90]
\draw (0,0) -- (-60:1) -- (-120:1) -- (0,0) ;
\draw[red, -latex] (0,0) -- (-75: 0.9);
\node at (-.02,-.13) {${\bf i}^\ast$}; 
\node[red] at (0,-.45)  {${\tiny w_0(\gamma})$};
\node[label=below: $k$] at (-75:.9) {};
\node[label=above: $p_1$] at (0,0) {};
\node at (-125:1.1) {$p_2$};
\node at (-55:1.1) {$p_3$};
\draw[thick, -stealth] (1.2, -.4) -- (1.8, -.4);
\node at (1.5, -.3) {$\ast$-involution};
\begin{scope}[xshift=3cm]
\draw (0,0) -- (-60:1) -- (-120:1) -- (0,0) ;
\draw[red, -latex] (0,0) -- (-75: 0.9);
\node at (-.02,-.12) {${\bf i}$}; 
\node[red] at (0,-.45)  {${\tiny -\gamma}$};
\node[label=below: $k$] at (-75:.9) {};
\node[label=above: $p_1$] at (0,0) {};
\node at (-125:1.1) {$p_2$};
\node at (-55:1.1) {$p_3$};
\end{scope}
\end{scope}
\end{tikzpicture}
\]
\newpage

{\bf Case 2: two punctures and one marked point.} See Figure 51, page 103 of \cite{GS16}
We ignore the frozen vertices on the side $m'm$ since we  focus on the space $\mathscr{X}_{\G, \bS}$.
\[
\begin{tikzpicture}
\node[minimum size=0pt,inner sep=0pt] (p1) at (1,1) {};
\node[minimum size=0pt,inner sep=0pt] (p2) at (-1,1) {};
\node[minimum size=0pt,inner sep=0pt] (m2) at (-1,-1) {};
\node[minimum size=0pt,inner sep=0pt] (m1) at (1,-1) {};
\node[minimum size=0pt,inner sep=0pt] (k) at (1,-.2) {};
\node at (-.65, .85) {$\bf i$};
\node at (1.2,1.2) {$p_1$};
\node at (-1.2,1.2) {$p_2$};
\node at (-1.2,-1.2) {$m'$};
\node at (1.2,-1.2) {$m$};
\node at (1.2,-.2) {$k$};
\node[red] at (.3,.4) {$\gamma$};
\foreach \from/\to in {p1/p2, p2/m2, m2/m1, m1/p1, p2/m1}
\draw (\from) -- (\to);
\draw[directed, red] (p2) -- (k);
\draw[thick, -stealth] (2.2, 0) -- (3.8, 0);
\node at (3, .3) {$w ~@ ~p_2$};

\begin{scope}[xshift=6cm]
\node[minimum size=0pt,inner sep=0pt] (p1) at (1,1) {};
\node[minimum size=0pt,inner sep=0pt] (p2) at (-1,1) {};
\node[minimum size=0pt,inner sep=0pt] (m2) at (-1,-1) {};
\node[minimum size=0pt,inner sep=0pt] (m1) at (1,-1) {};
\node[minimum size=0pt,inner sep=0pt] (k) at (1,-.2) {};
\node at (1.2,1.2) {$p_1$};
\node at (-1.2,1.2) {$p_2$};
\node at (-1.2,-1.2) {$m'$};
\node at (1.2,-1.2) {$m$};
\node[red] at (1.2,0) {$\gamma_-$};
\node[red] at (.2,.2) {$w(\gamma)$};
\foreach \from/\to in {p1/p2, p2/m2, m2/m1}
\draw (\from) -- (\to);
\foreach \from/\to in {p2/m1, p1/m1}
\draw[directed, red] (\from) -- (\to);
\draw[thick, -stealth] (2.2, 0) -- (3.8, 0);
\node at (3, .3) {$w_0 ~@ ~p_1$};
\end{scope}

\begin{scope}[yshift=-120]
\node[minimum size=0pt,inner sep=0pt] (p1) at (1,1) {};
\node[minimum size=0pt,inner sep=0pt] (p2) at (-1,1) {};
\node[minimum size=0pt,inner sep=0pt] (m2) at (-1,-1) {};
\node[minimum size=0pt,inner sep=0pt] (m1) at (1,-1) {};
\node[minimum size=0pt,inner sep=0pt] (k) at (1,-.2) {};
\node at (1.2,1.2) {$p_1$};
\node at (-1.2,1.2) {$p_2$};
\node at (-1.2,-1.2) {$m'$};
\node at (1.2,-1.2) {$m$};
\node[red] at (0,1.2) {$w_0(\gamma_-)$};
\node[red] at (.2,.2) {$\gamma_+$};
\foreach \from/\to in {p2/m2, m2/m1, m1/p1}
\draw (\from) -- (\to);
\foreach \from/\to in {p2/m1, p1/p2}
\draw[directed, red] (\from) -- (\to);
\draw[thick, -stealth] (2.2, 0) -- (3.8, 0);
\node at (3, .3) {Flip};
\end{scope}

\begin{scope}[yshift=-120, xshift=6cm]
\node[minimum size=0pt,inner sep=0pt] (p1) at (1,1) {};
\node[minimum size=0pt,inner sep=0pt] (p2) at (-1,1) {};
\node[minimum size=0pt,inner sep=0pt] (m2) at (-1,-1) {};
\node[minimum size=0pt,inner sep=0pt] (m1) at (1,-1) {};
\node[minimum size=0pt,inner sep=0pt] (k) at (1,-.2) {};
\node at (1.2,1.2) {$p_1$};
\node at (-1.2,1.2) {$p_2$};
\node at (-1.2,-1.2) {$m'$};
\node at (1.2,-1.2) {$m$};
\node[red] at (0,1.2) {$w(\gamma)$};
\node[red] at (-.2,.2) {$-\gamma_+$};
\foreach \from/\to in {p2/m2, m2/m1, m1/p1}
\draw (\from) -- (\to);
\foreach \from/\to in {m2/p1, p2/p1}
\draw[directed, red] (\from) -- (\to);
\draw[thick, -stealth] (2.2, 0) -- (3.8, 0);
\node at (3, .3) {$w' @ ~p_2$};
\end{scope}

\begin{scope}[yshift=-240]
\node[minimum size=0pt,inner sep=0pt] (p1) at (1,1) {};
\node[minimum size=0pt,inner sep=0pt] (p2) at (-1,1) {};
\node[minimum size=0pt,inner sep=0pt] (m2) at (-1,-1) {};
\node[minimum size=0pt,inner sep=0pt] (m1) at (1,-1) {};
\node[minimum size=0pt,inner sep=0pt] (k) at (-.2,-.2) {};
\node at (-0.8, 0.8) {${\bf i}^*$};
\node at (1.2,1.2) {$p_1$};
\node at (-1.2,1.2) {$p_2$};
\node at (-1.2,-1.2) {$m'$};
\node at (1.2,-1.2) {$m$};
\node at (0,-.2) {$k$};
\node[red] at (.1,.4) {$w_0(\gamma)$};
\foreach \from/\to in {p1/p2, p2/m2, m2/m1, m1/p1, p1/m2}
\draw (\from) -- (\to);
\foreach \from/\to in {p2/k}
\draw[directed, red] (\from) -- (\to);
\draw[thick, -stealth] (2.2, 0) -- (3.8, 0);
\node at (3, .3) {Rotation};
\end{scope}

\begin{scope}[yshift=-240, xshift=6cm]
\node[minimum size=0pt,inner sep=0pt] (p1) at (1,1) {};
\node[minimum size=0pt,inner sep=0pt] (p2) at (-1,1) {};
\node[minimum size=0pt,inner sep=0pt] (m2) at (-1,-1) {};
\node[minimum size=0pt,inner sep=0pt] (m1) at (1,-1) {};
\node[minimum size=0pt,inner sep=0pt] (k) at (1,-.2) {};
\node at (-0.65, 0.85) {${\bf i}^*$};
\node at (1.2,1.2) {$p_1$};
\node at (-1.2,1.2) {$p_2$};
\node at (-1.2,-1.2) {$m'$};
\node at (1.2,-1.2) {$m$};
\node at (1.2,-.2) {$k$};
\node[red] at (.3,.4) {$w_0(\gamma)$};
\foreach \from/\to in {p1/p2, p2/m2, m2/m1, m1/p1, p2/m1}
\draw (\from) -- (\to);
\foreach \from/\to in {p2/k}
\draw[directed, red] (\from) -- (\to);
\draw[thick, -stealth] (2.2, 0) -- (3.8, 0);
\node at (3, .3) {$\ast$-involution};
\end{scope}

\begin{scope}[yshift=-360]
\node[minimum size=0pt,inner sep=0pt] (p1) at (1,1) {};
\node[minimum size=0pt,inner sep=0pt] (p2) at (-1,1) {};
\node[minimum size=0pt,inner sep=0pt] (m2) at (-1,-1) {};
\node[minimum size=0pt,inner sep=0pt] (m1) at (1,-1) {};
\node[minimum size=0pt,inner sep=0pt] (k) at (1,-.2) {};
\node at (-0.65, 0.85) {${\bf i}$};
\node at (1.2,1.2) {$p_1$};
\node at (-1.2,1.2) {$p_2$};
\node at (-1.2,-1.2) {$m'$};
\node at (1.2,-1.2) {$m$};
\node at (1.2,-.2) {$k$};
\node[red] at (.3,.4) {$-\gamma$};
\foreach \from/\to in {p1/p2, p2/m2, m2/m1, m1/p1, p2/m1}
\draw (\from) -- (\to);
\foreach \from/\to in {p2/k}
\draw[directed, red] (\from) -- (\to);
\end{scope}
\end{tikzpicture}
\]
\newpage

{\bf Case 3: one puncture and two marked points.}

\[
\begin{tikzpicture}[scale=1.3]
\node[minimum size=0pt,inner sep=0pt] (p) at (90:1) {};
\node[minimum size=0pt,inner sep=0pt] (n2) at (162:1) {};
\node[minimum size=0pt,inner sep=0pt] (m2) at (234:1) {};
\node[minimum size=0pt,inner sep=0pt] (n1) at (306:1) {};
\node[minimum size=0pt,inner sep=0pt] (m1) at (18:1) {};
\node[minimum size=0pt,inner sep=0pt] (k) at (-72: .31) {};
\node[red] at (0.1, 0.3) {$\gamma$};
\node at (85:0.8) {$\bf i$};
\node at (90:1.2) {$p$};
\node at (162:1.2) {$m_2'$};
\node at (234:1.2) {$m_2$};
\node at  (306:1.2) {$m_1'$};
\node at (18:1.2) {$m_1$};
\node at (-72: .51) {$k$};
\foreach \from/\to in {p/n2, n2/m2, m2/n1, n1/m1, m1/p, p/m2, m1/m2}
\draw (\from) -- (\to);
\foreach \from/\to in {p/k}
\draw[directed, red] (\from) -- (\to);
\draw[thick, -stealth] (2.2, 0) -- (3.8, 0);
\node at (3, .3) {$w ~@~ p$};

\begin{scope}[xshift=6cm]
\node[minimum size=0pt,inner sep=0pt] (p) at (90:1) {};
\node[minimum size=0pt,inner sep=0pt] (n2) at (162:1) {};
\node[minimum size=0pt,inner sep=0pt] (m2) at (234:1) {};
\node[minimum size=0pt,inner sep=0pt] (n1) at (306:1) {};
\node[minimum size=0pt,inner sep=0pt] (m1) at (18:1) {};
\node[minimum size=0pt,inner sep=0pt] (k) at (-72: .31) {};
\node[red] at (-0.1, 0.2) {$w(\gamma)$};
\node[red] at (-72:0.4) {$\gamma_-$};
\node at (90:1.2) {$p$};
\node at (162:1.2) {$m_2'$};
\node at (234:1.2) {$m_2$};
\node at  (306:1.2) {$m_1'$};
\node at (18:1.2) {$m_1$};
\foreach \from/\to in {p/n2, n2/m2, m2/n1, n1/m1, m1/p}
\draw (\from) -- (\to);
\foreach \from/\to in {p/m2, m1/m2}
\draw[directed, red] (\from) -- (\to);
\draw[thick, -stealth] (2.2, 0) -- (3.8, 0);
\node at (3, .3) {Flip $@~m_1m_2$};
\end{scope}

\begin{scope}[yshift=-120]
\node[minimum size=0pt,inner sep=0pt] (p) at (90:1) {};
\node[minimum size=0pt,inner sep=0pt] (n2) at (162:1) {};
\node[minimum size=0pt,inner sep=0pt] (m2) at (234:1) {};
\node[minimum size=0pt,inner sep=0pt] (n1) at (306:1) {};
\node[minimum size=0pt,inner sep=0pt] (m1) at (18:1) {};
\node[minimum size=0pt,inner sep=0pt] (k) at (-72: .31) {};
\node[red] at (-0.4, -0.1) {$\gamma_+$};
\node[red] at (.4,-.1) {$-\gamma_-$};
\node at (90:1.2) {$p$};
\node at (162:1.2) {$m_2'$};
\node at (234:1.2) {$m_2$};
\node at  (306:1.2) {$m_1'$};
\node at (18:1.2) {$m_1$};
\foreach \from/\to in {p/n2, n2/m2, m2/n1, n1/m1, m1/p}
\draw (\from) -- (\to);
\foreach \from/\to in {p/m2, p/n1}
\draw[directed, red] (\from) -- (\to);
\draw[thick, -stealth] (2.2, 0) -- (3.8, 0);
\node at (3, .3) {Flip $@~pm_2$};
\end{scope}

\begin{scope}[yshift=-120, xshift=6cm]
\node[minimum size=0pt,inner sep=0pt] (p) at (90:1) {};
\node[minimum size=0pt,inner sep=0pt] (n2) at (162:1) {};
\node[minimum size=0pt,inner sep=0pt] (m2) at (234:1) {};
\node[minimum size=0pt,inner sep=0pt] (n1) at (306:1) {};
\node[minimum size=0pt,inner sep=0pt] (m1) at (18:1) {};
\node[minimum size=0pt,inner sep=0pt] (k) at (-72: .31) {};
\node[red] at (-0.4, -0.1) {$-\gamma_+$};
\node[red] at (.4,-.1) {$w(\gamma)$};
\node at (90:1.2) {$p$};
\node at (162:1.2) {$m_2'$};
\node at (234:1.2) {$m_2$};
\node at  (306:1.2) {$m_1'$};
\node at (18:1.2) {$m_1$};
\foreach \from/\to in {p/n2, n2/m2, m2/n1, n1/m1, m1/p}
\draw (\from) -- (\to);
\foreach \from/\to in {p/n1, n2/n1}
\draw[directed, red] (\from) -- (\to);
\draw[thick, -stealth] (2.2, 0) -- (3.8, 0);
\node at (3, .3) {$w' ~@~ p$};
\end{scope}

\begin{scope}[yshift=-240]
\node[minimum size=0pt,inner sep=0pt] (p) at (90:1) {};
\node[minimum size=0pt,inner sep=0pt] (n2) at (162:1) {};
\node[minimum size=0pt,inner sep=0pt] (m2) at (234:1) {};
\node[minimum size=0pt,inner sep=0pt] (n1) at (306:1) {};
\node[minimum size=0pt,inner sep=0pt] (m1) at (18:1) {};
\node[minimum size=0pt,inner sep=0pt] (k) at (-108: .31) {};
\node at (100:0.8) {$\bf i^\ast$};
\node[red] at (-0.4, 0.2) {$w_0(\gamma)$};
\node at (90:1.2) {$p$};
\node at (162:1.2) {$m_2'$};
\node at (234:1.2) {$m_2$};
\node at  (306:1.2) {$m_1'$};
\node at (18:1.2) {$m_1$};
\foreach \from/\to in {p/n2, n2/m2, m2/n1, n1/m1, m1/p, p/n1, n1/n2}
\draw (\from) -- (\to);
\foreach \from/\to in {p/k}
\draw[directed, red] (\from) -- (\to);
\draw[thick, -stealth] (2.2, 0) -- (3.8, 0);
\node at (3, .3) {$\ast$-involution};
\end{scope}

\begin{scope}[yshift=-240, xshift=6cm]
\node at (100:0.8) {$\bf i$};
\node[minimum size=0pt,inner sep=0pt] (p) at (90:1) {};
\node[minimum size=0pt,inner sep=0pt] (n2) at (162:1) {};
\node[minimum size=0pt,inner sep=0pt] (m2) at (234:1) {};
\node[minimum size=0pt,inner sep=0pt] (n1) at (306:1) {};
\node[minimum size=0pt,inner sep=0pt] (m1) at (18:1) {};
\node[minimum size=0pt,inner sep=0pt] (k) at (-108: .31) {};
\node[red] at (-0.4, 0.2) {$-\gamma$};
\node at (90:1.2) {$p$};
\node at (162:1.2) {$m_2'$};
\node at (234:1.2) {$m_2$};
\node at  (306:1.2) {$m_1'$};
\node at (18:1.2) {$m_1$};
\foreach \from/\to in {p/n2, n2/m2, m2/n1, n1/m1, m1/p, p/n1, n1/n2}
\draw (\from) -- (\to);
\foreach \from/\to in {p/k}
\draw[directed, red] (\from) -- (\to);
\draw[thick, -stealth] (2.2, 0) -- (3.8, 0);
\node at (3, .3) {Rotation};
\end{scope}

\begin{scope}[yshift=-360]
\node at (85:0.8) {$\bf i$};
\node[minimum size=0pt,inner sep=0pt] (p) at (90:1) {};
\node[minimum size=0pt,inner sep=0pt] (n2) at (162:1) {};
\node[minimum size=0pt,inner sep=0pt] (m2) at (234:1) {};
\node[minimum size=0pt,inner sep=0pt] (n1) at (306:1) {};
\node[minimum size=0pt,inner sep=0pt] (m1) at (18:1) {};
\node[minimum size=0pt,inner sep=0pt] (k) at (-72: .31) {};
\node[red] at (0.2, 0.3) {$-\gamma$};
\node at (90:1.2) {$p$};
\node at (162:1.2) {$m_2'$};
\node at (234:1.2) {$m_2$};
\node at  (306:1.2) {$m_1'$};
\node at (18:1.2) {$m_1$};
\node at (-72: .51) {$k$};
\foreach \from/\to in {p/n2, n2/m2, m2/n1, n1/m1, m1/p, p/m2, m1/m2}
\draw (\from) -- (\to);
\foreach \from/\to in {p/k}
\draw[directed, red] (\from) -- (\to);
\end{scope}
\end{tikzpicture}
\]

{\bf Case 4: three marked points.} It reduces to the case when $\bS=D_6$ is a hexagon:
\[
\begin{tikzpicture}[scale=1.3]
\draw (30:1) --(90:1) --(150:1)--(210:1)--(270:1)--(330:1)--(30:1);
\draw (90:1) --(210:1) --(330:1)--(90:1);
\draw[red, -latex] (90:1)--(0.2,-.5);
\node[red] at (0,0) {$\gamma$};
\node at (90:1.2) {$m_1$};
\node at (30:1.2) {$m_1'$};
\node at (330:1.2) {$m_2$};
\node at (270:1.2) {$m_2'$};
\node at (210:1.2) {$m_3$};
\node at (150:1.2) {$m_3'$};
\node at (0,.8) {${\bf i}$};
\draw[thick, -stealth] (2,0)--(3,0);
\node at (2.5, 0.2) {$\mathcal{DT}^t_{\G, D_6}$};
\begin{scope}[xshift=5cm]
\draw (30:1) --(90:1) --(150:1)--(210:1)--(270:1)--(330:1)--(30:1);
\draw (90:1) --(210:1) --(330:1)--(90:1);
\draw[red, -latex] (90:1)--(0.2,-.5);
\node[red] at (0,0) {$-\gamma$};
\node at (90:1.2) {$m_1$};
\node at (30:1.2) {$m_1'$};
\node at (330:1.2) {$m_2$};
\node at (270:1.2) {$m_2'$};
\node at (210:1.2) {$m_3$};
\node at (150:1.2) {$m_3'$};
\node at (0,.8) {${\bf i}$};
\end{scope}
\end{tikzpicture}
\]
It is a special case of the DT transformations for double Bott-Samelson cells,  proven by  the second author and D. Weng in \cite[Theorem 4.8]{SW}.


 
\medskip

\section{Quantum groups via quantum moduli spaces} \la{SSEECC11.3}

\medskip

\subsection{Basics of quantum groups} \la{SSECC11.1}
\medskip

Let $\U_q(\g)$ be the quantum group of finite type  with generators $\{\bE_i, \bF_i, \bK_i^{\pm}\}_{i\in {\rm I}}$, 
satisfying  relations  (\ref{Re1})-(\ref{Re2}) below. 
The generators  are rescaled version of the ones $\{E_i, F_i, K_i\} $ used in   \cite{L}:
\be \la{RESC}
\bE_i= q_i^{-\frac{1}{2}}(q_i-q_i^{-1}) E_i, \hskip 7mm \bF_i= q_i^{\frac{1}{2}}(q_i^{-1}-q_i) F_i, \hskip 7mm \bK_i =K_i;\hskip 10mm \mbox{where $q_i=q^{1/d_i}$.}
\ee
Therefore we have the modified relation
\be \la{Re1}
\bE_i\bF_j - \bF_j\bE_i = \delta_{ij}(q_i^{-1} - q_i)(\bK_i - \bK^{-1}_i).
\ee
Let us introduce the following notations
\[
[n]_q=\frac{q^n-q^{-n}}{q-q^{-1}}, \hskip 7mm [n]_q^{!}=\prod_{s=1}^n [s]_q, \hskip 7mm  {n\brack s}_q=\frac{[n]_q^!}{[s]_q^![n-s]_q^!}.
\] 
The other relations are:
\be \la{Re2}
\begin{split}
&\bK_i \bK_j = \bK_j \bK_i, ~\qquad~ \bK_i\bK_i^{-1}=1.\\
&\bK_i \bE_j = q^{{\rm C}_{ij}}\bE_j \bK_i,~\qquad~\bK_i \bF_j = q^{-{\rm C}_{ij}}\bF_j \bK_i,\\
&\sum_{s=0}^{1-{\rm C}_{ij}} (-1)^s {1-c \brack s}_{q_{i}} {\bf E}_i^s {\bf E}_j {\bf E}_i^{1-{\rm C}_{ij}-s}=0,\\
&\sum_{s=0}^{1-{\rm C}_{ij}} (-1)^s {1-c \brack s}_{q_{i}} {\bf F}_i^s {\bf F}_j {\bf F}_i^{1-{\rm C}_{ij}-s}=0.\\
\end{split}
\ee
The first two relations in (\ref{Re2}) tell that the elements $\bK_i$ give rise to  an action of the Cartan group $\H$ on $\U_q(\g)$. 
The last two relations are called the {\it quantum Serre relations.}

\subsubsection{The antipode $S$ and the counit $\varepsilon$.} They are given by
 \be \la{ACU}
 \begin{split}
 &S(\bK_i) = \bK^{-1}_i, ~\qquad~S(\bE_i) =  - \bK^{-1}_i\bE_i, ~\qquad~S(\bF_i) =  - \bF_i \bK_i.\\
 &\varepsilon(\bK_i) = 1, ~\qquad~\varepsilon(\bE_i) = 0, ~\qquad~\varepsilon(\bF_i) = 0.\\
 \end{split}
 \ee
 
\subsubsection{The coproduct  on $\U_q(\mathfrak{g})$.} It is given by 
\be \la{383+}  \begin{split}
 &\Delta ({\bf K}_i) = {\bf K}_i \otimes  {\bf K}_i,\\
  &\Delta({\bf E}_i) = {\bf E}_i \otimes 1 + {\bf K}_i \otimes {\bf E}_i,\\
   &\Delta ({\bf F}_i) = 1 \otimes {\bf F}_i + {\bf F}_i  \otimes {\bf K}^{-1}_i.\\
    \end{split}
  \ee
  
  \subsubsection{The $\ast$-algebra structure on $\U_q(\mathfrak{g})$.} Consider the following anti-involution $\ast$ of $\U_q(\mathfrak{g})$:
\be \la{ASST}
\ast: \U_q(\g) \lra \U_q(\g)^{\rm op}, \hskip 9mm \bE_i \lms \bE_i, \hskip 6mm \bF_i \lms \bF_i, \hskip 6mm \bK_i\lms \bK_i, \hskip 6mm q\lms q^{-1}.
\ee
 Thanks to the rescaling (\ref{RESC}) of the generators  of $\U_q(\mathfrak{g})$, it preserves the relations.

\subsubsection{Lusztig's braid group action  on  $\U_q(\g)$.} For simplicity, let us assume that  $\g$ is simply-laced. 
The braid group generators $s_i$ act by the automorphisms $T_i$ $(i\in {\rm I})$ of the algebra $\U_q(\g)$,  
given  by 
\[
\begin{split}
&\bE_j \mapsto  q^{-1}\bK_j^{-1}\bF_j,    \hskip 14mm\qquad ~\qquad~ \bF_j \mapsto q\bE_j\bK_j ,  \hskip 21mm\qquad ~\qquad~ \bK_j \mapsto \bK_j^{-1}, \hskip 9mm   \mbox{if } j=i; \\
\\
&\bE_j \mapsto  \bE_j,   \hskip 14mm\qquad \qquad \qquad \qquad ~ \bF_j \mapsto \bF_j ,  \hskip 21mm\qquad \qquad \qquad\bK_j \mapsto \bK_j, \hskip 12mm  \mbox{if } {\rm C}_{ij}=0; \\
\\
&\bE_j \mapsto  {\displaystyle \frac{q^{1/2}\bE_j \bE_i- q^{-1/2}\bE_i\bE_j}{q-q^{-1}}},  ~\qquad~
  \bF_j \mapsto {\displaystyle \frac{q^{1/2}\bF_j \bF_i- q^{-1/2}\bF_i\bF_j}{q-q^{-1}}} , ~\qquad~  \bK_j \mapsto \bK_i\bK_j, \hskip 6mm  \mbox{if } {\rm C}_{ij}=-1. \\
\end{split}
\]
The maps $T_i$ commute with the anti-involution $\ast$. 
 They satisfy  the  braid relations \cite{L}: 
\begin{itemize}
\item $T_i T_j T_i= T_j T_i T_j$ if ${\rm C}_{ij}=-1$, and $T_i T_j= T_j T_i$ otherwise. 
\end{itemize}
One assigns to  each Weyl group element $w$  a map $$
T_w= T_{i_1}T_{i_2}\ldots T_{i_k}, ~\qquad~\mbox{where $w=s_{i_1}s_{i_2}\ldots s_{i_k}$ is a reduced decomposition}.
$$

\begin{proposition}[{\cite[Proposition 1.8]{L}}] Let $w\in W$ such that $l(ws_i)=l(w)+1$. Then 
\[
T_w(\bE_i)\in \U_q(\mathfrak{n}_+), \qquad  T_w(\bF_i)\in \U_q(\mathfrak{n}_-).
\]
If, in addition,  $w(\alpha_i)=\alpha_k$ is a simple positive root, then $T_w(\bE_i)=\bE_k$, $T_w(\bF_i)=\bF_k$.
\end{proposition}

Let ${\bf i}=(i_1, \ldots, i_N)$ be a reduced word for the longest element $w_0\in W$. 
The following elements were used crucially by Lusztig in the construction of a PBW basis
\be \la{ELLWin}
\bE_{{\bf i}, k} := T_{s_{i_1}\ldots s_{i_{k-1}}}(\bE_{i_k}) \in \U_q(\mathfrak{n}_+).
\ee

\vskip 2mm
 The quasi-classical $q\to 1$ limits of these elements, divided by $q-q^{-1}$,  give the Chevalley basis $\{E_\alpha\}$ of the Lie algebra $\n_+$. It is parametrized by the positive simple roots, 
 ordered by the reduced decomposition ${\bf i}$ of $w_0$. In the quantum case the set of elements $\{\bE_{{\bf i}, k}\}$  depends 
on  ${\bf i}$. 

The compatibility of the geometric braid group action with Lusztig's led  to  
a geometric realization of the elements $\bE_{{\bf i},k}$ given in Theorem  \ref{REALI}. 
On the other hand, establishing the geometric interpretation of the elements $\bE_{{\bf i},k}$ is sufficient 
to prove that the geometric braid group action is compatible with Lusztig's. 
 
\medskip

\subsection{Mapping ${\U}_q(\mathfrak{b})$ to functions on quantum cluster varieties  related to a braid semigroup} \la{SECT16.2} 
\medskip

 In Section \ref{SECT16.2} we start from  a root system $\Delta$ associated with any 
skewsymmetrizable Cartan matrix ${\rm C}_{ij}$, $i,j\in {\rm I}$. 
We denote by  $\mathfrak{g}_\Delta$ 
 the corresponding  Kac-Moody Lie algebra, and by  $\mathfrak{b}_\Delta$ its Borel subalgebra.  
 The Lie algebra $\mathfrak{b}_\Delta$ is infinite dimensional unless $\Delta$ is a finite type. 
We recall
 the braid group  ${\Bbb B}_\Delta$,  and  its positive braid semigroup ${\Bbb B}^+_\Delta$, see Section \ref{SECT14.1.2}.
 
 Let us pick any element $b \in {\Bbb B}_\Delta^+$. Then there is the associated cluster Poisson variety 
 $\mathscr{X}_{[b]}$. 
 Namely, let  ${\bf i}=(i_1, \ldots, i_m)$ be any  word representing  $b$. The construction of the elementary quivers ${\bf J}(i_k)$ from Section \ref{11.2.1} is extended verbatim to the case of an arbitrary root system $\Delta$. Amalgamating the elementary quivers according to the word ${\bf i}$, we get a quiver
\[
{\bf J}({\bf i})= {\bf J}(i_1)\ast \cdots \ast {\bf J}(i_m).
\]
In particular, it provides a lattice   $\Lambda$   
with a skew-symmetric bilinear form $(\ast, \ast)$. The quantum torus algebra $\T_{\bf i}= \mathcal{O}_q({\rm T}_{\bf i})$ is generated by $\X_v, v \in \Lambda$, satisfying $\X_v \X_w= q^{(v,w)} \X_{v+w}$. 
Equivalently, the cluster Poisson torus ${\rm T}_{\bf i}$   is defined by the amalgamation of elementary 
 cluster Poisson varieties $\mathscr{X}_{s_i}$ assigned to the quivers ${\bf J}(s_i)$:
 $$
{ \rm T}_{\bf i} = \mathscr{X}_{s_1}  \ast \ldots \ast \mathscr{X}_{s_m}.
 $$  
Next, as proved in \cite[\S 3.6, 3.7]{FG05}, the cluster seeds corresponding to different word decompositions of $b$ are all related by cluster mutations.  Therefore any elementary braid relation ${\bf i} \sim {\bf i}'$ leads to a birational cluster Poisson transformation between the associated cluster tori: ${\rm T}_{\bf i} \lra {\rm T}_{\bf i'}$. 
 The cluster Poisson variety  $\mathscr{X}_{[b]}$ is obtained by gluing the cluster tori ${\rm T}_{\bf i}$ by these birational transformations. 
  When the root system $\Delta$ is of finite type, it was introduced in \cite{FG05}.  For an arbitrary root system $\Delta$ the non-trivial result, obtained in   \cite[Theorem 2.18]{SW}, 
  is  that 
 if a sequence of braid relations leads to the same reduced decomposition, then the corresponding cluster transformation is trivial. Therefore the gluing process is well defined. 
 The variety $\mathscr{X}_{[b]}$ is isomorphic to decorated double Bott-Samelson cell ${\rm Conf}_{b}^e(\mathcal{A}_{\rm ad})$ in {\it loc.cit.}

 \bt \la{AnyDelta} For any element $b \in {\Bbb B}_\Delta^+$, there is a canonical map of algebras
 \be
 \kappa_b: {\U}_q(\mathfrak{b}) \lra {\cal O}_q(\mathscr{X}_{[b]}).
 \ee
 \et
 \begin{proof} Let us explain first the architecture of the proof. After fixing a word ${\bf i}$ of $b$, we obtain a cluster seed and a quantum torus algebra $\mathbb{T}_{\bf i}$ defined above. Proposition \ref{Prop11.11} below shows that there is homomorphism from $\U_q(\mathfrak{b})$ to $\mathbb{T}_{\bf i}$. The element ${\rm K}_{\alpha}$ is a Laurent monomial for every cluster seed of $\mathscr{X}_{[b]}$, and therefore can be uniquely promoted to the quantized cluster algebra ${\cal O}_q(\mathscr{X}_{[b]})$ by the quantum lift theorem in Section \ref{sec18}. For every $ {\rm W}_{\alpha}$ given by a local quiver as in Figure \ref{local.quiver.alpha}, let us mutate at the vertices $e_{n-1}^\alpha$, ... $e_3^\alpha$, $e_2^\alpha$. In this way, we obtain an optimized quiver for $e_1^\alpha$, in the sense that all the arrows from mutable vertices are coming into $e_1^\alpha$. In this quiver, the quantized cluster variable associated with $e_1^\alpha$ is precisely ${\rm W}_\alpha$, and therefore belongs to ${\cal O}_q(\mathscr{X}_{[b]})$. 
 Let us elaborate  this program. 
 \vskip 2mm

Recall the level map from the set of vertices of ${\bf J}({\bf i})$ to {\rm I}.  The full subquiver formed by vertices at level $\alpha$ is a type $A_n$ quiver in Figure \ref{local.quiver.alpha}.
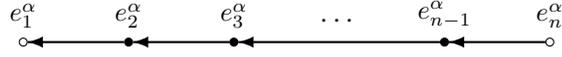
\begin{figure}
\begin{tikzpicture}[scale=1.4]
\node [circle,draw,fill=white,minimum size=3pt,inner sep=0pt, label=above:{\small $e_1^\alpha$}] (a) at (0,0) {};
\node [circle,draw,fill,minimum size=3pt,inner sep=0pt, label=above:{\small $e_2^\alpha$}] (b) at (1,0) {};
 \node [circle,draw,fill, minimum size=3pt,inner sep=0pt, label=above:{\small $e_3^\alpha$}]  (c) at (2,0){};
  \node [circle,draw,fill, minimum size=3pt,inner sep=0pt, label=above:{\small $e_{n-1}^\alpha$}]  (d) at (4,0){};
 \node [circle,draw,fill=white,minimum size=3pt,inner sep=0pt, label=above:{\small $e_n^\alpha$}] (e) at (5,0) {};
 \foreach \from/\to in {b/a, c/b, e/d}
                  \draw[directed, thick] (\from) -- (\to);
\foreach \from/\to in {d/c}                   
              \draw[directed, thick] (\from) -- (\to);      
  \node at (3,0.2) {$\cdots $};                  
     \end{tikzpicture}
     \caption{Local $A_n$ quiver on level $\alpha$.}
       \label{local.quiver.alpha}
    \end{figure}
Here $\{e_1^\alpha, ..., e_{n}^\alpha \}$ are the basis vectors of  $\Lambda$ for  the  level $\alpha$ vertices, counted from the left to the right. Note that $e_1^\alpha$ and $e_n^\alpha$ are frozen. Let us set 
\be
\la{def.347}
\begin{split}
&v_i^\alpha := e_1^\alpha + \cdots + e_i^\alpha.\\
&{\rm W}_\alpha := \sum_{i=1}^{n-1} \X_{v_i^\alpha}, \qquad {\rm K}_{\alpha}:= \X_{v_n^\alpha}.\\
\end{split}
\ee

\bp 
\la{Prop11.11}
There is a well defined map of $\ast$-algebra
\be
\begin{split}
&\U_q(\mathfrak{b}) \lra \mathbb{T}_{\bf i},\\
&{\bf E}_{\alpha} \lms {\rm W}_{\alpha}, \qquad {\bf K}_\alpha \lms {\rm K}_{\alpha}.\\
\end{split}
\ee
\ep

\begin{proof} Let $\alpha, \beta \in {\rm I}$. Recall in Figure \ref{pin30} the subquiver ${\bf Q}_{\alpha \beta}$ of vertices at levels $\alpha$ and $\beta$.

Assume that ${\bf Q}_{\alpha \beta}$ consist of  $n$   vertices at level  $\alpha$ and  $m$  vertices at level $\beta$. Let   $\{v_1^\alpha, \ldots, v_{n}^\alpha\}$ and $\{v_1^\beta, \ldots, v_{m}^\beta\}$
be the collections of vectors  in \eqref{def.347}. They  enjoy the following properties:

\vskip 1mm 1. $(v_i^\alpha, v_j^\alpha)=1$ whenever $ i <j\leq n $;

\vskip 1mm2. $(v_m^\beta, v_n^\alpha)=0$;

\vskip 1mm3. for every $k< m$, there exists an $i_k <n $ such that, setting $c= {C}_{\alpha\beta}$:  
\[
(v_k^\beta, v_i^\alpha)=    \left\{\begin{array}{ll}  
   -c/2 & \mbox{if }  i\leq i_k, \\
     c/2 & \mbox{if } i > i_k. \\
   \end{array}\right.
\]
Set $q_\alpha:= q^{1/d_\alpha}$. It is easy to check that 
\[
{\rm K}_\alpha {\rm W}_\beta = q_\alpha^{c}{\rm W}_\beta {\rm K}_\alpha, \hskip 7mm {\rm K}_\alpha {\rm W}_\alpha = q_\alpha^2  {\rm W}_\alpha {\rm K}_\alpha, \hskip 7mm {\rm K}_\alpha{\rm K}_\beta = {\rm K}_\beta K_{\alpha}.
\]
Precisely, the first identity follows since $(v_k^\beta, v_n^\alpha)=c/2$ for every  $k<m$, the second identity follows from Property 1, and the third identity follows from Property 2.

It remains to prove  the quantum Serre relations:  
\be \la{373}
\sum_{s=0}^{1-c} (-1)^s {1-c \brack s}_{q_{\alpha}} {\rm W}_\alpha^s {\rm W}_\beta {\rm W}_\alpha^{1-c-s}=0.
\ee
By Properties 1 and 3, for every $k<m$, we may write ${\rm W}_\alpha=M+N$ such that 
\be
\begin{split}
&M:= \sum_{  i\leq i_k} \X_{v_i^\alpha}, \hskip 7mm N:= \sum_{i_k< i < n} \X_{v_i^\alpha};\\
&MN=q_\alpha^2NM, \hskip 7mm \X_{v_k^\beta} M =q_\alpha^{-c}M \X_{v_k^\beta}, \hskip 7mm \X_{v_k^\beta} N= q_\alpha^c N\X_{v_k^\beta}.\\
\end{split}
\ee
By Lemma \ref{qsr.lemma}, ${\rm W}_\alpha$, $X_{v_k^\beta}$ satisfy the quantum Serre relation. Taking sum over $k$,  \eqref{373} follows.
\end{proof}

\bl 
\la{qsr.lemma}
Let $n\geq 0$ be an integer. Let $E_1=M+N$ and $E_2=L$ be such that 
\[
MN=q^2NM, \hskip 7mm LM=q^{n}ML, \hskip 7mm LN=q^{-n}NL.
\]
Then $E_1$ and $E_2$ satisfy the quantum Serre relation
\be \la{366}
\sum_{s=0}^{1+n} (-1)^s {1+n \brack s}_q E_1^s E_2 E_1^{1+n-s}=0.
\ee
\el
\begin{proof} Set $T:=MN^{-1}$. We have
$
NT=q^{-2}TN$,  $LT=q^{2n}TL$. 
Note that $E_1 =(1+T)N$.  

Let us set 
$(x;q)_s:= (1+x) (1+q^2x)   \ldots  (1+q^{2s-2}x)$. Then
\[
E_1^s=(M+N)^s= (T; q^{-1})_s N^s.
\]
Hence
\begin{align}
E_1^sE_2E_1^{1+n-s} &= (T; q^{-1})_s ~N^s L ~(T; q^{-1})_{1+n-s}~ N^{1+n-s} \nonumber\\
&= (T; q^{-1})_s ~(q^{2n-2s}T; q^{-1})_{1+n-s} ~N^s L N^{1+n-s} \nonumber\\
&=(T; q^{-1})_s ~(q^{2n-2s}T; q^{-1})_{1+n-s} ~q^{ns} LN^{1+n} \nonumber\\
&=q^{ns}(q^{-2s+2}T;q)_n ~ (1+T) L N^{1+n}. \nonumber
\end{align}
Therefore
\[
{\rm LHS} (\ref{366}) = \left(\sum_{s=0}^{1+n}(-1)^s q^{ns}{1+n \brack s}_q (q^{-2s+2}T;q)_{n}\right)\cdot (1+T)LN^{1+n}.
\]
It remains to show that
\be
\label{2019.1.14.2.25ts}
f_n(T):= \sum_{s=0}^{1+n}(-1)^s q^{ns}{1+n \brack s}_q (q^{-2s+2}T;q)_{n}=0.
\ee
Clearly
$f_0(T)=0$.
Note that the quantum binomial coefficients satisfy the identity
\be
\la{2019.1.14.2.32.ts}
{n+1 \brack s}_q={n \brack s}_q q^{-s} +{n \brack s-1}_q q^{1+n-s}.
\ee
Therefore we get
\begin{align}
q^{ns}{1+n \brack s}_q (q^{-2s+2}T;q)_{n}
=\left(q^{(n-1)s}{n \brack s}_q+ q^{2n}q^{(n-1)(s-1)}{n \brack s-1}_q\right)  (q^{-2s+4}T; q)_{n-1} (1+T)\nonumber\\
=\left(q^{(n-1)s}{n \brack s}_q (q^{-2s+2} q^2T; q)_{n-1}+ q^{2n}q^{(n-1)(s-1)}{n \brack s-1}_q(q^{-2s+4}T; q)_{n-1}\right)   (1+T). \nonumber
\end{align}
Plugging it into \eqref{2019.1.14.2.25ts}, we get
$
f_n(T)= \left(f_{n-1}(q^2T)- q^{2n} f_{n-1}(T)\right) (1+T).
$ 
By induction, $f_{n-1}(T)=0$ implies $f_n(T)=0$.
\end{proof}
Theorem \ref{AnyDelta} is proved. \end{proof}

Let us explore   basic properties of the patterns  $W_\alpha$ and $K_\alpha$ in  \eqref{def.347} for future use.

For a type $A_n$ quiver ${\bf p}$ with $e_1$ and $e_n$ frozen
 
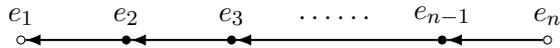
\begin{center}
\begin{tikzpicture}[scale=1.4]
\node [circle,draw,fill=white,minimum size=3pt,inner sep=0pt, label=above:$e_1$] (a) at (0,0) {};
\node [circle,draw,fill,minimum size=3pt,inner sep=0pt, label=above:$e_2$] (b) at (1,0) {};
 \node [circle,draw,fill, minimum size=3pt,inner sep=0pt, label=above:$e_3$]  (c) at (2,0){};
  \node [circle,draw,fill, minimum size=3pt,inner sep=0pt, label=above:$e_{n-1}$]  (d) at (4,0){};
 \node [circle,draw,fill=white,minimum size=3pt,inner sep=0pt, label=above:$e_n$] (e) at (5,0) {};
 \foreach \from/\to in {b/a, c/b, e/d}
                  \draw[directed, thick] (\from) -- (\to);
\foreach \from/\to in {d/c}                   
              \draw[directed, thick] (\from) -- (\to);      
  \node at (3,0.2) {$\cdots \cdots$};                  
     \end{tikzpicture}
    \end{center}
we define
 \be \la{378}
  \begin{split}
  &{\rm W}_{\bf p}:= \X_{e_1}+ \X_{e_1 + e_2}+ \cdots + \X_{e_1 + \ldots + e_{n-1}},\\
  &{\rm K}_{\bf p}:= \X_{e_1 + \ldots + e_{n}}.\\
  \end{split}
  \ee 
  \bl 
  \label{mutation.patterns.}
  After mutation at $2, 3, \ldots, n-1$, the elements ${\rm W}_{\bf p}$ and ${\rm K}_{\bf p}$ become monomials.
  \el
  \begin{proof}
 Let us mutate at $2$, the rest vertices forms a type $A_{n-1}$ quiver ${\bf p}'$:
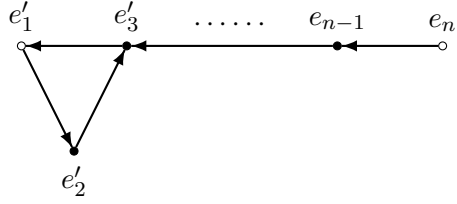
  \begin{center}
\begin{tikzpicture}[scale=1.4]
\node [circle,draw,fill=white,minimum size=3pt,inner sep=0pt, label=above:$e_1'$] (a) at (1,0) {};
\node [circle,draw,fill,minimum size=3pt,inner sep=0pt, label=below:$e_2'$] (b) at (1.5,-1) {};
 \node [circle,draw,fill, minimum size=3pt,inner sep=0pt, label=above:$e_3'$]  (c) at (2,0){};
  \node [circle,draw,fill, minimum size=3pt,inner sep=0pt, label=above:$e_{n-1}$]  (d) at (4,0){};
 \node [circle,draw,fill=white,minimum size=3pt,inner sep=0pt, label=above:$e_n$] (e) at (5,0) {};
 \foreach \from/\to in {a/b, b/c, c/a, e/d}
                  \draw[directed, thick] (\from) -- (\to);
\foreach \from/\to in {d/c}                   
              \draw[directed, thick] (\from) -- (\to);      
  \node at (3,0.2) {$\cdots \cdots$};                  
     \end{tikzpicture}
    \end{center}
    
It is easy to check that
$\X_{e_1'} = \X_{e_1}+ \X_{e_1+e_2}$ and $\X_{e_1'+e_3'} = \X_{e_1+e_2+e_3}$. 

Therefore
$
{\rm W}_{\bf p} = {\rm W}_{{\bf p}'}$ and ${\rm K}_{\bf p}= {\rm K}_{{\bf p}'}.  
$  

Mutating at 3 and so on until $n-1$, we obtain a type $A_2$ subquiver.  

An alternative geometric proof is given in the proof of Lemma \ref{MLe}. 
\end{proof}  
 
 Take a similar quiver ${\bf q}$  with  vertices labeled by a basis $\{f_1, ..., f_m\}$. Denote by ${\bf p}\ast {\bf q}$ the quiver obtained by gluing the right vertex of  ${\bf p}$ 
  to the left one of ${\bf q}$. Equivalently,  it  is the amalgamation of quivers ${\bf q}$ and ${\bf p}$,   described by a basis $\{e_1, ..., e_{n-1}, e_n+f_1, f_2, ..., f_m\}$. The amalgamation gives rise to a $\ast$-algebra embedding:
  \[
  \mathcal{O}_q(\mathscr{X}_{{\bf p}\ast {\bf q}}) \hlra   \mathcal{O}_q(\mathscr{X}_{{\bf p}} \times   \mathscr{X}_{ {\bf q}}).
  \]
Under this embedding, we get
   \be \la{383a}  \begin{split}
  &{\rm W}_{{\bf p}\ast {\bf q}}=  {\rm W}_{\bf p} + {\rm K}_{\bf p}{\rm W}_{\bf q},\\
  &{\rm K}_{{\bf p}\ast {\bf q}}:= {\rm K}_{\bf p}{\rm K}_{\bf q}.\\
  \end{split}
  \ee
\subsubsection{Proof of Theorem \ref{QWEaa}.} \la{16.2.1} The construction of the map $\kappa_b$ follows from Theorem \ref{AnyDelta}. The proof of the part (3) of Theorem \ref{QWEaa} follows immediately from this. 
\medskip

\subsection{Proof of Theorems \ref{QWE}, \ref{UEAB}, and \ref{5.25.24q}}
\la{SEC12.1}
 \medskip

Let us first prove the case when $\bS$ is a triangle $t$. Pick a special point $s$ of $t$. Let us show that  the functions $\mathcal{W}_{s, i}$ and $\mathcal{K}_{s, i}$ are standard monomials in the sense of Definition \ref{quantum.promotion.cls}. Then it follows from Quantum Lift Theorem \ref{quantum.promotion.f} that $\mathcal{W}_{s, i}$ and $\mathcal{K}_{s, i}$ admit natural quantum lifts 
\[
{\rm W}_{s, i} , ~\qquad~{\rm K}_{s, i} \in \mathcal{O}_q(\mathscr{P}_{\G, t}).
\]

By Lemma \ref{L5.10} ii),  $\mathcal{W}_{s, i}$ and $\mathcal{K}_{s, i}$ can be expressed as  monomials in one cluster seed of $\mathscr{P}_{\G', t}$.  It is easy to check that the monomial 
$P_{{\bf i},1} X_{i \choose 0}$ in \eqref{W.dist2} Poisson commutes with all the unfrozen cluster variables, and $P_{{\bf i},1}$   Poisson commutes with 
all the unfrozen cluster variables, except  that
\[
\left\{ \log P_{{\bf i},1}, \log X_{i \choose 1} \right\} =1/ d_i \geq 0. 
\]
By  Proposition \ref{standard monomial}, they are standard monomials. Theorem \ref{QWE} for a triangle $t$ is proved.

\vskip 2mm

Let ${\bf s}$ be a seed of $\mathscr{P}_{\G, t}$ associated to a special point $s$ of $t$ and a reduced word ${\bf i}$ of $w_0$. The full subquiver of ${\bf s}$ at level $\alpha_i$ is of type $A$. Recall the patterns ${\rm W}_{\alpha_i}$ and ${\rm K}_{\alpha_i}$ in \eqref{def.347}. By Lemma \ref{L5.10} i), their specialization at $q=1$ are  $\mathcal{W}_{s, i}$ and $\mathcal{K}_{s, i}$. By Lemma \ref{mutation.patterns.}, ${\rm W}_{\alpha_i}$ and ${\rm K}_{\alpha_i}$ can be presented as monomials in certain seed. Thus they are quantum lifts of $\mathcal{W}_{s, i}$ and $\mathcal{K}_{s, i}$. Theorem \ref{UEAB} i) for a triangle $t$ follows from Proposition \ref{Prop11.11}.

\vskip 2mm

Now let $s$ be a special point of a general surface ${\bS}$. Let us assume that the boundary component carrying $s$ has at least 2 special points. We may cut off a triangle $t$ from $\bS$ such that $t$ contains the special point $s$ and the two boundary intervals adjacent to $s$. Then Theorem  \ref{QWE}  and the first part of Theorem  \ref{UEAB}  is reduced to the triangle case, which has been proved above. 
By Theorem \ref{THH6.4} iii), the group  $W^n$ acts on $\mathscr{P}_{\G, \bS}$ as cluster automorphisms. The functions $\mathcal{W}_{s, i}$ and $\mathcal{K}_{s, i}$ are standard monomials invariant under the  $W^n$-action. By Quantum Lift Theorem \ref{quantum.promotion.f}, their quantum lifts ${\rm W}_{s, i}, ~{\rm K}_{s, i}$  are  $W^n$-invariants.
\vskip 2mm

If $s$ is the only special point at a boundary component, the potential $\mathcal{W}_{s, i}$ is not necessarily a standard monomial. In this case, let us give an alternative definition of the quantum lift of $\mathcal{W}_{s, i}$.

 Let $\mathcal{T}$ be an ideal triangulation of $\bS$. Recall the gluing map
$$
 \gamma_{\cal T} :~ \prod_{t \in {\cal T}} \mathscr{P}_{\G, t} \lra \mathscr{P}_{\G, \bS}. 
$$
 It gives rise to natural embeddings:
 \be
 \la{embed.TTT}
 \mathcal{O}\left(  \mathscr{P}_{\G, \bS} \right) \hlra \mathcal{O}\Bigl( \prod_{t \in {\cal T}} \mathscr{P}_{\G, t}\Bigr).
 \ee  
A special point $s$ of $\bS$  determines   an ordered sequence of special points $s_1, \ldots, s_n$ - the vertices of the triangles of   ${\cal T}$ containing $s$. 

\bl Under the embedding \eqref{embed.TTT}, we have
\be
\la{382b}
\begin{split}
&\mathcal{K}_{s, i} =  \prod_{k=1}^{n} \mathcal{K}_{s_k, i},\\
&\mathcal{W}_{s, i} = \mathcal{W}_{s_1, i} + \mathcal{K}_{s_1, i} \mathcal{W}_{s_2, i} + \ldots + \Bigl( \prod_{k=1}^{n-1} \mathcal{K}_{s_k, i}\Bigr) \cdot   \mathcal{W}_{s_n, i}. \\
\end{split}
\ee
\el
 \begin{proof} It suffices to prove the case when $n=2$. As shown on Figure \ref{pin31}, let  ${\rm q}, {\rm q}'$ be the pinnings near  $s_1$, with   ${\rm q}$  on the left. 
Let ${\rm p}, {\rm p}'$ be pinnings near $s_2$, with  ${\rm p}$ on the left. 

\begin{figure}[ht]
\centerline{
\begin{tikzpicture}[scale=1.2]                                 
 \draw (-1.5,-1.2) -- (0,0) -- (0,-2);      
 \draw (2.9,-1.2) -- (1.4,0) -- (1.4,-2);     
\node   (a) at (-0.3,0) {};
\node   (b) at (0.2,-.2) {};
\node   (c) at (1.2,-.2) {};
\node   (d) at (1.7,0) {};
\node [circle,draw=red,fill=red,minimum size=3pt,inner sep=0pt, label=above: {\scriptsize $s_1$}]   at (0,0) {};
\node [circle,draw=red,fill=red,minimum size=3pt,inner sep=0pt, label=above: {\scriptsize $s_2$}]   at (1.4,0) {};
\node at (-.5,-1) {{\small $\mathcal{W}_{s_1, i}$}};
\node at (1.9,-1) {{\small $\mathcal{W}_{s_2, i}$}};
\node[blue] at (.45,-1) {{\small $q'$}};
\node[blue] at (.95,-1.05) {{\small $p$}};
\node[blue] at (-1,-.3) {{\small $q$}};
\node[blue] at (2.4,-.3) {{\small $p'$}};
\draw[<-,>=latex', blue] (-1.6,-1.04) -- (a);
\draw[<-,>=latex', blue] (3,-1.04) -- (d);
\draw[->,>=latex', blue] (b) -- (.2,-2);
\draw[->,>=latex', blue] (c) -- (1.2,-2); 
\draw[->, >=stealth] (3.5,-.5) -- (4.9, -.5); 
\node at (4.2, -.3) {{\small \bf glue}};
\begin{scope}[shift={(7,0)}]
\begin{scope}[shift={(-1.4,0)}]
 \draw (2.9,-1.2) -- (1.4,0);   
 \node[blue] at (2.4,-.3) {{\small $p'$}};
 \node [inner sep=0pt]  (A) at (1.7,0) {};
 \draw[->,>=latex', blue] (A) -- (3,-1.04); 
\end{scope}
 \draw (-1.5,-1.2) -- (0,0);
 \draw[dashed, line width=0.07mm] (0,0) -- (0,-2);        
\node [inner sep=0pt]  (a) at (-0.3,0) {};
\node [circle,draw=red,fill=red,minimum size=3pt,inner sep=0pt, label=above: {\scriptsize $s$}]   at (0,0) {};
\node at (0.1,-1) {{\small $\mathcal{W}_{s, i}$}};
\node[blue] at (-1,-.3) {{\small $q$}};
\draw[<-,>=latex', blue] (-1.6,-1.04) -- (a);
\end{scope}
     \end{tikzpicture}}
\caption{Gluing   the moduli spaces  by identifying the pinnings: ${\rm q}'={\rm p}$.}
\label{pin31}
\end{figure}

Let us identify $q$ with a standard pinning. There is a unique element $b_{s_1}\in \B$ mapping ${\rm q}$ to ${\rm q'}$. Similarly, let $b_{s_2} $ map ${\rm p}$ to ${\rm p}'$, and  $b_s$ map ${\rm q}$ to ${\rm p}'$.
Here  $b_{s} = b_{s_1} b_{s_2} $. Consider the decomposition $b=uh$, where $u \in \U$ and $h \in {\rm H}$. We have
\[
h_s= h_{s_1} h_{s_2}, \hskip 7mm u_s = u_{s_1} \cdot {\rm Ad}_{h_{s_1}} (u_{s_2}).
\]
They prove the identities \eqref{382b} respectively.

\end{proof}

 \paragraph{\bf Examples.}   1. According to (\ref{EXFO}), 
 the cluster Poisson coordinate $X_{13}$ on the  space ${\mathscr P}_{ {\rm PGL_2}, t}$ assigned to the side $(13)$ is  potential   ${\cal W}_1$ at the vertex $1$: 
 $X_{13} =  {\cal W}_{1}$.

\vskip 1mm 2. Let $\bS = {\rm P}_n$ be an $n$-gon.  
Pick its triangulation ${\cal T}$. 
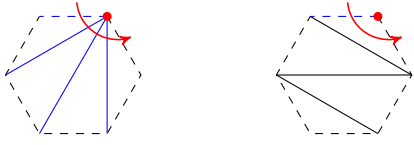
\begin{center}
\begin{tikzpicture}[scale=0.9]
\draw[dashed] (60:1) -- (0:1) -- (-60:1) -- (-120:1) -- (-180:1) --(-240:1); 
\draw[dashed, blue] (60:1) -- (120:1); 
\draw[blue] (60:1) -- (-60:1); 
\draw[blue] (60:1) -- (-120:1); 
\draw[blue] (60:1) -- (-180:1); 
\node [circle,draw=red,fill=red,minimum size=3pt,inner sep=0pt] (b) at (60:1) {};
\begin{scope}[shift={(60:1.3)}]
\draw[->,>=stealth',semithick, red] (185:.6cm) arc (185:290:.6cm);
\end{scope}
\begin{scope}[shift={(4,0)}]
\draw[dashed] (60:1) -- (0:1) -- (-60:1) -- (-120:1) -- (-180:1) --(-240:1); 
\draw[dashed, blue] (60:1) -- (120:1); 
\draw (120:1) -- (0:1); 
\draw (0:1) -- (180:1); 
\draw (-60:1) -- (180:1); 
\node [circle,draw=red,fill=red,minimum size=3pt,inner sep=0pt] (b) at (60:1) {};
\begin{scope}[shift={(60:1.3)}]
\draw[->,>=stealth',semithick, red] (185:.6cm) arc (185:290:.6cm);
\end{scope}
\end{scope}
\end{tikzpicture}
\end{center}

Take  a vertex $v$ and a counterclockwise 
 arc $\alpha_v$ near $v$. Denote by $E_1, ..., E_m$   
the  sides and diagonals of  ${\cal T}$   intersecting the arc, 
counted along the arc. 
  It follows from (\ref{382b}) that the potential ${\cal W}_v$ on the space ${\mathscr P}_{{\rm PGL_2}, \bS}$ is given by 
\be \la{FPW}
{\cal W}_v:= X_{E_1}  + X_{E_1} X_{E_2}  + \cdots + X_{E_1} X_{E_2} \ldots X_{E_{m-1}}.
\ee

\bl \la{MLe}The potential  ${\cal W}_v$ is a   standard monomial in certain cluster Poisson coordinates. \el

\begin{proof}  Mutate any triangulation  to the one which has no diagonals incident to the vertex $v$. \end{proof}

As an quantum analog  of \eqref{embed.TTT}, we get an embedding
\be
\la{emb.auzn.2}
  \mathcal{O}_q\left(  \mathscr{P}_{\G, \bS} \right) \hlra \mathcal{O}_q\Bigl( \prod_{t \in {\cal T}} \mathscr{P}_{\G, t}\Bigr).
\ee
We define   quantum lifts of $\mathcal{W}_i$ and $\mathcal{K}_i$   as   for triangle cases. Identity \eqref{382b} provides an alternative definition of quantum lifts of $\mathcal{W}_{s, i}$ and $\mathcal{K}_{s, i}$ for a general $\bS$. Using  embedding \eqref{emb.auzn.2}, we set
\be
\la{quantum.prom.def2}
\begin{split}
&{\rm K}_{s, i} := {\rm K}_{s_1, i}\cdots {\rm K}_{s_n, i}\\
&{\rm W}_{s, i}:= {\rm W}_{s_1, i} + \mathcal{\rm K}_{s_1, i} {\rm W}_{s_2, i} + \ldots +  {\rm K}_{s_1, i}\cdots {\rm K}_{s_{n-1}, i} \cdot   {\rm W}_{s_n, i}\\
\end{split}
\ee

Let us show that \eqref{quantum.prom.def2} is independent of the triangulation $\mathcal{T}$ chosen. If $\bS$ is a quadrilateral, then it admits two triangulations, corresponding to two diagonals of $\bS$. 
For the diagonal  disconnecting the special point $s$,  the functions in \eqref{quantum.prom.def2} are standard polynomials in $\mathcal{O}_q(\mathscr{P}_{\G, \bS})$ since $n=1$. 
For the diagonal connecting $s$, we use the angles near $s$ and reduced words ${\bf i}, {\bf i}'$ to obtain a seed for $\mathcal{O}_q(\mathscr{P}_{\G, \bS})$. By \eqref{383a}, the functions in \eqref{quantum.prom.def2} are still of the basic pattern \eqref{378}. By Lemma \ref{mutation.patterns.}, they are standard monomials. The very fact that it works for quadrilateral allows us to flip diagonals of $\mathcal{T}$ for arbitrary surfaces. So it is independent of $\mathcal{T}$.

The definition using \eqref{quantum.prom.def2} works for every special point $s$.
When $s$ is not the only special point at its boundary component, this definition  coincides with the previous definition using standard monomials of quantum cluster algebras. 
It completes the proof of Theorem \ref{QWE}.

The proof of the second part of Theorem \ref{UEAB} requires a few preparations.

\subsubsection{The outer monodromy is in the    center.}
Let $\pi$ be a boundary component of $\bS$ with  special points $s_1, \ldots, s_d$.  Recall the ``outer monodromy" map $\mu_{(\pi)}: \mathscr{P}_{\G, \bS} \lra {\rm H}_{(\pi)}$, see (\ref{DEF1.15}). 
For every $i \in {\rm I}$, we get a regular function
\[
\mathcal{K}_{i}^{(\pi)}:=\alpha_i\circ \mu^*_{(\pi)}.
\]
\bp \la{14.3.1} 1) The function $\mathcal{K}_{i}^{(\pi)}$ Poisson commutes with all the functions of ${\cal O}^{\rm cl}(\mathscr{P}_{\G, \bS})$.

\vskip 1mm2) The quantum lifts ${\rm K}_{i}^{(\pi)}$ of the functions $\mathcal{K}_{i}^{(\pi)}$ lie in the center of ${\cal O}_q(\mathscr{P}_{\G, \bS})$.
\ep  \la{P11.8}
 \begin{proof} 1) We give a proof for $d$ being even.
The proof for odd $d$ is similar. By definition, we have
 \[
 \mathcal{K}_{i}^{(\pi)}= \mathcal{K}_{s_1, i} \mathcal{K}_{s_2, i^\ast} \ldots \mathcal{K}_{s_d, i^\ast}. 
\]
By Lemma  \ref{L5.10} ii), every  $\mathcal{K}_{s_k, i}$ commutes with all the cluster Poisson variables except  the frozen ones placed on the two boundary intervals adjacent to $s_k$. Therefore $ \mathcal{K}_{i}^{(\pi)}$ commutes with all cluster Poisson variables that are not  on the boundary component $\pi$. For a frozen variable on the boundary interval connecting $s_{k-1}$ and $s_k$, by Lemma \ref{L5.10} i) and ii), it commutes with the product $\mathcal{K}_{s_{k-1}, i} \mathcal{K}_{s_{k}, i^\ast}$. So it  commutes with $ \mathcal{K}_{i}^{(\pi)}$.   

2) Follows immediately from 1). \end{proof} 
For example, if $\pi$ contains exactly two special point $e$ and $f$, we get central elements
 \[
 {\rm K}^{(\pi)}_i:= {\rm K}_{e, i} {\rm K}_{f, i^\ast} \in {\rm Center} \left(\mathcal{O}_q(\mathscr{P}_{\G, \bS})\right).
 \]

\subsubsection{Quantum relations.} Let us label the vertices of a triangle $t$ counterclockwise by $1, 2, 3$. For every $k\in \{1,2,3\}$ and every $i \in {\rm I}$, we obtain a pair 
 \[
 {\rm W}_{k, i}, \qquad{\rm K}_{k, i}  \in \mathcal{O}_q(\mathscr{P}_{\G, t}).
 \]

 \bl 
 \la{quantum.mota.rule.1}
 We have the following relations, where  $\delta_{ij}$ be the Kronecker symbol and $q_i= q^{1/d_i}$:
  \[
 \begin{split}
 &[ {\rm W}_{2, i},  {\rm W}_{1, j^\ast}] = \delta_{ij} (q_{i}^{-1} - q_i) {\rm K}_{2, i},\\
 &[{\rm W}_{2, i}, {\rm K}_{1,j^\ast}]=0, \\
 &{\rm K}_{2, i} {\rm W}_{1, j^\ast} = q_{j}^{-{\rm C}_{ij}} {\rm W}_{1, j^\ast}  {\rm K}_{2, i},\\
 & {\rm K}_{2, i} {\rm K}_{1, j^\ast} = q_j^{-{\rm C}_{ij}} {\rm K}_{1, j^\ast}{\rm K}_{2, i}. \\
\end{split}
 \]
Due to  cyclic invariance of the cluster structure on $\mathscr{P}_{\G, t}$, the same relations hold when the indices $(1,2)$ are replaced by $(2, 3)$ or $(3,1)$.
 \el
 \begin{proof}
 Consider the seed of $\mathscr{P}_{\G, t}$ associated to the vertex $1$ and a reduced word ${\bf i}$ of $w_0$ starting with $i^\ast$. The vertices at level $i^\ast$ and the first frozen vertex on the side $\{23\}$ form a full subquiver
 
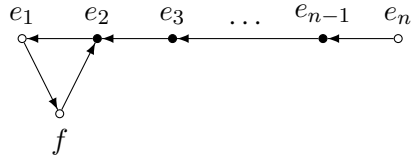
\begin{center}
\begin{tikzpicture}
\node [circle,draw,fill=white,minimum size=3pt,inner sep=0pt, label=above:$e_1$] (a) at (0,0) {};
\node [circle,draw,fill=white,minimum size=3pt,inner sep=0pt, label=below:$f$] (A) at (0.5,-1) {};
\node [circle,draw,fill,minimum size=3pt,inner sep=0pt, label=above:$e_2$] (b) at (1,0) {};
 \node [circle,draw,fill, minimum size=3pt,inner sep=0pt, label=above:$e_3$]  (c) at (2,0){};
  \node [circle,draw,fill, minimum size=3pt,inner sep=0pt, label=above:$e_{n-1}$]  (d) at (4,0){};
 \node [circle,draw,fill=white,minimum size=3pt,inner sep=0pt, label=above:$e_n$] (e) at (5,0) {};
 \foreach \from/\to in {b/a,a/A, A/b, c/b, e/d}
                  \draw[directed] (\from) -- (\to);
\foreach \from/\to in {d/c}                   
              \draw[directed] (\from) -- (\to);            
  \node at (3,0.2) {$\cdots$};            
     \end{tikzpicture}
    \end{center} 
 By Lemma  \ref{L5.10}, we get 
  \[
  \begin{split}
  &{\rm K}_{1, i^\ast}= \X_{e_1 + \ldots + e_{n}},\qquad ~\qquad~ {\rm W}_{1, i^\ast}= \X_{e_1}+ \X_{e_1 + e_2}+ \ldots + \X_{e_1 + \ldots + e_{n-1}},\\
  & {\rm K}_{2, i}= X_{f+ e_1},\qquad \qquad \qquad{\rm W}_{2, i}= \X_{f}.\\
  \end{split}
  \]
  Therefore
  \[
 \begin{split}
 & [{\rm W}_{2,i}, {\rm W}_{1, i^\ast}] = [\X_f, \X_{e_1}] = (q_i^{-1}-q_i){\rm K}_{2, i}.\\
&  [{\rm W}_{2, i},  {\rm K}_{1, i^\ast}] = [\X_f, \X_{e_1+\ldots+e_n}] = 0.\\
  \end{split}
    \]
  For any vertex $v$ at level $j\in {\rm I}-\{i\}$, there is no arrow between $f$ and $v$. Therefore ${\rm W}_{2,i}$ commutes with ${\rm W}_{1, j^\ast}$ and ${\rm K}_{1, j^\ast}$. The last two identities follow from the relation between ${\rm K}_{2, i}$ and frozen cluster variables at the side $\{12\}$.
 \end{proof}

 \vskip 2mm \subsubsection{\bf Proof of Theorem \ref{UEAB} iii).} Let $\pi$ a boundary component of $\bS$ consists of exactly two special points $e$ and $f$. It remains to show that 
 \be
 \la{quam.equa.c}
 [{\rm W}_{e, i}, {\rm W}_{f, j^\ast}] = \delta_{ij} (q_i^{-1}-q_i) ({\rm K}_{e, i} -{\rm K}_{f, i^\ast}).
 \ee
 The rest relations follows from part i) of Theorem \ref{UEAB}. 
 
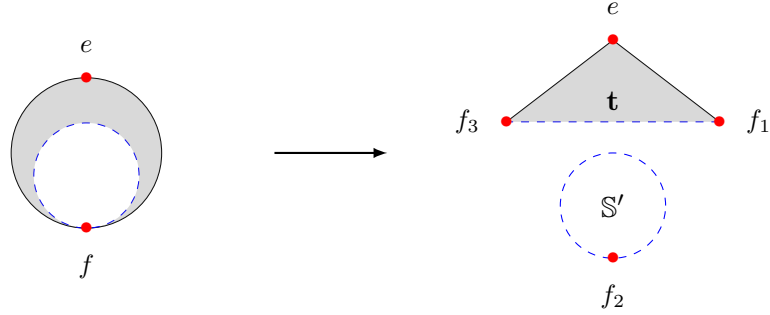
  \begin{center}
\begin{tikzpicture}
\draw (0,0) circle (1cm);
\draw[dashed, blue] (0,-.3) circle (7mm);  
\fill [gray,opacity=0.3] (0,-1) arc (-90:270:1) -- (0, -1) arc (270:-90:0.7);
\node[red, label=above: {\small $e$}] at (0, 1) {\small $\bullet$};
\node[red, label=below: {\small $f$}] at (0, -1) {\small $\bullet$};
\draw[-latex,thick] (2.5,0) -- (4,0);
\begin{scope}[shift={(7,-1)}]
\draw (45:2) -- (90:2.5) -- (135:2);
\draw[blue, dashed] (45:2) -- (135:2);
\fill [gray,opacity=0.3] 
(45:2) -- (90:2.5) -- (135:2) -- (45:2);
\draw[dashed, blue] (0,0.3) circle (7mm);  
\node[red, label=above: {{\small $e$}}] at (90:2.5) {\small $\bullet$};
\node[red, label=right: {{\small $f_1$}}] at (45:2) {\small $\bullet$};
\node[red, label=left: {{\small $f_3$}}] at (135:2) {\small $\bullet$};
\node[red, label=below: {{\small $f_2$}}] at (0, -.4) {\small $\bullet$};
\node at (90:1.7) {\bf t};
\node at (0, .3) {$\bS'$};
\end{scope}
\end{tikzpicture}
 \end{center}

Cutting $\bS$ along the dashed edge, we obtain a triangle $t$ and a surface $\bS'$. It gives rise to an embedding
\[
\mathcal{O}_q\left(\mathscr{P}_{\G, \bS}\right) \hlra \mathcal{O}_q\left(\mathscr{P}_{\G, t}\times \mathscr{P}_{\G, \bS'}\right).
\]
By \eqref{quantum.prom.def2}, we get
\[
 \begin{split}
 &{\rm K}_{f,j^\ast}={\rm K}_{f_1,j^\ast}{\rm K}_{f_2,j^\ast}{\rm K}_{f_3,j^\ast},\\
 &{\rm W}_{f, j^\ast} = {\rm W}_{f_1,j^\ast}+ {\rm K}_{f_1,j^\ast}{\rm W}_{f_2,j^\ast}+ {\rm K}_{f_1,j^\ast}{\rm K}_{f_2,j^\ast}{\rm W}_{f_3,j^\ast}.\\
 \end{split}
  \]
 By Lemma \ref{quantum.mota.rule.1}, we get
 \[
 \begin{split}
 & [{\rm W}_{e, i}, ~{\rm W}_{f_1, j^\ast}] = \delta_{ij} (q_{i}^{-1} - q_i) {\rm K}_{e, i},\\
 & [{\rm W}_{e, i}, ~ {\rm K}_{f_1,j^\ast}{\rm W}_{f_2,j^\ast}]=0,\\
& [{\rm W}_{e, i}, ~ {\rm K}_{f_1,j^\ast}{\rm K}_{f_2,j^\ast}{\rm W}_{f_3,j^\ast}]= {\rm K}_{f_1,j^\ast}{\rm K}_{f_2,j^\ast} [{\rm W}_{e, i}, {\rm W}_{f_3,j^\ast}] = -  \delta_{ij} (q_{i}^{-1} - q_i) {\rm K}_{f, i^\ast}.
\end{split}
 \]
 Adding them together, we obtain \eqref{quam.equa.c}.

 \subsubsection{\bf Proof of Theorem \ref{5.25.24q}.} \la{24.5.26.10.19}
 Consider the ideal triangle containing the boundary intervals connected to the marked point $t$. When we cut along the dashed edge, the angle at $s$ is divided into two parts, denoted $s_1$ and $s_2$, as shown in the following figure.
 
  \begin{center}
 \begin{tikzpicture}
 \fill [gray,opacity=0.3] (0,0) --(2,0)--+(-60:2)--(0,0);
 \draw [thick] (0,0)--(2,0)--+(-60:2);
 \draw [thick, dashed] (0,0)--(3, -1.732);
 \node at (1.7,-.5) {${\bf t}$};
 \node at (0,0.3) {$s$};
 \node at (2, 0.3) {$t$};
 \node[red] at (0,0) {$\bullet$};
 \node[red] at (2,0) {$\bullet$};
 \draw[thick] (-120:2)--(0,0);
 \draw[-latex] (4,-.5)--(6,-.5);
 \end{tikzpicture}
 \hskip 1cm
 \begin{tikzpicture}
 \fill [gray,opacity=0.3] (1,0) --(3,0)--+(-60:2)--(1,0);
 \draw [thick] (1,0)--(3,0)--+(-60:2);
 \draw [thick, dashed] (1,0)--(4, -1.732);
 \draw [thick, dashed] (0,0)--(3,-1.732);
 \node at (2.7,-.5) {${\bf t}$};
 \node at (0,0.3) {$s_1$};
 \node at (1,0.3) {$s_2$};
 \node at (3, 0.3) {$t$};
 \node[red] at (0,0) {$\bullet$};
 \node[red] at (1,0) {$\bullet$};
 \node[red] at (3,0) {$\bullet$};
 \draw[thick] (-120:2)--(0,0);
 \end{tikzpicture}
\end{center}
 We have
 \[
 {\rm K}_{s,i}= {\rm K}_{s_1, i}{\rm K}_{s_2, i}, \qquad 
 {\rm W}_{s,i}= {\rm W}_{s_1,i}+{\rm K}_{s_1,i} {{\rm W}_{s_2,i}}.
 \]
 Here ${\rm W}_{s_1,i}, ~{\rm K}_{s_1,i}$ commute with ${\rm W}_{t,j^\ast},~ {\rm K}_{t,j^\ast}$. The relations among $\{{\rm W_{s_2,i}}, ~{\rm K}_{s_2, i}\}$ and $\{{\rm W}_{t, j^\ast},~ {\rm K}_{t, j^\ast}\}$ are described by Lemma \ref{quantum.mota.rule.1}. Therefore
 \[
 [{\rm W}_{s,i}, ~{\rm W}_{t,j^\ast}]=[{\rm K}_{s_1,i}{\rm W}_{{s_2},i}, ~{\rm W}_{t,j^\ast}]= {\rm K}_{s_1,i}[{\rm W}_{s_2, i}, ~{\rm W}_{t,j^\ast}]= \delta_{ij}(q_i^{-1}-q_i){\rm K}_{s, i}.
 \]
 The rest formulas follow by a similar calculation. 
 \medskip

 \subsection{Geometry of the {PBW} bases and Lusztig's braid group action} \la{SSECC8}
\medskip

\bl 
\la{181}
The following three spaces are canonically isomorphic to each other
\begin{enumerate}

\vskip 1mm\item ${\rm U}_\ast := {\rm U} \cap \B^-w_0 \B^-$. 

\vskip 1mm\item the moduli space ${\rm Conf}^\times_3(\mathcal{A})_{e,e}$  parametrizes $\G$-orbits of triples $(\A_1, \A_2, \A_3)\in (\G/\U)^3$ such that the pairs $(\A_1, \A_2)$,  $(\A_2,\A_3)$, 
$(\A_1,\A_3)$ are generic and
\be \la{condd}
h(\A_1, \A_2)=h(\A_2, \A_3)=1.
\ee
\vskip 1mm\item the moduli space  ${\mathscr P}_{\G, t}^{(1)}$ parametrizes $\G$-orbits of the data $(\B_1, \B_2, \B_3; p_{23})$, where the flags $(\B_1, \B_2, \B_3)$ are pairwisely generic, and 
$p_{23}$ is a pinning over $(\B_2, \B_3)$. 
\end{enumerate}
\el
\begin{proof}
The isomorphism from (1) to (2) is given by 
\[ \label{isomo.2018.7.11}
\U_\ast \stackrel{\sim}{\lra} {\rm Conf}^\times_3(\mathcal{A})_{e,e}, \qquad ~ u \lms ( [\U], [\U^-],  \overline{w}_0u \cdot [\U^-])=(\A_1, \A_2, \A_3).
\]
Here the pair $(\A_1, \A_3)$ is generic is due to the very fact that $u\in \U_\ast$.

The isomorphism from (2) to (3) is 
\be \la{MAPii}
\begin{split}
&\rho: {\rm Conf}^\times_3(\mathcal{A})_{e,e} \stackrel{\sim}{\lra} {\mathscr P}_{\G, t}^{(1)},\\
&(\A_1, \A_2, \A_3; h_{12}=h_{23}=1) \lms (\B_1, \A_2, \A_3) = (\B_1, \B_2, \B_3; p_{23}).
\\
\end{split}
\ee
It forgets the decoration of the flag  $\A_1$, and interprets the pair of decorated flags $(\A_2, \A_3)$ satisfying the condition $h(\A_2, \A_3)=1$ as a pinning $p_{23}$ over $(\B_2, \B_3)$. 
The inverse map assigns to   $(\B_1, \A_2, \A_3)$ the unique decorated flag $\A_1$ over   $\B_1$ such that $h(\B_1, \B_2)=1$. 
\end{proof}

Recall   the space ${\mathscr P}_{\G, t}$ parametrizing generic triples of flags with three pinnings $(p_{12}, p_{23}, p_{13})$.  
Forgetting the   pinnings $(p_{12},   p_{13})$ we get  a   surjection:
\[  
 j:~{\mathscr P}_{\G, t}  \lra 
{\mathscr P}_{\G, t}^{(1)}. 
\]
The space ${\mathscr P}_{\G, t}^{(1)}$ inherits a cluster Poisson structure from the former by deleting frozen variables on the sides $\{12\}$ and $\{13\}$. By the quantization of cluster Poisson varieties \cite{FG03b}, we obtain a $\ast$-algebra $\mathcal{O}_q\left( {\mathscr P}_{\G, t}^{(1)} \right)$.

Dropping  conditions (\ref{condd}), we get an embedding:
  \be
i: {\rm Conf}^\times_3(\mathcal{A})_{e, e} \subset {\rm Conf}^\times_3(\mathcal{A}).
\ee 
Hence  ${\rm Conf}^\times_3(\mathcal{A})_{e, e}$ inherits a cluster $K_2$  structure from the latter by restricting the frozen $K_2$ variables on the sides $\{12\}$ and $\{23\}$ to be $1$.
 
 The space ${\rm Conf}_3(\mathcal{A})_{e, e}$, just as an arbitrary cluster $K_2$-variety, {\it a priori}   
does not have neither a canonical Poisson structure, nor a canonical quantization. Instead, after imposing  extra data called {\it compatible pairs}, Bereinstein-Zelevinsky \cite{BZ} provide a family of $q$-deformations of upper cluster algebras, and thus,  in our terminology, cluster $K_2-$varieties. In Section \ref{Sec.quant.s} we explain in what sense  the quantization of \cite{BZ} is compatible with the quantization of \cite{FG03b}.

\bp \la{1.10.19.5} 
The isomorphism (\ref{MAPii}) provides a natural  compatible pair for  ${\rm Conf}^\times_3(\mathcal{A})_{e, e}$. 
 Therefore the  quantum version of the isomorphism (\ref{MAPii})  provides  
  a $\ast$-algebra isomorphism
\be \la{ALGa}
\rho_q^\ast:~ \mathcal{O}_q\left({\mathscr P}_3^{(1)}\right) \stackrel{=}{\lra}
\mathcal{O}_q({\rm Conf}^\times_3\left(\mathcal{A})_{e,e}\right)
\ee
\ep

\begin{proof}  There is a commutative diagram, where $i$  is the map (\ref{condd}), and the map $\pi_3$ in \eqref{2018.9.22.18.42hhX}:
\begin{equation}\label{1.10.19.1}
\begin{gathered}
\xymatrix{
{\rm Conf}^\times_3({\cal  A} )  \ar[d]_{\pi_3}& \ar[l]_{i} {\rm Conf}^\times_3({\cal  A})_{e, e}  \ar[d]^{\rho}_{\sim}&\\
{\mathscr P}_{\G, t} \ar[r]_{j} & {\mathscr P}_{\G, t}^{(1)} &}
\end{gathered}
\end{equation}
Recall that the pair $\left({\rm Conf}^\times_3({\cal  A} ), {\mathscr P}_{\G, t}\right)$ forms a {\it cluster ensemble} of \cite{FG03a}. Therefore the map $\rho$ is a monomial map in arbitrary cluster chart  of the spaces.  

Precisely, let ${\bf s}$ be  a seed of the cluster ensemble $\left({\rm Conf}^\times_3({\cal  A} ), {\mathscr P}_{\G, t}\right)$. Recall that the vertices of ${\bf s}$ are parametrized by the set
\[
J:= J_{uf} \cup J_{12} \cup J_{23} \cup J_{31}.
\]
Each $J_{ab}$ parametrizes the $r$ many frozen vertices of the the side $\{ab\}$. The exchange matrix $\varepsilon$ of ${\bf s}$ is  a $J\times J$ matrix. Let $p$ be the integral $\left( J_{uf} \cup J_{31}\right) \times  \left(J_{uf}  \cup J_{23}\right)$ submatrix of $\varepsilon$. Following Proposition \ref{11.24.18.1}, the monomial map $\rho$ in terms the coordinates associated to the seed ${\bf s}$ is exactly given by the matrix $p$. Since $\rho$ is an isomorphism, we get $\det(p)=\pm 1$. Using $p^{-1}$ in \eqref{comp.BZ.pair}, we obtain a compatible pair of \cite{BZ} following Lemma \ref{lem.comp.BZ.pair}.
\end{proof}

Following the proof of Theorems \ref{QWE}, we get a homomorphism 
\be
\la{split.kappa.map.n}
\kappa: ~\U_q (\n) \lra \mathcal{O}_q\left( {\mathscr P}_{\G, t}^{(1)} \right)=\mathcal{O}_q({\rm Conf}^\times_3\left(\mathcal{A})_{e,e}\right) , \hskip 7mm {\bf E}_i \lms {\rm W}_{2, i}.
\ee
In the rest of this Section,  we present a geometric interpretation of Lusztig's {PBW} bases of $\U_q(\n)$, in terms of  specific sets of quantized cluster $K_2$ variables in $\mathcal{O}_q({\rm Conf}^\times_3\left(\mathcal{A})_{e,e}\right)$.

\begin{figure}[ht]
\epsfxsize 200pt
\center{
\begin{tikzpicture}[scale=1.5]
\node [circle,draw=red,fill=red,minimum size=3pt,inner sep=0pt, label=above:{\small ${\A}_1$}] (a) at (0,0) {};
\node [circle,draw=red,fill=red,minimum size=3pt,inner sep=0pt, label=below:{\small ${\A}_2$}] (b) at (-1.8,-3) {};
 \node [circle,draw=red,fill=white, minimum size=3pt,inner sep=0pt, label=below:{\small ${\A}_2''$}]  (c) at (-0.36,-3){};
  \node [circle,draw=red,fill=white, minimum size=3pt,inner sep=0pt, label=below:{\small ${\A}_1''$}]  (g) at (-1.08,-3){};
 \node [circle,draw=red,fill=white, minimum size=3pt,inner sep=0pt, label=below:{\small ${\A}_{n-1}''$}] (h) at (1.08,-3) {};
 \node [circle,draw=red,fill=red,minimum size=3pt,inner sep=0pt, label=below:{\small ${\A}_3$}] (e) at (1.8,-3) {};
  \node [circle,draw=red,fill=white, minimum size=3pt,inner sep=0pt, label=left:{\small ${\A}_1'$}]  (c1) at (-0.36,-0.6){};
  \node [circle,draw=red,fill=white, minimum size=3pt,inner sep=0pt, label=left:{\small ${\A}_2'$}]  (g1) at (-0.72,-1.2){};
 \node [circle,draw=red,fill=white, minimum size=3pt,inner sep=0pt, label=left:{\small ${\A}_{n-1}'$}] (h1) at (-1.44,-2.4) {};
 \foreach \from/\to in {b/g, g/c, h/e, a/e, a/c1, c1/g1, h1/b}
                  \draw[thick] (\from) -- (\to);
\foreach \from/\to in {c/h, g1/h1}
                  \draw[dotted, thick] (\from) -- (\to);
\foreach \from/\to in {c1/g, g1/c, h1/h}                   
              \draw[dashed] (\from) -- (\to);            
 \node[blue] at (-1.44,-3.14) {\small $s_{i_1^\ast}$};
  \node[blue] at (-0.72,-3.14) {\small ${s}_{i_2^\ast}$}; 
   \node[blue] at (1.44,-3.14) {\small $s_{i_N^\ast}$};          
    \node[blue] at (-0.34,-0.3) {\small $s_{i_1}$};
  \node[blue] at (-0.72,-0.9) {\small ${s}_{i_2}$}; 
   \node[blue] at (-1.81,-2.7) {\small $s_{i_N}$};        
    \end{tikzpicture}
 }
 \caption{Decomposition of the   rotation map $ \overline w_0 u: \{\A_1, \A_2\}\to \{\A_2, \A_3\}$.}
 \label{2018.2.23.9.46sh}
 \end{figure}
 
Let $(\A_1, \A_2, \A_3)\in {\rm Conf}_3(\mathcal{A})_{e, e}$. Recall that  each reduced word ${\bf i}=(i_1,\ldots, i_N)$ gives rise to two sequences of intermediate decorated flags, as illustrated on Figure \ref{2018.2.23.9.46sh}:
\[
\begin{split}
&\A_1= \A_0' \stackrel{s_{i_1}}{\lra} \A_1' \stackrel{s_{i_2}}{\lra} \A_2' \stackrel{s_{i_3}}{\lra}\cdots \stackrel{s_{i_{N}}}{\lra} \A_{N}' = \A_2;\\
&\A_2= \A_0'' \stackrel{s_{i_1}}{\lra} \A_1'' \stackrel{s_{i_2}}{\lra} \A_2'' \stackrel{s_{i_3}}{\lra}\cdots \stackrel{s_{i_{N}}}{\lra} \A_{N}'' = \A_3.\\
\end{split}
\]
Let us set  
\be \la{DEFA}
 a_{{\bf i},k}:= \Delta_{i_k}(\A_{k-1}', \A_k'').
 \ee

 \bl \la{PPPa} Every $a_{{\bf i},k}$ in \eqref{DEFA}  is a cluster $K_2$ variable 
 of ${\rm Conf}^\times_3(\mathcal{A})_{e,e}$.
  \el
  
  \begin{proof} 
  Let $\mathcal{A}_k$ be the partial configuration space parametrizing $\G$-diagonal orbits of quadruples of decorated flags $(\A_1, \A'_k, \A_{k}'', \A_3)$ such that 
  $h(\A_1, \A'_k) = h(\A'_k, \A''_k) =  h(\A''_k, \A_3) =   1$. Recall the sequence of isomorphisms given by reflections in \eqref{seq.of.refs.a}. See Figure \ref{pin20}.  Since the variables assigned to the side $(\A'_k, \A_{k}'')$ are equal to $1$, the monomial transformation associated with the reflection map is trivial. Therefore, when restricted to the subspaces $\mathcal{A}_k\subset {\rm Conf}_{u_k}^{v_k}(\mathcal{A})$, we obtain a chain of cluster $K_2$ isomophisms, see Figure \ref{pin20}: 
 \be
 \la{refl.confaee.ak}
 {\rm Conf}^\times_3(\mathcal{A})_{e,e} =  {\mathscr A}_0 \stackrel{r_1}{\lra} {\mathscr A}_{1} \stackrel{r_2}{\lra} {\mathscr A}_{2} \stackrel{r_3}{\lra} \ldots \stackrel{r_{N}}{\lra} {\mathscr A}_{N}.
  \ee
 Clearly $a_{{\bf i},k}$ is a cluster variable on ${\mathscr A}_{k}$. 
  Using $(r_1\circ r_2 \circ \ldots \circ r_{k})^{-1}$, we conclude that $a_k$ is a cluster variable on  ${\rm Conf}^\times_3(\mathcal{A})_{e,e}$.  
      \begin{figure}[bt]
\epsfxsize150pt
\center{
 \begin{tikzpicture}[scale=0.6]
  \draw[thick] (-14,-3)  --(-14,0)-- (-12.5,0) -- (-12.5,-3)-- (-14,-3); 
  \node [circle,draw=red,fill=red,minimum size=3pt,inner sep=0pt,label=below:${\A}_{k}''$]  at (-14,-3) {};
  \node [label=below:${\A}_{3}$]  at (-12.5,-3) {};
   \node [label=below:${\A}_{1}$]  at (-12.5,1.3) {};
    \node [label=below:${\A}_{k}'$]  at (-14,1.3) {};
  \path [fill=red!30] (-8,-3) -- (-7.33, -1) --(-7,0)-- (-7,-3) -- (-8,-3);         
   \path [fill=red!30] (-1,0) --(0,0)-- (0,-3) -- (-1,0);         
 \node [circle,draw=red,fill=red,minimum size=3pt,inner sep=0pt,label=below:{\small ${\A}_{k}''$}] (a) at (-8,-3) {};
\node [circle,draw=red,fill=red,minimum size=3pt,inner sep=0pt,label=above:{\small ${\A}_{k}'$}] (e) at (-7,0) {};
 \node [circle,draw=red,fill=white,minimum size=3pt,inner sep=0pt,label=left:{\small ${\A}_{k+1}'$}]  (f) at (-7.33,-1){};
\node [circle,draw=red,fill=red,minimum size=3pt,inner sep=0pt,label=above:{\small ${\A}_1$}] (b) at (-5.5,0) {};
 \node [circle,draw=red,fill=white,minimum size=3pt,inner sep=0pt,label=below:{\small ${\A}_{k+1}''$}]  (d) at (-7,-3){};
 \node [circle,draw=red,fill=red,minimum size=3pt,inner sep=0pt,label=below:{\small ${\A}_3$}] (c) at (-5.5,-3) {};
 \foreach \from/\to in {e/b, b/c, d/c, f/a, e/f, a/d}
                  \draw[thick] (\from) -- (\to);
\foreach \from/\to in {e/d}
                  \draw[dashed] (\from) -- (\to);                     
  \draw[thick,->] (-4,-1.5) -- (-2,-1.5);   
               \node at (-3,-1.9) {\small $\mbox{reflection}$};         
\node [circle,draw=red,fill=red,minimum size=3pt,inner sep=0pt,label=above:{\small ${\A}_{k+1}'$}] (a) at (-1,0) {};
\node [circle,draw=red,fill=white,minimum size=3pt,inner sep=0pt,label=above:{\small ${\A}_{k}'$}] (e) at (0,0) {};
\node [circle,draw=red,fill=red,minimum size=3pt,inner sep=0pt,label=above:{\small ${\A}_1$}] (b) at (1.5,0) {};
 \node [circle,draw=red,fill=red,minimum size=3pt,inner sep=0pt,label=below:{\small ${\A}_{k+1}''$}]  (d) at (0,-3){};
 \node [circle,draw=red,fill=red,minimum size=3pt,inner sep=0pt,label=below:{\small ${\A}_3$}] (c) at (1.5,-3) {};
 \foreach \from/\to in {a/e, e/b, b/c, d/c, d/a}
                  \draw[thick] (\from) -- (\to);
\foreach \from/\to in {e/d}
                  \draw[dashed] (\from) -- (\to);                               
 \end{tikzpicture}
  }
\caption{The configuration space ${\mathscr A}_k$, and the reflection map $r_{k+1}: {\mathscr A}_k \lra {\mathscr A}_{k+1}$.}
\label{pin20}
\end{figure}
  \end{proof}

\begin{remark}
Under  isomorphism \eqref{isomo.2018.7.11}, the element $\overline{w}_0 u$ takes the pair $\{\A_1, \A_2\}=([\U], [\U^-])$ to $\{\A_2, \A_3\}$. We obtain a decomposition of $\overline{w}_0 u$ by considering the maps
\[
\{\A_1, \A_2\} = \{\A_0',\A_0''\} \lra \{\A_1', \A_1''\}\lra \cdots \lra \{\A_N', \A_N''\}=\{\A_2, \A_3\}.
\]
Then, using  $a_{{\bf i}, k}$,  each   step corresponds to a group element 
$
z_{i_k}(a_{{\bf i}, k}) := \overline{s}_{i_k} x_{i_k}(a_{{\bf i}, k}). 
$
They provide a decomposition of the group element $\overline{w}_0 u$:
\be
\la{natural.dcom.u.zz}
\overline{w}_0 u = z_{i_1}(a_{{\bf i},1})\ldots z_{i_N}(a_{{\bf i}, N}).
\ee
These functions  give rise to an open embedding
\be
\la{open embedding.e}
\alpha_{\bf i}: ~\qquad~ {\rm Conf}_3(\mathcal{A})_{e,e}=\U_{\ast} {\hlra} \mathbb{A}^N,  \qquad  (\A_1, \A_2, \A_3) \lms (a_{{\bf i},1}, \ldots, a_{{\bf i},N}).
\ee
Each $a_{{\bf i},k}$ is a cluster variable by Proposition \ref{PPPa}. Note that $a_{{\bf i},k}$'s are not from the same seed. Therefore $\alpha_{\bf i}$ is not a cluster coordinate system.

One easily checks that 
\be
\label{z-tran0} 
\left.
\begin{array}{ll}
z_i(a)z_j(b)=z_{j}(b)z_i(a), \qquad &\mbox{if } {\rm C}_{ij}=0.\\
z_i(a) z_j(b) z_i(c) = z_j(c) z_i (ac- b) z_j(a), \qquad ~\qquad~& \mbox{if } {\rm C}_{ij}={\rm C}_{ji}=-1.\\
z_i(a)z_j(b)z_i(c)z_j(d)=z_j(d)z_i(ad^2-2bd+c)z_j(ad-b)z_i(a), \qquad ~\qquad~ & \mbox{if } {\rm C}_{ij}=2{\rm C}_{ji}=-2.\\
z_i(a)z_j(b)z_i(c)z_j(d)z_i(e)z_j(f)=z_j(f)z_i(p)z_j(q)z_i(r)z_j(s)z_i(a), \qquad ~\qquad~ & \mbox{if } {\rm C}_{ij}=3{\rm C}_{ji}=-3.\\
\end{array}
\right.
\ee
where
\[
p=af^3-3bf^2+3df-e, \hskip 3mm q= af^2-2bf+d, \hskip 3mm r=a^2f^3-3abf^2+3b^2f+c-3bd+ae, \hskip 3mm s=af-d.
\]
For any pair of reduced decompositions ${\bf i}$ and ${\bf i}'$ of $w_0$, the transition map $\alpha_{\bf i}\circ\alpha_{\bf i'}^{-1}$ is a bijection expressed as a composition of  
  elementary transition maps provided by (\ref{z-tran0}).
\end{remark}

 The cluster variables $\alpha_{\bf i}=\{a_{{\bf i},1}, \ldots, a_{{\bf i},N}\}$ give rise to quantum cluster variables
\be \la{1.10.19.6}
\{{\rm A}_{{\bf i},1}, \ldots, {\rm A}_{{\bf i},N}\} \in \mathcal{O}_q({\rm Conf}^\times_3(\mathcal{A})_{e,e}).
\ee
By definition, under the isomorphism \eqref{ALGa}, they are quantum standard monomials in $\mathcal{O}_q\left(\mathscr{P}_{\G, t}^{(1)}\right)$.

\bt \la{REALI} 
Recall the elements ${\bf E}_{{\bf i}, k}$ in \eqref{ELLWin} and map $\kappa$ in \eqref{split.kappa.map.n}. Then

i) If $\G$ is simply-laced,  then
\be\la{EAa}
\kappa(\bE_{{\bf i}, k} )= {\rm A}_{{\bf i}, k}.
\ee

ii) If  (\ref{EAa}) holds, then $\kappa$ is injective.  So the map $\kappa$ is injective for simply laced $\G$. 
\et

\bl \la{LLLa0} Let ${\bf i}=(..., i_{k-1}, i_k,   ...) = (..., i, j,  ...)$ with ${\rm C}_{ij}=0$. Let  ${\bf i}'=  (..., j, i,  ...)$ be  obtained  from ${\bf i}$ by switching $i\leftrightarrow j$. Then 
\be \la{EQQQ10}
{\rm A}_{{\bf i}',{k-1}}\stackrel{}{=} {\rm A}_{{\bf i}, {k}} 
\ee
\el
\begin{proof}
By the first identity of \eqref{z-tran0}, the cluster variables
\[
a_{{\bf i}',{k-1}}=a_{{\bf i},k}.  
\]
Therefore their quantum lifts coincide.
\end{proof}

\bl \la{LLLa}
 Let ${\bf i}=(..., i_{k-1}, i_k, i_{k+1}, ...) = (..., i, j, i, ...)$ with ${\rm C}_{ij}={\rm C}_{ji}=-1$, and  ${\bf i}'=(..., j, i, j, ...)$. Then 
\be
\label{local.check220}
{\rm A}_{{\bf i}', {k-1}}\stackrel{}{=} {\rm A}_{{\bf i}, {k+1}}.  
\ee
and 
\be
\label{local.check22}
{\rm A}_{{\bf i},k} = \frac{q^{1/2}{\rm A}_{{\bf i}, k-1}{\rm A}_{{\bf i}, k+1}- q^{-1/2}{\rm A}_{{\bf i}, k+1}{\rm A}_{{\bf i}, k-1}}{q-q^{-1}}.
\ee
\el
\begin{proof}
By the second identity of  \eqref{z-tran0}, we get $a_{{\bf i}',{k-1}}=a_{{\bf i},k+1}$.  
The formula \eqref{local.check22} follows.

Let us identify ${\rm Conf}_3(\mathcal{A})_{e, e}$ with $\mathcal{A}_{k-2}$ via the reflections \eqref{refl.confaee.ak}.
For the latter, we have the following cluster $K_2$-coordinates from one seed
\[
\begin{split}
&{\Delta_i(\A_{k-2}', \A_{k-1}'')}:= x_{v_1}, \qquad {\Delta_j(\A_{k-2}', \A_{k}'')}:= x_{v_2}, \qquad {\Delta_i(\A_{k-2}', \A_{k+1}'')}:= x_{v_3},\\
\end{split}
\]
where $x_{v_i}$ are the corresponding standard monomials of $\mathscr{P}_{\G, t}^{(1)}$ in the same seed.
We have
\[
a_{k-1}^{\bf i} = x_{v_1}, \hskip 5mm a_{{\bf i},k}={\Delta_j(\A_{k-1}', \A_{k}'')}={\Delta_j(\A_{k-2}', \A_{k}'')} =x_{v_2}.
\]
The switching $(i, j, i)\leftrightarrow(j,i,j)$ from ${\bf i}$ to ${\bf i'}$ is equivalent to a mutation at $x_{v_1}$. 
Therefore 
\[
a_{{\bf i}, k+1}=a_{{\bf i}', k-1}=\frac{x_{v_2}+x_{v_3}}{x_{v_1}}= x_{v_2-v_1} + x_{v_3-v_1}.
\]
 
After quantization, we get
\[
{\rm A}_{{\bf i}, k-1} = \X_{v_1}, \hskip 5mm {\rm A}_{{\bf i},k}=\X_{v_2}, \hskip 5mm {\rm A}_{{\bf i}, k+1}= \X_{v_2-v_1}+ \X_{v_3-v_1}.
\]
Recall the Poisson bracket
$
\left\{ \log x_{v_i}, ~\log x_{v_j}\right\} = 2 (v_i, v_j).
$
\begin{lemma} \la{L11.16} We have
$
(v_1,v_2)=(v_3, v_1)=1/2, \hskip 5mm (v_2, v_3)=0.
$
\end{lemma}
\begin{proof} Recall the primary coordinates of $\mathscr{P}_{\G, t}$
\[
P_{k-1}=\frac{{\Delta_i(\A_{1}, \A_{k-1}'')}}{{\Delta_i(\A_{1}, \A_{k-2}'')}},
\hskip 7mm
P_k=\frac{{\Delta_j(\A_{1}, \A_{k}'')}}{{\Delta_j(\A_{1}, \A_{k-1}'')}},
\hskip 7mm
P_{k+1}=\frac{{\Delta_i(\A_{1}, \A_{k+1}'')}}{{\Delta_i(\A_{1}, \A_{k}'')}}.
\]
By Theorem \ref{TH9.5}, we have
\[
\left\{ \log P_{k-1}, \log P_{k}\right\}= \left\{ \log P_k, \log P_{k+1}\right\} = \left\{ \log P_{k+1}, \log P_{k-1}\right\} =1.
\]
Now let us consider 
\[
\begin{split}
&P_{k-1}'=\frac{{\Delta_i(\A_{k-2}', \A_{k-1}'')}}{{\Delta_i(\A_{k-2}', \A_{k-2}'')}}= x_{v_1}, \qquad P_k'=\frac{{\Delta_j(\A_{k-2}', \A_{k}'')}}{{\Delta_j(\A_{k-2}', \A_{k-1}'')}}= x_{v_2},\qquad P_{k+1}'=\frac{{\Delta_i(\A_{k-2}', \A_{k+1}'')}}{{\Delta_i(\A_{k-2}', \A_{k}'')}}= x_{v_3-v_1}.\\
\end{split}
\]
As illustrated below, we get these   coordinates by a sequence of moves from $(..., i, \overline{i},...)$ to $(..., \overline{i}, {i},...)$.
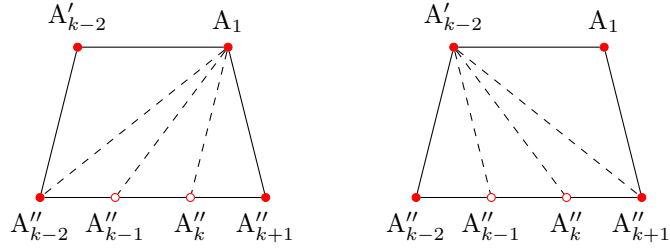
\begin{center}
\begin{tikzpicture}
 \node [circle,draw=red,fill=red,minimum size=3pt,inner sep=0pt,label=below:{\small ${\A}_{k-2}''$}] (a) at (0,0) {};
 \node [circle,draw=red,fill=white,minimum size=3pt,inner sep=0pt,label=below:{\small ${\A}_{k-1}''$}] (b) at (1,0) {};
 \node [circle,draw=red,fill=white,minimum size=3pt,inner sep=0pt,label=below:{\small ${\A}_{k}''$}] (c) at (2,0) {};
 \node [circle,draw=red,fill=red,minimum size=3pt,inner sep=0pt,label=below:{\small ${\A}_{k+1}''$}] (d) at (3,0) {};
  \node [circle,draw=red,fill=red,minimum size=3pt,inner sep=0pt,label=above:{\small ${\A}_{k-2}'$}] (a0) at (0.5,2) {};
 \node [circle,draw=red,fill=red,minimum size=3pt,inner sep=0pt,label=above:{\small ${\A}_{1}$}] (b0) at (2.5,2) {};
 \draw (a0) -- (b0) -- (d) -- (c) -- (b) --(a) -- (a0);
 \draw[dashed] (b0)--(a);
  \draw[dashed] (b0)--(b);
   \draw[dashed] (b0)--(c);
   \begin{scope}[shift={(5,0)}]
 \node [circle,draw=red,fill=red,minimum size=3pt,inner sep=0pt,label=below:{\small ${\A}_{k-2}''$}] (a) at (0,0) {};
 \node [circle,draw=red,fill=white,minimum size=3pt,inner sep=0pt,label=below:{\small ${\A}_{k-1}''$}] (b) at (1,0) {};
 \node [circle,draw=red,fill=white,minimum size=3pt,inner sep=0pt,label=below:{\small ${\A}_{k}''$}] (c) at (2,0) {};
 \node [circle,draw=red,fill=red,minimum size=3pt,inner sep=0pt,label=below:{\small ${\A}_{k+1}''$}] (d) at (3,0) {};
  \node [circle,draw=red,fill=red,minimum size=3pt,inner sep=0pt,label=above:{\small ${\A}_{k-2}'$}] (b0) at (0.5,2) {};
 \node [circle,draw=red,fill=red,minimum size=3pt,inner sep=0pt,label=above:{\small ${\A}_{1}$}] (a0) at (2.5,2) {};
 \draw (b0) -- (a0) -- (d) -- (c) -- (b) --(a) -- (b0);
 \draw[dashed] (b0)--(d);
  \draw[dashed] (b0)--(b);
   \draw[dashed] (b0)--(c);   
 \end{scope}  
\end{tikzpicture}
\end{center}
Following Corollary  \ref{local.primary} and the proof of Proposition \ref{switchi,ibar.}, we still have
\[
\left\{ \log P_{k-1}', \log P_k'\right\}= \left\{ \log P_k', \log P_{k+1}'\right\} = \left\{ \log P_{k+1}', \log P_{k-1}'\right\} =1.
\]
The Lemma is proved.
\end{proof}
Using  Lemma \ref{L11.16}, we get 
\be
\begin{split}
& \frac{q^{1/2}{\rm A}_{{\bf i}, k-1} {\rm A}_{{\bf i}, k+1}- q^{-1/2} {\rm A}_{{\bf i}, k+1}{\rm A}_{{\bf i}, k-1}}{q-q^{-1}}= \\
&\frac{q^{1/2}(q^{1/2} \X_{v_2}+ q^{-1/2}\X_{v_3})-q^{-1/2}(q^{-1/2} \X_{v_2}+q^{1/2}\X_{v_3})}{q-q^{-1}}= \X_{v_2} = {\rm A}_{{\bf i}, k+1}.\\
\end{split}
\ee
\end{proof}

\begin{lemma} \la{1.10.19.2} Recall the potential $\mathcal{W}_{2, i}$ of $\mathscr{P}_{\G, t}^{(1)}$.
If $s_{i_1}\ldots s_{i_{k-1}}(\alpha_k)=\alpha_i$ is a simple positive root, then 
 \[
\mathcal{W}_{2, i}= a_{{\bf i},k}:= \Delta_{i_k}(\A_{k-1}', \A_k'').
\]
\end{lemma}
\begin{proof} 
It is clearly true when $k=1$. The rest follows from  the identities  \eqref{z-tran0}. 
\end{proof}

  \paragraph{\bf Proof of Theorem \ref{REALI}.}

i)   Lemma \ref{1.10.19.2} tells that if  $s_{i_1}\ldots s_{i_{k-1}}(\alpha_k)=\alpha_i$ is a simple positive root, then 
\be \la{1.11.19.4}
\kappa({\bf E}_{{\bf i},k})= \kappa({\bf E}_i) = {\rm W}_{2, i}= {\rm A}_{{\bf i},k}.
\ee
In particular,    it holds for $k=1$.

We prove \eqref{EAa} by induction on $k$.  Assume it holds for $k <n$ for any ${\bf i}=(i_1, \ldots, i_N)$. 
Let $ w_{{\bf i}, l}=s_{i_1}\ldots s_{i_l}$.  To simply notation, we  denote  $w=w_{{\bf i}, n-2}$, $w'=w_{{\bf i}, n-1}$,   $w''= w_{{\bf i}, n}$, and $(i_{n-1}, i_n)=(i,j)$.
We divide the problem into the following   three cases.

\begin{itemize}
\vskip 1mm\item \underline{${\rm C}_{ij}=0$}.    
Let ${\bf i}'$ be  obtained from ${\bf i}$ by switching $i_{n-1}$ and $i_n$. 
Since $T_i\bE_j = \bE_j$, we have 
\be \la{INNd}
T_{w'}\bE_j = T_{w}T_i\bE_j = T_{w}\bE_j.
\ee
 Since $w'' = ws_js_i = ws_is_j$, we get
\be \la{eeeQ}
\kappa({\bE_{{\bf i},n}})\stackrel{\rm def}{=} \kappa(T_{w'}\bE_j)\stackrel{(\ref{INNd})}{=} \kappa(T_{w}\bE_j)\stackrel{\rm def}{=}  \kappa({\bE_{{\bf i}',n-1}})
\stackrel{\rm induction}{=} {\rm A}_{{\bf i}', n-1}\stackrel{ \ast}{=} {\rm A}_{{\bf i}, n}.
\ee
The last step is due to Lemma \ref{LLLa0}.
  
\vskip 2mm

\item \underline{${\rm C}_{ij}={\rm C}_{ji}=-1$ and $l(ws_{j})= l(w)-1$}. 
There is a reduced word $(i_1', \ldots, i_{n-2}')$ of $w$ with    $i_{n-2}'=j$. Let ${\bf i}'=(i_1',\ldots, i_{n-3}', j, i, j, i_{n+1}, \ldots)$. Applying the braid relation $s_is_js_i=s_js_is_j$ to ${\bf i}'$, we get  ${\bf i}''=(i_1',\ldots, i_{n-3}', i, j, i, \ldots)$. 
Therefore
\[
\kappa({\bE_{{\bf i},n}})=\kappa({\bE_{{\bf i}',n}})= \kappa({\bE_{{\bf i}'',n-2}})\stackrel{\rm induction}{=} {\rm A}_{{\bf i}'', n-2}{=}{\rm A}_{{\bf i}', n}= {\rm A}_{{\bf i}, n}.
\]
For the second equality, we note that ${\bE_{{\bf i}'',n-2}} = T_{w'}\bE_j$,  ${\bE_{{\bf i}',n}} = T_{w''}T_iT_j\bE_i$, and  $\bE_j=T_iT_j\bE_i$  due to Lusztig. 
The fourth equality  follows from \eqref{local.check220}. The last one holds since $a_{{\bf i}', n}=a_{{\bf i}, n}$.

\vskip 1mm\vskip 1mm\item \underline{${\rm C}_{ij}={\rm C}_{ji}=-1$ and $l(ws_{j})= l(w)+1$}.
We claim that  $l(ws_i) = l(w)+1$ and $l(ws_j) = l(w)+1$ imply
  that $l(ws_{i}s_j s_i)= l(w)+3$. Indeed, 
the first one is equivalent to $w\alpha_i>0$, and the second is equivalent to $w\alpha_j>0$. We have $l(ws_is_j) = l(ws_i)+1$ since $ws_i(\alpha_j) = w(\alpha_i+\alpha_j) = w(\alpha_i) +w(\alpha_j)>0$. Similar argument proves the claim. 

Without loss of generality, let us assume that $(i_{n-1}, i_{n}, i_{n+1})=(i,j,i)$.    Using the braid relation we move the $(n+1)^{st}$ generator  to the $(n-1)^{st}$.  
Then we get
\[
\kappa(\bE_{{\bf i},n-1})= {\rm A}_{{\bf i}, n-1},  \hskip 7mm \kappa(\bE_{{\bf i},n+1})= {\rm A}_{{\bf i}, n+1}.\]
  Now using similar arguments to the ones from previous steps, plus \eqref{local.check22}, we get 
\[\begin{split}
{\rm A}_{{\bf i}, n} &= \frac{q^{1/2}{\rm A}_{{\bf i}, n-1}{\rm A}_{{\bf i}, n+1}- q^{-1/2}{\rm A}_{{\bf i}, n+1}{\rm A}_{{\bf i}, n-1}}{q-q^{-1}} \\
&=\kappa \left( \frac{q^{1/2}\bE_{{\bf i},n-1} \bE_{{\bf i},n+1}- q^{-1/2}\bE_{{\bf i}, n+1} \bE_{{\bf i},n-1}}{q-q^{-1}} \right) \\
&= \kappa(\bE_{{\bf i},n}).\\
\end{split}
\]
\end{itemize}

ii) We use crucially the fact that ${\rm A}_{{\bf i}, n} = \kappa(\bE_{{\bf i},n})$. Let ${\bf m}:= (m_1, ..., m_N)\in {\Bbb N}^N$. Set
$$
\bE_{\bf i}^{\bf m}:= \bE^{m_1}_{{\bf i},1} \bE^{m_2}_{{\bf i},2}  \cdots \bE^{m_N}_{{\bf i},N}.
 $$
The set $ \{\bE_{\bf i}^{\bf m}~|~ {\bf m}\in {\Bbb N}^N\}$ is  a $\Q(q)$-linear  basis of   $\U_q(\mathfrak{n})$. It is Lusztig's PBW basis. Set
$$
{\bf A}_{\bf i}^{\bf m}:= \kappa(\bE_{\bf i}^{\bf m})={\rm A}^{m_1}_{{\bf i},1} {\rm A}^{m_2}_{{\bf i},2}  \cdots {\rm A}^{m_N}_{{\bf i},N}.
$$
Let us show that $\{{\bf A}_{\bf i}^{\bf m}\}$ is   $\Q(q)$-linearly independent. 
If not, there exists  relation  
\[f_1(q) {\bf A}_{\bf i}^{\bf m_1} + \cdots + f_k(q) {\bf A}_{\bf i}^{\bf m_k} =0,\hskip 7mm \mbox{where } f_i(q) \in \Q(q). 
\]  
By suitable rescaling, let us assume that  $f_i(q)$ are Laurent polynomials, and $f_i(1)\not = 0$ for at least one of $i$. 
By taking specialization at $q=1$, we get   
\be 
\la{rel.lin.anca}
f_1(1) \bar {\bf A}_{\bf i}^{\bf m_1} + \ldots + f_k(1) \bar {\bf A}_{\bf i}^{\bf m_k} =0,
\ee
where $\bar {\bf A}_{\bf i}^{\bf m_j}$ is a monomial of  $a_{{\bf i}, 1}, \ldots , a_{{\bf i}, N} \in {\cal O}({\rm Conf}_3({\cal A})_{e,e})$. 
The open embedding \eqref{open embedding.e}  provides an injective homomorphism from the polynomial ring $\Q[a_{{\bf i}, 1}, \ldots , a_{{\bf i}, N}]$ into $ {\cal O}({\rm Conf}_3({\cal A})_{e,e})$. 
Therefore $\{\bar {\bf A}_{\bf i}^{\bf m_j}\}$ are linearly independent, which contradicts \eqref{rel.lin.anca}. 
Theorem \ref{REALI} is proved. 

\medskip
\subsection{Proof of Theorem \ref{Th2.21}} \la{10.5}

\medskip

\noindent
 1)  This is a special case of Theorem \ref{UEAB},  proved in Section \ref{SEC12.1}.

 2)  
 Recall that the geometric coproduct $\Delta^{\rm g}$   is induced by the gluing of two punctured discs with a pair of special points on each.

\bl One has 
\be \la{383}  \begin{split}
 &\Delta^{\rm g} ({\bf K}_i) = {\bf K}_i \otimes  {\bf K}_i,\\
  &\Delta^{\rm g} ({\bf E}_i) = {\bf E}_i \otimes 1 + {\bf K}_i \otimes {\bf E}_i,\\
   &\Delta^{\rm g} ({\bf F}_i) = 1 \otimes {\bf F}_i + {\bf F}_i  \otimes {\bf K}^{-1}_i.\\
    \end{split}
  \ee
  \el
  \begin{proof} 
The first two formulas (\ref{383}) for the coproduct $\Delta^{\rm g}$ are dual to the amalgamation formulas (\ref{383a}) and \eqref{quantum.prom.def2}. The last is obtained similarly to the second one.
\end{proof}

3) Using (\ref{ACU}), the claim is evident for the counit, and easy to see for the antipode map. 

4)  Let us deduce it from Theorem \ref{REALI}. Let us cut a punctured disk along the dashed edge:
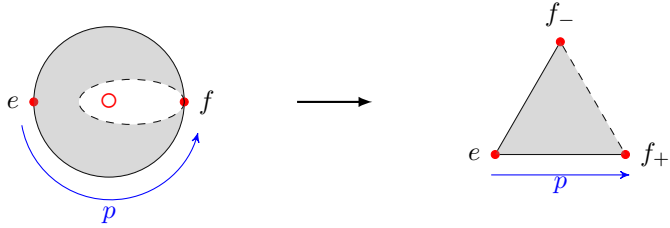
\begin{center}
\begin{tikzpicture}
\draw (0,0) circle (1cm);
\node[red, thick]  at (0,0) {\Large $\circ$}; 
\node [circle,draw=red,fill=red,minimum size=3pt,inner sep=0pt, label=left: {\small $e$}]  at (-1,0) {};
\node [circle,draw=red,fill=red,minimum size=3pt,inner sep=0pt, label=right: {\small $f$}] at (1,0) {};
\node [blue] at (0,-1.5) {\small  ${p}$};
\draw[-stealth', blue] (195:1.2) arc (190:345:1.2cm);
\draw[dashed] (0.3, 0) ellipse (0.7 and 0.3);   
\fill [gray,opacity=0.3] (1,0) arc (0:360:1) arc (360:0:0.7 and 0.3); 
\draw[-latex, thick] (2.5,0) -- (3.5,0);
\begin{scope}[shift={(6,-.2)}]
\node [circle,draw=red,fill=red,minimum size=3pt,inner sep=0pt, label=left: {\small $e$}] (a) at (210:1) {};
\node [circle,draw=red,fill=red,minimum size=3pt,inner sep=0pt, label=right: {\small $f_+$}] (b) at (330:1) {};
\node [circle,draw=red,fill=red,minimum size=3pt,inner sep=0pt, label=above: {\small $f_-$}] (c) at (90:1) {};
\node [blue] at (0,-.9) {\small  ${p}$};
\draw[-stealth', blue] (220:1.2) -- (320:1.2);
\draw (c) -- (a) -- (b);
\draw[dashed] (c) -- (b);
\fill [gray,opacity=0.3] (90:1) -- (210:1) -- (330:1)-- (90:1); 
\end{scope}
     \end{tikzpicture}
    \end{center}
The shadowed part together with its bottom pinning provides a map from $\mathscr{P}_{\G, \odot}$ to $\mathscr{P}_{\G, t}^{(1)}$. 

Fix a reduced word ${\bf i}=(i_1, \ldots, i_N)$, the functions $a_{{\bf i}, k}$ in \eqref{open embedding.e} can be pull back to $\mathscr{P}_{\G, \odot}$. 
In particular, assume that $i_N=i^\ast$ and set ${\bf i}'=(i, i_1, \ldots, i_N)$. By the very definition of  our geometric braid group action on $\mathscr{P}_{\G, \odot}$, for $k<N$, we get
\[
T({s_i})^* a_{{\bf i}, k} = a_{{\bf i}', k+1},\hskip 7mm \mbox{if } k<N, \hskip 10mm   T(s_i)^* a_{{\bf i}, N} = \mathcal{W}_{f, i^\ast} \mathcal{K}_{f, i^\ast}.
\]
Since $T(s_i)$ is a quasi-cluster transformation,  this formula  admits a quantum lift
\be
\la{classical.as}
T({s_i})^* {\rm A}_{{\bf i}, k} = {\rm A}_{{\bf i}', k+1};\hskip 7mm \mbox{if }k<N; \hskip 10mm   T(s_i)^* {\rm A}_{{\bf i}, N} =  q^{-1} {\rm K}_{f, i^\ast} {\rm W}_{f, i^\ast}.
\ee
Plugging \eqref{EAa} into \eqref{classical.as}, for $k<N$, we get
\[
T({s_i})^* \kappa({\bE}_{{\bf i}, k}) = \kappa({\bE}_{{\bf i}', k+1}) = \kappa(T_{s_i} {\bE}_{{\bf i}, k}).
\]
Note that $\bE_{{\bf i}, N}= \bE_i$. For $k=N$, we get
\[
T({s_i})^* \kappa({\bE}_{{\bf i}, N}) = q^{-1} {\rm K}_{f, i^\ast} {\rm W}_{f, i^\ast} =\kappa (T_{s_i} \bE_{{\bf i}, N}).
\]
The same argument works for the generators ${\bf F}_i$. For the generators ${\bf K}_i$, it is due to the quantization of the last Formula in Lemma \ref{L2.8.19.1}.

Let us  show  the map $\kappa$ is injective. By Theorem \ref{4}, we get an isomorphism
\[
\mathcal{L}_{\G, \odot}={\rm G}^\ast = {\rm U}\times {\rm U}^- \times {\rm H}\stackrel{\sim}{\lra} \mathbb{A}^{2(l(w_0))}\times \mathbb{G}_m.
\]Recall the rescaled PBW basis $\{\bE_{\bf i}^{{\bf m}} {\bf K}^{\nu} {\bf F}_{\bf i}^{\bf n}\}$of $\U_q(\mathfrak{g})$.  Applying   decomposition \eqref{natural.dcom.u.zz} to unipotent elements associated to both special points,   the   $q=1$ specialization of $\kappa (\bE_{\bf i}^{{\bf m}} {\bf K}^{\nu} {\bf F}_{\bf i}^{\bf n})$  provides a linear basis for $\mathcal{O}(\mathcal{L}_{\G, \odot})$.  Following the   argument in the proof of Theorem  \ref{REALI}, we prove the injectivity of $\kappa$. 
Theorem \ref{Th2.21} is proved.

 \medskip

\subsection{An example: $\G = {\rm PGL}_2$} \la{SEC12.4}

\medskip

The quiver below describes a  cluster for the quantum space  ${\mathscr P}_{{\rm PGL_2}, \odot}$,  given by  a quantum torus algebra with the generators $\X_{\pm e_i}$ and relations $\X_{e_i} \X_{e_{i+1}} = q^2\X_{e_{i+1}}\X_{e_i}$, $i \in \Z/4\Z$.

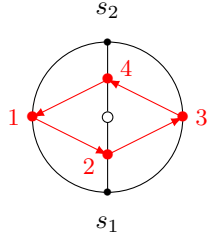
\begin{center}
\begin{tikzpicture}[scale=1]
\draw (0,0) circle (1cm);
\node[circle,draw, fill=white, minimum size=4pt,inner sep=0pt] (a) at (0,0) {};
\draw (0,1) -- (a); 
\draw (0,-1) -- (a); 
\node[label=above: {\small $s_2$}] at (0,1) {{\tiny $\bullet$}};
\node[label=below: {\small $s_1$}] at (0,-1) {{\tiny $\bullet$}};
\node[red] (a) at (0,-.5) {{\small $\bullet$}};
\node[red] (c) at (0,.5) {{\small $\bullet$}};
\node[red] (b) at (-1,0) {{\small $\bullet$}};
\node[red] (d) at (1,0) {{\small $\bullet$}};
\draw[red, directed] (-1,0)--(0,-.5);
\draw[red, directed] (1,0)--(0,.5);
\draw[red, directed] (0,-.5)--(1,0);
\draw[red, directed] (0,.5)--(-1,0);
\node[red] at (-.25,-.65) {{\footnotesize 2}};
\node[red] at (.25,.65) {{\footnotesize 4}};
\node[red] at (-1.25,0) {{\footnotesize 1}};
\node[red] at (1.25,0) {{\footnotesize 3}};
\end{tikzpicture}
 \end{center}

 The Cartan elements ${\rm K}_{s_i}$ at the special points $s_1, s_2$ are given by 
\[
 {\rm K}_{s_1} = \X_{e_1+e_2+e_3}, \qquad  {\rm K}_{s_2} = X_{e_1+e_3+e_4}.
\]

The monodromy around the puncture is   $\mu = \X_{e_2+e_4}$.

There is an action of the Weyl group $W=\Z/2\Z$, whose  generator $\sigma$ acts by 
\[
 \sigma:~\qquad~\X_{e_2} \lms \X_{-e_4}, \qquad\X_{e_4} \lms \X_{-e_2}, \qquad\X_{e_1} \lms \X_{e_1}\frac{(1+\X_{e_4})}{(1+q\X_{-e_2})}, \qquad\X_{e_3} \lms \X_{e_3} \frac{(1+q\X_{e_2})}{(1+q\X_{-e_4})}.
 \] 
 
 The generator $\beta$ of the braid group acts as follows: 
 $$
 \beta: \qquad \X_{e_1} \lms \X_{-e_1-e_4},   ~\qquad~ \X_{e_2} \lms \X_{e_4}, ~\qquad~ \X_{e_3} \lms \X_{-e_2-e_3},~\qquad~\X_4 \lms \X_2. 
 $$
Note that $\beta^2$ acts trivially on  unfrozen  variables, but on the   frozen ones  by 
$$
\beta^2: \qquad \X_1 \lms \X_1 {\rm K}_{s_2}{\rm K}^{-1}_{s_1},~\qquad~\X_3 \lms \X_3{\rm K}_{s_1}{\rm K}^{-1}_{s_2}.
$$ 
 
 The reduced quantum subspace $\widetilde {\mathscr L}_{{\rm PGL_2}, \odot}$ is described by the relation 
\be  \la{CENT}
 {\rm K}_{s_1} {\rm K}_{s_2} = \X_{2e_1+e_2+2e_3+e_4} = 1.
\ee
Note that   ${\rm K}_{s_1} {\rm K}_{s_2}$ belongs to the center. So the equation ${\rm K}_{s_1} {\rm K}_{s_2}=1$  makes sense.
 The Weyl group   acts by  $\sigma: {\rm K}_{s_i}   \mapsto  {\rm K}_{s_i}$.\footnote{Geometrically this is clear: altering the invariant flag  at the puncture we do not touch the pinnings, and thus  
  do not change the Cartan elements ${\rm K}_{s_i}$.}   The braid group   acts by $\beta: {\rm K}_{s_i} \mapsto {\rm K}^{-1}_{s_i}$. 
    So both actions  descend  to ${\cal O}_q(\widetilde {\mathscr L}_{{\rm PGL_2}, \odot})$.   
 
 \vskip 2mm
  The quantum group ${\U}_q(sl_2)$ is realized  in the  subalgebra ${\cal O}_q({\mathscr L}_{{\rm PGL_2}, \odot}):= {\cal O}_q(\widetilde {\mathscr L}_{{\rm PGL_2}, \odot})^W$ as follows:
 \be \la{REALQU}
 \begin{split} 
&\kappa: {\U}_q(sl_2) \lra {\cal O}_q({\mathscr L}_{{\rm PGL_2}, \odot}),\\
& {\bf E} \lms \X_{e_3} + \X_{e_2+e_3}, ~\qquad~{\bf F}\lms \X_{e_1} + \X_{e_1+e_4}, ~\qquad~{\bf K} \lms \X_{e_1+e_2+e_3}.\\
\end{split}
\ee
 Let    $\widetilde \mu = \X_{e_1+e_3}$. Then $\mu =  \X_{e_2+e_4}$. So $\mu = \widetilde \mu^{-2}$, thanks to the relation (\ref{CENT}). The center of  ${\cal O}_q(\widetilde {\mathscr L}_{{\rm PGL_2}, \odot})$ is generated by  $\widetilde \mu$. 
  The center of   ${\cal O}_q({\mathscr L}_{{\rm PGL_2}, \odot})$ is generated by the quantum Laplacian 
  $$
  \Delta = \frac{1}{2}\left({\bf E}{\bf F} + {\bf F}{\bf E}  - (q+q^{-1})({\bf K} + {\bf K}^{-1})\right)= \X_{e_1+e_3} + \X_{-e_1-e_3} = \widetilde \mu + \widetilde \mu^{-1}.
  $$ 
Therefore
  $ \mu + \mu^{-1}$ generates an index two subalgebra of the center of the algebra ${\cal O}_q({\mathscr L}_{{\rm PGL_2}, \odot})$.  \vskip 2mm
  
  Note that the variable $q$ which we use  in the cluster set-up, and in particular in Section \ref{SEC12.4},  match the $q$ in the quantum group set-up  in Section \ref{SSECC11.1}. 
  For example, we have ${\bf K}{\bf E}{\bf K}^{-1} = q^2{\bf E}$ and ${\bf K}{\bf F}{\bf K}^{-1} = q^{-2}{\bf F}$. 
 
 \vskip 2mm
 Realization (\ref{REALQU}) coincides with Faddeev's realization \cite[ (27)-(28)]{Fa2} of the quantum group ${\U}_q(sl_2)$  after a monomial transformation  
 $\omega_1 = \X_{e_3}, ~\omega_2 = \X_{e_2+e_3}, ~\omega_3 = \X_{e_1}, ~\omega_4 =\X_{e_1+e_4}$.  
   

  \medskip          
  
   \section{Realizations and quasiclassical limits of the principal series}  \la{Sec8}
 
  \medskip
  
      Below $\G$ is a split semi-simple  adjoint group,    $\G'$   its universal cover, 
 and  $\H$  and $ \H'$ are the Cartan groups of $\G$  and     $\G'$, respectively. 
  The kernel of the    isogeny 
$\psi:  \H' \rightarrow \H$  is identified with the center  $\G'$.  The dual map of   character lattices is    the 
embedding   of the root lattice into the weight lattice. 

In Section \ref{Sec8} we use extensively the quantum symplectic double construction, which is reviewed in Section \ref{Sec8.1}. The reader should consult Section \ref{Sec8.1} before reading this Section. 

\subsection{Cluster symplectic double and the Poisson moduli space for the   double of $\bS$} 

 \medskip

 \vskip 1mm
Let $\bS$ be a decorated surface without punctures, and $\bS^\circ$   the same  surface    with the opposite orientation. 
Gluing     $\bS$ and $\bS^\circ$   along the matching boundary intervals we get a surface $\bS_{\cal D}$ - 
the topological double of $\bS$. 
It comes with a set of { punctures} $\{p_1, \ldots , p_n\}$ matching the special points of $\bS$.  
Namely, each  oriented boundary interval ${\rm I_k}$ of $\bS$ gives rise to a    puncture $p_{k}$ on the surface $\bS_{\cal D}$,  corresponding to the endpoint of 
 ${\rm I_k}$.    

\vskip 1mm
There are  the following three moduli spaces / cluster varieties:

   \vskip 2mm 
  
\begin{enumerate}

\item   The   cluster Poisson variety ${\mathscr P}_{\G, \bS}$.    \vskip 1mm  

\item  The  cluster symplectic double ${\cal D}_{\G, \bS}$ of the cluster   Poisson variety ${\mathscr P}_{\G, \bS}$. \vskip 1mm  
\item   The moduli space $ {\mathscr P}_{  \G, \bS_{\cal D}}$ of  framed $\G$-local systems 
  on the topological   double $\bS_{\cal D}$.    
  
  \end{enumerate}
  \vskip 2mm 

Let us fix a quiver ${\bf q}$ providing   cluster coordinate systems $\{Y_i, \B_i\}$, $\{X_i\}$,   $\{X_{i^\circ}\}$ for the spaces $ \mathcal{D}_{\G, \bS}$, $\mathscr{P}_{\G, \bS}$,   $\mathscr{P}_{\G, \bS^\circ}$ respectively. 
We denote by ${\rm V}_{\G, \bS}$ the set parametrising the vertices of the quiver ${\bf q}$. 

To distinguish the   cluster Poisson coordinates on   $\mathscr{P}_{\G, \bS}$   from the ones on   $ {\mathscr P}_{  \G, \bS_{\cal D}}$, we denote the latter by  
$\{{\rm Z}_f, {\rm Z}_j, {\rm Z}_{j^\circ}\}$, where $\{{\rm Z}_f\}$ are the frozen variables matching the ones $\{X_f\}$ on ${\mathscr P}_{\G, \bS}$, 
$\{{\rm Z}_j\}$ match  the unfrozen variables $\{X_j\}$ on ${\mathscr P}_{\G, \bS}$, and $\{{\rm Z}_{j^\circ}\}$ match   the    unfrozen  ones $\{X_{j^\circ}\}$ on ${\mathscr P}_{\G, \bS^\circ}$, see Figure \ref{amalgam2}. 

Any cluster Poisson transformation $c^{\mathscr X}$ of $\mathscr{P}_{\G, \bS}$ gives rise naturally to a cluster Poisson transformation $c^{\cal Z}$ of     $ {\mathscr P}_{  \G, \bS_{\cal D}}$ and 
cluster symplectic double transformation $c^{\cal D}$ of   $ {\cal D}_{  \G, \bS}$. See Section \ref{SEC9.1} for the details.

 \bd  \la{16.1}  Given a decorated surface $\bS$ with boundary intervals $\{{\rm I}_k\}$, $k=1, ..., n$, the subalgebra 
 $$
 {\bf B }_{\rm I_k}\subset {{\cal O}_q(\cal D}_{\G, \bS})
 $$
  is generated by the frozen variables at  the boundary interval ${\rm I_k}$.  
  These subalgebras generate the subalgebra 
$$
{\bf B } = \bigotimes_{1 \leq k \leq n}{\bf B}_{\rm I_k} \subset {{\cal O}_q(\cal D}_{\G, \bS}).
$$
The centraliser of the subalgebra ${\bf B}$ in the algebra ${{\cal O}_q(\cal D}_{\G, \bS})$ is the subalgebra
\be \la{SUBQSD} 
{\cal O}_q({\cal D}_{\G, \bS})^{\bf B }\subset {{\cal O}_q(\cal D}_{\G, \bS}).
\ee
  \ed
The subalgebra ${\bf B }$ is the center   ${\cal O}_q({\cal D}_{\G, \bS})^{\bf B }$ for generic $q$.  One has a canonical    isomorphism 
\be \la{CaniB}
{\bf B }_{\rm I_k} = {\cal O}({\H'}).
\ee   

\vskip 2mm
 
 Our  goal is to compare the quantized algebra  ${\cal O}_q({\mathscr P}_{\G, \bS_{\cal D}} )$ for the topological double $\bS_{\cal D}$ of the decorated surface $\bS$, and 
 the subalgebra (\ref{SUBQSD}) of the quantum cluster symplectic double algebra ${{\cal O}_q(\cal D}_{\G, \bS})$. Let us start from the geometric counterpart of the subalgebra 
 ${\bf B}$. 
 \vskip 2mm

The  monodromy around the   puncture $p_{k}$ provides a  map 
 $ \mu_{p_{k}}: {\mathscr P}_{\G, \bS_{\cal D}} \lra \H$.  
 
 The  center of  the algebra ${\cal O}_q({\mathscr P}_{  \G, \bS_{\cal D}})$ for generic $q$ is given by  
 $
 \bigotimes_{ k } \mu_{p_{k}}^*{\cal O}_{\H}. 
 $ 

Following    (\ref{CORT}), we set
\be \la{CORTa}
\mathbb{B}_{k}:=  \prod_{j\in{\rm I}} B_{j}^{-C_{kj}}.
\ee

 \bt  \la{11.26.18.1} Let $\bS$ be a decorated surface without punctures. 
Then there is an  injective  map  
\be \la{zeta}
\xi^*:     {\cal O}_q({\mathscr P}_{\G, \bS_{\cal D}} ) \hra {\cal O}_q({\cal D}_{\G, \bS})^{\bf B }. 
\ee
 \vskip 2mm
 
  \begin{enumerate}
\item 
 It is given in   each cluster coordinate system     by the following  monomial map:
  \begin{equation} \label{12.12.xcv}
\xi_{\bf i}^*: {\rm Z}_i \lms
\left\{ \begin{array}{lll} {\rm Y}_i
& \mbox{ if $i$ is unfrozen, from from $\bS$},\\
\widetilde {\rm Y}_i& \mbox{ if $i$ is unfrozen, from $\bS^\circ$},\\
{\Bbb B}^{-1/2}_i \prod_{j \in {\rm V}_{\G, \bS}}\B_j^{-\varepsilon_{ij}}& \mbox{ if $i$ is frozen.}\\
\end{array}\right.
\end{equation}
It  intertwines  cluster transformations ${\bf c}^{\cal Z}$ and ${\bf c}^{\cal D}$ generated by unfrozen   mutations. 
  \vskip 1mm

\item  
For each puncture $p_{k}$, the map $\xi^*$ induces a canonical $W-$equivariant  injective  map  
\be \la{MONo}
  \mu_{p_{k}}^*{\cal O}_{\H} \stackrel{ }{\hlra}  {\bf B_{{\rm I}_k}} .
\ee
   Via  the isomorphism  ${\bf B_{{\rm I}_k}}  = {\cal O}_{{\H'}}$, see (\ref{CaniB}), it is induced by the isogeny $\psi: {\H'} \lra \H$. 
\end{enumerate}
 \et
 
  Theorem \ref{11.26.18.1}  
   is the precise form of Theorem \ref{11.26.18.1Xaaa}. In Section \ref{Sec12.2n} we elaborate its simplest example. 
   In Section \ref{SEC9.1} we  prove its general cluster counterpart,  Theorem  \ref{THH8.9}. Theorem \ref{11.26.18.1} is proved in Section \ref{SEC9.2}.
   In Sections \ref{sec16.5:RMSP},  \ref{SEC9.4},  
   we give  the main applications to the cluster representation theory of quantum groups: 
   realizations of the multiplicity spaces,  and 
    realizations of the principal series of representations.

\medskip

    \subsection{An example: $\G = {\rm PGL}_2$ and $\bS$ is a polygon} \la{Sec12.2n}
 
 \medskip
 
  Let  $\bS = {\Bbb P}_n$ be an $n-$gon.   Its  double $\bS_{\cal D}$ is a sphere  with $n$ punctures 
 $S^2 - \{p_1, ..., p_n\}$ $\&$ an ideal $n-$gon ${\mathscr P}_n$. Below we restrict to the simplest case $n=4$ which has all the features of the general case. 
 
 Pick a diagonal $E$   of the square ${\Bbb P}_4$. 
 Then its double     is triangulated 
by  the  sides 
 of   ${\mathscr P}_4$ matching the sides of the square, and  the  edges  $E,  E^\circ$ corresponding to the diagonal   of ${\Bbb P}_4$, see the middle  of Figure \ref{amalgam2}. 

 In each of the  algebras below, the corresponding cluster is generated by the following variables: 
\vskip 2mm
 
 ${\cal O}_q({\mathscr P}_{{\rm PGL}_2,  \bS_{\cal D}})$:  $6$  quantum variables 
  $({\rm Z}_1, ..., {\rm Z}_4, {\rm Z}_E, {\rm Z}_{ E^\circ})$ 
  shown in the middle of Figure \ref{amalgam2}.    \vskip 1mm
    
    ${\cal O}_q({\cal D}_{{\rm PGL}_2,  \bS})$:  $ \ 10$ quantum   variables  
  $({\rm Y}_1, \B_1, ..., {\rm Y}_4, \B_4, {\rm Y}_E,  \B_E)$, and  the subalgebra ${\bf B}=\Z[\B_1, ..., \B_4]$.   \vskip 1mm
    
  ${\cal O}_q({\cal D}_{{\rm PGL}_2, \bS})^{\bf B }$:   $6$  quantum variables 
  $(\B_1, \B_2, \B_3, \B_4, {\rm Y}_E, {\rm Y}_{E^\circ})$, shown on the right of Figure \ref{amalgam2}.    
  
  \vskip 1mm
  \vskip 1mm

   \begin{figure}[ht]
   \begin{center}
\begin{tikzpicture}
\draw[red] (-1,-1) -- (-1,1) -- (1,1) -- (1,-1) -- (-1,-1);
\draw (-1,-1)--(1,1);
\node[red] at (90:1.5) {${\rm X}_1$}; 
\node[red] at (0:1.5) {${\rm X}_2$}; 
\node[red] at (270:1.5) {${\rm X}_3$}; 
\node[red] at (180:1.5) {${\rm X}_4$}; 
\node at (.3,-.3) {${\rm X}_E$}; 
\begin{scope}[shift={(5,0)}]
\node at (1,1) {\footnotesize $\bullet$};
\node at (-1,1) {\footnotesize $\bullet$};
\node at (1,-1) {\footnotesize $\bullet$};
\node at (-1,-1) {\footnotesize $\bullet$};
\draw (-1,1) .. controls (0,1.2) .. (1,1);
\draw (1,1) .. controls (1.2,0) .. (1,-1);
\draw (1,-1) .. controls (0,-1.2) .. (-1,-1);
\draw (-1,-1) .. controls (-1.2,0) .. (-1,1);
\draw[dashed] (-1,-1) .. controls (-.6,0) and (0,.6) .. (1,1);
\draw (-1,-1) .. controls (0,-.6) and (.6,0) .. (1,1);
\node[blue] at (1.3,1.3) {$p_1$};
\node[blue] at (-1.3,1.3) { $p_4$};
\node[blue] at (1.3,-1.3) { $p_2$};
\node[blue] at (-1.3,-1.3) {$p_3$};
\node at (90:1.5) {${\rm Z}_1$}; 
\node at (0:1.5) {${\rm Z}_2$}; 
\node at (270:1.5) {${\rm Z}_3$}; 
\node at (180:1.5) {${\rm Z}_4$}; 
\node at (.5,-.5) {${\rm Z}_E$}; 
\node at (-.5,.5) {${\rm Z}_{E^\circ}$}; 
\end{scope}
\begin{scope}[shift={(10,0)}]
\node at (1,1) {\footnotesize $\bullet$};
\node at (-1,1) {\footnotesize $\bullet$};
\node at (1,-1) {\footnotesize $\bullet$};
\node at (-1,-1) {\footnotesize $\bullet$};
\draw (-1,1) .. controls (0,1.2) .. (1,1);
\draw (1,1) .. controls (1.2,0) .. (1,-1);
\draw (1,-1) .. controls (0,-1.2) .. (-1,-1);
\draw (-1,-1) .. controls (-1.2,0) .. (-1,1);
\draw[dashed] (-1,-1) .. controls (-.6,0) and (0,.6) .. (1,1);
\draw (-1,-1) .. controls (0,-.6) and (.6,0) .. (1,1);
\node[blue] at (1.3,1.3) {$p_1$};
\node[blue] at (-1.3,1.3) { $p_4$};
\node[blue] at (1.3,-1.3) { $p_2$};
\node[blue] at (-1.3,-1.3) {$p_3$};
\node at (90:1.5) {${\rm B}_1$}; 
\node at (0:1.5) {${\rm B}_2$}; 
\node at (270:1.5) {${\rm B}_3$}; 
\node at (180:1.5) {${\rm B}_4$}; 
\node at (.5,-.5) {${\rm Y}_E$}; 
\node at (-.5,.5) {${\rm B}_{E}$}; 
\end{scope}
\end{tikzpicture}
\end{center}
\caption{On the left:   cluster   coordinates for the square ${\Bbb P}_4$, the frozens are red. In the middle:    cluster coordinates for $S^2 - \{p_1, ..., p_4\}$. On the right:   coordinates for $S^2 - \{p_1, ..., p_4\}$ from   the double of ${\Bbb P}_4$.}
\label{amalgam2}
\end{figure}

The map $\xi^*$ in (\ref{zeta}) is given by:
   \be
   \begin{split}
 &{\rm Z}_1 \lms \B_1\B_4\B_E^{-1}, \qquad   {\rm Z}_2 \lms \B_2\B_3^{-1}\B_E, \qquad  {\rm Z}_3 \lms \B_3\B_2\B_E^{-1}, \qquad   {\rm Z}_4 \lms \B_4\B_1^{-1}\B_E,\\
 &{\rm Z}_E \lms {\rm Y}_E, \qquad {\rm Z}_{E^\circ} \lms  {\rm Y}^{-1}_E \B_1\B^{-1}_2\B_3\B^{-1}_4. \\
 \end{split}
   \ee
The monodromy $\mu_{p_i}$ around the puncture $p_i$  is given in  $Z-$coordinates by the product of  coordinates on the edges sharing the puncturte, e.g.  $\mu_{p_1} = {\rm Z}_1{\rm Z}_2{\rm Z}_E{\rm Z}_{E^\circ}$. In the $\B-$coordinates  it  is given by $\B_i^2$, that is:
$$
\xi^*(\mu_{p_i}) = \B_i^2.
$$ 
The $\B-$variables are related to the Cartan group $\H_{{\rm SL}_2}$ of ${\rm SL}_2$. The ${\rm Z}-$variables are related to the 
Cartan group $\H_{{\rm PGL}_2}$ of ${\rm PGL}_2$. The map $\H_{{\rm SL}_2} \to \H_{{\rm PGL}_2}$ is given by $t \lms t^2$. This is why we get $\B_i^2$. 

 \subsection{Quantum double of a cluster Poisson variety with frozen variables}  \la{SEC9.1}
 
 \medskip
 
 Let ${\bf q}$ be a quiver. Let ${\rm V} = {\rm V}^u \cup  {\rm V}^f$ be the set of vertices of ${\bf q}$, where ${\rm V}^f$  consists of  frozen ones. 
The chiral dual  quiver ${\bf q}^\circ$ is obtained from ${\bf q}$ by reversing arrows.   Amalgamating  ${\bf q}$ and ${\bf q}^\circ$  along the subset ${\rm V}^f$, we get a {\it double quiver} 
 with the set of   vertices  ${\rm V}^d = {\rm V}^u \cup  {\rm V}^f \cup  \widetilde {\rm V}^u$:
   $$
{\bf q}_{\bf d}:=  {\bf q}^\circ \cup_{{\rm V}^f} {\bf q}.
$$
There are several cluster varieties associated with the quivers ${\bf q}$ and  ${\bf q}_{\bf d}$:
\vskip 2mm
\begin{enumerate}
 
\item   The cluster Poisson variety ${\mathscr X}_{{\bf q}}$ associated with the  quiver  ${\bf q}$.  
It has  a cluster Poisson coordinate system $\{X_i, X_f\}$, where 
 $i \in {\rm V}^u, f \in {\rm V}^f$. Similarly, we have  ${\mathscr X}^\circ_{{\bf q}} = {\mathscr X}_{{\bf q^\circ}}$.\vskip 1mm
 
\item   The cluster Poisson variety ${\mathscr X}_{{\bf q}_{\bf d}}$ associated with the  quiver ${\bf q}_{\bf d}$. 
 It has a cluster Poisson coordinate system $\{{\rm Z}_i, {\rm Z}_{i^\circ}, {\rm Z}_f\}$, where 
 $i \in {\rm V}^u, f \in {\rm V}^f$. \vskip 1mm
 
 \item    The cluster symplectic double\footnote{See Theorem \ref{8.15.05.10} and \cite{FG07} for the definition of symplectic double. To distinguish from the coordinates $X_i$ of cluster Poisson varieties, we shall use $Y_i$ for the coordinates of symplectic doubles.} ${\cal D}_{{\bf q}}$ of ${\mathscr X}_{{\bf q}}$. It has a cluster   coordinate system $\{Y_i, Y_f, B_{i}, B_f\}$.\vskip 1mm
 
  \item A  new cluster variety ${\cal D}^{\bf b}_{\bf q}$,   defined below.   It has a cluster   coordinate system $\{Y_i, B_{i}, B_f\}$. 
\end{enumerate}
 \vskip 2mm

Let  ${\bf i}$ be a quiver  obtained from the initial quiver ${\bf q}$ by mutations.  Recall the quantum torus algebra ${\cal O}_q({\cal D}_{\bf i})$ assigned to ${\bf i}$, see \eqref{4.28.03.11x}. 
 Let ${\bf B}$ be a commutative subalgebra of   ${\cal O}_q({\cal D}_{\bf i})$ generated by the frozen variables ${\rm B}_f$, $ f \in {\rm V}^f$.    The centralizer  ${\cal O}_q({\cal D}_{\bf i})^{\bf B }$ of ${\bf B }$ is a quantum torus algebra generated by $\{{\rm Y}_i, {\rm B}_{i}, {\rm B}_f\}$. Note that ${\bf B }$ is the center  of ${\cal O}_q({\cal D}_{\bf i})^{\bf B }$. 
 
 The quantum cluster transformations are generated by  mutations at the unfrozen subset ${\rm V}^u$. They act trivially on the variables ${\rm B}_f$, and hence 
 on the subalgebra ${\bf B}$. Therefore they provide isomorphisms of the fraction fields of   algebras  ${\cal O}_q({\cal D}_{\bf i})^{\bf B }$. 
 
 \bd  \la{12.28.18.1} The quantum space ${\cal D}_{\bf q}^{\bf b}$ is obtained by gluing the quantum  torus algebras ${\cal O}_q({\cal D}_{\bf i})^{\bf B }$  
 by the quantum cluster transformations. Its algebra of regular functions is  
 $$
 {\cal O}_q({\cal D}_{\bf q}^{\bf b}) = {\cal O}_q({\cal D}_{\bf q})^{\bf B}.
  $$
  In particular, setting $q=1$ we get a  cluster variety ${\cal D}_{\bf q}^{\bf b}$.   
  \ed
  
  The variety ${\cal D}_{\bf q}^{\bf b}$ is a quotient space of  ${\cal D}_{\bf q}$. The quotient map
from ${\cal D}_{\bf q}$ to   ${\cal D}_{\bf q}^{\bf b}$
  is a Poisson map. The Hamiltonian flows of  Hamiltonians $B_f$ act on the fibers, generating them. 
   
   \vskip 2mm
 Every  $k \in {\rm V}^{u}$
gives rise  to mutations of different flavors, both classical and quantum: 

\begin{enumerate}

\item   A mutation $\mu^{\mathscr X}_k$ of the cluster Poisson variety ${\mathscr X}_{{\bf q}}$. \vskip 1mm
\item    A cluster Poisson transformation $\mu_k^{\cal Z}:= \mu^{\mathscr X}_k\circ \mu^{\mathscr X}_{k^\circ} $ of  ${\mathscr X}_{{\bf q}_{\bf d}}$. Note that  $\mu^{\mathscr X}_k $ and $\mu^{\mathscr X}_{k^\circ} $ commute. \vskip 1mm

\item   A mutation $\mu_k^{\cal D}$ of the cluster symplectic double   ${\cal D}_{{\bf q}}$, as well as of ${\cal D}^{\bf b}_{{\bf q}}$. 
 \end{enumerate}

\vskip 2mm
 
A  cluster transformation ${\bf c}^{\mathscr X}$ of ${\mathscr X}_{{\bf q}}$  is a sequence of mutations at unfrozen vertices. It defines:

  \begin{enumerate}
   
\item  A cluster transformations ${\bf c}^{\cal Z}$ of the cluster Poisson variety ${\mathscr X}_{{\bf q}_{\bf d}}$ and its q-deformation.   \vskip 1mm

\item  A cluster transformation ${\bf c}^{\cal D}$ of the symplectic double 
  ${\cal D}_{{\bf q}}$,  its quotient ${\cal D}^{\bf b}_{{\bf q}}$, and their q-deformations. 
 
\end{enumerate}

\vskip 2mm

Let us set 
\[
\widetilde {\rm Y}_i:= {\rm Y}_i^{-1}\prod_{j \in {\rm V}}{\rm B}_j^{-\varepsilon_{ij}}. 
\]  
Recall that $\varepsilon_{ff'} \in \frac{1}{2}\Z$ if $f,f' \in {\rm V}^f$.  We assign to each $f  \in {\rm V}^f$ an arbitrary   monomial  in frozen variables with  exponents in $\frac{1}{2}\Z$:
\be \la{12.28.18.5}
 {\Bbb B}_f :=\prod_{k \in {\rm V}^f}{\rm B}^{-\beta_{fk}}_{k}, \qquad \beta_{fk}\in \frac{1}{2}\Z,
\ee   
 such that  ${\Bbb B}_f\prod_{j \in {\rm V}}{\rm B}_j^{-\varepsilon_{fj}}$ is a monomial with integral exponents. 
So there is an integral ${\rm V}\times {\rm V}$  matrix  $p=(p_{ij})$   such that 
 $$
 \prod_{j \in {\rm V}}{\rm B}_j^{-p_{fj}} = {\Bbb B}_f\prod_{j \in {\rm V}}{\rm B}_j^{-\varepsilon_{fj}}.
   $$
So  $p_{ij}= \varepsilon_{ij}+\beta_{ij}$ if $i, j \in {\rm V}^f$, and $p_{ij}= \varepsilon_{ij}$ otherwise.

 \bt \la{THH8.9}There is a    map of   algebras:
$$
\xi^*:     {\cal O}_q({\mathscr X}_{{\bf q}_{\bf d} } ) \lra {\cal O}_q({\cal D}_{ {\bf q}})^{\bf B }, 
$$
induced in   each cluster coordinate system     by the following  monomial map:
  \begin{equation} \label{12.12.xcvr}
\xi^*: {\rm Z}_v \lms
\left\{ \begin{array}{lll} {\rm Y}_i
& \mbox{ if $v=i\in {\rm V}^u$},\\
\widetilde {\rm Y}_i& \mbox{ if $v = i^\circ$, $i \in  {\rm V}^u$},\\
 \prod_{j \in {\rm V}}{\rm B}_j^{-p_{fj}}& \mbox{ if $v=f\in {\rm V}^f$.}\\
\end{array}\right.
 \end{equation} 
 It intertwines  cluster transformations ${\bf c}^{\cal Z}$ and ${\bf c}^{\cal D}$ generated by  unfrozen   mutations. 
 
 If the integral matrix $p$ is of full rank, then the map $\xi^*$ is injective. 
   \et

 \begin{proof} The injectivity of $\xi^*$ follows by   easy   linear algebra.
We present a proof of  the rest by  calculation. There is a less computational proof using the quantum variant of the comment after Theorem \ref{8.15.05.10}. 
 
 If $i, j\in {\rm V}^u$, then   ${\rm Z}_i$   commutes with   ${\rm Z}_{j^\circ}$, and ${\rm Y}_i$ commutes with  $\widetilde {\rm Y}_j$. Therefore 
 the $\xi^*$ is an algebra homomorphism when restricted on the subalgebra generated by ${\rm Z}_i$ and   ${\rm Z}_{j^\circ}$.   
Now let $i\in {\rm V}^u$   and $f \in   {\rm V}^f$. Then by commutation relations (\ref{4.28.03.11x}), ${\Bbb B}_f$ and $X_i$ commute, and we have:
\be \nonumber
\begin{split}
&{\rm Y}_i \cdot \Bigl(\prod_{j \in {\rm V}}{\rm B}_j^{-p_{fj}}\Bigr) =q_i^{-2\varepsilon_{fi} } \Bigl(\prod_{j \in {\rm V}}{\rm B}_j^{-p_{fj}}\Bigr) \cdot{\rm Y}_i, \\
&{\rm Z}_i  {\rm Z}_f  = q^{-2   \widehat \varepsilon_{fi}}{\rm Z}_f {\rm Z}_i= q_i^{-2\varepsilon_{fi} } {\rm Z}_f {\rm Z}_i.\\
\end{split}
\ee
So $\xi^\ast$ is an algebra map.   
 Let $k\in {\rm V}^u$. By definition, $\mu_k^{\cal Z} = \mu^{\mathscr X}_k\circ \mu^{\mathscr X}_{k^\circ}$. It remains to check that 
 \be
 \label{commuta,ac}
 \xi^\ast\circ (\mu_k^{\cal Z} )^* =(\mu_k^{\cal D})^* \circ \xi^*.
 \ee
Note that if $i, j\in {\rm V}^u$ then mutations at   $i$ do not change   ${\rm Z}_{j^\circ}$, and vice versa. So the maps in \eqref{commuta,ac} are equal to each other when restricted on the subalgebra generated by ${\rm Z}_i$ and ${\rm Z}_{j^\circ}$. 
For a frozen ${\rm Z}'_f$, we get
$$
(\mu_k^{\cal Z} )^*({\rm Z}'_f) = {\rm Ad}_{\Psi_{q_k}({\rm Z}_k) \Psi_{q_k}({\rm Z}_{k^\circ})}\Bigl(q^{-\varepsilon_{fk }^2}{\rm Z}_f {\rm Z}_{k^\circ}^{[\varepsilon_{fk^\circ}]_+}{\rm Z}_{k }^{[\varepsilon_{fk }]_+}\Bigr).
$$
Under the map $\xi^*$ it goes to the following expression: 
\be \la{12.5.18.1}\begin{split}
&
q^{-\varepsilon_{fk }^2}\cdot{\rm Ad}_{\Psi_{q_k}({\rm Y}_k) \Psi_{q_k}(\widetilde {\rm Y}_{k})}\Bigl( \mathbb{B}_f\prod_j\B^{-\varepsilon_{fj}}_j  \widetilde {\rm Y}_{k }^{[-\varepsilon_{fk}]_+}
{\rm Y}_{k }^{[\varepsilon_{fk }]_+}\Bigr)  \\
&=
\mathbb{B}_f\prod_{j\not = k}\B^{-\varepsilon_{fj} - [\varepsilon_{fk}]_+  \varepsilon_{kj}}_j  \cdot  q^{-\varepsilon_{fk }^2}\cdot {\rm Ad}_{\Psi_{q_k}({\rm Y}_k) 
\Psi_{q_k}(\widetilde {\rm Y}_{k})} \Bigl(\B_k ^{-\varepsilon_{fk}}\cdot  \widetilde {\rm Y}_{k }^{-\varepsilon_{fk }}\Bigr) .\\
\end{split}
\ee
On the other hand, 
\be \la{12.5.18.2}
\begin{split}
&(\mu_k^{\cal D})^*\circ \xi^* ({\rm Z}'_f)=  \mathbb{B}_f\prod_{j \not = k}{\B}^{-\varepsilon'_{fj}}_j \cdot (\mu_k^{\cal D})^*({\B_k'})^{-\varepsilon'_{fk}}       \\
&\qquad \qquad=\mathbb{B}_f \prod_{j \not = k}{\B}^{-\varepsilon'_{fj}}_j \cdot {\rm Ad}_{\Psi_{q_k}({\rm Y}_k) /\Psi_{q_k}(\widetilde {\rm Y}_{k}^{-1})}\Bigl( {\Bbb B}_k^-/{\B}_k  \Bigr)^{\varepsilon_{fk}}\\
&\qquad \qquad={\Bbb B}_f \prod_{j \not = k}{\B}^{-\varepsilon'_{fj} + [-\varepsilon_{kj}]_+ \varepsilon_{fk}}_j \cdot {\rm Ad}_{\Psi_{q_k}({\rm Y}_k) /\Psi_{q_k}(\widetilde {\rm Y}_{k}^{-1})} \left(\B_k ^{-\varepsilon_{fk}}\right).\\
\end{split}
\ee
 
 By  the mutation formula for $\varepsilon'_{fj}$, we get
$
-\varepsilon_{fj} - [\varepsilon_{fk}]_+  \varepsilon_{kj} = -\varepsilon'_{fj} + [-\varepsilon_{kj}]_+ \varepsilon_{fk}. 
$
Therefore the two  left $\B$-terms in the last line of (\ref{12.5.18.1})  and  (\ref{12.5.18.2}) are equal. 
Note also that ${\rm Ad}_{\Psi_{q_k}({\rm Y}_k)}$ in those lines give the same result, and can be cancelled. 
So it remains to check that 
 \be
  \la{cahciashochdow}
\begin{split}
&q^{-\varepsilon_{fk }^2}\cdot 
  {\rm Ad}_{\Psi_{q_k}(\widetilde {\rm Y}_{k})} 
 \left(\B_k ^{-\varepsilon_{fk}}\cdot  \widetilde {\rm Y}_{k }^{-\varepsilon_{fk }}\right)\stackrel{?}{=} 
     {\rm Ad}^{-1}_{ \Psi_{q_k}(\widetilde {\rm Y}_{k}^{-1})} (\B_k^{-\varepsilon_{fk}}).\\
\end{split}
\ee
Let $u=\widetilde {\rm Y}_{k}, v= \B_k^{-\varepsilon_{fk}}$ and $n = -\varepsilon_{fk}$. Formula \eqref{cahciashochdow} follows from Lemma \ref{cnsadibvcsoa}. 
\end{proof}

\bl  
\la{cnsadibvcsoa}
Let $n \in \Z$. Then if $uv=q^{-2n}vu$, then 
$$
 q^{-n^2}{\rm Ad}_{\Psi_q(u)}(v  u^n ) = {\rm Ad}^{-1}_{\Psi_q(u^{-1})}(v).
$$
\el

\begin{proof} Assume first that $n \geq 0$. Then calculating the left   and the right hand sides, we get:
\be \nonumber
\begin{split}
&{\rm LHS}= q^{-n^2}\Psi_q(u)v \Psi_q(u)^{-1} u^n= q^{-n^2}v \Psi_q(q^{-2n}u)\Psi_q(u)^{-1}u^n \\ 
& \qquad ~\qquad~ = q^{-n^2}v \left( \prod_{k=1}^n(1+q^{-2k+1}u)^{-1} \right) u^n = v  \prod_{k=1}^n(1+q^{2k-1}u^{-1})^{-1}  .\\
&{\rm RHS}= \Psi_q(u^{-1})^{-1}v\Psi_q(u^{-1}) = v \Psi_q(q^{2n}u^{-1})^{-1}\Psi_q(u^{-1})= {\rm LHS}.\\
\end{split}
\ee
If $n<0$,  similar arguments prove the claim.  
\end{proof}

   \medskip

\subsection{Proof of Theorems \ref{11.26.18.1} and \ref{11.26.18.1Xaaa} }  \la{SEC9.2}

\medskip

\begin{proof}  1. Recall the map $\pi_\bS$, defined in   \eqref{CORT1X.general.nvadfo}:
$$
\pi_{\bS}: \mathscr{A}_{{\G'}, \bS} \lra \mathscr{P}_{\G, \bS}.
$$     
 By Theorem \ref{11.24.18.1},  
 in any cluster coordinate system we have:
 \be \la{CORT1}
\pi_\bS^*X_k = \prod_j A_j^{p_{kj}} =  \left\{    \begin{array}{ll}  
      {\prod_{j\in{\rm I}_{\G, \bS}} A_j^{\varepsilon_{kj}}} & \mbox{if $k$ is unfrozen}, \\
      &\\
 \mathbb{A}_{k}^{1/2} \cdot {\prod_{j\in {\rm V}_{\G, {\bS}}} A_j^{ \varepsilon_{kj}}} & \mbox{if $k$ is frozen}. \\
    \end{array} \right.
\ee
The right hand side of this formula  defines the matrix $p=(p_{ij})$. This matrix is integral by  Theorem \ref{11.24.18.1}.

Here is an example of the map $\pi^*_\bS$ when $\bS=t$ is  a triangle, and $\G'=SL_3$.   
\begin{figure}[ht]
 \begin{center}
\begin{tikzpicture}[scale=2]
\node (a1) at (120:1) {$A_1$};
\node (a2) at (60: 1) {$A_2$};
\node (b1) at (180: 1) {$A_3$};
\node (b2) at (0,0) {$A_4$};
\node (b3) at (0: 1) {$A_5$};
\node (c1) at (240: 1) {$A_6$};
\node (c2) at (300: 1) {$A_7$};
\foreach \from/\to in {a2/a1, b3/b2, b2/b1, b1/c1, a1/b2, b2/c2, c2/b3, c1/b2, b2/a2}
                 \draw[directed] (\from) -- (\to);   
\foreach \from/\to in {b1/a1, a2/b3, c2/c1}
                 \draw[directed, dashed] (\from) -- (\to);        
      
\begin{scope}[shift={(5,0)}]
\node (a1) at (120:1) {$\frac{A_4}{A_1 A_2}$};
\node (a2) at (60: 1) {$\frac{A_1A_5}{A_2A_4}$};
\node (b1) at (180: 1) {$\frac{A_1A_6}{A_3A_4}$};
\node (b2) at (0,0) {$\frac{A_2A_3A_7}{A_1A_5A_6}$};
\node (b3) at (0: 1) {$\frac{A_4}{A_5A_7}$};
\node (c1) at (240: 1) {$\frac{A_4}{A_3A_6}$};
\node (c2) at (300: 1) {$\frac{A_5A_6}{A_4A_7}$};
\foreach \from/\to in {a2/a1, b3/b2, b2/b1, b1/c1, a1/b2, b2/c2, c2/b3, c1/b2, b2/a2}
                 \draw[directed] (\from) -- (\to);   
\foreach \from/\to in {b1/a1, a2/b3, c2/c1}
                 \draw[directed, dashed] (\from) -- (\to);     
\end{scope} 
\end{tikzpicture}
\end{center}
\caption{The map $\pi^*_\bS$ for the triangle $\bS=t$, and $\G'=SL_3$.}
\label{FigA2}
\end{figure}

Since   $\bS$ has no punctures, the map $\pi_\bS$ is surjective, and  $p$ is of full rank.  
 
 There is a   surjective map
\be
\la{double.maps.advfs}
\begin{split}
&\alpha: \mathcal{D}_{\G, \bS} \lra \mathscr{P}_{\G, \bS} \times \mathscr{P}_{\G, \bS^\circ}\\
& \alpha^\ast(X_i) = Y_i, \qquad \alpha^*(X_{i^\circ}) = Y_i^{-1} \prod_j B_j^{-p_{ij}}.\\
\end{split}
\ee
It is compatible with cluster mutations, and so is globally well defined. However it may  not be Poisson. 

Next, recall the   gluing map, inducing  a birational isomorphism after taking the quotient by the diagonal action of $\H^n$:
   \be \la{GLUEINGM}
   \begin{split}
&\gamma:  \mathscr{P}_{\G, \bS} \times \mathscr{P}_{\G, \bS^\circ} \lra \mathscr{P}_{\G, \bS_{\cal D}}.\\
& \gamma:  ({\mathscr P}_{\G, \bS} \times  {\mathscr P}_{\G, \bS^\circ} )/ \H^n    \stackrel{\sim}{\lra}      {\mathscr P}_{\G, \bS_{\cal D}}.\\
 \end{split}
  \ee  
We get a dominant map 
$
  \gamma\circ \alpha: \mathcal{D}_{\G, \bS} \lra \mathscr{P}_{\G, \bS_{\cal D}}. 
$
It factorizes through   a dominant map 
$\xi: \mathcal{D}^{\bf b}_{\G, \bS} \lra \mathscr{P}_{\G, \bS_{\cal D}}$:

 \begin{displaymath}
    \xymatrix{
        \mathcal{D}_{\G, \bS}  \ar[dr]^{\gamma\circ \alpha}  \ar[d]_{} &       \\
          \mathcal{D}^{\bf b}_{\G, \bS} \ar[r]^{\xi} &   \mathscr{P}_{\G, \bS_{\cal D}} }
         \end{displaymath}
 
 The map $\alpha$ becomes Poisson  after adding a cross-term Poisson bivector on ${\mathscr P}_{ \G, \bS} \times {\mathscr P}_{ \G, \bS^\circ}$, 
which  maps to zero under the map $\delta$. So the composition 
$\alpha \circ \delta$ is a Poisson monomial map.

  By Theorem \ref{THH8.9}, the map $\xi$ delivers   map 
(\ref{zeta}). Part 1) of Theorem \ref{11.26.18.1} is proved. 

\vskip 2mm
2. {\it Monodromy around punctures of $\mathscr{P}_{\G, \bS_{\cal D}}$.} 
Consider the following surjective map
\be \la{11.24.18.2}
\begin{split}
&\delta:~ \mathscr{A}_{\G', \bS} \times \mathscr{A}_{\G', \bS^\circ} \lra \mathcal{D}_{\G, \bS},\\
&\delta^\ast (Y_i) := \prod_j A_j^{p_{ij}},  \qquad  \delta^\ast (B_j) := \frac{A_j^\circ}{A_j}.\\
\end{split}
\ee

Formulas (\ref{11.24.18.2}) define the map $\delta^*$ in each  cluster coordinate system. 
It is compatible with mutations. Indeed,    by Theorem \ref{8.15.05.10} the ratio $A^\circ_i/A_i$   
 mutates  as the $B_i$-coordinate  on the double,    
 and, up to a   non-mutating frozen coefficient ${\Bbb A}_k^{1/2}$,  $\delta^*(Y_i)$ is the standard     $p-$map. 
The composition $\alpha\circ \delta$ is a surjection
\[
\alpha\circ \delta=(\pi, \pi^\circ): ~  \mathscr{A}_{\G', \bS} \times \mathscr{A}_{\G', \bS^\circ}  \lra  \mathscr{P}_{\G, \bS} \times \mathscr{P}_{\G, \bS^\circ}.
\]
By definition, the map $\pi$ coincides with the map $\pi_{\bS}$, the map $\pi^\circ$ is defined by
\[
(\pi^\circ)^\ast (X_{i^\circ})= \prod_j (A^\circ_{j })^{-p_{ij}}. 
\]
The map $\pi^\circ$ is different from the map $\pi_{\bS^\circ}$ since it shifts the decorated flag to a different side. 

\vskip 2mm \paragraph{\bf Remarks.}  1. The maps $\alpha$ and $\delta$ are   relatives of the maps $\varphi$ and $\pi$   in the definition of the cluster symplectic double, see 
Theorem   \ref{8.15.05.10}  or {\cite[Theorem 2.3]{FG07}}. The notable difference is that the map $\alpha$ and $\delta$ are 
surjective if $\bS$ has no punctures. Indeed, we use the integral matrix $p_{ij}$   instead of the possibly non-integral exchange matrix $\varepsilon_{ij}$. 
The matrix $p_{ij}$ is of full rank if $\bS$ has no punctures.

2.  The map $\xi$ can be defined   geometrically as follows. We start with
 the geometric map  
 \be \la{pppp}
  \pi_\bS\times \pi^\circ_\bS: {\mathscr A}_{\G', \bS} \times  {\mathscr A}_{\G', \bS^\circ} \lra  {\mathscr P}_{\G', \bS} \times  {\mathscr P}_{\G', \bS^\circ}.
\ee 
 It is a surjective map. It factorizes as the composition 
$$
 \pi_\bS\times \pi^\circ_\bS: {\mathscr A}_{\G', \bS} \times  {\mathscr A}_{\G', \bS^\circ} \lra {\cal D}_{\G', \bS}  
 \lra  {\mathscr P}_{\G', \bS} \times  {\mathscr P}_{\G', \bS^\circ}.
 $$
Combining    (\ref{pppp}) with the gluing map $\gamma$, we get a map ${\cal D}_{\G', \bS}  \lra {\mathscr P}_{\G', \bS_{\cal D}}$, which   evidently 
factorizes through the quotient 
${\cal D}_{\G', \bS}/ {\H'}^n  = {\cal D}^{\bf b}_{\G', \bS}$.   
Therefore we get a commutative diagram:
   
\begin{equation}\nonumber
\begin{gathered}
\xymatrix{
 {\mathscr A}_{\G', \bS} \times  {\mathscr A}_{\G', \bS^\circ}  \ar[r]^{\qquad\delta} \ar[dr]_{\pi_\bS \times \pi_\bS^\circ} &{\cal D}_{\G', \bS}       
    \ar[r]^{ \rho}  \ar[d]^{\alpha}&  \ar[d]^{\xi}  {\cal D}_{  \G, \bS}^{{\bf b}}\\
  &  {\mathscr P}_{\G, \bS} \times  {\mathscr P}_{\G, \bS^\circ}     \ar[r]^{\quad \gamma}   &  {\mathscr P}_{\G, \bS_{\cal D}} }
\end{gathered}
\end{equation}

\vskip 2mm

As shown on Figure \ref{gluing1}, under the map $\pi=\pi_{\bS}$, we obtain pinnings
\[
p_k= (\A_k, \A_{k+1}'), \qquad ~ p_{k-1}= (\A_{k-1}, \A_k').
\]
Under the map $\pi^\circ$, we obtain pinnings
\[
p_k^\circ = (\A_k^\circ, \A^{\circ'}_{k+1}), \qquad ~\qquad~ p_{k-1}^\circ= (\A_{k-1}^\circ, \A_k^{\circ'}).
\]
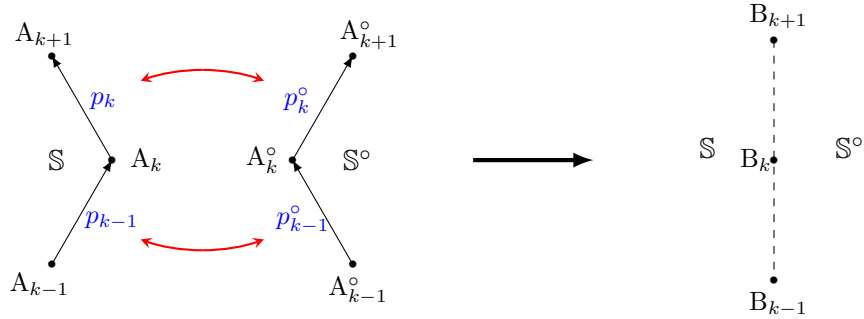
\begin{figure}[ht]
\begin{center}
\begin{tikzpicture}[scale=0.8]
\node at (0:0.5) {{ \small $\A_k$}};
\node[blue] at (98:1)  {{\small $p_{k}$}};
\node[blue] at (270:1)  {{\small $p_{k-1}$}};
\node at (0:-1) {{ $\bS$}};
\node at (120:2.4)  {{ \small $\A_{k+1}$}};
\node at (240:2.4)  {{\small $\A_{k-1}$}};
\draw[red, thick, -stealth] (1.5, 1.5) arc (90:70:3);
\draw[red, thick, -stealth] (1.5, 1.5) arc (90:110:3);
\draw[red, thick, -stealth] (1.5, -1.5) arc (-90:-70:3);
\draw[red, thick, -stealth] (1.5, -1.5) arc (-90:-110:3);
\node at (0:0) {{\tiny $\bullet$}};
\node at (120:2) {{\tiny $\bullet$}};
\node at (240:2) {{\tiny $\bullet$}};
\draw[directed] (240:2)--(0:0);
\draw[directed] (0,0)--(120:2);
\begin{scope}[shift={(3,0)}]
\node at (180:0.5) {{\small $\A_k^\circ$}};
\node[blue] at (280:1)  {{\small $p_{k-1}^\circ$}};
\node[blue] at (90:1)  {{ \small $p_{k}^\circ$}};
\node at (180:-1) {{ $\bS^\circ$}};
\node at (295:2.4)  {{ \small $\A_{k-1}^\circ$}};
\node at (60:2.4)  {{ \small $\A_{k+1}^\circ$}};
\node at (180:0) {{\tiny $\bullet$}};
\node at (300:2) {{\tiny $\bullet$}};
\node at (60:2) {{\tiny $\bullet$}};
\draw[directed] (-60:2)--(0:0);
\draw[directed] (0,0)--(60:2);
\end{scope}
\draw[ultra thick, -latex] (6,0)--(8,0);
\begin{scope}[shift={(11,0)}]
\node at (180:0.3) {{\small $\rF_k$}};
\node at (10:1.2) {{ $\bS^\circ$}};
\node at (170:1.2) {{ $\bS$}};
\node at (270:2.4)  {{ \small $\rF_{k-1}$}};
\node at (90:2.4)  {{ \small $\rF_{k+1}$}};
\node at (180:0) {{\tiny $\bullet$}};
\node at (270:2) {{\tiny $\bullet$}};
\node at (90:2) {{\tiny $\bullet$}};
\draw[dashed] (90:2)--(180:0)--(270:2);
\end{scope}
\end{tikzpicture}
\end{center}
 \caption{Gluing $\mathscr{P}_{\G, \bS}$ with $\mathscr{P}_{\G, \bS^\circ}$ along the pinnings.}
\label{gluing1}
\end{figure}
After gluing $\mathscr{P}_{\G, \bS}$ with $\mathscr{P}_{\G, \bS^\circ}$ along the pinnings, the semi-simple part of the monodromy 
$h_k$ around the puncture which carries the flag $\rF_k$ is given by
\[
h_k= h(\A_k', \A_k) h (\A_k^{\circ}, \A_k^{\circ'}) = h(\A_{k}, \A_{k-1})^{-1} h(\A_k^\circ, \A_{k-1}^\circ ) = \prod_{j\in {\rm I}} \alpha^\vee_j(\B_j). 
\]
Here $\B_i$ are the $\B$-coordinates of the symplectic double associated to frozen vertices on the boundary interval $\{k-1, k\}$ of $\bS$. In particular, we have
\be \la{MM}
M_{i}= \alpha_i(h_k) = \prod_{j\in {\rm I}} B_j^{C_{ij}}. 
\ee
The monodromy  around the puncture $p_k$ is a  standard monomial in  $\mathscr{P}_{\G, \bS_{\cal D}}$, and therefore admits the unique  quantum lift ${\rm M}_i$ in $\mathcal{O}_q(\mathscr{P}_{\G, \bS_{\cal D}})$.
Moreover, it is clear from (\ref{MM}) that we have
\[
\xi^\ast({\rm M}_i)= \prod_{j\in {\rm I}} {\rm B}_j^{C_{ij}}.
\]

Theorem \ref{11.26.18.1} and hence Theorem \ref{11.26.18.1Xaaa} are proved. 
\end{proof}

  \subsection{Complimentary results to Theorem \ref{11.26.18.1}.}  
\medskip

The group $ {\H'}^n$   acts  diagonally on ${\mathscr A}_{\G', \bS} \times  {\mathscr A}_{\G', \bS^\circ}$.

 \bp \la{Lemma9.7}There is a birational isomorphism
\be \la{mptheta}
\theta:  {\mathscr P}_{\G', \bS_{\cal D}} \stackrel{\sim}{\lra} ({\mathscr A}_{\G', \bS} \times  {\mathscr A}_{\G', \bS^\circ} )/{\H'}^n.
\ee
  \ep
 \begin{proof} Let us define the map $\theta$.  First,  cut out a  little disc around each puncture $p_{k}$. Let $\alpha_k$ be   the boundary loop, 
  oriented   against the surface orientation. We picture  on Figure \ref{pin10} the  surface orientation   as  clockwise.

\begin{figure}[ht]
\begin{center}
\begin{tikzpicture}[scale=1.4]
\node at (90:0.5) {{ $\widehat{\rm F}_1$}};
\node at (210:.5)  {{ $\widehat{\rm F}_2$}};
\node at (330:.5)  {{ $\widehat{\rm F}_3$}};
\node at (90:1.5) {{ $\widehat{\rm F}_1^\circ$}};
\node at (210:1.5)  {{ $\widehat{\rm F}_2^\circ$}};
\node at (330:1.5)  {{ $\widehat{\rm F}_3^\circ$}};
\draw[red, thick, -stealth] (90:.75) arc (270:90:0.25);
\draw[red, thick, -stealth, rotate=120] (90:.75) arc (270:90:0.25);
\draw[red, thick, -stealth, rotate=240] (90:.75) arc (270:90:0.25);
\node at (90:1) {{\tiny $\bullet$}};
\node at (210:1) {{\tiny $\bullet$}};
\node at (330:1) {{\tiny $\bullet$}};
\draw[dashed] (90:1)--(210:1)--(330:1)--(90:1);
\end{tikzpicture}
\end{center}
\caption{The map $\theta$ for a pair of pants glued from two triangles. The $\bS-$triangle is shown.}
\label{pin10}
\end{figure}

Take a framed $\G'$-local system  on 
$\bS_{\cal D}$. Near each puncture $p_k$ there is  an invariant  flag ${\rm F}_k$. 
We assign it to a point $x_k$ at the $\bS-$part of the loop $\alpha_k$.  
Pick  a   decorated flag  $ \widehat {\rm F}_k $ over   
$ {\rm F}_k $, and assign  it  to a positively oriented   tangent vector to the loop $\alpha_k$ at  $x_k$.\footnote{Positively oriented is depicted  on Figure \ref{pin10} as  oriented clockwise around the puncture $p_k$.} Decorated flags $(\widehat  {\rm F}_1, \ldots , \widehat  {\rm F}_n)$ provide a point of the  space ${\mathscr A}_{\G', \bS}$. 
Moving  $\widehat  {\rm F}_k$ along the   {oriented} loop $\alpha_k$ to $\bS^\circ$ we get 
a decorated flag $\widehat  {\rm F}^\circ_k$, see Figure \ref{pin10}.   Decorated flags $(\widehat  {\rm F}^\circ_1, \ldots , \widehat  {\rm F}^\circ_n)$ define a point of 
  ${\mathscr A}_{\G', \bS^\circ}$. So we get   a map
 $$
   \theta':  {\mathscr P}_{\G', \bS_{\cal D}} \stackrel{}{\lra}  {\mathscr A}_{\G', \bS} \times  {\mathscr A}_{\G', \bS^\circ}. 
   $$
It depends on a choice of  decorated flags $\widehat  {\rm F}_k$. Another choice gives  configurations   
$( \widehat  {\rm F}_1 \cdot  h_1, \ldots ,    \widehat  {\rm F}_n\cdot  h_n)$  and $(  \widehat  {\rm F}^\circ_1\cdot  h_1, \ldots ,   \widehat  {\rm F}^\circ_n\cdot  h_n)$, 
where $h_k \in {\H'}$. Therefore the image in the quotient by the diagonal action of $ {\H'}^n$ is well defined.  So we get the map $\theta$.  

Let us define now the inverse map 
\be \nonumber
\eta:  ({\mathscr A}_{\G', \bS} \times  {\mathscr A}_{\G', \bS^\circ} )/{\H'}^n    \stackrel{\sim}{\lra} {\mathscr P}_{ \G', \bS_{\cal D}}.
\ee
Take a    generic point 
$
 x\in ({\mathscr A}_{\G', \bS} \times  {\mathscr A}_{\G', \bS^\circ} )/ {\H'}.
 $ 
  Then the   ratio of the $h$-invariants of the   pairs of decorated flags at  the side 
 $(p_{k}, p_{k+1})$ is well defined: 
\be \nonumber
 {\bf b}_k:= h( {\rm F}^\circ_k,  {\rm F}^\circ_{k+1} ) / h(  {\rm F}_k,  {\rm F}_{k+1} )   \in {\H'}. 
\ee
 
We can glue      trivial $\G'$-local systems with decorated flags  $( \widehat  {\rm F}_1, \ldots ,    \widehat  {\rm F}_n)$ and  $( \widehat  {\rm F}^\circ_1, \ldots ,    \widehat  {\rm F}^\circ_n)$ on    $\bS$ and $\bS^\circ$ by   identifying the pairs 
 $(\widehat  {\rm F}_k \cdot {\bf b}_k,  \widehat {\rm F}_{k+1} )$ and  $( \widehat {\rm F}^\circ_k, \widehat  {\rm F}^\circ_{k+1} )$ 
since  their $h$-invariants coincide. 
 Let $\eta(x)$ be the obtained framed $\G'$-local system on  
 $\bS_{\cal D}$. Then  by the construction  $\eta = \theta^{-1}$.  
 \end{proof}  
 
 Composing    $\delta$ from (\ref{11.24.18.2}) with the projection ${\cal D}_{\G, \bS} \to {\cal D}_{  \G, \bS}^{{\bf b}}$ we get a  map 
$$
\bar \delta: ({\mathscr A}_{\G', \bS} \times  {\mathscr A}_{\G', \bS^\circ} )/ {\H'}^n \to {\cal D}_{  \G, \bS}^{{\bf b}}.
$$

 Therefore there is  a commutative diagram 
of rational maps: 
  \begin{equation}\nonumber
\begin{gathered}
\xymatrix{
  {\mathscr P}_{\G', \bS_{\cal D}} \ar[r]_{\sim\qquad \qquad }^{\theta\qquad \qquad }  &({\mathscr A}_{\G', \bS} \times  {\mathscr A}_{\G', \bS^\circ} )/ {\H'}^n     
    \ar[r]^{\qquad \qquad \bar \delta}  \ar[d]_{\pi_\bS \times \pi_\bS^\circ}& \ar[d]^{\xi}  {\cal D}_{  \G, \bS}^{{\bf b}}\\
    &({\mathscr P}_{\G, \bS} \times  {\mathscr P}_{\G, \bS^\circ} )/ \H^n    \ar[r]^{\qquad \qquad \sim}_{\qquad \qquad \gamma}   &  {\mathscr P}_{\G, \bS_{\cal D}} }
\end{gathered}
\end{equation}
They are monomial maps in each  cluster coordinate system.  The map $\xi$ is given by (\ref{12.12.xcv}). The composition $$
\gamma\circ (\pi_\bS,  \pi_\bS^\circ)\circ \theta: {\mathscr P}_{\G', \bS_{\cal D}} \to  {\mathscr P}_{\G, \bS_{\cal D}} 
$$ 
  assigns to a  $\G'$-local system  its quotient by   ${\rm Cent}(\G')$.

 Let  $h_k \in {\H'}$ be the monodromy  of the $\G'$-local system $\eta(x)$ around the puncture $p_{k}$.
 It  is the inverse of the monodromy around the loop $\alpha_k$.    
 
    \bl \la{MMONN} The element  ${\bf b}_k$ describes the monodromy around the puncture $p_{k}$:  ${\bf b}_k =  {h_k}\in {\H'}$. 
  \el
  
  \begin{proof}
  
Moving   $\widehat {\rm F}_{k}$ to $\bS^\circ$ along the   loop $\alpha_{k}$ we  get
 a decorated flag $\widehat {\rm F}^\circ_{k}$. Moving   $\widehat {\rm F}_{k}$ to the same point    
 against the orientation of the loop $\alpha_{k}$ we get a decorated flag $ \widehat {\rm G}^\circ_{k}$. 
 Then we have 
\be \la{11.27.18.2}
h(\widehat{\rm F}_k, \widehat{\rm F}_{k+1}) = h(\widehat{\rm G}^\circ_k,   {\widehat{\rm F}}^\circ_{k+1}) \in {\H'}.
\ee
 Indeed, if we put  $\widehat {\rm F}_k$ and $\widehat {\rm F}_{k+1}$    at the intersection   of the loops $\alpha_k$ and $\alpha_{k+1}$ 
 with  the side $(p_k, p_{k+1})$, 
 then  $\widehat {\rm G}^\circ_{k} = \widehat {\rm F}_{k}$ and $\widehat {\rm F}^\circ_{k+1}= \widehat {\rm F}_{k+1}$. 
Since $ \widehat {\rm G}^\circ_{k}$ is obtained from   $ \widehat {\rm F}^\circ_{k}$ by the monodromy 
along the  loop $\alpha_k$, we have $\widehat {\rm G}^\circ_{k}   =  \widehat {\rm F}^\circ_{k}\cdot h^{-1}_{k}$. Therefore
\be \nonumber
\begin{split}
&{\bf b}_k = h(\widehat{\rm F}^\circ_k,   {\widehat{\rm F}}^\circ_{k+1})/ h(\widehat{\rm F}_k, \widehat{\rm F}_{k+1})  =h_k  \cdot h(\widehat{\rm G}^\circ_k,   {\widehat{\rm F}}^\circ_{k+1})/ h(\widehat{\rm F}_k, \widehat{\rm F}_{k+1}) \stackrel{({\ref{11.27.18.2}})}{=}  h_k .\\
\end{split}
\ee
\end{proof}

 \bc \la{12.8} Any cluster coordinate system $(A_1, \ldots , A_N)$ on the space ${\mathscr A}_{\G', \bS}$  
  gives rise to  
   a maximal collection of Poisson commuting rational functions $B_i:= \theta^*(A^\circ_i / A_i)$ on the space ${\mathscr P}_{\G', \bS_{\cal D}}$. 
   The frozen coordinates $\{B^{(k)}_1, ..., B^{(k)}_r\}$ at the boundary interval ${\rm I}_k$ describe the monodromy around the puncture $p_k$ on $ \bS_{\cal D}$. \ec

  \medskip 
  
  \subsection{Realizations of  the multiplicity spaces for the quantum group principal series} \la{sec16.5:RMSP}

\medskip
Recall  that cluster varieties come in   triples $({\mathscr  A}, {\mathscr  X}, {\mathscr  D})$. Here ${\mathscr  X}$ is a cluster Poisson variety,  ${\mathscr  A}$ is the related cluster $K_2-$variety such that the pair $({\mathscr  A}, {\mathscr  X})$ is a cluster ensemble, and  ${\mathscr  D}$ is the cluster symplectic double of  the 
cluster Poisson variety ${\mathscr  X}$. 

Recall 
 the canonical representation of the  $\ast-$algebra $\mathcal{A}_\hbar({\cal D})$    in the cluster Schwartz space\footnote{As always,  "cluster vector space" refers not to a single 
 vector space, but rather to a collection of vector spaces assigned to different clusters, related by the canonical intertwiners.} 
\be \la{}
{\cal S}\left({\mathscr A}(\R_{>0})\right) \subset {\rm L}_2\left({\mathscr A}(\R_{>0}), \mu_{\mathscr A}\right)
\ee
 of the cluster Hilbert space ${\rm L}_2\left({\mathscr A}(\R_{>0}), \mu_{\mathscr A}\right)$ assigned to the  cluster $K_2-$variety  ${\mathscr A}$, see Section \ref{Sec8a}.  
 \vskip 2mm
 
 Let $\bS$ be an $n$-gon ${\Bbb P}_n$. Since the pair  $({\mathscr A}_{  \G, \P_n}, {\mathscr P}_{\G, \P_n})$ is a cluster ensemble, 
  there is a canonical      representation of the $\ast$-algebra   $\mathcal{A}_\hbar({\cal D}_{\G, \P_n})$    in the cluster Schwartz space  
\be \la{HDela+}
 {\cal S}_{\G, \P_n} \subset  
   {\rm L_2}\left({\mathscr A}_{  \G, \P_n}(\R_{>0}), \mu_{\mathscr A}\right).
\ee   
\vskip 2mm

Recall the commutative subalgebra ${\bf B} \subset \mathcal{A}_\hbar({\cal D}_{\G, \P_n})$ from the  Definition \ref{16.1}. 
 By the   isomorphism (\ref{CaniB}),  ${\bf B} = {\cal O}(\H'^n)$. 
According to (\ref{handyZ}),  the subalgebra ${\bf B}$ acts in the representation of the $\ast-$algebra 
$\mathcal{A}_\hbar({\cal D}_{\G, \P_n})$ in (\ref{HDela+})  by multiplication by the  
frozen coordinates 
on the space ${\mathscr A}_{  \G, \P_n}(\R_{>0})$.  The latter are the coordinates on $ {\H'}^n(\R_{>0}) =   \H^n(\R_{>0})$. 
So the spectral decomposition for the action of subalgebra ${\bf B}$   is on the nose the decomposition into an integral of subspaces 
obtained by fixing the values $\lambda\in   \H^n(\R_{>0})$ of the frozen variables. 

\be \la{QREF}
 {\rm L_2}({\mathscr A}_{  \G, \P_n}(\R_{>0}), \mu_{\mathscr A}) = \int {\cal M}_{\lambda} d\lambda, \qquad \lambda\in   \H^n(\R_{>0}).
\ee
\vskip 2mm

 We claim  that the   $\ast$-representations ${\cal M}_\lambda$   are  the  multiplicity spaces  for the tensor products of the 
  principal series  representations of the $\ast-$algebra given by the modular double   
  of the quantum group 
  $$
\mathcal{A}_\hbar(\mathfrak{g}) = {\U}_q(\mathfrak{g}) \otimes {\U}_{q^\vee}(\mathfrak{g}^\vee).
$$  
By Theorem \ref{MFCQIX}, which assumes  the Modular Functor Conjecture \ref{MK}, we have:
\vskip 2mm

(1a) {\it The multiplicity space  for the tensor product of the principal series $\ast-$representations of the quantum group modular double 
$\ast-$algebra $\mathcal{A}_\hbar(\mathfrak{g})$  is the canonical $\ast-$representation of the $\ast-$algebra $\mathcal{A}_\hbar({\mathscr X}_{\G, S^2_n})$.}
\vskip 2mm

Next, Theorem \ref{11.26.18.1} provides the following  maps of $\ast-$algebras 
 $ \xi^*_q$ and  $\xi_\hbar^*:=  \xi^*_q\otimes \xi^*_{q^\vee}$:
 \be \la{16.6.1}
 \begin{split}
\xi_q^*: ~ &{\cal O}_q({\mathscr X}_{\G, S^2_n}) \stackrel{}{\lra} {\cal O}_q({\cal D}_{\G, \P_n})^{\bf B},\\
\xi_\hbar^*:  ~&\mathcal{A}_\hbar({\mathscr X}_{\G, S^2_n}) \stackrel{}{\lra} \mathcal{A}_\hbar({\cal D}_{\G, \P_n})^{\bf B}.\\
\end{split}
 \ee
Therefore, combining the map $\xi_\hbar^*$  with the canonical representation of the $\ast-$algebra $\mathcal{A}_\hbar({\cal D}_{\G, \P_n})$ 
in the cluster Schwartz space  (\ref{HDela+}), 
 we arrive at: \vskip 2mm
 
(1b)  {\it Realization of the canonical representation of the $\ast-$algebra $\mathcal{A}_\hbar({\mathscr X}_{\G, S^2_n})$ 
in the   space (\ref{HDela+})}.

Therefore, combining (1a) and (1b) we conclude that
\vskip 2mm

(1)   {\it Assuming  the Modular Functor Conjecture \ref{MK},  the space (\ref{HDela+}) is the quantum multiplicity space}. 
\vskip 2mm

 Now let us
   compare the realizations of the multiplicity spaces for the classical and quantum principal series representations,  discussed in Sections \ref{SECT5.2.1}  and  \ref{SC2.7}.

 For any group $\G$, not necessarily adjoint, 
the group 
 $\G \times \H$ acts   on the principal affine space ${\cal A}_\G=\G/\U$: the group $\G$ acts from the left, and the Cartan group $\H$ acts from the right, since normalises the 
 maximal unipotent subgroup $\U$. So the group $\G(\R) \times \H(\R)$  acts on the vector 
 space ${\cal F}({\cal A}_\G(\R))$ of complex valued functions on the space ${\cal A}_\G(\R)$,  providing the 
 Gelfand-Naimark  principal series of representations. For each quasicharacter $\chi: \H(\R) \lra \C^\times$ there is a representation ${\cal V}_\chi$ 
 in the space of functions   homogeneous of the degree $\chi$ under the action of   $\H(\R)$:
 $$
 {\cal V}_\chi = \{f \in {\cal F}({\cal A}_\G(\R)) \ | \ f(ah) = \chi(h) f(a) \ \ \ \forall  h \in \H(\R)\}.
 $$
 The principal affine space ${\cal A}_\G$ carries a canonical volume form $\mu$. It is $\G(\R)-$invariant,   and covariant 
 with the weight $2\rho$ under the action of the group $\H(\R)$. So there is a canonical half-density $|\mu|^{1/2}$, and we get 
  a unitary representation ${\cal V}^{\rm un}$ of the group $\G(\R)$ in 
  the Hilbert space of half-densities  $L_2({\cal A}_\G(\R), |\mu|^{1/2})$, called the unitary principal series.\footnote{Using the canonical half-density $|\mu|^{1/2}$, it is
   identified with the space of square integrable functions on ${\cal A}_\G(\R)$. } 
  Using the Mellin transform for the group $\H(\R)$ we decomposes it  into an integral of irreducible unitary representations 
  ${\cal V}_{\alpha+\rho}$, where  $\alpha$ are the unitary characters of the group $\H(\R)$:
    $$
  {\rm L_2}({\cal A}_\G(\R);  |\mu|^{1/2}) = \int {\cal V}_{\alpha+\rho} d\alpha.
    $$
 The $n$-th tensor power   of the principal series representation of the group $\G(\R)$ is realized  in the space of functions 
 on the $n-$th power of the space ${\cal A}_\G(\R)$:
   $$
  {\cal V} \otimes \ldots \otimes {\cal V} \stackrel{\sim}{=}   {\cal F}({\cal A}_\G(\R) \times \ldots \times {\cal A}_\G(\R)).
    $$ 
  So the  space of $\G(\R)$-invariants in the $n$-th tensor power  of the principal series representation ${\cal V}$  is realized in the space of functions on the coinvariants 
  of the action of  the group $\G(\R)$ on the space ${\cal A}_\G(\R) \times \ldots \times {\cal A}_\G(\R)$, which is the   space  ${\rm Conf_n}({\cal A}_\G)(\R)$: 
   \be \la{MSPX+}
    ({\cal V} \otimes \ldots \otimes {\cal V})^{\G(\R)} \stackrel{\sim}{=}    {\cal F}({\rm Conf_n}({\cal A}_\G)(\R)). 
\ee

  
  \bl  \cite[Section 8]{FG03a} \la{8.20.18.1a} The moduli space ${\mathscr A}_{  \G, \P_n}$ is   identified with  the configuration space ${\rm Conf}_n^{\rm t}({\cal A}_{\G'})$ of  twisted 
  cyclic configurations of decorated flags for the simply-connected group $\G'$. 
  Picking a vertex of the polygon $\P_n$ we identify it with ${\rm Conf}_n({\cal A}_{\G'})$:
$$
  {\mathscr A}_{  \G, \P_n} =   {\rm Conf}_n^{\rm t}({\cal A}_{\G'}) \stackrel{\sim}{=} {\rm Conf}_n({\mathcal A}_{\G'}).
$$
   \el

Using Lemma \ref{8.20.18.1a} we present the space on the right of (\ref{MSPX+})  for the simply-connected group $\G'$ as follows:
 $$
 {\cal F}({\rm Conf_n}({\cal A}_{\G'})(\R)) =  {\cal F}( {\mathscr A}_{  \G, \P_n}(\R)).   
 $$
Therefore the  space  ${\cal F}( {\mathscr A}_{  \G, \P_n}(\R))$ is  the {classical multiplicity space}. Its ${\rm L_2}-$variant  looks as follows:
 \vskip 2mm

(2). {\it The classical unitary multiplicity space}: 
\be \la{cunmas}
{\rm L_2}( {\mathscr A}_{  \G, \P_n}(\R), |\mu|^{1/2}).   
\ee
 \vskip 2mm

\subsubsection{\it Conclusion} The spaces (\ref{cunmas}) and (\ref{HDela+}) are very similar, but yet different. 
The latter is the space of functions on the positive   component of the  space underlying the former. 
And the cluster volume form $\mu_{\mathscr A}$ is the Lebesgue form  in the logarithmic coordinates, which make sense only on the positive component.

Therefore the quasiclassical limit of the quantum multiplicity space (1) can be realized as the classical multiplicity space (2). 
This illustrates an 
 essential feature of the cluster representation theory: all vector spaces which appear naturally are not just vector spaces, but carry a rich extra structure 
 of representations of the quantized 
algebras of functions - e.g. the algebra $\mathcal{A}_\hbar({\mathscr X}_{\G, S^2_n})$ for the quantum multiplicity spaces.  

\medskip

  \subsection{Realizations of quantum group  principal series representations}  \la{SEC9.4}

 \medskip
 
 Let $\Delta$ be a triangle with two colored sides $(12)$ and $(23)$. 

\bd \la{DEF13.12} The moduli space   ${\mathscr P}_{\G, \Delta}= {\mathscr P}_{ \Delta}$ 
parametrizes the following  data
$
(\B_1, \B_2, \B_3; p_{12}, p_{13}), 
$ 
 considered modulo the $\G$-action, where $(\B_1, \B_2, \B_3)$ is a triple of Borel subgroups in $\G$ which are pairwisely generic,  and   $p_{12}$ and 
$p_{23}$ are pinnings over $(\B_1, \B_2)$ and $(\B_2, \B_3)$ respectively, illustrated on the left:
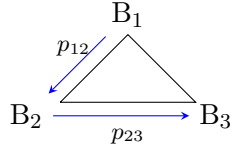
\ed

\begin{figure}[ht]
\epsfxsize 200pt
\center{
\begin{tikzpicture}[scale=.9]
\node (A) at (-.5, -0.2) {$\B_2$}; 
\node (B) at (1, 1.3) {$\B_1$}; 
\node (C) at (2.3, -0.2) {$\B_3$}; 
\draw[blue, stealth-] (A) --(B);
\draw[blue, -stealth] (A) --(C);
\draw (0,0) -- (1,1) -- (2,0) -- (0,0);   
\node at (1,-.5) {{\footnotesize $p_{23}$}};
\node at (0.2,.8) {{\footnotesize $p_{12}$}}; 
\end{tikzpicture} }
\label{pin14}
\end{figure}

 Forgetting the pinnings over $(\B_1, \B_3)$ we get  a projection  ${\mathscr P}_{\G, t} \lra {\mathscr P}_{\Delta}$. It induces a cluster Poisson structure on ${\mathscr P}_{\Delta}$.
 The set of   frozen vertices of   its  quiver ${\bf Q}_\Delta$  is  a union of two subsets,  
each containing $r$  elements, corresponding to the two special sides   of the triangle.  The quiver for ${\rm PGL_3}$ is pictured below. 
\begin{figure}[ht]
\epsfxsize 200pt
\center{
\begin{tikzpicture}[scale=.9]
\node (a1) at (120:1) {$\circ$};
\node (a2) at (60: 1) {};
\node (b1) at (180: 1) {$\circ$};
\node (b2) at (0,0) {$\bullet$};
\node (b3) at (0: 1) {};
\node (c1) at (240: 1) {$\circ$};
\node (c2) at (300: 1) {$\circ$};
\foreach \from/\to in { b2/b1, b1/c1, a1/b2, b2/c2,   c1/b2}
                 \draw[directed] (\from) -- (\to);   
\foreach \from/\to in {b1/a1,  c2/c1}
                 \draw[directed, dashed] (\from) -- (\to);        
\end{tikzpicture}}
\label{pin14*}
\end{figure}
 It has one unfrozen vertex, obtained by taking the  quiver   shown  after (\ref{CORT1}), and deleting the right side.

Denote by $\mathscr{D}_{\G, \Delta} = \mathscr{D}_{\Delta} $  the cluster symplectic double of the cluster Poisson variety ${\mathscr P}_{\Delta}$, and by  ${\cal O}_q(\mathscr{D}_{\Delta})$ its  quantized algebra. Let ${\bf B}_{12}$ and ${\bf B}_{23}$ be the subalgebras of  ${\cal O}_q(\mathscr{D}_{\Delta})$ generated by the frozen $\B$-variables at the sides $(12)$ and $(23)$ respectively. Note that cluster mutations keep ${\bf B}_{12}$ and ${\bf B}_{23}$ intact. Let  ${\cal O}_q(\mathscr{D}_{\Delta})[\sqrt{{\bf B}_{12}}]$ be the extension of the algebra ${\cal O}_q(\mathscr{D}_{\Delta})$ by adding square roots $\B_i^{1/2}$ for all frozen $\B_i$ at $(12)$. We further consider the centralizer of ${\bf B}_{23}$ in ${\cal O}_q(\mathscr{D}_{\Delta})[\sqrt{{\bf B}_{12}}]$:
\[
{\cal O}_q(\widetilde{\mathscr{D}}_{\Delta}^{\bf b}) :=\left({\cal O}_q(\mathscr{D}_{\Delta})[\sqrt{{\bf B}_{12}}]\right)^{{\bf B}_{23}}.
\]
Geometrically, let $\widetilde{\mathscr{D}}_{\Delta}$ be the $2^r$-to-$1$ covering space of $\mathscr{D}_{\Delta}$ by allowing taking square roots of frozen $B$-coordinates at the side $(12)$. The modified cluster symplectic double $\widetilde{\mathscr{D}}_{\Delta}^{\bf b}$ is obtained from $\widetilde{\mathscr{D}}_{\Delta}$ by forgetting the frozen $Y$-coordinates at the side $(23)$. 

Recall the moduli space $\widetilde {\mathscr L}_{\odot}:= \widetilde {\mathscr L}_{\G, \odot}$, defined as the fiber over $e$ of the outer monodromy map:
$$
\widetilde {\mathscr L}_{\G, \odot}:= \mu_{\rm out}^{-1}(e), \qquad \mu_{\rm out}: {\mathscr P}_{\G, \odot} \lra \H.
$$ 
 It admits a quantization $\mathcal{O}_q(\widetilde {\mathscr L}_{\odot})$ containing the quantum group $\mathcal{U}_q(\mathfrak{g})$ as a subalgebra.

  \vskip 2mm

Any unfrozen vertex $k$ of the  quiver   ${\bf Q}_{\Delta}$ 
gives rise a mutation of several  different flavors:

\vskip 1mmi) A mutation $\mu_k^{\mathscr X}$ of the cluster Poisson variety ${\mathscr P}_\Delta$.  

\vskip 1mmii)  A mutation $\mu_k^{\cal Z}:= \mu_k^{\mathscr X}\circ \mu_k^{{\mathscr X}^\circ}$ of the cluster Poisson variety ${\cal R}_\odot$.

\vskip 1mmiii) A mutation $\mu_k^{\mathscr{D}}$ of the  modified cluster symplectic double   $\widetilde{\mathscr{D}}_\Delta^{\bf b}$.

\bt  \la{Th9.13} There is a canonical    injective map of   algebras:
\be \la{12.28.18.4}
\xi_q^*: {\cal O}_q(\widetilde {\mathscr L}_{\odot} ) \lra {\cal O}_q({\widetilde{\mathscr{D}}}^{\bf b}_{\Delta}).
\ee

\vskip 1mm1. It is   a monomial map in each cluster coordinate system for ${\mathscr P}_\Delta$. 

\vskip 1mm2. It intertwines  cluster transformations ${\bf c}^{\cal Z}$ and ${\bf c}^\mathscr{D}$ generated by  the  mutations $\mu_k^{\cal Z}$ and $\mu_k^\mathscr{D}$.

\vskip 1mm 3.  It induces  canonical $W$-equivariant map
$
  \mu^*{\cal O}_\H   \stackrel{}{\lra} {\bf B}_{23}. 
   $
\et

\begin{proof}  Following the construction of \eqref{double.maps.advfs}, there is a surjective map
\be
\la{adsovbd}
\begin{split}
&\alpha:~\mathscr{D}_\Delta \lra {\mathscr P}_{\Delta} \times {\mathscr P}_{\Delta^\circ} \\
&\alpha^*(X_i)= Y_i, \qquad\alpha^*(X_{i^\circ}) = Y_i^{-1}\prod_j B_j^{-p_{ij}}.\\
\end{split}
\ee
Its image is denoted by the top and bottom triangles on the left of Figure \ref{pin13}. 
Let $B_i^{(12)}$ and $B_i^{(23)}$ be   coordinates on  $\mathscr{D}_\Delta$ assigned to  frozen vertices on the sides $(12)$ and $(23)$ of  $\Delta$ respectively. Set
$$
h:= \left(\prod_{i=1}^r\alpha_i^\vee(B_i^{(12)})\right)^{-1/2},\hskip 7mm h':= \prod_{i=1}^r\alpha_i^\vee(B_i^{(23)}).
$$
Let $p_{12}=(\A_1, \B_2)$ and $p_{12}'=(\A_1', \B_2')$ be the pinnings as shown on the figure. We construct two extra pinnings $p_{13}$ and $p_{13}'$ such that
\be
\la{bscaobsaon}
p_{13} =(\A_1\cdot h, ~\B_3), \hskip 7mm p_{13}'=(\A_1'\cdot h', ~ \B_3').
\ee

\begin{figure}[ht]
\epsfxsize 200pt
\center{
\begin{tikzpicture}[scale=.9]
\begin{scope}[shift={(0,1)}]
\node (A) at (-.5, -0.2) {$\B_2$}; 
\node (B) at (1, 1.3) {$\B_1$}; 
\node (C) at (2.5, -0.2) {$\B_3$}; 
\draw[blue, stealth-] (A) --(B);
\draw[blue, -stealth] (A) --(C);
\draw[blue, dashed, -stealth] (B) --(C);
\draw (0,0) -- (1,1) -- (2,0) -- (0,0);   
\node at (1,-.5) {{\footnotesize $p_{23}$}};
\node at (0.2,.8) {{\footnotesize $p_{12}$}}; 
\node at (1.8,.8) {{\footnotesize $p_{13}$}}; 
\end{scope}
\begin{scope}[shift={(0,-1)}]
\node (A) at (-.5, 0.2) {$\B_2'$}; 
\node (B) at (1, -1.3) {$\B_1'$}; 
\node (C) at (2.5, 0.2) {$\B_3'$}; 
\draw[blue, stealth-] (A) --(B);
\draw[blue, -stealth] (A) --(C);
\draw[blue, dashed, -stealth] (B) --(C);
\draw (0,0) -- (1,-1) -- (2,0) -- (0,0);   
\node at (1,.5) {{\footnotesize $p_{23}'$}};
\node at (0.2,-.8) {{\footnotesize $p_{12}'$}}; 
\node at (1.8,-.8) {{\footnotesize $p_{13}'$}}; 
\end{scope}
\draw (5,-.02) -- (6,-1.02) -- (7,-.02) -- (5,-.02);  
\draw (5,.02) -- (6,1.02) -- (7,.02) -- (5,.02); 
\node[red] at (6, -1.02) {$\bullet$};
\draw[blue, stealth-] (4.8, .15) --(5.7, 1.05);
\draw[blue, stealth-] (4.8,-.15) --(5.7,-1.05);
\draw[blue, dashed, -stealth] (6.3, 1.05) --(7.2, .15);
\draw[blue, dashed, -stealth] (6.3, -1.05) --(7.2, -.15);
\draw[ultra thick, -stealth] (3,0)--(4.2,0);
\draw[ultra thick, -stealth] (8,0)--(9,0);
\begin{scope}[shift={(0,-1.5)}]
\draw (11,1.5) circle (8mm);
\draw (11, 1.3) ellipse (3mm and 6mm);
\draw (11,1.5) -- (11, .7);
\node[red] at (11, 1.5) {{\small $\bullet$}};
\node[black] at (11, .7) {{\small $\bullet$}};
\node[black] at (11, 2.3) {{\small $\bullet$}};
\draw[blue, stealth-] (10.95, .4) arc (270:90:1.1);
\draw[blue, dashed, -stealth] (11.05, 2.6) arc (90:-90:1.1);
\end{scope}
\end{tikzpicture}
 }
\caption{Gluing a pair of triangles with two pinnings into a punctured disc with two pinnings.}
\label{pin13}
\end{figure}

Let us glue two triangles $\Delta$ on the left of Figure \ref{pin13} along the pinnings $p_{23}$ and $p_{23}'$, obtaining a moduli space ${\mathscr P}_{\square_4}$ assigned to the rectangle $\square_4$ with four pinnings  in the middle of Figure \ref{pin13}. We further glue the two bottom sides of ${\rm R}$ along the pinnings $p_{12}'$ and $p_{13}'$, obtaining a moduli space $\mathscr{P}_{\odot}$. A direct calculation shows that the outer ${\rm H}$-monodormy is 1 in this construction.
So our    construction gives a dominant map
\[
\xi: \widetilde{\mathscr{D}}^{\bf b}_{\Delta} \lra \widetilde {\mathscr L}_{\odot}.
\]
By the constructions \eqref{adsovbd} and \eqref{bscaobsaon}, it is easy to show that $\xi$ is a monomial map in each cluster coordinate system for ${\mathscr P}_{\Delta}$. A further check shows that $\xi$ is a Poisson map. Therefore its quantization gives rise to the map \eqref{12.28.18.4}. 
Properties 2)-3) follow by construction.
\end{proof}

\subsubsection{Example: realization of ${\U}_q(\mathfrak{sl}_2)$.}
The quiver associated with the triangle is illustrated by the  left on Figure \ref{T3a}, and the quiver for its symplectic double is shown by the right. The latter encodes  a rank $4$ lattice with a basis $e_1, e_2, b_1, b_2$ with the following skew-symmetric form
\[
(e_1, e_2)=1, \qquad (b_1, b_2)=0, \qquad (e_i, b_j)=\delta_{ij}.
\]
\begin{figure}
\begin{tikzpicture}
\draw (90:1)--(210:1)--(330:1)--(90:1);
\node[red] (a) at (150:0.5) {$\bullet$};
\node[red] (b) at (270:0.5) {$\bullet$};
\draw[blue, thick, -latex] (a)--(b);
\node at (150:0.8) {$e_1$};
\node at (270:0.8) {$e_2$};

\node[red] (e1) at (3,1) {$\bullet$};
\node[red] (e2) at (3, -.5) {$\bullet$};
\node[red] (b1) at (4.5, 1) {$\bullet$};
\node[red] (b2) at (4.5, -.5) {$\bullet$};
\draw[blue, thick, -latex] (e1)--(e2);
\draw[blue, thick, -latex] (e1)--(b1);
\draw[blue, thick, -latex] (e2)--(b2);
\node at (2.8, 1.2) {$e_1$};
\node at (2.8, -.7) {$e_2$};
\node at (4.7, 1.2) {$b_1$};
\node at (4.7, -.7) {$b_2$};
\end{tikzpicture}
\caption{On the left: the quiver for the   ${\rm PGL}_2-$triangle with two pinnings. On the right:  the quiver for its cluster symplectic double.}
\label{T3a}
\end{figure}

The lattice gives rise to a quantum torus algebra, with generators
\[
{\rm Y}_1={\rm T}_{e_1}, \qquad {\rm Y}_2={\rm T}_{e_2}, \qquad {\rm B}_1= {\rm T}_{b_1}, \qquad {\rm B}_2={\rm T}_{b_2}. 
\]
Let us apply the twisted double construction. It first gives rise to two triangles with the variables shown on left of Figure \ref{T3}. 
After gluing, we obtain the quiver as the right picture. Note that the vector $e_2$ is cancelled after gluing. 
The generators of ${\U}_q(\mathfrak{sl}_2)$  are presented as
\be \la{SL}
{\bf E}= {\rm T}_{e_1}, \qquad {\bf K}={\rm T}_{-b_1}, \qquad {\bf F}= {\rm T}_{-e_1-b_1}+{\rm T}_{-e_1+b_2}+ {\rm T}_{-e_1-b_2}+{\rm T}_{-e_1+b_1}.
\ee
In particular, we get the following central element of ${\U}_q(\mathfrak{sl}_2)$:
\[
{\bf C}:= {\bf E}{\bf F}- q^{-1}{\bf K}-q{\bf K}^{-1}= {\rm T}_{b_2}+{\rm T}_{-b_2}.
\]
\begin{figure}[h]
\begin{tikzpicture}
\draw (90:1)--(210:1)--(330:1)--(90:1);
\node[red] (a) at (150:0.5) {$\bullet$};
\node[red] (b) at (270:0.5) {$\bullet$};
\node[red] (c) at (30:0.5) {$\bullet$};
\node at (140:0.8) {$e_1$};
\node at (270:0.8) {$e_2$};
\node at (1.1, 0.5) {$-e_1-b_1$};
\begin{scope}[shift={(0,-2.3)}]
\draw (-90:1)--(-210:1)--(-330:1)--(-90:1);
\node[red]  at (-150:0.5) {$\bullet$};
\node[red]  at (-270:0.5) {$\bullet$};
\node[red]  at (-30:0.5) {$\bullet$};
\node at (90:0.8) {$-e_2+b_1+b_2$};
\node at (-1.6,-.25) {$-e_1+b_1-b_2$};
\node at (1.6,-0.25) {$e_1-b_1-b_2$};
\end{scope}
\begin{scope}[shift={(5,0)}]

\node[red]  (d) at (150:0.8) {$\bullet$};
\node[red]  (e) at (270:0.8) {$\bullet$};
\node[red]  (f) at (30:0.8) {$\bullet$};
\draw[thick, blue, -latex] (f) --(d);
\draw[thick, blue, -latex] (d)--(e);
\draw[thick, blue, -latex] (e)--(f);
\node[red] at (0, -1.9) {$\bullet$};
\node at (150:1.1) {$e_1$};
\node at (270:1.1) {$b_1+b_2$};
\node at (1.5, 0.55) {$-e_1-b_1$};
\node at (0, -2.3) {$-2b_2$};
\end{scope}

\begin{scope}[shift={(-8.5,-3)}, scale=1.7]
\draw (11,1.5) circle (8mm);
\draw (11, 1.2) ellipse (4mm and 5mm);
\draw (11,1.5) -- (11, .7);
\node[red] at (11, 1.5) {{\small $\times$}};
\node[black] at (11, .7) {{\small $\bullet$}};
\node[black] at (11, 2.3) {{\small $\bullet$}};
\node at (10, 1.5) {$e_1$};
\node at (12.3, 1.5) {$-e_1-b_1$};
\node at (11, 2.6) {${\bf E}$};
\node at (11, .4) {${\bf F}$};
\node at (11, 1.2) {$-2b_2$};
\node at (11, 1.8) {$b_1+b_2$};
\draw[red, thick] (11, 0.9) arc (90:70:1.3);
\draw[red, thick, -latex] (11, 0.9) arc (90:110:1.3);
\draw[red, thick,-latex] (11, 2.1) arc (270:290:1.3);
\draw[red, thick] (11, 2.1) arc (270:250:1.3);
\end{scope}

\end{tikzpicture}
\caption{The twisted double construction for the ${\rm PGL}_2-$triangle with two pinnings. The right picture is dual to the middle graph.}
\label{T3}
\end{figure}

Using notation  (\ref{handy}) and representation (\ref{5.26.08.1aa}) of the Heisenberg algebra in ${\rm L}_2({\cal A}_{\rm log}(\R))$, we get: 
  \begin{equation} \label{5.26.08.1aa}
\begin{split}
& y_k:=   2\pi i  \hbar \cdot  \frac{\partial}{\partial a_k}
-\alpha_k^+,\\ 
&b_k:= a_k.\\
\end{split}
\end{equation} 
So the  exponentials $ {\rm Y}_p:= {\rm exp}( y_p)$ and  $ \B_p:= {\rm exp}( b_p)$ act in  ${\rm L}_2({\cal A}(\R_{>0}))$ as  
 difference operators: 
\begin{equation} \nonumber
 {\rm Y}_p f (A_1, ..., A_n) = ({\mathbb A}^{+}_p)^{-1} \cdot f(A_1, ..., q^2A_p, 
\ldots A_n), ~\qquad~ 
\B_pf= A_pf.
\end{equation}

Then the quantum group acts on the  space ${\rm L}_2(\mathcal{A}_{\rm SL_2}(\mathbb{R}_{>0}))={\rm L}_2(\mathbb{R}_{>0}^2)$ as follows:
\be \la{4eq}
\begin{array}{ll}
{\bf E}\, f(A_1, A_2)&= A_2^{-1} f(q^2 A_1, A_2),\\
{\bf K}\, f(A_1, A_2)&= A_1^{-1}f(A_1, A_2),\\
{\bf F}\, f(A_1, A_2)&= (qA_1^{-1}+A_2+A_2^{-1}+ q^{-1}A_1) A_2 f(q^{-2}A_1, A_2),\\
{\bf C}\, f(A_1, A_2)&= (A_2+A_2^{-1}) f(A_1, A_2).\\
\end{array}
\ee

  \vskip 2mm \subsubsection{\it Realization of   principal series   $\ast-$representations of  the modular double  $ \mathcal{A}_\hbar({\mathfrak{g}})$ of the quantum group.}       Consider the following compositions,  
where $\kappa_\hbar:= \kappa_q \otimes \kappa_{q^\vee}$ and $\xi_\hbar:= \xi_q \otimes \xi_{q^\vee}$:
\be \la{16.7.1}
\begin{split}
&{\cal U}_q({\mathfrak{g}}) \stackrel{\kappa_q}{\hlra} {\cal O}_q(\widetilde
 {\mathscr L}_{\G, \odot}) \stackrel{\xi_q}{\lra} {\cal O}_q(\widetilde{\mathscr{D}}_{\G, \Delta})^{{\bf B}_{23}};\\
 &\mathcal{A}_\hbar({\mathfrak{g}}) \stackrel{\kappa_\hbar}{\hlra} 
 \mathcal{A}_\hbar({\cal R}_{\G, \odot}) \stackrel{\xi_\hbar}{\lra} \mathcal{A}_\hbar(\widetilde{\mathscr{D}}_{\G, \Delta})^{{\bf B}_{23}}.\\
 \end{split}
 \ee
 Let ${\mathscr A}_{\G, \Delta}\subset {\rm Conf}_3(\mathcal{A})$ be the space parametrizing $\G$-orbits of pairwisely generic triples $(\A_1, \A_2, \A_3)$ of decorated flags such that 
  \[
  h(\A_1, \A_3)=1.
  \]
The space ${\mathscr A}_{\G, \Delta}$ inherits a cluster $K_2-$structure from ${\rm Conf}_3(\mathcal{A})$ by deleting the frozen vertices on side $(13)$. The pair $({\mathscr A}_{\G, \Delta}, {\mathscr P}_{\G, \Delta})$ form a cluster ensemble. 
  Recall the canonical     representation of  the $\ast$-algebra $\mathcal{A}_\hbar(\widetilde{\mathscr{D}}_{\G, \Delta})$    in the Hilbert space  
\be \la{HDel}
   {\rm L_2}\left( {\mathscr A}_{  \G,\Delta}(\R_{>0}), \mu_{\mathscr A}\right).
\ee 
  The action of the commutative subalgebra ${\bf B}_{23}$ 
 commutes with the action of $\kappa_\hbar(\mathcal{A}_\hbar({\mathfrak{g}}))$. 
Decomposing the representation ${\cal H}_{\Delta}$ according to the points $\lambda\in \H(\R_{>0})$ we get the principal series of  representations ${\cal V}_\lambda$ of 
the  $\ast$-algebra $ \mathcal{A}_\hbar({\mathfrak{g}})$.  
  Lemma \ref{8.20.18.1} relates    Hilbert space 
  (\ref{HDel})  with the representation space  of the  classical principal series for the group $\G(\R)$.  

\begin{figure}[ht]
\epsfxsize 200pt
\center{
\begin{tikzpicture}[scale=.9]
\node at (-.3, -0.1) {$\mathcal{B}$}; 
\node at (1, 1.2) {$\mathcal{B}$}; 
\node at (2.3, -0.1) {$\mathcal{B}$}; 
\node at (1,-0.8) {{ \small Poisson variety}};
\draw[blue, stealth-] (-.1, .2) --(.7,1);
\draw[blue, -stealth] (0, -.2) --(2,-.2);
\draw (0,0) -- (1,1) -- (2,0) -- (0,0);     
\draw[thick, latex-latex] (3,0.5) -- (4,0.5);
\node at (3.5,0.8) {{ \small dual}};
\begin{scope}[shift={(5,0)}]
\node at (-.3, -0.1) {$\mathcal{A}$}; 
\node at (1, 1.2) {$\mathcal{A}$}; 
\node at (2.3, -0.1) {$\mathcal{A}$}; 
\node at (1,-0.8) {{ \small $K_2$ variety}};
\draw[red, stealth-stealth, dashed] (2.1, .2) --(1.3,1);
\node[red] at (1.8,0.8) {{ \small 1}};
\draw (0,0) -- (1,1) -- (2,0) -- (0,0); 
\node[thick] at (3.5,0.5) {$\sim$};
\begin{scope}[shift={(5,0)}]
\node at (-.3, -0.1) {$\mathcal{A}$}; 
\node at (1, 1.2) {$\mathcal{B}$}; 
\node at (2.3, -0.1) {$\mathcal{A}$}; 
\draw (2,0) -- (0,0) -- (1,1) ;  
\draw[dashed] (2,0) -- (1,1);
\node[thick] at (3.5,0.5) {$=$};
\node[thick] at (4.5,0.5) {$\mathcal{A}$};
\end{scope}
\end{scope}
\end{tikzpicture}
 }
\caption{ The cluster dual ${\mathscr A}_{\G, \Delta}$ to the cluster Poisson variety ${\mathscr P}_{\G, \Delta}$ is identified with  ${\mathscr A}_\G$.}
\label{pin15}
\end{figure}

  \bl \la{8.20.18.1}The cluster $K_2$-variety ${\mathscr A}_{\G, \Delta}$      is naturally identified, by a birational isomorphism respecting the cluster structure,  with  the principal affine space ${\mathscr A}_\G= \G/\U_-$. 
  \el
  
  \begin{proof}     
  
  The space ${\mathscr A}_{\G, \Delta}$ - illustrated by the second triangle on Figure \ref{pin15} - parametrizes 
  triples of decorated flags   $(\A_1, \A_2, \A_3)$    with $h(\A_1, \A_3)=1$    modulo the action of $\G$. 
  This space is identified with the space parametrizing triples 
  $(\B_1, \A_2, \A_3)$ modulo $\G$-action,  where the pair $(\B_1, \A_3)$ is generic - illustrated by the third triangle on Figure \ref{pin15}. 
 The group $\G$ acts simply transitively on the generic  pairs $(\B_1, \A_3)$. So assigning to a configuration $(\B_1, \A_2, \A_3)$  the decorated flag $\A_2$ we get the identification with   ${\mathscr A}_\G$:
\be \nonumber
\begin{split}
 & {\mathscr A}_{\G, \Delta} \stackrel{\sim}{\lra} \G/\U^-, \qquad (\B^+, \A_2, \U^-) \lms \A_2,\\
  \end{split}
 \ee 
 \end{proof} 
 
 Now we can compare the realizations of the classical and quantum principal series representations. 
 \vskip 2mm
 
 (1). {\it The  Gelfand-Naimark unitary principal series  representations of the group $\G(\R)$ are realized   in  
  the Hilbert space of half-densities  on the principal affine space}:
\be \la{PAS+}
 { \rm L_2}({\cal A}_\G(\R), |\mu|^{1/2}).
\ee
  \vskip 2mm

 On the other hand, recall  the map of $\ast-$algebras given by (\ref{16.7.1}):   
 \be 
 \mathcal{A}_\hbar({\mathfrak{g}}) \stackrel{ }{\hlra}  \mathcal{A}_\hbar(\widetilde{\mathscr{D}}_{\G, \Delta})^{{\bf B}_{23}}. 
 \ee 
Combining it with the canonical     $\ast-$representation of  the $\ast$-algebra $\mathcal{A}_\hbar(\widetilde{\mathscr{D}}_{\G, \Delta})$    in  
  (\ref{HDel}), we arrive at   \vskip 2mm
 
(2).  {\it The principal series $\ast-$representation of the modular double $ \mathcal{A}_\hbar({\mathfrak{g}})$ of the quantum group  can be realized   
in the cluster Schwartz space of the cluster Hilbert space} 
\be \la{HDel++}
{\rm L_2}\left( {\mathscr A}_{  \G,\Delta}(\R_{>0}), \mu_{\mathscr A}\right).
\ee 
\vskip 2mm

Lemma \ref{8.20.18.1}   identifies ${\cal A}_{\G'}$ for the simply-connected group $\G'$ with the cluster $K_2-$variety ${\mathscr A}_{\G', \Delta}$:
 $$
  {\cal A}_{\G'} =   {\mathscr A}_{\G', \Delta}.
 $$

\subsubsection{\it Conclusion} The classical and quantum principal series representations are realized in the space of functions   (\ref{PAS+}) and (\ref{HDel++})   on essentially the same space. Precisely, in the classical case this is the space of real points of the principal affine space, and in the quantum case the space of the positive points of a cluster $K_2-$variety, isomorphic to the principal affine space by Lemma \ref{8.20.18.1}, just as we have seen  in Section \ref{sec16.5:RMSP} for the multiplicity spaces. The cluster variety structure is essential for the definition of the  positive points.  
Therefore  the classical unitary principal series  (1) can be realized as the quasiclassical limit of the principal series $\ast-$representation of the quantum group (2). 

This shows again that the vector spaces in the cluster representation theory always carry an extra structure 
 of representations of the quantized 
algebras of functions on the underlying spaces.

 \subsection{A version of the cluster double realization} \la{SECTa17.8}

Let ${\bf q}$ be a quiver with vertices labelled by a set ${\rm V}$. Let $\widetilde{\bf q}$ be an extension of ${\bf q}$ by adding more frozen variables labelled by the set $S$.  Let $\widetilde {\bf q}^{\circ}$ be its opposite quiver. We glue $\widetilde {\bf q}$ with $\widetilde {\bf q}^\circ$ along the newly added frozen vertices, obtaining a new quiver
\[
\widetilde{\bf q}_{\bf d}:=\widetilde{\bf q}\cup_{S}\widetilde{\bf q}^\circ.
\]  

We consider the following two quantum algebras
\begin{itemize}
\item $\mathcal{O}_q(\mathscr{D}_{|{\bf q}|})$ is the quantum double cluster algebra associated with ${\bf q}$, 
\item $\mathcal{O}_q(\mathscr{X}_{|\widetilde{\bf q}_{\bf d}|})$ is the quantum cluster Poisson algebra associated with the quiver $\widetilde{\bf q}_{\bf d}$.
\end{itemize}

Recall that $\mathcal{O}_q(\mathscr{D}_{|{\bf q}|})[\sqrt{{\rm B}_f}]$ is an extension 
of $\mathcal{O}_q(\mathscr{D}_{|{\bf q}|})$ that adds  square roots of frozen ${\rm B}_f$'s. Precisely, let ${\rm V}^f$ and ${\rm V}^u$ be the set of frozen and unfrozen vertices of ${\bf q}$ respectively. Consider the double lattice
\[
  \Lambda_{\bf q}^\mathscr{D}:= \Lambda_{\bf q} \bigoplus (\oplus_{i\in {\rm V}^u}\mathbb{Z}b_i) \bigoplus (\oplus_{i\in {\rm V}_f}\frac{1}{2}\mathbb{Z} {b_i})
\]  
such that 
\[
(b_i, b_j)=0, \qquad (b_i, e_j)=-\delta_{ij}.
\]  
Locally, that is for each quiver ${\bf q}$, it gives rise to a quantum torus algebra  $\mathbb{T}_{\Lambda_{\bf q}^\mathscr{D}}$ generated by 
\[{\rm Y}_i:= {\rm T}_{e_i},~~ \forall i \in {\rm V};
\qquad {\rm B}_{i}:= {\rm T}_{b_i},  ~~ \forall i \in {\rm V}^u;
\qquad {\rm B}_i^{\frac{1}{2}}= {\rm T}_{\frac{b_i}{2}}, ~~ \forall i \in {\rm V}^f. 
\]
Recall that 
\[
\widetilde{\rm Y}_i= {\rm T}_{-e_i-\sum_k \varepsilon_{ik}b_k}= {\rm Y}_i^{-1}\prod_{k}{\rm B}_k^{-\varepsilon_{ik}}.
\]
By definition, the algebra $\mathcal{O}_q(\mathscr{D}_{|{\bf q}|})[\sqrt{{\rm B}_f}]$ is the algebra of universally Laurent polynomials in every $\mathbb{T}_{\Lambda_{\bf q}^\mathscr{D}}$ under the  mutations \eqref{zx1qt}.

\bp 
There is a canonical algebra homomorphism
\[
\xi^\ast: ~ \mathcal{O}_q(\mathscr{X}_{|\widetilde{\bf q}_{\bf d}|}) \longrightarrow \mathcal{O}_q(\mathscr{D}_{|{\bf q}|})[\sqrt{{\rm B}_f}].
\]
It is induced in each cluster coordinate system by the following monomial map: 
 \begin{equation} 
 \la{24.527.922}
\xi^*: {\rm Z}_v \lms
\left\{ \begin{array}{lll} {\rm Y}_i
& \mbox{ if $v=i\in {\rm V}^u$},\\
\widetilde {\rm Y}_i& \mbox{ if $v = i^\circ$, $i \in  {\rm V}^u$},\\
 \prod_{j \in {\rm V}}{\rm B}_j^{-\varepsilon_{fj}}& \mbox{ if $v=f\in {\rm V}^f$.}\\
\end{array}\right.
 \end{equation} 
 It intertwines  cluster transformations ${\bf c}^{\cal Z}$ and ${\bf c}^\mathscr{D}$ generated by  unfrozen   mutations. 
\ep
 
 \begin{proof}
 This is a special case of Theorem \ref{THH8.9} where the newly added ${\rm B}$-variables associated with the vertices in $S$ are set to 1. 
\end{proof}  

Now we investigate the kernel of $\xi^\ast$.
Locally, the monomial map \eqref{24.527.922} corresponds to a linear map
\be
\label{24.5.26.15.55}
\xi_{\bf q}: ~ \Lambda_{\widetilde{\bf q}_d} \lra \Lambda_{\bf q}^\mathscr{D}.
\ee
Recall the Casimir sub-lattice
\[
{\rm Ker}~\Lambda_{\widetilde{\bf q}_{\bf d}}:= \left\{ v ~\middle |~ (v,u)=0, \forall u \in \Lambda_{\widetilde{\bf q}_{\bf d}}\right\}.
\]
\bl We have 
\[
{\rm Ker} ~\xi_{\bf q} \subseteq {\rm Ker}~\Lambda_{\widetilde{\bf q}_{\bf d}}.
\]
These two sub-lattices are equal if 
\[
{\rm dim}~{\rm Ker}\ \Lambda_{\widetilde{\bf q}_{\bf d}}= {\#}S.
\]
\el
\begin{proof} Let $v\in {\rm Ker} ~\xi_{\bf q}$ and let $u\in \Lambda_{\widetilde{\bf q}_{\bf d}}$. We have
\[
(v, u) = (\xi_{\bf q}(v), \xi_{\bf q}(u)) = (0, \xi_{\bf q}(u))=0.
\]
Therefore $v\in {\rm Ker}~\Lambda_{\widetilde{\bf q}_{\bf d}}.$  
If ${\rm dim}~{\rm Ker}\ \Lambda_{\widetilde{\bf q}_{\bf d}}= {\#}S$, then 
\[
{\rm dim} ~{\rm Ker} \ \xi_{\bf q} \geq \dim \Lambda_{\widetilde{\bf q}_{\bf d}} - \dim  \Lambda_{\bf q}^\mathscr{D} = \# S = {\rm dim}{\rm Ker}~\Lambda_{\widetilde{\bf q}_{\bf d}}. 
\]
The quotient 
$ 
{\rm Ker}~\Lambda_{\widetilde{\bf q}_{\bf d}} ~{\big \slash}~ {\rm Ker} ~\xi_{\bf q}
$ 
can be identified with a sub-lattice of $\Lambda_{\bf q}^{\mathscr{D}}$, and therefore is torsion free. The equality of the two lattices follows follows from the equality of their dimensions.  
\end{proof}  
  
  As a direct consequence, we get
  \bp
  \la{5.24.26.521}
   If the dimension of the Casimir torus of  $\mathcal{O}_q(\mathscr{X}_{\widetilde{\bf q}_{\bf d}})$ is equal to $\#S$, then we get a natural injective homomorphism
  \[
  \xi: \mathcal{O}_q(\mathscr{X}_{|\widetilde{\bf q}\cup_{S}\widetilde{\bf q}^\circ |}) {\big \slash} {\mathcal{I}} \longrightarrow \mathcal{O}_q(\mathscr{D}_{|{\bf q}|})[\sqrt{\rm B}_f],
  \]
  where the quotient modulo $\mathcal{I}$ corresponds to imposing the relation ${\rm C}=1$ for every Casimir element ${\rm C}$ of $\mathcal{O}_q(\mathscr{X}_{\widetilde{\bf q}_{\bf d}})$.
  \ep

  \subsubsection{A version of the quantum double construction for  surfaces.} Let $\bS$ be a partially colored surface. We further assume that 
   \begin{itemize}
   \item $\bS$ has no punctures,
   \item every boundary circle of $\bS$ has at least one non-colored boundary interval.
   \end{itemize}
  Associated with $\bS$ is the symplectic double $\mathscr{D}_{{\G, \bS}}$. In this way, we obtain the algebra
  $\mathcal{O}_q(\mathscr{D}_{{\G, S}})[\sqrt{{\rm B}_f}]$  
  
  Let $\widetilde{\bS}$ be the surface with every boundary interval colored.
 
  Let $\mathcal{D}_{\bS}$ be the topological double of $\bS$ obtained by gluing of $\bS$ with $\bS^{\circ}$ along non-colored intervals. 
  
  \begin{figure}
  \begin{tikzpicture}
  \draw[dashed] (18:1)--(90:1)--(162:1)--(234:1)--(306:1)--(18:1);
   \draw[thick, red] (18:1)--(90:1)--(162:1);
  \draw[thick, red] (234:1)--(306:1);
    \node at (0,0) {$\bS$};
   \begin{scope}[shift={(0,-3)}]
   \draw[dashed] (-18:1)--(-90:1)--(-162:1)--(-234:1)--(-306:1)--(-18:1);
    \draw[thick, red] (-18:1)--(-90:1)--(-162:1);
  \draw[thick, red] (-234:1)--(-306:1);
    \node at (0,0) {$\bS^\circ$};
   \end{scope}
   
   \begin{scope}[shift={(-2,0)}]
   \draw[thick,red] (7, -1.5) ellipse (1cm and 0.5cm);
    \draw[thick,red] (7, -1.5) ellipse (2cm and 1.5cm);
    \draw[dashed] (5, -1.5) --(6, -1.5);
    \draw[dashed] (8, -1.5) --(9, -1.5);
    \node at (7,0) {$\bullet$};
    \node at (5, -1.5) {$\bullet$};
    \node at (7, -3) {$\bullet$};
    \node at (9, -1.5) {$\bullet$};
    \node at (6,-1.5) {$\bullet$};
    \node at (8, -1.5) {$\bullet$};
    \node at (7, -0.5) {$\bS$};
    \node at (7, -2.5) {$\bS^\circ$};
    \end{scope}
    
     \begin{scope}[shift={(10,-1.5)}]
   \draw[thick, red] (18:1)--(90:1)--(162:1)--(234:1)--(306:1)--(18:1);
    \node at (0,0) {$\widetilde{\bS}$};
   \end{scope}
  \end{tikzpicture}
  \caption{The surface in the middle is the topological double  ${\cal D}_\bS$ of the  colored surface $\bS$  on the left. It is obtained by gluing the surfaces on the left via matching dashed, that is uncolored, boundary  intervals. The colored boundary intervals carry pinnings. 
  The decorated surface $\widetilde \bS$ 
  on the right is  obtained by coloring all boundary intervals of $\bS$. Note that the number of the boundary components of the middle surface coincides with the number of dashed  boundary intervals of $\bS$. The set $S$ parametrised frozen variables at the dashed intervals. So $\#S$ is the rank $r$ times the number of dashed boundary intervals on $\bS$. }
   \label{pin110++}
    \end{figure}

\bt 
\la{24.527.1013} 
(i) There is a natural injective homomorphism
\be \la{5.27.24.b}
\mathcal{O}_q(\mathscr{P}_{\G, {\cal D}_{\bS}}) {\big \slash} {{\rm C}=1} \longrightarrow \mathcal{O}_q(\mathscr{D}_{\G, \bS})(\sqrt{{\rm B}_f}),
\ee
where the left hand side is modulo the relation ${\rm C}=1$ for every Casimir element ${\rm C}$ of $\mathcal{O}_q(\mathscr{P}_{\G, {\cal D}_{\bS}})$.

(ii) Specialising $q=1$, we get  an  isomorphism of  spaces  of real positive points:
\be \la{isompoa}
  \mathscr{D}_{\G, \bS}(\R_{>0}) \stackrel{\sim}{\lra} \mathscr{P}_{\G, {\cal D}_{\bS}}(\R_{>0})_{\mu_{\rm out} = 1}.
\ee\et  
  
 \begin{proof} 
 (i) Let $k$ be the number of uncolored boundary intervals of $\bS$ and let $r$ be the rank of ${\rm G}$. 
 Let $\widetilde{\bf q}$ be a quiver of $\mathscr{P}_{\G, \widetilde{\bS}}$. It is an extension of a quiver ${\bf q}$ of $\mathscr{P}_{\G,\bS}$ by adding $kr$ more frozen vertices. 
 Meanwhile, we note that the topological double surface $\mathcal{D}_{\bS}$ has $k$ many boundary circles, where every boundary circle is either a puncture or consists of even number of marked points.  We see that the dimension of the Casimir torus of $\mathscr{P}_{\G, {\cal D}_{\bS}}$ is $kr$. Therefore the Theorem is a special case of Proposition 
 \ref{5.24.26.521}.

(ii) Follows from (i). Note that over real positive  numbers there is unique way to define $\sqrt{{\rm B}_f}$. \end{proof}

  \subsubsection{A version of the quantum double construction for  punctured surfaces.} Let $\bS$ be a partially colored surface which can have punctures. 
  We further assume that 
   \begin{itemize}
   \item every boundary circle of $\bS$ with special points has at least one non-colored boundary interval.
   \end{itemize}
As above, we have  the symplectic double $\mathscr{D}_{{\G, \bS}}$ and the algebra
  $\mathcal{O}_q(\mathscr{D}_{{\G, S}})[\sqrt{{\rm B}_f}]$ . 
  
  Let $\widetilde{\bS}$ be the surface with every boundary interval colored.
 
  Let $\mathcal{D}_{\bS}$ be the topological double of $\bS$ obtained by gluing  $\bS$ and $\bS^{\circ}$ along all matching  non-colored intervals,  and 
   all matching punctures. Note that for each puncture $p$ on $\bS$ there is a canonical isotopy class $\nu_p$ 
   of a simple loop on the surface ${\cal D}_\bS$, which we call {\it neck loop}.   
    It is given by gluing the boundaries of the holes corresponding to the punctures. 
  
Denote by $\widehat{\mathscr{P}}_{\G, {\cal D}_{\bS}}$ a variant of the  moduli space ${\mathscr{P}}_{\G, {\cal D}_{\bS}}$ where 
   in addition the data $({\cal L}, ...)$ of a $\mathscr{P}-$space we add a flat section of the associated local system ${\cal L}_{\cal B}$ of flag varieties over every neck loop.
 When $S$ is a surface with punctures, it was considered in \cite{FG14}. Combining the construction of Section \ref{SECTa17.8} with the   main construction of \cite{FG14}, 
we can define the   cluster symplectic  double coordinates on 
the subspace of the space ${\widehat{\mathscr{P}}_{\G, {\cal D}_{\bS}}}$  given by the condition that the outer monodromy along every boundary component is equal to $1$. 
In particular, we have the following. 

\bt 
\la{5.27.24}
There is a natural isomorphism of the sets of real positive points
  $$
  \mathscr{D}_{\G, \bS}(\R_{>0}) \stackrel{\sim}{\lra} \widehat{\mathscr{P}}_{\G, {\cal D}_{\bS}}(\R_{>0})_{\mu_{\rm out} = 1}.
  $$
\et

 \subsection{Realizations of the quantum group analog of the  regular representation}  \la{SEC16.8} 
  \medskip

The cylinder ${\rm Cl}_{2,2}$  with $2+2$ ordered special boundary points  is the topological double of the colored rectangle  
 $\square$ with  the two opposite colored sides,  shown by   the   vertical sides  on 
  Figure \ref{pin110+}.  

 The canonical representation of the quantum double from Section \ref{Sec8a} gives  
 the  cluster realization of the $\ast-$algebra $\mathcal{A}_\hbar({\cal D}_{\G, \square})$  
 in the cluster Schwarz subspace of the cluster Hilbert space for the    cluster variety ${\mathscr A}_{\G', \square}$:
\be \la{CLSCH}
 {\cal S}({\mathscr A}_{\G', \square}(\R_{>0}))\subset {\rm L}_2({\mathscr A}_{\G', \square}(\R_{>0})).
\ee

 On the other hand, applying Theorem \ref{24.527.1013}, we have the following.  
   \bt  \la{11.26.18.1Xaaaa}  
There are  injective  maps of   $\ast-$algebras, where $\xi^*_{{\mathscr A}_\hbar}: = \xi^*_{{\cal O}_q} \otimes \xi^*_{{\cal O}_{q^\vee}}$: 
\be \la{MAPZEaaa+}
\begin{split}
&\xi^*_{{\cal O}_q}:     {\cal O}_q({\mathscr P}_{\G,{\rm Cl}_{2,2}  } ){\big \slash} {{\rm C}=1}\hra {\cal O}_q({\cal D}_{\G,  \square})(\sqrt{{\rm B}_f}),\\
&\xi^*_{{\mathscr A}_\hbar}:     \mathcal{A}_\hbar({\mathscr P}_{\G, {\rm Cl}_{2,2} } ){\big \slash} {{\rm C}=1} \hra \mathcal{A}_\hbar({\cal D}_{\G, \square})(\sqrt{{\rm B}_f}).\\ 
\end{split}
\ee
\et  
 The composition of the map $\xi^*_{{\mathscr A}_\hbar}$ with the cluster realization of the algebra $\mathcal{A}_\hbar({\cal D}_{\G, \square})$ 
 provides a representation of the $\ast-$algebra 
 $\mathcal{A}_\hbar({\mathscr P}_{\G, {\rm Cl}_{2,2} } ) $ in the cluster Schwartz space (\ref{CLSCH}). 
 This representation kills the outer monodromy and the extra central elements, and is exact on the quotient 
 $\mathcal{A}_\hbar({\mathscr P}_{\G, {\rm Cl}_{2,2} } ) {\big \slash} {{\rm C}=1}$.
  
   \begin{figure}[ht]
 \begin{center}
 \begin{tikzpicture}[scale=0.6]
   \draw[red, thick] (0,0) -- (0,1);
   \draw[red, thick] (3,0) -- (3,1);
   \draw[dashed] (0,0)-- (3,0);
   \draw[dashed] (0,1) -- (3,1);
    \node  at (0,0) {\tiny $\bullet$}; 
 \node  at (0,1) {\tiny $\bullet$}; 
  \node  at (3,0) {\tiny $\bullet$}; 
   \node  at (3,1) {\tiny $\bullet$}; 
   \node at (4,.5) {$\ast$};
   \begin{scope}[shift={(5,0)}]
      \draw[red, thick] (0,0) -- (0,1);
   \draw[red, thick] (3,0) -- (3,1);
   \draw[dashed] (0,0)-- (3,0);
   \draw[dashed] (0,1) -- (3,1);
    \node  at (0,0) {\tiny $\bullet$}; 
 \node  at (0,1) {\tiny $\bullet$}; 
  \node at (3,0) {\tiny $\bullet$}; 
   \node  at (3,1) {\tiny $\bullet$}; 
   \end{scope}
   \draw[-latex, thick] (9, .5) -- (10, .5);
   \begin{scope}[shift={(11,0.5)}, scale=.5]
    \draw[thick, red] (0,0) ellipse (0.5 and 1.25);
\draw (0, -1.25) -- (5, -1.25);
\draw[thick, red] (5, -1.25) arc (-90:90:0.5 and 1.25);
\draw [thick, red, dashed] (5, -1.25) arc (270:90:0.5 and 1.25);
\draw (5, 1.25) -- (0, 1.25);  
\node  at (-0.5,0) {\tiny $\bullet$};
\node  at (0.5,0) {\tiny $\bullet$};
\node  at (5.5,0) {\tiny $\bullet$};
\node at (4.5,0) {\tiny $\bullet$};
\end{scope}
 \end{tikzpicture}
 \end{center}
\caption{Similarly to Figure \ref{pin110++}, the topological double of a  rectangle  $\square$, with the two (red) colored sides, is a cylinder ${\rm Cl}_{2,2}$.}
\label{pin110+}
\end{figure}
\vskip 2mm

  The cluster Poisson moduli space $  {\mathscr P}_{\G, \square} $ 
 is identified with the  Poisson Lie group $\G$ by Theorem \ref{Th4.3.3}.
 
 Similarly,  the cluster $K_2-$moduli space $  {\mathscr A}_{\G', \square} $ 
 is identified with $\G'$.

The pair $({\mathscr A}_{\G', \square}, {\mathscr P}_{\G, \square})=(\G', \G)$ is a cluster ensemble by Theorem \ref{MTHa}.  \vskip 2mm
  
 There is a finite cover 
 $$
 \pi: {\mathscr A}_{\G', \square} \lra {\mathscr P}_{\G, \square}.
 $$ 
 After the above identifications it reduces to the  canonical group projection $\pi: \G'\to \G$.
So there are positive structures on both groups $\G$ and $\G'$, and the projection $\pi$ 
  induces an isomorphism $\G'({\Bbb R}_{>0}) \to \G({\Bbb R}_{>0})$ on the space of the positive points.  This way we get canonical identifications:
 $$
 \G'({\Bbb R}_{>0}) = {\mathscr A}_{\G', \square}(\R_{>0})  =   {\mathscr P}_{\G, \square}({\Bbb R}_{>0}) =  \G({\Bbb R}_{>0}).
 $$

 Therefore we can define  the cluster Schwarz space 
 $  {\cal S}({\G}({\R}_{>0}))$ for the Hilbert space $ {\rm L}_2({\G}({\R}_{>0}))$ as follows: 
  $$
  {\cal S}({\G}({\R}_{>0})):= {\cal S}({\mathscr A}_{\G', \square}(\R_{>0}))  \subset {\rm L}_2({\mathscr A}_{\G', \square}(\R_{>0})) = {\rm L}_2({\G}({\R}_{>0})).
 $$ 
  Then  the $\ast-$algebra $\mathcal{A}_\hbar(\G, {\rm Cl}_{2,2}){\big \slash {\rm C}=1}$ has a cluster  realization in  $  {\cal S}({\G}({\R}_{>0}))$.
  \vskip 2mm
  
  The outer monodromies around the boundary components 
   provide a map  
   \be
 \mu_{\rm out}:  {\mathscr P}_{\G, {\rm Cl}_{2,2}}   \lra \H \times \H.
   \ee
   Denote by   ${\cal R}_{\G, {\rm Cl}_{2,2}}$ its fiber over the unit $e\in \H \times \H$.  
 By  part 3) of Theorem \ref{UEAB}, for each  boundary component $\pi$ of ${\rm Cl}_{2,2}$ there is a canonical embedding  
  \be
 \kappa_\pi:  {\U}_q(\mathfrak{g}) \hra {\cal O}_q({\mathscr R}_{\G,{\rm Cl}_{2,2}  } ):= {\cal O}_q({\mathscr P}_{\G,{\rm Cl}_{2,2}  } )_{{\rm out}=1}. 
  \ee
  Let $\pi$ and $\pi'$ be  the two boundary components of the cylinder. By the last part of Theorem \ref{5.25.24q}, the subalgebras $ \kappa_\pi(\mathcal{U}_q(\mathfrak{g}))$ and $\kappa_{\pi'}(\mathcal{U}_q(\mathfrak{g}))$ commute. Therefore we get a map of algebras: 
    \be \la{5.28.24.1}
 \kappa_\pi \otimes \kappa_{\pi'}:   {\U}_q(\mathfrak{g})  \otimes_{\Bbb C} {\U}_q(\mathfrak{g}) \lra {\cal O}_q({\mathscr R}_{\G,{\rm Cl}_{2,2}  } ).
  \ee  
It provides  a representation of the   $\ast-$algebra  $\mathcal{A}_\hbar(\mathfrak{g}) \otimes_\C \mathcal{A}_\hbar(\mathfrak{g}) $ 
  in the cluster Schwarz space of 
 ${\rm L}_2({\G}({\R}_{>0}))$.  As we explain below, this  is an analog of the  regular representation of the group $\G(\R) \times \G(\R)$ in $L_2(\G(\R))$. 
 \vskip 2mm

 \subsubsection{Quantum monodromy} The map   (\ref{5.28.24.1}) is no longer injective. To describe its kernel, recall the center $Z_{\mathfrak{g}}:= Z( {\U}_q(\mathfrak{g}))$  of the quantum group. 

   \bcon \la{5.28.24.3} The following two  subalgebras of the algebra $ {\cal O}_q({\mathscr P}_{\G,{\rm Cl}_{2,2}  } )$ coincide:
  \be
   \kappa_\pi({Z_{\mathfrak{g}}}) =   \kappa_{\pi'}({Z_{\mathfrak{g}}}). 
   \ee
   \econ
     Assuming Conjecture \ref{5.28.24.3}, denote by ${Z_{\mathfrak{g}} }$ the common subalgebra. This subalgebra provides an effective definition of the quantum monodromy 
     of a framed $\G-$local system over the middle loop $\alpha$ on the cylinder. Indeed, in this definition we can calculate the quantum monodromy explicitly knowing the center   
     of ${\U}_q(\mathfrak{g})$. In fact we can do it in two different ways, using the component $\pi$ or $\pi'$. 
      Conjecture \ref{5.28.24.3} tells that they lead to the same element - the quantum monodromy element - 
and  the map (\ref{5.28.24.1}) factors through the map 
\be
  {\U}_q(\mathfrak{g}) \otimes_{Z_{\mathfrak{g}} }  {\U}_q(\mathfrak{g}) \longrightarrow {\cal O}_q({\mathscr R}_{\G,{\rm Cl}_{2,2}  } ).
\ee

  \subsubsection{Example: the $\mathfrak{sl}_2$ case.}
 Figure \ref{6.2.24.23.53} shows the quantum double realization of the $\mathfrak{sl}_2$ case, obtained applying formulas (\ref{24.527.922}).

 For the inner circle, we have an embedding of ${\rm U}_q(\mathfrak{sl}_2)$ such that
 \[
 {\bf E}= {\rm T}_{-e_1+b_2}+ {\rm T}_{-e_1+b_1}; \qquad
 {\bf K}={\rm T}_{b_1}; \qquad
 {\bf F}={\rm T}_{e_1}+ {\rm T}_{e_1+e_2}+{\rm T}_{e_1+e_2-b_2+b_3}+{\rm T}_{e_1-b_1-b_2}
 \]
 For the outer circle, we have
 \[
 {\bf E}'= {\rm T}_{-e_3+b_2} + {\rm T}_{-e_3+b_3}, \qquad {\bf K}'= {\rm T}_{b_3}, \qquad {\bf F}'= {\rm T}_{e_3}+{\rm T}_{e_2+e_3}+{\rm T}_{e_2+e_3+b_1-b_2}+{\rm T}_{e_3-b_2-b_3}.
 \]
 The quantum monodromy
 \begin{align*}
 {\bf C}:&={\bf E}{\bf F}-q^{-1}{\bf K}-q{\bf K}^{-1}= {\bf E}'{\bf F}'-q^{-1}{\bf K}'-q{\bf K'}^{-1}\\
 &= {\rm T}_{b_2}+ {\rm T}_{-b_2}+{\rm T}_{e_2+b_2}+ {\rm T}_{e_2+b_1}+  {\rm T}_{e_2+b_3}+  {\rm T}_{e_2+b_1-b_2+b_3}.
 \end{align*}
 In particular, we verified Conjecture \ref{5.28.24.3} for $\mathfrak{sl}_2$.
 \begin{figure}[h]
 \begin{tikzpicture}[scale=1.3]
 \draw[dashed] (0,0)--(2,0);
 \draw (2,0)--(2,1);
 \draw[dashed] (2,1)--(0,1);
 \draw (0,1)--(0,0);
 \draw (0,1)--(2,0);
 \node[red] (a) at (0,0.5) {$\bullet$};
 \node[red] (b) at (1, 0.5) {$\bullet$};
 \node[red] (c) at (2, 0.5) {$\bullet$};
 \draw[blue, thick, -latex] (b)--(a);
 \draw[blue, thick, -latex] (b)--(c);
 
  \draw[dashed] (0,-2)--(2,-2);
  \draw (2,-2)--(2,-3);
  \draw[dashed] (2,-3)--(0,-3);
  \draw (0,-3)--(0,-2);
 \draw (0,-3)--(2,-2);
 \node[red] (d) at (0,-2.5) {$\bullet$};
 \node[red] (e) at (1, -2.5) {$\bullet$};
 \node[red] (f) at (2, -2.5) {$\bullet$};
  \draw[blue, thick, -latex] (d)--(e);
 \draw[blue, thick, -latex] (f)--(e);

 \node at (-0.2, 0.7) {$e_1$};
 \node at (1, 0.7) {$e_2$};
 \node at (2.2, 0.7) {$e_3$};
 
 \node at (-.7, -2.7) {$-e_1+b_2$};
 \node at (1, -2.7) {$-e_2-b_1-b_3$};
 \node at (2.7, -2.7) {$-e_3+b_2$};
 
 \draw[thick] (7, -1) ellipse (1cm and 0.5cm);
    \draw[thick] (7, -1) ellipse (2cm and 1.5cm);
    \draw (7.5, -1) ellipse (1.5cm and 1cm);
    \draw (5, -1) --(6, -1);
    \draw (8, -1) --(9, -1);
    \node[red] at (7,0.5) {$\bullet$};
    \node at (7, 0.7) {$e_3$};
      \node[red] at (7,-0.5) {$\bullet$};
      \node at (7, -0.3) {$e_1$};
        \node[red] at (7,-1.5) {$\bullet$};
        \node at (7, -1.7) {$-e_1+b_2$};
        \node[red] at (7.5, 0) {$\bullet$};
        \node at (7.5, 0.2) {$e_2$};
         \node[red] at (7.5, -2) {$\bullet$};
          \node at (7.5, -2.2) {$-e_2-b_1-b_3$};
    \node[red] at (5.5, -1) {$\bullet$};
    \node at (5.5, -.8) {$b_3-b_2$};
    \node[red] at (7, -2.5) {$\bullet$};
    \node at (7, -2.7) {$-e_3+b_2$};
    \node[red] at (8.5, -1) {$\bullet$};
    \node at (8.5,-1.2) {$b_1-b_2$};
 \end{tikzpicture}
 \caption{The quantum double realization for the $\mathfrak{sl}_2$ case. We glue the two rectangles  on the left along the  horisontal sides into an annulus on the right. Cutting the annulus along the  horisontal edges carrying variables $b_3-b_2$ and $b_1-b_2$,  we recover the  rectangles. }
 \label{6.2.24.23.53}
 \end{figure}
 
 \subsubsection{The quantum analog of the regular representation}

Cutting the cylinder ${\rm Cl}_{2,2}$ by a loop $\alpha$, see the Figure \ref{pin100+}, we get a map of  moduli spaces
$$
{\rm Res}_\alpha: {\cal R}_{\G, {\rm Cl}_{2,2}} \lra {\mathscr L}_{\G, \odot}  \times {\mathscr L}_{\G, \odot}.
$$   
Then one should have the map of quantum algebras induced by this map 
 \be \la{SDMR}
  {\rm Res}^*_{\alpha}: {\cal O}_q({\mathscr L}_{\G, \odot})  \otimes_\C {\cal O}_q({\mathscr L}_{\G, \odot})\stackrel{?}{\lra} {\cal O}_q({\cal R}_{\G, {\rm Cl}_{2,2}} ).
  \ee
Recall that, by  for any $\G$,  there are   canonical injective maps   
\be
\begin{split}
&\kappa_{\odot}: {\U}_q(\mathfrak{g})\lra {\cal O}_q({\mathscr L}_{\G, \odot}).\\
&\kappa_{{\rm C}_{2,2}}: {\U}_q(\mathfrak{g}) \otimes_\C {\U}_q(\mathfrak{g})\lra {\cal O}_q({\mathscr R}_{\G, {\rm C}_{2,2}})\\
\end{split}
\ee
 The map $  {\rm Res}^*_{\alpha}$ in (\ref{SDMR}) is supposed to make the following diagram commutative:
  \begin{displaymath}
    \xymatrix{
        {\cal O}_q({\cal R}_{\G, {\rm Cl}_{2,2}} ) & \ar[l]_{{\rm Res}_\alpha^* \qquad}  {\cal O}_q({\mathscr L}_{\G, \odot})  \otimes_\C {\cal O}_q({\mathscr L}_{\G, \odot}) \\
          \U_q(\mathfrak{g}) \otimes   \U_q(\mathfrak{g})  \ar[u]^{\kappa_{{\rm C}_{2,2}}}& \ar[l]_{=} \U_q(\mathfrak{g}) \otimes   \U_q(\mathfrak{g})  \ar[u]_{\kappa_{\odot \times \odot}} }
         \end{displaymath}
So, rather than  define the map  $  {\rm Res}^*_{\alpha}$, we  define its composition with $\kappa_{\odot \times \odot}$, forcing  the diagram commute:
$$
{\rm Res}^*_{\alpha}\circ \kappa_{\odot \times \odot}:= \kappa_{{\rm C}_{2,2}}.
$$
Note that when $\G$ is simply laced, we know thanks to \cite{S22} that the map $\kappa_{\odot}$, and hence $\kappa_{\odot \times \odot}$, are  isomorphisms, and 
presumably the same is true in general. So this forces the definition of the map $  {\rm Res}^*_{\alpha}$.      
    
          \begin{figure}[ht]
 \begin{center}
 \begin{tikzpicture}[scale=0.4]
 \draw (0,0) ellipse (0.5 and 1.25);
\draw (0, -1.25) -- (5, -1.25);
\draw (5, -1.25) arc (-90:90:0.5 and 1.25);
\draw [dashed] (5, -1.25) arc (270:90:0.5 and 1.25);
\draw [blue] (2.5, -1.25) arc (-90:90:0.5 and 1.25);
\draw [blue, dashed] (2.5, -1.25) arc (270:90:0.5 and 1.25);
\draw (5, 1.25) -- (0, 1.25);  
\node [red] at (-0.5,0) {\small $\bullet$};
\node [red] at (0.5,0) {\small $\bullet$};
\node [red] at (5.5,0) {\small $\bullet$};
\node [red] at (4.5,0) {\small $\bullet$};
\node [blue] at (2.5,0) {\small $\alpha$};
 \draw (13,0) circle (1.25);
  \draw (17,0) circle (1.25);
  \node [red] at (11.75,0) {\small $\bullet$};
\node [red] at (14.25,0) {\small $\bullet$};
\node [red] at (15.75,0) {\small $\bullet$};
\node [red] at (18.25,0) {\small $\bullet$};
\node [blue, thick] at (13,0) {$\circ$};
\node [blue, thick] at (17,0) {$\circ$};
\draw [-latex, thick] (8,0) -- (10, 0);
\end{tikzpicture}
\end{center}
\caption{Cutting a cylinder with $2+2$ special points by a loop $\alpha$.}
\label{pin100+}
\end{figure}
   Therefore    Modular Functor Conjecture \ref{MK} predicts an isomorphism of Hilbert spaces, compatible with the action of the corresponding $\ast-$algebras on the Schwarz spaces:
   \be \la{GRD}
 {\cal H}({\cal R}_{\G, {\rm Cl}_{2,2}}) \stackrel{\sim}{\lra} \int_{\lambda\in \H(\R_{>0})}  {\cal H}({\cal R}_{\G, \odot})_\lambda  \otimes {\cal H}({\cal R}_{\G, \odot})_{\lambda^{-1}}d\lambda.
 \ee   
Isomorphism  (\ref{GRD}) should intertwine  the action 
of the algebra $\mathcal{A}_\hbar(\mathfrak{g}) \otimes  \mathcal{A}_\hbar(\mathfrak{g}) $ on   ${\cal S}({\cal R}_{\G, {\rm Cl}_{2,2}})$ via the map 
$\kappa_{{\rm Cl}_{2,2}}$ with 
     its action 
     on  ${\cal S}({\cal R}_{\G, \odot})_\lambda  \otimes {\cal S}({\cal R}_{\G, \odot})_{\lambda^{-1}}$  via the map $\kappa_{\odot \times \odot}$.  
     The decomposition (\ref{GRD}) is the spectral decomposition for the quantum monodromy algebra, realized as the action of the common center $Z_{\mathfrak{g}}$.

Therefore decomposition (\ref{GRD}) is a quantum analog of the decomposition  
 \be
 {\rm L}_2(\G(\R)) = \int_{\chi}  {\cal V}_\chi  \otimes {\cal V}_{\chi^{*}}d\chi 
  \ee
of the regular representation of  the group $\G(\R) \times \G(\R)$  via the unitary representations ${\cal V}_\chi$  of  $\G(\R)$.   The decomposition of the latter   into an integral of irreducibles was the main problem of the classical representation theory, 
 solved by Gelfand and Naimark and by Harish-Chandra. Note that in the classical setting besides the unitary principal series there are other representations entering into 
 the decomposition, e.g. the discrete series for ${\rm SL}_2(\R)$.   The Modular Functor Conjecture \ref{MK} predicts that in the quantum setting the decomposition goes through the principal series 
 representations of the modular double quantum group.  We conclude that \vskip 2mm

{\it Modular Functor Conjecture \ref{MK} for the cylinder ${\rm Cl}_{2,2}$ is equivalent to the quantum analog of the   decomposition of  the regular representation for the Lie group $\G(\R)$.}

{\it The crucial point is that  the cylinder ${\rm C}_{2,2}$ can be presented as the topological double of the rectangle $\square$, and at the same time the loop $\alpha$ cuts the cylinder into 
two copies of the decorated surface $\odot$ providing the cluster realization of the quantum group $\U_q(\mathfrak{g})$.}

 
\medskip

 \section{Principal series of $\ast$-representations of quantized cluster varieties}  \la{Sec8a}
 
 \medskip

In Section \ref{Sec8.1}  we recall  the construction of the quantum symplectic double from \cite{FG07}. 
Then we  define several $\ast$-algebra structures on ${\cal O}_q(\mathscr{D})$ and its  modular double, and 
  construct     principal series of  representations of these $\ast$-algebras.  
 In Section \ref{SEC8.4}, starting from an arbitrary cluster Poisson variety $\mathscr{X}$,  we define   the principal series of representations of the   $\ast$-algebra ${\cal O}_q({\mathscr X})$ and its   modular double. The case $\beta\in \R$ was done in \cite{FG07}. 
 
 \subsection{The quantum double} \la{Sec8.1}
 
\medskip

\subsubsection{Doubling a  quiver.} Recall the notion 
$$
\mbox{quiver = (lattice, skew-symmetrizable  form, basis, multipliers)}= 
\left({\Lambda}, ( \ast, \ast), \{e_i\}, \{d_i\}\right)
$$ from Definition \ref{DQUIV}: 
We often use the coordinate description, setting 
\be
\varepsilon_{ij}:=(e_i, e_j), \qquad \widehat \varepsilon_{ij} := \varepsilon_{ij}d_j^{-1}; ~\qquad~\widehat \varepsilon_{ij}  = -\widehat \varepsilon_{ji} .
\ee
Let  $\Lambda^\circ \subset {\rm Hom}(\Lambda, \Q)$ be the quasi-dual lattice  generated by the quasi-dual basis $f_j:= d_j^{-1}e^*_j$, where $e^*_j(e_i)=\delta_{ij}$. 
The double $
\Lambda_\mathscr{D}:= \Lambda\oplus  \Lambda^\circ 
 $ of   $\Lambda$ 
carries  a skew-symmetrizable bilinear form $(\ast, \ast)_\mathscr{D}$: 
 \be
(e_i, e_j)_\mathscr{D} := (e_i, e_j), ~\qquad~(e_i, f_j)_\mathscr{D} := d_i^{-1}\delta_{ij},  ~\qquad~(f_i, f_j)_\mathscr{D} := 0.
\ee
 
\subsubsection{ The symplectic double.} 
We assign to a quiver ${\bf q}$ a split torus   with   coordinates  
$(B_i, Y_i)$, $i \in {V}$:
\be \la{splitD}
\mathscr{D}_{\bf q}:= {\rm Hom}(\Lambda_\mathscr{D}, {\Bbb G}_m).
\ee
It has a Poisson structure $\{, \}$ and a symplectic $2$-form $\Omega_{\bf q}$:
\begin{equation}\label{zx1}
\begin{split}
&\{\B_i, \B_j\}= 0, ~\qquad~ \{Y_i, \B_j\}= d_i^{-1}\delta_{ij}Y_i\B_j, ~\qquad~ \{Y_i, Y_j\} = 
\widehat \varepsilon_{ij}Y_iY_j.\\
&\Omega_{\bf q}:=  -\frac{1}{2}\sum_{i,j\in V} d_i \cdot 
\varepsilon_{ij} d\log \B_i \wedge d\log \B_j -
 \sum_{i \in V} d_i \cdot d\log \B_i \wedge d\log Y_i.\\
\end{split}
\end{equation} 
The Poisson structure  coincides with the one defined by the 
$2$-form $\Omega_{\bf q}$.  
Set 
$$
{\mathbb B}_k^+:= \prod_{i|\varepsilon_{ki}>0}\B_i^{\varepsilon_{ki}}, \qquad 
{\mathbb B}_k^-:= \prod_{i|\varepsilon_{ki}<0}\B_i^{-\varepsilon_{ki}}.
$$
Given an element $k \in {V}$, we define a birational transformation 
$\mu_k: \mathscr{D}_{\bf q} \to \mathscr{D}_{\bf q'}$ acting  
on the coordinates $\{Y'_i, \B'_i\}$ on the torus $\mathscr{D}_{\bf q}$ by 
\begin{equation}\label{zx1qt}
\begin{split}
&\mu_k^*: Y' _{i} \lms \left\{\begin{array}{lll} Y_k^{-1}& \mbox{ if $i=k$}  \\
 Y_i(1+Y_k^{-{\rm sgn}(\varepsilon_{ik})})^{-\varepsilon_{ik}}, & \mbox{ if $i \not = k$}.
\end{array} \right.\\
&\mu_k^*: \B'_{i} \lms 
\left\{\begin{array}{ll} 
\B_i & \mbox{ if $i\not =k$}, \\
 \frac{{\mathbb B}^-_k + Y_k {\mathbb B}_k^+}{\B_k(1+Y_k)} &
\mbox{ if $i = k$}.
\end{array} \right.\\
\end{split}
\end{equation}

The symplectic double $\mathscr{D}$ is  obtained by gluing the   tori $\mathscr{D}_{\bf q}$ using 
formulas (\ref{zx1qt}).  
Let us set
\begin{equation}\label{zx1q1}
 Y^\circ_i: =Y^{-1}_i
\prod_{j\in V}\B_j^{-\varepsilon_{ij}}. 
\end{equation}

Given a  variety ${\cal Y}$  with  a Poisson bracket/2-form,  denote by ${\cal Y}^{\rm op}$ 
the same variety with the opposite Poisson bracket/2-form.  
If ${\cal Y}$ is a cluster variety, we denote by ${\cal Y}^{\circ}$ its chiral dual, obtained by setting 
$\varepsilon^\circ_{ij} := -\varepsilon_{ij}$.  
The Poisson structure/2-form on ${\cal Y}^{\circ}$ is   the negative of the original one. Note that ${\cal Y}^\circ$ is a different cluster variety then ${\cal Y}$, 
while   ${\cal Y}^{\rm op}$ is the same variety with the negative Poisson bracket/2-form. 

 
Denote by $(A_i, A_i^\circ)$ the coordinates on 
the cluster torus ${\mathscr A}_{\bf q}\times {\mathscr A}^\circ_{\bf q}$.  
Let 
$p_-, p_+$ be the projections of ${\mathscr A}\times {\mathscr A}$ onto the two 
factors. 
The key
 properties of the symplectic double are the following: 
\begin{theorem} [{\cite[Theorem 2.3]{FG07}}]   \label{8.15.05.10} The 
   symplectic 
space     $\mathscr{D}$ has the following features: 
\begin{enumerate}

\vskip 1mm\item  There is a  Poisson map $\pi: \mathscr{D} \to {\mathscr X}\times {\mathscr X}^{\circ}$, given in any cluster 
coordinates by\footnote{Our definition of $\widetilde Y_i$ in (\ref{zx1q1}) is the inverse of the original one used in \cite[Formula (25)]{FG07}.} 
$$
\pi^*(X_i \otimes 1) = Y_i, ~\qquad~ \pi^*(1 \otimes X^\circ_j) = 
 Y^\circ_j.
$$

\vskip 1mm\item  There is a map $\varphi: {\mathscr A}\times {\mathscr A}^\circ \to \mathscr{D}$, such that $
\varphi^*\Omega_\mathscr{D} = 
p_-^*\Omega_{\mathscr A} -  p_+^*\Omega_{{\mathscr A}}$, given   in any cluster   by  
\be \la{11.28.18.11}
\varphi^*(Y_i) = \prod_{j\in {V}} A_j^{\varepsilon_{ij}}, ~\qquad~ 
\varphi^*(\B_i) = \frac{A_i^\circ}{A_i}.
\ee

\vskip 1mm\item  There are  commutative diagrams, where  $j$ is a Lagrangian  embedding, and $\Delta_{\mathscr X}$ the diagonal:  \begin{displaymath}
   \xymatrix{
        {\mathscr A}\times {\mathscr A}^\circ  \ar[dr]_{\varphi} &\stackrel{p \times p^\circ}{\lra}
         &\ar[dl]^{\pi} {\mathscr X}\times {\mathscr X}^{\circ}    \\
        &  \mathscr{D} &    }   
         \qquad \qquad ~\qquad~  \xymatrix{
        {\mathscr X} \ar[r]^{j} \ar[d]_{} &\mathscr{D}         \ar[d]^{\pi} \\
         \Delta_{\mathscr X}\ar[r]^{ } &  {\mathscr X}\times {\mathscr X}^{\circ} }
         \end{displaymath}
 The  intersection of the image of $j$ with each cluster torus is 
 given by equations $B_i =1$, $i \in V$. 

\vskip 1mm\item There is an  isomorphism $\zeta: \mathscr{D}\to \mathscr{D}^{\circ}$ 
interchanging the   
components   map $\pi$,  acting   any cluster 
  by 
$$
\zeta^*(\B_i) = \B^{-1}_i, ~\qquad~ \zeta^*(Y_i) =  {Y^\circ_i}^{-1}.
$$

\vskip 1mm\item  The map $\varphi$ is the quotient by the diagonal action of the 
torus $H_{\mathscr A}$. The map $\pi$ is the quotient by a free Hamiltonian action of the 
torus $H_{\mathscr A}$ on $\mathscr{D}$, whose 
 commuting Hamiltonians 
are provided by the Poisson composition map 
$
\mathscr{D} \stackrel{\pi}{\lra} {\mathscr X} \times {\mathscr X}^{\circ} {\lra} {\mathscr X} 
\stackrel{\theta }{\lra} H_{\mathscr X}.  
$  
\end{enumerate}

\end{theorem}

  It is instructive to give   a key argument for  the proof of 1), clarifying the second formula in (\ref{zx1qt}).  Mutating variables $A_k$ and $A^\circ_k$, we get new variables $A_k'$ and ${A_k^\circ}'$ such that:
   \be \nonumber
    \begin{split}
 &A_kA_k'= {\Bbb A}_k^+ + {\Bbb A}_k^-,\qquad ~\qquad~ A^\circ_k  {A^\circ_k}'=   {{\Bbb A}^\circ_k}^+ + \  {{\Bbb A}^\circ_k}^-.\\
 \end{split}
  \ee
 We have:
 \be \nonumber
  \begin{split}
 & {\Bbb B}_k^\pm:=      {{\Bbb A}^\circ_k}^\pm / {\Bbb A}_k^\pm, 
 \qquad Y_k:=  {\Bbb A}_k^+ /  {\Bbb A}_k^-,
 \qquad \widetilde Y_k:=   {{\Bbb A}^\circ_k}^- /    {{\Bbb A}^\circ_k}^+.\\
 \end{split}
  \ee 
  Then  the mutation map $\mu_k^*$ acts on $  A^\circ_k/A_k$ as follows
   \be \nonumber
    \begin{split}
 & \mu_k^*: \frac{{A^\circ_k}'}{A'_k} \lms \frac{ { {\Bbb A}^\circ_k}^+ +  { {\Bbb A}^\circ_k}^-}{\B_k({\Bbb A}_k^+ +  {\Bbb A}_k^-)}  =
  \frac{  {{\Bbb A}^\circ_k}^+/   {\Bbb A}_k^- +    {\Bbb B}_k^- }{\B_k(1+Y_k)} =  \frac{ {\Bbb B}_k^- + {Y}_k   {\Bbb B}_k^+}{\B_k(1+Y_k)}. \\
 \end{split}
  \ee

\subsubsection{ The quantum double.}    
The lattice 
$\Lambda_\mathscr{D}$ with the skew-symmetric form $\langle\ast, \ast\rangle_\mathscr{D}$  gives rise to a quantum torus algebra   
${\cal O}_q(\mathscr{D}_{\bf q}) $. Set $q_i:= q^{1/d_i}$.  A   basis $\{e_i, f_i\}$  of $\Lambda_\mathscr{D}$ 
provides the generators 
$\B_i, Y_i$ of the algebra ${\cal O}_q(\mathscr{D}_{\bf q})$, satisfying   
the relations 
\begin{equation} \label{4.28.03.11x}
q_i^{-1} Y_i\B_i  = 
q_i \B_iY_i, \qquad   \B_i Y_j = 
Y_j\B_i ~\qquad~ \mbox{if $i \not = j$}, \qquad   
q^{-\widehat \varepsilon_{ij}} Y_i Y_j = 
q^{-\widehat \varepsilon_{ji}} Y_jY_i.
\end{equation}
 
Denote by $ {\cal F}_q(\mathscr{D}_{\bf q}) $ the (non-commutative) fraction field of ${\cal O}_q(\mathscr{D}_{\bf q})$.

Recall the quantum dilogarithm power series
\be \la{PSIps}
{\bf \Psi}_q(x)= 
\prod_{k=1}^{\infty} (1+q^{2k-1}x)^{-1}.
\ee

Let $\mu_k: {\mathbf q} \to {\mathbf q}'$ 
be a mutation. We define 
a {\it quantum mutation} map  as an isomorphism 
\be \la{MMQD}  
\mu^\ast_{k, q}: 
{\cal F}_q(\mathscr{D}_{\bf ic})     \lra {\cal F}_q(\mathscr{D}_{\bf q}). 
\ee

The mutation map (\ref{MMQD})  
is the composition 
$
\mu^\ast_{k, q} := (\mu^{\sharp}_k)^*\circ (\mu'_k)^*,  
$ where: 

\begin{itemize} 

\item  The map 
 $
(\mu_k^{\sharp})^*: {\cal F}_q(\mathscr{D}_{\bf q})\lra {\cal F}_q(\mathscr{D}_{\bf q})
 $ 
is the conjugation by 
${\bf \Psi}_{q_k}(Y_k) / {\bf \Psi}_{q_k}(\widetilde Y^{-1}_k)$:
\be \la{MUSH}
(\mu_k^{\sharp})^*:= {\rm Ad}_{{\bf \Psi}_{q_k}(Y_k)} \circ {\rm Ad}^{-1}_{{\bf \Psi}_{q_k}(\widetilde Y^{-1}_k)}.
\ee
  Note that since    $Y_i$ and $\widetilde Y_j$ 
commute, 
${\bf \Psi}_q(Y_k)$ commutes with ${\bf \Psi}_q(\widetilde Y_k^{-1})$.

\vskip 1mm\item Denote by $(B_i', Y_i')$ the coordinates of the quiver torus 
$\mathscr{D}_{\bf i'}$. 
Then:
\begin{equation} \label{11.18.06.10hr}
\begin{split}
&(\mu'_{k})^*: 
{\cal F}_q(\mathscr{D}_{\bf i'})     \lra {\cal F}_q(\mathscr{D}_{\bf q}); \\
&\B'_{i} \lms \left\{\begin{array}{lll} \B_i& \mbox{ if } & i\not =k, \\
    {\Bbb B}_k^-/\B_k
 & \mbox{ if } &  i= k; \\ 
\end{array} \right. \qquad 
Y'_{i} \lms \left\{\begin{array}{lll} Y_k^{-1}& \mbox{ if } & i=k, \\
   q^{-[\varepsilon_{ik}]^2_+} Y_i(Y_k)^{[\varepsilon_{ik}]_+} & \mbox{ if } &  i\neq k. \\
\end{array} \right. \\
\end{split}
\end{equation} 
 \end{itemize}

Few  comments   are in order.

\vskip 1mm 1. The map $\mu_k^{\sharp}$ acts on the coordinates $\{\B_i, Y_i\}$ 
conjugating them by 
${\bf \Psi}_{q_k}(Y_k)$.

\vskip 1mm 2. Although  
$\mu_k^{\sharp}$ is   the conjugation by 
  quantum dilogarithms, it is  
  a birational map.

3. We define a quantum space $\mathscr{D}_q$ by using gluing 
isomorphisms (\ref{MMQD}). They satisfy the generalized pentagon relations. 
We can talk about them geometrically, 
saying that gluing isomorphisms correspond
to birational maps of non-commutative space 
$ \mathscr{D}_{{\bf q}, q} \to \mathscr{D}_{{\bf q'}, q}$, and  the 
quantum space $\mathscr{D}_q$ is obtained by gluing them via these maps. 
However all we use is the compositions of  gluing maps  (\ref{MMQD}). 
We denote by $\mathscr{D}^{\rm op}_q$ the quantum space obtained by taking the opposite product $a\circ b = ba$ on   quantum torus algebras.

The main properties of the quantum double are summarized in the following theorem: 
\begin{theorem}[{\cite[Theorem 3.3]{FG07}}]
\label{8.15.05.10gg} 
There is a 
quantum 
space $\mathscr{D}_q$ together with:

\begin{enumerate}

\vskip 1mm\item  An  involutive 
isomorphism 
$\ast : \mathscr{D}_{q}\to \mathscr{D}_{q^{-1}}^{\rm op}$,  
given in any cluster coordinate system by 
$$
\ast(q) = q^{-1}, ~\qquad~ \ast (Y_i) = Y_i, ~\qquad~ \ast (\B_i) = \B_i.
$$

\vskip 1mm\item  A map  $\pi:\mathscr{D}_q  \lra {\mathscr X}_q\times {\mathscr X}^{\circ}_q $ 
given in a cluster coordinate system by 
\begin{equation} \nonumber
\pi^*(X_i \otimes 1) = Y_i, ~\qquad~ \pi^*(1 \otimes X^\circ_i) = 
 Y^\circ_i:= Y^{-1}_i\prod_{j}\B_j^{-\varepsilon_{ij}}.
\end{equation}

\vskip 1mm\item An involutive map $\zeta: \mathscr{D}_q \stackrel{\sim}{\lra} \mathscr{D}_{q}^{\rm op}$ 
interchanging the two components of  $\pi$,   
given in any cluster 
  by 
$$
\zeta^*\B_i = \B^{-1}_i, ~\qquad~ \zeta^*Y_i = {Y_i^{\circ}}^{-1}. 
$$

\vskip 1mm\item  A canonical map of 
quantum spaces $\theta_q: {\mathscr X}_q \to H_{\mathscr X}$. 
\end{enumerate}
\end{theorem}

\subsubsection{Connections between quantum  varieties.} 
 There are three ways to alter the space ${\mathscr X}_{q}$: 

(i) change $q$ to $q^{-1}$, 

(ii) change the quantum space ${\mathscr X}_{ q}$ to its chiral dual ${\mathscr X}^\circ_{ q}$ - that is the space with $\varepsilon^\circ_{ij}:=-\varepsilon_{ij}$, 

(iii) change the quantum space ${\mathscr X}_{ q}$ to the opposite quantum space 
${\mathscr X}^{\rm op}_{ q}$.

\noindent

\noindent
The resulting three quantum  spaces are canonically isomorphic (\cite{FG07}, \cite[Lemma 2.1]{FG03b}):
\begin{lemma} \label{6.9.03.11} 
There are canonical isomorphisms of quantum spaces
\be
\begin{split}
&\alpha_{\mathscr X}^q: {\mathscr X}_{q} \lra {\mathscr X}^{\rm op}_{ q^{-1}}, \qquad  (\alpha_{\mathscr X}^q)^*
: X_i \lms X_i. \\
&i^q_{\mathscr X}: {\mathscr X}_{ q} \lra {\mathscr X}^\circ_{ {q^{-1}}}, \qquad  (i^q_{\mathscr X})^*: X_i^\circ\lms X_i^{-1}.\\
&\beta_{\mathscr X}^q:= \alpha_{\mathscr X}^q\circ i^q_{\mathscr X}: 
{\mathscr X}^\circ_{ q} \lra {\mathscr X}^{\rm op}_{ q}, \qquad  X_i \lms {X_i^\circ}^{-1}. \\
\end{split}
\ee
\end{lemma}

 Similarly there are   three ways to alter the space $\mathscr{D}_{q}$, and three  
isomorphisms acting the same way on the $X$-coordinates, 
and identically on the $B$-coordinates. They are compatible with the projection 
$\pi$, while the involution $i$ in Theorem \ref{8.15.05.10gg} 
interchanges the two components of  $\pi$.

 \subsection{Canonical $\ast$-representation of the quantum double:   skew-symmetric case} \la{Sec8.2a}
 
\medskip

 In Section \ref{Sec8.2a} we consider  arbitrary {\it simply-laced} quivers ${\bf q}$, i.e. quivers with   $\varepsilon_{ij} = - \varepsilon_{ji}$.

\subsubsection{The modular quantum dilogarithm.} 
Let ${\rm sh}(t):= (e^t-e^{-t})/2$. Recall the modular  quantum dilogarithm function:
\be \la{12.12.18.8}
\Phi_\hbar(z) := {\rm exp}\Bigl(-\frac{1}{4}\int_{\Omega}\frac{e^{-i pz}}{ {\rm sh} (\pi   p)
{\rm sh} (\pi   \hbar p) } \frac{dp}{p} \Bigr).
\ee
Here the contour $\Omega$ goes along the real axes 
from $- \infty$ to $\infty$ bypassing the origin 
from above. 

Let us pick a complex number $\beta$ such that ${\rm Re}(\beta) > 0$.\footnote{The traditional notation for $\beta$ is $b$. We change it to $\beta$ to avoid confusion with the variables $b_k$ on the double.} Set 
\be \la{12.12.18.9}
\varphi_\beta(z) := {\rm exp}\Bigl(-\frac{1}{4}\int_{\Omega}\frac{e^{-i pz}}{ {\rm sh} (\pi \beta^{-1} p)
{\rm sh} (\pi  \beta p) } \frac{dp}{p} \Bigr).
\ee
If $z \in \R$, the integral converges since ${\rm Re}(\beta) \not = 0$ and the contour $\Omega$ is going along the real line. Indeed,  in this case 
  the denominator of the integrand have no zeros, and grows exponentially. 

If $\beta$ is a positive number   such that $\beta^2 = \hbar$, then  making the substitution $p=q\beta$ we get 
\be  \la{12.12.18.3}
\Phi_\hbar(z ) = \varphi_\beta(z/\beta).
\ee

Recall the  relation with the power series $\Psi_q(z)$, see (\ref{PSIps}), where $q = e^{i\pi \beta^2}$ and $q^\vee = e^{i\pi \beta^{-2}}$:
\be \la{KYREL}
\varphi_{\beta}(z)  = \frac{\Psi_q(e^{\beta z})}{\Psi_{1/q^\vee}(e^{z/\beta })}.
\ee

The following result on  the unitarity of the function $\varphi_\beta(z)$  is crucial  in our story. 

\bl \la{12.15.18.1} Assume that  $z\in {\mathbb R}$, and $(\beta+\beta^{-1})^2 >0$,  i.e.  $\beta$ satisfies one of the following   conditions:
\be \la{12.21.18.1}
{\rm R)}~ \beta \in \R,   \qquad {\rm U)}~ |\beta|=1, ~\qquad~\beta \not \in \pm i.  
\ee
Then one has   $|\varphi_\beta(z)|=1$.
\el 

\begin{proof} Condition (\ref{12.21.18.1}) on $\beta$ just means that  the complex conjugation acts as follows:
$$
{\rm R)}~    (\beta, \beta^{-1}) \lms (\beta, \beta^{-1}),   \qquad {\rm U)}~ (\beta, \beta^{-1}) \lms  (\beta^{-1}, \beta).
$$
Each of these maps preserves the function ${\rm sh} (\pi   \beta^{-1} p){\rm sh} (\pi   \beta p)$. 
 Using this  with  $\overline z =z$, and  changing   the integration
variable $q=-\overline{p}$, we get:
\be \nonumber 
\begin{split}
\overline{\varphi_\beta(z)}=~&{\rm exp}\Bigl(-\frac{1}{4}\int_{\overline \Omega}\frac{e^{i p\overline z}}{ {\rm sh} (\pi \overline \beta^{-1} p)
{\rm sh} (\pi  \overline \beta p) } \frac{dp}{p} \Bigr) \stackrel{q=-\overline p}{=}   \\
&{\rm exp}\Bigl(-\frac{1}{4}\int_{-\Omega}\frac{e^{-i p  z}}{ {\rm sh} (\pi   \beta^{-1} p)
{\rm sh} (\pi    \beta p) } \frac{dp}{p} \Bigr) =   \\
&{\rm exp}\Bigl(\frac{1}{4}\int_{\Omega}\frac{e^{-i p  z}}{ {\rm sh} (\pi   \beta^{-1} p)
{\rm sh} (\pi    \beta p) } \frac{dp}{p} \Bigr) =    {\varphi_\beta(z)}^{-1}.\\
\end{split}
\ee 
  \end{proof}
  
  The function $\varphi_\beta(z)$ satisfies the difference equations
  \be \nonumber
\begin{split}
&\varphi_\beta(z+2\pi i \beta) = (1+qe^{\beta z}) \varphi_\beta(z),\\
&\varphi_\beta(z+2\pi i \beta^{-1}) = (1+q^\vee e^{  z/\beta}) \varphi_\beta(z).\\
\end{split}
\ee

\subsubsection{Hilbert spaces.}   Denote by    $\{A_i\}$   cluster coordinates on a cluster $K_2$-variety ${\mathscr A}$ assigned to      a quiver ${\bf q}$.  
 The set of   positive points ${\mathscr A}(\R_{>0})$ has   
a vector  space structure with  logarithmic   coordinates 
    $a_i:= \log A_i$.  So  it has a   defined up to a sign  volume form $\mu_{\mathscr A}= da_1 \ldots da_n$. 
Consider  the related Hilbert space:     
  \be \nonumber
{\cal H}_{{\mathscr A}} := {\rm L}_2\left({\mathscr A}(\R_{>0}), \mu_{\mathscr A}\right).
\ee
   The volume form changes  the sign under  mutations. So  the 
    Hilbert spaces 
 assigned to   quivers related by cluster transformations are canonically identified.

\subsubsection{Representations.} 
We use the following notation: 
  \begin{equation} \label{handy}
\begin{split}
&\alpha_k^+:= \sum_{j\in {\rm V}} 
[\varepsilon_{kj}]_+a_j, \qquad  
\alpha_k^-:= \sum_{j\in {\rm V} }
[-\varepsilon_{kj}]_+a_j.\\
&\hbar:= \beta^2, \qquad  q:= e^{\pi i \beta^2}, 
\qquad  \hbar^{\vee}:= \beta^{-2},   \qquad  q^{\vee}:= e^{\pi i \beta^{-2}}.\\
\end{split}
\end{equation}
Consider the following  differential 
operators on the manifold $ {\mathscr A}(\R_{>0})$: 
\begin{equation} \label{5.26.08.1aa}
\begin{split}
& p_k:=   2\pi i  \cdot  \frac{\partial}{\partial a_k}
-\alpha_k^+,\\ 
& {{ p}^\circ_k} :=   
 -  2\pi i  \cdot \frac{\partial}{\partial a_k}
+ \alpha_k^-,\\
&a_l:= a_l.\\
\end{split}
\end{equation}
Since $ \alpha_k^+ - \alpha_k^- = 
\sum_{j }\varepsilon_{kj}a_j$,  we have
\begin{equation} \nonumber
\begin{split}
& {{ p}^\circ_k} := - p_k 
- \sum_{j\in {\rm V} }\varepsilon_{kj}   a_j, \qquad [ {{ p}^\circ_k},  p_l]=0.\\ 
\end{split}
\end{equation}
 Operators (\ref{5.26.08.1aa})    satisfy   relations of the Heisenberg $\ast$-algebra:
\be \nonumber
\begin{split}
&[ p_k,  p_l] = 2 \pi i   \cdot   \varepsilon_{kl}, 
~\qquad~ \qquad ~[ p_k,  a_l] = 2 \pi i   \cdot \delta_{kl}, \qquad ~\qquad~
\ast  p_k =  p_k, \qquad \ast a_l = a_l.\\
&[ {{ p}^\circ_k},  {{{ p}^\circ_l} }] 
= -2 \pi i   \cdot   \varepsilon_{kl}, \qquad  
[ {{ p}^\circ_k},  a_l]  = -2 \pi i    \cdot \delta_{kl}, 
\qquad 
\ast  {{ p}^\circ_k} =  {{ p}^\circ_k}. \\
\end{split}
\ee

 Rescaling  operators (\ref{5.26.08.1aa}) by $\beta $, we get the following operators: 
 \begin{equation} \label{handyZ}
\begin{split}
&b_k := \beta  \cdot a_k,\\
& y_k:= \beta  \cdot p_k, \\ 
& {{ y}^\circ_k} := \beta \cdot p_k^\circ.
\\
\end{split}
\end{equation}

Operators $b_k, y_k$ in (\ref{handyZ})  satisfy   relations of the Heisenberg $\ast$-algebra 
$ {\cal H}_{{\hbar},{\bf q}}$:
\be \label{5.26.08.2}\begin{split}
&[ y_k,  y_l] = 2 \pi i  \beta^2  \cdot \varepsilon_{kl}, 
\qquad  ~\qquad~[ y_k,  b_l] = 2 \pi i \beta^2\cdot \delta_{kl}, \qquad 
\ast  y_k =  y_k, ~\qquad~ \ast  b_l =  b_l.\\
\end{split}
\ee

  The exponentials $ Y_p:= {\rm exp}( y_p)$ and  $ B_p:= {\rm exp}( b_p)$ are 
 difference operators. 
The operators $\{ {Y}_p,  {B}_p\}$ satisfy the commutation relations  
of the quantum torus algebra ${\cal O}_q(\mathscr{D}_{\bf q})$.

 The operators $b_k, y^\circ_k$ in (\ref{handyZ})  satisfy the  relations of the    opposite  Heisenberg $\ast$-algebra   ${\cal H}^{\rm op}_{{\hbar},{\bf q}}$:
\be \nonumber
\begin{split}
&[ {{ y}^\circ_p},  {{{ y}^\circ_q} }] 
= -2 \pi i \beta^2 \cdot   \varepsilon_{kl}, \qquad ~ 
[ {y^\circ_k},  b_l]  = -2 \pi i \beta^2\cdot \delta_{kl}, 
\qquad \ast y^\circ_k =  {{ y}^\circ_k}. \\
\end{split}
\ee
There is a canonical isomorphism ${\cal H}^{\rm op}_{{\hbar},{\bf q}} \to {\cal H}_{\hbar, \bf c^o}$ 
acting as the identity on the generators.

\subsubsection{The modular dual representation.} Let us introduce the  modular dual operators
\begin{equation} \nonumber
\begin{split}
& {y}^{\vee}_k :=   {y}_k/\beta^2, ~\qquad~ 
 {b}^{\vee}_k :=   {b}_k/\beta^2, \\
\end{split}
\end{equation} 
Equivalently,  we can write 
 \begin{equation} \label{handyZDW}
\begin{split}
&b^{\vee}_k := \beta^{-1}  \cdot a_k,\\
& y^{\vee}_k:= \beta^{-1}  \cdot p_k, \\ 
& {{ y}^\circ_k}^{\vee} := \beta^{-1} \cdot p_k^\circ.
\\
\end{split}
\end{equation}

The operators $\{ {y}^{\vee}_k,  {b}^{\vee}_k\}$  
satisfy   commutation relations of the Heisenberg $\ast$-algebra 
${\cal H}_{\hbar^{\vee},{\bf q}}$: 
\begin{equation} \nonumber
\begin{split}
&[{y}^{\vee}_k, {y}^{\vee}_l] = 2 \pi i \beta^{-2}   \varepsilon_{kl},\\
&[{y}^{\vee}_k, {b}^{\vee}_l] = 2 \pi i \beta^{-2} \delta_{kl}, \qquad [{b}^{\vee}_k, {b}^{\vee}_l] =0.\\
\end{split}
\end{equation}
Their exponents are the difference operators  satisfying   the relations of the algebra  
${\cal O}_{q^{\vee}}(\mathscr{D}_{\bf q})$:
 \begin{equation} \nonumber
\begin{split}
& Y^{\vee}_k:= {\rm exp}( y^{\vee}_k), ~\qquad~ 
 B^{\vee}_k:= {\rm exp}( b^{\vee}_k)\\
\end{split}
\ee

The commutators of the operators $\{ {y}_k,  {b}_k\}$ with the ones $\{ {y}^{\vee}_l,  {b}^{\vee}_l\}$ lie in $2\pi i \Z$:
\begin{equation} \nonumber
\begin{split}
&[{y}_k, {y}^{\vee}_l] =   2 \pi i  \cdot \varepsilon_{kl},\\
&[{y}_k, {b}^{\vee}_l] = 2 \pi i  \cdot  \delta_{kl}, \qquad [{b}_k, {b}^{\vee}_l] =0.\\
\end{split}
\end{equation}
So operators $\{ {Y}^{\vee}_k,  {B}^{\vee}_k\}$ commute with the ones $\{ {Y}_k,  {B}_k\}$. 
 
\subsubsection{A subspace $W_{\bf q}$.} Let $W_{\bf q} $ be the subspace spanned by the   functions 
\begin{equation} \label{9.4.07.1}
P(a_1, \ldots, a_n){\rm exp}{(- \sum_{i } (a_i^2/2+ \lambda_ia_i))}, ~\qquad~\mbox{
where  $\lambda_i\in \C$, and $P$ is a polynomial.} 
\end{equation}
It is dense in ${\cal H}_{{\mathscr A}}$, and invariant under the Fourier transform. 
The operators $ \{B_p, Y_p, B^{\vee}_p,   Y^{\vee}_p\}$ are symmetric unbounded operators 
   on  $W_{\bf q}$, 
 satisfying  the standard  commutation relations.  Set
\be \nonumber
\mathcal{A}_\hbar(\mathscr{D}_{\bf q}) = {\cal O}_{q}(\mathscr{D}_{\bf q}) \otimes {\cal O}_{q^{\vee}}(\mathscr{D}_{\bf q}).
\ee 
So we get a representation of the algebra $\mathcal{A}_\hbar(\mathscr{D}_{\bf q})$ in the subspace $W_{\bf q}$ of the Hilbert space ${\cal H}_{{\mathscr A}}$.

 \vskip 2mm
 
By Lemma \ref{12.15.18.1}, the  modular quantum  dilogarithm function $\varphi_\beta(z)$ is unitary when $z \in \R$ and $\beta$ is real  or unitary. 
 Below we show that then the algebra $\mathcal{A}_\hbar(\mathscr{D}_{\bf q})$ has a $\ast$-algebra structure. 
  When $\beta$ is real, it is induced by a $\ast$-algebra structure on ${\cal O}_q(\mathscr{D}_{\bf q})$.
 
\subsubsection{The involution $\ast_{\rm R}$.} Let $\beta \in \R$.   
Then  there are   antiholomorphic involutive antiautomorphisms $\ast_\R$ of the algebras ${\cal O}_q(\mathscr{D}_{\bf q})$ and $\mathcal{A}_\hbar(\mathscr{D}_{\bf q})$, for which 
 operators (\ref{handyZ}) $\&$ (\ref{handyZDW}) are real: 
\be \nonumber
\begin{split}
&\ast_{\rm R}(y_k) = y_k, \qquad \  ~\ast_{\rm R}(y^\circ_k) = {y^\circ_k}, \qquad \quad \ast_{\rm R}(b_k) = b_k.\\
&\ast_{\rm R}(y^\vee_k) = y^\vee_k, \qquad \ast_{\rm R}({y^\circ_k}^\vee) = {y^\circ_k}^\vee, \qquad  \ast_{\rm R}(b^\vee_k) = b^\vee_k.\\
\end{split}
\ee

\subsubsection{The involution $\ast_{\lr}$.} Let   $\overline \beta  = \beta^{-1}$. So $\beta \in \U(1)= \{z\in \C^*~|~|z|=1\}$.  Then  there is a similar involution $\ast_{\lr}$ of the algebra  $\mathcal{A}_\hbar(\mathscr{D}_{\bf q})$ which   interchanges the operators
  (\ref{handyZ}) and  (\ref{handyZDW}):
$$
\ast_{\lr}(y_k) = y_k^\vee, \qquad \ast_{\lr}(y^\circ_k) = {y^\circ_k}^\vee, \qquad  \ast_{\lr}(b_k) = b_k^\vee.
$$

Summarising, we arrive at the following result. 

\bl Let ${\bf q}$ be any  skew-symmetric quiver. Then 

\vskip 1mm$\rm R$) If $ \beta  \in \R$, 
then  there is a $\ast$-algebra structure  on the algebra ${\cal O}_q(\mathscr{D}_{\bf q})$, and hence on 
 $\mathcal{A}_\hbar(\mathscr{D}_{\bf q})$.

\vskip 1mm$\U$) If  $|\beta|=1$,  then  there is a $\ast$-algebra  structure on the algebra $\mathcal{A}_\hbar(\mathscr{D}_{\bf q})$. 
 
 These  $\ast$-algebras     
 have   representations in the subspace $W_{\bf q}\subset {\cal H}_{\mathscr A}$, provided  by  the operators 
$\{ {Y}_k,  {B}_k, {Y}^{\vee}_k,  {B}^{\vee}_k\}$ for $\mathcal{A}_\hbar(\mathscr{D}_{\bf q})$,  and $\{{Y}_k,  {B}_k\}$ for ${\cal O}_q(\mathscr{D}_{\bf q})$.
\el

\vskip 2mm
 
 Recall    the  algebra of universally Laurent polynomials ${\cal O}_q(\mathscr{D})$.  Set 
\be \nonumber
\mathcal{A}_\hbar(\mathscr{D}):= {\cal O}_q(\mathscr{D})\otimes {\cal O}_{ q^\vee}(\mathscr{D}).
\ee
 
 Recall the mutation $
\mu^\ast_{k} := (\mu^{\sharp}_k)^*\circ (\mu'_k)^*$, where   $(\mu_k^{\sharp})^*:= {\rm Ad}_{{\bf \Psi}_{q}(Y_k)} \circ {\rm Ad}^{-1}_{{\bf \Psi}_{q}({Y^\circ_k}^{{-1}})}$, see (\ref{MUSH}). 
 
 Proposition \ref{MUMAP} is crucial for our story. 
 
 \bp \la{MUMAP} The   mutation $\mu_k^*: \mathcal{A}_\hbar(\mathscr{D}_{{\bf q'}}) \lra \mathcal{A}_\hbar(\mathscr{D}_{\bf q})$ commutes with the involution $\ast_{\rm L}$. 
So  there is a $\ast$-algebra structure on the algebra $\mathcal{A}_\hbar(\mathscr{D})$,  given in any cluster coordinates by the involution 
 $\ast_{\rm L}$. \ep
 
 \begin{proof} This is clear for the mutation map $(\mu'_k)^*$. So the claim follows from 
 the next Lemma. 
 
 \bl The   automorphism $(\mu_k^{\sharp})^*$ of the  algebra  $\mathcal{A}_\hbar(\mathscr{D}_{\bf q})$ commutes with the involution $\ast_{\rm L}$.
 \el
 
 \begin{proof}  Indeed, we have  
 $$
 \ast_{\lr}(\Psi_q(Y_k)) = \Psi_{1/q^\vee}(Y_k^\vee).
 $$
So, since $\ast_{\lr}$ is an antiautomorphism, and using the  power series identity $\Psi_{1/q}(X)^{-1} = \Psi_{q}(X)$, 
$$
\ast_{\lr}\circ {\rm Ad}_{\Psi_q(Y_k)} =    {\rm Ad}_{\Psi_{1/q^\vee}(Y_k^\vee)^{-1}}\circ \ast_{\lr} = {\rm Ad}_{\Psi_{q^\vee}(Y_k^\vee)}\circ \ast_{\lr} .
$$
Since $\ast_{\lr}(\widetilde Y_k) =\widetilde Y_k^\vee$, there is a similar identity for ${\rm Ad}^{-1}_{{\bf \Psi}_{q}(\widetilde Y^{-1}_k)}$.
\end{proof}
\end{proof}

\subsubsection{Remark.}    Let $\beta \in i\R$.   
Then  there is an   antiholomorphic involutive antiautomorphism $\ast_{\rm I}$ of the   algebra ${\cal O}_q(\mathscr{D}_{\bf q})$, and hence of $\mathcal{A}_\hbar(\mathscr{D}_{\bf q})$, so that    operators (\ref{handyZ}) $\&$ (\ref{handyZDW}) are purely imaginary: 
\be \nonumber
\begin{split}
&\ast_{{\rm I}}(y_k) = -y_k, \qquad ~\ast_{\rm I}(y^\circ_k) = -{y^\circ_k}, \qquad \qquad \ast_{\rm I}(b_k) = -b_k.\\
&\ast_{\rm I}(y^\vee_k) = -y^\vee_k, \qquad \ast_{\rm I}({y^\circ_k}^\vee) = -{y^\circ_k}^\vee, \qquad  \ast_{\rm I}(b^\vee_k) = -b^\vee_k.\\
\end{split}
\ee 
There is a similar involution $\ast_{{\rm I}}$ of any quantum torus algebra. 
However the involution $\ast_{{\rm I}}$ does not commute with mutations: the analog of Proposition \ref{MUMAP} fails. So it does not define an involution of the quantum cluster variety. 
And the function $\varphi_{\beta}(z)$ is problematic when $\beta \in i\R$.  

\subsubsection{Intertwiners.} The  operators $p_k$ and 
$ {p^\circ_k}$ defined in (\ref{5.26.08.1aa}) are self-adjoint. Indeed, $a_k$ are real variables, 
so it is a standard fact about the operator $i\frac{\partial}{\partial a}$ and the operator of multiplication by $a$ acting on functions $f(a)$ in ${\rm L}_2(\R)$. 
Therefore one can apply to them any continuous function on the real line. 
If the function  is unitary, we get a unitary operator. 
We use below the   modular quantum dilogarithm function $\varphi_{\beta}(x)$ defined in (\ref{12.12.18.9}).

\begin{definition} \label{intertwiner-beta}
Let $\hbar\in \C^*$ such that $\hbar + \hbar ^{-1} >-2$. 

Given a mutation $\mu_k: {\bf q} \to {\bf q}'$, we consider 
the  intertwining operator 
\begin{equation} \label{12.01.06.1}
{\cal I}_{{\bf q}'\to {\bf q}}:= {\bf K}^{\sharp}\circ 
{\bf K}^{'} : L^2({\mathscr A}(\R_{>0})) \lra 
L^2({\mathscr A}(\R_{>0})), \qquad \mbox{where}
\end{equation}
\begin{itemize}
\vskip 1mm\item
{\it The operator ${\bf K}^{\sharp}$} is  the 
 ratio  of quantum dilogarithms of the operators $p_k$ and $ -{p}^\circ_k$: 
\begin{equation} \label{*RUu}
\begin{split}
&{\bf K}^{\sharp}:=  \varphi_{\beta}( p_k)
\varphi_{\beta} (  -{p}^\circ_k)^{-1}. \\
\end{split}
\end{equation}

\vskip 1mm\item
{\it The operator ${\bf K}^{'}$} is induced by  the  linear map 
acting on the coordinates as follows:\footnote{This is just the coordinate transformation corresponding to the mutation  of the quasi-dual basis $f_i$, see (8) of \cite{FG03b}.}
\begin{equation} \label{11.18.06.10d}
{a}_i'\lms \left\{\begin{array}{lll} 
a_i& \mbox{ if } & i\not =k, \\
    \alpha_k^- - a_k
 & \mbox{ if } &  i= k. \\
\end{array} \right.
\end{equation}
\end{itemize}
\end{definition} 

Consider the Fourier transform ${\cal F}_{a_k}$
along the $a_k$-coordinate: 
$$
{\cal F}_{a_k}(f)(c) =   \frac{1}{2\pi i} \cdot \int e^{a_kc/2\pi i }f(a_k)da_k;~\qquad~
f(a_k)=\frac{1}{2\pi i }\int e^{-a_kc/2\pi i }\widehat{f}(c)dc.
$$
Then:
\be
\nonumber
{\bf K}^{\sharp} = {\cal F}_{a_k}^{-1} \circ \varphi_{\beta }\left(-  c -\alpha_k^+\right)\varphi_{\beta }\left( -  c -\alpha_k^-\right)^{-1} \circ {\cal F}_{a_k}. 
\end{equation}

The inverse ${\cal I}^{-1}_{{\bf q'}\to {\bf q}}$   has a  simpler  integral presentation:
\begin{equation}\nonumber 
\begin{split}
&({\cal I}^{-1}_{{\bf q}'\to {\bf q}}f)(a_1,\ldots, a'_k, \ldots, a_n) := \int G(a_1,\ldots, 
a'_k+a_k, \ldots, a_n)f(a_1,\ldots, a_k,\ldots, a_n)da_k,\\
&G(a_1,\ldots,a_n):= \frac{1}{(2\pi i )^2 }\int \varphi_{\beta }(- c
-\alpha_k^+)^{-1}\varphi_{\beta_k}(- c
-\alpha_k^-)^{-1}\exp\left(c
\frac{a_k-\alpha_k^-}{2\pi i  }\right)dc.\\
\end{split}
\end{equation}

\subsubsection{The Schwartz space.} \la{S17.2.10}

\begin{definition} \label{7.10.07.1} The Schwartz subspace ${\cal S}_{{\bf q}}  \subset 
L^2({\mathscr A}(\R_{>0}))$ consists 
of all vectors $f$ such that for any $ A \in \mathcal{A}_\hbar(\mathscr{D})$, the functional $w \to (f,   Aw)_{L^2 }$ on $W_{\bf q}$
is continuous.
\end{definition}
The subspace ${\cal S}_{\bf q}$
is
 the common domain of definition of
operators from $\mathcal{A}_\hbar(\mathscr{D})$ in ${\cal H}_{{\mathscr A}}$. Indeed,
since $W_{\bf q}$ is dense in  ${\cal H}_{{\mathscr A}}$,
by the Riesz theorem  for any $f \in {\cal S}_{\bf q}$ there is a
unique $g \in {\cal H}_{{\mathscr A}}$ 
such that $(g,w) = (f,   Aw)$. We set $ {A}^\ast f := g$. 
The  space ${\cal S}_{{\bf q}}$ has a   topology given by  seminorms
$$
\rho_B(f):= ||Bf||_{L^2}, \qquad  \mbox{$B$ runs through a basis in $\mathcal{A}_\hbar(\mathscr{D})$}.
$$
So for a given quiver ${\bf q}$ we get a representation of the $\ast$-algebra $\mathcal{A}_\hbar(\mathscr{D})$ in the   subspace ${\cal S}_{{\bf q}} \subset {\cal H}_{\mathscr A}$.

\subsubsection{General intertwiners.} A 
  cluster transformation  
${\bf c}: {\bf q} \to {\bf q'}$ provides a unitary operator 
$$
{\cal I}_{\bf c^\circ}: {\rm L}^2({\mathscr A}(\R_{>0})) \lra {\rm L}^2({\mathscr A}(\R_{>0})).  
$$
Indeed,  we assigned in (\ref{12.01.06.1}) to   a   
 mutation ${\bf q} \to {\bf q'}$ 
an intertwiner 
${\cal I}_{{\bf q'} \to {\bf q}}$. An  automorphism $\sigma$ 
of a quiver ${\bf q}$  gives rise 
to 
a unitary operator given by a permutation of  
coordinates   $\{a_i\}$. The  
 cluster transformation  
 ${\bf q}$  is a composition of mutations and automorphisms. 
Taking the reverse 
composition of the corresponding intertwiners, we get the unitary operator 
${\cal I}_{\bf c^\circ}$.

A  cluster transformation ${\bf c}$ gives rise to 
a   cluster transformation    
${\bf c}^\mathscr{D}_q$ of the quantum space   $\mathscr{D}$. 
Denote by $\gamma_{\bf t^o}$ the induced    automorphism of the algebra 
$\mathcal{A}_\hbar(\mathscr{D})$.

\begin{theorem} \la{12.23.18.1}
 \label{K0} Let $\hbar\in \C^*$ such that $\hbar >-2$. Then the operator ${\cal I}_{\bf c^\circ}$ induces  a map of 
 Schwartz spaces intertwining the automorphism $\gamma_{\bf t^o}$ of $\mathcal{A}_\hbar(\mathscr{D})$:
\be   \label{tu31}
\begin{split}
&{\cal I}_{\bf c^\circ}: {\cal S}_{\bf q'} \lra {\cal S}_{\bf q},\\
&{\cal I}_{\bf c^\circ}  A  (s) =  {\gamma_{\bf c^\circ}(A)} {\cal I}_{\bf c^\circ}(s)\qquad \forall A \in \mathcal{A}_\hbar(\mathscr{D}), ~\forall s \in {\cal S}_{\bf i'}.
\\
\end{split}
\ee
  If the cluster transformation  
${\bf c}^\mathscr{D}_q$  
is identity map, then   
 $
{\cal I}_{\bf c^\circ} = \lambda_{\bf c^\circ}{\rm Id}$,  where 
$|\lambda_{\bf c^\circ}|=1.$

So we get a representation of the $\ast$-algebra $\mathcal{A}_\hbar(\mathscr{D})$ in the Schwarz  subspace of the Hilbert space ${\cal H}_{\mathscr A}$, equivariant under a unitary projective representation of the cluster modular group $\Gamma$ in   ${\cal H}_{\mathscr A}$. \end{theorem} 

\begin{proof} In the case $\beta \in \R$ this is equivalent to   \cite[Theorem 7.7]{FG07} thanks to the identity (\ref{12.12.18.3}).  

The second line in (\ref{tu31}) follows   by analytic continuation from the 
case $\beta \in \R$.  

Here is the  point. Relation (\ref{KYREL}) plus the   power series identity  $\Psi_{1/q}(X)^{-1} = \Psi_{q}(X)$ imply 
$$
\varphi_{\beta}(p_k)  = \frac{\Psi_q(Y_k)}{\Psi_{1/q^\vee}(Y_k^\vee)} =  \Psi_q(Y_k) \Psi_{q^\vee}(Y_k^\vee).
$$
Since the variables $\{Y_k, B_k\}$ and $\{Y^\vee_l, B_l^\vee\}$ commute, the   intertwiner (\ref{*RUu}) acts on the algebra ${\cal O}_q(\mathscr{D})$ as the conjugation by $\Psi_q(Y_k)/\Psi_q(\widetilde Y_k^{-1})$. The latter coincides, by the  
 definition from Section \ref{Sec8.1}, with the mutation $(\mu_k^\sharp)^*$ of the quantum symplectic double. Similarly story for the action of the intertwiner (\ref{*RUu})   on the modular dual algebra 
 ${\cal O}_{q^\vee}(\mathscr{D})$.
 
 The claim in the first line in (\ref{tu31}) is proved  similarly to  \cite[Theorem 7.7]{FG07}.

The last claim   is also deduced from the $\beta \in \R$ case 
 by analytic continuation. \end{proof}
 
\subsubsection{Conclusion.}  The modular double  of the  algebra ${\cal O}_q(\mathscr{D})$ of the symplectic double of a  cluster Poisson 
variety ${\mathscr X}$ is a $\ast$-algebra $\mathcal{A}_\hbar(\mathscr{D})$, defined for any Planck constant $\hbar= \beta^2 \in \C^*$ and $\hbar + \hbar^{-1} \in \R$. 
 It has a canonical 
   representation. By this we mean the following data:

\begin{itemize}

\vskip 1mm\item The Hilbert $L_2$-space ${\cal H}_{{\mathscr A}}$ 
   assigned to the canonical   measure 
$\mu_{\mathscr A}$ on the manifold ${\mathscr A}(\R_{>0})$.

\vskip 1mm\item Representations $\{\rho_{\bf q}\}$  of the   $\ast$-algebra $\mathcal{A}_\hbar(\mathscr{D})$ 
 in  the  Hilbert space  ${\cal H}_{{\mathscr A}}$, assigned 
 to   quivers ${\bf q}$. 
Each representation $\rho_{\bf q}$ is given by unbounded operators,  defined on  a dense   subspace
 \be \nonumber
{\cal S}_{\bf q}\subset{\cal H}_{{\mathscr A}}.
\ee

\vskip 1mm\item Representations $\rho_{\bf q}$  are related by 
   unitary  intertwiners preserving  the 
 subspaces ${\cal S}_{\bf q} $. Precisely,  each cluster transformation  ${\bf c}: {\bf q} \lra {\bf q}'$ gives rise to a unitary operator 
 \be \la{OPI}
\begin{split}
 & {\cal I }_{{\bf c}^\circ}: {\cal H}_{{\mathscr A}} \lra {\cal H}_{{\mathscr A}}, \qquad 
 {\cal I}_{{\bf c}^\circ} ({\cal S}_{\bf q'})= {\cal S}_{\bf q}, \\
  \end{split}
 \ee 
  such that 
 $\forall A\in \mathcal{A}_\hbar$ and $\forall s\in {\cal S}_{\bf q'}$ one has 
 $$
{\cal I}_{{\bf c}^\circ}  \circ  \rho_{\bf q'}  (A) (s)= \rho_{\bf q}  (A)\circ {\cal I}_{{\bf c}^\circ} (s).
 $$
 Trivial cluster transformations   amounts to  unitary scalar intertwiners $\lambda \cdot {\rm Id}$.

 \end{itemize}
 
 Few remarks are in order. 

 \begin{enumerate} \vskip 1mm\item Although the  
  Hilbert space ${\cal H}_{\mathscr A}$ does not depend on  a choice of the quiver ${\bf q}$, the  representations $\rho_{\bf q}$ 
 of the   $\ast$-algebra $\mathcal{A}_\hbar(\mathscr{D})$ do.
  
\vskip 1mm  i) One uses crucially the cluster coordinate system on the space ${\mathscr A}$ assigned to a quiver ${\bf q}$ to define the representation $\rho_{\bf q}$: the  
  generators act as exponents of   first order differential operators in the logarithmic cluster coordinates assigned to   ${\bf q}$. 
 
\vskip 1mm ii)  Intertwiners (\ref{OPI}) involve the quantum dilogarithm function, and  do not reduce to coordinate transformations between cluster coordinate systems. 
 
\vskip 1mm \item To specify a vector in the canonical  representation one needs to specify  a vector in ${\cal H}_{\mathscr A}$ {\bf and} a cluster coordinate system. 
 This is in a sharp contrast with the usual intuition, where one just needs  a    function on a manifold.

\vskip 1mm \item The principal series   $\ast$-representations of the algebra $\mathcal{A}_\hbar({\mathscr X}) $, discussed in Section \ref{SEC8.4}, have similar features, with one notable exception: 
  Hilbert spaces of representations are not as canonical as the one   ${\cal H}_{\mathscr A}$. Indeed, there is no space, like the space ${\mathscr A}$, related to them;  we must choose a polarization of the underlying lattice $\Lambda$ to construct them. Since there is no natural choice for the polarization, and since it is hard to make such choices compatible for quivers
   related by cluster transformation, 
  it is   difficult to construct   a vector in this set-up. The only exception is ${\rm PGL_m}$, where we have the special cluster coordinate systems constructed in \cite{FG03a}. 
 \end{enumerate}

\medskip

\subsection{Canonical $\ast$-representation of the quantum double: the general  case} \la{Sec8.2}
 
\medskip 

In Section \ref{Sec8.2},   we extend  results of Section \ref{Sec8.2a} 
  to the   general  skew-symmetrizable 
  case. The case $\beta \in \R$ was done in \cite{FG07}.  We use the same set-up as in Section \ref{Sec8.2a}.

\subsubsection{Representations.} 
Consider the following  differential 
operators on the manifold $ {\mathscr A}(\R_{>0})$: 
\begin{equation} \label{5.26.08.1aaX}
\begin{split}
& p_k:=   2\pi i  d_k^{-1}\cdot  \frac{\partial}{\partial a_k}
-\alpha_k^+,\\ 
& {{ p}^\circ_k} :=   
 -  2\pi i  d_k^{-1}\cdot \frac{\partial}{\partial a_k}
+ \alpha_k^-,\\
&a_k:= a_k.\\
\end{split}
\end{equation}
Note that   $[ {{ p}^\circ_k},  p_l]=0$. Operators (\ref{5.26.08.1aaX})    satisfy   relations of the Heisenberg $\ast$-algebra with $\hbar=1$:
\be \nonumber
\begin{split}
&[ p_k,  p_l] = 2 \pi i   \cdot   \widehat\varepsilon_{kl}, 
~\qquad~ \qquad ~[ p_k,  a_l] = 2 \pi i   d_k^{-1}\cdot \delta_{kl}, \qquad ~\qquad~
\ast  p_k =  p_k, \qquad \ast a_k = a_k.\\
&[ {{ p}^\circ_k},  {{{ p}^\circ_l} }] 
= -2 \pi i   \cdot   \widehat \varepsilon_{kl}, \qquad  
[ {{ p}^\circ_k},  a_l]  = -2 \pi i  d_k^{-1}  \cdot \delta_{kl}, 
\qquad 
\ast  {{ p}^\circ_k} =  {{ p}^\circ_k}. \\
\end{split}
\ee

 Rescaling  operators (\ref{5.26.08.1aaX}) by $\beta $, we get the following operators: 
 \begin{equation} \nonumber
\begin{split}
&b_k := \beta  \cdot a_k,\\
& y_k:= \beta  \cdot p_k, \\ 
& {{ y}^\circ_k} := \beta \cdot p_k^\circ.
\\
\end{split}
\end{equation}
The operators $\{b_k, y_k\}$   satisfy  the relations of the Heisenberg $\ast$-algebra 
$  {\cal H}_{{\hbar},{\bf q}}$:
\be \nonumber
\begin{split}
&[ y_k,  y_l] = 2 \pi i  \beta^2  \cdot \widehat \varepsilon_{kl}, 
\qquad  ~\qquad~[ y_k,  b_l] = 2 \pi i \beta^2d_k^{-1}\cdot \delta_{kl}, \qquad 
\ast  y_k =  y_k, ~\qquad~ \ast  b_l =  b_l.\\
\end{split}
\ee
 The operators $\{b_k, y^\circ_k\}$   satisfy the  relations of the    opposite  Heisenberg $\ast$-algebra   ${\cal H}^{\rm op}_{{\hbar},{\bf q}}$.

 \vskip 2mm
Let us introduce the Langlands modular dual collection of operators
 \begin{equation} \label{7.26.07.1}
\begin{split}
&b^{\vee}_k := d_k {a}_k/\beta,\\
&y^{\vee}_k:=d_k {p}_k/\beta .\\ 
\end{split}
\end{equation}
They 
satisfy   commutation relations of the Heisenberg $\ast$-algebra 
${\cal H}_{\hbar^{\vee},{\bf q}^{\vee}}$:
\begin{equation} \nonumber
\begin{split}
&[{y}^{\vee}_k, {y}^{\vee}_l] = 2 \pi i \beta^{-2} \cdot d_k d_l  \widehat \varepsilon_{kl} = 2 \pi i \beta^{-2}  \cdot \widehat \varepsilon^\vee_{kl},\\
&[{y}^{\vee}_k, {b}^{\vee}_l] = 2 \pi i \beta^{-2} \cdot d_k \delta_{kl} = 2 \pi i \beta^{-2} \cdot(d_k^{\vee})^{-1}\delta_{kl},\\
&[{b}^{\vee}_k, {b}^{\vee}_l] =0.\\
\end{split}
\end{equation}

\vskip 2mm
The commutators of the operators $\{ {y}_k,  {b}_k\}$ with the ones $\{ {y}^{\vee}_l,  {b}^{\vee}_l\}$ lie in $2\pi i \Z$:
\begin{equation} \label{12.11.18.1j}
\begin{split}
&[{y}_k, {y}^{\vee}_l] = 2 \pi i \cdot \widehat \varepsilon_{kl}d_l  = 2 \pi i  \cdot  \varepsilon_{kl},\\
&[{y}_k, {b}^{\vee}_l] = 2 \pi i  \cdot  \delta_{kl}, \qquad [{b}_k, {b}^{\vee}_l] =0.\\
\end{split}
\end{equation}

  The exponentials $ Y_k:= {\rm exp}( y_k)$ and  $ B_k:= {\rm exp}( b_k)$ are 
 difference operators.  
The operators $\{ {Y}_k,  {B}_k\}$ satisfy the commutation relations  
of the quantum torus algebra ${\cal O}_q(\mathscr{D}_{\bf q})$. 
  Set 
\begin{equation} \label{7.26.07.1y}
\begin{split}
& Y^{\vee}_k:= {\rm exp}( y^{\vee}_k), ~\qquad~ 
 B^{\vee}_k:= {\rm exp}( b^{\vee}_k).\\
\end{split}
\end{equation} 
The operators $\{ {Y}^{\vee}_k,  {B}^{\vee}_k\}$  
satisfy   commutation relations of the algebra  
${\cal O}_{q^{\vee}}(\mathscr{D}_{\bf i^{\vee}})$. 
Thanks to (\ref{12.11.18.1j}) they  commute with the ones $\{ {Y}_k,  {B}_k\}$. 
Therefore we get a representation of the algebra
$$
\mathcal{A}_\hbar(\mathscr{D}_{\bf q})= {\cal O}_{q}(\mathscr{D}_{\bf q}) \otimes {\cal O}_{q^{\vee}}(\mathscr{D}_{\bf q^{\vee}}).
$$ 

\subsubsection{The $\ast$-algebra structures.}  
If  $\beta \in \R$, we get a $\ast$-algebra structure on the algebras ${\cal O}_{q}(\mathscr{D}_{\bf q})$ and    $\mathcal{A}_\hbar(\mathscr{D}_{\bf q})$, denoted 
by $\ast_{\rm R}$.  The above formulas define representations of these $\ast$-algebras in the Hilbert space ${\cal H}_{\mathscr A}$.

If $|\beta|=1$,  the Hilbert space structure provides an involution $\ast_\U$ such that
\be \la{12.23.18.100}
\ast_{\rm U}(Y_k^{d_k}) = Y_k^\vee, \qquad \ast_{\rm U}(B_k^{d_k}) = B_k^\vee.
\ee
Therefore there is a subalgebra ${\mathscr A}'_\hbar(\mathscr{D}_{\bf q}) \subset \mathcal{A}_\hbar(\mathscr{D}_{\bf q})$ on which  $\ast_\U$ induces a $\ast$-algebra structure. 
However the intertwiner operators (\ref{INTOP}) seem to be non-unitary.

\subsubsection{Intertwiners.} 
Recall the skew-symmetrizer $d_k$. Consider the function 
\be \nonumber
\begin{split}
&\widetilde \varphi_{\beta_k}(z) := {\rm exp}\Bigl(-\frac{1}{4}\int_{\Omega}\frac{e^{-i pz}}{ {\rm sh} (\pi  \beta^{-1}\cdot  p)
{\rm sh} (\pi   \beta d^{-1}_k \cdot p) } \frac{dp}{p} \Bigr).\\
\end{split}
\ee
It is related to the quantum dilogarithm function $\Phi_{\hbar_k}(z)$ in (\ref{12.12.18.8}),  used in \cite{FG07},  by 
\be \nonumber
\begin{split}
&\Phi_{\hbar_k}(z) = \widetilde \varphi_{\beta_k}(z/\beta).\\
\end{split}
\ee
Same argument as in the proof of Lemma \ref{12.15.18.1} shows that if $\beta \in \R$, then 
\be \la{RELA}
\overline {\widetilde \varphi_{\beta_k}(z)} =\widetilde \varphi_{\beta_k}(\overline z)^{-1}.
\ee

The operators $p_k$ and 
$ {p^\circ_k}$ are self-adjoint. So there are     operators $\widetilde \varphi_{\beta_k}(p_k)$ and $\widetilde \varphi_{\beta_k}(p^\circ_k)$.

\begin{definition} \label{intertwiner}
Given a mutation $\mu_k: {\bf q} \to {\bf q}'$, 
the    intertwining operator 
\begin{equation} \label{12.01.06.1XX}
{\cal I}_{{\bf q}'\to {\bf q}}: L^2({\mathscr A}(\R_{>0})) \lra 
L^2({\mathscr A}(\R_{>0}))
\end{equation}
is defined as   
$$
{\cal I}_{{\bf q'}\to {\bf q}}:= {\bf K}^{\sharp}\circ 
{\bf K}^{'}, 
$$ 
where 
 the  unitary operator ${\bf K}^{'}$  is induced by    linear map (\ref{11.18.06.10d}), and 
  \be \la{INTOP}
 {\bf K}^{\sharp}:=  \widetilde \varphi_{\beta_k}( p_k)
\widetilde \varphi_{\beta_k}(  -{p}^\circ_k)^{-1}.
\ee
\end{definition} 

If $\beta \in \R$ and $x \in \R$, then (\ref{RELA}) implies that $|\widetilde \varphi_{\beta_k}(x)|=1$, and so 
    operators $\widetilde \varphi_{\beta_k}(p_k)$ and $\widetilde \varphi_{\beta_k}(p^\circ_k)$  are unitary. 
    Therefore in this case the intertwiner (\ref{12.01.06.1XX}) is unitary. 
  
\vskip 2mm

  Recall the dense subspace   $W_{\bf q} $ in (\ref{9.4.07.1}).  
 Then $\{B_k, Y_k,  B^{\vee}_k,  Y^{\vee}_k\}$ are symmetric unbounded operators 
   preserving $W_{\bf q}$, 
and satisfying on it the standard  commutation relations. 
 Recall
\be \nonumber
\mathcal{A}_\hbar(\mathscr{D}_{\mathscr X}):= {\cal O}_q(\mathscr{D}_{\mathscr X})\otimes  {\cal O}_{q^\vee}(\mathscr{D}_{  {\mathscr X}^\vee}).
\ee

Following Definition \ref{7.10.07.1}, we define  
 the common domain ${\cal S}_{\bf q}$ of definition of
operators from $\mathcal{A}_\hbar(\mathscr{D}_{\mathscr X})$ in the realization assigned to a quiver ${\bf q}$.   

Just as in Section \ref{Sec8.2a}, a 
  cluster transformation  
${\bf c}: {\bf q} \to {\bf q'}$ provides a unitary operator 
$$
{\cal I}_{\bf c^\circ}: {\rm L}^2({\mathscr A}(\R_{>0})) \lra {\rm L}^2({\mathscr A}(\R_{>0})),  
$$
and  gives rise to an 
 automorphism   $\gamma_{\bf c^\circ}$   of the algebras  ${\cal O}_q(\mathscr{D}_{\mathscr X})$ and $\mathcal{A}_\hbar(\mathscr{D}_{\mathscr X})$.

\begin{theorem}
 \label{K0} The operator ${\cal I}_{\bf c^\circ}$ induces  a map  
 $ {\cal S}_{\bf q'} \to {\cal S}_{\bf q}$ intertwining the automorphism $\gamma_{\bf c^o}$:
\be   \nonumber
\begin{split}
&{\cal I}_{\bf c^\circ}  A  (s) =  {\gamma_{\bf c^\circ}(A)} {\cal I}_{\bf c^\circ}(s)\qquad \forall A \in {\cal O}_q(\mathscr{D}_{\mathscr X}), ~\forall s \in {\cal S}_{\bf q'}.
\\
\end{split}
\ee
  If the cluster transformation  
${\bf c}$ induces the identity map of the   quantum space   $\mathscr{D}_q$, then   
$$
{\cal I}_{\bf c^\circ} = \lambda_{\bf c^\circ}{\rm Id}, \qquad
|\lambda_{\bf c^\circ}|=1.
$$

Therefore if $\beta \in \R \cup i\R$, we get a $\Gamma$-equivariant representation of the $\ast$-algebra $\mathcal{A}_\hbar(\mathscr{D}_{\mathscr X})$ in the Schwarz  subspace of the Hilbert space ${\cal H}_{\mathscr A}$.
\end{theorem} 

 If  if $\beta \in \R$, this is \cite[{\rm Theorem 7.7}]{FG07}.

   \medskip
  
  \subsection{Principal series of $\ast$-representations of quantized cluster varieties } \la{SEC8.4}

\medskip

In Section \ref{SEC8.4} we start from an arbitrary cluster Poisson variety ${\mathscr X}$, assigned to a quiver ${\bf q}$.

Definition \ref{DEFINV} provides   a $\ast$-algebra structure on the modular double $\mathcal{A}_\hbar$ of the quantum torus algebra assigned to a lattice with a 
skew-symmetric  bilinear form. We   apply it to a quiver ${\bf q}$. Recall the corresponding split Poisson torus ${\mathscr X}_{\bf q}$. Let us  summarize the result. 

Let  ${\cal H}_{\bf q}$ be the Hilbert space given by the Weil representation construction in Theorem \ref{Weil}.  

\bl \la{12.23.18.2} Let ${\bf q}$ be any  skew-symmetric quiver, and   $\hbar \in \C^*$ such that $\hbar  + \hbar ^{-1}\in \R$. 

\vskip 1mm1) If $ \hbar  \in \R$,  
then  there is a $\ast$-algebra structure   on the algebra  ${\cal O}_q({\mathscr X}_{\bf q})$ and  hence on 
 $\mathcal{A}_\hbar({\mathscr X}_{\bf q})$.

\vskip 1mm 2) If  $|\hbar |=1$, then  there is a $\ast$-algebra  structure on the algebra $\mathcal{A}_\hbar({\mathscr X}_{\bf q})$. 
 
\vskip 1mm 3)  These  $\ast$-algebras   
 have   representations in the subspace $W_{\bf q}$ of the Hilbert space $ {\cal H}_{\bf q}$.
 \el

The basis $\{e_i\}$ of the lattice $\Lambda$ is represented by selfadjoint operators $\{p_i\}$ in the Hilbert space $ {\cal H}_{\bf q}$.    
So one can apply to the operator $p_k$ any continuous function $\varphi(t)$ on $\R$. 
If $|\varphi(t)|=1$, we get a unitary operator.

\begin{definition} \label{intertwiner-betaX}
Let $\beta\in \C^*$ such that $(\beta + \beta^{-1})^2 >0$. 

Given a mutation $\mu_k: {\bf q} \to {\bf q}'$, we consider 
the  intertwining operator 
\begin{equation} \label{12.01.06.1X}
{\cal I}_{{\bf q}'\to {\bf q}}:= {\bf K}^{\sharp}\circ 
{\bf K}^{'} : {\cal H}_{\bf q'} \lra 
{\cal H}_{\bf q}, \qquad \mbox{where}
\end{equation}
\begin{itemize}
\vskip 1mm\item
{\it The operator ${\bf K}^{\sharp}:=  \varphi_{\beta}( p_k) : {\cal H}_{\bf q} \lra 
{\cal H}_{\bf q}$} is  the 
  quantum dilogarithm of the operator $p_k$.

\vskip 1mm\item
{\it The operator ${\bf K}^{'} : {\cal H}_{\bf q'} \lra 
{\cal H}_{\bf q}$} is induced by  the    mutation  of the  basis $e_i$ of $\Lambda$. 
\end{itemize}
\end{definition}

  Recall    the  algebra of universally Laurent polynomials ${\cal O}_q({\mathscr X})$ and its modular double 
\be \nonumber
\mathcal{A}_\hbar( {\mathscr X}):= {\cal O}_q( {\mathscr X})\otimes {\cal O}_{ q^\vee}( {{\mathscr X}}).
\ee
 We define the subspace ${\cal S}_{\bf q}$
as 
 the common domain of definition of
operators from $\mathcal{A}_\hbar({\mathscr X})$ in ${\cal H}_{\bf q}$.

\vskip 2mm

  A 
  cluster transformation  
${\bf c}: {\bf q} \to {\bf q'}$ provides a unitary operator 
$$
{\cal I}_{\bf c^\circ}: {\cal H}_{\bf q'} \lra {\cal H}_{\bf q}.  
$$
Indeed,  we assigned in (\ref{12.01.06.1X}) to   a   
 mutation ${\bf q} \to {\bf q'}$ 
an intertwiner 
${\cal I}_{{\bf q'} \to {\bf q}}$. An  automorphism $\sigma$ 
of a quiver ${\bf q}$  gives rise 
to 
a unitary operator. The  
 cluster transformation  
 ${\bf c}$  is a composition of mutations and automorphisms. 
  The map 
${\cal I}_{\bf c^\circ}$ is the reverse 
composition of the  intertwiners.

A  cluster transformation ${\bf c}$ gives rise to 
a   cluster transformation    
${\bf c}^{\mathscr X}_q$ of the quantum space   ${\mathscr X}_q$. 
Denote by $\gamma_{\bf c^\circ}$ the induced    automorphism of the algebra 
$\mathcal{A}_\hbar( {\mathscr X})$.

\begin{theorem} \la{12.23.18.11}
 \label{K0} Let $\hbar \in \C^*$ such that $\hbar  + \hbar ^{-1}\in \R$. The operator ${\cal I}_{\bf c^\circ}$ induces  a map of 
 Schwartz spaces ${\cal I}_{\bf c^\circ}: {\cal S}_{\bf q'} \lra {\cal S}_{\bf q}$ intertwining the automorphism $\gamma_{\bf c^\circ}$ of $\mathcal{A}_\hbar( {\mathscr X})$:
\be  \nonumber
\begin{split}
&{\cal I}_{\bf c^\circ}  A  (s) =  {\gamma_{\bf c^\circ}(A)} {\cal I}_{\bf c^o}(s)\qquad \forall A \in \mathcal{A}_\hbar( {\mathscr X}), ~\forall s \in {\cal S}_{\bf i'}.
\\
\end{split}
\ee

  If the cluster transformation  
${\bf c}^{\mathscr X}_q$  
is identity map, then   
$
{\cal I}_{\bf c^\circ} = \lambda_{\bf c^\circ}{\rm Id}, \qquad
|\lambda_{\bf c^\circ}|=1.
$

So we get a $\Gamma$-equivariant representation of the $\ast$-algebra $\mathcal{A}_\hbar( {\mathscr X})$ in the  subspace  ${\cal S}_{\bf q} \subset {\cal H}_{\bf q}$. \end{theorem} 

\begin{proof} In the case $\beta \in \R$ this is equivalent to   \cite[Theorem 7.7]{FG07} thanks to the identity (\ref{12.12.18.3}).  
The general case is deduced   from this by 
 the analytic continuation, as was spelled in the proof of Theorem \ref{12.23.18.1}. 
Here are two key  points. 

 1. Just like in the proof of Proposition \ref{MUMAP},   the  mutation map $\mu_k^*: \mathcal{A}_\hbar({\mathscr X}_{\bf ic}) \lra {\mathscr A}_{\hbar^\vee}({\mathscr X}_{\bf q})$ commutes with the involution $\ast_{\rm U}$.  Therefore there is a well defined $\ast$-algebra structure on the algebra $\mathcal{A}_\hbar({\mathscr X})$, which in any 
 cluster coordinate system is given by the involution $\ast_\U$.   
 
 The similar claim for the involution $\ast_{\rm R}$ was proved in \cite{FG03b}. 
 
 2. Relation (\ref{KYREL}) plus the identity $\Psi_{1/q}(X)^{-1} = \Psi_{q}(X)$ imply 
$$
\varphi_{\beta}(p_k)  = \frac{\Psi_q(X_k)}{\Psi_{1/q^\vee}(X_k^\vee)} =  \Psi_q(X_k) \Psi_{q^\vee}(X_k^\vee).
$$
Variables $\{X_k\}$ and $\{X^\vee_l\}$ commute. So  the  intertwiner  $\varphi_{\beta}( p_k)$ acts on the algebra ${\cal O}_q({\mathscr X})$ as   ${\rm Ad}_{\Psi_q(X_k)}$, which is    
 just the mutation $\mu_k^\sharp$ of the quantum cluster variety ${\mathscr X}_q$. 
 Similarly   $\varphi_{\beta}( p_k)$ acts on the  algebra 
 ${\cal O}_{q^\vee}({\mathscr X})$ as   ${\rm Ad}_{\Psi_{q^\vee}(X^\vee_k) }$. 
  \end{proof}

 {The $\hbar \lms - \hbar$ involution.} Given a complex vector space $V$, the complex conjugate vector space $\overline V$ is  the abelian group $V$ with 
   the new action $\circ$ of $\C$ given by $\lambda \circ v:= \overline \lambda v$. 
   If $V$ had a non-degenerate hermitian form, there is a canonical isomorphism of complex vector spaces $V^* = \overline V$. 
   
   By Lemma \ref{6.9.03.11} there are canonical antiholomorphic isomorphisms 
   $$
  i_{\cal O}: {\cal O}_{q^{-1}}({\mathscr X}^\circ)   \stackrel{\sim}{\lra} {\cal O}_q({\mathscr X}), \qquad i_{\mathscr A}: {\mathscr A}_{-\hbar}({\mathscr X}^\circ)  \stackrel{\sim}{\lra} \mathcal{A}_\hbar({\mathscr X}).   
   $$
   Therefore taking a $\ast-$representation of the algebra $\mathcal{A}_\hbar({\mathscr X})$   in a Hilbert space ${\cal H}_{{\mathscr X}, \hbar}$ and composing the antiholomorphic isomorphism $i_{\mathscr A}$  
   with this representation, we get a $\ast-$representation of the algebra $ {\mathscr A}_{-\hbar}({\mathscr X}^\circ) $ in the complex conjugate Hilbert space $\overline {{\cal H}_{{\mathscr X}, \hbar}}$. 
   Note that $\overline {{\cal H}_{{\mathscr X}, \hbar}} =  {{\cal H}_{{\mathscr X}, \hbar}}^*$. Since the intertwiners are defined uniquely up to a unitary constant by their 
   intertwining property, for generic $\hbar$ it must be isomorphic to $ {{\cal H}_{{\mathscr X}^\circ, -\hbar}}$. In the case $\hbar = -1$ when the representation  of 
    ${\mathscr A}_{\hbar = -1}({\mathscr X}^\circ) $      is not defined, we can use the representation $\overline {\cal H}_{{\mathscr X}, \hbar=1}$ of ${\mathscr A}_{\hbar = 1}({\mathscr X})$   instead. 
   
 {The general  case.} Let ${\bf q}$ be an arbitrary skew-symmetrizable quiver. Then there is a $\ast$-algebra structure on ${\cal O}_q({\mathscr X})$  and hence on $\mathcal{A}_\hbar({\mathscr X})$    for any $\beta \in \R$, but not  for $|\beta|=1$. 
   
   As in   Lemma \ref{12.23.18.2}, there is a representation    of the $\ast$-algebra $\mathcal{A}_\hbar({\mathscr X})$ in a Hilbert space 
   ${\cal H}_{\bf q}$.
   
The intertwiner (\ref{12.01.06.1X}) is defined the same way,  but using the function $\varphi_{\beta_k}$ rather then $\varphi_{\beta}$:  
   $$
   {\bf K}^{\sharp}:=  \varphi_{\beta_k}( p_k) : {\cal H}_{\bf q} \lra 
{\cal H}_{\bf q}.
   $$
  The analog of Theorem \ref{12.23.18.11} holds for any  $\beta \in \R$. 
  However the function $\varphi_{\beta_k}(x)$ is not unitary if $|\beta|=1$.

 \subsubsection{Conclusion.}  Given a cluster Poisson variety ${\mathscr X}$, the modular double $\mathcal{A}_\hbar({\mathscr X})$ of the  algebra ${\cal O}_q( {\mathscr X})$  is a $\ast$-algebra, defined for any Planck constant $\hbar \in \C^*$ and $\hbar  + \hbar ^{-1}\in \R$. 
 It has a principal series of 
   representation. By this we mean the following data:

\begin{itemize}

\vskip 1mm\item A Hilbert $L_2$-space ${\cal H}_{{\bf q}}$ 
   assigned to a quiver ${\bf q}$ by the Weil representation construction. 

\vskip 1mm\item A representation $\rho_{\bf q}$  of the   $\ast$-algebra $\mathcal{A}_\hbar( {\mathscr X})$ 
 in  the  Hilbert space  ${\cal H}_{{\bf q}}$,  given by unbounded operators,  defined on  a dense   subspace $
{\cal S}_{\bf q}\subset{\cal H}_{{\bf q}}.
$

\vskip 1mm\item Representations $\rho_{\bf q}$  are related by 
   unitary  intertwiners preserving  the 
 subspaces ${\cal S}_{\bf q} $. Precisely,  each cluster transformation  ${\bf c}: {\bf q} \lra {\bf q}'$ gives rise to a unitary isomorphism 
 \be \nonumber
\begin{split}
 & {\cal I}_{{\bf c}^\circ}: {\cal H}_{{\bf q}'} \lra {\cal H}_{{\bf q}}, \qquad 
  {\cal I}_{{\bf c}^\circ} ({\cal S}_{\bf q'} ) = {\cal S}_{\bf q}, \\
  \end{split}
 \ee 
  such that 
 $\forall A\in \mathcal{A}_\hbar$ and $\forall s\in {\cal S}_{\bf q'}$ one has 
 $$
{\cal I}_{{\bf c}^\circ}  \circ  \rho_{\bf q'}  (A) (s)= \rho_{\bf q}  (A)\circ {\cal I}_{{\bf c}^\circ} (s).
 $$
 Trivial cluster transformations   amounts to  unitary scalar intertwiners $\lambda \cdot {\rm Id}$.

 \end{itemize}
 
 Finally  we want to stress again the following.   
\be \nonumber
\begin{split}
&\mbox{\it To specify a vector in a representation of the quantized cluster Poisson variety 
 one needs}\\ 
 &\mbox{\it  to specify a cluster  coordinate system assigned to a quiver ${\bf q}$, {\bf and} a vector in ${\cal H}_{\bf q}$.} 
\end{split}
\ee
  
        

\medskip

  \section{The Quantum Lift Theorem} \la{sec18} 
  \la{Sec.quant.s}
 
\medskip

\subsubsection{The set-up.} \la{sec18.1}We shall fix the following notations.
 
\vskip 1mm 1. $m\leq n$ are positive integers.  Set $[m]:=\{1,\ldots, m\}$, and  $[n]:=\{1, \ldots, n\}$. 

\vskip 1mm 2.  $\Lambda$ is a free $\mathbb{Z}$-module of dimension $n$.

\vskip 1mm 3.  $\mathcal{E}=\{e_1,\ldots, e_n\}$ is a basis of $\Lambda$.

\vskip 1mm 4.  $(\ast, \ast): ~\Lambda\times \Lambda \rightarrow \mathbb{Q}$ is a skew-symmetric bilinear map such that $\widehat{\varepsilon}_{ij}=(e_i, e_j)$.

\vskip 1mm 5.  $\{d_1, \ldots, d_m\}$ is a collection of positive integers such that  
 \[\varepsilon_{ij}:= \widehat{\varepsilon}_{ij}d_j \in \mathbb{Z}, \hskip 7mm \forall (i,j )\in [n]\times[m]. 
 \]

\vskip 1mm 6.   $\mathcal{F}=\{f_1, \ldots, f_n\}$ is a basis of $\Lambda_{\mathbb{Q}}= \Lambda\otimes \mathbb{Q}$ such that

\vskip 1mmi) $\Lambda$ is contained in the $\mathbb{Z}$-linear span of $\mathcal{F}$,

\vskip 1mmii) $(f_i, e_j)= d_j^{-1}\delta_{ij}$, where $(i, j)\in [n]\times [m]$.
 
In this case, we say that the pair $(\mathcal{F}, \mathcal{E})$ is  compatible.

 \vskip 2mm
 
Let us fix a compatible pair and write $e_i=\sum_{j \in [n]} p_{ij}f_j$. 
\begin{lemma} \la{Lem18.1} For $(i,j)\in [n]\times[m]$, we have $p_{ij}= \varepsilon_{ij}$.
\end{lemma}
\begin{proof}
It follows from the identity $(e_i, e_j)= d_j^{-1}p_{ij}=\widehat{\varepsilon}_{ij}$.
\end{proof}

We define a mutation at $j\in [m]$,   leading to  new bases $\mathcal{E}'=\{e_1', \ldots, e_n'\}$ and  $\mathcal{F}'=\{f_1', \ldots, f_m'\}$: 
\[
e_i':= \left\{ \begin{array}{ll} e_i+[\varepsilon_{ij}]_+e_j &\mbox{if } i\neq j,\\
-e_j &\mbox{if } i=j. \end{array} \right.  \qquad f_i':= \left\{ \begin{array}{ll} 
f_i &\mbox{if } i\neq j, \\
-f_i+\sum_{k\in [n]}[-p_{jk}]_+f_k &\mbox{if } i= j.\\ \end{array} \right.\]

\begin{lemma} The pair $(\mathcal{F}',\mathcal{E}')$ is compatible.
\end{lemma}
\begin{proof} It is due to the following direct calculations.

1.  $(f_j', e_j')=(-f_i+\sum_{k\in [n]}[-p_{jk}]_+f_k, -e_j)= (-f_j, -e_j)=d_j^{-1}$.

2.  For $i\neq j$, we have $(f_i', e_j')=(f_i, -e_j)=0$.

3. For $k\in [m]-\{j\}$ and $i\in [n]-\{j\}$, we have 
$(f_i', e_k')= (f_i, e_k+[\varepsilon_{kj}]_+ e_j) = (f_i, e_k)=d_k^{-1}\delta_{ik}.$

4. For $k\in [m]-\{j\}$, we have
\be \nonumber
\begin{split}
(f_j', e_k')&= \Bigl(-f_j+\sum_{s\in [n]}[-p_{js}]_+ f_s, ~ e_k+ [\varepsilon_{kj}]_+ e_j \Bigr)   \\
&= -\left(f_j, [\varepsilon_{kj}]_+ e_j\right)+ \left([-p_{jk}]_+ f_k, e_k\right)  
=-[\varepsilon_{kj}]_+d_j^{-1}+ [-p_{jk}]_+d_k^{-1} =0.  \\
\end{split}
\ee
\end{proof}
 Let $p'=(p_{ik}')$ such that $e_i'= \sum_{k\in [n]} p_{ik}' f_k'$. Then a direct calculation shows that
 \[
 p_{ik}'= \left\{ \begin{array}{ll} -p_{ik} &  \mbox{if } j\in \{i,k\}, \\
 p_{ik}+[p_{ij}]_+p_{jk}+p_{ij}[-p_{jk}]_+ & \mbox{if } j\notin \{i,k\}. \\
 \end{array}\right.
 \]
\subsubsection{Quantization.} We assign to each $v\in \Lambda_{\mathbb{Q}}$   a variable $T_v$ such that
$
T_v T_w = q^{(v, w)}T_{v+w} 
$ 
and set 
\[
X_i :=T_{e_i}, \qquad A_i:= T_{f_i}.
\]
Recall the quantum dilogarithm $\Psi_q(x) = \prod_{n=0}^\infty (1+q^{2n+1}x)^{-1}$. 
Let $q_j= q^{1/d_j}$. Define
\[
X_i'= {\rm Ad}_{\Psi_{q_j}(X_j)} T_{e_i'}, \hskip 5mm A_i'= {\rm Ad}_{\Psi_{q_j}(X_j)} T_{f_i'}.
\]
Note that the map $X_i \lms X_i'$ is the quantized cluster transformation introduced in \cite{FG03b}.
\bt \la{quan.equi.s}
We have
\[
A_i'=\left\{ \begin{array}{ll} T_{f_i} & \mbox{if } i\neq j, \\
T_{f_j'}+T_{f_j'+e_j}  &\mbox{if } i=j.\\
\end{array}\right.
\]
\et
\begin{proof} Note that $(e_j, f_i')= d_j^{-1}\delta_{ij}$. Therefore
\[
A_j' =T_{f_j'}\Psi_{q_j}(q_j^2 X_j) \Psi_{q_j}(X_j)^{-1}= T_{f_j'}(1+q_j T_{e_j}) = T_{f_j'}+ T_{f_j'+ e_j}.
\]
If $i\neq j$, then $A_{i}'= T_{f_i'}=A_i$.
\end{proof}

\subsubsection{Comparison with  Berenstein-Zelevinsky's quantization \cite{BZ}.} \la{sec18.3}
Recall the matrix $p=(p_{ij})$ such that $e_i= \sum_{j\in [n]}p_{ij}f_j$, the skewsymmetric matrix $\widehat{\varepsilon}$, and the  $n\times n$ diagonal matrix $D$ given by 
$D={\rm diag}(d_1,\ldots, d_m, 1,\ldots, 1)$. We have $\varepsilon=\widehat{\varepsilon}D$.

Let $\widetilde{B}=(b_{ij})$ be an $m\times n$ matrix such that $b_{ij}=p_{ji}$. So $\widetilde{B}^t$ is given by the first $m$ rows of $p$.

Let $d$ be a negative integer such that the following matrix has integral coefficients: 
\be \la{comp.BZ.pair}
\Pi = (\pi_{ij}):= d p^{-1}\widehat{\varepsilon} (p^{-1})^t. 
\ee

\begin{lemma} 
\la{lem.comp.BZ.pair}
$(\widetilde{B}, \Pi)$ is a compatible pair in the sense of \cite[\S 3]{BZ}, i.e.,

\vskip 1mm1.  The matirx $\Pi$ is skew-symmetric.

\vskip 1mm2. The $m\times n$ matrix $\widehat{B}^t \Pi =  ( D', 0)$, when $D'= {\rm diag}(-d/d_1, \ldots, -d/d_n)$.
 
\end{lemma}
\begin{proof}
The first claim is clear. 
By definition, $p=   \begin{pmatrix} 
      \widetilde{B}^t \\
      C \\
   \end{pmatrix}$. 
   The first $n$ columns of $\varepsilon$ and $P$ coincide. So
   \[
   p \Pi = p dp^{-1} \widehat{\varepsilon}(p^{-1})^t = -d \widehat{\varepsilon}^t (p^{-1})^t = - d (p^{-1}\widehat{\varepsilon})^t.  
   \]
   Comparing the first $n$ rows of LHS and RHS, we get the second claim.
\end{proof}

By definition $(f_i, f_j)= \lambda_{ij}/d$.   We have
$
T_{f_i} T_{f_j} = q^{\pi_{ij}/d} T_{f_i+f_j}.
$

Recall that $A_i= T_{f_i}$. We consider mutation at $j$. Let us set 
\be \nonumber
\begin{split}
&f_j':= -f_j +\sum_{k} [-p_{jk}]_+ f_k = - f_j - \sum_{k ~|~ b_{kj}<0} b_{kj} f_k,\\
&f_j'' = - f_j'+e_j =- f_j + \sum_{k ~|~ b_{kj}>0} b_{kj} f_k.\\
\end{split}
\ee
Therefore $A_j'= T_{f_j'} + T_{f_j''}$. 
It recovers formula (4.23)  of \cite{BZ}. 

To summarize, the data 1)-5) at the beginning of this section gives rise to a cluster Poisson variety $\mathscr{X}$. Its coordinate ring $\mathcal{O}(\mathscr{X})$ admits a quantization $\mathcal{O}_q(\mathscr{X})$ introduced in \cite{FG03b}. Adding data $6)$, we obtain a compatible pair $(\tilde{B}, \Lambda)$ of \cite{BZ}. It further gives rise to a quantized upper cluster algebra $\mathcal{O}_q(\mathcal{A})$ following the procedure of \cite{BZ}. Theorem \ref{quan.equi.s} asserts that these two quantizations are compatible, i.e., there is a natural embedding 
$
p^*: \mathcal{O}_q(\mathscr{X}) \hlra \mathcal{O}_q(\mathcal{A}).
$
 In particular, if  $\mathcal{F}$ is a basis of  $\Lambda$, then $p^\ast$ is an isomorphism. Its classical limit generalizes the $p$-map in \cite{FG03a}.

\subsubsection{Standard monomials  and the quantum lift.} 
Every $v=\sum_{i=1}^m \lambda_i e_i \in \Lambda$ gives rise to a rational function $x_v := x_1^{\lambda_1}\ldots x_m^{\lambda_m}$ of $\mathscr{X}$.\footnote{We use the notation $x_i$ for cluster Poisson coordinates, reserving $X_i$ for their quantum counterparts.}
 We say $x_v$ is a {\it standard monomial} if $x_v \in \mathcal{O}(\mathscr{X})$.

\begin{proposition} 
\la{standard monomial}
The function $x_v$ is a standard monomial  if and only if $(v, e_i)\geq 0$ for every unfrozen basis vector $e_i\in \Lambda$.
\end{proposition}
  
  \begin{proof} Let us double the lattice $\Lambda$, setting  $\widetilde{\Lambda}:= \Lambda \bigoplus \oplus_{i=1}^n \mathbb{Z}e_i'$, and extend the form $(\ast, \ast)$ so that   
  \[
  (e_i', e_j)=  d_j^{-1}\delta_{ij}, \hskip 7mm (e_i', e_j') =0. 
  \]
  The lattice $\widetilde{\Lambda}$ with $\widetilde{\mathcal{E}}=\{e_1, \ldots, e_n, e_1', \ldots, e_n'\}$ determines a cluster Poisson variety $\widetilde{\mathscr{X}}$.  There is a  surjection $j: \widetilde{\mathscr{X}} \to {\mathscr{X}}$.
Note that $\widetilde{\Lambda}$ admits another basis $\widetilde{\mathcal{F}}$ such that  the pair  $(\widetilde{\mathcal{F}}, \widetilde{\mathcal{E}})$ is compatible. It determines a cluster $K_2-$variety $\widetilde{\mathscr A}$ and an isomorphism $p: \widetilde{\mathscr A}\to \widetilde{\mathscr X}$.
Consider the composition: 
$$
\kappa: \widetilde{\mathscr A} \stackrel{p}{\lra} \widetilde{\mathscr X} \stackrel{j}{\lra} {\mathscr X}. 
$$
Hence we get an injection $\kappa^*: \mathcal{O}(\mathcal{\mathscr X})\hra  \mathcal{O}(\widetilde{\mathscr A})$, and $x_v \in \mathcal{O}(\mathcal{\mathscr X})$ if and only if $\kappa^* (x_v) \in \mathcal{O}(\widetilde{\mathscr A})$.
Let $(A_1, \ldots, A_n, A_1', \ldots, A_n')$ be cluster variables of $\widetilde{\mathscr A}$ associated to $\widetilde{\Lambda}$. By definition, we have
\[
\kappa^*(x_v) = \mathbb{A} \prod_{i=1}^m A_i^{(v, e_i)d_i}, \qquad  \mbox{{where } $\mathbb{A}$ is a Laurent monomial of frozen variables}.
\]

If $(v, e_i)\geq 0$ for every unfrozen $e_i\in \Lambda$, then $\kappa^* (x_v)$ is a   monomial. By  Laurent Phenomenon theorem,   $\kappa^* (x_v) \in \mathcal{O}(\widetilde{\mathscr A})$. 
If $(v, e_i)<0$ for some unfrozen $e_i$, then after mutation at $i$, $\kappa^*(x_v)$ is not a Laurent polynomial, and therefore $\kappa^* (x_v) \notin \mathcal{O}(\widetilde{\mathscr A})$.
  \end{proof}
  
Let ${\rm X}_v$ be the quantum Laurent monomial associated to $v\in \Lambda$. By the Laurent Phenomenon for quantum cluster variables and the same argument as in the proof of  Proposition  \ref{standard monomial}, we conclude that ${\rm X}_v \in \mathcal{O}_q(\mathscr{X})$ if and only if $(v, e_i)\geq 0$ for every unfrozen $e_i\in \Lambda$.

\bd
\la{quantum.promotion.cls}
A regular function $f \in \mathcal{O}(\mathscr{X}) $ is   a {\em standard monomial} if it is a standard monomial in one cluster chart of $\mathscr{X}$. Let $\mathcal{S}_{\mathscr{X}}$ be the collection of standard monomials, and  $\mathcal{S}_{\mathscr{X}}^q$ its  quantum counterpart. 
\ed

The specialization at $q=1$ provides a  projection ${\rm Sp}: \mathcal{S}_{\mathscr{X}}^q \rightarrow \mathcal{S}_{\mathscr{X}}$, called the {\em specialization map}. 
It is evidently  equivariant under the action of the   cluster modular group $\mathcal{G}_{\mathscr{X}}$. 

\btc 
\la{quantum.promotion.f}
The specialization map is a bijection. 

The inverse map ${\rm Sp}^{-1}: \mathcal{S}_{\mathscr{X}} \to \mathcal{S}_{\mathscr{X}}^q$ is called the {\em quantum lift}.
\etc

\begin{proof}
  
Let ${\bf c}: {\bf s}'\rightarrow {\bf s}$ be a quasi-cluster transformation. It induces an isomorphism of the fraction fields of the quantum torus algebras
$ 
{\bf c}_q^{\mathscr X}:~ {\rm Frac} ~ {\cal O}_q({\mathscr X}_{\bf s})  \lra  {\rm Frac} ~{\cal O}_q({\mathscr X}_{\bf s'}). 
$

Let ${\bf c}^{\mathscr X}$ be its classical counterpart. Denote by $\Lambda$ the lattice underlying quivers  ${\bf s}$  and ${\bf s}'$. 
\bp
\label{quantum.promotion.lemma}
Let $x_\lambda$, $\lambda \in \Lambda$ be a standard monomial associated to the seed ${\bf s}$. Suppose that ${\bf c}^{\mathscr X}(x_{\lambda})= x_{\mu}'$ remains a Laurent monomial in ${\bf s}'$. Then ${\bf c}_q^{\mathscr X}({\rm X}_\lambda)= {\rm X}_\mu'$.
\ep
\begin{proof} Every quantum quasi-cluster transformation can be expressed as  a  monomial transformation followed by a conjugation of a product of quantum dilogs  $\psi:=\Psi_q(Y_1^{\varepsilon_1})^{\varepsilon_1} \ldots \Psi_q(Y_n^{\varepsilon_n})^{\varepsilon_n}$, where $\varepsilon_k \in \{\pm 1\}$. By the sign-coherence of $c$-vectors, there is a unique sequence of $\varepsilon_k$ such that every $Y_k^{\varepsilon_k}$ is a monomial of $\{\X_i\}$ with positive exponents. Therefore
$
\psi:= 1+  \mbox{higher order terms}.
$ 
We refer the readers to \cite[\S 2]{GS16} for more details. In particular, it follows that
\be \nonumber
{\bf c}_q^{\mathscr X}({\rm X}_\lambda) = \X'_\mu(1+ \mbox{higher order terms}) \qquad \mbox{for some $\mu \in \Lambda$.}
\ee
Recall  the  isomorphism of quantum cluster Poisson varieties ${\mathscr X}_q \lra {\mathscr X}_{q^{-1}}^\circ$ given in any cluster   coordinate system by 
$\X_i \lms \X_i^{-1}$. Combining it with the argument above we conclude that 
\be \nonumber
\X_\lambda = \X'_\nu(1+ \mbox{lower order terms}). 
\ee

Since $x_\lambda$ is a standard monomial, the quantum ${\rm X}_\lambda \in \mathcal{O}_q(\mathscr{X})$. Therefore ${\bf c}_q^{\mathscr X}({\rm X}_\lambda)$ is a Laurent polynomial in the seed ${\bf s}'$. 
Combining  the last two formulas,  
we get 
\[
\X_\lambda = \X'_\mu + \mbox{terms strictly between} + \X'_\nu.
\]
 Specializing at $q=1$, and using the fact that the lowest and the top terms are standard monomials, we get 
$\mu = \nu$. The Proposition follows. 
\end{proof}
Theorem \ref{quantum.promotion.f} is a direct consequence of   Proposition \ref{quantum.promotion.lemma}. 
\end{proof} 
      
\medskip

 \section{Appendix} \la{Appen12}
 
   \medskip  

  \subsection{Amalgamation   of cluster varieties} \la{sec2}

 \medskip

We recall the amalgamation of cluster Poisson varieties  \cite{FG05}, and extend it to cluster   $K_2-$varieties. 
 
Recall  quivers from Definition \ref{DQUIV}.  A quiver ${\bf c}$ gives rise to    tori with   extra structures,    described below.

\subsubsection{The quiver Poisson torus   
$
{\mathscr X}_ {\bf c}:= {\rm Hom}(\Lambda, {\Bbb G}_m) 
$} 
It has  a  Poisson structure  
provided by the form $(\ast, \ast)$, given by 
 $
\{X_v, X_w\} = 2 (v,w) X_v X_w.
 $ 
The basis $\{e_i\}$ gives rise to {\it cluster Poisson coordinates $\{X_i\}$}.

\subsubsection{The quiver $K_2$-torus ${\mathscr A}_ {\bf c}$. }  

Take the dual lattice 
 $
\Lambda^*:= {\rm Hom}(\Lambda, \Z).
$ 
 The basis $\{e_i\}$ provides a dual basis $\{e^{\ast}_i\}$ of 
  $\Lambda^*$.
Consider the vectors 
 $\{f_i\} \in  \Lambda^*\otimes \Q$ given by 
$
f_i = d_i^{-1}e^{\ast}_i. 
$ 
Let $\Lambda^{\circ}$ be 
the sublattice spanned by the vectors $f_i$. 
The {\it quiver $K_2$-torus}  
is  
 $
{\mathscr A}_ {\bf c}:= {\rm Hom}(\Lambda^\circ, {\Bbb G}_m).
 $ 
The basis 
 $\{f_i\}$ provides {\it cluster $K_2$-coordinates $\{A_i\}$}. 
 Denote by $\Q({\mathscr A})$ the field of rational functions on a space ${\mathscr A}$, and by $\Q({\mathscr A})^*$ 
 its multiplicative group. The 2-form $\Omega_ {\bf c}:= \sum_{i,j}  d_id_j\langle e_i, e_j\rangle d\log A_i \wedge d\log A_j$ on the torus lifts to an element
 $$
  W_ {\bf c}:= \sum_{i,j}  d_id_j\langle e_i, e_j\rangle  A_i \wedge A_j \in  \Q({\mathscr A}_{ {\bf c}})^*\wedge \Q({\mathscr A}_{ {\bf c}})^*.
  $$

\begin{definition}
A mutation of a quiver $ {\bf c}$ 
in the direction of a basis vector $e_k$ is a new quiver $  {\bf c}'$.
Its lattice and the form  
are the same as of $ {\bf c}$. 
The basis $\{  e'_i\}$ of $\ {\bf c}'$ is defined by
\begin{equation} \label{12.12.04.2a}
  e'_i := 
\left\{ \begin{array}{lll} e_i + (e_{i}, e_k)_+e_k
& \mbox{ if } &  i\not = k\\
-e_k& \mbox{ if } &  i = k.\end{array}\right.
\end{equation}
The mutation transports 
the frozen basis vectors and the multipliers isomorphically. 
\end{definition}

A mutation of quivers $\mu_{e_k}:  {\bf c} \lra   {\bf c}'$ gives rise to birational isomorphisms ${\mathscr X}_{ {\bf c}} \stackrel{\sim}{\lra} {\mathscr X}_{  {\bf c'}}$ and ${\mathscr A}_{ {\bf c}} \stackrel{\sim}{\lra}  {\mathscr A}_{ {\bf c}'}$ of the cluster tori, called cluster mutations,   preserving the Poisson/$K_2$-structure:
We denote by   ${\mathscr X}_{|\mathbf c|}$ and ${\mathscr A}_{|\mathbf c|}$   the cluster varieties assigned to a quiver $\mathbf c$, see details in \cite{FG03b}.

\subsubsection{Frozen and  non-frozen variables.} For any cluster Poisson variety ${\mathscr X}$ there is a map
 \be \la{FRX}
 F_{\mathscr X}: {\mathscr X} \lra  {\mathscr X}^{\rm nf}. 
 \ee
 Here  ${\mathscr X}^{\rm nf}$ is the cluster Poisson variety provided by the non-frozen part of the quiver. Indeed,  cluster Poisson   mutations of the unfrozen variables do not depend on the   frozen ones - see (\ref{pois.clust.mut}).
 
 For any cluster $K_2-$variety ${\mathscr A}$ there is a dual kind of map on the torus given by the frozen  ${\mathscr A}-$variables:  
 \be  \la{FRA}
  F_{\mathscr A}: {\mathscr A} \lra  {\mathscr A}^{\rm fr} = {\Bbb G}_m^f. 
 \ee
 Indeed, the frozen ${\mathscr A}-$variables do not change under the mutations of the unfrozen ones.  
 
\subsubsection{Amalgamation of quivers.} 
Let $\{ {\bf c}_{s}\}$ be a collection of quivers 
parametrised by a finite set $S$:
$$
\mathbf c_{s} = \Bigl(\Lambda_{s}, (\ast, \ast)_{s}, \{e_{s, i}\}, \{f_{s, j}\}, \{d_{s, i}\}\Bigr). 
$$
 
 Denote by $V_s$ the sets parametrising the basis vectors 
of the quiver $ {\bf c}_{s}$. 

Denote by $F_s$ the subset of $V_s$ parametrizing the frozen vectors.

\bd An {\it   amalgamation data} for a collection of quivers $\{ {\bf c}_{s}\}$   is  an epimorphism 
\be \la{AMD}
\varphi: 
\coprod_{s \in S}V_s \lra K, \qquad ~\qquad~\mbox{such that:}
\ee

\vskip 1mm (i) Images of any two of the subsets may intersect only at the frozen elements. 

\vskip 1mm (ii)  The multiplier function on 
$\coprod_{s \in S}V_s$ descends to a function $d$ on $K$: 

for any   $i \in V_s$ and $i' \in V_{s'}$   mapping to the same element of   $K$ we have 
$d_i=d_{i'}$. 
\ed

Consider  the  image of the 
 frozen elements: 
 $ 
 K_0 := \varphi(\coprod_{s \in S}F_s) \subset K. 
 $  
 
The surjectivity of   $\varphi$ plus  condition i) just mean that there is an isomorphism
$$
\varphi: \coprod_{s \in S}(V_s -  F_s) \stackrel{\sim}{\lra} K - K_0. 
$$
The lattice $ \oplus_{s\in S}\Lambda_s$     
carries  a form $(\ast, \ast)':=  \oplus_{s\in S}(\ast, \ast)_s$, and 
  a basis $\{e_u\}$ where  $u\in \coprod_{s \in S}V_s$.  
\begin{definition} \label{Def3.4} 
The amalgamation of quivers  
$\mathbf c_{s}$ for the amalgamation   data (\ref{AMD}) is a
 quiver $ {\bf c}$ 
whose lattice $\Lambda$ is a sublattice of  $\oplus_{s\in S}\Lambda_s$ 
generated by the vectors 
$$
e_t:= \sum_{\varphi(u)= t}e_{u}, \qquad  t\in K.
$$
The sum is over $u\in \coprod_{s \in S}V_s$. The form $(\ast, \ast)$ on 
$\Lambda$ is induced by the form $(\ast, \ast)'$.   
Frozen vectors are parametrised by the set $K_0$. 
The multiplier function is the function $d$ from  condition (ii). 
\end{definition} 
 
Basis vectors of $\Lambda$ corresponding to   
non-frozen elements of different quivers $ {\bf c}_s$ are orthogonal for the form 
$(\ast, \ast)$.

\subsubsection{Amalgamation of cluster tori.} 
  The amalgamation data provides an embedding of lattices 
$$
\mu: \Lambda \hra \oplus_{s\in S}\Lambda_s.
$$ 
It gives rise to maps of the corresponding cluster tori:

\begin{equation} \label{2.14.05.10}
\begin{split}
&a_{\mathscr X}: \prod_{s \in S}{\mathscr X}_{\mathbf c_s} \lra
{\mathscr X}_{\mathbf c}, \qquad  
a_{\mathscr X}^*X_t = \prod_{
\varphi(u)= t}X_u.\\
&a_{\mathscr A}: {\mathscr A}_{\mathbf c} \lra \prod_{s \in S}{\mathscr A}_{\mathbf c_s}, \qquad  
a^*_{\mathscr A} A_u = A_{\varphi(u)}. 
\end{split}
\end{equation}
Here $u\in \coprod_{s \in S}V_s$. The map $a_{\mathscr X}$ is surjective, and the map $a_{\mathscr A}$ is injective. 
Since the map $\mu$ 
respects the forms $( \ast, \ast)$, the map $a_{\mathscr X}$ provides a 
similar map 
of the quantum tori.

\subsubsection{Amalgamation of cluster varieties.} 
It is easy to check that the amalgamation of quivers   commutes with quiver cluster transformations. So we arrive to the following Lemma. 

\begin{lemma} \label{9.20.04.11}  The amalgamation maps of cluster tori (\ref{2.14.05.10}) 
 commute with 
mutations, and thus gives rise to   maps of  cluster  varieties, called the   amalgamation maps:
\be \nonumber
\begin{split}
&a_{\mathscr X}: \prod_{s \in S}{\mathscr X}_{|\mathbf c_s|} \to 
{\mathscr X}_{|\mathbf c|}, \qquad \qquad a_{\mathscr A}: {\mathscr A}_{|\mathbf c|} \lra \prod_{s \in S}{\mathscr A}_{|\mathbf c_s|}.\\
\end{split}
\ee
The map   $a_{\mathscr X}$ is a surjective map of cluster Poisson varieties. The map $a_{\mathscr A}$ is an injective  map of  cluster $K_2$-varieties.
  There is a  quantum analog of the Poisson amalgamation map given by  a map of  
quantum 
cluster  varieties, understood as  an injective map of algebras
\be \nonumber
\begin{split}
&a_{{\mathscr X}_q}^*: {\cal O}_q({\mathscr X}_{|\mathbf c|}) \lra \bigotimes_{s \in S} {\cal O}_q( {\mathscr X}_{|\mathbf c_s|}).\\
\end{split}
\ee   \end{lemma}

\subsubsection{Defrosting.} 
Let 
$K_0' \subset K_0$. Assume 
that the restriction of the form $(\ast, \ast)$  to 
$K- K'_0$ 
takes values in $\Z$. 
Then there is a new quiver  $\mathbf c':= (K, K_0', \varepsilon, d)$. We say that 
{\it the quiver  $\mathbf c'$ 
is obtained from $\mathbf c$ by defrosting of  $K_0-K_0'$}. 
Abusing notation, one may refer to amalgamation followed by defrosting 
 simply as amalgamation. 
However then the defrosted subset must be specified.

\medskip

  \subsection{Quasi-cluster transformation} \la{sec.quasicl}
  
  \medskip
  
    Let $\widehat{\varepsilon}=(\widehat{\varepsilon}_{ij})$ be an  $n \times n$ skewsymmetric matrix over $\mathbb{Q}$ and let $d=\{d_{1}, \ldots, d_m\}$ be a set of positive 
   numbers such that $m\leq n$ and 
   \begin{itemize}
\item $\varepsilon_{ik}:= \widehat{\varepsilon}_{ik} d_k\in \mathbb{Z}$ for $1\leq i\leq n$ and $1\leq k \leq m$.
\end{itemize}
One assigns to $\widehat{\varepsilon}$  a triple $\left( {\rm T},\alpha,  \omega \right)$, where ${\rm T}$ is a split algebraic torus,  $\alpha=\{X_1, \ldots, X_n\}$   a basis of the character lattice ${\rm Hom}({\rm T}, \mathbb{G}_m)$, and $\omega$   a Poisson bivector:
$
\omega= \sum_{i,j=1}^n \widehat{\varepsilon}_{ij} X_i  \partial_{X_i}\wedge X_j \partial_{ X_j}. 
$


The Poisson bivector $\omega$  determines a Poisson bracket   $\{*,*\}$ such that $\{X_i, X_j\}= 2 \widehat{\varepsilon}_{ij} X_i X_j$.  

The datum ${\bf s}=\{{\rm T}, \alpha, \omega, d \}$ is called a {\it seed}. 
The {\it exchange matrix} $\varepsilon$ is determined by $\omega$ and $d$.
 
 \vskip 2mm \noindent
Let $\tau$ be a permutation of $\{1,\ldots, m\}$. A seed ${\bf s}'=\{{\rm T}', \alpha', \omega', d' \}$ is \emph{$\tau$-linearly equivalent} to ${\bf s}$ if 
\begin{itemize}
\item ${\rm T}'={\rm T}$, $\omega'=\omega$, and $d=\{d_{\tau(i)}\}$; 
\vskip 1mm\item $\alpha'=\{X_i'\}$  is another basis of ${\rm Hom}({\rm T}, \mathbb{G}_m)$ such that  $X_k'=X_{\tau{(k)}}$ for $1\leq k\leq m$.
\end{itemize}
In other words, $X_k'=x_{\tau{(k)}}$ $1 \leq k \leq m$, while $X_{m+i}$ are 
any monomials of coordinates $X_j'$, $1 \leq j \leq n$,  such  that the induced monomial transformation 
is Poisson.

\subsubsection{Cluster Poisson mutation.} 
Let ${\rm T}'$ be a  split  rank $n$ torus and  $\alpha'=\{X_i'\}$   a basis of its character lattice. For $ k\in\{1,\ldots, m\}$,   consider a birational map
$\mu_k:{\rm T}\rightarrow {\rm T}'$ given by
\be 
\label{pois.clust.mut}
\mu_k^\ast X_i'=   \left\{ \begin{array}{ll} 
      X_{k}^{-1} & \mbox{if } i=k \\
      X_i\left(1+X_k^{-{\rm sgn}(\varepsilon_{ik})}\right)^{-\varepsilon_{ik}} & \mbox{if } i\neq k \\
   \end{array}\right.
\ee
Let $\omega'=(\mu_k)_\ast \omega$ and $d'=d$. The seed $\mu_k({\bf s}):= \{{\rm T}', \alpha',\omega', d' \}$ is called a mutated  at the direction $k$ seed.
An easy calculation shows that

 \begin{lemma} 
\label{linear.lemm} 
If  seeds ${\bf s}_1$ and ${\bf s}_2$ are $\tau$-linearly equivalent, then $\mu_k({\bf s})$ and $\mu_{\tau(k)}({\bf s})$ are $\tau$-linearly equivalent.
 \end{lemma}
 
Given  a seed ${\bf s}$, let  $|{\bf s}|$ be the collection of seeds obtained from ${\bf s}$ by a finite number of mutations and basis changes.
The {\it cluster Poisson variety} $\mathscr{X}_{|{\bf s}|}$ is defined by gluing the tori in  all seeds of $|{\bf s}|$ via the corresponding sequences of  birational maps \eqref{pois.clust.mut} and basis change maps. 
We   often skip the subscript $|{\bf s}|$,  writing $\mathscr{X}$  for $\mathscr{X}_{|{\bf s}|}$. The bivector $\omega$ defines a   Poisson bracket $\{\ast, \ast\}$ on   $\mathscr{X}$.  

By the very definition,  every ${\bf s}' \in |{\bf s}|$ induces  an open embedding 
$$
\chi_{{\bf s}'}: \mathbb{G}_m^n \hra \mathscr{X},
$$
called a {\it quasi-cluster Poisson chart}. The atlas $\mathcal{C}_{\mathscr{X}}$ consists of all such charts.
If ${\bf s}'$ and ${\bf s}''$ are $\tau$-linear equivalent, then the images of $\chi_{{\bf s}'}$ and $\chi_{{\bf s}''}$ coincide.

\begin{definition}
A {\it quasi-cluster automorphism} $\phi$ of $\mathscr{X}$ is a biregular isomorphism of $\mathscr{X}$ which 
 
\vskip 1mm1.  Preserves the Poisson bracket $\{\ast, \ast \}$;

\vskip 1mm2.  Preserves the atlas $\mathcal{C}_{|{\bf s}|}$.  
 
The quasi-cluster modular group $\mathcal{G}_{\mathscr{X}}$ consists of quasi-cluster automorphims of $\mathscr{X}$.
\end{definition}
Thanks to Condition 1, if $\phi\circ \chi\in \mathcal{C}_{{\mathscr{X}}}$ for one $\chi\in \mathcal{C}_{\mathscr{X}}$, then Condition 2 holds.

The group $\mathcal{G}_{\mathscr{X}}$ acts on the ring $\mathcal{O}(\mathscr{X})$ of regular functions of $\mathscr{X}$. 
Let $\mathcal{O}_q(\mathscr{X})$ be the quantization of $\mathcal{O}(\mathscr{X})$ introduced in \cite{FG03b}. The action of $\mathcal{G}_{\mathscr{X}}$  on $\mathcal{O}(\mathscr{X})$ can be lifted to an action on $\mathcal{O}_q(\mathscr{X})$.

  \medskip

\subsection{The Casimir Center} \la{SECT20.4}

\medskip

Recall that every quiver ${\bf c}$ contains a lattice $\Lambda_{\bf c}$ with a basis $\{e_i\}$, together with a skewsymmetric bilinear form $(\ast, \ast): \Lambda_{\bf c}\times \Lambda_{\bf c}\rightarrow \mathbb{Q}$, see Definition \ref{DQUIV}. Consider the sub-lattice
\[
{\rm Ker}\,\Lambda_{\bf c}:= \left\{ v\in \Lambda_{\bf c} ~\middle|~ ( v, w) =0, ~\forall w \in \Lambda_{\bf c} \right\}.
\]
Every $v\in {\rm Ker}\,\Lambda_{\bf c}$ gives rise to a {\it Casimir monomial} $X_{v}$ in the quantum torus algebra $\mathcal{O}_q({\rm T}_{\Lambda_{\bf c}})$. 
Note that a Casimir monomial for one quiver remains a Casimir monomial for every quiver in the mutation equivalence class. We define the {\it Casimir torus} ${\rm C}_{\mathscr{X}}$  by setting
$$
{\rm C}_{\mathscr{X}}:= {\rm Hom}({\rm Ker}\ \Lambda_{\bf c}, {\Bbb G}_m).
$$ 

\bd The Casimir center ${\rm Cas}(\mathscr{X})$ is an algebra given by the linear span of the Casimir monomials. 
\ed

Recall that the quantized  algebra of functions on a cluster Poisson variety ${\mathscr X}$
\[
\mathcal{O}_q(\mathscr{X})=\bigcap_{{\bf c}'} \mathcal{O}_q({\rm T}_{\Lambda_{{\bf c}'}}),
\]
where the intersection is over all the quivers ${\bf c}'$ that are mutation equivalent to ${\bf c}$. The Casimir center ${\rm Cas}(\mathscr{X})$ is a subalgebra of $\mathcal{O}_q(\mathscr{X})$. 
The lattice $\Lambda_{\bf c}$ is equipped with a basis $\{e_i\}$. Let 
\[\Lambda_{\bf c}^+:=\bigoplus \mathbb{Z}_{\geq 0} e_i\]
 be the positive cone inside $\Lambda_{\bf c}$. For $u, w\in \Lambda_{\bf c}$, we say that $u>w$ if $u\neq w$ and $u-w\in \Lambda_{\bf c}^+$. 
 
 Let $n=\dim \Lambda_{\bf c}$.
 We say that $\mathcal{O}_q(\mathscr{X})$ is {\it large enough} if there exist $n$ many linearly independent vectors $w_1, \ldots, w_n \in \Lambda_{\bf c}$, which  give rise to the following functions  in $ \mathcal{O}_q(\mathscr{X})$
 \[
 \mathbb{I}(w_i) = X_{w_i}+ \sum_{u> w_i} c_{\bf c}(u, w_i) X_{u}.
 \]
 Note that the quantum Duality Conjecture implies that $\mathcal{O}_q(\mathscr{X})$ is large enough.

\bp Assume that $\mathcal{O}_q(\mathscr{X})$ is large enough. When $q$ is not a root of unity, then ${\rm Cas}(\mathscr{X})$ is the center of $\mathcal{O}_q(\mathscr{X})$.
\ep

\begin{proof} By the definition of the quantum torus algebra $\mathcal{O}_q({\rm T}_{\Lambda_{\bf c}})$, the algebra ${\rm Cas}(\mathscr{X})$ is the center of  $\mathcal{O}_q({\rm T}_{\Lambda_{\bf c}})$. Therefore, ${\rm Cas}(\mathscr{X})$ belongs to the center of $\mathcal{O}_q(\mathscr{X})$. 
Conversely, let $f$ be a central element of $\mathcal{O}_q(\mathscr{X})$. Let us fix a quiver ${\bf c}$. We have the decomposition $f=g+h,$
where $g\in {\rm Cas}(\mathscr{X})$, and
\[
h=\sum_{v\notin {\rm Ker} \Lambda_{\bf c}} c(v, h) X_v.
\]
Note that $h$ is in the center of $\mathcal{O}_q(\mathscr{X})$. 
Suppose that $h\neq 0$. There exists a minimal vector $v$ such that the coefficient $c(v,h)\neq 0$ and $c(u,h)=0$ whenever $v>u$.
Assume that $\mathcal{O}_q(\mathscr{X})$ is large enough. We have 
\[
\mathbb{I}(w_i) \cdot h = h \cdot \mathbb{I}(w_i).
\]
By comparing the coefficients of the term $X_{w_i+v}$, we get 
\[
q^{(w_i, v)}c(v,h)= q^{(v, w_i)}c(v,h).
\]
If $q$ is not a root of unity, then $(v, w_i)=0$ for $1\leq i \leq n$. Hence, $v\in {\rm Ker} \Lambda_{\bf c}$, which is a contradiction. Therefore $h=0$ and $f\in {\rm Cas}(\mathscr{X})$.
\end{proof}

\medskip

\subsubsection{The cluster $p$ map}

For every $t\in {\rm H}$, we set $t^\ast:= w_0(t^{-1})$. Let $e\in {\rm H}$ be the identity.
\begin{lemma} 
Let $(\A_1, \A_2)$ be a pair of decorated flags of generic position. 
There exists a $t\in {\rm H}$ such that 
\be
\label{hcond=e}
h(\A_1\cdot t, \A_2 \cdot t^\ast)=e.
\ee
\end{lemma}

\begin{proof}
Note that 
\[
h(\A_1\cdot t, \A_2\cdot t^\ast)= t^2 \cdot h(\A_1, \A_2).
\]
We take $t$  such that $t^2=h(\A_1, \A_2)^{-1}$.
\end{proof}
Note that the $t$ in \eqref{hcond=e} is not unique. However, if $h(\A_1, \A_2)\in {\rm H}(\mathbb{R}_{>0})$, then there is a unique $t\in {\rm H}(\mathbb{R}_{>0})$ that satisfies the condition \eqref{hcond=e}. In this case, we write
\be
\label{average.t}
t=h(\A_1, \A_2)^{-\frac{1}{2}},
\ee
and we denote the corresponding pinning by 
\be
\overline{p}(\A_1, \A_2):= (\A_1\cdot t, \A_2 \cdot t^\ast).
\ee

Let $\G'$ be the universal cover of $\G$.
Recall the cluster $p$-map
\[
p: \mathscr{A}_{\G', \bS} \longrightarrow \mathscr{P}_{\G, \bS}
\]
such that $p^\ast X_i =\prod_j A_j^{\varepsilon_{ij}}$. 
When $\bS$ contains marked points, $\varepsilon_{ij}$ are not necessarily integers. Therefore $p$ may  not be a rational map. However, $p$ is well-defined when restricted to the totally positive parts:
\be
\label{p-positive-map}
p: \mathscr{A}_{\G', \bS}(\mathbb{R}_{>0})  \longrightarrow \mathscr{P}_{\G, \bS} (\mathbb{R}_{>0}).
\ee
A modified rational map $\pi_{\bS}$ has been investigated in Section \ref{SEC9.2x}. 

\bp The map \eqref{p-positive-map} is defined in the following way
\begin{itemize}
\item For every puncture $s$ of $\bS$, the map $p$ takes its associated decorated flag $\A_s$ to $\B_s=\pi(\A_s)$.
\item For every boundary interval $I$ of $\bS$, the map $p$ assigns to $I$ the pinning
\[
\overline{p}_{I}=\overline{p}(\A_1, \A_2)
\]
as in the following picture. 

\begin{center}
\begin{tikzpicture}
\draw[thick] (-135:1)--(0,0)--(2,0)--+(-45:1);
\node at (0,0) {$\bullet$};
\node at (2,0) {$\bullet$};
\node at (0,0.3) {\small $\A_1$};
\node at (2, 0.3) {\small $\A_2$};
\node at (1, -.3) {$I$};
\draw[-latex, thick] (4,0)--(6,0);
\end{tikzpicture}
\hskip 1cm
\begin{tikzpicture}
\draw[thick] (-135:1)--(0,0)--(2,0)--+(-45:1);
\node at (0,0) {$\bullet$};
\node at (2,0) {$\bullet$};
\node at (1,0.4) {$\overline{p}_I=\overline{p}(\A_1, \A_2)$};
\end{tikzpicture}
\end{center}
\end{itemize}
\ep
\begin{proof} The proof is similar to the proof of Theorem \ref{11.24.18.1}. 
Indeed, the map $\pi_{\bS}$ in \eqref{CORT1X.general.nvadfo} associates to each boundary interval the pinning determined by the first decorated flag $\A_1$, while the map $p$ uses the ``normalized" pinning $\overline{p}(\A_1, \A_2)$, adjusted by $t$  as in \eqref{average.t}. A direct calculation shows that the adjustment by $t$ cancels the term $\mathbb{A}_k^{\frac{1}{2}}$ as in \eqref{CORT1X.general}, which concludes the proof of this Proposition.
\end{proof}

\medskip
\subsubsection{The monodromy map $\mu_\bS$.}

 Recall the map $\mu_{\bS}$ in Theorem \ref{Th1.14}, part (6).  
We restrict it to the totally positive part
\[
\mu_{\bS}:= (\mu_p, \mu^\pi_{{\rm out}}): \mathscr{P}_{\G, \bS}(\mathbb{R}_{>0}) \longrightarrow {\rm C}_{\G, \bS}(\mathbb{R}_{>0}).
\]
We prove that $\mu_{\bS}$ is a symplectic fibration. 
Let us consider the fiber of $\mu_{\bS}$ over the identity $e$:
\[
\mathscr{P}_{\G, \bS}(\mathbb{R}_{>0})_e: = \mu_{\bS}^{-1}(e).
\]

\bt 
\label{image.of.p.24}
The image of $p$ in \eqref{p-positive-map} coincides with the fiber $\mathscr{P}_{\G, \bS}(\mathbb{R}_{>0})_e$.
\et
\begin{proof}
As illustrated in the following figure, under the map $p$, the decorated flag ${\rm A}_i$ is rescaled by $t_i$ and $t_{i-1}^\ast$, giving rise to the pinnings to its neighbored boundary intervals.
\begin{center}
\begin{tikzpicture}[scale=1.5]
\draw[thick] (-150:2)--(0,0)--(-30:2);
\node at (0,0) {$\bullet$};
\node at (0, 0.3) {$\A_i$};
\node[red] at (-.5,0.2) {$\A_i \cdot t_i$};
\node[red] at (.6,0.2) {$\A_i \cdot t_{i-1}^\ast$};
\node[red] at (0, -.3) {$\rho_i$};
\draw[latex-latex, red] (-1.7,-0.8) --+ (30:1.5);
\draw[latex-latex, red] (1.7,-0.8) --+ (150:1.5);
\end{tikzpicture}
\end{center}
Here $\rho_i=t_{i-1}^\ast t_i^{-1}$. Let $\pi$ be a boundary circle with $d$ many marked points. If $d$ is even, then
\[
\mu^\pi_{\rm out}= \rho_1 \cdot \rho_2^\ast\cdot \ldots \cdot \rho_d^\ast = (t_d^{\ast} t_1^{-1}) \cdot (t_1^\ast t_2^{-1})^\ast \cdots (t_{d-1}^\ast t_d^{-1})^\ast =e \in {\rm H}.
\]
If $d$ is odd, then 
\[
\mu^\pi_{\rm out}=\left(\rho_1 \cdot \rho_2^\ast\cdot \ldots \cdot \rho_d\right)_\ast =\left( (t_d^{\ast} t_1^{-1}) \cdot (t_1^\ast t_2^{-1})^\ast \cdots (t_{d-1}^\ast t_d^{-1}) \right)_\ast = (t_d^\ast t_d^{-1})_\ast=e \in {\rm H}_\ast.
\]
Meanwhile, the monodromy of $\mathscr{A}_{{\rm G}', \bS}$ surrounding each puncture is unipotent. Under the map $p$, we get $\mu=e$. To summarize, we get 
\[
p\left(\mathscr{A}_{\G', \bS}(\mathbb{R}_{>0})\right) \subset \mathscr{P}_{\G, \bS}(\mathbb{R}_{>0})_{e}.
\]

Now we prove that $\mathscr{P}_{\G, \bS}(\mathbb{R}_{>0})_e \subset p(\mathscr{A}_{\G', \bS}(\mathbb{R}_{>0}))$. Pick a local system in $\mathscr{P}_{\G, \bS}(\mathbb{R}_{>0})_e$, whose data associated with a marked point $i$ is illustrated by the following picture. We shall assign to $i$ a decorated flag ${\rm A}_i$.
\begin{center}
\begin{tikzpicture}[scale=1.5]
\draw[thick] (-150:2)--(0,0)--(-30:2);
\node at (0,0) {$\bullet$};
\node[red] at (-1,-.3) {$p_i$};
\node[red] at (1, -.3) {$p_{i-1}$};
\node[red] at (-.5,0.2) {$\A_i '$};
\node[red] at (.6,0.2) {$\A_i ''$};
\node[red] at (0, -.3) {$\rho_i$};
\draw[latex-latex, red] (-1.7,-0.8) --+ (30:1.5);
\draw[latex-latex, red] (1.7,-0.8) --+ (150:1.5);
\end{tikzpicture}
\end{center}
Note that $\rho_i\in {\rm H}(\mathbb{R}_{>0})$.
Let $d$ be the number of marked points on a boundary circle $\pi$. If $d$ is even, then
\[
\rho_1\rho_2^\ast \cdots \rho_d^\ast =e.
\]
 There exists $t_1, t_2, \ldots, t_n \in {\rm H}(\mathbb{R}_{>0})$ such that 
\[
\rho_i=t_i^{-1}\cdot t_{i-1}^\ast.
\]
If $d$ is odd, then we have
\[
(\rho_1\rho_2^\ast \cdots \rho_d)^\ast = (\rho_1\rho_2^\ast \cdots \rho_d)^{-1}
\]
Let us take $t_d\in {\rm H}(\mathbb{R}_{>0})$ such that 
$
t_d^2= (\rho_1\rho_2^\ast \cdots \rho_d)^{-1}$. Note that $t_d^\ast =t_d^{-1}$.
We obtain $t_i$ recursively via $\rho_i=t_i^{-1}t_{i-1}^\ast$. 
Let us assign to $i$ the decorated flag 
\[ \A_i= \A_{i}'\cdot t_i^{-1} = \A_i''\cdot (t_{i-1}^\ast)^{-1} \in {\rm G} /{\rm U}.
\]
Similarly, for every puncture $s$ of $\bS$, the monodromy of every local system in  $\mathscr{P}_{\G, \bS}(\mathbb{R}_{>0})_e$ is unipotent. Hence we may promote the associated flag ${\rm B}_s$ to a decorated $\A_s$ such that every partial potential $\mathcal{W}_{s, i}$ is positive. 
In this way, we lift an element from $\mathscr{P}_{\G, \bS}(\mathbb{R}_{>0})_e$ to $\mathscr{A}_{\G, \bS}(\mathbb{R}_{>0})$.  
By the constructions of positive structures on $\mathscr{A}$-spaces as in \cite{FG03a},  $\mathscr{A}_{\G, \bS}(\mathbb{R}_{>0})$ coincides with  $\mathscr{A}_{\G', \bS}(\mathbb{R}_{>0})$, which concludes the proof of the Theorem.
\end{proof}

\smallskip

\subsubsection{Proof of Theorem \ref{TH1.16a}, part 1.} \la{20.3.3p}

\bl 
\label{sur.M-maps}
The map $\mu_{\bS}$ in \eqref{mubs-tcas} is surjective. 
\el
\begin{proof}
Let ${\bf H}_{\bS}$ be the product of ${\rm H}$ over all the punctures and special points of $\bS$. Gathering the maps $\mu_p$ in \eqref{mup} and  $\rho_s$ in  \eqref{sspm}, we define the map
\be
\label{MbS.maps}
M_{\bS}=\left(\{\mu_{p}\}, \{\rho_s\}\right):~ \mathscr{P}_{\G, \bS}\longrightarrow {\bf H}_{\bS}.
\ee
Note that the surjectivity of $\mu_{\bS}$ follows from the surjectivity of $M_{\bS}$.

The map $M_{\bS}$ is compatible with the gluing map $\gamma_{{\rm I}_1, {\rm I}_2}$ in \eqref{AMAX}, i.e., the following  diagram commutes:
\[ \begin{tikzcd}
\mathscr{P}_{\G,\bS} \arrow{r}{\gamma_{{\rm I}_1, {\rm I}_2}} \arrow[swap]{d}{M_{\bS}} & \mathscr{P}_{\G, \bS'} \arrow{d}{M_{\bS'}} \\%
{\bf H}_{\bS} \arrow{r}{\tilde{\gamma}_{{\rm I}_1,{\rm I}_2}}& {\bf H}_{\bS'}
\end{tikzcd},
\]
where  $\tilde{\gamma}_{{\rm I}_1,{\rm I}_2}$ is defined by taking the products of the corresponding Cartan elements that are merged together, and is therefore surjective. Hence the surjectivity of $M_{\bS}$ implies the surjectivity of $M_{\bS'}$. Without loss of generality, we assume that $\bS$ is a disjoint union of triangles. 

It suffices to prove the surjectivity of $M_t$ for a triangle $t$.
Let ${\rm U}^\ast= {\rm U}\cap \B^-w_0\B^-$. 
Let $u_1\in {\rm U}^\ast$ be the unipotent invariant associated with the top angle of $t$. Recall the isomorphism 
\be
\label{iso.p.uhh}
\mathscr{P}_{\G, t}\stackrel{\sim}{\longrightarrow} {\rm H}\times {\rm H} \times {\rm U}^\ast,
\ee
which takes the value $(\rho_1, \rho_2, u_1)$. 
\begin{center}
\begin{tikzpicture}
\draw[thick] (0,0)--(60:2)--(2,0)--(0,0);
\node[red] at (0.4, 0.2) {$\rho_3$};
\node[red] at (1.6, 0.2) {$\rho_2$};
\node[red] at (1, 1.4) {$\rho_1$};
\end{tikzpicture}
\end{center}
Let ${\rm T}$ be the Cartan subgroup of the universal cover $\G'$. Let $\pi$ be the projection from ${\rm T}$ to ${\rm H}$.
Recall the positive map $\beta$ from ${\rm U}^\ast$ to {\rm T} in (161) of \cite{GS13}. By Lemma 6.2 (1) and Lemma 6.4 (4) in {\it loc.cit.}, we have
\be
\label{betta,ms}
\pi\left(\beta(u_1)\right)= \rho_1^{-1}\rho_3^{-1} \rho_2^{*}.
\ee
 The map $M_t$ can be factorized as 
\[
\mathscr{P}_{\G, t}\stackrel{\sim}{\longrightarrow} {\rm H}\times {\rm H} \times {\rm U}^\ast \longrightarrow {\rm H}^3,
\qquad 
(\rho_1, \rho_2, u_1) \longmapsto (\rho_1, \rho_2, \pi(\beta(u_1)^{-1})\rho_1^{-1}\rho_2^\ast).
\]
Since the $\beta$ map is surjective, the map $M_t$ is also surjective, which concludes the proof of the Lemma.  
\end{proof}

We now complete the proof of Theorem \ref{TH1.16a}, part 1, and hence of Theorem \ref{TH1.16a}. Let $X^\ast({\rm C}_{\G, \bS})$ be the character lattice of ${\rm C}_{\G, \bS}$. Since every character of  ${\rm C}_{\G, \bS}$ is Casimir element of  $\mathscr{P}_{\G, \bS}$, we obtain a homomorphism from $X^\ast({\rm C}_{\G, \bS})$ to the Casimir lattice $X^\ast({\rm C}_{\mathscr{P}_{\G, \bS}})$.
The surjectivity of  $\mu_{\bS}$ implies that this homomorphism is injective, and therefore the map $i_{\rm C}$ in  \eqref{ISOIC}  is surjective.  By Theorem \ref{image.of.p.24} and the definition of the cluster $p$ map, the ranks of $X^\ast({\rm C}_{\G, \bS})$ and $X^\ast({\rm C}_{\mathscr{P}_{\G, \bS}})$ are both equal to the corank of the cluster exchange matrix of $\mathscr{P}_{\G, \bS}$. Therefore, $X^\ast({\rm C}_{\G, \bS})$ is of finite index in $X^\ast({\rm C}_{\mathscr{P}_{\G, \bS}})$, and the isomorphism \eqref{ISORPP} follows. 
As we explained right after the statement of Theorem \ref{TH1.16a}, the  proof of  Theorem \ref{TH1.16a} follows from this. 
\vskip 2mm
Therefore  Theorem \ref{TH1.16a} is  proved.

\smallskip

\subsubsection{Extra Casimir elements of $\mathscr{P}_{\G, \bS}$.} \la{20.3.4}

Below we assume that $\bS$ is connected and fully colored.
We construct extra Casimir elements of $\mathscr{P}_{\G, \bS}$ which do not  necessarily belong to $\mu_\bS^*\mathcal{O}({\rm C}_{\G, \bS})$.
Assume that each boundary interval of $\bS$ is colored. We consider the map
\be
\label{rhos.global}
\rho_{\bS}=  \prod_{p} \mu_p \times \prod_{s} \rho_s : \, \mathscr{P}_{\G, \bS} \longrightarrow {\rm H}, 
\ee
where the products are taken over all  punctures $p$ and all  special points $s$. 

Below we use the notations in the proof of Lemma \ref{sur.M-maps}.
\bl
The map $\rho_{\bS}$ can be naturally lifted to a regular map
\[
\tilde{\rho}_{\bS}: \mathscr{P}_{\G, \bS} \longrightarrow {\rm T}
\]
such that $\rho_{\bS}=\pi\circ \tilde{\rho}_{\bS}$.
\el
\begin{proof} 
Note  that \eqref{rhos.global} is compatible with the gluing map $\gamma_{{\rm I}_1, {\rm I}_2}$:
$\rho_{\bS}=\rho_{\bS'}\circ \gamma_{{\rm I}_1, {\rm I}_2}. $
Therefore it suffices to prove the Lemma for the case when $\bS=t$ is a triangle. 
By \eqref{betta,ms}, we get
\[
\rho_t=\rho_1\rho_2\rho_3= (\rho_2\rho_2^*) \cdot \pi\left(\beta(u_1)^{-1}\right).
\]
Note that $\rho_2\rho_2^*=\rho_2w_0(\rho_2^{-1})$ can be naturally lifted to an element $\widetilde{\rho_2\rho_2^*}$ in ${\rm T}$. Hence we get
\[
\tilde{\rho}_t= \widetilde{\rho_2\rho_2^*} \cdot \beta(u_1)^{-1} \in {\rm T}. \qedhere
\]
\end{proof}

As a consequence, every character $\Lambda\in X^*({\rm T})$ gives rise to a {\it positive} regular function
\be
\Lambda(\tilde{\rho}_{\bS}): ~\mathscr{P}_{\G, \bS} \lra \mathbb{G}_m.
\ee
\begin{example} Let $\G$ be of type $A_{m-1}$ and let $\Lambda_1$ be its first fundamental weight. We have
\[
m \Lambda_1 =\sum_{k=1}^{m-1} k \alpha_{m-k},
\]
where $\alpha_{m-k}$ are the simple positive roots. 
Then 
\[
\Lambda_1(\tilde{\rho}_t)^m= \prod_{k=1}^{m=1} \alpha_{m-k}(\rho_1\rho_2\rho_3)^k.
\]
Recall the cluster Poisson structure on $\mathscr{P}_{{\rm PGL}_m, t}$ in Section \ref{SSECC3.2.2}.
Note that the cluster variable $X_{a_1, a_2, a_3} $ only appears once in $\alpha_{m-a_1}(\rho_1)$, $\alpha_{m-a_2}(\rho_2)$, and $\alpha_{m-a_3}(\rho_3)$. Therefore   
\[
 \prod_{k=1}^{m=1} \alpha_{m-k}(\rho_1\rho_2\rho_3)^k =\left (\prod_{(a_1, a_2, a_3)\in {\rm I}^{(m)}} X_{a_1, a_2, a_3}\right)^m.
\]
The function $\Lambda_1(\tilde{\rho}_t)$ is positive. Therefore 
\be
\Lambda_1(\tilde{\rho}_t) = \prod_{(a_1, a_2, a_3)\in {\rm I}^{(m)}} X_{a_1, a_2, a_3}.
\ee
For a general surface $\bS$, we use a cluster Poisson chart $\{X_f\}_{f\in I}$ of $\mathscr{P}_{{\rm PGL}_m, \bS}$ defined in Section \ref{SSECC3.2.3}. By the gluing property of $\rho_{\bS}$, we get
\be
\label{mon.tsist}
\Lambda_1(\tilde{\rho}_{\bS}) = \prod_{f\in I} X_f.
\ee
\end{example}

If $\bS$ has only punctures, the function $\Lambda(\tilde{\rho}_{\bS})$ is a Casimir element. For instance, in Section   \ref{SSECC3.2.3}, we observe that  every quiver of $\mathscr{P}_{{\rm PGL}_m, \bS}$  is balanced, meaning each vertex has an equal number of incoming and outgoing arrows. Therefore $\Lambda_1(\tilde{\rho}_{\bS}) = \prod_{f\in I} X_f$ is a Casimir function. Furthermore, $\Lambda_1(\tilde{\rho}_{\bS})$ cannot be obtained by pullback from $\mathcal{O}({\rm C}_{\G, \bS})$. 

If $\bS$ has boundary intervals, then for any weight $\Lambda \in X^*({\rm T})$ with $\Lambda=-w_0(\Lambda)$, the function $\Lambda(\tilde{\rho}_{\bS})$ is a Casimir element.

\medskip

Let $X^\ast({\rm H})_{\mathbb{R}}:=X^\ast({\rm H})\otimes \mathbb{R}$. Let ${\varpi}_{\bS}=\left\{ \omega_p, \omega_s\right\}$ be an assignment of $\omega_p\in X^\ast({\rm H})_{\mathbb{R}}$ to every puncture $p$ of $\bS$ and an assignment of $\omega_{s}$ to every special point.  Composing with \eqref{MbS.maps}, we define the following map on the positive part of $\mathscr{P}_{\G, \bS}$:
\[
\tau_{\varpi_{\bS}}= \prod_p \omega_p \cdot \prod_s \omega_s(\rho_s):~ \mathscr{P}_{\G, \bS} (\mathbb{R}_{>0}) \longrightarrow \mathbb{R}_{>0}.
\]
For $c\in \mathbb{R}$, the floor function $\lfloor c \rfloor$ is the greatest integer less than or equal to $c$. Let $\alpha_1, \ldots, \alpha_r \in X^\ast({\rm H})$ be simple positive roots. For $\omega=\sum_i c_i \alpha_i \in X^\ast({\rm H})_{\mathbb{R}}$, we set
\[
\lfloor \omega \rfloor :=\sum_{i} \lfloor c_i \rfloor \alpha_i, \qquad \overline{\omega}:= \omega-\lfloor \omega \rfloor.
\]
Similarly, we set 
\[
\lfloor \varpi_{\bS} \rfloor:=\left\{\lfloor \omega_p\rfloor, \lfloor \omega_s\rfloor\right\}, \qquad  \overline{\varpi}_{\bS}:=\left\{\overline{\omega}_p, \overline{\omega}_s\right\},
\]
Note that 
\[
\tau_{\varpi_{\bS}}= \tau_{\lfloor \varpi_{\bS}\rfloor}\cdot \tau_{\overline{\varpi}_{\bS}}.
\]

\bl 
\label{regularity.lemma.center}
The function $\tau_{\varpi_{\bS}}$ can be extended to a regular function on $\mathscr{P}_{\G, \bS}$ if and only if there exists a weight $\Lambda\in X^\ast ({\rm T})$ such that \[
\tau_{\overline{\varpi}_{\bS}}=\Lambda(\tilde{\rho}_{\bS}).\]
\el

\begin{proof} The ``if part" follows since both $\tau_{\lfloor \varpi_{\bS}\rfloor}$ and $\Lambda(\tilde{\rho}_{\bS})$ are regular functions on $\mathscr{P}_{\G, \bS}$.

For the ``only if part", we first prove the case when $\bS=t$. Let $\varpi_t=\{\omega_1, \omega_2, \omega_3\}$ be an assignment of vectors to the vertices of $t$. Let $(\rho_1, \rho_2, u_1)\in {\rm H}\times {\rm H}\times {\rm U}^\ast$. By the isomorphisms \eqref{iso.p.uhh} and \eqref{betta,ms}, we get
\be
\label{taumaopit}
\tau_{\varpi_t}= \omega_1(\rho_1)\omega_2(\rho_2)\omega_3(\rho_3)=  (\omega_1-\omega_3)(\rho_1) \cdot (\omega_2+\omega_3^\ast)(\rho_2)\cdot \omega_3(\beta(u_1)^{-1}) .
\ee
Recall Lusztig's coordinates of ${\rm U}^\ast$ associated with a reduced word ${\bf i}=(i_1, \ldots, i_n)$ of $w_0$:
\[
\mathbb{G}_m^n \longrightarrow {\rm U}^\ast, \hskip 10mm u_1=x_{i_1}(a_1)x_{i_2}(a_2)\cdots x_{i_n}(a_n). 
\]
By Lemma 5.3 of \cite{GS13}, we get 
\be
\label{u1beta,ap}
\beta(u_1)=\prod_{k=1}^n \beta_{k}^{\bf i}(a_k^{-1}),
\ee
where $\beta_k^{\bf i}$ are the positive coroots of $\G$ defined in \eqref{reduced.word.positive.roots}. Plugging \eqref{u1beta,ap} into \eqref{taumaopit}, we get 
\[
\tau_{\varpi_t}= (\omega_1-\omega_3)(\rho_1) \cdot (\omega_2+\omega_3^\ast)(\rho_2)\cdot \prod_{k=1}^n a_k^{\langle \beta_{k}^{\bf i}, \omega_3\rangle }.
\]
The regularity of $\tau_{\varpi_t}$ implies that  $\omega_1-\omega_3$ and $\omega_2+\omega_3^\ast$ are in $X^\ast ({\rm H})$ and $\omega_3$ is in $X^\ast({\rm T})$. 
Or equivalently, 
\[
\overline{\omega}_1= \overline{\omega}_2= \overline{\omega}_3:=\Lambda \in X^\ast({\rm T}).\] Therefore $\tau_{\overline{\varpi}_{t}}=\Lambda(\tilde{\rho}_{t}).$

Now let $\bS$ be a connected fully colored surface. Let ${\cal T}$ be an ideal triangulation of $\bS$. Let $\varpi_t$ be the restriction of  $\varpi_{\bS}$  to every ideal triangle $t\in {\cal T}$. The map $\tau_{\varpi_{\bS}}$ is compatible with the gluing map. Therefore, the following diagram commutes
\[ \begin{tikzcd}
\prod_{t\in \mathcal{T}}\mathscr{P}_{\G,t} \arrow{r}{\rm gluing} \arrow[swap]{dr}{\prod_{t}\tau_{\varpi_t}} & \mathscr{P}_{\G, \bS} \arrow{d}{\tau_{\varpi_{\bS}}} \\%
& \mathbb{G}_m
\end{tikzcd}.
\]
If $\tau_{\varpi_{\bS}}$ is regular, then $\tau_{\varpi_t}$ is regular for every ideal triangle $t$. Since $\bS$ is connected, we get $\overline{\omega}_s=\overline{\omega}_t:=\Lambda$ for every $s$ and $p$. Therefore $\tau_{\overline{\varpi}_\bS}=\Lambda(\tilde{\rho}_\bS)$. 
 \end{proof}

 \bp Every Casimir element $\tau$ of $\mathscr{P}_{\G, \bS}$ can be decomposed as 
 \[
 \tau=\lfloor \tau \rfloor \cdot \Lambda(\tilde{\rho}_{\bS}),
 \]
 where 
 \begin{itemize}
 \item $\lfloor \tau \rfloor \in \mu_{\bS}^\ast\mathcal{O}({\rm C}_{\G, \bS})$,
 \item $\Lambda \in X^\ast({\rm T})$. If $\bS$ has boundary intervals, then $\Lambda^\ast = \Lambda$.
 \end{itemize}
 \ep
 
 \begin{proof} By theorem \ref{TH1.16a}, part 1, we have $\tau= \tau_{\varpi_{\bS}}$. 
 By Lemma \eqref{regularity.lemma.center}, we have
 \[
 \tau_{\varpi_{\bS}}= \tau_{\lfloor \varpi_{\bS}\rfloor}\cdot \Lambda(\tilde{\rho}_{\bS}).
 \]
 
 If $\bS$ has only punctures, then $\tau_{\lfloor \varpi_{\bS}\rfloor}\in \mu_{\bS}^\ast\mathcal{O}({\rm C}_{\G, \bS})$ and $\Lambda \in X^\ast({\rm T})$. The proposition follows.
 
If $\bS$ has a boundary circle with $n$ many special points, then denote by $\omega_1, \ldots, \omega_n$ the restriction of $\varpi_{\bS}$ to these special points. We have 
 \[
 \omega_k=\omega_{k+1}^\ast, \qquad \forall k \in \mathbb{Z}/ n \mathbb{Z}. 
 \]
 Therefore $\lfloor \omega_k \rfloor = \lfloor \omega_{k+1} \rfloor^\ast$ and  we get   $\tau_{\lfloor \varpi_{\bS}\rfloor}\in \mu_{\bS}^\ast\mathcal{O}({\rm C}_{\G, \bS})$. Meanwhile
 \[
 \Lambda= \overline{\omega}_k= \overline{\omega_{k+1}^\ast}= \Lambda^\ast \in X^\ast({\rm T}),
 \]
 which concludes the proof of the proposition.
 \end{proof}


\begin{thebibliography}{Lus2}

\bibitem[AF91]{AF91} A.Yu. Alekseev, L.D. Faddeev:
\newblock  Commun. Math. Phys. 141 (1991) 413-422.

\bibitem[AGT]{AGT}  L. F. Alday, D. Gaiotto,   Y. Tachikawa:
\newblock  Liouville Correlation Functions from
Four-dimensional Gauge Theories, 
\newblock {\em Lett. Math. Phys.} 91 (2010), 167-197. 
\newblock \href{https://arxiv.org/abs/0906.3219}{arXiv:0906.3219}.

\bibitem[Ba]{Ba}  E.W. Barnes: 
\newblock  The genesis of the double gamma function.  
\newblock {\em Proc. London Math. Soc}. 31 (1899), 358-381. 

\bibitem[Bax]{Bax}  R. Baxter: 
\newblock {Exactly solved models in statistical mechanics. Integrable systems in statistical mechanics}. 
\newblock {\em Ser. Adv. Statist. Mech., }1, World Sci. Publishing, Singapore (1985).

\bibitem[BPZ]{BPZ} A. Belavin, A. Polyakov, A. Zamolodchikov:
\newblock {Infinite conformal symmetry in two-dimensional quantum field theory}.
\newblock {\em Nucl. Phys.} B241 (1984), no. 2, 333-380.


\bibitem[BK]{BK} J. Bernstein, B. Kr\"otz: 
\newblock Smooth Frechet globalizations of Harish-Chandra modules.
Israel J. Math, 199(1):45 - 111, 2014.

\bibitem[BZ]{BZ} 
A. Berenstein,  A. Zelevinsky: 
\newblock Quantum cluster algebras. 
\newblock {\em Adv. Math.} 195 (2005), n. 2, 405-455,
 \newblock \href{https://arxiv.org/abs/0404446}{arXiv:0404446}.

\bibitem[BFZ]{BFZ}
A. Berenstein, S. Fomin,  A. Zelevinsky:
\newblock Cluster algebras {III}: Upper bounds and double {B}ruhat cells.
\newblock {\em Duke Math. J.} 126 (2005), no. 1, 1--52.
\newblock \href{https://arxiv.org/abs/math/0305434}{arXiv:math/0305434}.


\bibitem[BT]{BT} A.G. Bytsko, J. Teschner:
\newblock R-operator, coproduct and Haar measure for the modular
double of $U_q(sl(2, R))$, 
\newblock {\em Comm. Math. Phys.}, 240 (2003), 171-196.
\newblock \href{https://arxiv.org/abs/math/0208191}{arXiv:math/0208191}.

\bibitem[C]{C} W. Casselman: \newblock Canonical extensions of Harish-Chandra modules to representations
of G. Canad. J. Math, 41(3):385 - 438, 1989.

\bibitem[CF]{CF}  L. Chekhov,  V. Fock: 
\newblock Quantum Teichm\"uller spaces.
\newblock \href{https://arxiv.org/abs/math/9908165}{arXiv:math/9908165}.


\bibitem[D]{D} V.~Drinfeld:  
\newblock Quantum groups.  
\newblock {\rm Proc. ICM-1986}, Vol. 1, 
 798--820, AMS, Providence, 1987.

 
 \bibitem[DO]{DO} H. Dorn, H.-J. Otto:
 \newblock Two- and three-point functions in {L}iouville theory.
 \newblock  {\em Nucl. Phys.} B429 (1994), 375-388.
 \newblock \href{https://arxiv.org/abs/hep-th/9403141}{arXiv:hep-th/9403141}.
 
 \bibitem[EL]{EL}  S. Evens, J.-H. Lu. 
\newblock Poisson geometry of the Grothendieck resolution of a complex
semisimple group. Mosc. Math. J., 7(4):613-642, 766, 2007.


   
 \bibitem[Fa1]{Fa1}  L.D.  Faddeev: 
\newblock Discrete Heisenberg-Weyl group and modular group.
\newblock  {\em Lett. Math. Phys}. 34 (1995), no. 3, 249--254.
 
 \bibitem[Fa2]{Fa2}  L.D.  Faddeev: 
\newblock Modular double of Quantum Group. 
\newblock \href{https://arxiv.org/abs/math/9912078}{arXiv:math/9912078}.


 
\bibitem[FK]{FK}  L.D.  Faddeev,  R.  Kashaev:  
\newblock Quantum dilogarithm. 
\newblock {\em Mod.Phys.Lett}. A 9 427 (1994).  

 
\bibitem[FKV]{FKV}  L.D. Faddeev,  R.M. Kashaev,  A.Yu. Volkov:
\newblock {Strongly coupled quantum discrete Liouville theory. I: Algebraic approach and duality}.
\newblock {\em Comm. Math. Phys.} 219 (2001) 199-219.
\newblock \href{https://arxiv.org/abs/hep-th/0006156}{arXiv:hep-th/0006156}


\bibitem[FRT]{FRT}  L.D. Faddeev,     N.Yu. Reshetikhin,  L.A. Takhtajan:
\newblock { Quantization of Lie groups and Lie algebras}.    
\newblock {\em St. Petersburg Math. J.} 1:1 (1990) 193-225.


\bibitem[F]{F} J. Fei.: \newblock {Tensor Product Multiplicities via Upper Cluster Algebras}.    
\newblock {\em Ann. Sci. \'Ec. Norm. Sup\'er}. 54, no. 6 (2021).
\newblock \href{https://arXiv:1603.02521}{arXiv:1603.02521} 



\bibitem[FF1]{FF1} B. Feigin, E. Frenkel: 
\newblock Quantization of the Drinfeld Sokolov reduction. 
\newblock {\em Phys. Lett.} B Vol 246, n 1-2,   75-81.


\bibitem[FF2]{FF2} B. Feigin, E. Frenkel: 
\newblock Integrals of Motion and Quantum Groups.
\newblock {\em Lecture Notes in Math.}, 1620, Fond. CIME/CIME Found. Subser., Springer, Berlin, 1996. 
\newblock \href{https://arxiv.org/abs/hep-th/9310022}{arXiv:hep-th/9310022}.   


\bibitem[FFu]{FFu} B. Feigin, D.B. Fuchs: 
\newblock   Verma modules over the Virasoro algebra. 
\newblock  {\em Functional Analysis and Its Applications}
July 1983, Volume 17, Issue 3, p. 241-242.
 
\bibitem[F94]{F94}  E. Frenkel:
\newblock Free field realizations in representation theory and conformal field theory. 
\newblock {\em Proc. of the ICM-94,} Vol. 1, 2 (Z\"urich, 1994), 1256-1269, Birkh\"auser, Basel, 1995.
\newblock \href{https://arxiv.org/abs/math/9408109}{arXiv:hep-th/ 9408109}.  


\bibitem[FBZ]{FBZ} E. Frenkel, D. Ben Zvi: Vertex Algebras and Algebraic Curves: Second Edition. Mathematical Surveys and Monographs
Volume: 88; 2004; 400 pp.


\bibitem[FG1]{FG03a}
V.V. Fock, A.B. Goncharov:
\newblock Moduli spaces of local systems and higher {Teichm\"{u}ller} theory.
\newblock {\em Publications Math\'{e}matiques   IHES,} No. 103 (2006),  1--211.
\newblock \href{https://arxiv.org/abs/math/0311149}{arXiv:math/0311149}.

\bibitem[FG2]{FG03b}
V.V. Fock,  A.B. Goncharov:
\newblock Cluster ensembles, quantization and the dilogarithm.
\newblock {\em Ann. Sci. \'{E}c. Norm. Sup\'{e}r.} (4) 42 (2009),  no.6,
  865--930.
\newblock \href{https://arxiv.org/abs/math/0311245}{arXiv:math/0311245}.

\bibitem[FG3]{FG05}
V.V. Fock, A.B. Goncharov:
\newblock Cluster {X}-varieties, amalgamation and {P}oisson-{L}ie groups.
\newblock {\em Algebraic Geometry and Number Theory,  
  Drinfeld's 50th birthday volume, Birkh{\"{a}}user}, (2006),  27--68.
\newblock \href{https://arxiv.org/abs/math/0508408}{arXiv:math/0508408}.

\bibitem[FG4]{FG07}
V.V. Fock,  A.B. Goncharov:
\newblock The quantum dilogarithm and representations of quantum cluster
  varieties.
\newblock {\em Invent. Math.} 175 (2009),  223--286.
\newblock \href{https://arxiv.org/abs/math/0702397}{arXiv:math/0702397}.

\bibitem[FG5]{FG14}
V.V. Fock, A.B. Goncharov:
\newblock Symplectic double for moduli spaces of {G}-local systems on surfaces.
\newblock {\em Adv. Math.} 300 (2016), 505--543.
\newblock \href{https://arxiv.org/abs/math/1410.3526}{arXiv:1410.3526}.

\bibitem[FG6]{FG05b}
V.V. Fock, A.B. Goncharov:
\newblock Dual Teichmuller and lamination spaces.
\newblock {\em  Handbook on Teichmuller theory}. IRMA Lect. Math. Theor. Phys., 11
European Mathematical Society (EMS), Zurich, 2007, 647–684.
\newblock \href{https://arxiv.org/abs/math/0510312 }{arXiv:0510312}.

\bibitem[FR]{FR} V.V. Fock, A.B. Rosly: 
\newblock Flat connections and polyubles. 
\newblock {\em Theor.and Math. Physics} 95 (2), 526-534.

\bibitem[FZ98]{FZ98} S.~Fomin,  A.~Zelevinsky:
\newblock  Double Bruhat cells and total positivity.
 \newblock {\em Journal AMS.} 12  (1999), no. 2, 335-380;  
\newblock \href{https://arxiv.org/abs/math/9802056}{arXiv:math/9802056}.

\bibitem[FZ]{FZI}
S.~Fomin, A.~Zelevinsky:
\newblock Cluster algebras {I}: Foundations.
\newblock {\em Journal AMS.} 15 (2002),  no. 2, 497--529.
\newblock \href{https://arxiv.org/abs/math/0104151}{arXiv:math/0104151}.

\bibitem[FZIV]{FZIV}
S.~Fomin,  A.~Zelevinsky:
\newblock Cluster algebras {IV}: Coefficients.
\newblock {\em Compos. Math.} 143 (2007), no. 1,  112--164.
\newblock \href{https://arxiv.org/abs/math/0602259}{arXiv:math/0602259}.


\bibitem[FI]{FI} I. Frenkel, I. Ip: 
\newblock Positive representations of split real quantum groups and future perspectives.
\newblock {\em Int. Math. Res. Not. IMRN} 2014, no. 8, 2126--2164.
\newblock \href{https://arxiv.org/abs/1111.1033}{arXiv:1111.1033}.

\bibitem[G]{G} D. Gaitsgory: 
\newblock  
\newblock {\em Quantum Langlands Correspondence}  
\newblock \href{https://arXiv:1601.05279}{arXiv:1601.05279}.
\bibitem[G]{G} D. Gaiotto, 
\newblock ${\cal N} = 2$ dualities.  
\newblock \href{https://arxiv.org/abs/0904.2715}{arXiv:0904.2715}.


\bibitem[GG73]{GG73}  S.I. ~Gelfand,  M.I.~Graev: 
\noindent Fourier-Weyl operators on the principal affine space of a Chevalley group.
\noindent Preprint N79,   Institute of Applied Mathematics, Moscow,  1973. 
 
 \bibitem[GK]{GK} Goncharov A., Kontsevich M.: 
 \noindent Spectral description of non-commutative local systems on surfaces and non-commutative cluster varieties.  
 \newblock \href{https://arxiv.org/abs/2108.04168 }{arXiv:2108.04168}.  
 
 \bibitem[GY]{GY} Goodearl, M. Yakimov: 
\noindent The Berenstein-Zelevinsky quantum cluster algebra conjecture.
\newblock {\em J. Eur. Math. Soc.} 22(2020), no.8, 2453-2509,
 \newblock \href{https://arxiv.org/abs/1602.00498}{arXiv:1602.00498}. 
 
 \bibitem[GKL]{GKL} A. Gerasimov, S. Kharchev, D. Lebedev.:
 \newblock Representation theory and quantum integrability.
 \newblock {\em Infinite dimensional algebras and quantum integrable systems}, 133-156, 
Progr. Math., 237, Birkh\"auser, Basel, 2005.  
 \newblock \href{https://arxiv.org/abs/0402112}{arXiv:math/0402112}.
 
  \bibitem[GKL2]{GKL2}  A. Gerasimov, S. Kharchev, D. Lebedev, S. Oblezin:
  \newblock  On a Class of Representations of Quantum Groups, 
   \newblock \href{https://arxiv.org/abs/math/0501473}{arXiv:math/0501473}.
 
  
 
 \bibitem[GS1]{GS13}
A.B. Goncharov, L.~Shen:
\newblock Geometry of canonical bases and mirror symmetry.
\newblock {\em Invent. Math.} 202 (2015), no. 2, 487--633.
\newblock \href{https://arxiv.org/pdf/1309.5922.pdf}{arXiv:1309.5922}.

\bibitem[GS2]{GS16}
A.B. Goncharov, L.~Shen:
\newblock {D}onaldson-{T}homas transformations for moduli spaces of {G}-local
  systems.
\newblock {\em Adv. Math.} 327 (2018),  225--348.
\newblock \href{https://arxiv.org/abs/1602.06479}{arXiv:1602.06479}.
 

\bibitem[GHKK]{GHKK} M.~Gross, P.~Hacking, S.~Keel, M.~Kontsevich:  
\newblock  Canonical bases for cluster algebras.  
\newblock {\em J. Amer. Math. Soc.} 31 (2018), no. 2, 497--608.
\newblock \href{https://arxiv.org/abs/1411.1394}{arXiv:1411.1394.}

\bibitem[H]{H} N.~Hitchin: 
\newblock  Lie groups and Teichmuller space. 
\newblock Topology 31 (1992), no. 3, 449--473.

\bibitem[JN]{JN}  S. Jeong, N. Nekrasov: 
\newblock  Opers, surface defects, and Yang-Yang functional.
\newblock {\em Adv. Theor. Math. Phys.} 24(2020) 7, 1789-1916,
 \newblock  \href{https://arxiv.org/abs/1806.08270v3}{arXiv:1806.08270} 
  
    
\bibitem[J]{J} M. Jimbo: 
\noindent A q-difference analogue of U(g) and the Yang-Baxter equation, {\em Lett. Math. Phys.}, 10 (1985), 63 - 69.

  \bibitem[KRR]{KRR}  V. Kac,  A.Raina,  N. Rozhkovskaya: 
\newblock {\it Bombay lectures on highest weight representations of infinite dimensional Lie algebras}.
\newblock   2013.

  \bibitem[K97]{K97} M. Kapranov: 
\newblock   Heisenberg doubles and derived categories. 
 \newblock  \href{https://arxiv.org/abs/9701009}{arXiv:9701009} 
 
   
   \bibitem[Kas95]{Kas95}   R. Kashaev: 
\newblock     Heisenberg double and the pentagon relation, 
  \newblock  \href{https://arxiv.org/abs/9503005}{arXiv:9503005}    
 
 
  \bibitem[Ke]{Ke}   S. Kerckhoff, 
  \newblock The Nielsen realization problem, Ann. of Math. (2) 117 (1983), no. 2, 235–265.    
    
  
  \bibitem[I1]{II1} I. Ip: 
\newblock     
Positive representations of split real simply-laced quantum groups.
 \newblock  \href{https://arxiv.org/abs/1203.2018 }{arXiv:1203.2018}.  
 \newblock   Positive representations of non-simply-laced split real quantum groups.
  \newblock  \href{https://arxiv.org/abs/1205.2940}{arXiv:1205.2940}. 
  
  
\bibitem[I2]{II2} I. Ip: 
\newblock Cluster realization of $U_q(\mathfrak{g})$ and factorization of  universal $R$-matrix.
 \newblock  \href{https://arxiv.org/abs/1612.05641}{arXiv:1612.05641}.  


\bibitem[IIO]{IIO} R. Inoue, T. Ishibashi, H. Oya:  
\newblock Cluster realizations of Weyl groups and higher Teichm\"uller theory. 
\newblock {\em Sel. Math. New Ser.} 27, 37(2021),
 \newblock  \href{https://arxiv.org/abs/1902.02716}{arXiv:1902.02716}.  
 
\bibitem[K1]{K1}  R. Kashaev: 
\newblock Quantization of Teichm\"uller spaces and quantum dilogarithm. 
\newblock {\em Letters Math. Phys.} 43 (1998) 2, 105-115.

\bibitem[KL]{KL} D. Kazhdan, G. Lusztig:  
\newblock Tensor structures arising from affine Lie algebras, I, 
JAMS 6 (1993), 905 --947, II, JAMS
6 (1993), 949 - 1011, IV JAMS 7 (1994), 383-453.

\bibitem[K13]{K13} B. Keller: 
\newblock Quiver mutation and combinatorial DT-invariants. 
\newblock {\em Contribution to the FPSAC} 2013. 
\newblock  http://webusers.imj-prg.fr/~bernhard.keller/publ/index.html. 

\bibitem[KLS]{KLS} S. Kharchev, D. Lebedev, M. Semenov-Tian-Shansky:  
\newblock Unitary representations of $U_{q}({sl}(2,\R))$, the modular double, and  multiparticle q-deformed Toda chains.  
\newblock {\em Comm. Math. Phys.} 225 (2002), no. 3, 573--609.
\newblock \href{https://arxiv.org/abs/0102180}{arXiv:0102180}. 

\bibitem[KR]{KR}  A.N. Kirillov, N. Reshetikhin: 
\newblock q-Weyl group and a multiplicative formula for universal R-matrices. \newblock {\em Comm. Math. Phys.} 134 (1990), no. 2, 421 - 431.

\bibitem[KS]{KS} Kontsevich M., Soibelman Y.:  
 \newblock    Stability structures, motivic Donaldson-Thomas invariants and cluster transformations. 
  \href{https://arxiv.org/abs/0811.2435}{arXiv:0811.2435}.   
    

\bibitem[Ku]{Ku}  S. Kumar: 
\newblock Proof of the Parthasarathy-Ranga Rao-Varadarajan conjecture.
\newblock  {\em Inv. Math.} 93, 117-130(1988). 
\href{http://www.unc.edu/math/Faculty/kumar/papers/kumar14.pdf}{Kumar's webpage}.


\bibitem[Le1]{IL1}
I.~Le:
\newblock Cluster structures on higher {T}eichm{\"u}ller spaces for classical
  groups.
  \newblock {\em Forum of Mathematics, Sigma.} 2019;7:e13
 \href{https://arxiv.org/abs/1603.03523}{ arXiv:1603.03523}.

  

 

\bibitem[Le2]{IL2}
I.~Le:
\newblock  An Approach to Cluster Structures on Moduli of Local Systems for General Groups. 
\newblock {\em International Mathematics Research Notices,} Volume 2019, Issue 16, August 2019, Pages 4899?4949,
\newblock \href{https://arxiv.org/abs/1606.00961}{arXiv:1606.00961}.




\bibitem[LS]{LS}  S. Levendorskii, Y. ~Soibelman:
\newblock Some applications of quantum Weyl groups. 
{\em J. Geom. and Phys.} 7 (1990), 241-254.


\bibitem[LS1]{LS1} S. Levendorskii,  Ya. Soibelman: 
\newblock Algebras of functions on compact quantum groups, Schubert cells and quantum tori, {\it Comm. Math. Phys.} 139, 141 - 170 (1991).


\bibitem[LV]{LV}
G. Lion, M. Vergne: 
\newblock {\em The Weil representation, Maslov index and Theta series}.
\newblock Progress in Mathematics, Vol. 6 (1980). .


\bibitem[Lus1]{L}
G.~Lusztig:
\newblock {\em Introduction to {Q}uantum {G}roups}.
\newblock Birkh{\"a}user, Boston, MA, 1993.

\bibitem[Lus2]{L1}
G.~Lusztig:
\newblock Total positivity in reductive groups.
\newblock {\em Progr. Math., 123},   531--568, 1994.


\bibitem[Lus3]{L4}
G.~Lusztig:
\newblock Total positivity and canonical bases.
\newblock {\em Algebraic groups and Lie groups. Austral. Math. Soc. Lect. Ser}. vol. 9,   281-295.  Cambridge Univ. Press, Cambridge (1997). 


\bibitem[Lus4]{L6} G.~Lusztig: 
\newblock  Quantum deformations of certain simple modules over enveloping algebras. 
{\em Adv in Math.} 70 (1988), 237-249.

\bibitem[Lus5]{L7} G.~Lusztig: 
\newblock Quantum groups at roots of 1. {\em Geom. Dedicata} 35 (1990), 89-114. 

\bibitem[Mir]{Mir} Mirzakhani M.: 
\newblock Ergodic theory of the earthquake flow. {\it IMRS}, Vol 2008, 39 pages. 

\bibitem[N]{N}   N.  Nekrasov: 
\newblock  Seiberg-Witten Prepotential From Instanton Counting, 
 {\em Adv.Theor. Math. Phys}. 7 (2004) 831--864, 
\newblock \href{https://arxiv.org/abs/math/0206161}{arXiv:hep-th/0206161}.

\bibitem[NS]{NS} N.~Nekrasov, S.~Shatashvili:
 \newblock Quantization of Integrable Systems and Four Dimensional Gauge Theories.
 \newblock {\em XVIth Int. Congress on Math. Phys.}, 265--289, World Sci. Publ., 
 \newblock \href{https://arxiv.org/abs/math/0908.4052}{arXiv:0908.4052}. 
 

\bibitem[Nov92]{Nov92} S.P. Novikov: 
\newblock Various doublings of Hopf algebras. Operator algebras on quantum groups, complex cobordisms. 
{\em Russian Mathematical Surveys} (1992),47(5):198.  
   
\bibitem[Pa]{Pa} Papadopoulos A.: 
\newblock Geometric intersection functions and hamiltonian flows on the space of measured foliations of a surface.
{\it Topology and its Applications} 41 (1991):147-177.

\bibitem[PT1]{PT1}  B. Ponsot, J. Teschner:  
\newblock Liouville bootstrap via harmonic analysis on a noncompact quantum group. 
\newblock \href{https://arxiv.org/abs/math/9911110}{arXiv:hep-th/9911110}. 

\bibitem[PT2]{PT2} B. Ponsot, J. Teschner:  
\newblock Clebsch-Gordan and Racah-Wigner coefficients for a continuous series of representations of $U_q(sl(2, R))$, 
\newblock {\em Comm. Math. Phys.}, 224, (2001), p. 613-655.

\bibitem[RS90]{RS90} N.Yu. Reshetikhin, M.A. Semenov-Tian-Shansky: 
\newblock  Lett. Math. Phys. 19
(1990)

\bibitem[Se95]{Se95}  N. Seiberg: 
\newblock Electric-magnetic duality in supersymmetric non-abelian gauge theories. 
\newblock{\em Nuclear Phys}. B 435(1-2):129-146, 1995. 
\newblock \href{https://arxiv.org/abs/math/9411149}{arXiv:hep-th/9411149}. 

\bibitem[SS1]{SS1}  G. Schrader, A. Shapiro:   
\newblock A cluster realization of $U_q(\mathfrak{sl_n})$ from quantum character varieties.
\newblock {\em Invent. math.} (2019)
\newblock  \href{https://arxiv.org/abs/1607.00271}{arXiv:1607.00271}.   

\bibitem[SS2]{SS2}  G. Schrader, A. Shapiro:   
\newblock   Continuous tensor categories from quantum groups I. 
\newblock  \href{https://arxiv.org/abs/1708.08107}{arXiv:1708.08107}.  
 
\bibitem[SS3]{SS3}  G. Schrader, A. Shapiro:     
\newblock On $b$-Whittaker functions.
\newblock  \href{https://arxiv.org/abs/1806.00747}{arXiv:1806.00747}. 

\bibitem[SS4]{SS4}  G. Schrader, A. Shapiro:   
\newblock To appear
\newblock  

\bibitem[STS]{STS} M. Semenov-Tian-Shansky: 
\newblock Dressing transformations and Poisson Lie groups actions. {\em Publ. Res. Inst. Math. Sci.} 21(6): 1237-1260. 1985.


\bibitem[STS92]{STS92} M. Semenov-Tian-Shansky: 
\newblock Poisson-Lie groups. The quantum duality principle and the twisted quantum double. 
{\em Theor and Math Physics} vol. 93, pp.1292-1307 (1992)

\bibitem[SW]{SW} L. Shen, D. Weng:
\newblock Cluster structures on Double Bott-Samelson cells,
{\em Forum of Mathematics, Sigma}, 9 E66
\newblock \href{https://arxiv.org/abs/1904.07992}{arXiv:1904.07992}  

\bibitem[S20]{S20} L. Shen:   
\newblock  Duals of semisimple Poisson-Lie groups and cluster theory of moduli spaces of G-local systems.
\newblock  {\em International Mathematics Research Notices,} Volume 2022, Issue 18, September 2022, Pages 14295-14318,
\newblock \href{https://arxiv.org/abs/arXiv: 2003.07901}{arXiv: 2003.07901}. 


\bibitem[S22]{S22} L. Shen:   
\newblock  Cluster nature of quantum groups.  
\newblock \href{https://arxiv.org/abs/2209.06258}{arXiv: 2209.06258}


\bibitem[Sh]{Sh}  T. Shintani: 
 \newblock On a Kronecker limit formula for real quadratic fields. 
{\em J. Fac. Sci. Univ. Tokyo Sect} 1A Math 24 (1977) 167-199. 

\bibitem[Sch]{Sch} K. Schm\"udgen: 
\newblock Operator representations of $U_q(sl(2, R))$, 
\newblock {\em Lett. Math. Phys.} 37 (1996) 211-222.

\bibitem[So]{So} Y.~ Soibelman: 
\newblock Algebra of functions on a compact quantum group and its representations. 
{\em Algebra and Analysis}, 2, (1990), 190-212.


\bibitem[So1]{So1} Y. Soibelman: Quantum Weyl group and some of its applications. 
{\it Rend. Circ. Mat. Palermo},    volume Suppl. No 26, no. 2, (1991), p.155 - 181.


\bibitem[T]{T}  J. Teschner: 
\newblock An analog of a modular functor from quantized Teichm\"uller theory. 
 {\em Handbook of Teichm\"uller theory. Volume 1. IRMA Lect.   Math.  and Theor.  Phys.} 11. 2007. 685-760. 
\newblock \href{https://arxiv.org/abs/0510174}{arXiv:QA/0510174}. 

\bibitem[T10]{T10} J. Teschner:
\newblock Quantization of the Hitchin moduli spaces, Liouville theory, and the geometric
Langlands correspondence I, 
\newblock {\em Adv. Theor. Math. Phys.} 15 (2011) 471-564. 
 \href{https://arxiv.org/abs/1005.2846}{arXiv:1005.2846}.

\bibitem[TV]{VT} J. Teschner, G. Vartanov:
\newblock Supersymmetric gauge theories, quantization of {$\mathcal{M}_{\rm flat}$}, and conformal field theory, 
 \newblock {\em Adv. Theor. Math. Phys.} 19 (2015) 1-135.
\newblock \href{https://arxiv.org/abs/1302.3778}{arXiv:1302.3778}  

\bibitem[W]{W} D.~Weng: 
\newblock Donaldson-Thomas Transformation of Double Bruhat Cells in Semisimple Lie Groups. 
\newblock {\em Ann. Sci. \'Ec. Norm. Sup\'er}. 53:291-352, 2020.
\newblock \href{https://arxiv.org/abs/1611.04186}{arXiv:1611.04186}.

 \bibitem[ZZ]{ZZ}  A. B. Zamolodchikov, Al. B. Zamolodchikov: 
  \newblock Structure Constants and Conformal Bootstrap in Liouville Field Theory. 
  \newblock {\em Nucl. Phys. B} 477 (1996), 577-605.
 \newblock \href{https://arxiv.org/abs/hep-th/9506136}{arXiv:hep-th/9506136}. 




\end{thebibliography}
\end{document}